\documentclass[11pt,reqno]{amsart} 
\usepackage[margin=1in]{geometry}  
\usepackage[utf8]{inputenc}  
\usepackage{amsmath,amsthm,amsfonts,amssymb}
\usepackage{times}
\usepackage{enumitem}
\usepackage{mathrsfs}
\usepackage{comment}
\usepackage{graphicx}
\usepackage[font=footnotesize]{caption}
\usepackage{tikz}
\usepackage{amsmath,amssymb}
\usetikzlibrary{arrows.meta}

\def\p{\partial}
\renewcommand{\epsilon}{\varepsilon}
\def\eps{\varepsilon}
\def\les{\lesssim}
\renewcommand*{\div}{\ensuremath{\mathrm{div}}\,}
\newcommand{\brak}[1]{\langle #1 \rangle} 
\newcommand{\norm}[1]{\left \| #1 \right\|} 
 
\newcommand{\abs}[1]{\left|#1\right|}

\newcommand*{\supp}{\ensuremath{\mathrm{supp}}\,}
\newcommand*{\sgn}{\ensuremath{\mathrm{sgn}}\,}
\newcommand*{\Id}{\ensuremath{\mathrm{Id}}\,}
\newcommand*{\tr}{\ensuremath{\mathrm{Tr}}\,}
\def\Re{\ensuremath{\mathrm{Re}}\,}

\newcommand{\Reals}{\mathbb{R}}
\newcommand{\Naturals}{\mathbb{N}}

\newcommand{\Compl}{\mathbb{C}}

\newcommand{\cV}{\mathcal{V}}

\newcommand{\OO}{\mathcal{O}}

\newcommand{\DD}{\mathcal{D}}

\newcommand{\YY}{\mathcal{Y}}
\newcommand{\TT}{\mathcal{T}}

\newcommand{\NN}{\mathcal{N}}

\newcommand{\gamkin}{\gamma_{\mathsf{kin}}}

\newcommand{\LLL}{\mathsf{L}}
\newcommand{\MMM}{\mathsf{M}}

\newcommand{\NNN}{\mathsf{N}}

\newcommand{\III}{\mathsf{I}}
\newcommand{\JJJ}{\mathsf{J}}

\newcommand{\cx}{\mathsf{c_r}}
\newcommand{\cxstar}{\mathsf{c_r}^{\!\!*}}
\newcommand{\cxbar}{\mathsf{\bar c_r}}
\newcommand{\cxtilde}{\mathsf{\tilde{c}_r}}
\newcommand{\cb}{\mathsf{c_b}}
\newcommand{\cbstar}{\mathsf{c_b}^{\!\!*}}
\newcommand{\cbbar}{ \mathsf{\bar c_b}} 
\newcommand{\cbtilde}{ \mathsf{\tilde{c}_b}} 
\newcommand{\cu}{\mathsf{c_u}}

\newcommand{\cubar}{ \mathsf{\bar c_u}} 
\newcommand{\cutilde}{ \mathsf{\tilde{c}_u}}

\newcommand{\crho}{\mathsf{c_\varrho}}
\newcommand{\crhobar}{ \mathsf{\bar c_\varrho}} 
\newcommand{\crhotilde}{ \mathsf{\tilde{c}_\varrho}} 

\newcommand{\cbb}{\mathsf{c_B}}
\newcommand{\cbbbar}{ \mathsf{\bar c_B}} 
\newcommand{\cbbtilde}{ \mathsf{\tilde{c}_B}}

\newcommand{\Wring}{\mathring{W}}
\newcommand{\Wringbar}{\overline{\mathring{W}}}
\newcommand{\Wringtilde}{\tilde{\mathring{W}} }
\newcommand{\Zring}{\mathring{Z}}
\newcommand{\Zringbar}{\overline{\mathring{Z}}}
\newcommand{\Zringtilde}{\tilde{\mathring{Z}}}
\newcommand{\Aring}{\mathring{A}}
\newcommand{\Aringbar}{\overline{\mathring{A}}}
\newcommand{\Aringtilde}{\tilde{\mathring{A}}}

\newcommand{\Ybar}{\overline{\YY}}
\newcommand{\Ytilde}{\tilde{\YY}}
\newcommand{\Tbar}{\overline{\TT}}
\newcommand{\Ttilde}{\tilde{\TT}}
\newcommand{\Dbar}{\overline{\DD}}
\newcommand{\Dtilde}{\tilde{\DD}}

\newcommand{\bfa}{{\boldsymbol{\alpha}}}

\newcommand{\bfu}{{\boldsymbol{u}}}
\newcommand{\bfU}{{\boldsymbol{U}}}

\newtheorem{theorem}{Theorem}[section]
\newtheorem{lemma}[theorem]{Lemma}
\newtheorem{proposition}[theorem]{Proposition}
\newtheorem{corollary}[theorem]{Corollary}
\newtheorem{definition}[theorem]{Definition}

\newtheorem{example}[theorem]{Example}

\theoremstyle{definition}
\newtheorem{remark}[theorem]{Remark}

\numberwithin{equation}{section}

\usepackage{xcolor}

\usepackage[colorlinks=true, pdfstartview=FitV, linkcolor=blue,citecolor=blue, urlcolor=blue]{hyperref}

\definecolor{labelkey}{rgb}{0,0,1}
\definecolor{blue2}{cmyk}{.94,.11,0,0}

\let\pa=\partial
\let\al=\alpha
\let\bet=\beta
\let\del=\delta
\let\eps=\varepsilon
\let \kp = \kappa
\let\lam=\lambda
\let\sig=\sigma
\let \les = \lesssim
\let \gtr = \gtrsim
\let \th = \theta
\let \pr = \prime
\let \vp = \varphi
\let\Gam=\Gamma
\let\Del=\Delta
\let\Lam=\Lambda
\let\Om=\Omega
\let\td = \tilde

\let\pa=\partial
\def\cA{{\mathcal A}}
\def\cB{{\mathcal B}}
\def\cC{{\mathcal C}}
\def\cD{{\mathcal D}}
\def\cE{{\mathcal E}}

\def\cI{{\mathcal I}}
\def\cJ{{\mathcal J}}
\def\cL{{\mathcal L}}
\def\cN{{\mathcal N}}
\def\cR{{\mathcal R}}
\def\cS{{\mathcal S}}
\def\cT{{\mathcal T}}
\def\cX{{\mathcal X}}
\def\cW{{\mathcal W}}
\def\na{\nabla}
\def\la{\langle}
\def\ra{\rangle}
\def\one{\mathbf{1}}
\def\udb{\underbrace}

\newcommand{\ee}{\mathbf{e}}
\newcommand{\uu}{\bfu} %\newcommand{\uu}{\mathbf{u}}
\newcommand{\vv}{{\boldsymbol{v}}} %\newcommand{\vv}{\mathbf{v}}
\newcommand{\UU}{\bfU} %\newcommand{\UU}{\mathbf{U}}
\newcommand{\WW}{{\boldsymbol{W}}} % \newcommand{\WW}{\mathbf{W}}

\newcommand\ang[1]{ { \la {#1}\ra } } 

\let \mr = \mathring
\let \mw = \mathrm
\let \msf = \mathsf
\let \tf = \tfrac
\let \itl = \intercal
\def\HH{\mathsf{H}}
\def\HHH {\mathsf{H}^{(3)}}
\def\HHT {\mathsf{H}^{(2)}}
\def\VVs {\mathsf{V}}

\def\LLs{\mathsf{L}}

\def \ccs {\mathsf{c}}
\def \aaa { \bar \vrho_0 }
\def \cca {  \bar{\mathsf{C}}_1 }
\def \aan { \bar \vrho_{ 2\NNN} }
\def \bb {\bar B_2}
\def \bbn {\bar B_{2\NNN + 2}}
\def \uua{\bar U_1}
\def \uun {\bar U_{2\NNN + 1}}
\def \iin {\mathsf{in}}
\def \bu {\bar \UU }

\def \xu { \td Z }
\def \divu { \na\cdot \tu }
\def \bcr {\cxbar}
\def \bcu {\cubar}

\def \mcb {\cbb}
\def \tcb {\cbbtilde}
\def \bcb {\cbbbar}
\def \brho {\bar \rho}
\def \cvr {\crho}
\def \bcvr {\crhobar}
\def \tcvr {\crhotilde}
\def \cc {\mathsf C}
\def \bc {\bar \cc}
\def \tu {\td \UU }

\def \tcr {\cxtilde}

\def \tcu {\cutilde}
\def \tw { \td \WW }
\def \bw {\bar \WW}

\def \twm {\td \WW_{\mm}}
\def \tum {\td \UU_{\mm} }
\def \tbm {\tb_{\mm} }
\def \tvrm {\tvr_{\mm}}
\def \rs { R(\tau) }
\def \vrs {\Bar{\Bar{\vrho}}}
\def \vrl {\Bar{\Bar{\vrho_{\mathsf{l}}}}}
\def \vru {\Bar{\Bar{\vrho_{\mathsf{u}}}}}
\def \bbl { \mb_{ \mathsf{l} } }
\def \bbu { \mb_{ \mathsf{u} } }
\def \ccl {\cc_{ \mathsf{l} } }
\def \ccu {\cc_{ \mathsf{u} } }
\def \vrho {\varrho}
\def \bvr {\bar \vrho}
\def \tvr{\td \vrho}
\def \tl { \mathsf{I}}
\def \rl {\mathsf{P}}
\def \EE {\mathscr{E}}
\def \EEt {\mathscr{E}_{\msf{tot}}}
\def \FFs {\mathscr{F}}
\def \GGs {\mathscr{G}}
\def \NNs {\mathscr{N}}
\def \tlk {\mathsf{I}_k}
\def \rlk {\mathsf{P}_k}

\def \rlm {\mathsf{P}_{\mathsf{M}}}

\def \mm {\mathsf{M}}
\def \ddd { \del_{\mathsf{0} }}
\def \dda { \del_{\mathsf{1} }}
\def \ddb {\del_{\mb}}
\def \mb {  \mathsf{B} }
\def \barb {\bar {\mathsf{B}}}
\def \tb { \td{\mathsf{B}}}
\def \mbb {  \mathsf{b} }
\def \vvrho{\eta}

\def \kk { {k_*}}

\def \Algm {\mathsf{AM}}
\def \Algmpos {\mathsf{AM}_{\geq 0} }

\newcommand\Std{ \tilde{\mathbb{S}}}
\newcommand\Htd{ \tilde{\mathbb{H}}}
\newcommand\Sringbb{\mathring{\mathbb{S}}}
\newcommand\Hringbb{\mathring{\mathbb{H}}}
\newcommand\Aringbb{\mathbb{A}}

 \newcommand{\bseq}{\begin{subequations}}
 \newcommand{\eseq}{\end{subequations}}

\newcommand{\beq}{\begin{equation}}
\newcommand{\eeq}{\end{equation}}
 \newcommand{\bal}{\begin{aligned} }
 \newcommand{\eal}{\end{aligned}}

\newcommand{\PP}{\mathbf{P}}
\newcommand{\bmm}{\mathbf{m}}

\newcommand{\MMn}{\mathbf{M}_n}
\newcommand{ \VVN}{V_{ \Del^{\NNN } } }
\newcommand{ \VVNother}{ V_{ 2 \NNN \backslash \Del^{\NNN } } }
\newcommand\OOL[1]{\mathscr{L}_{{#1}}}

\title{Smooth and stable Euler implosions}

\author{Jiajie Chen}
\address{Department of Mathematics, University of Chicago, Chicago, IL 60637.}
\email{\href{jiajiechen@uchicago.edu}{jiajiechen@uchicago.edu}}
\thanks{The work of JC was partially supported by NSF Grant
DMS--2408098.}

\author{Steve  Shkoller}
\address{Department of Mathematics, University of California Davis, Davis, CA 95616.}
\email{\href{shkoller@math.ucdavis.edu}{shkoller@math.ucdavis.edu}}
\thanks{The work of SS was partially supported by the Collaborative NSF grant DMS--2307680.}

\author{Vlad Vicol}
\address{Courant Institute of Mathematical Sciences, New York University, New York, NY 10012.}
\email{\href{vicol@cims.nyu.edu}{vicol@cims.nyu.edu}}
\thanks{The work of VV was partially supported by the Collaborative NSF grant DMS--2307681 and a Simons Investigator Award.}
 
\begin{document}

\begin{abstract}
We construct a new class of self-similar implosion profiles for the multi-dimensional compressible Euler equations. These profiles are \emph{smooth}, \emph{genuinely non-isentropic}, radially/spherically symmetric, and have \emph{explicit (closed-form) similarity exponents}. We prove that the exact Euler solution corresponding to the ground state implosion profile is \emph{stable} to radially symmetric perturbations,  as a solution to the full nonlinear compressible Euler equations, modulo a one-dimensional compatibility condition on the initial data. For perturbations of the Euler solution corresponding to the ground state implosion profile of a monatomic or diatomic gas, that do not obey any symmetry assumptions, we provide a complete characterization of the set of initial data that yield nonlinear stability.
\end{abstract}

\maketitle

\setcounter{tocdepth}{1}
\hypersetup{bookmarksdepth=2}
\tableofcontents

 %%%%%%%%%%%%% START OF INTRO %%%%%%%%%%%%%

\section{Introduction}
The compressible Euler equations admit a remarkable class of singularities known as \emph{implosions}, in which at least one of the primary flow variables (typically density and/or pressure) becomes unbounded at a single point in spacetime. This is in contrast to shock singularities, where the primary flow variables remain bounded but their gradients become unbounded. Understanding the \emph{existence} and \emph{stability} of such singularities is a fundamental problem in fluid/gas dynamics.

\begin{figure}[htb!]
\centering
\includegraphics[width=0.55\textwidth]{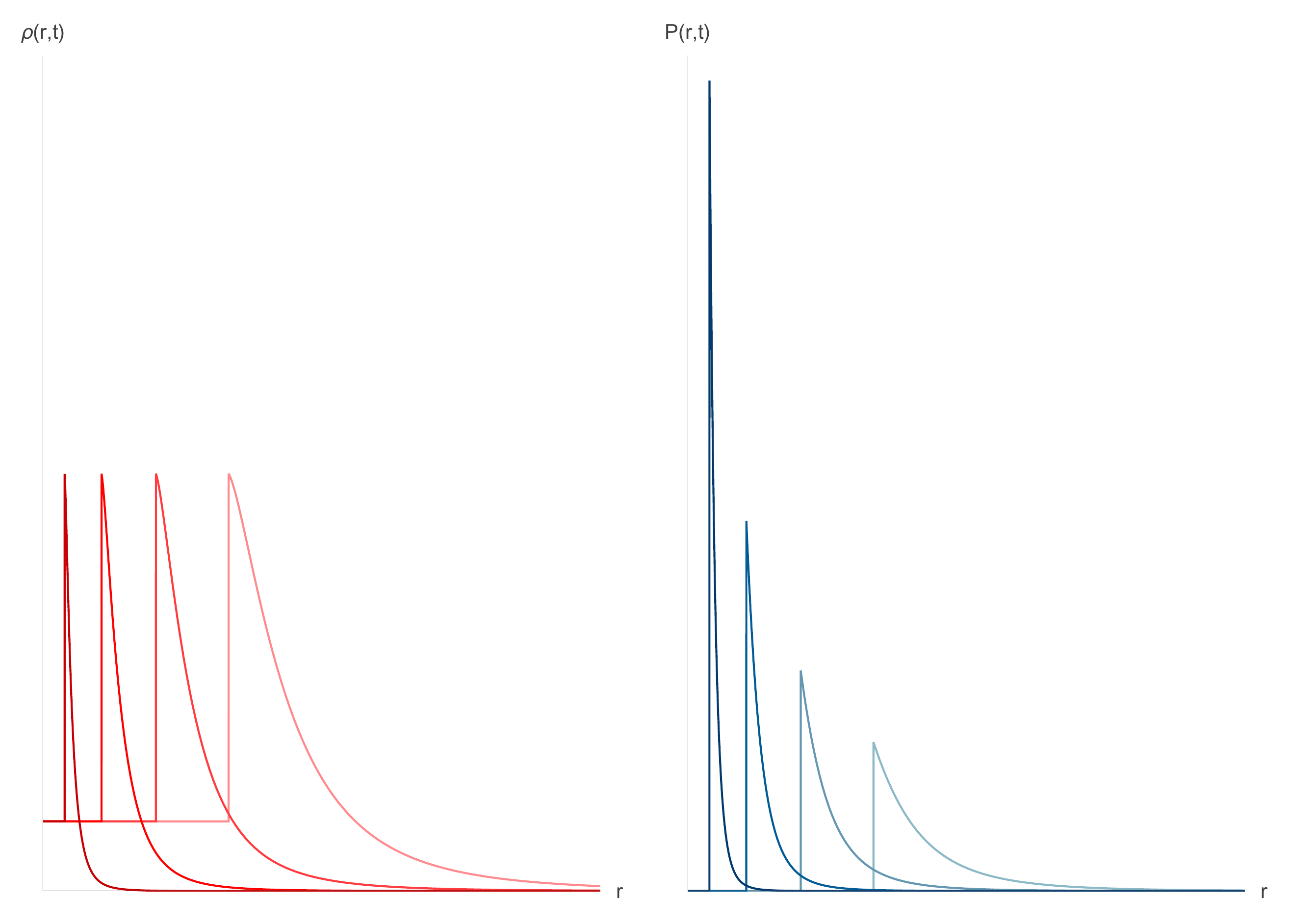}
\caption{Snapshots of the radial density $\rho(r,t)$ (red, left) and pressure $p(r,t)$ (blue, right) for the Guderley implosion (for $d=3$, $\gamma = 7/5$, and similarity exponent reported by Lazarus~\cite{Lazarus1981}), at four successive times, with darker shades indicating later times. Ahead of the imploding shock the density is a positive constant and the pressure vanishes identically. As $t \to t_*^-$, the time of implosion, the density remains bounded while the pressure diverges near the origin.}
\label{fig:Guderley}
\end{figure}

Prior to this work, implosions lying in two \emph{extreme} regularity classes were proven to exist for the compressible Euler equations: \emph{imploding shocks} and \emph{smooth isentropic implosions}. The main features and drawbacks of these two classes of solutions may be briefly summarized as follows (a full discussion and detailed bibliographic account is given in Section~\ref{sec:literature} below):
\begin{itemize}[leftmargin=1em]

\item \textbf{Imploding shocks}. 
The classical implosion scenario is given by Guderley's 1942 self-similar solution~\cite{Guderley1942}, re-discovered by Landau and Stanyukovich~\cite{Stanyukovich1960} in 1944--1945. It describes a radial/spherical shock wave converging inward into a quiescent medium (constant density, zero velocity, zero pressure), leading to a finite time collapse at the spatial origin. As the shock focuses, both the pressure and the radial velocity diverge at the moment of collapse, but the density remains finite; see Figure~\ref{fig:Guderley}. Guderley's solution applies to the full (non-isentropic) Euler system. 
Linear mode analysis suggests stability of the Guderley solution under radially symmetric perturbations~\cite{Morawetz1951,ChenZhangPanarella1995}, while linear instability to non-radial perturbations has been established in~\cite{BrushlinskiiKazhdan1963,WuRoberts1996PRE}. A significant limitation of Guderley's solution is that the similarity exponent can only be determined \emph{numerically}, to a specified precision, by solving a nonlinear eigenvalue problem at the sonic point~\cite{Lazarus1981, Ramsey2012}. Although Guderley's imploding shock was recently shown to arise dynamically from classical (shock-free) initial data~\cite{CialdShkVic2025}, it cannot arise from initial data that is globally $C^2$-smooth.  A final drawback is that Guderley's imploding shock solution assumes propagation into a quiescent core at zero pressure (and hence identically zero temperature); consequently, it is inconsistent with both the vanishing-viscosity limit of the compressible Navier--Stokes equations and the hydrodynamic limit of standard kinetic models such as the Boltzmann or Landau equations.
\end{itemize}

\begin{figure}[htb!]
\centering
\includegraphics[width=0.55\textwidth]{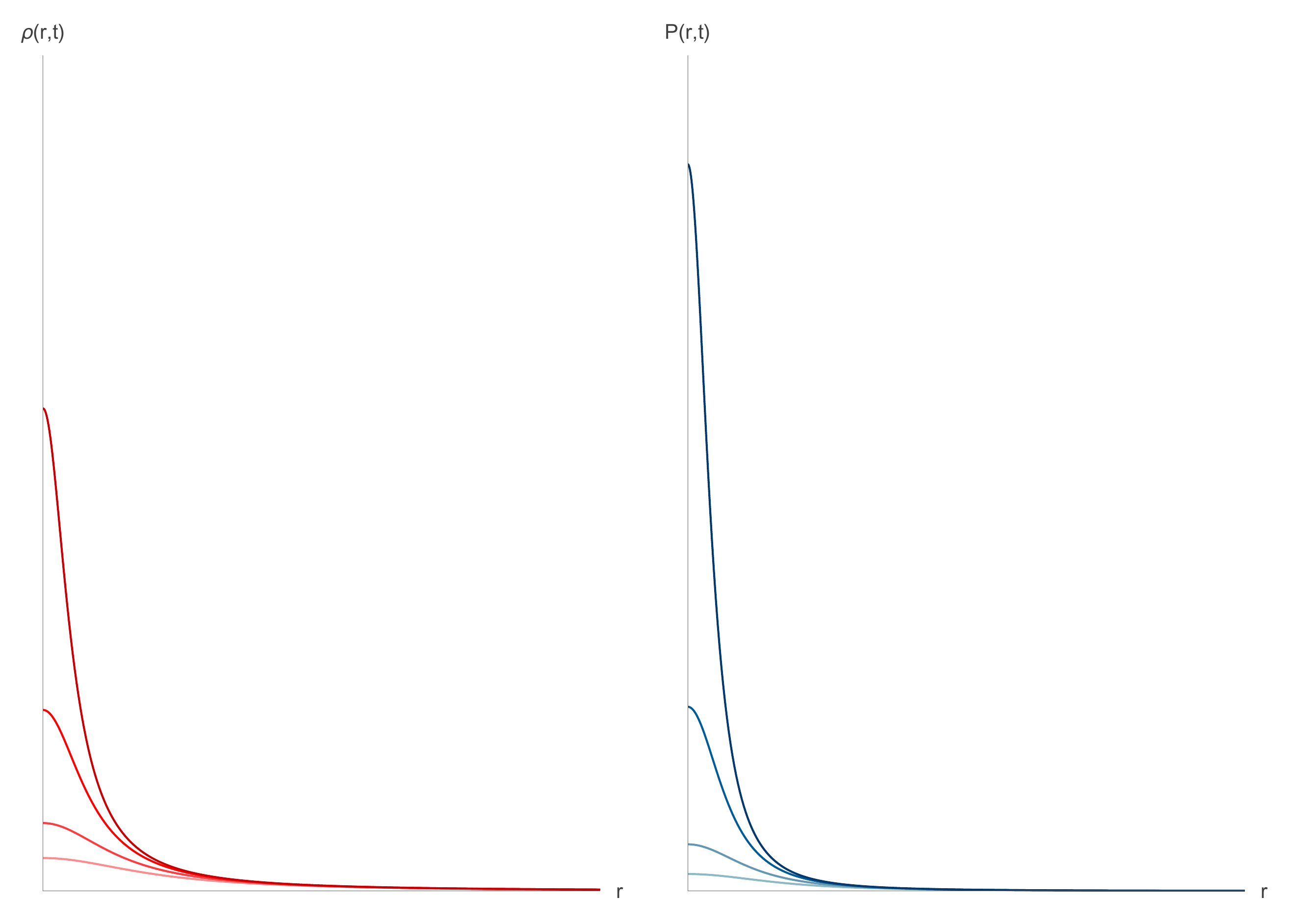}
\caption{Snapshots of the radial density $\rho(r,t)$ (red, left) and pressure $p(r,t)$ (blue, right) for the smooth isentropic implosion of Merle, Rapha\"el, Rodnianski, and Szeftel~\cite{MRRS2022a} (for $d=3$, $\gamma = 7/5$, $\NNN=1$, and similarity exponent as reported by Biasi~\cite{Biasi2021}), at four successive times, with darker shades indicating later times. The profiles are everywhere $C^\infty$-smooth and shock-free; in particular, no vacuum forms at any finite $r$, even though $\rho(r,t)$ and $p(r,t)$ decay rapidly in $r$. As $t \to t_*^-$, the relation $p = \tfrac{1}{\gamma}\rho^{\gamma}$ forces the density and pressure to diverge \emph{simultaneously} at $r=0$.}
\label{fig:MRSS}
\end{figure}

\begin{itemize}[leftmargin=1em]

\item \textbf{Smooth isentropic implosions}. 
A remarkable class of new Euler implosions was discovered by Merle, Rapha\"el, Rodnianski, and Szeftel~\cite{MRRS2022a}; see also the refinements in~\cite{BCG2025,ShaoWangWeiZhang2025}. These self-similar solutions contain no vacuum regions and remain $C^\infty$ \emph{smooth} until the singular time; in particular, they contain no shocks. These smooth implosions arise in the context of the \emph{isentropic} Euler system, with $p = \frac{1}{\gamma}\rho^\gamma$ and constant entropy; as such, as the density blows up, the pressure is forced to blow up as well, with the radial velocity also diverging at the time of collapse; see Figure~\ref{fig:MRSS}. They are constructed for \emph{quantized} (discrete) values of the similarity exponent, which are determined implicitly by requiring smooth passage through the sonic point. Thus, as in the Guderley problem, the similarity exponent of the smooth implosions constructed in~\cite{MRRS2022a,BCG2025,ShaoWangWeiZhang2025} is \emph{not explicit}, and can only be approximated numerically. A significant limitation of these smooth isentropic implosions is that they are \emph{conjectured to be unstable}, even within radial/spherical symmetry, on the basis of the numerical mode analysis of~\cite{Biasi2021}. Specifically, the studies in~\cite{Biasi2021} suggest that generic radial perturbations of the implosion profiles from~\cite{MRRS2022a} are deflected away from the smooth implosion, typically resulting in a gradient catastrophe (shock formation) before collapse. This conjectured instability (even in radial/spherical symmetry), combined with the fact that the similarity exponent can only be computed numerically (to a given precision), implies that smooth isentropic implosions cannot be observed in direct numerical simulations of the compressible Euler system, even when the scheme is restricted to the radially symmetric setting.
\end{itemize}

\begin{figure}[htb!]
\centering
\includegraphics[width=0.55\textwidth]{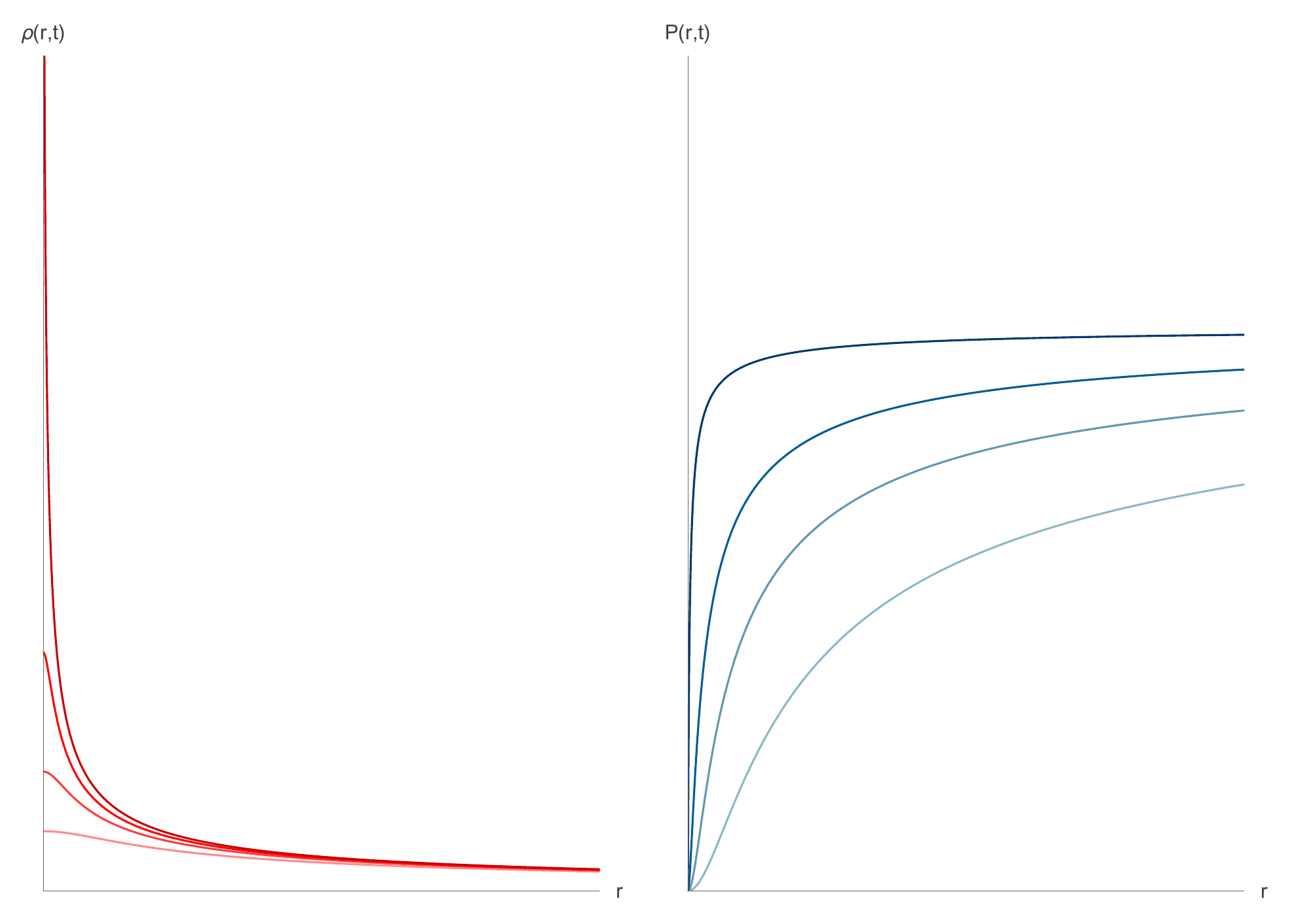}
\caption{Snapshots of the radial density $\rho(r,t)$ (red, left) and pressure $p(r,t)$ (blue, right) for the new ground-state implosion profile constructed in this paper, for $d=2$ with $\gamma = 2$ (2D monatomic gas), and $\NNN=1$, at four successive times, with darker shades indicating later times. As in the 3D monatomic case, the profiles are everywhere $C^\infty$-smooth and shock-free, the density is strictly positive for all~$r$, and the pressure vanishes \emph{only} at $r=0$. In contrast to the 3D case, however, as $t \to t_*^-$ \emph{only the density $\rho(r,t)$ diverges in a vicinity of the origin}: the pressure $p(r,t)$ remains uniformly bounded, with the spatial pressure profile steepening and saturating at a finite asymptote for large~$r$.}
\label{fig:New:Implosion:2}
\end{figure}

The goal of the present paper is to introduce a new class of smooth, self-similar implosion solutions for the \emph{full} (non-isentropic) compressible Euler equations (see Sections~\ref{sec:intro:existence} and~\ref{sec:intro:stability} for the main results in abbreviated form), which mitigates the limitations of all known previous constructions:

\begin{itemize}[leftmargin=1em]
\item \textbf{Smooth.} Given a dimension $d$ and an adiabatic exponent $\gamma$, we construct a sequence of globally self-similar imploding Euler solutions, indexed by an integer $\NNN \geq 1$, see~Section~\ref{sec:profiles}. 
These solutions have smooth (real analytic) density, velocity, and pressure; they are radially/spherically symmetric, and they have a strictly positive density (see Figure~\ref{fig:New:Implosion:2}). In particular, there is no shock discontinuity, and there is no quiescent core. See Remark~\ref{rem:globally:SS:exact} for the formula relating the Euler solution to the self-similar profiles.
\end{itemize}
\begin{figure}[htb!]
\centering
\includegraphics[width=0.55\textwidth]{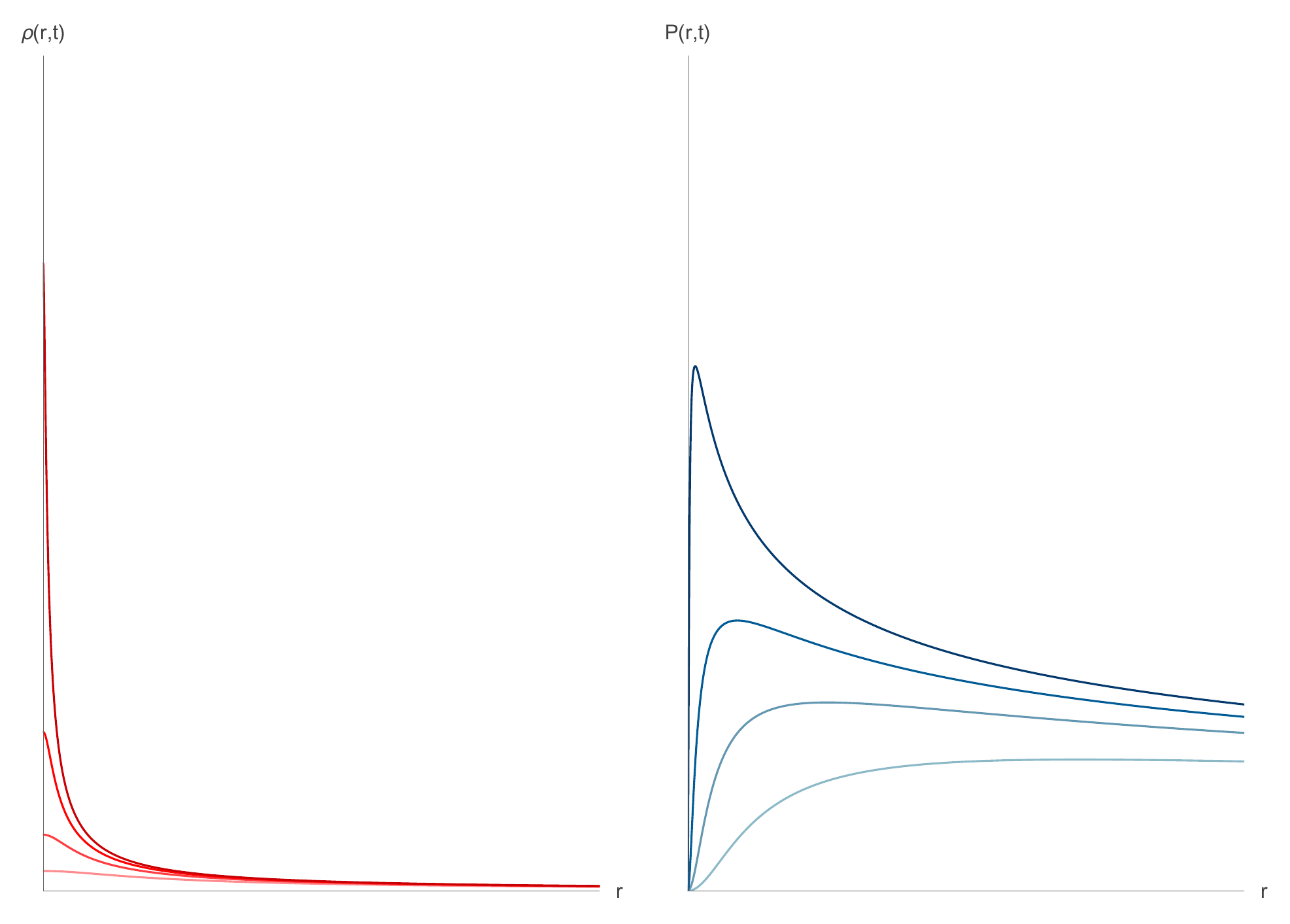}
\caption{Snapshots of the radial density $\rho(r,t)$ (red, left) and pressure $p(r,t)$ (blue, right) for the ground-state implosion profile constructed in this paper, for $d=3$, $\gamma = \tfrac 5 3$ (3D monatomic gas), and $\NNN=1$, at four successive times, with darker shades indicating later times. The profiles are everywhere $C^\infty$-smooth and shock-free, the density is strictly positive for all~$r\geq 0$, and the pressure vanishes \emph{only} at $r=0$; a one-point compatibility condition dynamically preserved by the Euler evolution. As $t \to t_*^-$, both the density and the pressure diverge in a vicinity of the origin, while $p(0,t) = 0$ is maintained for all~$t < t_*$.}
\label{fig:New:Implosion:1}
\end{figure}

\begin{itemize}[leftmargin=1em]
\item \textbf{Explicit similarity exponents.} The similarity exponents $\cxbar$, $\cubar$, and $\cbbar$ (see the self-similar ansatz in Section~\ref{sec:self:similar:intro}), are given by explicit closed-form expressions    in terms of the adiabatic constant $\gamma$, the spatial dimension $d$, and the integer parameter $\NNN \geq 1$, see~Definition~\ref{def:exponents}. These similarity exponents are largest for $\NNN=1$, and we refer to the associated self-similar profiles as the \emph{ground state} profiles. As opposed to Guderley's imploding shock, and as opposed to the smooth isentropic implosions in~\cite{MRRS2022a,BCG2025,ShaoWangWeiZhang2025}, no shooting method or numerical eigenvalue computation is required to determine the similarity exponents.

\item \textbf{Genuinely non-isentropic.} The solutions we analyze in this paper solve the full Euler system~\eqref{eq:euler0} and are genuinely non-isentropic. The pressure is taken to vanish at $x=0$ and be strictly positive elsewhere; together with $\uu(\cdot,0) = 0$, this is dynamically preserved by the Euler evolution~\eqref{eq:euler:primary}. The density, by contrast, is strictly positive everywhere, including at $x=0$. Consequently, since $p(0,t) = 0$ but $\rho(0,t) > 0$, the isentropic pressure law $p = \frac{1}{\gamma}\rho^\gamma$ cannot hold at $x=0$. This decoupling opens the possibility that the pressure and the density behave qualitatively differently at the time of collapse: the density always blows up (see Definition~\ref{def:implosion}, Remark~\ref{rem:SS:density:pressure}, and Section~\ref{sec:limiting:fields}), but the pressure need not (see Figure~\ref{fig:New:Implosion:2}).
\end{itemize}

\begin{itemize}[leftmargin=1em]
\item \textbf{Full stability of the ground state in radial symmetry.} For all $d$ and $\gamma$, the smooth implosion profile of the ground state at $\NNN = 1$ is asymptotically stable under perturbations in an \emph{open set of smooth radially/spherically symmetric} initial data, with pressure vanishing at the origin. See~Theorem~\ref{thm:stability:radial} and Remark~\ref{rem:ground:state:stable}, which provide the full stability picture under symmetry assumptions, for all $\NNN\geq 1$.  
Note that for Guderley's imploding shock, linear mode analysis~\cite{Morawetz1951,ChenZhangPanarella1995} supports stability within radial symmetry; however, no rigorous nonlinear PDE stability result for the Guderley solution is presently available, even under radial symmetry. We also note that the stability of our ground state profile (cf.~Theorem~\ref{thm:stability:radial}) stands in sharp contrast to the smooth isentropic implosions of~\cite{MRRS2022a,BCG2025,ShaoWangWeiZhang2025}, which are conjectured (on the basis of the numerical evidence in~\cite{Biasi2021}) to be unstable even within radial symmetry.

\item \textbf{Complete characterization of non-radial instabilities.} For the ground state profile at $\NNN=1$, and the physically most-relevant cases of a monatomic gas ($\gamma = \frac{d+2}{d}$) and diatomic gas ($\gamma = \frac{2d+1}{2d-1}$), we determine the \emph{exact} dimension of the unstable+center manifold for \emph{non-radial/non-spherical perturbations} and provide a complete characterization of this manifold.
Modulo these finitely many compatibility conditions on the initial data, the implosion is nonlinearly stable for general (non-radial) perturbations; see~Theorems~\ref{thm:nonrad:stab} and~\ref{thm:ODE_stab}. This provides, to the best of our knowledge, the first complete nonlinear PDE stability result for smooth imploding solutions of the multi-dimensional compressible Euler equations, with respect to perturbations that need not preserve any spatial symmetry of the underlying profile.

\item \textbf{A new numerical benchmark.} Due to their stability properties and the explicit nature of the similarity exponents, the smooth Euler implosions constructed in this paper (corresponding to the ground state at $\NNN=1$ in the case of a monatomic/diatomic gas) provide new exact analytic solutions of the multi-dimensional Euler equations in a regime of unbounded density growth, and may be used to \emph{benchmark the accuracy of existing numerical codes}.
\begin{itemize}
\item[$\circ$] For numerical codes which implement the Euler equations directly in radial symmetry~\eqref{eq:euler2}, the advantages our smooth implosions present over the Guderley solution (which has historically been used extensively as a benchmark~\cite{Ramsey2012}) are threefold: the similarity exponents are explicit (no shooting method required to compute them), the solutions do not contain a shock (no front-tracking required), and the pressure only vanishes at one point (a boundary/compatibility condition dynamically preserved by Euler cf.~\eqref{eq:euler:p}) instead of in the entire quiescent region (zero pressure in an open set is difficult to propagate numerically with the existing regularization schemes).

\item[$\circ$] For numerical codes which implement the Euler equations~\eqref{eq:euler:primary} without any symmetry assumptions, our solutions present significant advantages that go beyond those mentioned above. Chief among these is the fact that our smooth implosions are stable even outside of spherical symmetry, if the initial data for density, velocity, and pressure is chosen to satisfy \emph{finitely many, explicit,} compatibility conditions at $x=0$;\footnote{
Near $x=0$, the power-series coefficients of our smooth imploding profiles are given \emph{explicitly}; hence the power series can be derived \emph{analytically}, without numerical error, to \emph{arbitrary} order. 
We impose explicit compatibility conditions on the perturbation, which in turn yield explicit compatibility conditions on the initial data.}
see~\eqref{eq:vanishing:order:at:origin:default} and item (v) in~Theorem~\ref{thm:ODE_stab}. 
In particular, we provide large classes of initial data with nonzero vorticity and nontrivial entropy variation, which are constant outside a compact set (cf.~Remark~\ref{rem:compact:support}) and which lead to a finite-time implosion in a \emph{stable fashion}.

\end{itemize}

\item \textbf{The hydrodynamic limit for a broad range of kinetic collision operators.} 
It was recently shown in~\cite{bedrossian2026finite} that there exist smooth, strictly positive initial data for the inhomogeneous Landau equation with very hard potentials $\gamkin \in (\sqrt{3},2]$, which lead to finite time singularities. The Landau 
singular solution constructed in~\cite{bedrossian2026finite} is \emph{asymptotically hydrodynamic}, in the sense that the distribution function converges to a local Maxwellian centered on hydrodynamic fields (density, velocity, temperature) which develop an asymptotically
self-similar Euler implosion, with smooth profiles taken from~\cite{ShaoWangWeiZhang2025}.
The constraint $\gamkin > \sqrt{3}$, which does not cover three exponents $\gamkin = -3$ (Coulomb potential), $\gamkin=0$ (Maxwell molecules), and $\gamkin = 1$ (hard spheres), arises largely because the smooth implosions in~\cite{ShaoWangWeiZhang2025} are isentropic ($\cbbar=0$).
This restricts the admissible range of $\gamkin$ for satisfying the ``self-similar hydrodynamic limit'' constraint
$-3\cbbar + (\gamkin+3)\cubar + 1 < 0$
(see~\eqref{eq:hydro}). 
The self-similar implosions constructed in this paper for a 3D monatomic gas ($\gamma = 5/3$) allow for nontrivial entropy, with $\cbbar > 1/3$. Remarkably, condition~\eqref{eq:hydro}, which becomes $\gamkin \cubar < 1/2$ in this setting (see~\eqref{eq:hydro2}), may then be satisfied for \emph{every} $\gamkin \in [-3,2]$, by appropriately choosing the integer parameter $\NNN \geq 1$. Moreover, for the above three cases $\gamkin \in\{-3,0,1\}$, condition~\eqref{eq:hydro2} is achieved by the ground state profile at $\NNN = 1$, which enjoys enhanced stability properties. See~Section~\ref{sec:intro:kinetic} for a detailed discussion.
\end{itemize}

The remainder of the Introduction is dedicated to explaining the PDE setup considered in this paper, presenting our main results in abbreviated form, discussing the related literature, and outlining the main ideas of the proof.

\subsection{The Euler equations}
The full (non-isentropic) Euler equations, for the unknown density $\rho$, momentum $\rho \uu$, and total energy $E$ are given by the system of conservation laws 
\begin{subequations} 
\label{eq:euler0}
\begin{align} 
\p_t \rho +  \operatorname{div}(\rho  \uu ) & = 0  ,  \label{eq:rho0} \\
\p_t (\rho \uu) + \operatorname{div}(\rho \uu \otimes \uu  +  p \, \Id \!) & = 0 , \label{eq:u0} \\
\p_t E  + \operatorname{div} ((p+E) \uu) &=0 , \label{eq:E0}
\end{align} 
where the pressure $p$ is given by
\begin{equation}
p = (\gamma-1) \bigl( E - \tfrac 12 \rho |\uu|^2 \bigr) , \label{eq:p0}
\end{equation}
\end{subequations}
and $\gamma>1$ is the adiabatic exponent. Physically relevant values of $\gamma$ include $5/3$ (monatomic gas in 3D), $7/5$ (diatomic gas in 3D), etc. The equations are posed on $\Reals^d$, in dimensions $d \in \{1, 2, 3\}$. Throughout this paper we assume for convenience that $\gamma \leq 2 d+1$.\footnote{This range of adiabatic exponents far exceeds the  physically relevant values for gas dynamics, which correspond to $\gamma = \frac{d+2}{d}$ for a monatomic gas, $\gamma=\frac{2d+1}{2d-1}$ for a diatomic gas, with smaller values of $\gamma$ for gases whose molecules contain more than three atoms. Our results extend beyond $\gamma = 2d+1$; in this paper we do not pursue the maximal range of adiabatic exponents.}  

Within the class of classical solutions,\footnote{If the flow variables are $C^{1}_{x,t}$-smooth and the density is strictly positive, globally in spacetime, the Euler solution is referred to as being ``classical'' or ``genuine''; see e.g.~\cite{Lax1954}.} the Euler system for the conservation-law variables $(\rho, \rho \uu, E)$, namely~\eqref{eq:euler0}, may be recast as an evolution equation for the primary flow variables: density $\rho$, velocity $\uu$, and pressure $p$, as
\begin{subequations}
\label{eq:euler:primary}
\begin{align}
\p_t\rho +    \uu \cdot \nabla \rho  + \rho \, \div \uu & = 0  ,  \label{eq:euler:rho} \\
\p_t \uu +  \uu \cdot \nabla \uu  +  \tfrac{1}{\rho}  \nabla p & =0 , \label{eq:euler:u} \\
\p_t p + \uu \cdot \nabla p + \gamma p \, \div \uu&=0 . \label{eq:euler:p}
\end{align} 
\end{subequations}
We find it useful to further rewrite system~\eqref{eq:euler:primary} as an evolution equation for 
the  rescaled sound speed  $\sigma$, the fluid velocity  $\uu$, and the square root of the pseudo-entropy $b$, as follows. To do so, for $\gamma>1$  we define $\alpha: = \frac{\gamma-1}{2} >0$.  The constraint $\gamma \leq 2d+1$ amounts to $\alpha \leq d$. We may then use the ideal gas relation 
\begin{subequations}
\label{eq:constitutive:relations}
\begin{equation}
p = p(\rho,s) = \tfrac{1}{\gamma} \rho^\gamma e^s ,
\end{equation}
where $s$ is the specific entropy, to define the square-root of pseudo-entropy
\begin{equation}
b: = \sqrt{\tfrac{\gamma p}{\rho^\gamma}} = e^{\frac s2},
\label{eq:b:def}
\end{equation}
the sound speed
\begin{equation}
 c = \sqrt{\tfrac{dp}{d\rho}} = \rho^\alpha e^{\frac{s}{2}} = \rho^\alpha b ,
\end{equation}
and the rescaled sound speed as
 \begin{equation}
 \sigma = \tfrac{1}{\alpha} c = \tfrac{1}{\alpha} \rho^\alpha b .
 \label{eq:sigma:def}
\end{equation}
\end{subequations}
The Euler system for the conservation-law variables $(\rho, \rho \uu, E)$~\eqref{eq:euler0}, or for the primary flow variables $(\rho,\uu,p)$, is then \emph{equivalent} within the class of smooth solutions to the Euler system for the fundamental variables $(\sigma,\uu,b)$:
\begin{subequations} 
\label{eq:euler1}
\begin{align} 
\p_t\sigma +   \uu \cdot \nabla \sigma  +  \alpha  \sigma  \operatorname{div} \uu & = 0  ,  \label{eq:sigma1} \\
\p_t \uu +  \uu \cdot \nabla \uu  +   \alpha  \sigma \nabla \sigma -  \tfrac{\alpha}{\gamma} \sigma^2 \tfrac{\nabla b}{b} & =0 , \label{eq:u1} \\
\p_t b + \uu \cdot \nabla b &=0 . \label{eq:s1}
\end{align} 
\end{subequations} 
The density, sound speed, pressure, and total energy may be recovered from $(\sigma,\uu,b)$ via~\eqref{eq:p0}--\eqref{eq:constitutive:relations}, namely
\begin{equation}
\rho = \bigl( \tfrac{\alpha \sigma}{b} \bigr)^{\frac{1}{\alpha}},
\qquad
c = \alpha \sigma,
\qquad
p = \tfrac{\alpha^2}{\gamma} \rho \sigma^2,
\qquad
E= \tfrac 12  \rho (\tfrac{\alpha}{\gamma} \sigma^2 +  |\uu|^2) 
.
\label{eq:new:constitutive:relations}
\end{equation}
The Cauchy problem for the Euler system~\eqref{eq:euler1} (or equivalently~\eqref{eq:euler:primary} or~\eqref{eq:euler0}) is supplemented with smooth initial conditions at a given initial time, taken to be $t=-1$ in this paper. 

Throughout this paper, we will consider solutions which are \emph{smooth} (see Definition~\ref{def:smooth:solution}) on $\Reals^d \times [-1,t_*)$, where $t_* \approx 0$ is the time at which the implosion singularity develops.

\begin{definition}[\bf Smooth solution]
\label{def:smooth:solution}
A solution of the Euler system~\eqref{eq:euler0}, or equivalently~\eqref{eq:euler:primary}, is called \emph{smooth} if the conservation law variables $(\rho, \rho \uu, E)$, or equivalently~the primary flow variables $(\rho,\uu,p)$, are \emph{smooth functions of space and time},\footnote{In this paper, by ``smooth'' we mean \emph{very smooth}: the exact imploding solutions are real-analytic, while their perturbations are taken to lie in a H\"older space $C^\MMM$, with $\MMM$ sufficiently large.}  and the density is \emph{strictly positive} ($\rho > 0$) on $\Reals^d \times [-1,t_*)$.
\end{definition}

For smooth solutions of the Euler system on $\Reals^d \times [-1,t_*)$, the three formulations~\eqref{eq:euler0}, \eqref{eq:euler:primary}, and~\eqref{eq:euler1} are equivalent because the density is (strictly) bounded away from zero. 

\begin{remark}[\bf Smooth solutions need not have smooth sound speed]
We note that while smooth solutions have smooth density $\rho$, velocity $\uu$, and pressure $p$, the sound speed $c$, the rescaled sound speed $\sigma$, and the  square-root of pseudo-entropy $b$ need not be smooth at $x=0$. To see this, assume that the density is constant and that the pressure is given by a constant  times $|x|^2$, near $x=0$; these are both $C^\infty$ functions of $x$. Then, the fields $b$, $c$, and $\sigma$ all behave as a constant times $|x|$ near $x=0$, which is not a smooth function of $x$. 
\end{remark}

\subsection{What is an implosion?}
Informally, an implosion is a flow in which geometric focusing causes a gas to converge inward toward a point, compressing matter to infinite density or pressure in finite time. We now make this notion mathematically precise, in the context of the compressible Euler system.

\begin{definition}[\bf Euler implosion]
\label{def:implosion}
A solution\footnote{By \emph{solution}, we mean a \emph{weak solution} $(\rho,\rho \uu,E)$ of the Euler equation in conservation law form \eqref{eq:euler0}, which is at least $C_{x,t}^1$-smooth on the complement of finitely many co-dimension $1$ (orientable and smooth) hypersurfaces in the spacetime $\Reals^d \times [-1,t_*)$, with strictly positive density, and which satisfies the Lax geometric (entropy) conditions. This notion of solution is quite general; in practice we either encounter \emph{regular shock solutions}---which are smooth on the complement of a co-dimension $1$ spacetime shock front, across which the density, normal momentum, and energy satisfy the Rankine-Hugoniot jump conditions, see~\cite[Definition 1.1]{CialdShkVic2025}, or \cite[Definition 1.1]{BDSV}---or \emph{smooth solutions} in the sense of Definition~\ref{def:smooth:solution} above.} of the compressible Euler equations on $\Reals^d \times [-1,t_*)$
is called an \emph{implosion solution} with \emph{implosion center}~$x_* \in \Reals^d$
and \emph{implosion time}~$t_*>-1$, if the following hold:
\begin{itemize}[leftmargin=2em]
\item[(i)]  \textsl{(Blowup of a primary flow variable).} As one approaches the spacetime point of collapse, $(x_*, t_*)$, \emph{at least one of the primary thermodynamic variables} (density or
pressure) becomes unbounded; that is, for \emph{every} $\varepsilon \in (0,1]$ we have
\begin{subequations}
\label{eq:implosion:big:def}
\begin{equation}
  \limsup_{t \to t_*^-}\;
    \sup_{|x - x_*| < \varepsilon} \frac{\rho(x,t)}{\|\rho(\cdot,-1)\|_{L^\infty(B_{100}(0))}}
   +  
   \limsup_{t \to t_*^-}\;
    \sup_{|x - x_*| < \varepsilon} \frac{p(x,t)}{\|p(\cdot,-1)\|_{L^\infty(B_{100}(0))}}
   = +\infty.
   \label{eq:implosion:big:def:a}
\end{equation}
 
\item[(ii)] \textsl{(Inward-focusing character).}
For  $t$ sufficiently close
to~$t_*$, the flow is \emph{compressive} toward~$x_*$;  that is, there exists a time-dependent radius $\mathsf{r} \colon [-1,t_*) \to (0,1]$ such that
\begin{equation}
  \bigl(\uu(x,t) - \uu(x_*,t)\bigr) \cdot (x - x_*) \leq 0\,,
  \label{eq:implosion:big:def:b}
\end{equation}
for all $|x - x_*| \leq \mathsf{r}(t) $, and all $t \in [-1,t_*)$.

\item[(iii)] \textsl{(No amplitude blowup away from the spacetime point of collapse).}
The primary flow variables $(\rho,\uu,p)$ are bounded at all spacetime points which are  away from $(x_*,t_*)$; that is, for any $\eps \in (0,1]$  there exists $C_\eps>0$ such that
\begin{equation}
\sup_{\substack{(x,t) \in \Reals^d \times[-1,t_*) \\ \eps \leq |x-x_*| + |t-t_*| \leq 1/\eps} } 
\frac{\rho(x,t)}{\|\rho(\cdot,-1)\|_{L^\infty(B_{100}(0))}}
+
\frac{|\uu(x,t)|}{\|\uu(\cdot,-1)\|_{L^\infty(B_{100}(0))}}
+
\frac{p(x,t)}{\|p(\cdot,-1)\|_{L^\infty(B_{100}(0))}}
\leq C_\eps.
\label{eq:implosion:big:def:c}
\end{equation}
\end{subequations}
\end{itemize}
\end{definition}

\begin{remark}[\bf Velocity blowup at the implosion center?]\label{rem:velocity}
It is natural to ask whether the blowup of density or pressure at~$(x_*,t_*)$, in the sense of~\eqref{eq:implosion:big:def:a}, 
forces the velocity amplitude to blow up there as well. For all previously known implosions~\cite{Guderley1942,Stanyukovich1960,Lazarus1981,Sedov2018,JangLiuSchrecker2025,MRRS2022a,BCG2025,ShaoWangWeiZhang2025,jenssen2023radially,Jenssen2025} this is the case: $|\uu(x(t),t)|$ diverges as $t \to t_*^-$, for a suitably chosen $x(t) \to x_*$. On the other hand, integrating~\eqref{eq:euler:rho} and~\eqref{eq:euler:p} along the Lagrangian flow associated to $\uu$ shows that density or pressure blowup, cf.~\eqref{eq:implosion:big:def:a}, only implies that the compressive strain rate diverges, namely $\div \uu \to -\infty$ as one approaches $(x_*,t_*^-)$. There is no clear mechanism in~\eqref{eq:euler:primary} that forces the blowup of the velocity amplitude: the center of implosion may in principle drift at finite speed, while the strain rate around it diverges. 

The present paper settles this question: for some of the Euler implosions constructed in this work, the velocity amplitude $|\uu|$ \emph{does not blow up} in the vicinity of $(x_*,t_*)$,  in the sense of~\eqref{eq:implosion:big:def:a}. This situation occurs when $\cxbar \geq 1$, a condition which holds for instance when $d=3$, $\gamma = 5/3$, and $\NNN$ is 
small. 
On the other hand, for the  Euler implosions with $\cxbar < 1$ constructed in this work, the velocity amplitude $|\uu|$ does blow up in the vicinity of $(x_*,t_*)$. The condition $\cxbar < 1$ holds for instance when $d=3$, $\gamma = 5/3$, and $\NNN$ is taken to be sufficiently large. See Remark~\ref{rem:SS:density:pressure}, Section~\ref{sec:limiting:fields}, and Example~\ref{ex:exponents:monatomic:3D}. 
\end{remark}
  
\subsection{Different flavors of Euler implosions}
Definition~\ref{def:implosion} is broad enough to encompass every class of Euler implosions known to us; to emphasize this point, we next list several natural specializations.

Consider an Euler implosion on $\Reals^d \times [-1,t_*)$ with implosion center $x_*$ and implosion time $t_*$, in the sense of Definition~\ref{def:implosion}. We have the following flavors:
\begin{itemize}[leftmargin=1em]

\item \textbf{Shock implosion.}
\emph{An implosion solution which is a regular shock solution of the Euler system, which contains a smooth co-dimension $1$ shock front $\Sigma(t)$ that collapses onto $x_*$ as $t\to t_*^-$, and with the shock strength (the jump in~$\rho$ or~$p$ across~$\Sigma(t)$) diverging as $t\to t_*^-$, is called a shock implosion.} 

By ``regular shock solution'' we mean a weak solution $(\rho,\rho \uu,E)$ of the Euler equation~\eqref{eq:euler0}, which is $C_{x,t}^1$-smooth on the complement of a co-dimension $1$ spacetime shock front $\Sigma(t)$, across which the jumps in density, normal momentum, and energy, satisfy the Rankine-Hugoniot conditions, and which satisfies the Lax geometric (entropy) conditions (cf.~\cite[Definition 1.1]{CialdShkVic2025} or \cite[Definition 1.1]{BDSV}). The converging shock front is typically a sphere, i.e., $\Sigma(t) = \{|x - x_*| = r_s(t)\}$ with $r_s(t) \to 0$ as $t \to t_*^-$. The requirement of a diverging shock strength is what ensures~\eqref{eq:implosion:big:def:a}, which is then  realized by the  values of~$\rho$
or~$p$ downstream of (behind) the shock. This is exactly the Guderley~\cite{Guderley1942} and Landau-Stanyukovich~\cite{Stanyukovich1960} setting, which has been by now revisited multiple times, see e.g.~\cite{Lazarus1977,Lazarus1981,MeyerTerVehn1982,Barenblatt1996,Ramsey2012,Sedov2018,JenssenTsikkou2018,Giron2023,JangLiuSchrecker2025,CialdShkVic2025} and references therein. We note that for Guderley's imploding shock, the density does not blow up as one approaches $(x_*,t_*)$, while the pressure and velocity amplitude do; thus ~\eqref{eq:implosion:big:def:a} is attained due to the pressure, not the density.

Shock implosions may be similarly defined in the case of multiple (nested) converging shock fronts; their construction is for instance discussed in~\cite{Lazarus1977,Lazarus1981,Giron2023}.

\item \textbf{Smooth implosion.}
\emph{If the primary flow variables $(\rho, \uu, p)$ are $C^\infty$-smooth in space and time on $\Reals^d \times [-1,t_*)$, and if the density is strictly positive on this set, we refer to the  imploding solution as being smooth.}

The solutions constructed in~\cite{MRRS2022a,BCG2025,ShaoWangWeiZhang2025} and the ones constructed in this paper fall in this category. Note  that for the implosions in~\cite{MRRS2022a,BCG2025,ShaoWangWeiZhang2025} both density and pressure blow up as one approaches $(x_*,t_*)$, in the sense of~\eqref{eq:implosion:big:def:a}; in fact, the velocity amplitude blows up too. However, for the implosions considered in this paper the pressure need not blow up at the location of collapse (when $\gamma(\cxbar-1) \geq \cbbar$), and the velocity amplitude need not blow up either (when $\cxbar \geq 1$); the density does always blow up, implying~\eqref{eq:implosion:big:def:a}; see Remark~\ref{rem:SS:density:pressure} and Section~\ref{sec:limiting:fields}.

\item \textbf{$C^k$-smooth implosion.}
\emph{We refer to the imploding solution as being $C^k$-smooth, if there exists an integer $k\geq 0$ such that the primary flow variables $(\rho, \uu, p)$ are $C^k$-smooth, but not $C^{k+1}$-smooth, in space and time on $\Reals^d \times [-1,t_*)$, and the density is strictly positive on this set.}

The Euler implosions discussed in~\cite{jenssen2023radially} fall into this category, with $k=0$ (continuous but not differentiable amplitude blowup); similarly, the solutions discussed in~\cite{Jenssen2025} are continuous, but only have finite regularity. For these solutions all primary flow variables $(\rho,\uu,p)$ blow up at the spacetime point of collapse. We also mention that the constructions in~\cite{MRRS2022a,BCG2025} allow one to construct $C^k$-smooth implosions for generic values of the similarity exponent. 
 
\item \textbf{Radially/spherically symmetric implosion.}
\emph{If the  implosion solution is such that the density and pressure depend only on $r = |x - x_*|$ and on $t$, and if the velocity is purely radial, meaning that  $\uu(x,t) = \frac{x - x_*}{|x-x_*|} u^r(|x - x_*|,  t)$ for a scalar function of two variables $u^r$, we refer to the solution as being radially symmetric (when $d=2$) or spherically symmetric (when $d=3$).}

Note that for radially/spherically symmetric implosions we may always take $x_*=0$, by appealing to the translation symmetry of the Euler system; this identification is assumed tacitly throughout the paper, see e.g.~Section~\ref{sec:radial:symmetry}. Note that in this setting the continuity of $\uu$ at $(x_*,t)$ with $t<t_*$ implies that $u^r(0,t) = 0$, which in turn implies that $\uu(x_*,t) = 0$ (the only vector fixed by every $SO(d)$ rotation is the zero vector); in particular, condition~\eqref{eq:implosion:big:def:b} reduces to $u^r(r,t) \leq 0$ for $0 \leq r \leq \mathsf{r}(t) \to 0$ as $t\to t_*^-$.

All the implosion solutions discussed above, those containing shocks~\cite{Guderley1942,Stanyukovich1960,Lazarus1977,Lazarus1981,MeyerTerVehn1982,JenssenTsikkou2018,JangLiuSchrecker2025}, the $C^\infty$-smooth isentropic ones~\cite{MRRS2022a,BCG2025,ShaoWangWeiZhang2025}, the ones which are continuous but not smooth~\cite{jenssen2023radially,Jenssen2025}, and also the $C^\infty$-smooth non-isentropic solutions considered in Sections~\ref{sec:profiles} and~\ref{sec:stability} of this paper, they are all radially/spherically symmetric.\footnote{The rotational symmetry of the Euler equations (a part of the Galilean symmetry group) is what allows one to make this  reduction from $d+1$ dimensions to $1+1$ dimensions.}

The only known examples of \emph{non-radial implosions} are those constructed by perturbing the radially/spherically symmetric implosions, establishing that the implosion mechanism survives non-radial perturbations; see e.g.~the recent works~\cite{CGSS2025,CCSV2024,Chen2024} and the solutions constructed in Section~\ref{sec:nonradial} of this manuscript.
 
\item \textbf{Self-similar implosion.}
\emph{An implosion solution is called \emph{globally} self-similar, if the primary flow variables are invariant (globally) under the scaling symmetry of the Euler equations. In the radially/spherically symmetric case, this means that there exist \emph{similarity exponents} $\delta  \in \Reals$ and $\beta >0$, such that with the time-dependent scaling function $\lambda(t) = (t_*-t)/(t_*+1)$  we have
\[
\rho(r,t) = \lambda(t)^{\delta} \rho\bigl(\tfrac{r}{\lambda(t)^\beta},-1\bigr),
\quad
u^r(r,t) = \lambda(t)^{\beta-1} u^r\bigl(\tfrac{r}{\lambda(t)^\beta},-1\bigr),
\quad
p(r,t) = \lambda(t)^{\delta + 2 (\beta-1)} p\bigl(\tfrac{r}{\lambda(t)^\beta},-1\bigr),
\]
for all $t\in [-1,t_*)$ and all $r\geq 0$.\footnote{In the setting of this paper, we have denoted $\beta \mapsto \cxbar$, and $\delta \mapsto (\cxbar -1 - \cbbar)/\alpha$; see Remark~\ref{rem:SS:density:pressure} and Remark~\ref{rem:globally:SS:exact}. The notation used by Lazarus~\cite{Lazarus1981} and most subsequent  papers on Guderley imploding shocks is: $\beta \mapsto 1/\lambda$ and $\delta \mapsto \kappa/\lambda$. The notation used in the smooth isentropic implosion papers~\cite{MRRS2022a,BCG2025} is $\beta \mapsto 1/r$ and $\delta \mapsto (1-r)/(\alpha r)$.} The functions $(\rho,u^r,p)(\cdot,-1)$ are called \emph{similarity profiles}.}  

We have only precisely defined \emph{globally self-similar} Euler implosions in the radially/spherically symmetric setting, because this is the \emph{only setting} in which  globally self-similar implosions were shown to exist:~\cite{Guderley1942,Stanyukovich1960,Lazarus1977,Lazarus1981,MeyerTerVehn1982,JenssenTsikkou2018,JangLiuSchrecker2025}, \cite{MRRS2022a,BCG2025,ShaoWangWeiZhang2025},  \cite{jenssen2023radially,Jenssen2025}, and  Section~\ref{sec:profiles} of this paper. This is by now a classical idea, discussed at length in the books~\cite{Stanyukovich1960,LandauLifshitz1987,Barenblatt1996,Sedov2018,ZelDovich2002,EggersFontelos2015}: by combining radial/spherical symmetry with a globally self-similar ansatz, one may reduce the $d+1$ dimensional PDE to an ODE, where powerful construction tools are available. 
 
We note that as opposed to the Taylor--von Neumann--Sedov blast wave solution~\cite{Taylor1950,vonNeumann1947,Sedov1946}, which describes a spherically symmetric shock wave expanding from a point explosion, and for which the similarity exponents $\beta$ and $\delta$ may be determined from dimensional analysis, the similarity exponents of globally-self similar Euler implosions cannot be determined from dimensional analysis; instead, $\beta$ and $\delta$ are determined by solving a so-called ``nonlinear eigenvalue problem''.\footnote{This is a rather loose terminology used by Barenblatt~\cite{Barenblatt1996}. In practice, the similarity exponents for Guderley's imploding shock wave may be computed by a shooting method~\cite{JangLiuSchrecker2025}, while the similarity exponents of smooth isentropic implosions were obtained as the zeros of a complicated function~\cite{MRRS2022a,BCG2025}. The similarity exponents $\cxbar$ and $\cbbar$ in this paper truly are obtained by solving an eigenvalue problem; see Section~\ref{sec:similarity:expo:matrix}.}
 
We note that in order to construct non-radial implosions, or to consider Euler implosions whose far-field state is a non-vacuum constant (i.e.,~strictly positive density and pressure), one cannot work within the class of \emph{globally} self-similar solutions, and must instead consider \emph{asymptotically self-similar} solutions (see Section~\ref{sec:self:similar:intro}). The existence of such solutions is typically proven by showing that the globally self-similar profiles are asymptotically stable (as the self-similar time diverges) to carefully designed perturbations; see~\cite{MRRS2022b,BCG2025,CGSS2025,CCSV2024,Chen2024} and the solutions constructed in Sections~\ref{sec:stability} and~\ref{sec:nonradial} of this manuscript.
\end{itemize}

\subsection{Radial/Spherical symmetry}
\label{sec:radial:symmetry}
The \emph{globally self-similar} implosions discussed in this paper are radially/spherically symmetric.
As such, Sections~\ref{sec:profiles} and~\ref{sec:stability} of this paper consider  the Euler equations~\eqref{eq:euler1} within the class of solutions which posses radial symmetry (when $d=2$), or spherical symmetry (when $d=3$):
\[
\uu(x,t) = \vec{e}_r u^r(r,t), \quad
\sigma(x,t) = \sigma(r,t), \quad 
b(x,t) = b(r,t), \quad 
\rho(x,t) = \rho(r,t), \quad 
p(x,t) = p(r,t),
\] 
where $r= |x| \in [0,\infty)$ and $\vec{e}_r = \frac{x}{r}$. Using that $\operatorname{div} \uu = \p_r u^r + \frac{d-1}{r} u^r$, we may write~\eqref{eq:euler1} as
\begin{subequations} 
\label{eq:euler2}
\begin{align} 
\p_t\sigma +    u^r \p_r \sigma  +  \alpha  \sigma ( \p_r u^r + \tfrac{d-1}{r} u^r)& = 0  ,  \label{eq:sigma2} \\
\p_t u^r +   u^r \p_r   u^r  +   \alpha  \sigma \p_r \sigma - \tfrac{\alpha}{\gamma} \sigma^2 \tfrac{\p_r b}{b}& = 0  , \label{eq:u2} \\
\p_t b + u^r \p_r  b &=0 , \label{eq:s2}
\end{align} 
\end{subequations} 
while \eqref{eq:euler:primary} may be written as
\begin{subequations}
\label{eq:euler2:primary}
\begin{align}
\p_t\rho +   u^r \partial_r  \rho  + \rho \, ( \p_r u^r + \tfrac{d-1}{r} u^r) & = 0  ,    \\
\p_t u^r  + u^r \partial_r  u^r  +  \tfrac{1}{\rho}  \partial_r p & =0 ,   \\
\p_t p + u^r \partial_r  p + \gamma p \, ( \p_r u^r + \tfrac{d-1}{r} u^r)&=0 .  
\end{align} 
\end{subequations}
Since $\sigma$ has the same units as $u^r$ (they are velocities), and since $b$ is purely transported, in this symmetric setting it is convenient to work with the formulation~\eqref{eq:euler2}.

For simplicity, we shall henceforth refer to both radially symmetric ($d=2$) and to spherically symmetric ($d=3$) functions as being \emph{radially symmetric}. When $d=1$, by radially symmetric functions we mean either even functions (for density, pressure), or odd functions (for velocity). 

Note, however, that in Section~\ref{sec:nonradial} we will consider non-radially symmetric implosions, which are only \emph{asymptotically self-similar}; in this setting, we need to work with the Euler evolution from~\eqref{eq:euler:primary}.

\subsection{Self-similar ansatz and  Euler evolution in radial symmetry}
\label{sec:self:similar:intro}
Our initial aim is to analyze self-similar solutions of \eqref{eq:euler2}, which are smooth at the initial time $t =-1 $ and develop a first singularity at time $t = t_* \approx 0$. For this purpose, we will use \emph{modulated self-similar coordinates} and \emph{modulated self-similar variables}, which will be defined precisely in Section~\ref{sec:stability} below. 

For clarity of the presentation, here we only introduce self-similar coordinates and variables \emph{as if they were not modulated}, meaning, they correspond to \emph{globally self-similar solutions}. To emphasize the fact that the ``real'' definitions are given later in Section~\ref{sec:true:ss:ansatz}, here we use the $\approx$ symbol instead of the $=$ symbol; roughly speaking, the ``$\approx$'' symbol becomes an ``$=$'' symbol as $t\to t_*\approx 0$. We proceed as follows:
\begin{itemize}[leftmargin=1em]
\item The initial data is specified at time $t = -1$.
\item The self-similar time variable $\tau \in [0,\infty)$ is given by
\[
\tau \approx - \log(-t)
\quad \Leftrightarrow \quad 
(-t) \approx e^{-\tau}
.
\]
Note that $t = -1 \Leftrightarrow \tau \approx0$, and $(-t) \to 0^- \Leftrightarrow \tau \to \infty$.

\item The self-similar radial space variable $R$ is defined as 
\[
R \approx \tfrac{r}{(-t)^{\cx}}
\quad \Leftrightarrow \quad 
r \approx R e^{-\cx \tau},
\]
where $\cx > 0$ is a \emph{space similarity exponent}. 

\item The self-similar radial velocity $U$  and self-similar rescaled sound speed $\Sigma$ are defined as 
\begin{align*}
&u^r (r,t) 
\approx (-t)^{\cu} U\bigl(\tfrac{r}{(-t)^{\cx}}, - \log(-t)\bigr) 
\quad \Leftrightarrow \quad
u^r(r,t) 
\approx
e^{-\tau \cu} U(R,\tau), \\
&\sigma (r,t) 
\approx (-t)^{\cu} \Sigma\bigl(\tfrac{r}{(-t)^{\cx}}, - \log(-t)\bigr) 
\quad \Leftrightarrow \quad
\sigma(r,t) 
\approx
e^{-\tau \cu} \Sigma(R,\tau),
\end{align*}
where $\cu \in \Reals$ is a \emph{speed similarity exponent}. At a first reading, it is convenient to use 
\[
\cu \approx \cx - 1,
\]
which is an expression of the fact that we must balance the $\p_t$ and the $u^r \p_r$ terms in~\eqref{eq:euler2}.

\item The self-similar square-root of pseudo entropy $B$ is given by
\[
b(r,t) 
\approx (-t)^\cb  B\bigl(\tfrac{r}{ (-t)^{\cx}},- \log(-t)\bigr) 
\quad \Leftrightarrow \quad 
b(r,t) \approx e^{- \tau \cb} B(R,\tau), 
\]
where $\cb\in \Reals$ is an \emph{entropy similarity exponent}.
\end{itemize} 

Using the above notation, the Euler equation in radial symmetry~\eqref{eq:euler2} may be rewritten as 
\begin{subequations} 
\label{eq:euler3:old}
\begin{align} 
\p_\tau \Sigma 
-  \cu \Sigma 
+ \cx R \p_R \Sigma 
+   U \p_R \Sigma  
+ \alpha    \Sigma ( \p_R U + \tfrac{d-1}{R} U)
& \approx 0 ,  \\
\p_\tau U 
- \cu U 
+  \cx    R \p_R U 
+   U \p_R  U  
+ \alpha     \Sigma \p_R \Sigma 
- \tfrac{\alpha}{\gamma}  \Sigma^2 \tfrac{\p_R B}{B}
& \approx 0  ,  \\
\p_\tau B 
 - \cb  B 
+   \cx  R \p_R B 
+  U \p_R  B
&\approx 0 .
\end{align} 
\end{subequations} 
The self-similar evolution~\eqref{eq:euler3:old} with $(R,\tau) \in [0,\infty) \times [0,\infty)$ is equivalent to the Eulerian evolution~\eqref{eq:euler2} with $(r,t) \in [0,\infty) \times [-1,0)$, via the above transformations. We refer to the functions $(\Sigma,U,B)$ as \emph{self-similar profiles}.

\begin{remark}[\bf Modulation functions]
\label{rem:modulations}
At a first pass, the reader may regard  $\tau,\cx,\cu,\cb$ as being constant in $t$, and think of the $\approx$ symbols as $=$ signs. The fact that $\tau,\cx,\cu,\cb$ are in fact functions of time only plays a role in the stability analysis of Sections~\ref{sec:stability} and~\ref{sec:nonradial}; there, these functions of time are finely tuned to ``mod-out'' apparent instabilities that correspond to the scaling symmetries of the Euler system. 
\end{remark}

\begin{remark}[\bf Self-similar density, pressure, energy, and pseudo-entropy]
\label{rem:SS:density:pressure}
The self-similar profiles and scaling for density $\rho$, sound speed $c$, pressure $p$ (hence also of total energy $E$ and pseudo-entropy $e^s$), may be deduced from the above definitions; they  are given by 
\begin{align*}
\rho(r,t) 
&\approx (-t)^{\frac{\cu-\cb}{\alpha}}
\left( (\alpha \Sigma)^{\frac{1}{\alpha}} B^{-\frac{1}{\alpha}}\right)(R,\tau) 
,\\
c(r,t)
&\approx
(-t)^{\cu}
(\alpha \Sigma)(R,\tau)
,\\
p(r,t) 
&\approx
 (-t)^{\frac{\gamma \cu-\cb}{\alpha}} \left( \tfrac{1}{\gamma}  (\alpha \Sigma)^{\frac{\gamma}{\alpha}} B^{-\frac{1}{\alpha}}\right) (R,\tau) 
,\\
E(r,t)
&\approx
 (-t)^{\frac{\gamma \cu-\cb}{\alpha}}
\left(\tfrac{1}{2\alpha\gamma}  (\alpha \Sigma)^{\frac{\gamma}{\alpha}} B^{-\frac{1}{\alpha}} + \tfrac 12 (\alpha \Sigma)^{\frac{1}{\alpha}} U^2 B^{-\frac{1}{\alpha}} \right)(R,\tau) 
,\\
e^s(r,t)
&=
\tfrac{\gamma p}{\rho^\gamma}(r,t)
\approx 
(-t)^{2\cb} B^2(R,\tau),
\end{align*}
with $\cu \approx \cx - 1$.
\end{remark}

\subsection{Main results: existence of smooth globally self-similar Euler implosions}
\label{sec:intro:existence}
The following theorem is an \emph{abbreviated form} of Theorem~\ref{thm:main:profiles}, which is the main result of Section~\ref{sec:profiles}; it establishes the existence of a countable family of smooth, globally self-similar implosion profiles for the \emph{non-isentropic} compressible Euler equations, with similarity exponents  given by explicit closed-form expressions.

\begin{theorem}[\bf Smooth implosion profiles---abbreviated version]
\label{thm:intro:profiles}
Let $d \in \{1,2,3\}$ and $1 < \gamma \leq 2d+1$ be arbitrary. For each integer $\NNN \geq 1$, use the explicit closed form expressions in~\eqref{eq:cx:admissible}--\eqref{eq:cb:admissible} to define the similarity exponents $\cxbar = \cxstar(d,\gamma,\NNN)$,  $\cubar = \cxbar - 1$, and $\cbbar = \cbstar(d,\gamma,\NNN)$.\footnote{Crucially, $\cbbar \neq 0$, so that the entropy is non-trivial: these self-similar implosions are genuinely non-isentropic.} 
Then, there exist exact globally self-similar solutions $(\bar \sigma, \bar u^r, \bar b)(r,t)$ to the radially symmetric Euler equations~\eqref{eq:euler2}, given in terms of unique smooth similarity profiles $(\bar \Sigma, \bar U, \bar B)(R)$, via the transformation described in Section~\ref{sec:self:similar:intro}. The profiles $(\bar \Sigma, \bar U, \bar B)$ solve the system~\eqref{eq:euler4}, they are real-analytic functions of $R \in [0,\infty)$, they have explicit asymptotics as $R \to 0^+$ (cf.~\eqref{eq:values:at:zero:and:cb}, \eqref{eq:V:Q:R=0}, and~\eqref{eq:bar:H:R=0}),\footnote{It is precisely in the asymptotic behavior at $R=0$ where the parameter $\NNN$ plays a fundamental role. For any given $\NNN\geq 1$, the profiles $(\bar \Sigma, \bar U, \bar B)$ are constructed such that $\bar \Sigma(R) = \bar q_0 R + \OO(R^{2\NNN+1})$, $\bar U(R) = \bar v_0 R + \OO(R^{2\NNN+1})$, and $\bar B(R) = R + \OO(R^{2\NNN+1})$ as $R\to 0^+$, for suitable constants $\bar q_0>0$ and $\bar v_0<0$.} and they satisfy power-law behavior as $R \to \infty$ (cf.~\eqref{eq:V:Q:R=infty} and~\eqref{eq:bar:H:R=infty}). We moreover have $\bar \Sigma(R) >0$ and $\bar B(R) >0$ for all $R > 0$, and $\bar \Sigma(R)/R$ is strictly monotone decreasing in $R$. Finally, and crucially for the nonlinear stability analysis carried out in Sections~\ref{sec:stability} and~\ref{sec:nonradial}, the three self-similar wave speeds in the system obey a strictly positive lower bound which is uniform in $R\in[0,\infty)$ (cf.~\eqref{eq:Omega:invariant:c}); we refer to this property as the \emph{global outgoing property}.
\end{theorem}

The family of globally self-similar solutions to~\eqref{eq:euler2} constructed in Theorem~\ref{thm:intro:profiles} naturally defines (using the identifications in Remark~\ref{rem:SS:density:pressure}) globally self-similar solutions  $(\bar \rho,\bar \uu,\bar p)$ of the Euler equations  in terms of primary flow variables~\eqref{eq:euler:primary}; precisely, these are defined as
\begin{subequations}
\label{eq:globally:SS:exact:primary}
\begin{align}
\bar \rho(x,t) &:=  (-t)^{\frac{\cxbar -1-\cbbar}{\alpha}}\Bigl( (\alpha \bar \Sigma)^{\frac{1}{\alpha}}  \bar B^{- \frac{1}{\alpha}}\Bigr)  \Bigl(\tfrac{|x|}{(-t)^{\cxbar}}\Bigr) ,
\label{eq:globally:SS:exact:primary:a}
\\
\bar \uu(x,t) &:= (-t)^{\cxbar -1} \tfrac{x}{|x|} \bar U\Bigl(\tfrac{|x|}{(-t)^{\cxbar}}\Bigr) , 
\label{eq:globally:SS:exact:primary:b}
\\
\label{eq:globally:SS:exact:primary:c}
\bar p(x,t) &:= (-t)^{\frac{\gamma(\cxbar -1) - \cbbar}{\alpha}} \Bigl(\tfrac{1}{\gamma} (\alpha \bar \Sigma)^{\frac{\gamma}{\alpha}} 
\bar B^{- \frac{1}{\alpha}}\Bigr)\Bigl(\tfrac{|x|}{(-t)^{\cxbar}} \Bigr).
\end{align}
\end{subequations} 
We will show in Section~\ref{sec:profiles} that $(\alpha \bar \Sigma)/\bar B$ is a strictly positive function of $R \in [0,\infty)$, with limiting value $\frac{1}{1+\alpha d} \sqrt{\frac{\alpha \gamma d}{2}} > 0$ as $R\to 0^+$; thus, the density $\bar \rho$ is smooth and strictly positive on $\Reals^d \times [-1,0)$. Similarly, note that $\bar U(R) = R \bar V(R)$, where $\bar V$ is a smooth function of $R^2$  with $\bar V(0) = - \frac{1}{1+\alpha d} < 0$; thus, the velocity vector $\bar \uu$ is smooth on $\Reals^d \times [-1,0)$ and satisfies~\eqref{eq:implosion:big:def:b}. Moreover, for $0<R\leq 1$ we have that $ (\bar \Sigma ^{\gamma/\alpha}  \bar B^{-1/\alpha})(R) = R^2 \times$a smooth function of $R^2$, which is uniformly bounded from below (and above); since $R^2 = |x|^2 (-t)^{-2 \cxbar}$ and $|x|^2$ is a smooth function of $x$, we deduce that the pressure $\bar p$ is smooth on $\Reals^d \times [-1,0)$, and for each $t\in [-1,0)$ we have $\bar p(0,t)=0$ and $\bar p(x,t) >0$ for $|x|>0$. Hence, condition~\eqref{eq:implosion:big:def:c} also holds. Lastly, we note that the explicit expressions for $\cxbar$ and $\cbbar$ in~\eqref{eq:cx:admissible}--\eqref{eq:cb:admissible} give
\[
\tfrac{1}{\alpha} \bigl( \cxbar -1-\cbbar \bigr) = \tfrac{1}{\alpha} \bigl( \tfrac{1}{1+\alpha d} -1  \bigr) = - \tfrac{d}{1+\alpha d} < 0.
\]
Therefore, the power of $(-t)$ appearing in  ~\eqref{eq:globally:SS:exact:primary:a} is negative, showing that no matter what $d$, $\gamma$, or $\NNN$ we consider, we must have $\bar \rho(0,t) \to +\infty$ as $t \to 0^-$, and thus~\eqref{eq:implosion:big:def:a} holds. 

In summary, for any $d\in \{1,2,3\}$, $1<\gamma\leq 2d+1$ and $\NNN\in\Naturals$, the radially symmetric, smooth (real-analytic), globally self-similar solution $(\bar \rho,\bar \uu,\bar p)$ of the Euler system~\eqref{eq:euler:primary} (cf.~\eqref{eq:globally:SS:exact:primary}),  experiences an implosion singularity due to density blow up at $x_*=0$ and $t_*=0$, in the sense of Definition~\ref{def:implosion}.

\subsubsection{The limiting fields at the time of implosion}
\label{sec:limiting:fields}
For the solution $(\bar \rho,\bar \uu,\bar p)$ defined in~\eqref{eq:globally:SS:exact:primary}, not only do we know that it is an implosion solution of Euler~\eqref{eq:euler:primary} with implosion center $x_*=0$ and implosion time $t_*=0$; we may explicitly determine the limiting values of these primary variables at the time of implosion, namely $(\bar \rho,\bar \uu,\bar p)(r,0^-)$ for all $r>0$.

The argument is as follows. For every fixed $r>0$, the self-similar space variable $R = r / (-t)^{\cxbar} \to \infty$ as $t \to 0^-$. Moreover, in item (iii) of Theorem~\ref{thm:main:profiles} (cf.~Proposition~\ref{prop:power:series:infinity} and~\eqref{eq:bar:H:R=infty}) we prove that the self-similar profiles $(\bar \Sigma, \bar U, \bar B)$ satisfy the asymptotic power law behavior:
\[
\lim_{R\to \infty} R^{\frac{1}{\cxbar} -1} \bar \Sigma(R) = \underline{q}_1,
\quad
\lim_{R\to \infty} R^{\frac{1}{\cxbar} -1} \bar U(R) = \underline{v}_1,
\quad 
\lim_{R\to \infty} R^{- \frac{\cbbar}{\cxbar} } \bar H(R) = \underline{h}_1,
\]
for some constants $\underline{q}_1, \underline{h}_1 \in (0,\infty)$ and $\underline{v}_1 \in \Reals$ which depend on $\gamma, d$, and $\NNN$. In fact, for the ground state profile at $\NNN=1$, in Proposition~\ref{prop:sign:v1} we prove that $\underline{v}_1 < 0$.

Using the two facts mentioned above, and by appealing to the formulas in~\eqref{eq:globally:SS:exact:primary}, it then immediately follows that for \emph{any} $r>0$, we have the pointwise limits
\begin{subequations} 
\label{eq:fields:at:time:of:implosion}
\begin{align}
&\lim_{t\to 0^-} \bigl(\bar u^r,\bar \sigma,\bar b\bigr)(r,t) = \bigl( \underline{v}_1 r^{1-\frac{1}{\cxbar}}, \underline{q}_1 r^{1-\frac{1}{\cxbar}}, \underline{h}_1 r^{\frac{\cbbar}{\cxbar}}\bigr),
\label{eq:fields:at:time:of:implosion:a} \\
&\lim_{t\to 0^-} \bar \rho(r,t) = \bigl(\tfrac{\alpha \underline{q}_1}{\underline{h}_1} \bigr)^{\frac{1}{\alpha}} r^{\frac{\cxbar - 1 - \cbbar}{\alpha \cxbar}} 
= \bigl(\tfrac{\alpha \underline{q}_1}{\underline{h}_1} \bigr)^{\frac{1}{\alpha}} r^{-\frac{d}{(1+\alpha d) \cxbar}}
\label{eq:fields:at:time:of:implosion:b}
, \\
&\lim_{t\to 0^-} \bar p(r,t) = \tfrac{1}{\gamma} (\alpha \underline{q}_1)^{\frac{\gamma}{\alpha}}  \underline{h}_1^{-\frac{1}{\alpha}}   r^{\frac{\gamma (\cxbar - 1) - \cbbar}{\alpha \cxbar}},
\label{eq:fields:at:time:of:implosion:d} \\
&\lim_{t\to 0^-} \bar E(r,t) = \tfrac 12 \bigl( \tfrac{\alpha}{ \gamma}  \underline{q}_1^{2}     
+   \underline{v}_1^2 \bigr) ( \alpha \underline{q}_1)^{\frac{1}{\alpha}}  \underline{h}_1^{-\frac{1}{\alpha}} r^{\frac{\gamma (\cxbar - 1) - \cbbar}{\alpha \cxbar}} 
\label{eq:fields:at:time:of:implosion:f} .
\end{align}
\end{subequations}
In particular, \eqref{eq:fields:at:time:of:implosion:b} reiterates our previously established fact: for any $d,\gamma,\NNN$, the density $\bar \rho$ blows up at $(0,0)$. 
Moreover, since $d > \frac{d}{(1+\alpha d) \cxbar} \Leftrightarrow \cxbar > \tfrac{1}{1+\alpha d}$, which is a direct consequence of~\eqref{eq:cx:admissible}, we have that $\lim_{t\to 0^-} \int_0^1 \bar \rho(r,t) r^{d-1} dr < \infty$, and so \emph{at the time of implosion we have locally finite mass}. 

Similarly, since $d > \frac{\cbbar - \gamma (\cxbar - 1) }{\alpha \cxbar}
\Leftrightarrow \alpha d \cxbar + \gamma (\cxbar - 1) > \cxbar - \frac{1}{1+\alpha d} \Leftrightarrow \cxbar  > \frac{d + 2+ 2\alpha d}{(d+2)(1+\alpha d)} $, which holds for all $\alpha \in(0,d]$ by using ~\eqref{eq:cx:admissible}--\eqref{eq:cb:admissible} and Lemma~\ref{lem:properties:exponents}, we have $\lim_{t\to 0^-} \int_0^1 \bar E(r,t) r^{d-1} dr < \infty$, and hence \emph{at the time of implosion we have locally finite energy}. 

Lastly, since $d> \frac{d}{(1+ \alpha d) \cxbar} + \frac{1}{\cxbar} - 1 \Leftrightarrow \cxbar  > \frac{d+1+\alpha d}{(d+1)(1+ \alpha d)}$, which holds for all $\alpha \in(0,d]$ by using ~\eqref{eq:cx:admissible}--\eqref{eq:cb:admissible} and Lemma~\ref{lem:properties:exponents}, we deduce that $\lim_{t\to 0^-} \int_0^1( \bar \rho \bar u^r)(r,t) r^{d-1} dr < \infty$, and hence \emph{at the time of implosion we have locally finite momentum}. 

As discussed by Jenssen and Tsikkou~\cite{JenssenTsikkou2018} in the context of the Guderley problem, the fact that at the time of implosion the solutions in~\eqref{eq:fields:at:time:of:implosion} have locally finite mass, momentum, and energy, implies that we may hope to continue the Euler solution \emph{past the time of the implosion} as a weak solution of the full Euler system~\eqref{eq:euler0}. This continuation-post-implosion problem, sometimes called the ``reflection problem'', will be analyzed in a subsequent work.

\subsubsection{Values of similarity exponents and the regularity of kinetic equations}
\label{sec:intro:kinetic}
A long-standing open problem in kinetic theory is whether the Landau or Boltzmann equation admits global-in-time smooth solutions. In the spatially homogeneous setting, global regularity has been established for the Landau equation with soft potentials (including the Coulomb case $\gamkin = -3$) by Guillen and Silvestre \cite{GuillenSilvestreLandau}, and for the Boltzmann equation by Imbert, Silvestre, and Villani \cite{imbert2026monotonicity}. The regularity problem is much less understood in the inhomogeneous setting. In this setting, kinetic equations are closely related to the compressible Euler equations through the hydrodynamic limit and the Hilbert expansion \cite{caflisch1980fluid}. In view of this connection, and following the work \cite{MRRS2022a,MRRS2022b} on smooth implosions, imploding singularities provide a potential scenario for singularity formation in kinetic equations.

Very recently, Bedrossian, Chen, Gualdani, Ji, Vicol, and Yang~\cite{bedrossian2026finite} constructed finite-time singularities for the inhomogeneous Landau equation with very hard potentials \(\gamkin \in (\sqrt{3},2]\) from smooth, strictly positive initial data by proving nonlinear finite co-dimension stability of the local Maxwellian associated with the Euler imploding profile in self-similar variables.\footnote{\label{foot:Maxwell} In original $(x,t;v)$ variables, the local Maxwellian is $\mathcal{M}_{\bar \rho,\bar \uu,\bar \theta}(x,t;v) = \frac{\bar \rho(x,t)}{(2\pi \bar\theta(x,t))^{3/2}} \exp\bigl(- \frac{|v - \bar \uu(x,t)|^2}{2\bar \theta(x,t)} \bigr)$, where $\bar \theta = \bar{\rho}^{-1} \bar p$ is the temperature. Using~\eqref{eq:globally:SS:exact:primary}, and upon denoting $(\mathcal{X},\mathcal{V}) = (x/(-t)^{\cxbar},v/(-t)^{\cubar})$ and re-normalizing the profiles from Theorem~\ref{thm:intro:profiles} as $(\bar \Sigma,\bar U,\bar B)(R) = R \cdot (\bar Q,\bar V,\bar H)(R)$, this self-similar local Maxwellian becomes 
\[
\mathcal{M}_{\bar \rho,\bar \uu,\bar \theta}(x,t;v) = (-t)^{-3 \cbbar}\mathcal{M}(\mathcal{X};\mathcal{V}), 
\quad \mbox{where} \quad
\mathcal{M}(\mathcal{X};\mathcal{V}) :=  \frac{1}{c_0 |\mathcal{X}|^{3} \bar H^{3}(|\mathcal{X}|)} \exp\Bigl(- \frac{15 |\mathcal{V} -  \mathcal{X} \bar V(|\mathcal{X}|)|^2}{2 |\mathcal{X}|^2 \bar Q^2(\mathcal{X})} \Bigr),
\quad \! c_0 = (\frac{6\pi}{5})^{3/2}.
\]
} 
This work provides the first example of a collisional kinetic model that is globally well-posed in the homogeneous setting but admits finite-time singularities for inhomogeneous data. To carry out this lifting from the Euler level to the kinetic level, one must use an imploding profile specific to a 3D~\emph{monatomic gas} (\(\gamma = \tfrac{5}{3}\)). In addition, the authors of~\cite{bedrossian2026finite} identified the following condition for the \emph{hydrodynamic limit} in self-similar variables:\footnote{
The similarity exponents $(c_f, c_v)$ used in \cite[Section 2.2]{bedrossian2026finite} correspond to $(- 3 \cbbar ,  \cubar)$ in \eqref{eq:hydro}, respectively.
}
\begin{equation}\label{eq:hydro}
 - 3  \cbbar  + (\gamkin+3) \bcu + 1 < 0,
\end{equation}
where \(\cbbar, \bcu\) are the similarity exponents constructed in Theorem~\ref{thm:intro:profiles}, and \(\gamkin\) is the parameter of the collision kernel. There are a few important parameter values: 
$\gamkin = 1$  for  Boltzmann with hard-sphere potential, $\gamkin = 0$ for Maxwell molecules, and $\gamkin = -3$ for Landau with Coulomb potential.

In the isentropic setting, the pseudo-entropy profile is constant, \(\bar B \equiv \text{const} \neq 0\), and \(\cbbar = 0\). Within the class of smooth radial isentropic imploding profiles for the Euler equations with a monatomic gas, Shao, Wang, Wei, and Zhang~\cite{ShaoWangWeiZhang2025} constructed profiles with asymptotically minimal similarity exponents, with $\bcu \to (\frac{\sqrt{3}-3}{6})^{+}$. 
As a result, it follows from~\eqref{eq:hydro} that lifting such isentropic Euler profiles to the kinetic level requires \(\gamkin > \sqrt{3}\). 
To reach smaller values of \(\gamkin\) using radially symmetric profiles, one \emph{must introduce a nontrivial entropy}, and one has to construct Euler implosion profiles with exponents compatible with~\eqref{eq:hydro}.

By introducing nontrivial entropy, the profiles for a 3D monatomic gas (\(\gamma = \tfrac{5}{3}\)) constructed in Theorem~\ref{thm:intro:profiles} significantly relax the constraint on \(\gamkin\). For $d=3$ and $\gamma = \tfrac{5}{3}$, Example~\ref{ex:exponents:monatomic:3D} gives $\cbbar = \cxbar - \tfrac12$ and $\cubar = \cxbar - 1$, so that~\eqref{eq:hydro} reduces to
\begin{equation}\label{eq:hydro2}
  \gamkin (\cxbar - 1) - \tfrac12 < 0.
\end{equation}
Since $\cxstar(3,\tfrac{5}{3},1) = \frac 12 + \frac 12 \sqrt{\frac{47}{12}} \approx 1.48953$ and $\cxstar(3,\tfrac{5}{3},\infty) = \frac 12 + \frac 12 \sqrt{\frac{5}{6}} \approx 0.956435$ (cf.~Example~\ref{ex:exponents:monatomic:3D}), for \emph{any} $\gamkin \in [-3, 2]$ one can choose an integer $\NNN = \NNN(\gamkin) \geq 1$ such that the profile with similarity exponent $\cxbar := \cxstar(3,\tfrac{5}{3},\NNN)$ satisfies~\eqref{eq:hydro2}. In particular, the three physically important cases $\gamkin \in \{-3, 0, 1\}$ are already covered by the ground state profile at $\NNN = 1$. Consequently, the Euler implosion profiles constructed in this work provide candidates for lifting implosions to the kinetic level, to construct finite time singularities for the Landau/Boltzmann equations at much smaller values of $\gamkin$.

We also point out that, due to the vanishing of the temperature at the origin (see Footnote \ref{foot:Maxwell}), the Maxwellians associated with these new profiles are singular at the spatial origin ($x=0$), where the local Maxwellians become Dirac masses (in $v$). 
Thus, any potential lift of the Euler implosion profiles constructed in this work to the kinetic level would require substantial new ideas beyond the framework of~\cite{bedrossian2026finite}, in order to overcome this one-point singularity. Likely, these new ideas would need to take into account the regularizing effect of the collision term.
Whether the full nonlinear stability analysis of~\cite{bedrossian2026finite} can be extended to these non-isentropic profiles is a separate question, not addressed here.

We also want to highlight a recent work of Golding and Henderson~\cite{GoldingHenderson25}, which establishes an $L^1_t L^\infty_{x,v}$ continuation criterion for the Landau equation with Coulomb potential $\gamkin = -3$. For an Euler imploding singularity to potentially be lifted to the kinetic (Landau) level when $\gamkin = -3$, the result of~\cite{GoldingHenderson25} suggests that the $L^1_t L^\infty_{x,v}$ norm of the local Maxwellian associated with the Euler imploding solution must become unbounded. For the profiles of a 3D monatomic gas $\left(\gamma = \frac53\right)$ constructed in Theorem~\ref{thm:intro:profiles}, the amplitude of the associated local Maxwellian blows up at the rate $(-t)^{-3\cbbar}$ as $t\to 0^-$; see Footnote~\ref{foot:Maxwell}. Moreover, for all such profiles, we have $\cbbar = \cxbar - \tfrac12 > \cxstar(3,\frac 53,\infty) - \frac 12 \approx 0.456435 > \frac13$; thus,  the time integral $\int_{-1}^{t}(-t')^{-3\cbbar} dt'$ diverges as $t\to 0^-$. Therefore, these profiles could potentially be useful for studying blowup in the Landau equation with Coulomb potential.

\subsection{Main results: stability of the globally self-similar implosions}
\label{sec:intro:stability}
The bulk of the paper is dedicated to the PDE stability of the globally self-similar solutions $(\bar \rho,\bar \uu, \bar p)$ constructed in~\eqref{eq:globally:SS:exact:primary}, within the class of smooth solutions of the full Euler equations~\eqref{eq:euler:primary}. 

We emphasize that for Guderley's imploding shocks, a \emph{rigorous} nonlinear stability analysis is currently \emph{not available}. The only available results deal with mode/spectral analysis for the linearized problem---starting with Morawetz's PhD thesis~\cite{Morawetz1951} for linear stability in radial symmetry, and the work of Brushlinskii and Kazhdan~\cite{BrushlinskiiKazhdan1963} for linear instability outside of radial symmetry---or consider arguments which are semi-analytical, a combination of asymptotic arguments and supporting numerical computations~\cite{GardnerBrookBernstein1982,ChenZhangPanarella1995,WuRoberts1996PRE,WuRoberts1996}. 

For the $C^\infty$-smooth isentropic of Merle, Rapha\"el, Rodnianski, and Szeftel~\cite{MRRS2022a}, rigorous stability results are available. Nonlinear stability for perturbations which have \emph{radial symmetry} and which lie in an \emph{unquantified} finite co-dimension subspace of a weighted Sobolev space, was established in~\cite{MRRS2022b,MRRS2022}, and subsequently refined in~\cite{BCG2025,CCSV2024,Chen2024}.\footnote{The numerical work of Biasi~\cite{Biasi2021} makes interesting conjectures about the nature of radially symmetric instabilities.} Nonlinear stability outside of radial symmetry, but still for perturbations which lie in an \emph{unquantified} finite co-dimension subspace of a weighted Sobolev space was established in~\cite{CGSS2025}. 

The results in this paper improve upon the previously available stability theory in several ways. In radial symmetry, we establish the \emph{complete nonlinear stability picture} for small smooth perturbations of the globally self-similar solutions $(\bar\rho,\bar\uu,\bar p)$ constructed in~\eqref{eq:globally:SS:exact:primary}. Outside of radial symmetry, we establish the complete nonlinear stability picture for the ground state ($\NNN=1$) at the adiabatic exponents corresponding to a monatomic or diatomic gas; for general $\gamma$ and $\NNN$, we prove finite co-dimension stability, give an explicit upper bound on the dimension of the unstable manifold, and characterize this manifold in terms of the unstable modes of an explicit, finite-dimensional ODE system. The mechanism that enables these results---and which distinguishes our analysis from the prior work~\cite{MRRS2022b,MRRS2022,BCG2025,CCSV2024,Chen2024,CGSS2025}---is the \emph{global outgoing property} of the self-similar profiles asserted in Theorem~\ref{thm:intro:profiles}. Our main results are as follows.

\subsubsection{Stability in radial symmetry}
Let $\NNN \geq 1$. In radial symmetry, the $\NNN^{\rm th}$ globally self-similar solution $(\bar \rho, \bar u^r,\bar p)$ constructed in~\eqref{eq:globally:SS:exact:primary} (from the $\NNN^{\rm th}$ profile in Theorem~\ref{thm:intro:profiles}) is nonlinearly stable (to small and smooth perturbations) within the class of smooth solutions $(\rho,u^r,p)$ of~\eqref{eq:euler2:primary} for which the initial pressure vanishes at $r=0$, and for which the initial data satisfies $3(\NNN-1)$ explicit compatibility conditions at $r=0$. In particular, for the ground state solution at $\NNN=1$, no compatibility conditions are required beyond $p(0,-1)=0$.
The detailed statement is given in Theorem~\ref{thm:stability:radial}, which is stated in terms of self-similar coordinates and variables. The implosion blowup in physical variables is recorded in Corollary~\ref{cor:radial:blowup:physical}. Here we give an informal summary of the results in Section~\ref{sec:stability}.

\begin{theorem}[\bf Radial stability---abbreviated version]
\label{thm:intro:radial:stab}
Fix $d \in \{1,2,3\}$, $1<\gamma\leq 2d+1$, and $\NNN \geq 1$. Fix the similarity exponents $\cxbar, \cbbar$, and the self-similar profiles $(\bar \Sigma, \bar U, \bar B)$ from Theorem~\ref{thm:intro:profiles}.  Then, there exists $\underline{\eps}^* = \underline{\eps}^*(d,\gamma,\NNN)> 0$ such that for any $0< \underline{\eps} \leq \underline{\eps}^*$,  the following holds. 

Consider initial data $(\rho_{\sf in},\uu_{\sf in},p_{\sf in})(x)$ for the Euler system in terms of primary variables~\eqref{eq:euler:primary}. Assume that this initial datum is radially symmetric and  belongs to $C^{17 \NNN+1}(\Reals^d)$.\footnote{Radial symmetry and regularity in $x$ near $x=0$ imply that $\rho_{\sf in}$ and $p_{\sf in}$ are ``even'' functions of $r=|x|$ (meaning, they depend on $r^2$) and that $u^r_{\sf in}$ is an ``odd'' function of $r=|x|$ (meaning that it is given by $r$ times a function which depends on $r^2$).} Denote the initial perturbation from the exact self-similar solution $(\bar \rho,\bar \uu, \bar p)$ defined in  \eqref{eq:globally:SS:exact:primary}, at time $t=-1$,  by 
\begin{equation*}
(\tilde \rho_{\sf in},\tilde u^r_{\sf in},\tilde p_{\sf in})(r)
= \bigl(\rho_{\sf in}(x) - \bar \rho(x,-1), \tfrac{x}{|x|} \cdot (\uu_{\sf in}(x) - \bar \uu(x,-1) ), p_{\sf in}(x) - \bar p(x,-1) \bigr),
\end{equation*}
where $r=|x|$. Assume that $\tilde p_{\sf in}(0) = 0$, and for $\NNN\geq 2$ assume the following $3(\NNN-1)$ compatibility conditions on the initial density, velocity, and pressure:\footnote{This is a rewriting of the compatibility condition in item (ii) of Theorem~\ref{thm:stability:radial}. Recalling that $\rho, u^r/r$, and $p$ are functions of $r^2$ (due to regularity in $x\in\Reals^d$), this compatibility condition together with the properties of the profiles $(\bar \Sigma,\bar U,\bar B)$ imply that $\rho_{\sf in} = \rho_0 + \OO(r^{2\NNN})$, $u^r_{\sf in} = u_1 r + \OO(r^{2\NNN+1})$, and $p_{\sf in} = p_2 r^2 + \OO(r^{2\NNN+2})$ for $r\ll 1$, for suitable constants $\rho_0,u_1,p_2$.}
\begin{equation}
\label{eq:thm:intro:radial}
\p_r^{2k+2} \tilde \rho_{\sf in}(0)
= 
\p_r^{2k+3} \tilde u^r_{\sf in}(0)
=
\p_r^{2k+4} \tilde p_{\sf in}(0)
=
0,
\qquad
\mbox{for all}
\qquad
0 \leq k \leq \NNN-2.
\end{equation}
Furthermore, assume that the Taylor coefficients of $(\tilde \rho_{\sf in},\p_r \tilde u^r_{\sf in}, \p_{rr} \tilde p_{\sf in})$ at $r=0$, of order $\leq 17 \NNN$, are small in terms of $\underline{\eps}$ (cf.~\eqref{eq:thm:init:Taylor}), and assume that certain weighted $L^\infty$-norms of the quantities $\p_r \tilde u^r_{\sf in} - \frac{1}{r} \tilde u^r_{\sf in}$, $\bar \rho(\cdot,-1)^{-1} r \p_r \tilde \rho_{\sf in}$, and $\bar p(\cdot,-1)^{-1} r \p_r \tilde p_{\sf in}$, are small with respect to $\underline{\eps}$ (cf.~\eqref{eq:thm:init:JM}).

Then, there exists a time $t_* = \OO(\underline{\eps}^{\frac 45})$ and a unique, smooth, radially symmetric, asymptotically self-similar implosion solution $(\rho,\uu,p)$ of the Euler equations~\eqref{eq:euler:primary} (equivalently,~\eqref{eq:euler2:primary}) on $\Reals^d \times [-1,t_*)$, with implosion center $x_*=0$, and implosion time $t_*$. 

In self-similar coordinates  (see~Section~\ref{sec:self:similar:intro}), the corresponding fields $(\Sigma,U,B)(R,\tau)$ which solve~\eqref{eq:euler3:old}, exist globally in time $\tau$, remain as smooth as the initial data, and converge exponentially fast (in a suitable topology) to the stationary profiles $(\bar \Sigma,\bar U,\bar B)(R)$ found in Theorem~\ref{thm:intro:profiles}; see~\eqref{eq:thm:conc:Taylor}, \eqref{eq:thm:conc:JM}, and item (d) of Theorem~\ref{thm:stability:radial}. The associated modulation functions $(\cx,\cu,\cb)(\tau)$ converge exponentially fast in $\tau$ to their stationary values $(\cxbar,\cubar,\cbbar)$.
\end{theorem}

\begin{remark}[\bf Ground state is stable to radially symmetric perturbations]
\label{rem:intro:ground:state}
For $\NNN = 1$, condition~\eqref{eq:thm:intro:radial} is \emph{empty}. Since the $\underline{\eps}$-smallness conditions mentioned in Theorem~\ref{thm:intro:radial:stab} form an open set in a suitable topology, we deduce that in radial symmetry the stability of the ground state self-similar profile/globally self-similar solution~\eqref{eq:globally:SS:exact:primary} holds for all permissible (meaning, $p_{\sf in}(0)=0$) initial perturbations in an \emph{open set}.  See also Remark~\ref{rem:ground:state:stable}.
\end{remark}

\subsubsection{Stability outside of radial symmetry}
\label{sec:intro:stability:nonradial}
Let $\NNN \geq 1$. We next consider smooth perturbations of the $\NNN^{\rm th}$ globally self-similar solution $(\bar \rho,\bar \uu,\bar p)$ constructed in~\eqref{eq:globally:SS:exact:primary} which \emph{do not} obey any symmetry assumptions. Our main result is that, modulo finitely many explicit linear conditions on the initial data at the spatial point $x=0$, the globally self-similar implosion is nonlinearly stable, within the class of smooth solutions $(\rho,\uu,p)$ of the full Euler system~\eqref{eq:euler:primary}. More precisely, we prove \emph{finite co-dimension stability} with an \emph{explicit upper bound on the co-dimension}, which moreover becomes \emph{sharp} for the ground state ($\NNN=1$) of a monatomic or diatomic gas. The detailed statements are given in Theorem~\ref{thm:nonrad:stab} (PDE stability of the bulk perturbation) and Theorem~\ref{thm:ODE_stab} (finite-dimensional ODE analysis at the origin); here we record an informal summary.

\begin{theorem}[\bf Non-radial stability---abbreviated version]
\label{thm:intro:nonrad:stab}
Fix $d \in \{1,2,3\}$, $1<\gamma \leq 2d+1$, and $\NNN \geq 1$. Fix the similarity exponents $(\cxbar, \cubar, \cbbar) = (\cxbar, \cxbar-1, \cbbar)$ and the globally self-similar solution $(\bar \rho,\bar \uu,\bar p)$ of~\eqref{eq:euler:primary} defined in~\eqref{eq:globally:SS:exact:primary}. There exist an integer $\MMM_* = \MMM_*(d,\gamma,\NNN)\geq 1$, an explicitly computable constant $\bar C(d) \geq 1$, a linear subspace $\Sigma_{\mathsf{uns}} \subset \Reals^{d_{\leq \MMM_*}}$ with\footnote{For an integer $n\geq 1$, the notation $d_{\leq n}$ 
is used to denote the length of the vector which enumerates the elements of the set $\cup_{f \in \{ \rho, \nabla \uu, \nabla^2 p\}} \{\p^\bfa f \colon \bfa \in \Naturals_0^d, |\bfa| \leq n \}$; see Equation~\eqref{def:TL_nth_order} and Footnote~\ref{foot:d:leq:n} following it.}
\[
\dim(\Sigma_{\mathsf{uns}}) \leq \bar C(d)\,\NNN^{d},
\]
and a sufficiently small constant $\bar \eps_* = \bar \eps_*(d,\gamma,\NNN) > 0$, such that for any $0 < \bar\eps \leq \bar \eps_*$ the following holds.

Consider initial data $(\rho_{\mathsf{in}},\uu_{\mathsf{in}},p_{\mathsf{in}})(x)$ for the Euler system~\eqref{eq:euler:primary}, and denote the initial perturbation from the exact self-similar solution at time $t=-1$ by
\[
  (\tilde \rho_{\mathsf{in}},\tilde \uu_{\mathsf{in}},\tilde p_{\mathsf{in}})(x)
  =
  \bigl(\rho_{\mathsf{in}}(x) - \bar \rho(x,-1),\;
  \uu_{\mathsf{in}}(x)-\bar \uu(x,-1),\;
  p_{\mathsf{in}}(x) - \bar p(x,-1)\bigr).
\]
Assume that $(\rho_{\mathsf{in}},\uu_{\mathsf{in}},p_{\mathsf{in}}) \in C^{\MMM_*+3}(\Reals^d)$ 
is such that the density is strictly positive, the pressure is non-negative,  and attains a global minimum at a \emph{unique point in space}, where the pressure vanishes; without loss of generality\footnote{The Euler equations~\eqref{eq:euler:primary} are invariant under spatial translations.}  we may assume that this global minimum of the pressure is attained at $x=0$, namely
\begin{subequations}
\label{eq:thm:intro:nonrad:compat}
\begin{equation}
\label{eq:thm:intro:nonrad:compat:origin}
  p_{\mathsf{in}}(0) = 0 \,.
\end{equation}
Automatically\footnote{Since the pressure is non-negative, if $p_{\sf in}(0)=0$ then $p_{\sf in}$ attains a global minimum at the origin, and thus $\nabla p_{\sf in}(0)= 0$ too.} we also deduce that $\nabla p_{\mathsf{in}}(0) = 0$, and without loss of generality\footnote{The Euler equations~\eqref{eq:euler:primary} are invariant under Galilean boosts.} we have $\uu_{\mathsf{in}}(0) = 0$.

Define the vector of Taylor coefficients of the initial perturbation at $x=0$ up to order~$\MMM_*$ by
\[
  V_{\leq \MMM_*}(0) :=
  \bigl(
  \{\p^{\bfa}\tilde \rho_{\mathsf{in}}(0)\}_{|\bfa|\leq \MMM_*},
  \;
  \{\p^{\bfa}\tilde \uu_{\mathsf{in}}(0)\}_{1\leq |\bfa|\leq \MMM_*+1},
  \;
  \{\p^{\bfa}\tilde p_{\mathsf{in}}(0)\}_{2\leq |\bfa|\leq \MMM_*+2}
  \bigr)
  \subset
  \Reals^{d_{\leq \MMM_*}},
\]
and assume the \emph{finite co-dimension compatibility condition}
\begin{equation}
\label{eq:thm:intro:nonrad:compat:Sigma}
  V_{\leq \MMM_*}(0) \in \mathcal{M}_{\mathsf{stab}} \,,
\end{equation}
where $\mathcal{M}_{\mathsf{stab}}$ is the nonlinear stable manifold (at the origin in $\Reals^{d_{\leq \MMM_*}}$) of the closed ODE system governing the Taylor coefficients at $x=0$ (cf.~Theorem~\ref{thm:ODE_stab}, items~(ii)--(iii)) 
and is related to 
the linear subspace $\Sigma_{\mathsf{uns}}$. In particular, \eqref{eq:thm:intro:nonrad:compat:Sigma} amounts to 
$\dim(\Sigma_{\mathsf{uns}})$ scalar compatibility conditions, which is at most $\bar C(d)\NNN^{d}$. 

Finally, assume that $|V_{\leq \MMM_*}(0)| \leq \bar \eps$, and assume that the bulk part of the initial perturbation $(\tilde \rho_{\mathsf{in}}, \tilde \uu_{\mathsf{in}},\tilde p_{\mathsf{in}})$ (cf.~\eqref{eq:nonrad_var} and \eqref{eq:self-similar:new:variables}) is bounded in terms of $\bar \eps$, in a suitable weighted Sobolev norm  (cf.~\eqref{ass:nonrad_init}).
\end{subequations}

Then, there exists a time $t_* = \OO(\bar\eps^{\frac{4}{5}})$ and a unique, smooth, asymptotically self-similar implosion solution $(\rho,\uu,p)$ of the Euler equations~\eqref{eq:euler:primary} on $\Reals^d \times [-1, t_*)$, with implosion center $x_*=0$ and implosion time~$t_*$.

In self-similar coordinates (see Sections~\ref{sec:convenient:unknowns} and~\ref{sec:modulate:SS:ansatz:nonradial}), the corresponding fields $(\vrho,\UU,\mb)(y,\tau)$ which solve~\eqref{eq:sys}, are such that the perturbations $(\tilde \vrho,\tilde \UU,\tilde \mb)(y,\tau)$ from the stationary profile $(\bar \vrho,\bar \UU,\bar \mb)$ exist globally in $\tau$, remain as smooth as the initial data, and decay exponentially fast to $0$ in a suitable weighted energy norm; see~\eqref{eq:nonrad_decay} and Theorem~\ref{thm:nonrad:stab}. The associated modulation functions $(\cx,\cu,\cbb,\crho)(\tau)$ defined in~\eqref{eq:normal-1} converge exponentially fast in $\tau$ to their stationary values $(\cxbar,\cubar,\cbbbar,\crhobar)$.
\end{theorem}

\begin{remark}[\bf Example of an admissible set of initial data]
\label{rem:intro:nonrad:empty}
The assumptions on $(\tilde \rho_{\mathsf{in}},\tilde \uu_{\mathsf{in}},\tilde p_{\mathsf{in}})$ in Theorem~\ref{thm:intro:nonrad:stab} are \emph{open conditions} in the topology of a suitably weighted Sobolev space, except for~\eqref{eq:thm:intro:nonrad:compat:Sigma}, which induces a finite co-dimension constraint. Condition~\eqref{eq:thm:intro:nonrad:compat:Sigma} is a finite collection of at most $\bar C(d)\NNN^{d}$ linear equalities among the Taylor coefficients of the initial perturbation at $x=0$. 
A simple explicit subset of admissible data is obtained by requiring that the initial perturbation vanishes to sufficiently high order at the origin; if
\begin{equation}
\label{eq:thm:intro:nonrad:highvanish}
  \nabla^{\leq \MMM_*} \tilde \rho_{\mathsf{in}}(0) = 0,
  \qquad
  \nabla^{\leq \MMM_*+1} \tilde \uu_{\mathsf{in}}(0) = 0,
  \qquad
  \nabla^{\leq \MMM_*+2} \tilde p_{\mathsf{in}}(0) = 0,
\end{equation}
then~\eqref{eq:thm:intro:nonrad:compat:Sigma} holds trivially, and~\eqref{eq:thm:intro:nonrad:highvanish} is preserved by the Euler evolution. See~Remark~\ref{rem:data:for:nonrad:stability}. Large classes of initial data satisfying~\eqref{eq:thm:intro:nonrad:highvanish} and all the other assumptions of Theorem~\ref{thm:intro:nonrad:stab} may be constructed: these data have non-zero vorticity, non-trivial entropy variation, and may be taken to be constant outside a compact set; see~Remark~\ref{rem:compact:support}.
\end{remark}

\begin{remark}[\bf Complete stability picture for the monatomic/diatomic ground state]
\label{rem:intro:nonrad:ground:state}
When $\NNN=1$, item (iv) in Theorem~\ref{thm:ODE_stab} shows that the unstable subspace $\Sigma_{\mathsf{uns}}$ is the trivial lifting\footnote{By this we mean that if $V \in \Sigma_{\mathsf{uns}}$, then $V = (V_{\leq 2},0\ldots,0)$, with $V_{\leq 2} \in \Sigma_{\mathsf{uns},\leq 2}$.} of $\Sigma_{\mathsf{uns},\leq 2} \subset \Reals^{d_{\leq 2}}$, the unstable+center manifold associated to the ODE evolution of the Taylor coefficients at $x=0$ \emph{of order at most~$2$}. For the physically most relevant adiabatic exponents, the dimension of  the unstable+center manifold is explicit:
\begin{center}
\renewcommand{\arraystretch}{1.2}
\begin{tabular}{l@{\hspace{2em}}l@{\hspace{2em}}c c }
\hline
dimension $d$ & adiabatic exponent $\gamma$ & $\dim (\Sigma_{\mathsf{uns}})$ & $\dim(\Reals^{d_{\leq 2}})$\\
\hline
$d=3$ & monatomic: $\gamma = \tfrac{5}{3}$ & $11$ & 98 \\
$d=3$ & diatomic:   $\gamma = \tfrac{7}{5}$ & $18$ & 98 \\
$d=2$ & monatomic: $\gamma = 2$         & $5$ &  36 \\
$d=2$ & diatomic: $\gamma = \tfrac{5}{3}$  & $7$ & 36 \\
$d=1$ & any $\gamma \in (1,3]$           & $1$ &  9 \\
\hline
\end{tabular}
\end{center}
Moreover, in each of these cases, item (v) of Theorem~\ref{thm:ODE_stab} provides an \emph{explicit} class of admissible initial data, for which condition~\eqref{eq:thm:intro:nonrad:compat:origin} is automatically satisfied: it suffices that the Taylor coefficients up to order $2$ of the initial perturbation vanish; that is $V_{\leq 2}(0)=0$. The assumption that the higher-order Taylor coefficients $V_{\leq n}(0)$ are small, and that the bulk initial perturbation is small in a suitable weighted Sobolev space (cf.~\eqref{ass:nonrad_init}), are \emph{open} conditions.  To the best of our knowledge, this is the first complete nonlinear stability result for  imploding solutions of the multi-dimensional compressible Euler equations outside of symmetry considerations.
\end{remark}

\subsection{Related literature}
\label{sec:literature}
The works which are directly relevant\footnote{The mathematical literature on compressible flows is too vast to review here. For instance, we do not discuss the problem of shock formation and shock development for the compressible Euler equations. For the modern theory of hyperbolic systems of conservation laws, which, in one space dimension, provides a framework for shock formation, development, and interaction, we refer the interested reader to the book of Dafermos~\cite{dafermos2005hyberbolic}. For a detailed account of the literature on multi-dimensional shocks for the Euler equations, we refer the interested reader to the summary given by Shkoller and Vicol in~\cite{ShkollerVicol2024}.} to the present paper have already been discussed in the opening part of the Introduction. Here we complement that summary with a more complete account of the literature around implosion singularities: rigorous existence theorems that were not discussed above, stability analyses of Euler implosions, and self-similar imploding singularities in related PDEs.

\subsubsection{Rigorous mathematical treatment of Guderley's imploding shock}

For more than seven decades after its introduction, the Guderley imploding shock was studied mostly by semi-analytical and numerical means~\cite{Lazarus1977,Lazarus1981,MeyerTerVehn1982,Ramsey2012,Giron2023,Sedov2018}: the similarity exponent is determined by a shooting method from the sonic point and the resulting profile is continued numerically through the spacetime point of collapse.
A fully rigorous construction, as a genuine weak solution of the full (non-isentropic) compressible Euler system with the prescribed self-similar structure, was established only recently by Jang, Liu, and Schrecker~\cite{JangLiuSchrecker2025} for the full range of adiabatic exponents $\gamma \in (1,3]$. Shortly thereafter, Cialdea, Shkoller, and Vicol~\cite{CialdShkVic2025} established that the Guderley profile \emph{arises dynamically} from $C^{1,1/3}$-smooth, classical, shock-free initial data. In particular, the paper~\cite{CialdShkVic2025} provides the first explicit construction of solutions
to the full compressible Euler equations that evolve from classical, shock-free initial data into the \emph{strong
shock regime}. Complementary to these forward-in-time statements, Jenssen and Tsikkou~\cite{JenssenTsikkou2018} rigorously showed that Guderley's self-similar solution provides a weak solution of the Euler system across the implosion time, so that the diverging explosion wave emerging at $t = t_*^+$ is a genuine admissible continuation of the imploding phase.  

\subsubsection{Linear mode analysis of the Guderley implosion}

Well before the nonlinear PDE stability of imploding profiles came within reach, the question of stability was approached by linear mode analysis. For the Guderley imploding shock in radial symmetry, the earliest systematic treatment is that of Morawetz's PhD thesis under Friedrichs~\cite{Morawetz1951}, which addresses the eigenvalue problem arising from a linearization of the Euler equations around the self-similar Guderley profile. For non-radial perturbations, Brushlinskii and Kazhdan~\cite{BrushlinskiiKazhdan1963} identified unstable discrete eigenmodes associated with low spherical-harmonic sectors, and Wu and Roberts~\cite{WuRoberts1996,WuRoberts1996PRE} refined this analysis using spherical-harmonic decomposition and asymptotic expansions near the sonic point. Closer to the fully nonlinear regime, Gardner, Brook, and Bernstein~\cite{GardnerBrookBernstein1982} and Chen, Zhang, and Panarella~\cite{ChenZhangPanarella1995} combined asymptotic methods with numerical experiments to argue for linear stability of the Guderley solution in radial symmetry. Taken together, these studies provide a consistent picture: linear stability of the Guderley solution under radial perturbations, and a small number of discrete linear instabilities at low angular frequencies outside radial symmetry. To date, there is no rigorous nonlinear PDE analysis which confirms these expectations.

\subsubsection{Smooth isentropic implosions: existence, stability, and implications}

Three existence results concerning the smooth isentropic implosions deserve special mention, beyond what was already reviewed in the introductory bullet list.
First, the original construction of Merle, Rapha\"el, Rodnianski, and Szeftel~\cite{MRRS2022a} applies to 
all adiabatic exponents \(\gamma>1\) except for a countable set, excluding in particular the physically relevant monatomic case $\gamma = \tfrac 53$ in $d=3$. Buckmaster, Cao-Labora, and G\'omez-Serrano~\cite{BCG2025} extended the construction to \emph{every} $\gamma > 1$, while Shao, Wang, Wei, and Zhang~\cite{ShaoWangWeiZhang2025} specifically treated the (algebraically degenerate) three-dimensional monatomic case $\gamma = \tfrac 53$, producing an infinite sequence of implosion profiles.

On the dynamical side, the stability of the smooth isentropic implosions is only known for perturbations which lie on a  finite co-dimension subset of a weighted Sobolev space, which is unquantified. Within radial symmetry, nonlinear stability of the profiles from~\cite{MRRS2022a} for data on a finite-codimension subspace of a weighted Sobolev space
was first established by Merle, Rapha\"el, Rodnianski, and Szeftel~\cite{MRRS2022b,MRRS2022}, and subsequently refined in the course of the existence work~\cite{BCG2025}. Outside radial symmetry, the analogous finite-codimension stability statement was obtained by Cao-Labora, G\'omez-Serrano, Shi, and Staffilani~\cite{CGSS2025}. A distinctive line of work combines the implosion mechanism with vorticity production, in symmetry-reduced settings. Chen, Cialdea, Shkoller, and Vicol~\cite{CCSV2024} used a two-dimensional axisymmetric perturbation of the profiles from~\cite{MRRS2022a}, to construct smooth solutions of compressible Euler which exhibit \emph{vorticity blowup}. This result was subsequently extended to all dimensions \(d \geq 3\) by Chen~\cite{Chen2024}, using a higher-dimensional analogue of two-dimensional axisymmetry with additional error terms and controlling general non-radial perturbations. On the numerical side, the work of Biasi~\cite{Biasi2021} is particularly informative: by computing the spectrum of radial perturbations of the profiles from~\cite{MRRS2022a}, the author identifies a small number of genuinely  unstable directions, along which typical perturbations deflect away from the implosion, generically producing a shock before the singular time.  

We emphasize that the smooth implosion mechanism discovered for the compressible Euler system has been shown to have profound implications for singularity formation in related equations. The works~\cite{MRRS2022b,BCG2025,CGSS2025,ShaoWangWeiZhang2025} establish  implosion for barotropic compressible Navier--Stokes with linear viscosity\footnote{We also mention two very recent works on the barotropic compressible Navier--Stokes equations that partially clarify the role of the viscosity structure in permitting or suppressing implosions. When the shear and bulk viscosities scale as $\rho^\delta$ with $0 < \delta < \delta_*(\gamma)$ for an explicit threshold, Chen, Liu, and Zhu~\cite{ChenLiuZhu2025} use isentropic Euler implosions to construct smooth initial data for barotropic Navier--Stokes, that lead to finite-time implosion; this shows that a sufficiently degenerate viscous term fails to suppress the inviscid blowup mechanism. In contrast, when the viscosities scale linearly in $\rho$ (as in the shallow water system, corresponding to $\delta=1$), Chen, Zhang, and Zhu~\cite{ChenZhangZhu2025} prove that arbitrary large spherically symmetric initial data with bounded positive density generate globally regular solutions in dimensions $d \in \{2,3\}$, and therefore cannot develop any implosion or cavitation.} by showing that the Euler self-similar profile dominates the dynamics close to the blowup time in appropriate parameter regimes, with the viscous term treated as a perturbation. In another direction, Merle, Rapha\"el, Rodnianski, and Szeftel~\cite{MRRS2022} used their smooth radial isentropic Euler implosion profiles to prove finite time blowup for the energy-supercritical defocusing NLS equation, resolving a longstanding open problem;  an extension to the non-radial case was obtained recently in~\cite{CGSS_NLS2024}.\footnote{These NLS blowup results are made possible because via the Madelung transform the defocusing NLS equation maps to a system resembling compressible Euler, but with a quantum pressure.}

\subsubsection{Non-smooth shockless implosions}

Between the two extremes of an Euler implosion driven by a strong imploding shock and a $C^\infty$-smooth implosion, there lies an intermediate class of solutions: radially symmetric implosions that contain \emph{no} shock discontinuity, but that  only have finite regularity. Solutions of this type were recently constructed by Jenssen and Tsikkou~\cite{jenssen2023radially} and by Jenssen~\cite{Jenssen2025} for the full Euler system, in the regime where the primary flow variables remain continuous but are not everywhere differentiable as $t \to t_*^-$. The construction combines matched asymptotic expansions with phase-plane analysis of the associated self-similar ODE system, and relies on both rigorous arguments and careful numerical computations of phase-plane orbits.   The dynamical PDE stability of these continuous-but-non-smooth implosions has not yet been investigated.

\subsubsection{Imploding singularities in related models}

Smooth self-similar imploding singularities have also been constructed for systems closely related to compressible Euler, with applications to astrophysics and to supercritical nonlinear wave equations.

For the isothermal Euler--Poisson system, which governs the gravitational collapse of a self-gravitating gas, Guo, Had\v{z}i\'c, and Jang~\cite{guo2021larson} rigorously constructed the \emph{Larson--Penston implosion profile}, a radially symmetric, $C^\infty$-smooth, globally self-similar imploding solution first postulated on astrophysical grounds by Larson~\cite{larson1969numerical} and Penston~\cite{penston1969dynamics}. In  contrast to the smooth isentropic implosions of~\cite{MRRS2022a,BCG2025,ShaoWangWeiZhang2025}, the Larson--Penston profile is expected to be \emph{stable} under radial perturbations, a conjecture supported by formal and numerical work in the physics literature (see, e.g.,~\cite{hanawa1997stability,ori1988simple}, and the discussion in~\cite{guo2025nonlinear}). This expectation was recently turned into a theorem by Guo, Had\v{z}i\'c, Jang, and Schrecker~\cite{guo2025nonlinear}, who established nonlinear dynamical stability of the Larson--Penston profile under radially symmetric perturbations. The proof combines global monotonicity properties of the profile, high-order weighted energy estimates, and computer-assisted arguments (\emph{interval arithmetic}).

For the relativistic Euler equations, Shao, Wei, and Zhang~\cite{SWZ2024} constructed smooth imploding self-similar solutions and used them to prove finite-time singularity formation for the supercritical defocusing nonlinear wave equation with complex-valued solutions~\cite{shao2025blowup}. In a parallel direction, Buckmaster and Chen~\cite{BuckChen2024} constructed relativistic imploding profiles and used them to establish blowup in dimension $d=4$ for the nonlinearity $|u|^{p-1}u$ with $p=7$, the endpoint case of this blowup mechanism for radially symmetric complex-valued solutions. These results illustrate that smooth self-similar imploding profiles for compressible fluid systems are not merely of intrinsic interest, but also serve as a flexible tool for establishing blowup of supercritical wave and dispersive equations.

\subsection{Main ideas of the proof: existence}
\label{sec:intro:proof:existence}
%%%%

%%%%
The construction of the smooth self-similar profiles $(\bar U, \bar \Sigma, \bar B)$ asserted in Theorem~\ref{thm:intro:profiles} is carried out at the level of the renormalized profiles
\begin{equation*}
\bar U(R) =: R \bar V(R), 
\qquad 
\bar \Sigma(R) =: R \bar Q(R),
\qquad
\bar B(R) =: R \bar H(R),
\end{equation*}
introduced in~\eqref{eq:profile:normalize}. This renormalization both encodes the correct vanishing behavior of $\bar U, \bar \Sigma, \bar B$ at $R = 0$ (necessary for the underlying physical fields $\uu, c, b$ to be smooth at $x=0$) and produces unknowns $(\bar V, \bar Q, \bar H)$ which are smooth functions of $R^2$ near $R=0$. In these renormalized variables, the self-similar Euler system~\eqref{eq:euler4} becomes the system~\eqref{eq:euler5}. The entropy equation~\eqref{eq:s5} is a pure transport equation for $\bar H$ along the radial wave speed $\cxbar + \bar V$, which is non-vanishing due to the \emph{outgoing condition} $\cxbar + \bar V > 0$ proven in~\eqref{eq:Omega:invariant:c}; in turn, this allows us to solve (cf.~\eqref{eq:bar:H:def}) for $\bar H$ explicitly in terms of $\bar V$ as 
\begin{equation*}
\bar H(R)
= \bar H(0) 
\exp  \Bigl( - \int_0^R \tfrac{\bar V(R') - \bar V(0)}{R'(\cxbar + \bar V(R'))} \, d R' \Bigr).
\end{equation*}
Note that if $\bar H(0)>0$, the positivity of $\bar H(R)$ for all $R > 0$ is then automatically true, with no additional verification required. Since $\bar H$ is only defined up to a scaling factor, we may let $\bar H(0)=1$.

The original three-equation system~\eqref{eq:euler5} thus collapses to a $2\times 2$ closed system for $(\bar V, \bar Q)$, namely~\eqref{eq:euler6}. Equivalently, this may be written as an autonomous system~\eqref{eq:euler7:all} in which $R \p_R$ acts as the ``time derivative'' and the right side is a rational function of the unknowns. This ODE system is supplemented with boundary conditions at $R=0$, which are obtained as follows.

The behavior of $(\bar V, \bar Q)$ at $R=0$ is uniquely determined by evaluating the system~\eqref{eq:euler5} at $R=0$, which gives  $\bar v_0:= \bar V(0) = - \tfrac{1}{1+\alpha d}$, $\bar q_0:= \bar Q(0) = \tfrac{1}{1+\alpha d}\sqrt{\tfrac{d\gamma}{2\alpha}}$, and the algebraic relation $\cbbar = \cxbar +\bar v_0$ between the two similarity exponents (see~\eqref{eq:values:at:zero:and:cb}). 

In order to capture the higher-order behavior at $R=0$, we postulate a convergent power series ansatz $\bar V(R) = \bar v_0 + \sum_{n\geq 1} \bar v_n R^{2n}$ and $\bar Q(R) = \bar q_0 + \sum_{n\geq 1} \bar q_n R^{2n}$. Substituting into~\eqref{eq:euler6} and matching powers of $R$ produces, at each order $n \geq 1$, a $2\times 2$ linear system $\MMM_n (\bar q_n, \bar v_n)^\intercal = (\mathsf{F}_1(n), \mathsf{F}_2(n))^\intercal$, where the matrix $\MMM_n$ is given explicitly in~\eqref{eq:Mn:def} and the right-hand side $\mathsf{F}_1(n), \mathsf{F}_2(n)$ depends only on the lower-order coefficients $\{\bar v_j, \bar q_j\}_{1\leq j \leq n-1}$. The key observation is now that $\mathsf{F}_1(1) = \mathsf{F}_2(1) = 0$ identically, because when $n=1$ the set $\{j \colon 1\leq j \leq n-1\}$ is empty. Consequently, if the matrix $\MMM_1$ were invertible, then $(\bar q_1, \bar v_1) = (0,0)$, which in turn implies $\mathsf{F}_1(2) = \mathsf{F}_2(2) = 0$. By induction, as long as the matrices $\MMM_n$ are invertible, we are forced into $(\bar q_n, \bar v_n) = (0,0)$ for all $n \geq 1$. Avoiding the trivial solution thus forces $\det(\MMM_\NNN) = 0$ for some integer $\NNN \in \mathbb{N}$. Since $\det(\MMM_n)$ is, by direct computation, a quadratic polynomial in $\cxbar + \bar v_0$ with explicit coefficients (see~\eqref{eq:det:Mn:def}), the condition $\det(\MMM_\NNN) = 0$ reduces to a quadratic equation whose positive root determines $\cxbar = \cxstar(d,\gamma,\NNN)$ via the closed-form expression~\eqref{eq:cx:admissible}; the value of $\cbbar = \cbstar(d,\gamma,\NNN)$ then follows from the algebraic relation $\cbbar = \cxbar +\bar v_0$, see~\eqref{eq:cb:admissible}. At the resonant order $n = \NNN$ the matrix $\MMM_\NNN$ has a one-dimensional kernel, which allows us to pick $(\bar v_\NNN,\bar q_\NNN$) to be parallel to this kernel element. At every other order $n \neq \NNN$, the coefficients $(\bar q_n, \bar v_n)$ are uniquely and recursively determined, and the resulting power series converges for $R>0$ sufficiently small (cf.~Corollary~\ref{cor:power:series:V:Q}).

In order to extend $(\bar V, \bar Q)$ globally to $R \in [0,\infty)$, we exhibit (see Definition~\ref{def:Omega}) a positively invariant rectangular region for the autonomous system~\eqref{eq:euler7:all}, given by
$
\Omega := \bigl\{ (\bar V, \bar Q) : \bar V \in (\bar v_0, - \tfrac{d\gamma}{2} \bar v_0), \; \bar Q \in (0, \bar q_0) \bigr\}
$.
The positive invariance of $\overline{\Omega}$ is verified by checking the sign of the inward normal velocity on each of its four edges (Proposition~\ref{prop:Omega:invariant}); a key role in this verification is played by the closed-form lower bound $\cxbar > \cxstar(d,\gamma,\infty) = \tfrac{1}{1+\alpha d}(1 + \sqrt{\tfrac{\alpha\gamma d}{2}})$ from Lemma~\ref{lem:properties:exponents}. The local-in-$R$ trajectory constructed via power series at $R=0$ enters $\Omega$ along the unstable manifold of the saddle fixed point $(\bar v_0, \bar q_0)$. Within $\overline{\Omega}$ the only fixed points of~\eqref{eq:euler7:all} are $(\bar v_0, \bar q_0)$ and $(0,0)$, both of which lie on $\partial \Omega$, and the strict monotonicity $R \p_R \bar Q < 0$ inside $\Omega$ (Corollary~\ref{cor:Omega:invariant}, see~\eqref{eq:Omega:invariant:a}) prevents the trajectory from returning to $(\bar v_0, \bar q_0)$ or producing a periodic orbit. The Poincar\'e--Bendixson theorem then forces the trajectory to approach the sink $(0,0)$ as $R \to \infty$, and Proposition~\ref{prop:global:profile} concludes the global existence of the profile. The decay rates $(\bar V, \bar Q)(R) \sim R^{-1/\cxbar}$ as $R \to \infty$ are sharpened in Proposition~\ref{prop:power:series:infinity} by constructing convergent power series expansions in inverse powers of $R^{1/\cxbar}$.

The third profile $\bar H$ is recovered from $\bar V$ via the explicit integral formula displayed above. The asymptotic behavior $\bar H(R) \sim R^{-(\cxbar - \cbbar)/\cxbar}$ as $R \to \infty$ is a direct consequence of the asymptotics of $\bar V$ (see~\eqref{eq:bar:H:R=infty}), and the leading order behavior $\bar H(R) = 1 - \tfrac{1}{2 \NNN(\cxbar + \bar v_0)} R^{2\NNN} + \OO(R^{4\NNN})$ as $R \to 0^+$ (see~\eqref{eq:bar:H:R=0}) makes explicit the connection between the resonant order $\NNN$ in the construction and the order of vanishing of $\bar H - 1$ at $R = 0$.

\subsection{Main ideas of the proof: stability in radial symmetry}
\label{sec:intro:proof:stability}
The stability analysis of the self-similar profiles $(\bar V, \bar Q, \bar H)$, mentioned in Theorems~\ref{thm:intro:radial:stab} and~\ref{thm:intro:nonrad:stab} above, is carried out first in radial symmetry (Section~\ref{sec:stability}). The key new ingredients
needed to handle non-radial perturbations (Sections~\ref{sec:nonradial} and~\ref{sec:ODE}) are discussed in Section~\ref{sec:intro:proof:stability:2} below.

The stability analysis for the Euler system in radial symmetry~\eqref{eq:euler2} is done in modulated self-similar coordinates, using the time-dependent transformation~\eqref{eq:self-similar:ansatz}, which uses three scalar modulation functions $(\cx, \cu, \cb)(\tau)$, and makes precise the discussion in Section~\ref{sec:self:similar:intro}. After renormalizing at $R=0$ via $V(R,\tau) := U(R,\tau)/R$, $Q(R,\tau) := \Sigma(R,\tau)/R$, $H(R,\tau) := B(R,\tau)/R$, and $K := \log H$ (see~\eqref{eq:class:of:perturbations}), we obtain a closed evolution equation for $(V, Q, K)$ (see~\eqref{eq:euler:final}). The stability goal is then to prove that if we choose the modulation functions $(\cx,\cu,\cb)$ and the initial data $(V,Q,K)|_{\tau=0}$ appropriately, then $(V,Q,K) \to (\bar V, \bar Q, \log \bar H)$ and $(\cx,\cu,\cb) \to (\cxbar,\cxbar-1,\cbbar)$ as $\tau\to \infty$.

In order to achieve this, we write the solution as a perturbation of the stationary profile,
\[
(V, Q, K) = (\bar V + \tilde V, \, \bar Q + \tilde Q, \, \bar K + \tilde K), \qquad
(\cx, \cu, \cb) = (\cxbar + \cxtilde, \, \cubar + \cutilde, \, \cbbar + \cbtilde),
\]
with $\bar K := \log \bar H$. Substituting into~\eqref{eq:euler:final} yields the equation satisfied by the perturbations~\eqref{eq:euler:tilde}. The role of the modulation functions $(\cxtilde, \cutilde, \cbtilde)$ which enter~\eqref{eq:euler:tilde} through the source terms $(\mathcal{F}_{\tilde Q}, \mathcal{F}_{\tilde V}, \mathcal{F}_{\tilde K})$ and the nonlinear terms $(\mathcal{N}_{\tilde Q},\mathcal{N}_{\tilde V})$, is to ``mod out'' a small number of apparent instabilities in the linearization of~\eqref{eq:euler:tilde} at $R=0$, as we describe next.

By construction, for $0 < R \ll 1$ we have that $\bar V(R) = \bar v_0 + \bar v_{\NNN} R^{2\NNN} + h.o.t.$, $\bar Q(R) = \bar q_0 + \bar q_{\NNN} R^{2\NNN} + h.o.t.$, and $\bar K(R) = \bar k_{\NNN} R^{2\NNN} + h.o.t.$; see~\eqref{eq:Taylor:R=0:all}. The sole role of the $3(\NNN-1)$ compatibility conditions in~\eqref{eq:thm:intro:radial} (which is the same as assumption (ii) in Theorem~\ref{thm:stability:radial}), is to ensure that the initial perturbations $(\tilde V,\tilde Q,\tilde K)|_{\tau =0}$ have the same structure for $0<R \ll 1$: a constant $+$ a constant times $R^{2\NNN}$ $+$ higher order terms. We may then show that this structure (order of vanishing) is preserved for positive time $\tau$, as long as the solution remains sufficiently smooth; see Lemma~\ref{lem:vanishing:at:R=0}. This allows us to write
\begin{equation}
\label{eq:intro:Taylor:radial}
\bigl( \tilde V, \tilde Q, \tilde K \bigr)(R,\tau)  
= 
\bigl(
\tilde v_0(\tau) + \tilde v_{\NNN}(\tau) R^{2\NNN} ,  
\tilde q_0(\tau) + \tilde q_{\NNN}(\tau) R^{2\NNN} ,
\tilde k_0(\tau) + \tilde k_{\NNN}(\tau) R^{2\NNN} 
\bigr)
+ \OO(R^{2\NNN+1}),
\end{equation}
for $0<R\ll1$ and $\tau \geq 0$; see~\eqref{eq:Taylor:R=0:tilde}.

At zeroth order, by restricting~\eqref{eq:euler:tilde} to $R=0$ we obtain an ODE system for $\vec{\mathfrak{u}}_0(\tau):= (\tilde v_0,\tilde q_0,\tilde k_0)^\intercal(\tau)$ 
(cf.~\eqref{eq:euler:tilde:vq:zero} and~\eqref{eq:euler:tilde:h:zero}), which 
is of the type $\frac{d}{d\tau} \vec{\mathfrak{u}}_0 + \mathcal{A}_0 \vec{\mathfrak{u}}_0 + \vec{\mathfrak{m}}_0 = $ nonlinear terms. Here $\mathcal{A}_0$ is a $3\times 3$ matrix with a positive eigenvalue, a negative eigenvalue, and a zero eigenvalue,\footnote{Throughout this paper, we adopt the sign convention: in the linear evolution $\frac{d}{d\tau} \vec{\mathfrak{u}} + \mathcal{A} \vec{\mathfrak{u}} = 0$, an eigenvalue of $\mathcal{A}$ with positive real part corresponds to exponential \emph{decay} (a stable direction), while an eigenvalue of $\mathcal{A}$ with negative real part corresponds to exponential \emph{growth} (an unstable direction). We treat zero and imaginary eigenvalues of $\mathcal{A}$ in the same way as unstable ones.} while $\vec{\mathfrak{m}}_0$ is a modulation vector, determined in terms of $\cxtilde$, $\cutilde$, and $\cbtilde$. The specific structure of \eqref{eq:euler:tilde:vq:zero}--\eqref{eq:euler:tilde:h:zero} then allows us to choose two linear relations among the modulation functions $(\cxtilde,\cutilde,\cbtilde)$, see~\eqref{eq:modulation:1} and~\eqref{eq:modulation:2}, so that with the resulting $\vec{\mathfrak{m}}_0$ the previous system becomes $\frac{d}{d\tau} \vec{\mathfrak{u}}_0 + \mathcal{A}_0^* \vec{\mathfrak{u}}_0  = $ nonlinear terms, where the matrix $\mathcal{A}_0^*$ only has positive (stable) eigenvalues.\footnote{The following toy example is revealing. Let $\mathfrak{u} \in \Reals^3$ solve the ODE system $\frac{d}{d\tau} \mathfrak{u}  + \mathrm{diag}(-1,0,1)  \mathfrak{u} + (\mathfrak{m}_1,\mathfrak{m}_2,0)^\intercal = 0$, with modulation functions $\mathfrak{m}_1$ and $\mathfrak{m}_2$ which we are free to choose, but only if they decay exponentially fast to $0$ as $\tau\to \infty$. The natural choice is $\mathfrak{m}_1:= 2 \mathfrak{u}_1(\tau)$ and $\mathfrak{m}_2:= \mathfrak{u}_2(\tau)$; with this choice, the system becomes $\frac{d}{d\tau} \mathfrak{u}  + \mathrm{diag}(1,1,1)  \mathfrak{u} = 0$, and each component of $\mathfrak{u}$ decays to $0$ exponentially fast. A posteriori, this decay justifies the choice of $\mathfrak{m}_1$ and $\mathfrak{m}_2$.} Thus, if $|\vec{\mathfrak{u}}_0(0)|$ is sufficiently small (to treat the nonlinear terms perturbatively), then we may expect $|\vec{\mathfrak{u}}_0(\tau)|$ to decay exponentially fast to $0$ as $\tau \to \infty$.

At order $R^{2\NNN}$ we have a similar story. From~\eqref{eq:euler:tilde} we deduce an ODE system for $\vec{\mathfrak{u}}_{\NNN}(t):=(\tilde v_{\NNN},\tilde q_{\NNN}, \tilde k_{\NNN})^\intercal(\tau)$ (cf.~\eqref{eq:tilde:vqh:2N}), which is of the type $\frac{d}{d\tau} \vec{\mathfrak{u}}_\NNN + \mathcal{A}_\NNN \vec{\mathfrak{u}}_\NNN + \vec{\mathfrak{m}}_\NNN = $ nonlinear terms. This time, the matrix $\mathcal{A}_\NNN$ has two positive eigenvalues and a zero eigenvalue. Thus, upon choosing a third linear relation among the modulation functions $(\cxtilde,\cutilde,\cbtilde)$, see~\eqref{eq:modulation:3}, with the resulting $\vec{\mathfrak{m}}_{\NNN}$ the previous system becomes $\frac{d}{d\tau} \vec{\mathfrak{u}}_\NNN + \mathcal{A}_\NNN^* \vec{\mathfrak{u}}_\NNN  = $ nonlinear terms, where the matrix $\mathcal{A}_{\NNN}^*$ has only positive (stable) eigenvalues. Thus, if $|\vec{\mathfrak{u}}_{\NNN}(0)|$ is sufficiently small, then we may expect $|\vec{\mathfrak{u}}_{\NNN}(\tau)|$ to decay exponentially fast to $0$ as $\tau \to \infty$.

The three linear relations on the modulation functions, summarized in~\eqref{eq:modulation:final}, make $(\cxtilde,\cutilde,\cbtilde)$ linear functions of $(\tilde v_0,\tilde q_0,\tilde k_0,\tilde v_{\NNN}, \tilde q_{\NNN}, \tilde k_{\NNN})$. A careful analysis then yields the exponential decay (as $\tau \to \infty$) for the three modulation functions and for the leading order Taylor coefficients (corresponding to $R^0$ and $R^{2\NNN}$) of $(\tilde V,\tilde Q, \tilde K)$ at $R=0$; see bootstraps~\eqref{eq:boot:0}, \eqref{eq:boot:1}, and~\eqref{eq:boot:3}. 

At all higher orders $R^{2n}$ with $n > \NNN$, the matrix corresponding to $\mathcal{A}_n$ in the above discussion may be shown to have all three eigenvalues with \emph{strictly positive real part} (see~\eqref{eq:dummy:equation:0}), so that all directions are stable for the forward in time evolution; see Proposition~\ref{prop:Taylor:boot:closure:higher}. Thus, if the initial data for the Taylor coefficients for $(\tilde V,\tilde Q,\tilde K)|_{\tau =0}$ at $R=0$, of orders $\ell \in \{ 2\NNN+2, 2\NNN+4, \ldots, \MMM\}$, is taken to be sufficiently small (with respect to $\underline{\eps}$), then we may expect all these Taylor coefficients to decay exponentially fast to $0$ as $\tau \to \infty$; see bootstrap~\eqref{eq:boot:2}. Here and throughout this subsection, we denote $\MMM=17 \NNN$.

So far, we have shown that the Taylor coefficients of the perturbation $(\tilde V,\tilde Q,\tilde K)$ at $R=0$, of orders $\leq \MMM$, decay exponentially in time to $0$ as $\tau \to \infty$; and we have shown that the perturbation modulation functions $(\cxtilde,\cutilde,\cbtilde)$ also decay exponentially to $0$. We emphasize that the only information about the profiles $(\bar V,\bar Q,\bar K)$ which was used in this part of the proof are the \emph{explicit values} of  $(\bar v_0,\bar q_0,\bar k_0)$, $(\bar v_{\NNN},\bar q_{\NNN},\bar k_{\NNN})$, and of $(\cxbar,\cubar,\cbbar)$; we did not use any information about the similarity profiles for $R>0$.

It remains to analyze the fields $(\tilde V,\tilde Q,\tilde K)(R,\tau)$ for $R>0$, and to show that they decay to $0$ as $\tau \to \infty$. A natural first attempt would be to return to the evolution equation~\eqref{eq:euler:tilde} satisfied by these perturbations, ignore for the moment the nonlinear $\mathcal{N}_{\bullet}$ and forcing $\mathcal{F}_{\bullet}$ terms, and to try to extract enough damping from the linear terms $\mathcal{D}_{\bullet}$ defined in~\eqref{eq:D:Q:def}--\eqref{eq:D:H:def} to guarantee a \emph{coercive/dissipative estimate} for a suitably defined (weighted Sobolev) norm of $(\tilde V,\tilde Q,\tilde K)$. This is the classical ``linearization + coercivity/dissipativity'' strategy. This strategy is very powerful, but it requires detailed information about the (variable in $R$) coefficients of the linear terms $\mathcal{D}_{\bullet}$; that is, information about $\bar V(R)$ and $\bar Q(R)$ for all $R>0$. In the setting of this paper, the profiles $\bar V(R)$ and $\bar Q(R)$ are not explicit. To overcome this fundamental difficulty 
and establish the global-in-$R$ decay of $(\tilde V,\tilde Q,\tilde K)(R,\tau)$ as $\tau \to \infty$
by a purely analytic argument, we proceed as follows:
\begin{itemize}[leftmargin=1em]

\item Our first observation is that the profiles $(\bar V,\bar Q)$ obey a \emph{global outgoing property}; see~\eqref{eq:Omega:invariant:c}, which we re-state here for convenience:
\begin{equation}
\label{eq:intro:outgoing}
\cxbar + \bar V(R) + \alpha \bar Q(R) > \cxbar + \bar V(R) > \cxbar +\bar V(R) - \alpha \bar Q(R) \geq \tfrac{1 + \frac 23 \alpha d}{4 \NNN(1+\alpha d)}>0, 
\qquad\forall R\geq0.
\end{equation}

\item Our second observation is that the radially symmetric setting presents a tremendous luxury: it allows us to perform the analysis of $(\tilde V,\tilde Q,\tilde K)(R,\tau)$ at the pointwise level, without having to resort to $L^2$-based energy estimates. The  idea is to introduce the differentiated Riemann-type variables\footnote{Differentiated Riemann-type variables, which linearize (to leading order) the differentiated Euler system, have played a crucial role in the recent developments on shock formation for the compressible Euler system~\cite{BDSV,NRSV24,NealShkollerVicol25,NSV25,CialdShkVic2025}.}  (see~\eqref{eq:ringW:ringZ:ringA:def} and~\eqref{eq:YY:def})
\begin{align*}
&\Wring := R \p_R V + R \p_R Q - \tfrac{1}{\gamma} Q R \p_R K, 
&&\Zring := R \p_R V - R \p_R Q + \tfrac{1}{\gamma} Q R \p_R K, \\
&\Aring := \tfrac{\alpha}{\gamma} Q R \p_R K, &&
\YY := (\Zring, \Aring, \Wring)^\intercal.
\end{align*}
The evolution of the vector $\YY$ is obtained by applying $R\p_R$ to ~\eqref{eq:euler:final}, and then diagonalizing the resulting system. This procedure shows that $\YY$ solves the quasilinear evolution equation (see~\eqref{eq:YY:evo})
\begin{equation}
\label{eq:intro:Y:evo}
\partial_\tau \YY + \TT  \, R \p_R \YY + \DD \, \YY + \NN(\YY, \YY) = 0,
\end{equation}
where $\TT= \TT(V,Q) := {\rm diag}(\cx + V - \alpha Q, \cx+ V, \cx+V + \alpha Q)$ is a diagonal \emph{transport matrix} whose entries are precisely the three self-similar wave speeds present in the system, $\DD = \DD(V,Q)$ is an explicit \emph{damping matrix} which is linear in $V$ and $Q$, and $\NN$ is a quadratic nonlinearity with constant coefficients. Crucially, the lower bound on the three self-similar wave speeds inherited from the global outgoing property of the profile (see~\eqref{eq:intro:outgoing} above) ensures that all three components of $\TT$ remain bounded below, uniformly in $R$ and $\tau$. Therefore, in~\eqref{eq:intro:Y:evo} \emph{information propagates outward} along the $R$-axis, in a uniform fashion, which emphasizes the importance of very sharp estimates at radii $0<R\ll 1$.

\item
The third observation concerns the evolution equation for the perturbation $\Ytilde = \YY-\Ybar$ (defined in the natural way) at radii $0<R\ll 1$; upon subtracting from~\eqref{eq:intro:Y:evo} the analogous stationary PDE for $\Ybar$, we obtain that $\Ytilde$ satisfies $\p_\tau \Ytilde + \bar \TT R\p_R \Ytilde + \bar \DD \Ytilde = \mathcal{F}$, for a suitable ``force'' $\mathcal{F}$. Upon multiplying this evolution by a singular weight, such as $R^{-
(\MMM+1)}$ for some $\MMM \in \Naturals$, we obtain that the quantity $R^{-(\MMM+1)} \Ytilde(R,\tau)$ solves  
\begin{equation}
\label{eq:intro:weighted:Ytilde:evo}
\p_\tau (R^{-(\MMM+1)} \Ytilde) + \bar \TT R \p_R (R^{-(\MMM+1)} \Ytilde) +  ( (\MMM+1) \bar \TT +\bar \DD)  (R^{-(\MMM+1)} \Ytilde) = R^{-(\MMM+1)} \mathcal{F},
\end{equation}
which is a damped and forced transport equation.  The emergence of the beneficial linear damping factor of strength $\MMM+1$ directly manifests the advantage of performing weighted-$L^\infty$ estimates with singular weights. This approach requires, however, that $R^{-(\MMM+1)} \Ytilde(R,\cdot)\in L^\infty_{\rm loc}([0,\infty))$, which is to say that $\Ytilde$ vanishes to high order (at least to order $\MMM+1$) at the origin.

\item The fourth observation is that although $\Ytilde$ \emph{does not vanish to high order} at $R=0$ (see~\eqref{eq:intro:Taylor:radial}), we may nonetheless apply the singular weight idea described in the previous bullet once we subtract from $\Ytilde$ its Taylor polynomial (at $R=0$) of degree $\MMM$. We denote by $\III_{\MMM} \Ytilde$ (see~\eqref{eq:Taylor:poly:1}) the Taylor polynomial of $\Ytilde$, at $R=0$, of degree $\MMM$, smoothly cutoff at radii $R \in [0,2 R_{\mathsf{in}}]$, where $R_{\mathsf{in}}\ll1$ is chosen suitably in the proof. Then, $\Ytilde - \III_{\MMM}\Ytilde$ vanishes to order $R^{\MMM+1}$ at the origin (assuming $\Ytilde \in C^{\MMM+1}_{\rm loc}$), and so we may consider the evolution equation
\eqref{eq:intro:weighted:Ytilde:evo} with $\Ytilde$ being replaced by $\Ytilde - \III_{\MMM}\Ytilde$. The important fact to note is that once this difference is estimated, we recover ``for free'' the desired information on $\Ytilde$, because $\III_{\MMM}\Ytilde$ was previously estimated globally in time (recall, we have shown that the Taylor coefficients of the perturbation $(\tilde V,\tilde Q,\tilde K)$ of orders $\leq \MMM$, decay exponentially in time to $0$ as $\tau \to \infty$).

\item Lastly, we observe that if we replace the singular weight $R^{-(\MMM+1)}$ appearing in~\eqref{eq:intro:weighted:Ytilde:evo} with $\psi(R)^{-(\MMM+1)}$, for a custom  designed weight function $\psi(R)$ (smooth, positive, linear at the origin, constant at infinity), then we may obtain global existence and decay as $\tau \to \infty$ for the quantity $\JJJ_{\MMM}\Ytilde:= (\Ytilde - \III_{\MMM}\Ytilde)\psi^{-(\MMM+1)}$ (cf.~\eqref{eq:J:M:def}), \emph{without using very precise information on $\bar \TT(R)$ and $\bar \DD(R)$ at all $R>0$}. 
The idea is that using the weight $\psi(R)^{-(\MMM+1)}$, we generate a strong damping coefficient $ (\MMM+1) \frac{R \pa_R \psi}{\psi} \bar \TT + \bar \DD$ in \eqref{eq:intro:weighted:Ytilde:evo}
instead of $(\MMM+1) \bar \TT + \bar \DD$. 
For any $R_1 \ll 1 \ll R_2$, we can design the weight $\psi(R)$ with $R \in [R_1, R_2]$ almost \emph{freely}. 
If we take $\MMM$ sufficiently large ($\MMM=17 \NNN$ suffices) and if $\psi$ is constructed very carefully (see~Proposition~\ref{prop:choice:of:M:psi}), we may get away with only knowing the explicit values of $\bar V$ and $\bar Q$ (hence $\bar \TT$ and $\bar \DD$) as $R\to 0^+$ and $R\to \infty$;\footnote{This point is important, because at $R=0$ and as $R\to \infty$ we know $\bar V$ and $\bar Q$ explicitly, while for intermediate values of $R$ we only know a few basic properties of these profiles. 
This allows us to complete the proof \emph{without} computer assistance.} in essence, this is made possible by the global outgoing property~\eqref{eq:intro:outgoing}.\footnote{Similar outgoing properties have played a crucial role in the stability analysis of self-similar blowup for the 3D incompressible Euler equations \cite{ChenHou2023a} (see also \cite{chen2025singularity}), 
and for the compressible Euler system~\cite{CCSV2024,Chen2024}.
}
That is, we are able to prove that the \emph{system of transport equations} obeyed by $\JJJ_{\MMM}\Ytilde$ (cf.~\eqref{eq:JM:tildeY:evo}) is such that the \emph{damping matrix is ``diagonally dominated''} (cf.~\eqref{eq:bazooka}), which allows us to establish global existence and exponential decay to $0$ as $\tau\to \infty$ via a Gershgorin-type argument.

\end{itemize}

The proof of Theorem~\ref{thm:stability:radial} then proceeds along the ideas described in the five points above. This is implemented technically by a bootstrap argument, with five mutually-reinforcing assumptions: exponential decay of the zeroth-order Taylor coefficients $(\tilde v_0, \tilde q_0, \tilde k_0)$ (cf.~\eqref{eq:boot:0}); exponential decay of the order-$\NNN$ Taylor coefficients $(\tilde v_\NNN, \tilde q_\NNN, \tilde k_\NNN)$ (cf.~\eqref{eq:boot:1}); exponential decay of all higher-order even Taylor coefficients up to order $\MMM = 17 \NNN$ (cf.~\eqref{eq:boot:2}); exponential decay of the modulation functions $(\cxtilde, \cutilde, \cbtilde)$ (cf.~\eqref{eq:boot:3}); and a weighted-pointwise decay $\langle R \rangle^\theta |\JJJ_\MMM \Ytilde| \lesssim \underline{\eps}^{1/5} e^{-\underline{\lambda} \tau}$, with $\theta := \tfrac{1}{2} \min\{1, 1/\cxbar\}$, see~\eqref{eq:boot:4}. Once these bootstrap bounds are closed, a standard continuity argument completes the proof.

\subsection{Stability outside of radial symmetry}
\label{sec:intro:proof:stability:2}

The stability analysis for non-radial perturbations of the globally self-similar solution $(\bar\rho,\bar\uu,\bar p)$, carried out in Sections~\ref{sec:nonradial} and~\ref{sec:ODE}, follows the same overall philosophy as in the radial case; the analysis is performed in modulated self-similar coordinates, and the perturbation is split into a finite collection of Taylor coefficients at $y=0$ (which obey a \emph{significantly} more complicated, but still closed, finite-dimensional system of ODEs), and a bulk remainder controlled by a singularly weighted PDE energy estimate. The two halves of the argument are coupled through a finite list of modulation functions and a bootstrap loop. Compared to the radial case, however, the absence of symmetry forces several new structural difficulties, which we describe in turn.

\subsubsection{Choice of variables} 
A first difficulty is at the level of unknowns. Smoothness in $x$ near the origin is no longer equivalent to smoothness in $|x|^2$ (mixed terms such as $x_i x_j$ with $i\neq j$ may now occur), so the radial renormalized variables $(V, Q, H)$ and the differentiated Riemann variables $\YY$ are no longer suitable. We therefore work directly in the unknowns 
\begin{equation}\label{eq:intro:nonrad:vars}
   (\vvrho,\uu,\mbb) := \bigl(\rho^{\gamma-1},\,\uu,\,e^s\bigr),
\end{equation}
introduced in~\eqref{eq:nonrad_var}--\eqref{eq:sys_UB}, whose smoothness is equivalent to that of the primary flow variables $(\rho,\uu,p)$ in~\eqref{eq:euler:primary}. A key advantage of $(\vvrho,\uu,\mbb)$ over $(\rho,\uu,p)$ is that the resulting evolution has purely \emph{quadratic nonlinearities}, with no occurrence of the singular factor $\rho^{-1}\nabla p$; this \emph{significantly} simplifies the higher-order energy estimates, which can be closed by directly applying the standard Leibniz rule.\footnote{In contrast, $\nabla^k(\rho^{-1}p)$ contains terms of the form $\rho^{-(k+1)}p\,(\nabla\rho)^k$ and $\rho^{-(k+1)}p\,(\nabla^{k-i}\rho)(\nabla^i\rho)$. Controlling these terms for large $|y|$ would require sharp \emph{decay} estimates on $\nabla^i \rho$ together with sharp \emph{lower bounds} on $\rho$; cf.~\cite{chen2025stability}.}

\subsubsection{Modulation functions and the perturbation system} 
The radial analysis used three scalar modulation functions, determined by the structure of the obstructions at orders $R^0$ and $R^{2\NNN}$. Outside of radial symmetry the structure of these obstructions is substantially richer, and we are forced to introduce four modulation functions $(\cx,\cu,\cbb,\crho)(\tau)$, see~\eqref{eq:globally:SS:exact:nonradial:a}, subject to the single algebraic constraint $\crho(\tau) + \cbb(\tau) = 2 \cu(\tau)$,
cf.~\eqref{eq:normal0}.
After passing to self-similar coordinates $(\vrho,\UU,\mb)$ for $(\vvrho,\uu,\mbb)$, the stability goal is to prove that, modulo finitely many obstructions, $(\vrho,\UU,\mb)\to(\bvr,\bu,\barb)$ and $(\cx,\cu,\cbb,\crho)\to(\cxbar,\cubar,\cbbbar,\crhobar)$ as $\tau\to\infty$, exponentially fast. 
Writing the solution as a perturbation of the stationary profile,
\[
(\vrho,\UU,\mb) = (\bvr+\tvr,\,\bu+\tu,\,\barb+\tb),
\quad
(\cx,\cu,\cbb,\crho) = (\cxbar+\cxtilde,\,\cubar+\cutilde,\,\cbbbar+\cbbtilde,\,\crhobar+\crhotilde),
\]
the perturbation $\tw := (\tvr,\tu,\tb)$ obeys a PDE system of the schematic form (cf.~\eqref{eq:lin:full}) 
\bseq
\label{eq:intro:nonradial:id1}
\begin{align}
\pa_\tau \tvr &= -( \cx y + \UU) \! \cdot  \! \nabla \tvr - 2 \alpha \vrho \, \div  \tu  
   && + \cvr \tvr  +  (\bar \DD  \tw)_1 + \OOL{1}( \tcvr, \tcu, \tcb)  ,  \\
\pa_\tau \tu &=  -  ( \cx y + \UU) \!  \cdot \! \nabla   \tu  - \tfrac{1}{2\alpha} \mb \nabla \tvr - \tfrac{1}{\gamma} \vrho \nabla \tb   && +  \cu  \tu  
+  (\bar \DD  \tw)_2 + \OOL{2}( \tcvr, \tcu, \tcb), \\ 
\pa_\tau \tb &=  - (\cx y + \UU)  \! \cdot  \! \nabla \tb   && +  \cbb \tb  +
  (\bar \DD  \tw)_3  + \OOL{3}( \tcvr, \tcu, \tcb) , 
 \end{align}
 \eseq
where $\bar \DD(y) \in \Reals^{3 \times 3}$ is a damping matrix, $ \OOL{i}( \tcvr, \tcu, \tcb)$ 
denotes perturbation terms linear in the modulation functions $\tcvr, \tcu, \tcb$, with 
coefficients depending on the profile $(\bvr,\bu,\barb)$. As in the radial case, our analysis of~\eqref{eq:intro:nonradial:id1} is partially nonlinear: we do not decompose the coefficients of the top-order derivatives into stationary profile plus perturbation, nor do we decompose the modulation functions appearing in the lower-order terms.

The argument is conditional on a bootstrap assumption for the Taylor coefficients of $\tw$ at $y=0$, encoded in~\eqref{eq:ass:EOM}; the verification of this assumption is carried out independently in Section~\ref{sec:ODE}. After subtracting a smooth Taylor cutoff, this bootstrap allows us to assume that the perturbation $(\tvr, \tu, \tb)$ vanishes to a sufficiently high order at the origin. We now describe the two halves of the argument.

\subsubsection{Bulk PDE stability via singularly weighted $H^k$ estimates}
For non-radial perturbations, the multi-D Euler system cannot be diagonalized, and the differentiated Riemann variable framework which we have employed in the radial case must be replaced. We exploit the symmetric hyperbolic structure of the Euler system, thereby avoiding derivative loss, and build on a singularly weighted $L^2$-energy framework inspired by~\cite{chen2019finite,chen2021HL,ChenHou2023a,CCSV2024,Chen2024}, which works directly with the system~\eqref{eq:intro:nonradial:id1} for  $(\tvr,\tu,\tb)$. The main result is recorded in Theorem~\ref{thm:nonrad:stab}. 

In~\eqref{energy:Hk} we introduce\footnote{Instead of $\tw$, we use $E_\al(\cdot)$ with the Taylor remainder term $\twm$, cf.~\eqref{eq:nonrad:Taylor}--\eqref{eq:pertb:main}, which vanishes to high order at $y=0$.} an $H^k$-energy density\footnote{The Euler equations form a symmetric hyperbolic system: after a suitable change of variables, akin to~\eqref{eq:intro:nonrad:vars}, they take the schematic form $\partial_\tau Y + \sum_{i\leq d} M_i\,\partial_{y_i} Y = \mathrm{l.o.t.}$ with $M_i = M_i^\top$. The cross term in $E_\al$ exploits this structure to allow integration by parts at top order without derivative loss. Similar energy structures, exploiting symmetric hyperbolicity, have proved effective in~\cite{CCSV2024,Chen2024,bedrossian2026finite}.}
\begin{equation}
\label{intro:energy:Hk}
E_\al(\tw) := \frac{|\partial^\al \tvr|^2}{(2\al\,\vrho)^2}
+ \frac{|\partial^\al \tu|^2}{\cc^2}
+ \kp_\mb \frac{|\partial^\al \tb|^2}{\mb^2}
+ \frac{2}{\gamma}\,\frac{\partial^\al \tvr}{2\al\,\vrho}\cdot\frac{\partial^\al \tb}{\mb},
\qquad
E_k(\tw) := \sum_{|\al|=k}E_\al(\tw),
\end{equation}
where $\cc := (\vrho\,\mb)^{1/2}$ is the sound speed, and $\kp_\mb$ is a large coupling parameter. For $\kp_\mb$ large enough (cf.~\eqref{eq:def:kp_mb}), the density $E_\al$ is positive definite. Since $\vrho,\cc,\mb$ may vanish at infinity, we modify the density profile (cf.~\eqref{eq:vrho_s}) and propagate two-sided pointwise bounds on $(\vrho,\mb)$ via the bootstrap~\eqref{eq:boot1}.

The $H^k$ estimate is closed by integrating $E_\al(\tw)$ against a carefully chosen, time-independent radial weight $\vp_k(y)$. After integration by parts in the symmetric hyperbolic top-order terms, the energy identity reads, schematically,
\begin{align}
&\frac{d}{d\tau}\!\int E_\al(\tw)\,\vp_k
\label{eq:intro:Hk}
\\
&\leq
\int \underbrace {\frac{(\cx y + \UU) \cdot \na \vp_k}{2\vp_k} E_{\al}(\tw) \vp_k  }_{:=I_1}
+ \underbrace{ \frac{\cc |\na \vp_k|}{\vp_k} \Big(  \frac{ |\pa^{\al} \tum| }{\cc} \frac{| \pa^{\al} \tvrm|}{2 \al \vrho}  \vp_k + 
\frac{1}{\gamma} \frac{|\pa^{\al} \tbm|}{\mb} \frac{|\pa^{\al} \tum|}{\cc} \Big) \vp_k }_{:=I_2}
 + \cR_k  \vp_k ,
\notag
\end{align}
where $\cR_k$ denotes the $\vp_k$-\emph{independent} contributions from the terms $ \bar \DD \tw, \OOL{i}$,
$\cvr \tvr ,  \cu  \tu  , \cbb \tb$ in \eqref{eq:intro:nonradial:id1}, and the contributions from 
the lower order terms, e.g. when derivatives acting on the denominators. 
The terms $I_1, I_2$ can be seen  as the \emph{damping effect} generated by the \emph{variation} of the weight $\vp_k$, which captures the outgoing property of the three self-similar wave speeds.
To see this, note that by the Cauchy-Schwarz inequality and using the definition of the energy density $E_{\al}$ in~\eqref{intro:energy:Hk}, we may estimate 
\begin{align*}
I_1 + I_2 
\leq \underbrace{\frac{ (\cx y + \UU) \cdot \na \vp_k + (1 + 3 \kp_{\mb}^{-1/2}  ) \cc |\na \vp_k|  }{ 2 \vp_k}  }_{=:\mathcal D_{\vp_k}(y) }  E_{\al}( \tw ) \vp_k .
\end{align*}

The aim is to obtain a sufficiently strong damping term in the weighted estimate~\eqref{eq:intro:Hk}, namely, to make $\mathcal D_{\vp_k}$ very negative. For this purpose, we choose a radially symmetric weight $\vp_k(y)$ that decreases in $|y|$ in a large domain. 
For sufficiently small perturbations, so that $\cx \approx \cxbar, U(R) \approx \bar U(R) = R \bar V(R),  \cc(R) \approx \bar \cc(R) = \alpha R \bar Q(R)$, the above damping coefficients satisfies
\[
\mathcal D_{\vp_k}(y) \approx 
- 
\Bigl[\, \cxbar + \bar V(R) - \alpha \bigl(1+3\kp_\mb^{-1/2}\bigr)\bar Q(R)\,\Bigr]\,
\frac{|R \partial_R \vp_k|}{\vp_k},
\qquad R = |y|.
\]
By taking $\kp_{\mb}$ sufficiently large (cf.~\eqref{eq:def:kp_mb}), from the \emph{global outgoing property} of the self-similar profiles (cf.~\eqref{eq:intro:outgoing}), we obtain that $\cxbar + \bar V(R) - \alpha (1 + 3 \kp_{\mb}^{-1/2}  ) \bar Q(R) \geq C_0 $ for some $C_0 = C_0(d,\gamma,\NNN)>0$, for all $y$. Thus $\mathcal D_{\vp_k}(y)< 0$ for all $y$.  
Consequently, by choosing $\vp_k = \psi_0^N$ with $\psi_0$ radially decreasing in a large compact annulus, using $ \frac{R \pa_R \vp_k }{ \vp_k } = N \frac{R \pa_R \psi_0}{\psi_0}$,  and taking $N$ large, we generate \emph{arbitrarily strong} damping in any prescribed annulus $R_1 \leq |y| \leq R_2$, with $R_1 \ll 1 \ll R_2$. Choosing $\vp_k \asymp |y|^{-2\MMM-d}$ near the origin produces the strong damping needed to control low-order terms (this plays the same role as the singular weight $R^{-(\MMM+1)}$ in~\eqref{eq:intro:weighted:Ytilde:evo}). For large values of $|y|$, we design $\vp_k$ to have sufficient decay to generate the necessary amount of damping.

With damping of arbitrary strength available, the remainder $\mathcal R_k$ may be treated \emph{perturbatively}, without any precise knowledge of the profile in the intermediate range of~$y$.\footnote{Such singularly weighted estimates have played an important role in the stability analysis of self-similar blowup for fluid equations and related models, see~\cite{ChenHou2023a,chen2025singularity,chen2020singularity,Chen2024,CCSV2024}.} 

We note importantly that the order $k$ in~\eqref{intro:energy:Hk} is chosen large enough to close the nonlinear estimates, but it is \emph{not too large}. In fact, the value of $k$ taken in our proof is \emph{explicit}, in terms  of $d$ (cf. \eqref{eq:k_star}); it is not an unquantified large constant as in~\cite{MRRS2022b,BCG2025}.

\subsubsection{Cutoff Taylor polynomials}
Due to the presence of the singular weight $\vp_k$, the energy appearing in~\eqref{eq:intro:Hk} requires that $\tw$ vanish to high order at $y=0$, which is not automatic for general perturbations. As in the radial case, this is handled by subtracting a smoothly cutoff Taylor polynomial $\tl_k\tilde f$ at the origin and propagating the estimates only on the remainder $\rl_k\tilde f$ (cf.~\eqref{eq:nonrad:Taylor}). A significant \emph{technical difficulty} arises because the three components $(\vrho,\uu,\mb)$ have \emph{different} vanishing orders near $y=0$, namely $\vrho \asymp 1$, $|\uu|\asymp |y|$, and $|\mb|\asymp|y|^2$ for smooth solutions close to the stationary profile; hence, the radial trick of dividing by $|y|^i$ and introducing the normalized variables $(V,Q,H)$ is not available. 

To overcome this difficulty, in the non-radial case we subtract Taylor polynomials at \emph{three different orders} (cf.~\eqref{eq:pertb:main}): 
\begin{equation}\label{eq:intro:Taylor:cutoff:new}
\tvr = \tl_\MMM\tvr + \rl_\MMM\tvr,\qquad
\tu = \tl_{\MMM+1}\tu + \rl_{\MMM+1}\tu,\qquad
\tb = \tl_{\MMM+2}\tb + \rl_{\MMM+2}\tb,
\end{equation}
where the cutoff Taylor polynomials $\tl_k \tilde f$ consist of the derivatives at $y=0$:
$\{\na^i \tvr(0)\}_{i\leq \MMM}$, $\{\na^{i} \tu\}_{i\leq \MMM+1}$, and $\{\na^{i} \tb\}_{i\leq \MMM+2}$.
Each remainder $(\rl_{\MMM} \tvr, \rl_{\MMM+1} \tu, \rl_{\MMM+2} \tb)$ vanishes to order $|y|^{\MMM+1}$ at the origin, which is the input required by~\eqref{intro:energy:Hk}.

\subsubsection{Closed ODE system for the Taylor coefficients at the origin}

As in the radial case, the above stability argument is conditional: the bulk decay estimate~\eqref{eq:nonrad_decay} is conditional on the bootstrap assumption~\eqref{eq:ass:EOM} on a finite collection of Taylor coefficients at $y = 0$, and on the modulation functions. The verification of this bootstrap assumption~\eqref{eq:ass:EOM} is the subject of Section~\ref{sec:ODE} and Theorem~\ref{thm:ODE_stab}.

For the $\NNN$-th profile, the modulation functions are determined, through a careful analysis paralleling the radial case at orders $R^0$ and $R^{2\NNN}$, by the Taylor coefficients at $y=0$ of 
\begin{equation}\label{eq:intro:ODE1:new}
\tvr \in \Reals,\quad 
\nabla\tu \in \Reals^{d\times d},\quad \nabla^2\tb\in\Reals^{d\times d},
\qquad
\Delta^\NNN\tvr \in \Reals,\quad \Delta^{\NNN+1}\xu \in \Reals,\quad \Delta^{\NNN+1}\tb\in \Reals,
\end{equation}
where $\xu(y,\tau) := y\cdot\tu(y,\tau)$ plays the role of a radial velocity (cf.~\eqref{eq:XU_V}). The four modulation functions $(\cxtilde,\cutilde,\cbbtilde,\crhotilde)$ and the coefficients in~\eqref{eq:intro:ODE1:new} together satisfy a closed ODE system (one of the modulation functions is redundant in view of the algebraic constraint~\eqref{eq:normal0}). Although this system is significantly more complicated than its radial counterpart, it admits a clean structural decomposition: at the linear level,
\begin{itemize}[leftmargin=1.5em]
\item the \emph{non-radial part}, consisting of the trace-free symmetric gradient $\nabla\tu + (\nabla\tu)^\top - \tfrac{2}{d}(\div\tu)\,\Id$, together with the trace-free Hessian $\nabla^2\tb - \tfrac{1}{d}(\Delta\tb)\,\Id$, and 
\item the \emph{radial part}, consisting of the symmetric scalar combinations in $(\tvr$, $\div\tu$, $\Delta\tb$, $\Delta^\NNN\tvr$, $\Delta^{\NNN+1}\xu$, $\Delta^{\NNN+1}\tb)$, together with the modulation functions $(\cxtilde,\cutilde,\cbbtilde,\crhotilde)$,
\end{itemize}
are \emph{decoupled}. The radial part is structurally analogous to the ODE system of the radial analysis, and the radial proof (in particular, the choice of modulation functions) extends to it directly. The non-radial part decouples further into uncoupled scalar and $2\times 2$ ODE blocks, which admit a \emph{complete} linear stability analysis; this analysis identifies the first batch of unstable directions associated with non-radial perturbations. 

\subsubsection{Higher order ODE system and unstable dimension counts}

The remainder of Theorem~\ref{thm:ODE_stab} concerns the higher mixed-derivative tower
$V_{=n}(\tau) = \bigl(\nabla^n\tvr,\,\nabla^{n+1}\tu,\,\nabla^{n+2}\tb\bigr)(0,\tau)$. 
For $n\geq 2 \NNN$, where $\NNN$ indexes the profile, the full vector $V_{\leq n} := \bigl(V_{=0},\dots,V_{=n}\bigr)$ obeys a closed linear-quadratic ODE,
\[
\frac{d}{d\tau} V_{\leq n} = M_{\leq n}\,V_{\leq n} + Q_n(V_{\leq n},V_{\leq n}),
\]
governed by an explicit constant matrix $M_{\leq n}\in\Reals^{d_{\leq n}\times d_{\leq n}}$ of size $O(n^d)$. The principal task is to count the eigenvalues of $M_{\leq n}$ with non-negative real part, as this count fixes the (nonlinear) unstable+center manifold $\Sigma_{\sf uns,\leq n}$. For general $n \geq 1$, this task is quite complicated.

The strategy is guided by the radial symmetry of the underlying profile $(\bvr,\bu,\barb)$: the more symmetric the derivative combination is, the cleaner the ODE block to which it belongs. To exploit this, for $k\geq 0$ and $\bet\in\Naturals_0^d$ we organize mixed derivatives by powers of the Laplacian,
\begin{equation}\label{intro:def:VVs}
\VVs_{\bet,k}(\tau) := \bigl(
\partial^\bet \Delta^k \tvr,\;
\partial^\bet \Delta^k(\div\tu),\;
\partial^\bet \Delta^{k+1}\xu,\;
\partial^\bet \Delta^{k+1}\tb
\bigr)^\intercal(0,\tau)\in\Reals^{4\times 1},
\end{equation}
where $\xu$ captures the velocity in the radial direction $\xu(y,\tau) = y \cdot \tu(y,\tau)$.

We analyze the system satisfied by the $\VVs_{\bet,k}$ in the order of decreasing $k$, for fixed $|\bet|$. In this ordering, the matrix $M_{\leq n}$ becomes \emph{lower block-triangular}, with diagonal blocks of size at most $4\times 4$ given explicitly by the matrices $\HHH_k$, $\HHT_k$, and $\HH_{|\bet|,k}$ of Proposition~\ref{prop:ODE_Nth}. The eigenvalue analysis then proceeds in two regimes:
\begin{itemize}[leftmargin=1.5em]
\item \emph{(Distinguished cases.)} For the ground state $\NNN=1$ together with the five physical pairs 
$(\gamma,d)\in\bigl\{(\tfrac{5}{3},3),(\tfrac{7}{5},3),(2,2),(\tfrac{5}{3},2)\bigr\}$ or $\gamma\in(1,3]$, $d=1$, the closed form of the diagonal blocks, combined with the Routh--Hurwitz criterion\footnote{
The Routh--Hurwitz criterion allows us to count the number of eigenvalues of a matrix with nonnegative real part by proving a few inequalities on its entries, \emph{without} deriving the eigenvalues \emph{explicitly}.
} (Theorem~\ref{thm:RH})
yields the \emph{exact} dimension counts of $\Sigma_{\sf uns}$ summarized in item~(iv) of Theorem~\ref{thm:ODE_stab}  (see also Remark~\ref{rem:intro:nonrad:ground:state}).
\item \emph{(General $\NNN,d,\gamma$.)} 
For the general case $\NNN \geq 1$, $d\in \{1,2,3\}$, and $1<\gamma\leq 2 d+1$, we use the classical result that a weighted diagonally dominated matrix does not have eigenvalues with non-negative real part (Lemma \ref{lem:diag_eig}). 
For high derivative orders $n \geq n_1 := (18d)^2\NNN$, such a condition holds 
for the ODE system of mixed derivatives $V_{=n}(\tau) =(\na^{n} \tvr,\na^{n+1} \tu, \na^{n+2} \tb)(0,\tau)$.
The analysis produces the stabilization $\dim(\Sigma_{\mathsf{uns}, \leq n}) = \dim(\Sigma_{\mathsf{uns}, \leq n_1})$ for $n \geq n_1$, and further provides the explicit upper bound for the 
unstable + center manifold: $\dim(\Sigma_{\sf uns}) \leq \bar C(d)\,\NNN^d$.

\end{itemize}
With the spectral structure of $M_{\leq n}$ in hand, a global solution of the nonlinear ODE on the unstable+center manifold is constructed by combining classical stable manifold theory with the splitting method of~\cite{ChenHou2023a}.

\subsubsection{Closing the loop}
The bulk PDE estimate of Theorem~\ref{thm:nonrad:stab} and the ODE analysis of Theorem~\ref{thm:ODE_stab} are coupled through the local well-posedness of the full system~\eqref{eq:sys}. On any time interval $[0,T]$ on which the PDE is regular, the Taylor coefficients of the PDE solution at the origin agree with those produced by Theorem~\ref{thm:ODE_stab}; this verifies the bootstrap assumption~\eqref{eq:ass:EOM} on $[0,T]$, which in turn permits the bulk decay estimate~\eqref{eq:nonrad_decay}, and finally the smallness afforded by~\eqref{eq:nonrad_decay} extends the local solution past $T$. A standard continuity argument in $\tau$ (see Remark~\ref{rem:relation}) closes the loop and yields the global-in-$\tau$ stability conclusion stated in Theorem~\ref{thm:intro:nonrad:stab}.

\subsection{Organization}
The paper is organized as follows. In Section~\ref{sec:profiles}, we construct the family of exact self-similar profiles and prove Theorem~\ref{thm:main:profiles}. In Section~\ref{sec:stability}, we establish the stability of these profiles within the class of radially symmetric perturbations and prove Theorem~\ref{thm:stability:radial}. In Section~\ref{sec:nonradial}, we extend the stability analysis to non-radial perturbations and prove Theorem~\ref{thm:nonrad:stab}. In Section~\ref{sec:ODE}, we analyze the ODE system governing Taylor coefficients at the origin and prove Theorem~\ref{thm:ODE_stab}. 
We emphasize that Sections~\ref{sec:stability}, \ref{sec:nonradial}, and~\ref{sec:ODE} may be read independently.
The appendices contain interpolation estimates and auxiliary algebraic computations.
%%%%%%%%%%%% 

\section{A family of exact self-similar profiles}
\label{sec:profiles}

The goal of this section is to show that for any $\NNN \in \Naturals$, we may find explicit similarity exponents
\[
\cxbar = \cxbar(\gamma,d,\NNN),
\qquad 
\cubar = \cxbar -1,
\qquad
\cbbar = \cbbar(\gamma,d,\NNN),
\] 
and exact globally self-similar solutions to~\eqref{eq:euler2}, given in terms of smooth profiles $(\bar \Sigma, \bar U, \bar B)$. That is, passing $\tau \to \infty$ in~\eqref{eq:euler3:old}, we seek solutions of
\begin{subequations} 
\label{eq:euler4}
\begin{align} 
(1-\cxbar) \bar\Sigma + \cxbar R \p_R \bar\Sigma +  \bar U \p_R \bar\Sigma  +  \alpha  \bar\Sigma ( \p_R \bar U + \tfrac{d-1}{R} \bar U)& = 0  ,  \label{eq:sigma4} \\
(1-\cxbar) \bar U + \cxbar R \p_R \bar U + \bar U \p_R  \bar U  +   \alpha \bar \Sigma \p_R \bar \Sigma - \tfrac{\alpha}{\gamma} \bar \Sigma^2 \tfrac{\p_R \bar B}{\bar B}& = 0  , \label{eq:u4} \\
- \cbbar \bar B + \cxbar R \p_R \bar B + \bar U \p_R  \bar B&=0 . \label{eq:s4}
\end{align} 
\end{subequations} 
We seek solutions $(\bar \Sigma, \bar U, \bar B)$ of \eqref{eq:euler4} which are smooth at $R=0$ and which grow sub-linearly as $R\to \infty$. Their existence and main properties are stated in Theorem~\ref{thm:main:profiles} below. This section is devoted solely to the analysis of system~\eqref{eq:euler4}, leading to the proof of Theorem~\ref{thm:main:profiles}.

\begin{remark}[\bf Globally self-similar solutions of the Euler equations]
\label{rem:globally:SS:exact}
Once a smooth solution $(\bar \Sigma,\bar U, \bar B)$ of~\eqref{eq:euler4} is found, for suitable values of $\cxbar$ and $\cbbar$, upon reversing the self-similar transformation discussed in Section~\ref{sec:self:similar:intro}, we directly obtain globally self-similar solutions $(\bar \sigma,\bar u^r, \bar b)$ of \eqref{eq:euler2}, defined as
\begin{subequations}
\label{eq:globally:SS:exact}
\begin{align}
\bar \sigma(r,t) &:=  (-t)^{\cxbar -1} \bar \Sigma\Bigl(\tfrac{r}{(-t)^{\cxbar}}\Bigr) ,
\\
\bar u^r(r,t) &:= (-t)^{\cxbar -1} \bar U\Bigl(\tfrac{r}{(-t)^{\cxbar}}\Bigr) , 
\\
\bar b(r,t) &:= (-t)^{\cbbar} \bar B\Bigl(\tfrac{r}{(-t)^{\cxbar}} \Bigr).
\end{align}
\end{subequations} 
The functions described in~\eqref{eq:globally:SS:exact} describe an exact solution of the Euler equations in radial symmetry~\eqref{eq:euler2} on $[-1,0) \times \Reals^d$, with initial datum $\bar \sigma(r,-1) = \bar \Sigma(r)$, $\bar u^r(r,-1) = \bar U(r)$, $\bar b(r,-1) = \bar B(r)$. This solution determines a smooth density $\bar \rho$ and smooth pressure $\bar p$  via~\eqref{eq:new:constitutive:relations}; the smooth radially symmetric velocity $\bar \uu = \vec{e}_r \bar u^r$  on $[-1,0) \times \Reals^d$. The solution $(\bar \rho,\bar \uu,\bar p)$ is the unique smooth solution of~\eqref{eq:euler:primary} with initial datum computed from $\mathsf {RHS}_{\eqref{eq:globally:SS:exact}}|_{t=-1}$ and~\eqref{eq:new:constitutive:relations}. This solution experiences an implosion singularity at $x_*=0$ and $t_*=0$, in the sense of Definition~\ref{def:implosion}.
\end{remark}

We now turn to the solvability of~\eqref{eq:euler4}. Two important remarks are in order.

\begin{remark}[\bf Non-isentropic implies $\cbbar \neq 0$]
\label{rem:cb:not:0}
When $\cbbar = 0$, then we observe that $\bar B = {\rm constant} $ is an explicit solution of \eqref{eq:s4}. Since for this solution $\p_R \bar B = 0$, the system~\eqref{eq:sigma4}--\eqref{eq:u4} becomes the self-similar \emph{isentropic} Euler system. As such, throughout this analysis we require that 
$\cbbar \neq 0$,
which precludes us from re-discovering isentropic implosion profiles.
\end{remark}

\begin{remark}[\bf $\cbbar \neq 0$ implies that smooth $\bar B$ and $\bar \Sigma$ must vanish at $R=0$]
\label{rem:B:Sigma:vanish:at:0}
Since $\bar U(0) = 0$,\footnote{This condition is necessary in order to ensure that $u(x,t)$ is smooth at $x=0$.} assuming that $|\partial_R \bar B(0)| < \infty$, we deduce from~\eqref{eq:s4} evaluated at $R=0$ that $\cbbar \bar B(0) = 0$. Since $\cbbar \neq 0$, we thus must have
$\bar B(0) = 0$.
The same argument applied to \eqref{eq:u4} evaluated at $R=0$,  and assuming that $|\partial_R \bar U(0)| + |\partial_R \bar \Sigma(0)| < \infty$, implies that $\lim_{R\to 0^+}\bar \Sigma^2(R) \frac{\p_R \bar B(R)}{\bar B(R)} = 0$. Assuming that $\bar B$ does not vanish to infinite order at $R=0$, meaning that $\bar B(R) \sim R^\nu$ for some $0 < \nu < \infty$, this implies that  
$ \bar \Sigma(0) = 0$.
\end{remark}

\subsection{Renormalization at \texorpdfstring{$R=0$}{R=0}}
We make three observations. First, the initial velocity $u(x,-1)$ is smooth at $|x|=0$ if and only if the Taylor series expansion of $u^r(r,-1)$ around $r=0$ contains only odd powers of $r$; as such, it is reasonable to expect that $\bar U(R) \sim R$ as $R \to 0^+$ and $U(R)/R$ is a function of $R^2$ for $|R| \ll 1$. Second, if we wish to ensure that the initial density $\rho(r,-1)$ does not vanish at $r=0$, by~\eqref{eq:new:constitutive:relations} and Remark~\ref{rem:B:Sigma:vanish:at:0} we obtain that $\bar \Sigma$ and $\bar B$ need to vanish \emph{at the same rate} as $R\to 0^+$; a quick inspection of \eqref{eq:euler4} suggests that it is reasonable to expect $\bar \Sigma(R), \bar B(R) \sim R$ as $R\to0^+$. Third, if we are to ensure that the initial density $\rho(r,-1)$ and the initial pressure $p(r,-1)$ are smooth around $r=0$, then by \eqref{eq:new:constitutive:relations} we obtain that $\bar \Sigma(R)/R$ and $\bar B(R)/R$ are functions of $R^2$ for $|R|\ll 1$. Based on these three observations, we make the following ansatz, whose goal is to renormalize the profiles $(\bar U,\bar \Sigma,\bar B)$ near $R=0$:
\begin{subequations}
\label{eq:profile:normalize}
\begin{equation}
\bar U(R) =: R  \bar V(R) ,
\qquad
\bar \Sigma(R) =: R \bar Q(R),
\qquad
\bar B(R) =: R \bar H(R) ,
\end{equation}
where the stationary renormalized profiles $(\bar V, \bar Q, \bar H)$ are smooth functions of $R^2$ for $|R|\ll 1$, and such that 
\begin{equation}
\bar V(0) < 0, \qquad \bar Q(0) > 0 , \qquad \bar H(0) > 0 .
\end{equation}
\end{subequations}
From~\eqref{eq:profile:normalize} and \eqref{eq:euler4} we obtain that $(\bar V, \bar Q, \bar H)$ must solve
\begin{subequations} 
\label{eq:euler5}
\begin{align} 
 \bigl( 1 + (1 + \alpha d) \bar V\bigr) \bar Q + ( \cxbar + \bar V)  R \p_R  \bar Q   +  \alpha    \bar Q  R \p_R  \bar V  & = 0  ,  \label{eq:sigma5} \\
\bigl(  \bar V + \bar V^2 + \tfrac{2\alpha^2}{\gamma} \bar Q^2 \bigr) + (\cxbar   +\bar V) R \p_R  \bar V   +   \alpha   \bar Q R \p_R  \bar Q  - \tfrac{\alpha}{\gamma}  \bar Q^2 \tfrac{R \p_R  \bar H }{  \bar H}& = 0  , \label{eq:u5} \\
\bigl(\cxbar - \cbbar + \bar V\bigr)   \bar H   + (\cxbar + \bar V) R \p_R  \bar H    &=0 . \label{eq:s5}
\end{align} 
\end{subequations} 

\begin{remark}[\bf The value of $\bar V(0)$, $\bar Q(0)$, and $\cbbar$ are determined]
Since we have assumed that $\bar Q(0) \neq 0$ and $|\partial_R \bar Q(0)| + |\partial_R \bar V(0)| < \infty$, upon evaluating \eqref{eq:sigma5} at $R=0$ we deduce that 
\begin{subequations}
\label{eq:values:at:zero:and:cb}
\begin{equation}
\bar V(0) = - \tfrac{1}{1 + \alpha d} .
\label{eq:bar:V:zero} 
\end{equation}
Since we have assumed that $\bar H(0) \neq 0$ and also $|\partial_R \bar H(0)| < \infty$, upon evaluating
\eqref{eq:u5} at $R=0$ we obtain
\begin{equation}
\bar Q(0)   = - \bar V(0) \sqrt{\tfrac{d \gamma}{2\alpha}} =    \tfrac{1}{1 + \alpha d}  \sqrt{\tfrac{d \gamma}{2\alpha}}
.
\label{eq:bar:Q:zero} 
\end{equation}
Lastly, by evaluating \eqref{eq:s5} at $R=0$ we deduce that 
\begin{equation}
\label{eq:cb:def}
\cbbar = \cxbar + \bar V(0) = \cxbar - \tfrac{1}{1 + \alpha d} .
\end{equation}
Since $\cbbar \neq 0$, we need to ensure that $\cxbar \neq - \bar V(0) = \frac{1}{1+\alpha d}$. Furthermore, we note that \eqref{eq:u5}--\eqref{eq:s5} remain the same if we multiply $\bar H$ by a constant, so it is natural that $\bar H(0) > 0$ is a free parameter in the problem; for simplicity, we let 
\begin{equation}
\label{eq:bar:H:zero} 
\bar H(0) = 1 \, .
\end{equation}
\end{subequations}
For compactness of notation, and for consistency with~\eqref{eq:Taylor:R=0:all} below, we shall henceforth denote
\[
\bar v_0 := \bar V(0) = - \tfrac{1}{1+\alpha d},
\qquad
\bar q_0 := \bar Q(0) = \tfrac{1}{1+\alpha d} \sqrt{\tfrac{d\gamma}{2\alpha}},
\qquad
\bar h_0 := \bar H(0) = 1.
\]
\end{remark}

\subsection{A closed system for \texorpdfstring{$\bar V$ and $\bar Q$}{bar V and bar Q}}
We may use \eqref{eq:s5} to express $\frac{R \partial_R \bar H}{\bar H}$ as a function of $\bar V$ alone, namely
\begin{equation}
\tfrac{R \p_R \bar H}{\bar H} = - \tfrac{\cxbar - \cbbar + \bar V}{\cxbar + \bar V} =  - \tfrac{\bar V - \bar v_0}{\cxbar + \bar V}
.
\label{eq:bar:H:RdR:bar:H}
\end{equation}
Inserting the above expression into the last term of \eqref{eq:u5} (which is justified as long as $\cxbar + \bar V(R) \neq 0$), and appealing to \eqref{eq:values:at:zero:and:cb}, leads us to a closed system obeyed by the unknowns $\bar V$ and $\bar Q$, namely 
\begin{subequations} 
\label{eq:euler6}
\begin{align} 
& \mathsf{P}_1[\bar V, \bar Q] 
+ ( \cxbar + \bar V)  R \p_R  \bar Q   
+  \alpha    \bar Q  R \p_R  \bar V   = 0  ,  \label{eq:sigma6} \\
&\mathsf{P}_2[\bar V, \bar Q]  
+   \alpha (\cxbar + \bar V)  \bar Q R \p_R  \bar Q
+ (\cxbar   +\bar V)^2 R \p_R  \bar V    = 0  . \label{eq:u6} 
\end{align} 
Here $\mathsf{P}_1[\bar V, \bar Q](R)$ and $\mathsf{P}_2[\bar V,\bar Q](R)$ are two explicit polynomials in $\bar V$ and $\bar Q$ (quadratic, respectively cubic), which vanish at $R=0$; specifically, we have
\begin{align}
\mathsf{P}_1[\bar V, \bar Q] 
&:= \bigl(1 + (1 + \alpha d)\bar V  \bigr) \bar Q  = (1 + \alpha d) \bigl(\bar V - \bar v_0 \bigr) \bar Q
,
\label{eq:P1:def}
\\
\mathsf{P}_2[\bar V, \bar Q] 
&:= 
(\cxbar + \bar V) \bigl(  \bar V + \bar V^2 + \tfrac{2\alpha^2}{\gamma} \bar Q^2 \bigr)
+ \tfrac{\alpha}{\gamma}  \bar Q^2 (\bar V - \bar v_0)
\notag\\
&\ = \bigl(\bar V - \bar v_0 \bigr)  \bigl( (\cxbar + \bar V)  (1 + \bar V + \bar v_0) 
+  \tfrac{\alpha}{\gamma}   \bar Q^2 \bigr)
+ \tfrac{2\alpha^2}{\gamma} \bigl(\bar Q - \bar q_0\bigr) (\cxbar + \bar V) \bigl(\bar Q + \bar q_0\bigr)  
.
\label{eq:P2:def}
\end{align}
\end{subequations}

\subsection{Power series for \texorpdfstring{$\bar V$ and $\bar Q$}{bar V and bar Q} near \texorpdfstring{$R=0$}{R=0} for a discrete family of similarity exponents}

We aim to find a power series solution for~\eqref{eq:euler6} near $R=0$, by postulating that 
\begin{equation}
\bar V(R) = \bar v_0 + \sum_{n\geq 1}  \bar v_n R^{2n},
\qquad
\bar Q(R) =  \bar q_0 + \sum_{n\geq 1}  \bar q_n  R^{2n},
\label{eq:V:Q:R=0}
\end{equation}
where $\bar v_0$ and $\bar q_0$ are given by~\eqref{eq:bar:V:zero}--\eqref{eq:bar:Q:zero}, and the coefficients $\{\bar v_n\}_{n\geq 1}$ and $\{\bar q_n\}_{n\geq 1}$ are to be determined from a recursion relation.

\subsubsection{Recursion relation for the power series coefficients}
Inserting~\eqref{eq:V:Q:R=0} into~\eqref{eq:sigma6}, leads to the recursion relation:
\begin{subequations}
\label{eq:recursion:1}
\begin{equation}
2 n (\cxbar + \bar v_0) \bar q_n
+ 
\bigl(1+ \alpha d + 2  \alpha  n \bigr) \bar q_0  \bar v_n 
= 
\mathsf{F}_1(n)
,
\label{eq:recursion:1:a}
\end{equation}
for all $n\geq 1$, where we have defined
\begin{equation}
\mathsf{F}_1(n):= - {\mathbf 1}_{n\geq 2} \sum_{m+j=n, m,j\geq 1}  \bar q_{m} \bar v_{j} \bigl(1+\alpha d + 2n - 2(1-\alpha) j   \bigr) 
.
\label{eq:recursion:1:b} 
\end{equation}
\end{subequations}
Similarly, inserting~\eqref{eq:V:Q:R=0} into~\eqref{eq:u6}, leads to the recursion relation:
\begin{subequations}
\label{eq:recursion:2}
\begin{equation}
\alpha (\cxbar + \bar v_0)  \bar q_0 \bigl( \tfrac{4\alpha}{\gamma}  + 2 n  \bigr) \bar q_n
+ 
\Bigl( 2n (\cxbar + \bar v_0)^2 +  \bigl(1 + 2 \bar v_0\bigr) (\cxbar + \bar v_0)  + \tfrac{\alpha}{\gamma} \bar q_0^2 \Bigr)  \bar v_n 
= 
\mathsf{F}_2(n)
,
\label{eq:recursion:2:a}
\end{equation}
for all $n\geq 1$, where we have defined
\begin{align}
\mathsf{F}_2(n)
&:= - {\mathbf 1}_{n\geq 2} \sum_{m+j=n, m,j\geq 1}   \bar v_{m} \bar v_{j} \bigl(1+ \cxbar  + 3 \bar v_0 + 4 j (\cxbar + \bar v_0)   \bigr) 
\notag\\
&\quad - {\mathbf 1}_{n\geq 2} \sum_{m+j=n, m,j\geq 1}   \bar v_{m} \bar q_{j}  2\alpha \bar q_0  (j+1)  
 - {\mathbf 1}_{n\geq 2} \sum_{m+j=n, m,j\geq 1}  \bar q_{m} \bar q_{j} (\cxbar + \bar v_0) \bigl( \tfrac{2\alpha^2}{\gamma} +  2 \alpha j   \bigr) 
\notag\\
&\quad - {\mathbf 1}_{n\geq 3} \sum_{m+\ell+j=n, m,\ell,j\geq 1} \bar v_{m}  \bigl( \bar v_{\ell} \bar v_j + \alpha \bar q_{\ell} \bar q_j\bigr) (2j+1).
\label{eq:recursion:2:b} 
\end{align}
\end{subequations}

\subsubsection{A discrete family of similarity exponents}
\label{sec:similarity:expo:matrix}
Upon introducing the matrix
\begin{equation}
\MMM_n :=
\begin{bmatrix}
2n(\cxbar + \bar v_0) 
&  (1+ \alpha d + 2  \alpha  n ) \bar q_0  
\\
\alpha (\cxbar + \bar v_0)  \bar q_0 \bigl( \tfrac{4\alpha}{\gamma} + 2 n  \bigr)  
& 2n (\cxbar + \bar v_0)^2 + \bigl(1 + 2 \bar v_0  \bigr) (\cxbar + \bar v_0) + \tfrac{\alpha}{\gamma} \bar q_0^2
\end{bmatrix}	
\label{eq:Mn:def}
\end{equation}
we may summarize \eqref{eq:recursion:1}--\eqref{eq:recursion:2} as
\begin{equation}
 \MMM_n \begin{bmatrix} \bar q_n \\ \bar v_n \end{bmatrix} = \begin{bmatrix}\mathsf{F}_1(n) \\ \mathsf{F}_2(n) \end{bmatrix}
 ,
 \label{eq:recursion:3}
\end{equation}
for all $n\geq 1$.

We notice at this stage that by definition $\mathsf{F}_1(1) = \mathsf{F}_2(1)  = 0$, and thus if $\MMM_1$ is invertible, then we must have $\bar q_1 = \bar v_1 = 0$. In turn, this implies $\mathsf{F}_1(2) = \mathsf{F}_2(2)  = 0$, and so if $\MMM_2$ is invertible, then we must have $\bar q_2 = \bar v_2= 0$, and so on. By induction, we would obtain that $\bar q_n = \bar v_n = 0$ for all $n\geq 1$, which is the \emph{trivial solution}. This situation may only be avoided if we have that $\det(\MMM_{\NNN}) = 0$ for some $\NNN \in \Naturals$. This imposes a constraint  on the \emph{admissible values} for the similarity parameter $\cxbar$, as follows. Using~\eqref{eq:bar:V:zero} and~\eqref{eq:bar:Q:zero} we may compute
\begin{align}
\det(\MMM_n) 
&= 4n^2 (\cxbar + \bar v_0) \Bigl(
 (\cxbar + \bar v_0)^2 
+ \tfrac{1}{2n} (\cxbar + \bar v_0) \bigl(1 + 2 \bar v_0 \bigr) - \bar v_0^2 \bigl( \tfrac{\alpha \gamma d}{2}   +  \mathsf{E}_n \bigr) \Bigr)
\label{eq:det:Mn:def}
,
\end{align}
where $\mathsf{E}_n = \mathsf{E}_n(\gamma,d)$ is given explicitly by the expression
\begin{equation}
\mathsf{E}_n
:=  
 \tfrac{\alpha \gamma d (d+2)}{4n} 
+ \tfrac{\alpha d (1+ \alpha d  ) }{2 n^2}   
.
\label{eq:En:def}
\end{equation}
We note that $\mathsf{E}_n \geq 0$, $\{\mathsf{E}_n\}_{n\geq 1}$ is monotone decreasing in $n$, and $\mathsf{E}_n  = \OO(\frac{1}{n})$ as $n\to \infty$.

At this stage we recall from \eqref{eq:cb:def} and the requirements $\cbbar \neq 0$, that $\cxbar + \bar v_0 \neq 0$. From \eqref{eq:det:Mn:def} and~\eqref{eq:bar:V:zero} it follows that $\det(\MMM_n) =0$ only if 
\begin{align}
\cxbar + \bar v_0 
&=  - \tfrac{1}{4n} (1 + 2 \bar v_0) \pm   \sqrt{\tfrac{1}{16 n^2} (1+2 \bar v_0)^2 +   \bar v_0^2 ( \tfrac{\alpha \gamma d}{2}  + \mathsf{E}_n) }
\notag\\
&= - \bar v_0 \Bigl( \tfrac{1-\alpha d}{4n}   \pm  \sqrt{ \tfrac{\alpha \gamma d}{2} + \mathsf{E}_n  + \tfrac{(1-\alpha d)^2}{16 n^2}  } \Bigr) .
\label{eq:cx:admissible:1}
\end{align}
The branch with the negative sign in front of the square root in~\eqref{eq:cx:admissible:1} leads to negative values of $\cxbar$ (e.g.~for $n\gg 1$) which is not a relevant parameter regime for implosion singularities. Therefore, we consider the branch with the positive sign in front of the square root in~\eqref{eq:cx:admissible:1}, and define:
\begin{definition}[\bf Admissible values of the similarity exponents]
\label{def:exponents}
Let $d \geq 1$, $\gamma >1$, and $\alpha =\frac{\gamma-1}{2}$. For each $\NNN\in \Naturals$, define the $\NNN^{\rm th}$ admissible similarity exponent by  
\begin{equation}
\cxstar(d,\gamma,\NNN)
:= \tfrac{1}{1+\alpha d} \Bigl(1  +  \sqrt{ \tfrac{\alpha \gamma d}{2}  + \mathsf{E}_\NNN  + \tfrac{(1-\alpha d)^2}{16 \NNN^2}  } + \tfrac{1-\alpha d}{4\NNN}   \Bigr),
\label{eq:cx:admissible}
\end{equation}
where $\mathsf{E}_{\NNN}$ is as defined in~\eqref{eq:En:def}.
Moreover, according to~\eqref{eq:cb:def} we define 
\begin{equation}
\cbstar(d,\gamma,\NNN)
:= \tfrac{1}{1+\alpha d} \Bigl( \sqrt{ \tfrac{\alpha \gamma d}{2} + \mathsf{E}_\NNN  + \tfrac{(1-\alpha d)^2}{16 \NNN^2}  } + \tfrac{1-\alpha d}{4\NNN}   \Bigr).
\label{eq:cb:admissible}
\end{equation}
We emphasize that the values of $\cxstar$ and $\cbstar$ defined above are explicit and exact.
\end{definition}

\begin{example}[The ``ground state'' at $\NNN = 1$]
When $\NNN=1$, by appealing to~\eqref{eq:En:def}  and~\eqref{eq:cx:admissible}, we arrive at
\begin{align}
 \cxstar(d,\gamma,1) &=  \tfrac{1}{1+\alpha d} \Bigl(1  +  \sqrt{   \tfrac{4 (\alpha d)^2 +  (2 +  d + 4 \gamma) \alpha d  }{4} + \tfrac{(1-\alpha d)^2}{16} } + \tfrac{1-\alpha d}{4} \Bigr) 
 .
 \label{eq:cxstar:ground:state}
\end{align}
\end{example}

\begin{example}[\bf Similarity exponents for a monatomic gas in three space dimensions]
\label{ex:exponents:monatomic:3D}
Consider the expression for $\cxstar$ in~\eqref{eq:cx:admissible} for $\gamma = \frac 53$ (a monatomic gas), so that $\alpha = \frac 13$, and $d=3$ (the three-dimensional case). Then, upon drastic simplification one may show that
\begin{subequations}
\begin{align}
 \cxstar(3,\tfrac{5}{3},\NNN) &=  \tfrac 12 + \tfrac{1}{2} \sqrt{ \tfrac{5}{6}  +  \tfrac{25 }{12 \NNN}  +   \tfrac{1}{\NNN^2} } 
 \\
 \cbstar(3,\tfrac{5}{3},\NNN) &= 
 \cxstar(3,\tfrac{5}{3},\NNN) - \tfrac{1}{2} = 
\tfrac{1}{2} \sqrt{ \tfrac{5}{6}  +  \tfrac{25 }{12 \NNN}   +  \tfrac{1}{\NNN^2}   }
 .
\end{align}
\end{subequations}
These are both monotonically decreasing sequences in $\NNN$, with  $ \cxstar(3,\frac{5}{3},1) \approx 1.48953$, $\cbstar(3,\frac 53,1) \approx 0.98953$, $\cxstar(3,\frac{5}{3},\infty) \approx 0.956435$, and $\cbstar(3,\frac{5}{3},\infty) \approx 0.456435$.
\end{example}

\begin{lemma}[\bf Properties of the admissible similarity exponents]
\label{lem:properties:exponents}
The exponents $\cxstar$ and $\cbstar$ defined in \eqref{eq:cx:admissible} and~\eqref{eq:cb:admissible} satisfy the following properties:
\begin{itemize}[leftmargin=1em]
\item Fix $\NNN \geq 1$ and let $\cxbar = \cxstar(d,\gamma,\NNN)$ in the definition of $\MMM_n$ in~\eqref{eq:Mn:def}. Then, we have that $\det(\MMM_n) = 0$ if and only if $n = \NNN$.  The eigenvector corresponding to the zero-eigenvalue of $\MMM_{\NNN}$ is parallel to the vector $(- \bar q_0 \frac{1+\alpha d + 2  \alpha  \NNN }{2\NNN (\cxbar + \bar v_0)},1)$. 

\item  For $d \in \{1,2,3\}$ and $1< \gamma \leq 2 d+1$,  the sequences $\cxstar(d,\gamma,\NNN)$ and $\cbstar(d,\gamma,\NNN)$ are monotone decreasing in $\NNN$, and accumulate  at $\cxstar(d,\gamma,\infty) = \tfrac{1}{1 +\alpha d} (1 +   \sqrt{\frac{\alpha \gamma d}{2}})$, respectively $\cbstar(d,\gamma,\infty) = \tfrac{1}{1 +\alpha d}  \sqrt{\frac{\alpha \gamma d}{2}}$. 

\item  For any $\NNN\geq 1$ and $1 < \gamma \leq 2d+1$, we have that $\cxstar(d,\gamma,\NNN) + \bar v_0 - \alpha \bar q_0 \geq \frac{1 + \frac 23 \alpha d}{4 \NNN(1+\alpha d)}$.
\end{itemize}
\end{lemma}
\begin{proof}[Proof of Lemma~\ref{lem:properties:exponents}]
We prove the claims in each bullet separately. 

\noindent \textit{Proof of the first bullet.}
Let $n\geq 1$, potentially different from $\NNN$. From~\eqref{eq:det:Mn:def}, setting $X := \cxbar + \bar v_0$ with $\cxbar = \cxstar(d,\gamma,\NNN)$, we have that $\det(\MMM_n) = 0$ (with $X \neq 0$) if and only if
\begin{equation}
\label{eq:det:zero:condition}
X^2 + \tfrac{1+2\bar v_0}{2n}X - \bar V^2(0)\Bigl(\tfrac{\alpha\gamma d}{2} + \mathsf{E}_n\Bigr) = 0.
\end{equation}
The two roots of this quadratic in $X$ are
\[
X_n^\pm = \tfrac{1}{1+\alpha d}\left(\tfrac{1-\alpha d}{4n} \pm \sqrt{\tfrac{(1-\alpha d)^2}{16n^2} + \tfrac{\alpha\gamma d}{2} + \mathsf{E}_n}\right).
\]
By definition~\eqref{eq:cx:admissible}, we have $X_n^+ = \cxstar(d,\gamma,n) + \bar v_0$. Moreover, $X_n^- < 0$ for all $n \geq 1$, since the square root term exceeds $\bigl|\frac{1-\alpha d}{4n}\bigr|$. For $\cxbar = \cxstar(d,\gamma,\NNN)$, we have $ \cxbar + \bar v_0 = X_\NNN^+ > 0$. Thus, we have $\det(\MMM_n) = 0$ if and only if $X_\NNN^+ = X_n^+$ or $X_\NNN^+ = X_n^-$. Since $X = X_\NNN^+ > 0$ and $X_n^- < 0$ for all $n$, $\det(\MMM_n) = 0$ if and only if $X_\NNN^+ = X_n^+$, which is equivalent to $\cxstar(d,\gamma,\NNN)=  \cxstar(d,\gamma,n)$. Due to the monotonicity claimed in the second bullet, this equality holds if and only if $n = \NNN$.

For the eigenvector computaton, when $\det(\MMM_\NNN) = 0$, the null space is determined by the first row of $\MMM_\NNN$:
\[
2\NNN(\cxbar + \bar v_0)\bar q_\NNN + (1+\alpha d + 2\alpha\NNN)\bar q_0\bar v_\NNN = 0,
\]
which gives $\bar q_\NNN = -\frac{(1+\alpha d + 2\alpha\NNN)\bar q_0}{2\NNN(\cxbar + \bar v_0)}\bar v_\NNN$. Thus the eigenvector $(\bar q_\NNN, \bar v_\NNN)$ is parallel to $\bigl(-\bar q_0\frac{1+\alpha d + 2\alpha\NNN}{2\NNN(\cxbar + \bar v_0)}, 1\bigr)$.

\noindent \textit{Proof of the second bullet.}
The statements about the limit as $\NNN \to \infty$ are immediate consequences of $\mathsf{E}_\NNN \to 0$ as $\NNN\to \infty$.
Moreover, the monotonicity of $\cbstar$ is equivalent to that of $\cxstar$. Let 
\[
h(\NNN) := \sqrt{f(\NNN)} + \tfrac{1-\alpha d}{4\NNN},
\qquad
f(\NNN) := \tfrac{\alpha\gamma d}{2} + \mathsf{E}_\NNN + \tfrac{(1-\alpha d)^2}{16\NNN^2}.
\]
In view of~\eqref{eq:cx:admissible}, in order to prove that $\cxstar$ is strictly decreasing in $\NNN$, it suffices to show that $h(\NNN)$ is strictly decreasing in $\NNN$. It is clear that $f(\NNN)$ is strictly decreasing in $\NNN$. Since for $\alpha d \leq 1$ we have $\frac{1-\alpha d}{4\NNN} \geq 0$ is decreasing in $\NNN$, in this case $h(\NNN)$ is the sum of two decreasing functions, hence decreasing.

It remains to consider the case $\alpha d > 1$. Since $\alpha \in(0,d]$, we only need to consider $d\in \{2,3\}$. For $\NNN \geq 1$, we have
\[
h(\NNN) - h(\NNN+1) = \tfrac{f(\NNN) - f(\NNN+1)}{\sqrt{f(\NNN+1)} + \sqrt{f(\NNN)}} -  \tfrac{\alpha d-1}{4 \NNN (\NNN+1)} .
\]
We need to show that the above expression is positive. Recalling the definition of $\mathsf{E}_{\NNN}$ in~\eqref{eq:En:def}, we have that 
\[
 \tfrac{f(\NNN) - f(\NNN+1)}{\sqrt{f(\NNN+1)} + \sqrt{f(\NNN)}} 
 >
 \tfrac{1}{2 \sqrt{f(\NNN)}} \left(  \tfrac{\alpha \gamma d (d+2)}{4\NNN (\NNN+1)} 
+ \tfrac{\alpha d (1+ \alpha d) (2\NNN+1)}{2 \NNN^2(\NNN+1)^2} + \tfrac{(\alpha d -1)^2}{16} \tfrac{2 \NNN+1}{\NNN^2 (\NNN+1)^2} \right) 
\]
and
\[
2 \sqrt{f(\NNN)} 
\leq \sqrt{2\alpha\gamma d } + \sqrt{ \tfrac{\alpha \gamma d (d+2)}{\NNN}} 
+ \tfrac{1}{\NNN} \sqrt{2 \alpha d (1+ \alpha d  )} + \tfrac{\alpha d-1}{2\NNN}
\]
Therefore, the inequality $h(\NNN)>h(\NNN+1)$ is true if we are able to show that 
\[
\tfrac{\alpha \gamma d (d+2)}{2} 
+ \tfrac{\alpha d (1+ \alpha d) (2\NNN+1)}{\NNN(\NNN+1)} 
+ \tfrac{(\alpha d -1)^2}{8} \tfrac{2 \NNN+1}{\NNN (\NNN+1)}
\geq  \tfrac{\alpha d-1}{2} \left( \sqrt{2\alpha\gamma d } + \sqrt{ \tfrac{\alpha \gamma d (d+2)}{\NNN}} 
+ \tfrac{1}{\NNN} \sqrt{2 \alpha d (1+ \alpha d  )} + \tfrac{\alpha d-1}{2\NNN} \right)
.
\]
In turn, the above inequality holds if we are able to show that 
\[
\tfrac{3 \alpha \gamma d (d+2)}{8} 
+ \tfrac{8 \alpha d (1+ \alpha d)- 5 (\alpha d-1)^2}{8 \NNN} 
+ \tfrac{8\alpha d (1+ \alpha d) - (\alpha d -1)^2 }{8 \NNN(\NNN+1)} 
\geq  
\tfrac{\alpha d-1}{2} \sqrt{2\alpha\gamma d }  
.
\]
Since $\alpha d > 1$, the above inequality does indeed hold for any $\NNN\geq 1$, if we are able to establish
\[
3\sqrt{d+2} (d+2) \sqrt{\alpha d} 
\geq  
4\sqrt{2d}(\alpha d-1)
.
\]
Since $x=\sqrt{\alpha d} \in (1,d]$ (recall that $\frac{1}{d} < \alpha \leq d$), the above inequality clearly holds true for both $d=2$ (where it becomes $24 x
\geq  8(x^2-1)$, for $x\in(1,2]$) and $d=3$ (where it becomes $15 \sqrt{5} x
\geq  4\sqrt{6}(x^2-1)$ , for $x\in(1,3]$).

\noindent \textit{Proof of the third bullet.}
Recall that $\cxstar(d,\gamma,\infty) = \frac{1}{1 +\alpha d} (1 + \sqrt{\frac{\alpha \gamma d}{2}})= \frac{1}{1 +\alpha d} ( 1 + \alpha \sqrt{\frac{ \gamma d}{2\alpha}}) = - \bar v_0 + \alpha \bar q_0$. The claimed lower bound follows by bounding from below the quantity $\cxstar(d,\gamma,\NNN) - \cxstar(d,\gamma,\infty)$. From~\eqref{eq:cx:admissible} this difference is explicitly given as
\begin{align*}
(1+\alpha d) \bigl(\cxstar(d,\gamma,\NNN) - \cxstar(d,\gamma,\infty)\bigr) 
&= \sqrt{\tfrac{\alpha \gamma d}{2} + \mathsf{E}_\NNN + \tfrac{1}{16 \NNN^2} (1-\alpha d)^2 } + \tfrac{1-\alpha d}{4\NNN} - \sqrt{\tfrac{\alpha \gamma d}{2} } \\
&= \tfrac{1-\alpha d}{4\NNN} + \tfrac{\mathsf{E}_\NNN + \frac{1}{16 \NNN^2} (1-\alpha d)^2}{\sqrt{\tfrac{\alpha \gamma d}{2} + \mathsf{E}_\NNN + \frac{1}{16 \NNN^2} (1-\alpha d)^2 } + \sqrt{\tfrac{\alpha \gamma d}{2} }} .
\end{align*}
Recalling that $\mathsf{E}_{\NNN}$ is monotone decreasing in $\NNN$ and positive, and using the explicit formula for the first term on the right side of~\eqref{eq:En:def}, we may thus lower bound 
\begin{align*}
&\tfrac{1}{\alpha d} \left( 4 \NNN(1+\alpha d) \bigl(\cxstar(d,\gamma,\NNN) - \cxstar(d,\gamma,\infty)\bigr) - 1 - \tfrac 23 \alpha d \right)\\
&= - \tfrac 53 + \tfrac{4 \NNN \mathsf{E}_\NNN + \frac{(1-\alpha d)^2}{4 \NNN} }{\alpha d(\sqrt{\frac{\alpha \gamma d}{2} + \mathsf{E}_\NNN + \frac{(1-\alpha d)^2}{16 \NNN^2}  } + \sqrt{\tfrac{\alpha \gamma d}{2} })}  
\geq - \tfrac 53  + \tfrac{ \gamma   (d+2) }{ \sqrt{\frac{\alpha \gamma d}{2} + \mathsf{E}_1 + \frac{(1-\alpha d)^2}{16}  } + \sqrt{\frac{\alpha \gamma d}{2} }} 
 =: F(\alpha,d),
\end{align*}
where $\mathsf{E}_1 =  \frac{\alpha \gamma d (d+2)}{4} 
+ \frac{\alpha d (1+ \alpha d  ) }{2} $ (cf.~\eqref{eq:En:def}).
A tedious but elementary analysis shows that the function $F(\alpha,d)$ is monotone decreasing with respect to $\alpha \in (0,d]$ for all $d\in\{1,2,3\}$; therefore, its minimum is attained at $\alpha =d$. An explicit computation shows that 
\[
F(d,d) = - \tfrac 53  + \tfrac{ (2d+1)   (1+\frac{2}{d}) }{ \sqrt{ \frac{(2d+1)(d+4)}{4} + \frac{(1 + 3 d^2)^2}{16 d^2} } + \sqrt{\frac{2d+1}{2} }} 
\geq \tfrac 14
\]
for $d\in\{1,2,3\}$, which concludes the proof of the lower bound claimed in the third bullet.
\end{proof}

\subsubsection{Solvability of the recursion relation}
We return to the solvability of the recursion relation~\eqref{eq:recursion:3} for the unknown sequences $\{\bar v_n\}_{n\geq 1}$ and $\{\bar q_n\}_{n\geq 1}$. As discussed earlier, \eqref{eq:recursion:3} has a nontrivial solution only if $\cx$ is chosen to be equal to one of the similarity exponents defined in~\eqref{eq:cx:admissible}; we may thus prove:

\begin{proposition}[\bf Solvability of the recursion corresponding to $R=0$]
\label{prop:solve:recursion:R=0}
Fix an $\NNN\geq 1$, and let $\cxbar := \cxstar(d,\gamma,\NNN)$\footnote{We will simply write $\cxbar$ and suppress the dependence on $d,\gamma,\NNN$, as these are now fixed.} in the matrix $\MMM_n$ (see~\eqref{eq:Mn:def}), in the forcing terms $\mathsf{F}_1(n)$ (see~\eqref{eq:recursion:1:b}), and $\mathsf{F}_2(n)$ (see~\eqref{eq:recursion:2:b}). Define:
\begin{itemize}[leftmargin=1em]
\item For $n \neq k\NNN$ with $k\in  \Naturals$, let $\bar q_n = \bar v_n = 0$.
\item For $n = \NNN$, define $(\bar q_{\NNN}, \bar v_{\NNN})$ to be an eigenvector of $\MMM_{\NNN}$ corresponding to the zero eigenvalue, more precisely, let 
\begin{equation}
\begin{bmatrix} 
\bar q_{ \NNN} \\ 
\bar v_{ \NNN} 
\end{bmatrix} 
:=  \begin{bmatrix} 
- \bar q_0 \frac{1+\alpha d + 2 \alpha \NNN }{2\NNN (\cxbar + \bar v_0)} \\ 
1
\end{bmatrix}
.
\label{eq:first:nontrivial:Taylor:coefficient}
\end{equation}
\item For $n = k\NNN$ with $k \in \Naturals$ and $k\geq 2$, the matrix $\MMM_{k \NNN}$ is invertible, we let 
\begin{equation}
  \begin{bmatrix} \bar q_{k \NNN} \\ \bar v_{k \NNN} \end{bmatrix} := \MMM_{k \NNN}^{-1} \begin{bmatrix}\mathsf{F}_1(k \NNN) \\ \mathsf{F}_2(k \NNN) \end{bmatrix}
 .
\end{equation}
\end{itemize}  
Then,  the above defined sequences $\{\bar q_n\}_{n \geq 1}$ and $\{\bar v_n\}_{n\geq 1}$ solve the recursion relation~\eqref{eq:recursion:3}, equivalently \eqref{eq:recursion:1}--\eqref{eq:recursion:2}.
\end{proposition}
\begin{proof}[Proof of Proposition~\ref{prop:solve:recursion:R=0}]
The proof proceeds by verifying that the recursion relation~\eqref{eq:recursion:3} is satisfied for each $n \geq 1$. We consider three cases.

\medskip
\noindent \textit{Case 1: $n$ is not a multiple of $\NNN$.}
We claim that $\mathsf{F}_1(n) = \mathsf{F}_2(n) = 0$. Indeed, from~\eqref{eq:recursion:1:b}, the sum defining $\mathsf{F}_1(n)$ runs over pairs $(m,j)$ with $m + j = n$ and $m, j \geq 1$. Since $n$ is not a multiple of $\NNN$, at least one of $m$ or $j$ is not a multiple of $\NNN$, and hence at least one of $\bar q_m$, $\bar v_j$ vanishes by definition. Therefore every term in the sum is zero, giving $\mathsf{F}_1(n) = 0$. The same argument applied to~\eqref{eq:recursion:2:b} shows $\mathsf{F}_2(n) = 0$.

Since $n \neq \NNN$, by the first bullet of Lemma~\ref{lem:properties:exponents}, we have $\det(\MMM_n) \neq 0$, so $\MMM_n$ is invertible. The recursion~\eqref{eq:recursion:3} then reads $\MMM_n (\bar q_n, \bar v_n)^T = (0,0)^T$, whose unique solution is $\bar q_n = \bar v_n = 0$. This matches the definition in the first bullet of the proposition.

\medskip
\noindent \textit{Case 2: $n = \NNN$.}
Since $\bar q_m = \bar v_m = 0$ for all $1 \leq m < \NNN$ (as these $m$ are not multiples of $\NNN$), every term in the sums defining $\mathsf{F}_1(\NNN)$ and $\mathsf{F}_2(\NNN)$ vanishes. Indeed, the constraint $m + j = \NNN$ with $m, j \geq 1$ forces $m, j < \NNN$. Thus $\mathsf{F}_1(\NNN) = \mathsf{F}_2(\NNN) = 0$.

By the first bullet of Lemma~\ref{lem:properties:exponents}, $\det(\MMM_\NNN) = 0$, and the null space of $\MMM_\NNN$ is one-dimensional, spanned by the vector $\bigl(-\bar q_0 \frac{1+\alpha d + 2\alpha\NNN}{2\NNN(\cxbar + \bar v_0)}, 1\bigr)$. The recursion~\eqref{eq:recursion:3} at $n = \NNN$ reads $\MMM_\NNN (\bar q_\NNN, \bar v_\NNN)^T = (0,0)^T$, which is satisfied by any vector in the null space of $\MMM_\NNN$. The choice~\eqref{eq:first:nontrivial:Taylor:coefficient} is precisely such a vector.

\medskip
\noindent \textit{Case 3: $n = k\NNN$ for some $k \geq 2$.}
We proceed by induction on $k$. Assume that $(\bar q_{m\NNN}, \bar v_{m\NNN})$ have been defined for all $1 \leq m < k$ and satisfy the recursion. The forcing terms $\mathsf{F}_1(k\NNN)$ and $\mathsf{F}_2(k\NNN)$ depend only on $\{\bar q_{m\NNN}, \bar v_{m\NNN}\}_{m=1}^{k-1}$ (since the sums in~\eqref{eq:recursion:1:b}--\eqref{eq:recursion:2:b} run over indices strictly less than $k\NNN$, and non-multiples of $\NNN$ contribute zero). Thus $\mathsf{F}_1(k\NNN)$ and $\mathsf{F}_2(k\NNN)$ are well-defined.

Since $k\NNN \neq \NNN$ for $k \geq 2$, by the first bullet of Lemma~\ref{lem:properties:exponents}, we have $\det(\MMM_{k\NNN}) \neq 0$, so $\MMM_{k\NNN}$ is invertible. The unique solution to~\eqref{eq:recursion:3} at $n = k\NNN$ is therefore
\[
\begin{bmatrix} \bar q_{k\NNN} \\ \bar v_{k\NNN} \end{bmatrix} = \MMM_{k\NNN}^{-1} \begin{bmatrix} \mathsf{F}_1(k\NNN) \\ \mathsf{F}_2(k\NNN) \end{bmatrix},
\]
which is precisely the definition given in the third bullet of the proposition. 
\end{proof}

\subsubsection{Bounds for the power series coefficients}
With Proposition~\ref{prop:solve:recursion:R=0}, we return to the power series ansatz~\eqref{eq:V:Q:R=0}, and establish a bound on the  radius of convergence. The heavy lifting is done by the following result:

\begin{proposition}[\bf Bounds for power series coefficients]
\label{prop:power:series:bounds}
Let $\{\bar q_n\}_{n \geq 1}$ and $\{\bar v_n\}_{n\geq 1}$ be as defined in Proposition~\ref{prop:solve:recursion:R=0}. Then, 
there exist constants $ \mathsf{C}_{\eqref{eq:power:series:bounds}} = \mathsf{C}_{\eqref{eq:power:series:bounds}}(\gamma,d,\NNN) \geq 1$ and $0 < \mathsf{C}^\prime_{\eqref{eq:power:series:bounds}} = \mathsf{C}^\prime_{\eqref{eq:power:series:bounds}}(\gamma,d,\NNN) < 1$, such that 
\begin{subequations}
\label{eq:power:series:bounds}
\begin{align}
 |\bar q_{k \NNN}| &\leq  \mathsf{C}^{\prime}_{\eqref{eq:power:series:bounds}}  (\mathsf{C}_{\eqref{eq:power:series:bounds}})^k  k^{-\frac 32}
 \label{eq:power:series:bounds:q}
 ,
 \\
 |\bar v_{k \NNN}| &\leq   \mathsf{C}^{\prime}_{\eqref{eq:power:series:bounds}}    (\mathsf{C}_{\eqref{eq:power:series:bounds}})^k  k^{-\frac 32}
 \label{eq:power:series:bounds:v}
 .
\end{align}
\end{subequations}
for each $k\in \Naturals$.  
\end{proposition}

\begin{proof}[Proof of Proposition~\ref{prop:power:series:bounds}]
For compactness we omit the subscript for the constants $\mathsf{C}$ and $\mathsf{C}^\prime$.
For $k=1$ the bound~\eqref{eq:power:series:bounds} holds as soon as $\mathsf{C}^{\prime} \mathsf{C}  \geq \max\{1,  |\bar q_0| \frac{1+\alpha d + 2\alpha \NNN }{2\NNN (\cxbar + \bar v_0)} \}$.
For $k\geq 2$, from~\eqref{eq:Mn:def} and~\eqref{eq:det:Mn:def} we have that 
\begin{align}
\MMM_{k \NNN}^{-1} 
&:=
\tfrac{1}{4 k^2 \NNN^2 (\cxbar + \bar v_0) \bigl(
 (\cxbar + \bar v_0)^2 
+ \tfrac{1}{2k \NNN} (\cxbar + \bar v_0) \bigl(1 + 2 \bar v_0 \bigr) - \bar V^2(0) \bigl(\tfrac{\alpha \gamma d}{2}   +  \mathsf{E}_{k \NNN} \bigr) \bigr)}
\notag\\
&\quad \times
\begin{bmatrix}
2 k \NNN (\cxbar + \bar v_0)^2 + \bigl(1 + 2 \bar v_0  \bigr) (\cxbar + \bar v_0) + \tfrac{\alpha}{\gamma} \bar Q^2(0)
&  -(1+ \alpha d + 2 \alpha k \NNN  ) \bar q_0  
\\
- \alpha (\cxbar + \bar v_0)  \bar q_0 \bigl( \tfrac{4\alpha}{\gamma} + 2 k \NNN  \bigr)  
& 2k \NNN (\cxbar + \bar v_0) 
\end{bmatrix}	.
\label{eq:Mn:inverse:def}
\end{align}
From \eqref{eq:Mn:inverse:def} we deduce that there exists a constant $\mathsf{C}_{\MMM} = \mathsf{C}_{\MMM}(\gamma,d,\NNN)$ such that 
\begin{equation}
|\MMM_{k \NNN}^{-1} | \leq \mathsf{C}_{\MMM} k^{-1}
.
\label{eq:Mn:inverse:bound}
\end{equation}

On the other hand, 
from~\eqref{eq:recursion:1:b}, for each $k\geq 2$ we have that 
\begin{equation}
\mathsf{F}_1(k \NNN):= -  \sum_{m+j= k, m,j\geq 1}  \bar q_{m\NNN} \bar v_{j \NNN} \bigl(1+\alpha d + 2k \NNN - 2(1-\alpha) j \NNN  \bigr) 
.
\label{eq:recursion:1:b:new} 
\end{equation}
Using the fact that $\sum_{j=1}^{k-1} j^{-\frac 32} (k-j)^{-\frac 32} \leq 6 k^{-\frac 32}$, and the inductive bounds~\eqref{eq:power:series:bounds}, we deduce that there exists a constant $\mathsf{C}_{1} = \mathsf{C}_{1}(\gamma,d,\NNN)$ such that 
\begin{equation}
 |\mathsf{F}_1(k \NNN)| 
 \leq \mathsf{C}_{1} k \sum_{j=1}^{k-1} |\bar q_{(k-j)\NNN}| \, |\bar v_{j \NNN}| 
 \leq 6 \mathsf{C}_{1} k^{-\frac 12} (\mathsf{C}^\prime)^2  \mathsf{C}^k .
 \label{eq:recursion:1:b:bound}
\end{equation}
On the other hand, from \eqref{eq:recursion:2:b} we have
\begin{align}
\mathsf{F}_2(k \NNN)
&:= -   \sum_{m+j= k, m,j\geq 1}   \bar v_{m\NNN} \bar v_{j\NNN} \bigl( 1 + \cxbar +  3 \bar v_0  + 4 j \NNN (\cxbar + \bar v_0)   \bigr) 
\notag\\
&\quad -  \sum_{m+j=k, m,j\geq 1}   \bar v_{m\NNN} \bar q_{j\NNN} 2\alpha \bar q_0  (j \NNN+1)  
-  \sum_{m+j=k, m,j\geq 1}  \bar q_{m\NNN} \bar q_{j\NNN} (\cxbar + \bar v_0) \bigl( \tfrac{2\alpha^2}{\gamma} +  2 \alpha j \NNN   \bigr) 
\notag\\
&\quad - {\mathbf 1}_{k\geq 3} \sum_{m+\ell+j=k, m,\ell,j\geq 1} \bar v_{m\NNN}  \bigl( \bar v_{\ell\NNN} \bar v_{j\NNN} + \alpha \bar q_{\ell\NNN} \bar q_{j\NNN}\bigr) (2j\NNN+1).
\label{eq:recursion:2:b:new} 
\end{align}
As before, using that $\sum_{j=1}^{k-1} j^{-\frac 32} (k-j)^{-\frac 32} \leq 6 k^{-\frac 32}$, the bound $\sum_{j=1}^{k-2} \sum_{\ell = 1}^{k-j-1} \ell^{-\frac 32} j^{-\frac 32} (k-j-\ell)^{-\frac 32} \leq 36 k^{-\frac 32}$, and the inductive bounds~\eqref{eq:power:series:bounds}, we deduce the existence of a constant $\mathsf{C}_{2} = \mathsf{C}_{2}(\gamma,d,\NNN)$ such that 
\begin{align}
 |\mathsf{F}_2(k \NNN)| 
 &\leq \mathsf{C}_{2} k \sum_{j=1}^{k-1} |\bar v_{(k-j)\NNN}| \, |\bar v_{j \NNN}| + |\bar v_{(k-j)\NNN}| \, |\bar q_{j \NNN}| + |\bar q_{(k-j)\NNN}| \, |\bar q_{j \NNN}|
 \notag\\
 &\quad 
 + {\mathbf 1}_{k\geq 3} \mathsf{C}_{2} k \sum_{j=1}^{k-2} \sum_{\ell=1}^{k-j-1} 
 | \bar v_{(k-j-\ell)\NNN}|  \bigl( |\bar v_{\ell\NNN}| \, |\bar v_{j\NNN}| + |\bar q_{\ell\NNN}| \, |\bar q_{j\NNN}| \bigr) 
 \notag\\
 &\leq 18 \mathsf{C}_{2} k^{-\frac 12} (\mathsf{C}^\prime)^2  \mathsf{C}^k + 72  \mathsf{C}_{2} k^{-\frac 12} (\mathsf{C}^\prime)^3  \mathsf{C}^k .
  \label{eq:recursion:2:b:bound}
\end{align}
From~\eqref{eq:Mn:inverse:bound}, \eqref{eq:recursion:1:b:bound}, and~\eqref{eq:recursion:2:b:bound}, we deduce that if $\mathsf{C}^\prime$ is chosen small enough to ensure $6 \mathsf{C}_{1} \mathsf{C}_{\MMM} \mathsf{C}^{\prime} \leq 1$ and $18 \mathsf{C}_{2} \mathsf{C}_{\MMM} \mathsf{C}^{\prime} + 72  \mathsf{C}_{2} \mathsf{C}_{\MMM} (\mathsf{C}^{\prime})^2 \leq 1$, then the inductive bounds hold.
\end{proof}

From Proposition~\ref{prop:solve:recursion:R=0}, Proposition~\ref{prop:power:series:bounds}, and ansatz \eqref{eq:V:Q:R=0} we obtain the immediate consequence.
\begin{corollary}[\bf Convergence of power series for $\bar V$ and $\bar Q$ at $R=0$]
\label{cor:power:series:V:Q}
Let $\NNN\geq 1$, and let $\cxbar = \cxstar(d,\gamma,\NNN)$. Let $\mathsf{C}_{\eqref{eq:power:series:bounds}} \geq 1$ be the constant from Proposition~\ref{prop:power:series:bounds}. Define the sequences $\{\bar q_n\}_{n \geq 1}$ and $\{\bar v_n\}_{n\geq 1}$ as in Proposition~\ref{prop:solve:recursion:R=0}. Then, the series  
\begin{equation}
\bar V(R) = \bar v_0 + \sum_{k \geq 1}  \bar v_{k\NNN} R^{2k\NNN},
\qquad
\bar Q(R) =  \bar q_0 + \sum_{k \geq 1}  \bar q_{k\NNN}  R^{2k \NNN},
\label{eq:power:series:V:Q:R=0}
\end{equation}
converge uniformly and absolutely for all $0 \leq R < \mathsf{C}_{\eqref{eq:power:series:bounds}}^{-\frac{1}{2\NNN}}$. The resulting real-analytic functions $\bar V$ and $\bar Q$ solve the system~\eqref{eq:euler6} for all $R \in (0,\mathsf{C}_{\eqref{eq:power:series:bounds}}^{-\frac{1}{2\NNN}})$. 
\end{corollary}

\begin{remark}[\bf No free parameters]
We emphasize that in the construction of $\bar V$ and $\bar Q$ in Corollary~\ref{cor:power:series:V:Q}, there are no free parameters; once we have chosen $d\in\{1,2,3\}$, $1<\gamma\leq 2d+1$, and we have chosen the integer $\NNN \geq 1$, the solution is uniquely determined.
\end{remark}

\subsection{The global solutions \texorpdfstring{$\bar V$ and $\bar Q$}{bar V and bar Q}}
Fix $\NNN\geq 1$. Corollary~\ref{cor:power:series:V:Q} provides a real-analytic solution $(\bar V,\bar Q)$ of \eqref{eq:euler6}, whose power series converges  for all $0 \leq R <  \mathsf{C}^{-\frac{1}{2\NNN}}_{\eqref{eq:power:series:bounds}}$. Our next goal is to extend this solution curve to  all $R<\infty$; see~Proposition~\ref{prop:global:profile} below.

In order to achieve this, we note that~\eqref{eq:euler6} is a $2\times2$ system of \emph{autonomous ODEs} (with respect to the $R \p_R$ derivative), whose coefficients that are rational functions of the unknowns. To see this, we may simply diagonalizing~\eqref{eq:euler6}, which leads to
\begin{subequations}
\label{eq:euler7:all}
\begin{equation} \label{eq:euler7}
R \p_R  \bar V   
=  \frac{\Delta_{\bar V}[\bar V,\bar Q]}{\Delta[\bar V,\bar Q]}  , \qquad
R \p_R  \bar Q   
= \frac{\Delta_{\bar Q}[\bar V,\bar Q]}{\Delta[\bar V,\bar Q]} , 
\end{equation}
where we have denoted
\begin{align}
\Delta[\bar V,\bar Q]
&:= (\cxbar + \bar V) (\cxbar + \bar V - \alpha \bar Q) (\cxbar + \bar V + \alpha \bar Q) ,
\label{eq:Delta:def}
\\
\Delta_{\bar V}[\bar V,\bar Q]
&:= (\cxbar + \bar V) \bigl( \alpha \bar Q \mathsf{P}_1[\bar V, \bar Q]  -  \mathsf{P}_2[\bar V, \bar Q] \bigr)
,
\label{eq:Delta:V:def}
\\
\Delta_{\bar Q}[\bar V,\bar Q]
&:=\alpha \bar Q \mathsf{P}_2[\bar V, \bar Q]  - (\cxbar + \bar V)^2 \mathsf{P}_1[\bar V, \bar Q]  
,
\label{eq:Delta:Q:def}
\end{align}
and we recall that 
\begin{align}
\mathsf{P}_1[\bar V, \bar Q] 
&:= \bigl(1 + (1 + \alpha d)\bar V  \bigr) \bar Q   
\label{eq:P1:recall}
,
\\
\mathsf{P}_2[\bar V, \bar Q] 
&:= 
(\cxbar + \bar V) \bigl(  \bar V + \bar V^2 + \tfrac{2\alpha^2}{\gamma} \bar Q^2 \bigr)
+ \tfrac{\alpha}{\gamma}  \bar Q^2 (\bar V - \bar v_0)
.
\label{eq:P2:recall}
\end{align}
\end{subequations} 
Note that $\Delta_{\bar V}$ and $\Delta_{\bar Q}$ are fourth order polynomials in the variables $(\bar V,\bar Q)$, while $\Delta$ is a third order polynomial. We may expand the expressions in~\eqref{eq:Delta:V:def} and~\eqref{eq:Delta:Q:def} as
\begin{subequations}
\label{eq:Delta:alt}
\begin{align}
\Delta_{\bar V}[\bar V,\bar Q]
&:= \tfrac{\alpha^2(2+\gamma d)}{\gamma} (\cxbar + \bar V) \bar Q^2  \bigl(  \bar V - \bar v_0 \bigr)  
- (\cxbar + \bar V)^2 \bigl(  \bar V + \bar V^2 + \tfrac{2\alpha^2}{\gamma} \bar Q^2 \bigr)
,
\label{eq:Delta:V:alt}
\\
\Delta_{\bar Q}[\bar V,\bar Q]
&:=\alpha \bar Q  (\cxbar + \bar V) \bigl(  \bar V + \bar V^2 + \tfrac{2\alpha^2}{\gamma} \bar Q^2 \bigr)
+ \tfrac{\alpha^2}{\gamma}  \bar Q^3 (\bar V - \bar v_0)
- (1+\alpha d) (\cxbar + \bar V)^2 \bar Q \bigl(\bar V - \bar v_0\bigr) 
.
\label{eq:Delta:Q:alt}
\end{align} 
\end{subequations}

\begin{remark}[\bf The fixed points of interest]
\label{rem:fixed:points}
We remark that the points $(\bar v_0,\bar q_0)$ (corresponding to $R=0$) and $(0,0)$ (corresponding to $R\to \infty$) are fixed points of the autonomous ODE  system~\eqref{eq:euler7}. This is because
by~\eqref{eq:Delta:alt}, \eqref{eq:values:at:zero:and:cb}, and 
Lemma~\ref{lem:properties:exponents} we have that
\begin{subequations}
\begin{align}
&\Delta_{\bar V}[\bar v_0,\bar q_0] = \Delta_{\bar Q}[\bar v_0,\bar q_0] = 0 , \\
&\Delta[\bar v_0,\bar q_0] = (\cxbar +  \bar v_0) (\cxbar +  \bar v_0 +\alpha \bar q_0) (\cxbar +  \bar v_0  - \alpha \bar q_0) > 0  ,
\end{align}
and also
\begin{align}
&\Delta_{\bar V}[0,0] = \Delta_{\bar Q}[0,0] = 0 , \\
&\Delta[0,0] =  \cxbar ^3 >0 .
\end{align}
\end{subequations}
Additionally, one may verify that $(\bar v_0,\bar q_0)$ is a saddle and $(0,0)$ is a sink, though this information is not explicitly used in our analysis. The local-in-$R$ solution we have constructed in Corollary~\ref{cor:power:series:V:Q} ``escapes'' the saddle point $(\bar v_0,\bar q_0)$ along its unstable manifold, in the direction of increasing $\bar V$ (this is because $\bar v_{\NNN}>0$, cf.~\eqref{eq:first:nontrivial:Taylor:coefficient}) and decreasing $\bar Q$ (this is because $\bar q_{\NNN}<0$, cf.~\eqref{eq:bar:Q:zero} and~\eqref{eq:first:nontrivial:Taylor:coefficient}).
\end{remark}

\begin{figure}[htb!]
\centering
\begin{minipage}{.48\textwidth}
  \centering
  \includegraphics[width=0.8\linewidth]{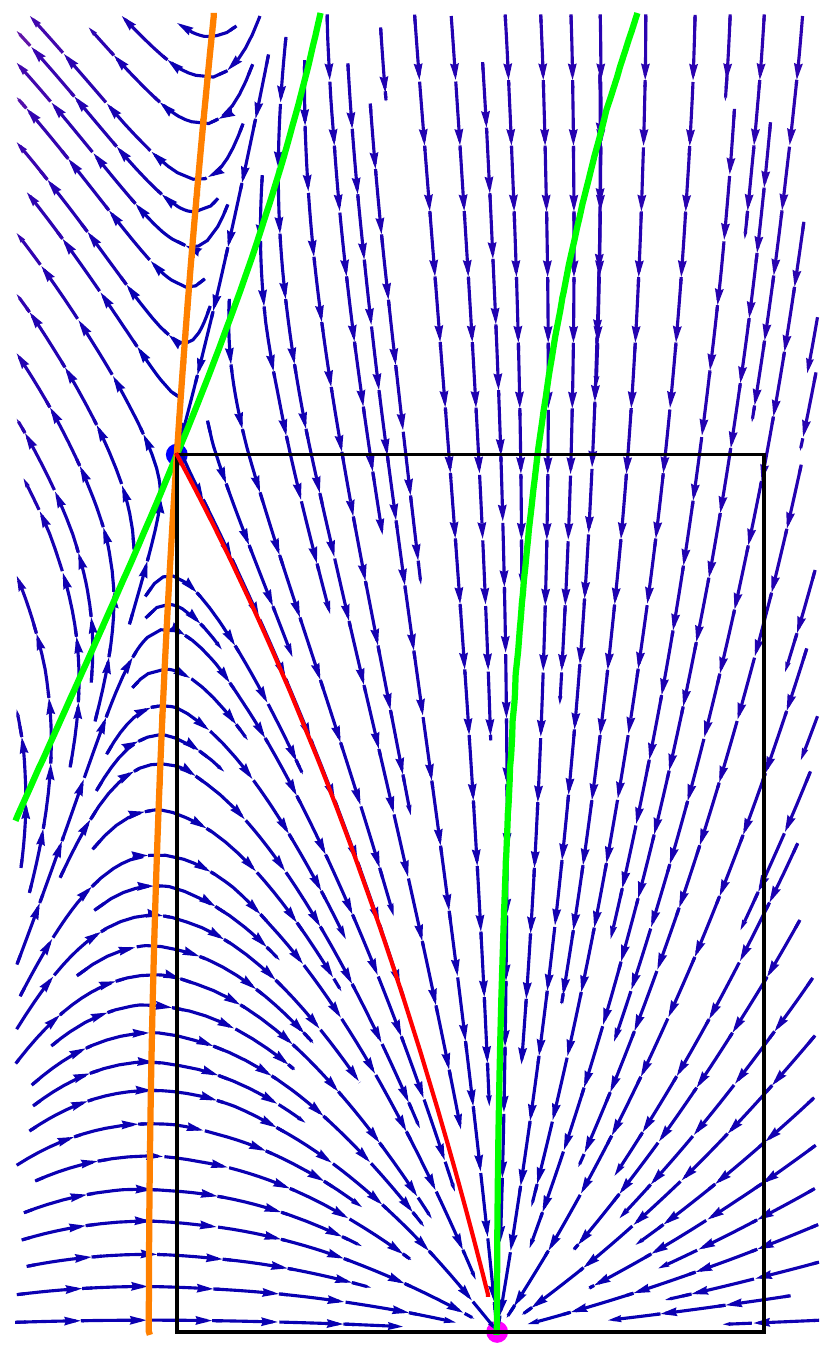}
\end{minipage}%
\begin{minipage}{.48\textwidth}
  \centering
  \includegraphics[width=0.63\linewidth]{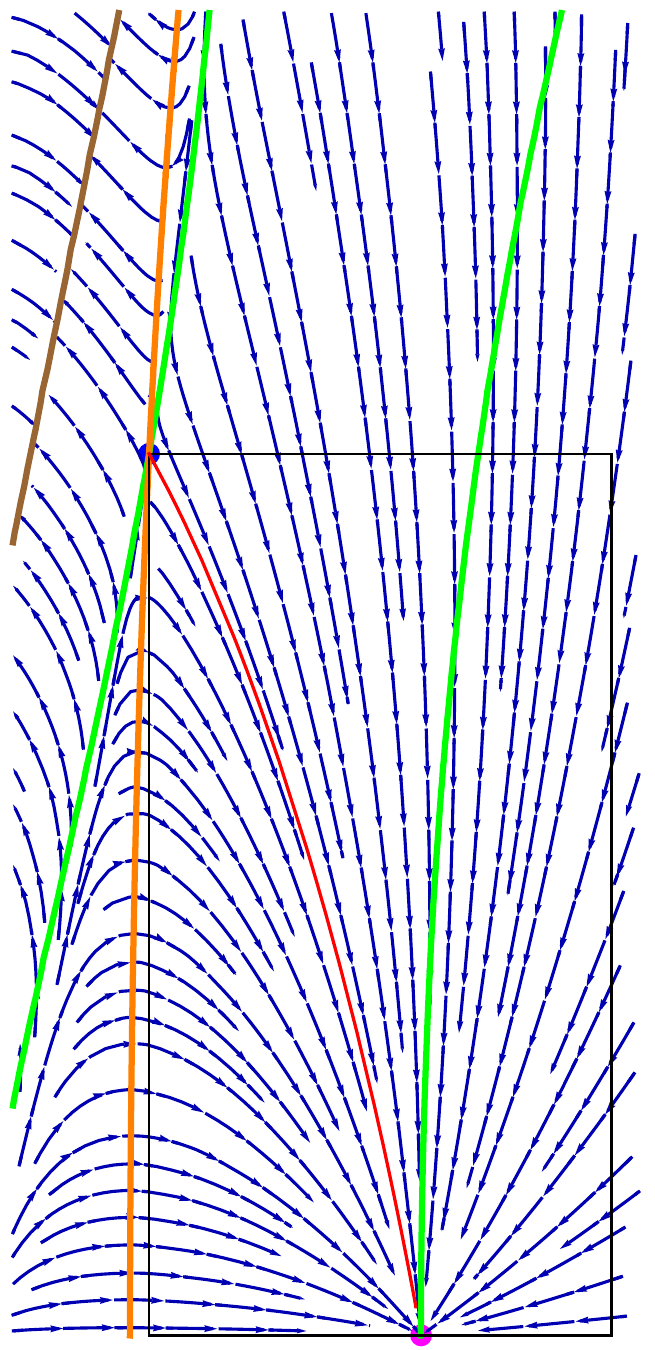}
\end{minipage}
\caption{The phase portrait of the $2\times2$ system of ODEs~\eqref{eq:euler7}, for the parameter choices $(d,\gamma,\NNN)=(3,\frac 53,1)$ for the left panel, and $(d,\gamma,\NNN)=(3,\frac 75,2)$ for the right panel. The vertical axis represents $\bar Q$, while the horizontal one is $\bar V$. The pink point has coordinates $(0,0)$. The blue point has coordinates $(\bar v_0,\bar q_0)$. The green curve represents one branch of $\{\Delta_{\bar V}[\bar V,\bar Q] = 0\}$, while the orange curve represents a branch of $\{\Delta_{\bar Q}[\bar V,\bar Q] = 0\}$. The brown curve (which only enters our plotting area for the right panel) represents a branch of $\{\Delta[\bar V,\bar Q] = 0\}$. In this phase portrait it is evident that if a trajectory $(\bar V(R),\bar Q(R))$ is able to ``escape'' the blue point in the direction $(+,-)$, which is precisely how we chose our eigenvector of $\MMM_{\NNN}$ in Proposition~\ref{prop:solve:recursion:R=0}, then this trajectory is confined to the rectangle $\Omega$ from Definition~\ref{def:Omega}, plotted here in black. By Poincar\'e-Bendixson, this trajectory is forced to approach (as $R\to \infty$) the stationary point $(0,0)$ which lies on the boundary of the confining region. The numerically computed trajectory $\{(\bar V(R),\bar Q(R)) \colon 0\leq R \leq \infty\}$ is represented by the thin red curve.}
\label{fig:phase:portrait:1}
\end{figure}

Figure~\ref{fig:phase:portrait:1} motivates us to define an \emph{invariant region} for the $2\times2$ autonomous ODE system~\eqref{eq:euler7}; from a technical viewpoint, it is  convenient to use a larger, ``rectangular'' invariant region, given by:

\begin{definition}[\bf Invariant region $\Omega$]
\label{def:Omega}
We let $\Omega = \bigl\{ (\bar V, \bar Q) \colon \bar Q \in (0,\bar q_0), \bar V \in (\bar v_0,- \tfrac{d\gamma}{2} \bar v_0) \bigr\}$.
\end{definition}

\begin{proposition}
\label{prop:Omega:invariant}
The closure of the domain $\Omega$ from Definition~\ref{def:Omega} is an invariant region for the ODE system~\eqref{eq:euler7}.
\end{proposition}
\begin{proof}[Proof of Proposition~\ref{prop:Omega:invariant}]
Let $\vec{n} = \vec{n}(\bar V,\bar Q)$ denote the outward unit normal to $\partial \Omega$; this vector is well-defined and smooth except for the corner points $(- \tfrac{d\gamma}{2} \bar v_0,0)$, $(\bar v_0,0)$, $(\bar v_0,\bar q_0)$, and $(- \tfrac{d\gamma}{2} \bar v_0,\bar q_0)$. The proof of the proposition amounts to showing that $(\Delta_{\bar V},\Delta_{\bar Q}) \cdot \vec{n} \leq 0$ for all $(\bar V,\bar Q) \in\partial \Omega$. We check this inequality individually on the four pieces of $\partial \Omega$:
\begin{itemize}[leftmargin=1em]
\item Let $\partial_{\rm bottom} \Omega := [\bar v_0,- \tfrac{d\gamma}{2} \bar v_0] \times \{0\}$, so that $\vec{n} = (0,-1)$. Since $\bar Q=0$, by~\eqref{eq:Delta:Q:alt} we have that $\Delta_{\bar Q}[\bar V,0] = 0$, and thus $(\Delta_{\bar V},\Delta_{\bar Q}) \cdot \vec{n} = 0$ on the ``bottom'' part of $\partial \Omega$. 

\item Let $\partial_{\rm left} \Omega: = \{\bar v_0\} \times [0,\bar q_0]$, so that $\vec{n} = (-1,0)$. By~\eqref{eq:Delta:V:alt}, and recalling~\eqref{eq:bar:V:zero}, \eqref{eq:bar:Q:zero}, 
we have that 
\begin{align*}
- \Delta_{\bar V}[\bar v_0,\bar Q] 
= (\cxbar + \bar v_0)^2 \bigl(  \bar v_0 + \bar v_0^2 + \tfrac{2\alpha^2}{\gamma} \bar Q^2 \bigr)
= \tfrac{2\alpha^2}{\gamma}  (\cxbar + \bar v_0)^2  \bigl( \bar Q^2 - \bar q_0^2 \bigr)
\leq 0 .
\end{align*}
Hence, $(\Delta_{\bar V},\Delta_{\bar Q}) \cdot \vec{n} \leq  0$ on the ``left'' part of $\partial \Omega$.

\item Let $\partial_{\rm top} \Omega: = [\bar v_0,- \tfrac{d\gamma}{2} \bar v_0] \times \{ \bar q_0 \}$, so that $\vec{n} = (0,1)$. By~\eqref{eq:Delta:Q:alt}, and recalling~\eqref{eq:bar:V:zero}, \eqref{eq:bar:Q:zero}, 
we have that 
\begin{align*}
\Delta_{\bar Q}[\bar V ,\bar q_0] 
&= \alpha   (\cxbar + \bar V) \bigl(  \bar V + \bar V^2 + \tfrac{2\alpha^2}{\gamma} \bar q_0^2 \bigr) \bar q_0
+ \bigl(\tfrac{\alpha^2}{\gamma}  \bar q_0^2 - (1+\alpha d) (\cxbar + \bar V)^2 \bigl) \bigl(\bar V - \bar v_0\bigr) \bar q_0
\\
&=-  \mathsf{Poly}[\bar V +\cxbar ] \bigl(\bar V - \bar v_0\bigr)\bar q_0 ,
\end{align*}
where
\begin{equation*}
\mathsf{Poly}[X]:= (1+\alpha (d-1)) X^2 + \alpha    \bigl( \cxbar- 1 - \bar v_0 \bigr) X
-  \tfrac{d \alpha }{2}  \bar v_0^2 
\end{equation*} 
is a quadratic polynomial in $X = \cxbar + \bar V \geq \cxbar + \bar v_0 >0$, whose leading order coefficient is positive.
Since $\bar q_0>0$ and $\bar V - \bar v_0 \geq 0$ for $\bar V \in[\bar v_0,-\frac{d\gamma}{2} \bar v_0]$, the sign of $\Delta_{\bar Q}[\bar V ,\bar q_0] $ is opposite to the sign of $\mathsf{Poly}[\bar V +\cxbar ]$. Thus, if we show that $\mathsf{Poly}[\bar V + \cxbar] \geq 0$ for all $\bar V \in[\bar v_0,-\frac{d\gamma}{2} \bar v_0]$, then $(\Delta_{\bar V},\Delta_{\bar Q}) \cdot \vec{n} = \Delta_{\bar Q}[\bar V,\bar q_0] \leq 0$ on the ``top'' part of $\partial \Omega$.

An explicit computation, which appeals to the second bullet in Lemma~\ref{lem:properties:exponents}, shows that for $X\geq \cxbar + \bar v_0$, we have
\begin{align*}
&\mathsf{Poly}[X] - \mathsf{Poly}[\cxbar + \bar v_0] 
\notag\\
&= (X - \cxbar - \bar v_0) \bigl( (1+\alpha (d-1)) X + (1+\alpha d) \cxbar + (1+ \alpha (d-2)) \bar v_0 - \alpha\bigr)
\notag\\
&\geq (X - \cxbar - \bar v_0) \tfrac{1}{1+\alpha d} \bigl( (2+2\alpha d-\alpha) (1 +   \sqrt{\tfrac{\alpha \gamma d}{2}}) - (2 + 2 \alpha d-3\alpha)  - \alpha (1+\alpha d) \bigr)
\notag\\
&\geq (X - \cxbar - \bar v_0) \tfrac{\alpha}{1+\alpha d} \bigl( (2+2\alpha d-\alpha)  \sqrt{d} + 1- \alpha d \bigr)
\geq 0,
\end{align*}
for $d\in \{1,2,3\}$ and $\alpha \in (0,d]$. Therefore, if we are able to show that $\mathsf{Poly}[\cxbar + \bar v_0]  \geq0$, then we automatically have $\mathsf{Poly}[X] \geq0$ for all $X\geq \cxbar + \bar v_0$. We are thus left to verify the positivity of the explicit expression
\begin{align*}
\mathsf{Poly}[\bar v_0 + \cxbar ]
&= (1+\alpha (d-1)) (\cxbar + \bar v_0)^2 + \alpha    \bigl( \cxbar- 1 - \bar v_0 \bigr) (\cxbar + \bar v_0)
-  \tfrac{d \alpha }{2}  \bar v_0^2
\\
&= 
(1+\alpha d)  \cxbar^2  
-  \tfrac{2 + \alpha (2d-1) + \alpha^2 d}{1+\alpha d} \cxbar 
+ \tfrac{1- \alpha + \frac{d \alpha }{2}  + \alpha^2 d}{(1+\alpha d)^2}  .
\end{align*}
We may now view the above expression as a quadratic polynomial in $\cxbar$, whose leading order coefficient is positive. Recall from Lemma~\ref{lem:properties:exponents} that $\cxbar>\cxstar(d,\gamma,\infty)= \frac{1}{1+\alpha d}(1+  \sqrt{\frac{\alpha \gamma d}{2}})$. Using the fact that $0<\alpha \leq d$, we may check that the above-displayed quadratic polynomial in $\cxbar$ is increasing for $\cxbar > \cxstar(d,\gamma,\infty)$; hence it attains its minimum at $\cxstar(d,\gamma,\infty)$. Therefore, for $\alpha \in (0,d]$ and $d\in\{1,2,3\}$, we conclude that 
\begin{align*}
\mathsf{Poly}[\bar v_0+\cxbar ]
&= \tfrac{1}{1+\alpha d} (1+  \sqrt{\tfrac{\alpha \gamma d}{2}})^2
-  \tfrac{2 + \alpha (2d-1) + \alpha^2 d}{(1+\alpha d)^2} (1+  \sqrt{\tfrac{\alpha \gamma d}{2}}) 
+ \tfrac{1- \alpha + \frac{d \alpha }{2}  + \alpha^2 d}{(1+\alpha d)^2} 
\notag\\
&=  
\tfrac{\alpha^2 d (2+  \gamma d)}{2(1+\alpha d)^2}   
+  \tfrac{\alpha - \alpha^2 d}{(1+\alpha d)^2} \sqrt{\tfrac{\alpha \gamma d}{2}}
\geq \tfrac{\alpha^2 d }{2(1+\alpha d)^2}   \Bigl( 2+  \gamma d 
+  \sqrt{\tfrac{2\gamma}{\alpha d}}
-   \sqrt{2\alpha \gamma d } \Bigr)
\geq 
\tfrac{\alpha^2 d }{ (1+\alpha d)^2}  
.
\end{align*}
The last inequality follows by considering separately the cases  $\alpha d \leq 1$ and  $1 < \alpha d \leq d^2$.
This concludes the proof of 
$(\Delta_{\bar V},\Delta_{\bar Q}) \cdot \vec{n}  \leq 0$ on the ``top'' part of $\partial \Omega$.

\item Let $\partial_{\rm right} \Omega: = \{- \tfrac{d\gamma}{2} \bar v_0\} \times [0,\bar q_0]$, so that $\vec{n} = (1,0)$. By~\eqref{eq:Delta:V:alt},  recalling~\eqref{eq:bar:V:zero}, \eqref{eq:bar:Q:zero}, and since Lemma~\ref{lem:properties:exponents} implies $\cxbar>\cxstar(d,\gamma,\infty)= (-\bar v_0)(1+  \sqrt{\frac{\alpha \gamma d}{2}})$, we obtain that 
\begin{align*}
\Delta_{\bar V}[- \tfrac{d\gamma}{2} \bar v_0,\bar Q] 
&=
\tfrac{\alpha^2(2+\gamma d)}{\gamma} (\cxbar - \tfrac{d\gamma}{2} \bar v_0) \bar Q^2   ( - \bar v_0  ) (\tfrac{d\gamma}{2}+1) 
- (\cxbar - \tfrac{d\gamma}{2} \bar v_0 )^2 \bigl( - \tfrac{d\gamma}{2} \bar v_0 + \tfrac{(d\gamma)^2}{4} \bar v_0^2 + \tfrac{2\alpha^2}{\gamma} \bar Q^2 \bigr)
\\
&\leq
 \Bigl(  \tfrac{d\alpha (2+\gamma d)}{2}(1 + \tfrac{d\gamma}{2}) - (1 + \tfrac{d\gamma}{2} + \sqrt{\tfrac{\alpha \gamma d}{2}}) (\tfrac{(d\gamma)^2}{4} + \alpha d) 
 \Bigr)   ( - \bar v_0  )^3 (\cxbar - \tfrac{d\gamma}{2} \bar v_0)    
 \\
&\leq
- \tfrac{d^2 (d+2)}{8}   ( - \bar v_0  )^3 (\cxbar - \tfrac{d\gamma}{2} \bar v_0)    
 <0 
 .
\end{align*}
In the second to last inequality we have used that the complicated expression appearing on the third line is maximized at $\gamma=1$, so that $\alpha = 0$.
Hence, $(\Delta_{\bar V},\Delta_{\bar Q}) \cdot \vec{n} \leq 0$ on the ``right'' part of $\partial \Omega$.
\end{itemize}
This concludes the proof of the Proposition.
\end{proof}

From Proposition~\ref{prop:Omega:invariant} we may immediately deduce:
\begin{corollary}
\label{cor:Omega:invariant}
For all points $(\bar V, \bar Q) \in \Omega$, we have the bounds  
\begin{subequations}
\begin{align}
\tfrac{1}{\bar Q} \Delta_{\bar Q}[\bar V, \bar Q] &\leq - \tfrac{\alpha^2 d }{ (1+\alpha d)^2}   (\bar V-\bar v_0) < 0 ,  
\label{eq:Omega:invariant:a}\\
\cxbar + \bar V + \alpha \bar Q > \cxbar + \bar V > \cxbar + \bar V - \alpha \bar Q &\geq \tfrac{1+ \frac 23 \alpha d}{4 \NNN (1+\alpha d)} .
\label{eq:Omega:invariant:c}
\end{align}
\end{subequations}
\end{corollary}
\begin{proof}[Proof of Corollary~\ref{cor:Omega:invariant}]
By definition, we have that $\bar Q>0$ in $\Omega$, which is an open set. Thus, the first two inequalities in~\eqref{eq:Omega:invariant:c} hold trivially. For the last inequality in~\eqref{eq:Omega:invariant:c}, due to the last bullet in Lemma~\ref{lem:properties:exponents} we have that for all $(\bar V,\bar Q) \in \Omega$ it holds that $\cxbar +\bar V - \alpha \bar Q > \cxbar + \bar v_0 -\alpha \bar q_0  \geq \frac{1+\frac 23 \alpha d}{4 \NNN (1+\alpha d)}$. This proves~\eqref{eq:Omega:invariant:c}. 

In order to prove~\eqref{eq:Omega:invariant:a}, we use~\eqref{eq:Delta:Q:alt} to rewrite
\begin{equation*}
\tfrac{1}{\bar Q}\Delta_{\bar Q}[\bar V,\bar Q]
=  \alpha   (\cxbar + \bar V) \bigl(  \bar V + \bar V^2 + \tfrac{2\alpha^2}{\gamma} \bar Q^2 \bigr)
+ \tfrac{\alpha^2}{\gamma}  \bar Q^2  (\bar V - \bar v_0)
- (1+\alpha d) (\cxbar + \bar V)^2   \bigl(\bar V - \bar v_0\bigr) 
.
\end{equation*}
This is a quadratic polynomial in $\bar Q$, which contains no linear term in $\bar Q$ and the coefficient of $\bar Q^2$ is strictly positive. Thus, the above expression is maximized when $\bar Q = \bar q_0$, and so 
\begin{equation*}
\tfrac{1}{\bar Q} \Delta_{\bar Q}[\bar V,\bar Q] \leq \tfrac{1}{\bar q_0}\Delta_{\bar Q}[\bar V,\bar q_0] .
\end{equation*}
Recall however that in the proof of Proposition~\ref{prop:Omega:invariant} (the third bullet in the proof) we have established $\frac{1}{\bar q_0} \Delta_{\bar Q}[\bar V,\bar q_0] = -  \mathsf{Poly}[\bar V +\cxbar ] \bigl(\bar V - \bar v_0\bigr) \leq - \tfrac{\alpha^2 d }{ (1+\alpha d)^2}  (\bar V - \bar v_0)$. This concludes the proof of~\eqref{eq:Omega:invariant:a}.
\end{proof}

We now are able to conclude the existence of a global solution $(\bar V,\bar Q)(R)$ of~\eqref{eq:euler7}.
\begin{proposition}[\bf The globally defined profile]
\label{prop:global:profile}
Fix $\NNN\geq 1$, and let $(\bar V,\bar Q)(R)$ be the solution of~\eqref{eq:euler7} constructed in Corollary~\ref{cor:power:series:V:Q}, which is real-analytic in $R$ for $0 \leq R <  \mathsf{C}^{-\frac{1}{2\NNN}}_{\eqref{eq:power:series:bounds}}$. This local-in-$R$ solution may be uniquely extended as a smooth (real-analytic) global-in-$R$ solution of~\eqref{eq:euler7}, with $(\bar V,\bar Q)(R) \to (0,0)$ as $R\to \infty$.
\end{proposition}
\begin{proof}[Proof of Proposition~\ref{prop:global:profile}]
As noted at the end of Remark~\ref{rem:fixed:points}, the local-in-$R$ solution constructed in Corollary~\ref{cor:power:series:V:Q} enters the interior of the invariant domain $\Omega$. Due to~\eqref{eq:Omega:invariant:c}, and recalling~\eqref{eq:Delta:def} and \eqref{eq:Delta:alt}, we have that the right side of the $2\times 2$ system of autonomous ODEs~\eqref{eq:euler7}, namely $\frac{\Delta_{\bar V}}{\Delta}[\bar V,\bar Q]$ and $\frac{\Delta_{\bar Q}}{\Delta}[\bar V,\bar Q]$, are analytic functions of $\bar V$ and $\bar Q$ in all of $\Omega$. Thus, the local-in-$R$ solution can be extended to a global-in-$R$ solution. Recall that $(\bar v_0,\bar q_0)$ and $(0,0)$ are the only fixed points of \eqref{eq:euler7} within the closure of the rectangle $\Omega$, and these lie on $\partial \Omega$. Moreover, since $\bar Q(R)$ is strictly decreasing in $R$ by ~\eqref{eq:Omega:invariant:a},
\eqref{eq:Omega:invariant:c}, and \eqref{eq:euler7:all}, the trajectory $(\bar V,\bar Q)(R)$ cannot return to $(\bar v_0,\bar q_0)$ as $R$ increases and  cannot produce a periodic orbit. By the Poincar\'e-Bendixson theorem, we obtain that as $R\to \infty$ the global-in-$R$ solution approaches the stable equilibrium $(0,0)$. Along this trajectory, $\bar Q(R)$ is monotone decreasing. 
\end{proof}

\subsection{Power series for \texorpdfstring{$\bar V$ and $\bar Q$}{bar V and bar Q} near \texorpdfstring{$R=\infty$}{R = infinity}}
Proposition~\ref{prop:global:profile} guarantees that the global-in-$R$ solution of~\eqref{eq:euler7} is real-analytic in $R$, which is a \emph{qualitative} statement. On the other hand, Proposition~\ref{prop:power:series:bounds} and Corollary~\ref{cor:power:series:V:Q} provide a~\emph{quantitative} description of the solution for $R\ll 1$. The goal of this section is to obtain a similar \emph{quantitative} description of the solution for $R\gg 1$. We have that:
\begin{proposition}[\bf Convergence of power series for $\bar V$ and $\bar Q$ near $R=\infty$]
\label{prop:power:series:infinity}
There exist real coefficients $\{\tilde v_n\}_{n\geq 1}$ and $\{\tilde q_n\}_{n\geq 1}$, defined recursively in~\eqref{eq:recursions:at:infinity} below, such that the solution $(\bar V,\bar Q)(R)$ from Proposition~\ref{prop:global:profile} may be described as the convergent power series expansion
\begin{equation}
\bar V(R) =  \sum_{n\geq 1} \underline v_n R^{- \frac{n}{\cxbar}},
\qquad
\bar Q(R) =  \sum_{n\geq 1}  \underline q_n  R^{- \frac{n}{\cxbar}},
\label{eq:V:Q:R=infty}
\end{equation}
where $\cxbar = \cxstar(d,\gamma,\NNN)$,  $ |\underline v_1| < \infty$, and $0 < \underline q_1 < \infty$. Moreover, there exists a constant $\mathsf{C}_{\eqref{eq:power:series:bounds:infty}} = \mathsf{C}_{\eqref{eq:power:series:bounds:infty}}(d,\gamma,\NNN)>1$ such that the power series in~\eqref{eq:V:Q:R=infty} converge uniformly and absolutely for $R > \mathsf{C}_{\eqref{eq:power:series:bounds:infty}}^{\cxbar} (|\underline v_1|^2 + |\underline q_1|^2)^{\frac{\cxbar}{2}}$.
\end{proposition}

\begin{remark}[\bf Velocity of the ground state profile is negative]
For the ground state at $\NNN=1$, the leading order coefficient $\underline{v}_1$ appearing in~\eqref{eq:V:Q:R=infty} is \emph{strictly negative} and we have that $\bar V(R) < 0$ for all $R \geq 0$. 
 This negativity of $\underline{v}_1$ and of $\bar V$ is not used anywhere in the bulk of the paper, so we have relegated the proof of this fact to  Appendix~\ref{sec:appendix:sign}. 
\end{remark}

\begin{proof}[Proof of Proposition~\ref{prop:power:series:infinity}]
Due to Proposition~\ref{prop:global:profile} we know that $(\bar V,\bar Q)(R) \to 0$ as $R\to \infty$, and that $\bar Q(R)>0$ for all $R>0$.  
From~\eqref{eq:Delta:def} and~\eqref{eq:Delta:alt}, a Taylor series computations shows that as $(\bar V,\bar Q)\to (0,0)$ we have
\begin{subequations}
\label{eq:detailed:R:dR:V:Q}
\begin{align}
R\p_R \bar V 
&= - \tfrac{1}{\cxbar} \bar V 
+ \tfrac{1-\cxbar}{\cxbar^2} \bar V^2 
- \tfrac{\alpha^2}{\gamma \cxbar^2} ( 2 \cxbar -\tfrac{2+\gamma d}{1+\alpha d} ) \bar Q^2 + \OO(|\bar V|^3 + \bar Q^3) ,
\\
R \p_R \bar Q 
&= - \tfrac{1}{\cxbar} \bar Q   
+ \tfrac{1+\alpha - \cxbar (1+\alpha d)}{\cxbar^2} \bar V \bar Q + \OO(|\bar V|^2 \bar Q + \bar Q^3).
\end{align}
\end{subequations}
From~\eqref{eq:detailed:R:dR:V:Q} and the qualitative statement $(\bar V,\bar Q)\to (0,0)$ as $R\to \infty$ we deduce that there exists a constant $C_0 = C_0(\alpha, d,\NNN)>0$ and $R_0 \geq 1$ sufficiently large, such that for all $R\geq R_0$ we have
\[
|R \p_R (R^{\frac{1}{\cxbar}} \bar V(R))| + |R \p_R (R^{\frac{1}{\cxbar}} \bar Q(R))|
\leq C_0 R^{-\frac{1}{\cxbar}} \bigl( R^{\frac{1}{\cxbar}} |\bar V(R)| +  R^{\frac{1}{\cxbar}} \bar Q(R)\bigr)^2
,
\]
and 
\[
|\bar V(R)| + \bar Q(R) \leq \tfrac{1}{8 C_0 \cxbar}
.
\]
From the two bounds above and the fundamental theorem of calculus we obtain that 
\[
R^{\frac{1}{\cxbar}} |\bar V(R)| +  R^{\frac{1}{\cxbar}} \bar Q(R)
\leq
2 R_0^{\frac{1}{\cxbar}} |\bar V(R_0)| + 2  R_0^{\frac{1}{\cxbar}} \bar Q(R_0) =: C_1,
\quad \mbox{for all} \quad
R \geq R_0.
\]
Therefore, we may  define 
\[
\underline v_1 := \lim_{R\to \infty} R^{\frac{1}{\cxbar}} \bar V(R) 
= R_0^{\frac{1}{\cxbar}} \bar V(R_0) + \int_{R_0}^{\infty}  
\underbrace{\p_{R^\prime} \bigl( (R^\prime)^{\frac{1}{\cxbar}} \bar V(R^\prime) \bigr)}_{|\cdot|\leq C_0 C_1^2 (R^\prime)^{-\frac{1}{\cxbar}-1} } dR^\prime,
\] 
and 
\[
\underline q_1 := \lim_{R\to \infty} R^{\frac{1}{\cxbar}} \bar Q(R)
= R_0^{\frac{1}{\cxbar}} \bar Q(R_0) + \int_{R_0}^{\infty}  \underbrace{  \p_{R^\prime} \bigl( (R^\prime)^{\frac{1}{\cxbar}} \bar Q(R^\prime) \bigr)}_{|\cdot|\leq C_0 C_1^2 (R^\prime)^{-\frac{1}{\cxbar}-1} }  dR^\prime,
\]
and  deduce that $|\underline v_1| < \infty $ and $0 \leq \underline q_1 < \infty$. Moreover, the second identity in~\eqref{eq:detailed:R:dR:V:Q} shows that there exists $C_2>0$ such that 
\[
R^{\frac{1}{\cxbar}} \bar Q(R) 
\geq 
R_0^{\frac{1}{\cxbar}} \bar Q(R_0) e^{- C_2 \int_{R_0}^R (R^\prime)^{-\frac{1}{\cxbar}-1} dR^\prime}
\geq 
R_0^{\frac{1}{\cxbar}} \bar Q(R_0) e^{- \cxbar C_2}
\]
for all $R\geq R_0$. Therefore, we have proven that $\underline q_1 >0$.

In order to recursively compute the coefficients~$\{\underline v_n\}_{n\geq 2}$ and $\{\underline q_n\}_{n\geq 2}$, formulation~\eqref{eq:euler6} of the system of $2\times2$ autonomous ODEs~\eqref{eq:euler7} is more convenient to use. Inserting the ansatz~\eqref{eq:V:Q:R=infty} into \eqref{eq:euler6} leads to the  recursion relations
\begin{subequations}
\label{eq:recursions:at:infinity}
\begin{align}
\cxbar (n-1) \underline q_n 
&=  \sum_{m+j=n, m,j \geq 1}    \bigl(\cxbar (1+\alpha d) -  n + (1- \alpha) j  \bigr) \underline q_{m}    \underline v_j  \label{eq:recursions:at:infinity:a}\\
\cxbar (n-1) \underline v_n
&=   \sum_{m+j=n, m,j \geq 1} \bigl( \tfrac{2\alpha^2 \cxbar}{\gamma} + \tfrac{\alpha}{\gamma (1+\alpha d)} - \alpha m \bigr) \underline q_{m} \underline q_j 
+   \bigl( \cxbar+1 - 2 m \bigr) \underline v_{m} \underline v_j 
\notag\\
&\qquad 
+ {\bf 1}_{n\geq 3} \sum_{m+\ell+j = n, m, \ell, j \geq 1} \bigl( 1 - \tfrac{j}{\cxbar} \bigr) \underline v_{m} \bigl( \underline v_\ell \underline v_{j} + \alpha \underline q_\ell \underline q_{j} \bigr),
\label{eq:recursions:at:infinity:b}
\end{align} 
\end{subequations}
for all $n\geq 2$. Recursion~\eqref{eq:recursions:at:infinity} uniquely defines the coefficients $\{\underline v_n\}_{n\geq 2}$ and $\{\underline q_n\}_{n\geq 2}$ from knowledge of $\underline v_1$ and $\underline q_1$.
Moreover, as in the proof of Proposition~\ref{prop:power:series:bounds}, we may prove the existence of a sufficiently large constant $\mathsf{C}_{\eqref{eq:power:series:bounds:infty}} = \mathsf{C}_{\eqref{eq:power:series:bounds:infty}}(d,\gamma,\NNN) \geq 1$ and a sufficiently small constant $0 < \mathsf{C}^\prime_{\eqref{eq:power:series:bounds:infty}} = \mathsf{C}^\prime_{\eqref{eq:power:series:bounds:infty}}(d,\gamma,\NNN)<1$ such that 
\begin{subequations}
\label{eq:power:series:bounds:infty}
\begin{align}
|\underline v_n| &\leq  \mathsf{C}^\prime_{\eqref{eq:power:series:bounds:infty}} (\mathsf{C}_{\eqref{eq:power:series:bounds:infty}})^n (|\underline v_1|^2 + |\underline q_1|^2)^{\frac n2} n^{-\frac 32} \\
|\underline q_n| &\leq  \mathsf{C}^\prime_{\eqref{eq:power:series:bounds:infty}} (\mathsf{C}_{\eqref{eq:power:series:bounds:infty}})^n (|\underline v_1|^2 + |\underline q_1|^2)^{\frac n2}  n^{-\frac 32}
\end{align} 
\end{subequations}
for all $n \geq 1$. The proof of the bounds in~\eqref{eq:power:series:bounds:infty} is identical to the proof of~\eqref{eq:power:series:bounds}, and we  omit it for the sake of brevity. In particular, \eqref{eq:power:series:bounds:infty} implies that the power series expansions in~\eqref{eq:V:Q:R=infty} converge uniformly and absolutely for all $R > \mathsf{C}_{\eqref{eq:power:series:bounds:infty}}^{\cxbar} (|\underline v_1|^2 + |\underline q_1|^2)^{\frac{\cxbar}{2}} $.
\end{proof}

\subsection{The globally defined profile \texorpdfstring{$\bar H$}{bar H}}
The globally defined profile $(\bar V,\bar Q)$ from Proposition~\ref{prop:global:profile} allows us to directly compute the profile $\bar H$ which solves \eqref{eq:s5}, because the factor $\cxbar + \bar V(R)$ does not vanish (see~\eqref{eq:Omega:invariant:c}). Using~\eqref{eq:s5} and~\eqref{eq:cb:def}, we have
\begin{equation*}
\tfrac{R \p_R \bar H(R)}{\bar H(R)} = - \tfrac{\bar V(R) - \bar v_0}{\cxbar + \bar V(R)} < 0 ,
\end{equation*}
and thus, since $\bar H(0) =1$, we have 
\begin{equation}
\label{eq:bar:H:def}
\bar H(R) :=  \exp\Bigl(- \int_0^R \tfrac{\bar V(R^\prime) - \bar v_0}{R^\prime (\cxbar + \bar V(R^\prime))} d R^\prime \Bigr)  .
\end{equation}
We note that $\bar H(R) >0$ for all $R>0$. Also, note that from Corollary~\ref{cor:power:series:V:Q}, \eqref{eq:first:nontrivial:Taylor:coefficient},  and~\eqref{eq:bar:H:def} we have that 
\begin{equation}
 \bar H(R) =  \exp\Bigl( - \tfrac{1}{2\NNN(\cxbar + \bar v_0)} R^{2\NNN} + \OO(R^{4\NNN}) \Bigr) ,
 \qquad \mbox{as} \qquad R\to 0^+ .
 \label{eq:bar:H:R=0}
\end{equation}
Similarly, using Proposition~\ref{prop:power:series:infinity} we may deduce that for $R\geq 1$ we have
\begin{equation}
\bar H(R) = \bar H(1) R^{- \frac{\cxbar - \cbbar }{\cxbar}}  \exp\Bigl( - \int_1^\infty \tfrac{(\cxbar + \bar v_0) \bar V(R^\prime)}{\cxbar   R^\prime (\cxbar + \bar V(R^\prime))} dR^\prime + \OO(R^{-\frac{1}{\cxbar}})\Bigr),
\qquad \mbox{as} \qquad R\to \infty.
\label{eq:bar:H:R=infty}
\end{equation}

\subsection{Main result: existence of globally self-similar profiles}
We summarize the analysis in this section by stating the precise existence theorem we have established:
\begin{theorem}[\bf Exact self-similar profiles for non-isentropic Euler]
\label{thm:main:profiles}
Let $d \in \{1,2,3\}$, $1 < \gamma \leq 2d+1$, and let $\alpha = \frac{\gamma-1}{2}$. For each $\NNN \in \Naturals$, define the similarity exponents $\cxbar := \cxstar(d,\gamma,\NNN)$ and $\cbbar := \cbstar(d,\gamma,\NNN)$, via~\eqref{eq:cx:admissible} and~\eqref{eq:cb:admissible}, respectively. Then, for each such fixed $d,\gamma,\NNN$ there exist unique \emph{real-analytic} solutions $\bar V, \bar Q, \bar H \colon [0,\infty) \to \Reals$, of the system~\eqref{eq:euler5}, with the following properties:
\begin{enumerate}[label=(\roman*), leftmargin=*]
\item \textsl{(ODE system).} The profiles $(\bar V, \bar Q)$ solve the $2\times2$ system of autonomous ODEs~\eqref{eq:euler6}, or equivalently~\eqref{eq:euler7:all}. The profile $\bar H$ is determined by~\eqref{eq:bar:H:def}.

\item \textsl{(Behavior at $R=0$).} The profiles $(\bar V, \bar Q, \bar H)$ have explicit values at $R=0$ given by~\eqref{eq:bar:V:zero}, \eqref{eq:bar:Q:zero}, and~\eqref{eq:bar:H:zero}. Near $R=0$, the profiles admit convergent power series expansions (see~\eqref{eq:power:series:V:Q:R=0} and~\eqref{eq:bar:H:R=0}).

\item \textsl{(Behavior as $R\to\infty$).} The profiles vanish asymptotically as $R\to\infty$ and have power law behavior given by  $(\bar V, \bar Q, \bar H)(R) \sim (R^{-\frac{1}{\cxbar}},R^{-\frac{1}{\cxbar}},R^{-\frac{\cxbar-\cbbar}{\cxbar}})$. The profiles $\bar V$ and $\bar Q$ admit a convergent power series expansion as $R \to \infty$ (see~\eqref{eq:V:Q:R=infty}), while the behavior of $\bar H$ is given by~\eqref{eq:bar:H:R=infty}.

\item \textsl{(Global properties).} The three self-similar wave speeds present in the system obey a global-in-$R$ lower bound~\eqref{eq:Omega:invariant:c}, which we refer to as a \emph{global outgoing property}. Moreover, $\bar Q(R)$ is strictly monotone decreasing in $R$, $\bar Q(R) > 0$ for all $R\geq 0$, and we have $\bar H(R) > 0$ for all $R \geq 0$.
\end{enumerate}
Upon defining $\bar U(R) := R \bar V(R)$, $\bar \Sigma(R) := R \bar Q(R)$, and $\bar B(R) := R \bar H(R)$ (cf.~\eqref{eq:profile:normalize}), the profiles $(\bar \Sigma, \bar U, \bar B)$ are smooth on $[0,\infty)$, real-analytic on $(0,\infty)$, and solve the self-similar Euler system~\eqref{eq:euler4}.
\end{theorem}
\begin{proof}[Proof of Theorem~\ref{thm:main:profiles}]
The existence and uniqueness of the similarity exponents $\cxbar$ and $\cbbar$ is established in Definition~\ref{def:exponents}, with properties given in Lemma~\ref{lem:properties:exponents}. 
The local existence of the real-analytic profiles $(\bar V, \bar Q)$ near $R=0$ is established in Corollary~\ref{cor:power:series:V:Q}, with the power series coefficients constructed in Proposition~\ref{prop:solve:recursion:R=0} and the quantitative bounds given in Proposition~\ref{prop:power:series:bounds}. The values at $R=0$ are determined by~\eqref{eq:values:at:zero:and:cb}.
The global extension of the solution to all $R \in [0,\infty)$ is established in Proposition~\ref{prop:global:profile}, using the invariance of the region $\Omega$ (Proposition~\ref{prop:Omega:invariant}) and the Poincar\'e--Bendixson theorem.  Corollary~\ref{cor:Omega:invariant} gives the global outgoing property. The monotonicity of $\bar Q$ follows from~\eqref{eq:Omega:invariant:a} and the fact that $\Delta > 0$ in $\Omega$ (cf.~\eqref{eq:Omega:invariant:c}).
The behavior of $(\bar V, \bar Q)$ as $R \to \infty$, including the convergent power series expansion, is established in Proposition~\ref{prop:power:series:infinity}.
The profile $\bar H$ is constructed via~\eqref{eq:bar:H:def}, with asymptotic expansions near $R=0$ and $R=\infty$ given by~\eqref{eq:bar:H:R=0} and~\eqref{eq:bar:H:R=infty}, respectively. The positivity $\bar H(R) > 0$ follows directly from~\eqref{eq:bar:H:def}.
Finally, the equivalence between the $(\bar V, \bar Q, \bar H)$ system and the original self-similar Euler system~\eqref{eq:euler4} for $(\bar \Sigma, \bar U, \bar B)$ is established  via the change of variables~\eqref{eq:profile:normalize}.
\end{proof}

%%%%%%%%%%%%%%%%%%%%%%

\section{Stability analysis in radial symmetry}
\label{sec:stability}

For any given $0<\alpha \leq d$, in the previous section we have established the existence of a family (indexed by $\NNN \in \Naturals$) of smooth implosion profiles (solutions of~\eqref{eq:euler4})
\[
\bar U(R) = R \bar V(R),
\qquad
\bar \Sigma(R) = R \bar Q(R),
\qquad
\bar B(R) 
= R  \bar H(R)
,
\]
for suitable values of $\cxbar, \cubar$, and $\cbbar$, which depend on $d$, $\alpha$, $\NNN$.
Via~\eqref{eq:globally:SS:exact}, these profiles define globally self-similar solutions $(\bar u^r, \bar \sigma, \bar b)$ of the Euler equations in radial symmetry~\eqref{eq:euler2}. In this section we establish the stability of these solutions $(\bar u^r, \bar \sigma, \bar b)$ within the class of radially symmetric solutions of~\eqref{eq:euler2}, for which the initial density, velocity, and pressure are smooth, and for which the initial pressure vanishes at $r=0$. The main result of this section is Theorem~\ref{thm:stability:radial} below.

\subsection{Modulated self-similar ansatz}
\label{sec:true:ss:ansatz}
The stability analysis is performed using carefully designed modulated self-similar coordinates, which we define next.
We let 
\[
\tau \in [0,\infty)
\]
denote the self-similar time variable. The fundamental modulation functions will be denoted by
\begin{equation}
\cx, \cu,\cb \colon [0,\infty) \to \Reals.
\label{eq:cr:cu:cb:first}
\end{equation}
These modulation functions account for invariances of the system~\eqref{eq:euler2} as follows: $\cx(\tau)$ -- rescaling of the space variable, $\cu(\tau)$ -- rescaling of the velocity, while $\cb(\tau)$ -- the rescaling of the pseudo entropy. The functions $(\cx, \cu,\cb)(\tau)$ will be chosen carefully in our analysis (cf.~\eqref{eq:modulation:final} below); for the moment, we only emphasize that $\cx(\tau) \to \cxbar$, $\cu(\tau) \to \cubar = \cxbar -1$, and $\cb(\tau) \to \cbbar$ as $\tau \to \infty$.

Assuming for the moment that the modulation functions in~\eqref{eq:cr:cu:cb:first} have been chosen, we define
\begin{equation}
\label{eq:Cr:Cu:Cb:def} 
C_r(\tau):= \exp\Bigl( - \! \int_0^\tau\! \cx(\tau') d\tau'\Bigr), 
\; 
C_u(\tau):= \exp\Bigl( - \! \int_0^\tau \! \cu(\tau') d\tau'\Bigr), 
\; 
C_b(\tau):= \exp\Bigl( - \! \int_0^\tau \!\cb(\tau') d\tau'\Bigr).
\end{equation}
The relation between the original time variable $t$ and the self-similar time variable $\tau$ is not given by the simple  $\tau = - \log(-t)$; we need to account for \emph{time modulation}. Since we aim to keep the initial time as $t=-1$ (corresponding to $\tau=0$), we define 
\begin{subequations}
\label{eq:self-similar:ansatz}
\begin{equation}
 t = t(\tau):= -1 + \int_0^\tau C_r(\tau') C_u^{-1}(\tau') d\tau' 
 .
 \label{eq:self-similar:ansatz:tau}
\end{equation}
The self-similar radial coordinate is then defined as 
\begin{equation}
R:=  r \; C_r^{-1}(\tau)
.
\label{eq:self-similar:ansatz:R}
\end{equation}
With this notation, the self-similar velocity $U$, the re-scaled sound speed $\Sigma$, and square-root of the pseudo-entropy are defined as 
\begin{align}
u^r(r,t) &:= C_u(\tau) U (R,\tau), 
\label{eq:self-similar:ansatz:U} \\
\sigma(r,t) &:= C_u(\tau) \Sigma (R,\tau), 
\label{eq:self-similar:ansatz:Sigma} \\
b(r,t) &:=C_b(\tau) B(R,\tau) .
\label{eq:self-similar:ansatz:B}
\end{align}
\end{subequations}

\begin{remark}[\bf Relation to the globally self-similar transformation and variables]
We detail the usage of the $\approx$ symbol in Section~\ref{sec:self:similar:intro} of the Introduction, by relating the modulated self-similar transformation defined in~\eqref{eq:self-similar:ansatz} above, and the globally self-similar transformation used in~\eqref{eq:globally:SS:exact}. Assuming that we replace the modulation functions from~\eqref{eq:cr:cu:cb:first} with their limiting values as $\tau\to \infty$, that is, we replace $(\cx,\cu,\cb) \mapsto (\cxbar,\cubar,\cbbar)$, then \eqref{eq:Cr:Cu:Cb:def} becomes $(C_r,C_u,C_b) \mapsto (e^{-\cxbar \tau}, e^{-\cubar \tau}, e^{-\cbbar \tau})$. Importantly, since $\cubar = \cxbar - 1$, we have that \eqref{eq:self-similar:ansatz} becomes $t \mapsto -1 + \int_0^\tau e^{-\cxbar \tau'} e^{(\cxbar -1) \tau'} d\tau' = -1 + \int_0^\tau e^{-\tau'} d\tau' = - e^{-\tau}$. Equivalently, $-\log(-t)\mapsto \tau$, and so $(C_r,C_u,C_b) \mapsto ((-t)^{\cxbar }, (-t)^{\cubar}, (-t)^{\cbbar})$. Then, \eqref{eq:self-similar:ansatz:R} becomes $R\mapsto r  e^{\cxbar \tau} = r /(-t)^{\cxbar}$. Lastly, upon replacing $(U,\Sigma,B) \mapsto (\bar U,\bar \Sigma,\bar B)$, the definitions~\eqref{eq:self-similar:ansatz:U}--\eqref{eq:self-similar:ansatz:B} become \eqref{eq:globally:SS:exact}.
\end{remark}

\begin{remark}[\bf The blowup time]
While the initial time for the both the globally self-similar and the modulated self-similar analysis is $t=-1$ (equivalently $\tau=0$), the time of blowup may differ. Indeed, in the globally self-similar coordinates $\tau \to \infty$ corresponds to $t\to 0$, while the in the modulated self-similar analysis the blowup time $\tau \to \infty$ corresponds to 
\[
t \to t_*:= -1 + \int_0^\infty C_r(\tau') C_u^{-1}(\tau') d\tau'.
\]
Note that if we were to replace $(C_r,C_u)$ with $(e^{-\cxbar \tau}, e^{-\cubar \tau})$, then the above definition directly gives $t_* = 0$.
\end{remark}

\subsection{Self-similar evolution}
After a short computation using the relations $\frac{d\tau}{dt} = C_u C_r^{-1}$ (which follows from~\eqref{eq:self-similar:ansatz:tau}), $\frac{dR}{dr} = C_r^{-1}$ and $\frac{dR}{dt} =  \cx R C_u C_r^{-1}$ (which follow from~\eqref{eq:Cr:Cu:Cb:def} and ~\eqref{eq:self-similar:ansatz:R}), one may show that with the self-similar ansatz~\eqref{eq:self-similar:ansatz}, the  Euler equations~\eqref{eq:euler2} may be rewritten as 
\begin{subequations} 
\label{eq:euler3}
\begin{align} 
\p_\tau \Sigma 
-  \cu \Sigma 
+ \cx R \p_R \Sigma 
+  U \p_R \Sigma  
+ \alpha  \Sigma ( \p_R U + \tfrac{d-1}{R} U)
& = 0  ,  \label{eq:sigma3} \\
\p_\tau U 
- \cu U 
+  \cx    R \p_R U 
+  U \p_R  U  
+ \alpha    \Sigma \p_R \Sigma 
- \tfrac{\alpha}{\gamma} \Sigma^2 \tfrac{\p_R B}{B}
& = 0  , \label{eq:u3} \\
\p_\tau B 
 - \cb  B 
+   \cx  R \p_R B 
+  U \p_R  B
&=0 . \label{eq:s3}
\end{align} 
\end{subequations} 
The self-similar evolution~\eqref{eq:euler3} with $(R,\tau) \in [0,\infty) \times [0,\infty)$ is equivalent to the Eulerian evolution~\eqref{eq:euler2} with $(r,t) \in [0,\infty) \times [-1,t_*)$, via the transformations in~\eqref{eq:self-similar:ansatz}. 

In analogy to~\eqref{eq:profile:normalize}, since in this section we are working with solutions of~\eqref{eq:euler3} for which the density is $\OO_{R\to 0}(1)$, the radial velocity profile is $\OO_{R\to0}(R)$, and the pressure is $\OO_{R\to 0}(R^2)$, it makes sense to normalize the  solution $(U,\Sigma,B)$ of~\eqref{eq:euler3} near $R=0$ by introducing  
\begin{equation}
V(R,\tau) = \tfrac{1}{R} U(R,\tau), \qquad 
Q(R,\tau) = \tfrac{1}{R} \Sigma(R,\tau) , \qquad
H(R,\tau) = \tfrac{1}{R} B(R,\tau) .
\label{eq:class:of:perturbations} 
\end{equation}
The equivalent unknowns $(V,Q,H)$ solve the system 
\begin{subequations} 
\label{eq:euler:final}
\begin{align} 
&\p_\tau Q 
+ \bigl( \cx - \cu +  (1 + \alpha d)  V\bigr) Q 
+ \bigl( \cx +  V \bigr) R \p_R Q   
+  \alpha   Q  R \p_R    V   = 0  
,  \label{eq:euler:final:a} \\
&\p_\tau V  
+ \bigl(\cx - \cu \bigr)  V 
+ \bigl( V^2 + \tfrac{2\alpha^2}{\gamma} Q^2 \bigr)
+ \bigl( \cx  +  V \bigr) R \p_R V   
+ \alpha  Q R \p_R Q  
- \tfrac{\alpha}{\gamma}  Q^2 \tfrac{R \p_R H }{H} = 0  , \label{eq:euler:final:b} \\
&\p_\tau H
+ \bigl( \cx  - \cb +  V\bigr) H 
+ \bigl( \cx  +  V \bigr) R \p_R H    =0  ,\label{eq:euler:final:c}
\end{align} 
where the $\tau$-dependent modulation functions $(\cx,\cu,\cb)$ are to be chosen suitably.

It is also convenient to denote\footnote{In writing $\log H$ we are implicitly assuming that $H > 0$. We know this positivity holds true for the steady state $\bar H$. Moreover, it is clear from the maximum principle that if $H|_{\tau=0} >0$ for all $R$, then the evolution equation~\eqref{eq:euler:final:c} preserves this positivity.}
\begin{equation}
K = \log H 
,
\label{eq:euler:final:d} 
\end{equation}
so that the last term on the left side of \eqref{eq:euler:final:b} may be rewritten as $- \frac{\alpha}{\gamma}  Q^2 R \p_R K$, and the evolution equation~\eqref{eq:euler:final:c} may be replaced by 
\begin{equation}
 \p_\tau K
+ \bigl( \cx  - \cb +   V\bigr)   
+ \bigl( \cx  +   V \bigr) R \p_R K   =0 \,.\label{eq:euler:final:e}
\end{equation}
\end{subequations} 
Henceforth we will interchangeably refer to the solution of~\eqref{eq:euler:final} as $(V,Q,H)$ or $(V,Q,K)$. 

\subsection{Definition of the perturbation and their evolution}
We write the solution $(V,Q,K)$ of~\eqref{eq:euler:final} as an additive perturbation of the stationary profile $(\bar V,\bar Q,\bar K)$ which solves~\eqref{eq:euler4}, as 
\begin{subequations}
\label{eq:L:infty:perturbation}
\begin{equation}
V = \bar V + \tilde V , \qquad
Q = \bar Q + \tilde Q , \qquad
K = \bar K + \tilde K = \log \bar H + \log (1+ \tilde H) .
\label{eq:L:infty:perturbation:a}
\end{equation}
Note that~\eqref{eq:L:infty:perturbation:a} implies $H = \bar H (1+ \tilde H)$; this  multiplicative nature of the perturbation for $H$, stems from the fact that  \eqref{eq:euler:final} truly is  an equation for $K = \log H$, and that an additive perturbation should be considered at the level of these logarithms. 

Similarly, we denote the similarity exponents $(\cx,\cu,\cb)$ as perturbations of their limiting values as $\tau \to \infty$, as 
\begin{equation}
\cx = \cxbar +  \cxtilde , \qquad
\cu = \cubar + \cutilde = \cxbar - 1 +  \cutilde , \qquad
\cb = \cbbar +  \cbtilde .
\label{eq:L:infty:perturbation:b}
\end{equation}
\end{subequations}

Our goal is to show that if $(\tilde V,\tilde Q,\tilde K)|_{\tau = 0}$ is sufficiently small (in a suitable sense), then by suitably choosing $( \cxtilde,  \cutilde,  \cbtilde)(\tau)$ to decay exponentially fast as $\tau \to \infty$, the unique local-in-$\tau$ smooth solution $(V,Q,K)$ of~\eqref{eq:euler:final} can be extended to a global-in-$\tau$ solution, and that for this solution we have $(\tilde V,\tilde Q,\tilde K) \to 0$  as $\tau \to \infty$, thereby establishing the asymptotic stability of the profile $(\bar V,\bar Q,\bar K)$; see Theorem~\ref{thm:stability:radial} below.

In order to achieve the aforementioned goal, we analyze the evolution equations satisfied by the perturbation $(\tilde V,\tilde Q,\tilde K)$, which are deduced using~\eqref{eq:euler5}, \eqref{eq:euler:final}, and~\eqref{eq:L:infty:perturbation}:
\begin{subequations} 
\label{eq:euler:tilde}
\begin{align} 
&\p_\tau \tilde Q 
+ \bigl( \cx +  V \bigr) R \p_R \tilde Q   
+  \alpha  Q  R \p_R \tilde V  
+ \mathcal{D}_{\tilde Q} 
+  \mathcal{N}_{\tilde Q} 
+ \mathcal{F}_{\tilde Q} 
= 0  
,  \label{eq:euler:tilde:a} \\
&\p_\tau \tilde V  
+ \bigl( \cx  +  V \bigr) R \p_R \tilde V   
+ \alpha  Q R \p_R \tilde Q  
- \tfrac{\alpha}{\gamma}  Q^2  R \p_R \tilde K
+  \mathcal{D}_{\tilde V} 
+  \mathcal{N}_{\tilde V} 
+ \mathcal{F}_{\tilde V} 
= 0  , \label{eq:euler:tilde:b} \\
&\p_\tau \tilde K
+ \bigl( \cx  +  V \bigr) R \p_R \tilde K    
+  \mathcal{D}_{\tilde K} 
+ \mathcal{F}_{\tilde K} 
=0 ,
\label{eq:euler:tilde:c}
\end{align} 
where we have kept the  transport terms explicitly, we have denoted the \emph{linear damping} terms  by
\begin{align}
\mathcal{D}_{\tilde Q}  
&:=   
\bigl(1 + (1 + \alpha d)  \bar V + \alpha R \p_R \bar V \bigr) \tilde Q
+
\bigl(  (1 + \alpha d)  \bar Q + R \p_R \bar Q \bigr) \tilde V 
,
\label{eq:D:Q:def}
\\
\mathcal{D}_{\tilde V} &:=  
(1 + 2 \bar V + R \p_R \bar V  ) \tilde V 
+ 
(\tfrac{4\alpha^2}{\gamma} \bar Q + \alpha R \p_R \bar Q - 
\tfrac{2\alpha}{\gamma}  \bar Q  R \p_R \bar K ) \tilde Q     
,
\label{eq:D:V:def}
\\
\mathcal{D}_{\tilde K} &:=
(1 +   R \p_R \bar K )  \tilde V 
,
\label{eq:D:H:def}
\end{align}
the \emph{nonlinear}  terms are given by
\begin{align}
\mathcal{N}_{\tilde Q} &:= 
( \cxtilde - \cutilde ) \tilde Q
+  (1 + \alpha d)  \tilde V \tilde Q
,
\label{eq:N:Q:def}
\\
\mathcal{N}_{\tilde V} &:= 
\bigl(\cxtilde - \cutilde \bigr)  \tilde V 
+
 (\tfrac{2\alpha^2}{\gamma} -
\tfrac{\alpha}{\gamma} R \p_R \bar K ) \tilde Q^2
,
\label{eq:N:V:def}
\end{align}
and the \emph{modulation forcing} terms are 
\begin{align}
\mathcal{F}_{\tilde Q} &:= 
(\cxtilde - \cutilde) \bar Q
+
\cxtilde  R \p_R \bar Q
,
\label{eq:F:Q:def}
\\
\mathcal{F}_{\tilde V} &:= 
\bigl(\cxtilde - \cutilde \bigr)  \bar V   
+ \cxtilde  R \p_R \bar V  
, 
\label{eq:F:V:def}
\\
\mathcal{F}_{\tilde K} &:= 
  - \cbtilde + \cxtilde \bigl( 1 +  R \p_R \bar K  \bigr)
.
\label{eq:F:H:def}
\end{align}
\end{subequations} 
 
\subsection{Taylor coefficients at \texorpdfstring{$R=0$}{R=0}}
We recall from~\eqref{eq:first:nontrivial:Taylor:coefficient}, ~\eqref{eq:power:series:V:Q:R=0} and~\eqref{eq:bar:H:R=0} that as $R\to 0$ we have 
\begin{subequations}
\label{eq:Taylor:R=0:all}
\begin{alignat}{2}
\bar V(R) &= \bar v_0 + \bar v_{\NNN} R^{2\NNN} +  \OO_{R\to 0}(R^{4\NNN}) ,  
\qquad && \bar v_0 = \bar V(0)\, , \qquad \bar v_{\NNN} = 1 , \\
\bar Q(R) &= \bar q_0 + \bar q_{\NNN} R^{2\NNN} +  \OO_{R\to 0}(R^{4\NNN}) , 
\qquad && \bar q_0 = \bar Q(0) , \qquad  \bar q_{\NNN} =  - \tfrac{\bar q_0(1 + \alpha d + 2 \alpha \NNN)}{2 \NNN (\cxbar + \bar v_0)} , \\
\bar H(R) &= \bar h_0  + \bar h_{\NNN} R^{2\NNN} +  \OO_{R\to 0}(R^{4\NNN}) ,
\qquad &&\bar h_0 = 1, \qquad  \bar h_{\NNN} = - \tfrac{1}{2 \NNN (\cxbar + \bar v_0)} .
\label{eq:Taylor:R=0:all:c}
\end{alignat}
Using that $\bar K = \log \bar H$, relation~\eqref{eq:Taylor:R=0:all:c} may be rewritten as
\begin{alignat}{2}
\bar K(R) &= \bar k_0  + \bar k_{\NNN} R^{2\NNN} + \mathcal{O}_{R\to 0}(R^{4\NNN}) ,
\qquad &&\bar k_0 = \log  \bar h_0 = 0, \qquad  \bar k_{\NNN} =  \bar h_{\NNN} = - \tfrac{1}{2 \NNN (\cxbar + \bar v_0)} .
\label{eq:Taylor:R=0:all:d}
\end{alignat}
\end{subequations}
In the stability analysis performed in this section, we only consider perturbations $(\tilde V,\tilde Q, \tilde K)$ which are consistent with the Taylor series expansions in~\eqref{eq:Taylor:R=0:all}, namely, such that 
\begin{subequations}
\label{eq:Taylor:R=0:tilde}
\begin{align}
\tilde V(R,\tau) &= \tilde v_0(\tau) + \tilde v_{\NNN}(\tau) R^{2\NNN} + \OO_{R\to 0}(R^{2\NNN+2}) , \\
\tilde Q(R,\tau ) &= \tilde q_0(\tau) + \tilde q_{\NNN}(\tau) R^{2\NNN}  +  \OO_{R\to 0}(R^{2\NNN+2}) , \\
\tilde K(R,\tau ) &= \tilde k_0(\tau) +  \tilde k_{\NNN}(\tau)  R^{2\NNN} + \OO_{R\to 0}(R^{2\NNN+2}) .
\end{align}
\end{subequations}
In writing~\eqref{eq:Taylor:R=0:tilde} we have implicitly assumed that the order of vanishing at $R=0$ of the functions $(V,Q,K)$, and hence also of $(\tilde V, \tilde Q, \tilde K)$, does not change in time. This is indeed the case, and it is a consequence of the following lemma.
\begin{lemma}[\bf Order of vanishing at $R=0$]
\label{lem:vanishing:at:R=0} Fix any integer $\LLL \geq 2 \NNN+1$.
Define $\mathcal{L}:= \{ \ell \in \Naturals \colon 1 \leq \ell \leq 2 \NNN - 1 \} \cup \{\ell \in 2 \Naturals+ 1\colon 2\NNN+1 \leq \ell \leq \LLL \}$. 
Assume that the initial conditions  of \eqref{eq:euler:tilde:a}--\eqref{eq:euler:tilde:c} are such that 
$(\p_R^\ell \tilde V, \p_R^\ell \tilde Q, \p_R^\ell \tilde K)(0,0) = 0$ for all integers $\ell \in \mathcal{L}$. Moreover, assume that the solution $(\tilde V,\tilde Q,\tilde K)$ of \eqref{eq:euler:tilde:a}--\eqref{eq:euler:tilde:c} belongs to $L^\infty([0,T]; W^{\LLL+1,\infty}_{\rm loc}(\Reals_+))$ for $T>0$. Then, we have that $(\p_R^\ell \tilde V, \p_R^\ell \tilde Q, \p_R^\ell \tilde K)(0,\tau) = 0$ for all integers  $\ell \in \mathcal{L}$, for all times $\tau \in (0,T]$.
\end{lemma}
\begin{proof}[Proof of Lemma~\ref{lem:vanishing:at:R=0}]
We consider the evolution equation \eqref{eq:euler:final:a}, \eqref{eq:euler:final:b}, \eqref{eq:euler:final:d} for the unknowns $(Q,  V,  K)$. Denote $(Q,V,K)(0,\tau) = (q_0,v_0,k_0)(\tau)$. 
For each integer $1 \leq \ell \leq \LLL$, we apply $\frac{1}{\ell!} \p_R^\ell$ to \eqref{eq:euler:final} and restrict the resulting at $R=0$ to obtain the system of ODEs
\begin{subequations}
\label{eq:coeff:evo}
\begin{align}
&\p_\tau
(\tfrac{\p_R^\ell Q}{\ell!} , \tfrac{\p_R^\ell V}{\ell!}  , \tfrac{\p_R^\ell K}{\ell!}  )|_{R=0}^\intercal
+ 
(E_0  + \ell  E_1  )
(\tfrac{\p_R^\ell Q}{\ell!} , \tfrac{\p_R^\ell V}{\ell!}  , \tfrac{\p_R^\ell K}{\ell!}  )|_{R=0}^\intercal
=
-  {\bf 1}_{\ell \geq 2}
 Z^{(\ell)}
.
\label{eq:coeff:evo:a}
\end{align}
where
\begin{align}
&E_0:= 
\begin{pmatrix}
\cx - \cu + (1+\alpha d)   v_0
&  (1+\alpha d)   q_0  
& 0\\
\frac{4\alpha^2}{\gamma}   q_0  
& \cx - \cu + 2   v_0  
& 0 \\
0
& 1
& 0
\end{pmatrix}
,
\qquad
E_1:=
\begin{pmatrix}
\cx  +   v_0
&\alpha    q_0 
& 0\\
\alpha    q_0 
&\cx +    v_0
& -\frac{\alpha}{\gamma}      q_0^2 \\
0
&0
& \cx  +    v_0
\end{pmatrix}
,
\label{eq:E0:E1:def}
\\
&
{\footnotesize 
Z^{(\ell)} :=
\begin{pmatrix}
\sum_{i=1}^{\ell-1}   (\tfrac{\p_R^i V}{i!}  \tfrac{\p_R^{\ell-i} Q}{(\ell-i)!} )|_{R=0}   \bigl( \ell- (1-\alpha) i   + (1+\alpha d) \bigr)
\\
\sum_{i=1}^{\ell-1} 
 (i + 1)  (\tfrac{\p_R^i V}{i!}  \tfrac{\p_R^{\ell-i} V}{(\ell-i)!} )|_{R=0}  
+ (\alpha i + \tfrac{2\alpha^2}{\gamma})  (\tfrac{\p_R^i Q}{i!}  \tfrac{\p_R^{\ell-i} Q}{(\ell-i)!} )|_{R=0}    
- \tfrac{\alpha}{\gamma} \sum_{i=1}^{\ell-1}  \sum_{j=0}^{\ell-i} 
i (\tfrac{\p_R^i K}{i!}  \tfrac{\p_R^j Q}{j!} \tfrac{\p_R^{\ell-i-j} Q}{(\ell-i-j)!} )|_{R=0}  
\\
\sum_{i=1}^{\ell-1} (\ell-i)  (\tfrac{\p_R^i V}{i!}  \tfrac{\p_R^{\ell-i} K}{(\ell-i)!} )|_{R=0}
\end{pmatrix}
}
\notag
.
\end{align}
\end{subequations}
We observe that the matrices $E_0$ and $E_1$ only depend on the Taylor coefficients $v_0$ and $q_0$, while the forcing functions $Z^{(\ell)}$ only depend on the Taylor coefficients with indices $\leq \ell-1$ and  $\geq 1$. It is also convenient to denote $Z^{(1)}= (0,0,0)^\intercal$, when the sums defining $Z^{(\ell)}$ are empty,

By construction, we have that $\tfrac{1}{\ell!} \p_R^\ell \bar Q |_{R=0}=  \tfrac{1}{\ell!} \p_R^\ell \bar V|_{R=0} =  \tfrac{1}{\ell!} \p_R^\ell \bar K|_{R=0} = 0$ for all $\ell \in \mathcal{L}$. Thus, proving the lemma is equivalent to showing that $\tfrac{1}{\ell!} \p_R^\ell   Q |_{R=0} (\tau) =  \tfrac{1}{\ell!} \p_R^\ell   V|_{R=0} (\tau) =  \tfrac{1}{\ell!} \p_R^\ell   K|_{R=0} (\tau) = 0$ for all $\ell \in \mathcal{L}$ and $\tau \in [0,T]$. Since by assumption the initial data is chosen such that $\tfrac{1}{\ell!} \p_R^\ell   Q |_{R=0} (0) =  \tfrac{1}{\ell!} \p_R^\ell   V|_{R=0} (0) =  \tfrac{1}{\ell!} \p_R^\ell   K|_{R=0} (0) = 0$ for all $\ell \in \mathcal{L}$, the Lemma directly follows from the evolution equation~\eqref{eq:coeff:evo:a}, if we are able to show that $Z^{(\ell)}(\tau) = (0,0,0)^\intercal$ for all $\ell \in \mathcal{L}$. The latter fact follows by induction on $\ell$. For the base step $\ell=1$, by definition it immediately follows that $Z^{(1)} (\tau)= (0,0,0)^\intercal$. For the induction step, consider two separate cases: $\ell \in \mathcal{L}$ with $\ell \leq 2\NNN-1$ and $\ell \in \mathcal{L}$ which is odd and $\ell \geq 2 \NNN+1$. For the first case, the induction assumption implies that $\tfrac{1}{i!} \p_R^i   Q |_{R=0} (\tau) =  \tfrac{1}{i!} \p_R^i   V|_{R=0} (\tau) =  \tfrac{1}{i!} \p_R^i   K|_{R=0} (\tau) = 0$ for all $1 \leq i \leq \ell -1$; hence, the sum defining $Z^{(\ell)}$ vanishes identically. For the second case, we note that if $\ell$ is odd and $1\leq i \leq \ell-1$, then either $i$ or $\ell-i$ is also odd; hence every summand in the sum defining $Z^{(\ell)}$ vanishes, concluding the proof of the induction step and hence of the Lemma.
\end{proof}

\subsection{Choice of modulation functions}
In this section we choose the modulation functions $(\cx,\cu,\cb)$ in~\eqref{eq:euler:tilde} to ensure that $(\tilde v_0,\tilde v_{\NNN},\tilde q_0, \tilde q_{\NNN},\tilde k_0, \tilde k_{\NNN})(\tau) \to 0 $ as $\tau \to \infty$; see~\eqref{eq:modulation:final} below.

\subsubsection{The coefficients of $R^0$}
Restricting the dynamics~\eqref{eq:euler:tilde:a}--\eqref{eq:euler:tilde:c} at $R=0$, appealing to~\eqref{eq:bar:V:zero}, \eqref{eq:Taylor:R=0:all}, noting that the transport-terms vanish at $R=0$, and using the notation in~\eqref{eq:Taylor:R=0:tilde}  results in 
\begin{subequations}
\label{eq:euler:tilde:zero} 
\begin{align}
\label{eq:euler:tilde:vq:zero}
&\tfrac{d}{d\tau}
\begin{pmatrix}
\tilde q_0 \\
\tilde v_0  
\end{pmatrix} 
+ 
\underbrace{
\begin{pmatrix}
0 & (1+ \alpha d) \bar q_0  \\
\tfrac{4\alpha^2}{\gamma} \bar q_0 & 1+ 2 \bar v_0  
\end{pmatrix}
}_{=: A_0}
\begin{pmatrix}
\tilde q_0 \\
\tilde v_0  
\end{pmatrix}
+ (\cxtilde - \cutilde)
\begin{pmatrix}
  \bar q_0 \\
 \bar v_0  
\end{pmatrix}
= - 
\begin{pmatrix}
\mathcal{N}_{\tilde Q}(0,\tau) \\
\mathcal{N}_{\tilde V}(0,\tau)  
\end{pmatrix}
,
\end{align}
and the decoupled evolution
\begin{align}
\tfrac{d}{d\tau} \tilde k_0
+ \tilde v_0
+  \cxtilde - \cbtilde  
= 0.
\label{eq:euler:tilde:h:zero}
\end{align}
\end{subequations} 

The matrix $A_0$ present on the left side of~\eqref{eq:euler:tilde:vq:zero} has:
\begin{itemize}[leftmargin=1em]
\item a positive (stable) eigenvalue of $\tfrac{2\alpha d}{1+\alpha d} = 2(1+\bar v_0)$ with eigenvector  $( \gamma, 4 \alpha^2 \bar q_0)^\intercal$;
\item a negative (unstable) eigenvalue of $-1$ with eigenvector $(\bar q_0,\bar v_0)^\intercal$.
\end{itemize} 
The first constraint on our modulation functions ensures that this unstable direction is ``removed''. For this purpose, we impose the constraint 
\begin{equation}
\cxtilde(\tau) - \cutilde(\tau) 
= 
- \tfrac{4\alpha^2}{\gamma} \tfrac{\bar q_0}{\bar v_0} \tilde q_0(\tau)
+ \tfrac{1}{\bar v_0} \tilde v_0(\tau)
.
\label{eq:modulation:1} 
\end{equation}
With~\eqref{eq:modulation:1}, \eqref{eq:euler:tilde:vq:zero}   becomes 
\begin{align}
\tfrac{d}{d\tau}
(
\tilde q_0 ,
\tilde v_0  
)^\intercal
+ 
\tfrac{2\alpha d}{1+\alpha d} 
(
\tilde q_0 ,
\tilde v_0  
)^\intercal
= - 
(
\mathcal{N}_{\tilde Q}(0,\tau) ,
\mathcal{N}_{\tilde V}(0,\tau)  
)^\intercal
.
\label{eq:euler:tilde:vq:zero:new}
\end{align}
It is now evident from~\eqref{eq:euler:tilde:vq:zero:new} that $(\tilde q_0,\tilde v_0)(\tau)$ is under control as $\tau \to \infty$; quantitative bounds are established in~\eqref{eq:boot:0} below.

Returning to~\eqref{eq:euler:tilde:h:zero}, we impose the second constraint on our modulation functions, whose goal is to ensure that $\tilde k_0(\tau)$ is ``damped'' as $\tau \to \infty$. The precise damping coefficient is free,  and for simplicity, and consistency with~\eqref{eq:euler:tilde:vq:zero:new}, we choose 
\begin{equation}
 \cxtilde(\tau) - \cbtilde (\tau) 
= \tfrac{2\alpha d}{1+\alpha d} \tilde k_0(\tau) - \tilde v_0(\tau)
.
\label{eq:modulation:2}
\end{equation}
With~\eqref{eq:modulation:2}, the evolution \eqref{eq:euler:tilde:h:zero} becomes
\begin{align}
\tfrac{d}{d\tau} \tilde k_0
+ \tfrac{2\alpha d}{1+\alpha d} \tilde k_0
= 0.
\label{eq:euler:tilde:h:zero:new}
\end{align}

Combining~\eqref{eq:euler:tilde:vq:zero:new} and~\eqref{eq:euler:tilde:h:zero:new}, we obtain
\begin{align}
\tfrac{d}{d\tau}
(
\tilde q_0 ,
\tilde v_0  ,
\tilde k_0
)^\intercal
+ 
\tfrac{2\alpha d}{1+\alpha d} 
(
\tilde q_0 ,
\tilde v_0 ,
\tilde k_0 
)^\intercal
= - 
(
\mathcal{N}_{\tilde Q}(0,\tau) ,
\mathcal{N}_{\tilde V}(0,\tau) , 
0 
)^\intercal
,
\label{eq:euler:tilde:vq:hq:zero:total}
\end{align}
From~\eqref{eq:euler:tilde:vq:hq:zero:total} it is now evident that we should expect global-in-time stability of the vector $(\tilde q_0,\tilde v_0,\tilde k_0)(\tau)$. Note that~\eqref{eq:modulation:1} and~\eqref{eq:modulation:2} only define two of the three modulation functions (namely $\cutilde$ and $\cbtilde$), and that one of the modulation functions (namely $\cxtilde$) is still free.

\subsubsection{The coefficients of $R^{2\NNN}$}
By either using~\eqref{eq:coeff:evo}, or by differentiating~\eqref{eq:euler:tilde} $2\NNN$ times with respect to $R$, appealing to~\eqref{eq:Taylor:R=0:all}, \eqref{eq:Taylor:R=0:tilde},  restricting the resulting dynamics at $R=0$, and dividing by $(2\NNN)!$, in analogy to~\eqref{eq:euler:tilde:zero} we obtain
\begin{subequations}
\label{eq:euler:tilde:vqh:2N}
\begin{equation}
\tfrac{d}{d\tau}
(
\tilde q_{\NNN} ,
\tilde v_{\NNN} ,
\tilde k_{\NNN}
)^\intercal
+ A_{1}^{(\NNN)} 
(
\tilde q_{\NNN} ,
\tilde v_{\NNN} ,
\tilde k_{\NNN}
)^\intercal
+ A_2^{(\NNN)}(\tau)
(
\bar q_{\NNN} ,
\bar v_{\NNN} ,
\bar k_{\NNN}
)^\intercal
= 
\mathsf{Nonlinear}_{\eqref{eq:tilde:vqh:2N}}(\tau)
\label{eq:tilde:vqh:2N}
,
\end{equation}
where we have denoted
\begin{align}
A_{1}^{(\NNN)}
&:=
\begin{pmatrix}
0 & (1+ \alpha d) \bar q_0 & 0 \\ 
 \frac{4\alpha^2}{\gamma} \bar q_0 &  1 + 2 \bar v_0& 0 \\ 
0 & 1 & 0
\end{pmatrix}
+ 
2\NNN
\begin{pmatrix}
\cxbar + \bar v_0  &  \alpha   \bar q_0 & 0 \\ 
 \alpha   \bar q_0  & \cxbar + \bar v_0 & -\frac{\alpha}{\gamma} \bar q_0^2 \\ 
0 & 0 & \cxbar + \bar v_0
\end{pmatrix} 
\label{eq:euler:tilde:vqh:2N:b}
\\
A_2^{(\NNN)}(\tau)
&:=
\begin{pmatrix}
\cxtilde - \cutilde + (1+\alpha d)  \tilde v_0  
&  (1+\alpha d)  \tilde q_0  
& 0 
\\
\tfrac{4\alpha^2}{\gamma}   \tilde q_0  
&  \cxtilde - \cutilde + 2  \tilde v_0  
& 0   
\\
0 
& 0 
& 0  
\end{pmatrix}
+ 2 \NNN 
\begin{pmatrix}
 \cxtilde  +  \tilde v_0 
& \alpha  \tilde q_0  
& 0 
\\
\alpha  \tilde q_0  
&  \cxtilde  + \tilde v_0 
& - \tfrac{2\alpha}{\gamma}  \bar q_0  \tilde q_0     
\\
0 
& 0 
& \cxtilde  + \tilde v_0  
\end{pmatrix}
\notag
\\
&=: 2 \NNN \bigl(\cxtilde(\tau)  + \tilde v_0(\tau) \bigr) \mathrm{I}_3 + \tilde A_2^{(\NNN)}(\tau) ,
\label{eq:euler:tilde:vqh:2N:c}
\end{align}
and the nonlinear terms are given by
\begin{align}
&\mathsf{Nonlinear}_{\eqref{eq:tilde:vqh:2N}}
:= - \tfrac{1}{(2\NNN)!}
\begin{pmatrix}
(\p_R^{2\NNN} \NN_{\tilde Q})(0,\tau)  \\
(\p_R^{2\NNN} \NN_{\tilde V})(0,\tau)  \\
(\p_R^{2\NNN} \NN_{\tilde K})(0,\tau)  
\end{pmatrix}
- 2 \NNN
\begin{pmatrix}
\cxtilde  + \tilde v_0  
&
\alpha  \tilde q_0 
&
0
\\
\alpha  \tilde q_0  
& \cxtilde  + \tilde v_0  
& - \tfrac{\alpha}{\gamma}  \tilde q_0 (2 \bar q_0 + \tilde q_0)\\
0 
&
0
& \cxtilde  + \tilde v_0  
\end{pmatrix}
\begin{pmatrix}
\tilde q_\NNN \\
\tilde v_\NNN \\
\tilde k_\NNN
\end{pmatrix}
\label{eq:euler:tilde:vqh:2N:d}
.
\end{align}
\end{subequations}
The terms on the right side of~\eqref{eq:euler:tilde:vqh:2N:d} contain solely  nonlinear terms. With the notation  from~\eqref{eq:E0:E1:def}, we recognize $A_1^{(\NNN)} = \bar E_0 + 2 \NNN \bar E_1$, $A_2^{(\NNN)}$ represents the components of $\tilde E_0 + 2 \NNN \tilde E_1$ which are linear in the the ``tilde terms''. 
In \eqref{eq:euler:tilde:vqh:2N:c}, we isolate the term 
$ \cxtilde(\tau)  + \tilde v_0(\tau)$ and will use it to determine the modulation function $\cxtilde(\tau) $ below; the ``tilde term" $\tilde A_2^{(\NNN)}(\tau)$ in \eqref{eq:euler:tilde:vqh:2N:c} depends only on 
the $0$-th order Taylor coefficients $ \td v_0, \td q_0, \td k_0$ (recalling \eqref{eq:modulation:1}
for $  \cxtilde - \cutilde$), which satisfy the ODE \eqref{eq:euler:tilde:vq:hq:zero:total} 
with a linear damping term and are expected to be stable global-in-time. 
The terms in~\eqref{eq:euler:tilde:vqh:2N:d} are the nonlinear terms; note that the second term in~$\mathsf{RHS}_{\eqref{eq:euler:tilde:vqh:2N:d}}$ is equal to $- 2 \NNN \tilde E_1 (\tilde q_\NNN,\tilde v_\NNN, \tilde k_\NNN)^\intercal$. We remark that the \emph{explicit form} of $\tilde A_2^{(\NNN)}(\tau)$ in~\eqref{eq:euler:tilde:vqh:2N:c} and the nonlinear terms $\mathsf{Nonlinear}_{\eqref{eq:tilde:vqh:2N}}$ 
are \emph{not used} in the stability analysis.

The matrix $A_1^{(\NNN)}$ defined in \eqref{eq:euler:tilde:vqh:2N:b} has the following eigensystem:
\begin{itemize}[leftmargin=1em]
\item a positive (stable) eigenvalue  
\begin{subequations}
\label{eq:A1:eigensystem}
\begin{equation}
\lambda_{\NNN}^\sharp := 2\NNN (\cxbar + \bar v_0),
\end{equation}
with eigenvector 
\begin{equation}
 {p}_{\NNN}^\sharp := 
\bigl(\bar q_0, 0, \gamma + \tfrac{2\alpha}{\NNN}\bigr)^\intercal \,;
\end{equation}
\item a positive (stable) eigenvalue\footnote{Here we have appealed to \eqref{eq:det:Mn:def}--\eqref{eq:cx:admissible}.}
\begin{equation}
\lambda_{\NNN}^\flat := 4 \NNN |\bar v_0| \sqrt{\tfrac{\alpha \gamma d}{2} + \mathsf{E}_{\NNN} + \tfrac{(1-\alpha d)^2}{16 \NNN^2}} = 4 \NNN (\cxbar + \bar v_0) + (1 + 2 \bar v_0) , 
\end{equation}
with eigenvector
\begin{equation}
 {p}_{\NNN}^\flat := \bigl(\tfrac{(1+\alpha d + 2 \NNN\alpha) \bar q_0}{2 \NNN ( \cxbar + \bar v_0) + 1 + 2\bar v_0}, 1 , \tfrac{1}{2 \NNN ( \cxbar + \bar v_0) + 1 + 2\bar v_0}\bigr)^\intercal \,;
\end{equation}
\item a vanishing (neutral) eigenvalue 
\begin{equation}
\lambda_{\NNN}^\dagger:=0 , 
\end{equation} 
with eigenvector
\begin{equation}
p_{\NNN}^\dagger := \bigl( - \tfrac{(1+\alpha d + 2 \NNN\alpha) \bar q_0}{2 \NNN (\cxbar + \bar v_0)}, 1 , -\tfrac{1}{2 \NNN (\cxbar + \bar v_0)} \bigr) = \bigl(\bar q_{\NNN},\bar v_{\NNN},\bar k_{\NNN} \bigr)^\intercal .
\end{equation}
\end{subequations}
\end{itemize}
Let $P$ be the matrix with columns $[p_{\NNN}^\sharp |  p_{\NNN}^\flat | p_{\NNN}^\dagger]$, so that 
\begin{equation*}
A_1^{(\NNN)} = P {\rm diag\,}(\lambda_{\NNN}^\sharp, \lambda_{\NNN}^\flat,0) P^{-1} .
\end{equation*}
It is also convenient to denote 
\begin{equation}
\label{eq:varphi:sharp:flat:def}
(\varphi^\sharp, \varphi^\flat, \varphi^\dagger)^\intercal
:= P^{-1}
(\tilde q_{\NNN}, \tilde v_{\NNN}, \tilde k_{\NNN})^\intercal
.
\end{equation}
With this notation, we multiply \eqref{eq:tilde:vqh:2N} from the left with $P^{-1}$, noting the decomposition in~\eqref{eq:euler:tilde:vqh:2N:c} and the fact that $P^{-1} p_{\NNN}^\dagger = e_3$, to obtain
\begin{equation}
\tfrac{d}{d\tau}
\begin{pmatrix}
\varphi^\sharp \\
\varphi^\flat \\
\varphi^\dagger
\end{pmatrix}
+ {\rm diag\,}(\lambda_{\NNN}^\sharp, \lambda_{\NNN}^\flat,0) 
\begin{pmatrix}
\varphi^\sharp \\
\varphi^\flat \\
\varphi^\dagger
\end{pmatrix}
+ 2 \NNN \bigl(\cxtilde + \tilde v_0 \bigr) 
\begin{pmatrix}
0\\
0 \\
1
\end{pmatrix}
= 
- P^{-1}  \tilde A_2^{(\NNN)} 
\begin{pmatrix}
\bar q_{\NNN} \\
\bar v_{\NNN} \\
\bar k_{\NNN}
\end{pmatrix}
+ P^{-1} 
\mathsf{Nonlinear}_{\eqref{eq:tilde:vqh:2N}} 
\label{eq:tilde:vqh:2N:alt}
.
\end{equation}

With~\eqref{eq:tilde:vqh:2N:alt} in mind, the third and last constraint on our modulation functions ensures that the  neutrally stable direction  corresponding to $\varphi^\dagger$ becomes stable. For this purpose, in analogy with~\eqref{eq:modulation:1} we impose
\begin{equation}
2 \NNN \bigl(\cxtilde(\tau) + \tilde v_0(\tau)  \bigr)
= \lambda_{\NNN}^\sharp \varphi^\dagger(\tau)
,
\label{eq:modulation:3} 
\end{equation}
so that \eqref{eq:tilde:vqh:2N:alt} becomes 
\begin{equation}
\tfrac{d}{d\tau}
\begin{pmatrix}
\varphi^\sharp \\
\varphi^\flat \\
\varphi^\dagger
\end{pmatrix}
+ {\rm diag\,}(\lambda_{\NNN}^\sharp, \lambda_{\NNN}^\flat,\lambda_{\NNN}^\sharp) 
\begin{pmatrix}
\varphi^\sharp \\
\varphi^\flat \\
\varphi^\dagger
\end{pmatrix}
= 
- P^{-1}  \tilde A_2^{(\NNN)} 
\begin{pmatrix}
\bar q_{\NNN} \\
\bar v_{\NNN} \\
\bar k_{\NNN}
\end{pmatrix}
+ P^{-1} 
\mathsf{Nonlinear}_{\eqref{eq:tilde:vqh:2N}} 
\label{eq:tilde:vqh:2N:alt:alt}
.
\end{equation}
Since the matrix $\tilde A_2^{(\NNN)}$ in \eqref{eq:euler:tilde:vqh:2N:c} depends only on $\td v_0, \td q_0, \td k_0$, which satisfy the ODE \eqref{eq:euler:tilde:vq:hq:zero:total} with a linear damping term, and since  $\lambda_{\NNN}^\sharp, \lambda_{\NNN}^\flat>0$, \eqref{eq:tilde:vqh:2N:alt:alt} shows that we may expect global-in-time stability for the vector $(\varphi^\sharp,\varphi^\flat,\varphi^\dagger)(\tau)$, and thus via~\eqref{eq:varphi:sharp:flat:def} we expect global-in-time stability for $(\tilde q_{\NNN}, \tilde v_{\NNN}, \tilde k_{\NNN})(\tau)$. Equation~\eqref{eq:modulation:3} provides the definition of the last modulation functions, namely $\cxtilde$.

\subsubsection{The definitions of the modulation functions}
Let $(\tilde q_0,\tilde v_0,\tilde k_0)(\tau)$ solve~\eqref{eq:euler:tilde:vq:hq:zero:total}, and let $(\varphi^\sharp, \varphi^\flat, \varphi^\dagger)(\tau)$ solve~\eqref{eq:tilde:vqh:2N:alt:alt}.
Combining~\eqref{eq:modulation:1}, \eqref{eq:modulation:2}, and \eqref{eq:modulation:3}, we have
\begin{subequations}
\label{eq:modulation:final}
\begin{align}
\cxtilde (\tau)
&=  \tfrac{\lambda_{\NNN}^\sharp}{2\NNN}  \varphi^\dagger(\tau) - \tilde v_0 (\tau) 
=  (\cxbar + \bar v_0)  \varphi^\dagger(\tau)   - \tilde v_0 (\tau)
\label{eq:modulation:2:final} ,
\\
\cutilde(\tau) 
&= \cxtilde(\tau)
+ \tfrac{4\alpha^2}{\gamma} \tfrac{\bar q_0}{\bar v_0} \tilde q_0(\tau)
- \tfrac{1}{\bar v_0} \tilde v_0(\tau)
=  (\cxbar + \bar v_0)  \varphi^\dagger(\tau) 
+ \tfrac{4\alpha^2}{\gamma} \tfrac{\bar q_0}{\bar v_0} \tilde q_0(\tau)
- \tfrac{1 + \bar v_0}{\bar v_0} \tilde v_0(\tau)
, \label{eq:modulation:1:final} \\
\cbtilde (\tau) 
&= \cxtilde(\tau) - \tfrac{2\alpha d}{1+\alpha d} \tilde k_0(\tau) + \tilde v_0(\tau)
= (\cxbar + \bar v_0)  \varphi^\dagger(\tau)    - \tfrac{2\alpha d}{1+\alpha d} \tilde k_0(\tau)  
.
\label{eq:modulation:3:final}
\end{align}
\end{subequations}

\subsection{Differentiated Riemann-type variables}
In order to avoid derivative losses at the level of (weighted) $C^1$ estimates,  and motivated by the propagation of regularity established in Proposition~\ref{prop:reg:propagation}, we define
\begin{align}
\label{eq:ringW:ringZ:ringA:def}
\Wring := R \partial_R V + R \partial_R Q - \tfrac{1}{\gamma  }Q R\partial_R K , \quad 
\Zring := R \partial_R V - R \partial_R Q + \tfrac{1}{\gamma }Q R \partial_R K , \quad
\Aring := \tfrac{\alpha }{\gamma} Q R \partial_R K .
\end{align}
In particular, we note that~\eqref{eq:ringW:ringZ:ringA:def} implies
\begin{equation}
R \p_R  Q 
=  \tfrac 12 \Wring - \tfrac 12 \Zring +  \tfrac{1}{\alpha} \Aring,
\qquad
R \p_R  V 
=  \tfrac 12 \Wring + \tfrac 12 \Zring  ,
\qquad
R \p_R K
= \tfrac{\gamma}{\alpha Q} \Aring
.
\label{eq:QVH:via:ring}
\end{equation}
For compactness of notation, it is convenient to introduce the vector
\begin{equation}
\YY: = (\Zring,\Aring,\Wring)^\intercal,
\label{eq:YY:def}
\end{equation}
and denote its components as $\YY_i$ for $i \in \{1,2,3\}$. 

\subsubsection{Evolution of the differentiated Riemann-type variables}

Let $(V,Q,K)$ solve~\eqref{eq:euler:final}. After applying $R \p_R$ to~\eqref{eq:euler:final},  and then diagonalizing the resulting system, we obtain the system of equations 
\begin{align}
\p_\tau \YY + \TT \, R \p_R \YY + \DD \, \YY +  \NN(\YY,\YY) 
= 0 , 
\label{eq:YY:evo}
\end{align}
where the diagonal ``transport matrix'' $\TT$ is given by
\begin{subequations}
\label{eq:YY:evo:aux}
\begin{equation}
\TT :=  (\cx +  V) {\rm Id} +  \alpha Q {\rm diag}(-1,0,1)  \, , 
\label{eq:YY:evo:TT}
\end{equation}
the ``damping matrix'' $\DD$ is given by 
\begin{equation}
\DD  := (\cx - \cu) {\rm Id} +     V {\rm D}^{(v)} +  \alpha Q {\rm D}^{(q)} , 
\label{eq:YY:evo:DD}
\end{equation}
where ${\rm D}^{(v)}, {\rm D}^{(q)}$ are two constant matrices which depend  only on $\alpha$ and $d$, and are explicitly given by
\begin{equation}
\label{eq:YY:evo:DD:vq}
{\rm D}^{(v)} 
:= 
\begin{pmatrix}
\tfrac{3+\alpha d}{2}  & 0 & \tfrac{1-\alpha d}{2}  
\\
0 & 1+\alpha d  & 0
\\
 \tfrac{1-\alpha d}{2}  & 0 & \tfrac{3 +\alpha d}{2} 
\end{pmatrix} 
,
\qquad 
{\rm D}^{(q)}
:=
\begin{pmatrix}
- \tfrac{d+2}{2}  & \frac{4}{\gamma}  & - \tfrac{(d+2)+2\alpha(d-2)}{2\gamma}
\\
\tfrac{1}{2\gamma}   & 0 & \tfrac{1}{2\gamma}
\\
\tfrac{(d+2)+2\alpha(d-2)}{2\gamma} & \tfrac{4}{\gamma}  & \tfrac{d+2}{2}
\end{pmatrix} 
,
\end{equation}
and where $\NN(\cdot,\cdot)$ is a  quadratic nonlinearity, written   in terms of its components  $\NN =(\NN_1,\NN_2,\NN_3)^\intercal$, each of which is a bilinear form $\Reals^3 \times \Reals^3 \to \Reals$ given explicitly as
\begin{equation}
\label{eq:YY:evo:NN}
\NN_i(\YY^\prime,\YY^{\prime\prime}) 
= (\tfrac{1 + (2-i)\alpha}{2} \YY^{\prime}_1 + (i-2) \YY^\prime_2 + \tfrac{1 + (i-2)\alpha}{2} \YY^\prime_3) \YY^{\prime\prime}_i 
+ \bigl( \tfrac{2-i}{2} + {\bf 1}_{i=2} \tfrac{\alpha}{2} \bigr) (\YY^\prime_1 + \YY^\prime_3) \YY^{\prime\prime}_2 
\, ,
\end{equation}
for each $i \in \{1,2,3\}$.
\end{subequations}

\begin{remark}[\bf Main terms for stability analysis]
While the system \eqref{eq:YY:evo}--\eqref{eq:YY:evo:aux} may appear complicated, in order to prove nonlinear stability
we mainly use the transport term $\cT R \pa_R \YY$ in \eqref{eq:YY:evo} 
and the constant damping terms $ (\cx - \cu) {\rm Id} $ in \eqref{eq:YY:evo:DD}. In particular, the damping terms associated with $ ( V {\rm D}^{(v)} +  \alpha Q {\rm D}^{(q)} ) \YY$ 
in \eqref{eq:YY:evo:DD} and the nonlinearity $\cN_i(\cdot, \cdot)$ 
are treated perturbatively in the stability analysis. \end{remark}
 
\subsubsection{Perturbations of the differentiated Riemann-type variables}
The steady state of the variable $\YY$ is
\begin{equation}
\Ybar 
:= \bigl(\Zringbar,\Aringbar,\Wringbar \bigr)^\intercal 
= \bigl(R \p_R \bar V - R \p_R \bar Q + \tfrac{1}{\gamma } \bar Q R\p_R \bar K, \tfrac{\alpha }{\gamma} \bar Q R\p_R \bar K,R \p_R \bar V + R \p_R \bar Q - \tfrac{1}{\gamma} \bar Q R\p_R \bar K \bigr)^\intercal 
.
\label{eq:ringW:ringZ:ringA:bar}
\end{equation}
We shall denote the perturbation to this steady state as 
\begin{equation}
\Ytilde 
= \YY - \Ybar  
= \bigl(\Zringtilde, \Aringtilde, \Wringtilde \bigr)^\intercal 
= \bigl(\Zring - \Zringbar, \Aring-\Aringbar,\Wring-\Wringbar \bigr)^\intercal .
\label{eq:ringW:ringZ:ringA:tilde}
\end{equation}
From~\eqref{eq:YY:evo}--\eqref{eq:YY:evo:aux} we deduce that the perturbation $\Ytilde$ solves
\begin{align}
\p_\tau \Ytilde 
+ \TT \, R \p_R \Ytilde  
+ \DD \, \Ytilde  
+ \Ttilde \, R\p_R \Ybar 
+ \Dtilde \, \Ybar 
+   \NN(\Ybar, \Ytilde)
+   \NN(\Ytilde,\Ybar)
+   \NN(\Ytilde,\Ytilde) 
= 0 , 
\label{eq:Ytilde:evo}
\end{align}
where in analogy to~\eqref{eq:YY:evo:TT}--\eqref{eq:YY:evo:DD} we may introduce
\begin{subequations}
\label{eq:Ytilde:evo:aux}
\begin{align}
\Tbar &:=  (\cxbar +  \bar V) {\rm Id} +  \alpha \bar Q {\rm diag}(-1,0,1)  \, , 
\label{eq:Ytilde:evo:Tbar}
\\
\Ttilde &:=  (\cxtilde  + \tilde V) {\rm Id} + \alpha  \tilde Q   {\rm diag}(-1,0,1)  \, , 
\label{eq:Ytilde:evo:Ttilde} 
\\
\Dbar  &:=  {\rm Id} +  \bar V {\rm D}^{(v)} +   \alpha \bar Q {\rm D}^{(q)} , 
\label{eq:Ytilde:evo:Dbar}
\\
\Dtilde  &:= (\cxtilde - \cutilde) {\rm Id} +    \tilde V  {\rm D}^{(v)} + \alpha  \tilde Q  {\rm D}^{(q)} .
\label{eq:Ytilde:evo:Dtilde}
\end{align}
\end{subequations}
In order to close the bootstrap~\eqref{eq:boot:4}, we perform $L^\infty$-type estimates for the forced transport equation~\eqref{eq:Ytilde:evo}, but only after subtracting a high enough Taylor polynomial at $R=0$ for $\Ytilde$.

\subsection{Subtracting the Taylor jet at \texorpdfstring{$R=0$}{R=0}}
Our goal is to introduce an operator $\JJJ_\MMM$, which takes as input a $C^{\MMM+1}$ smooth function, and outputs the same function minus its $\MMM^{\mbox{th}}$ order Taylor polynomial at $R=0$, divided by $R^{\MMM+1}$ (near $R=0$). In order to construct the operator $\JJJ_\MMM$, we first define a $C^\infty$ smooth non-increasing function $\phi = \phi(R)$  such that  
\begin{equation}
\phi(R) = 1  \;\; \mbox{for} \;\; R \leq R_{\sf in} ,
\qquad
\phi(R) = 0 \;\; \mbox{for} \;\; R \geq 2 R_{\sf in} ,
\qquad
 |R \p_R \phi| \leq 2 \;\; \mbox{for}\;\; R \in (R_{\sf in}, 2 R_{\sf in}),
\label{eq:phi:basic:properties}
\end{equation}
where $0 < R_{\sf in} \ll 1$ is a parameters to be chosen later (see Proposition~\ref{prop:choice:of:M:psi} below).
Then, for any $\MMM \geq 0$, and any function $f\colon \Reals_+^2 \to \Reals$ such that $f(\cdot,\tau) \in C^{\MMM}$, we define the cutoff Taylor polynomial 
\begin{subequations}
\label{eq:Taylor:poly:all}
\begin{equation}
\III_{\MMM} f (R,\tau) = \phi(R) \sum_{j=0}^{\MMM} \tfrac{1}{j!} \p_R^{j} f(0,\tau) 
R^j .
\label{eq:Taylor:poly:1}
\end{equation}
The Leibniz rule implies the useful relations
\begin{align}
\III_{\MMM} ( R \p_R f) (R,\tau) 
- R \p_R \bigl(\III_{\MMM} f(R,\tau)\bigr)
&= - \tfrac{R \p_R \phi(R)}{\phi(R)} \III_{\MMM} f(R,\tau)
\label{eq:Taylor:poly:3}
,
\\
\III_{\MMM} f(R,\tau) \III_{\MMM} g(R,\tau) -
\III_{\MMM} (f \, g)(R,\tau)
&= 
\III_\MMM^{\sf err}[f,g](R,\tau)
\label{eq:Taylor:poly:2}
,
\end{align}
where
\begin{align}
\III_\MMM^{\sf err}[f,g](R,\tau)
&:= R^{\MMM+1} \phi(R)^2 \sum_{1\leq j,k \leq \MMM; j+k \geq \MMM+1} \tfrac{1}{j! k!} \p_R^{j} f(0,\tau)  \p_R^{k} g(0,\tau) R^{j+k-\MMM-1}
\notag\\
&\qquad + 
\phi(R) (\phi(R) -1) \sum_{\ell=0}^{\MMM} \sum_{j=0}^{\ell} \tfrac{1}{j! (\ell-j)!} \p_R^j f(0,\tau) \p_R^{\ell-j} g(0,\tau) R^\ell . 
\end{align}
\end{subequations}
Second, we define a $C^\infty$ smooth  monotone increasing function $\psi = \psi(R)$ with the properties:
\begin{equation}
\psi(R) = R \;\; \mbox{for} \;\;  R \leq   2  R_{\sf in} ,
\qquad
\psi(R) = \psi(R_{\sf out})  \;\;  \mbox{for} \;\;  R \geq  R_{\sf out} ,
\qquad
\p_R \psi(R) \geq 0 \; \; \mbox{for all}\; \; R>0,
\label{eq:psi:basic:properties}
\end{equation}
where $0 < R_{\sf in} \ll 1 \ll R_{\sf out}$ are chosen in Proposition~\ref{prop:choice:of:M:psi} below. 
Lastly, we introduce the notation
\begin{equation}
\JJJ_{\MMM} f(R,\tau) := \frac{f(R,\tau) - \III_{\MMM} f (R,\tau)}{\psi(R)^{\MMM+1}}
,
\label{eq:J:M:def}
\end{equation}
to denote the Taylor remainder of $f(\cdot,\tau)$ at $R=0$, weighted by the correct power of  $R$ (when $f(\cdot,\tau) \in C^{\MMM+1}$). With our definitions, we have that 
\begin{subequations}
\label{eq:J:M:properties}
\begin{equation}
\lim_{R\to 0^+} \JJJ_{\MMM} f(R,\tau) = \tfrac{1}{(\MMM+1)!} \p_R^{\MMM+1} f(0,\tau) 
,
\qquad \mbox{and} \qquad 
\JJJ_{\MMM} f(R,\tau) =  \tfrac{1}{\psi(R)^{\MMM+1}}  f(R,\tau), \; \; \forall R\geq  2 R_{\sf in}.
\label{eq:J:M:properties:a}
\end{equation}
Moreover, the identities~\eqref{eq:Taylor:poly:2}--\eqref{eq:Taylor:poly:3} become:
\begin{align}
\JJJ_{\MMM}(f \, g) 
&= f \JJJ_{\MMM} g + \III_{\MMM} \, g \JJJ_{\MMM} f  + \III^{\sharp}_{\MMM}[f,g]
=  g \JJJ_{\MMM} f + \JJJ_{\MMM} g \, \III_{\MMM} \, f  + \III^{\sharp}_{\MMM}[f,g]
, \label{eq:J:M:properties:b} \\
\JJJ_{\MMM}(R\p_R f) 
&=
R\p_R ( \JJJ_{\MMM} f) 
+ (\MMM+1) \tfrac{R\p_R \psi}{\psi} \, \JJJ_{\MMM} f
+ \III^{\flat}_{\MMM} f 
\, ,\label{eq:J:M:properties:c}
\end{align}
where
\begin{align}
 \III^{\sharp}_{\MMM}[f,g](R,\tau) 
 &:= \tfrac{1}{\psi(R)^{\MMM+1}} \III^{\sf err}_{\MMM}[f,g](R,\tau) 
 =  \phi(R)^2 \sum_{1\leq j,k \leq \MMM; j+k \geq \MMM+1} \tfrac{1}{j! k!} \p_R^{j} f(0,\tau)  \p_R^{k} g(0,\tau) R^{j+k-\MMM-1}
\notag\\
&\qquad \qquad \qquad \qquad \qquad \qquad + 
\tfrac{\phi(R) (\phi(R) -1)}{R^{\MMM+1}} \sum_{\ell=0}^{\MMM} \sum_{j=0}^{\ell} \tfrac{1}{j! (\ell-j)!} \p_R^j f(0,\tau) \p_R^{\ell-j} g(0,\tau) R^\ell
 , \\
  \III^{\flat}_{\MMM} f (R,\tau) & := \tfrac{R \p_R \phi(R)}{\phi(R) \psi(R)^{\MMM+1}} \III_{\MMM} f(R,\tau)
=  \tfrac{\p_R \phi(R)}{R^{\MMM}} \sum_{j=0}^{\MMM} \tfrac{1}{j!} \p_R^{j} f(0,\tau) 
  .
\end{align}
The fact that $\phi$ and $\psi$ are independent of time implies that 
\begin{equation}
\p_t (\JJJ_{\MMM} f) = \JJJ_{\MMM}(\p_t f)
. \label{eq:J:M:properties:d}
\end{equation}
\end{subequations}
We emphasize that the terms $\III_{\MMM}^{\sharp}[f,g]$ and $\III_{\MMM}^{\flat} f$ do not contain any singularity at $R=0$ because
\begin{equation}
0 \leq \phi(R)^2 \leq {\bf 1}_{R \leq 2 R_{\sf in}}
,
\quad
0 \leq \tfrac{\phi(R) (\phi(R) -1)}{R^{\MMM+1}} \leq {\bf 1}_{R_{\sf in} \leq R \leq 2 R_{\sf in}} \tfrac{1}{R_{\sf in}^{\MMM+1}}
,
\quad 
- {\bf 1}_{R_{\sf in} \leq R \leq 2 R_{\sf in}} \tfrac{1}{R_{\sf in}^{\MMM}} \leq \tfrac{\p_R \phi(R)}{R^{\MMM}} \leq 0
.
\label{eq:psi:phi:quotient}
\end{equation}
Moreover, $\III_{\MMM}^{\sharp}[f,g]$ and $\III_{\MMM}^{\flat} f$ are uniquely computed from knowledge of the Taylor coefficients of $f$ and $g$ at $R=0$, up to order $\MMM$; because of this, it is fair to consider $\III_{\MMM}^{\sharp}$ and $\III_{\MMM}^{\flat}$ as exponentially decaying error terms (in $\tau$). The following result will be useful in the subsequent analysis.
\begin{lemma}
\label{lem:JM:aux:1}
Let $f \in C^{\MMM+1}$ and assume $R_{\sf in} \leq \frac 14$. For any $R \in (0,2 R_{\sf in}]$ we have the pointwise estimate
\begin{subequations}
\label{eq:J:M:R:pR:bound}
\begin{align}
|\JJJ_\MMM f(R)| 
&\leq \bigl( \tfrac{1}{\MMM+1} + {\bf 1}_{R \geq 2 R_{\sf in}} \log\tfrac{R}{2 R_{\sf in}} \bigr) \|\JJJ_\MMM (R \p_R f) \|_{L^\infty(0,R)} 
\notag\\
&\qquad + {\bf 1}_{R \geq R_{\sf in}}  \tfrac{4}{\psi(R)^{\MMM+1}}  \max_{0\leq j\leq\MMM} \tfrac{1}{(j+1)!} |\p_R^j f(0)|
.
\label{eq:J:M:R:pR:bound:a}
\end{align}
In addition, if there exists $\theta>0$ such that $\brak{R}^\theta R \p_R f \in L^\infty$, then for all $R\geq 2 R_{\sf in}$ we have
\begin{align}
 \brak{R}^{\theta} |\JJJ_\MMM f(R)| 
\leq \tfrac{2^\theta}{\theta} \Bigl( {\bf 1}_{R\geq R_{\sf out}} + {\bf 1}_{2 R_{\sf in} \leq R \leq R_{\sf out}} \bigl(\tfrac{\psi(R_{\sf out})}{\psi(R)}\bigr)^{\MMM+1}
 \Bigr) \| \brak{R}^\theta \JJJ_{\MMM}(R \p_R f)\|_{L^\infty(R,\infty)}
.
\label{eq:J:M:R:pR:bound:b}
\end{align}
\end{subequations}
The above inequalities hold for any choice of functions $\phi$ and $\psi$ which satisfy~\eqref{eq:phi:basic:properties} and~\eqref{eq:psi:basic:properties}.
\end{lemma}
\begin{proof}[Proof of Lemma~\ref{lem:JM:aux:1}]
Denote $R \p_R f = F$, then by~\eqref{eq:Taylor:poly:all} and~\eqref{eq:J:M:def} we have 
\begin{subequations}
\begin{equation}
 \JJJ_\MMM f(R) 
=  \int_0^R   \tfrac{\psi(R^\prime)^{\MMM+1}}{\psi(R)^{\MMM+1}} \JJJ_{\MMM} F(R^\prime)  \tfrac{d R^{\prime}}{R^\prime}
- \sum_{j=0}^{\MMM} \tfrac{1}{j!} \p_R^j f(0) \tfrac{1}{\psi(R)^{\MMM+1}}\int_0^R  \p_{R^\prime} \phi(R^\prime) {R^\prime}^{j}  d R^{\prime} 
,
\label{eq:J:M:R:pR:identity}
\end{equation}
while~\eqref{eq:J:M:properties:a} allows us to write 
\begin{equation}
\JJJ_{\MMM} f(R) =  - \int_R^{\infty} \tfrac{\psi(R^\prime)^{\MMM+1}}{\psi(R)^{\MMM+1}} \JJJ_\MMM F(R^\prime) \tfrac{dR^\prime}{R^\prime}
.
\label{eq:J:M:R:pR:identity:new}
\end{equation}
\end{subequations}
The bounds in~\eqref{eq:J:M:R:pR:bound} are now direct consequences of  the properties of $\phi$ and $\psi$ expressed in~\eqref{eq:phi:basic:properties} and~\eqref{eq:psi:basic:properties} (in particular the monotonicity of $\psi$ and the bound~\eqref{eq:psi:phi:quotient}), identity~\eqref{eq:J:M:R:pR:identity} for $R\leq 2 R_{\sf in}$, and identity~\eqref{eq:J:M:R:pR:identity:new} for $R\geq 2 R_{\sf in}$.
\end{proof}

\subsection{Global-in-time bootstraps} 
\label{sec:boot}
The first set of bootstrap assumptions are related to the behavior of the perturbations $(\tilde V,\tilde Q,\tilde K)$ at $R=0$; these are as follows:
\begin{itemize}[leftmargin=1em]

\item For the coefficients $(\tilde v_{\mathsf{0}}, \tilde q_{\mathsf{0}}, \tilde k_{\mathsf{0}})$ defined in~\eqref{eq:Taylor:R=0:tilde}, we make the bootstrap assumption
\begin{equation}
|\tilde q_0(\tau)| 
+ 
|\tilde v_0(\tau)|
+
|\tilde k_0(\tau)| 
\leq 
\underline{\eps} 
e^{- \underline{\lambda} \tau} ,
\label{eq:boot:0}
\end{equation}
for all $\tau\geq0$, where $\underline{\eps} \in (0,\frac 12]$ is taken to be sufficiently small with respect to $\gamma$, $d$, and $\NNN$, and 
\begin{align}
\underline{\lambda}  = \underline{\lambda} (\gamma,d,\NNN) \in (0, \tfrac{\alpha}{2(1+\alpha d)}]
\label{eq:underline:lambda:1}
\end{align}
is to be defined explicitly in~\eqref{eq:underline:lambda:def} below.

\item For the coefficients $(\tilde v_{\NNN}, \tilde q_{\NNN}, \tilde k_{\NNN})$ defined in~\eqref{eq:Taylor:R=0:tilde}, we assume that 
\begin{equation}
|\tilde q_{\NNN}(\tau)| 
+
|\tilde v_{\NNN}(\tau)| 
+
|\tilde k_{\NNN}(\tau)| 
\leq 
\underline{\eps}^{\frac 35}  e^{- \underline{\lambda} \tau} ,
\label{eq:boot:1}
\end{equation}
for all $\tau\geq0$.

\item Let $\MMM = 17 \NNN$ (see~Proposition~\ref{prop:choice:of:M:psi}). We assume that for all even integers\footnote{We recall that smooth radially-symmetric functions have power series expansions in powers of $R^2$ at $R=0$, so for all odd integers $\ell \in \{2\NNN+1,\ldots,\MMM\}$ and all $\tau>0$, Lemma~\ref{lem:vanishing:at:R=0} guarantees the vanishing of the $\ell^{\rm th}$ Taylor coefficient.}$\ell$ such that $2\NNN+1 \leq \ell \leq \MMM$ we have
\begin{equation}
\tfrac{1}{\ell!} |(\p_R^\ell \tilde Q)(0,\tau)| 
+
\tfrac{1}{\ell!} |(\p_R^\ell \tilde V)(0,\tau)| 
+
\tfrac{1}{\ell!} |(\p_R^\ell \tilde K)(0,\tau)| 
\leq 
\underline{\eps}^{\frac 25}  e^{- \underline{\lambda} \tau} ,
\label{eq:boot:2}
\end{equation}
for all $\tau\geq0$.

\item For the modulation functions defined in~\eqref{eq:modulation:final}, we make the bootstrap assumption assume that 
\begin{equation}
|\cxtilde(\tau)| + |\cutilde(\tau)| + |\cbtilde(\tau)|  
\leq \underline{\eps}^{\frac 45} e^{-\underline{\lambda} \tau}
\label{eq:boot:3}
\end{equation} 
for all $\tau \geq 0$.
\end{itemize}

We only make one bootstrap assumption on the behavior of  the functions $(R \p_R \tilde V, R\p_R \tilde Q, R\p_R \tilde K)$, for $R>0$; in fact, we consider bootstrap assumptions for the  differentiated Riemann-type variables  defined in~\eqref{eq:ringW:ringZ:ringA:def}:
\begin{itemize}[leftmargin=1em]
\item 
Let $\MMM = 17 \NNN$ and $\theta = \tfrac 12 \min\{1,\tfrac{1}{\cxbar}\}$. 
We assume that 
\begin{equation}
\brak{R}^\theta |\JJJ_{\MMM} \Ytilde (R,\tau)|
\leq 
\underline{\eps}^{\frac 15} e^{-\underline{\lambda} \tau} 
\label{eq:boot:4}
\end{equation}
for all $R>0$ and $\tau \geq 0$.
\end{itemize}
 
\subsection{Consequences of the bootstrap assumptions}
\label{sec:boot:consequence}
The bootstraps~\eqref{eq:boot:0}--\eqref{eq:boot:4} have a number of immediate useful consequences, which we list below. Before stating these consequences, we record a few bounds that pertain to the steady states $\bar V$, $\bar Q$, $\bar H$, and $\Ybar$.

\begin{remark}[\bf Estimates for $\bar V$, $\bar Q$, $\bar K$, and $\Ybar$]
\label{eq:bar:quantities:at:infinity}
The power law decay of $(\bar V,\bar Q)$ as $R \to \infty$ determined in~\eqref{eq:V:Q:R=infty}, implies that there is a sufficiently large constant $M_0>0$, which only depends on $d,\alpha$, and $\NNN$, such that 
\begin{subequations}
\begin{equation}
|\bar V(R)| + |\bar Q(R)|  + |R \p_R \bar V(R)| + |R \p_R \bar Q(R)|    
+ | (R \p_R)^2 \bar V(R)| + |(R\p_R)^2 \bar Q(R)|
\leq M_0 \brak{R}^{- \frac{1}{\cxbar}} , 
\label{eq:barV:barQ:at:infinity}
\end{equation}
and
\begin{equation}
\bar Q(R) \geq \tfrac{1}{M_0} \brak{R}^{- \frac{1}{\cxbar}} 
\label{eq:barQ:global:lower} 
\end{equation}
for all $R>0$. 
Identity~\eqref{eq:bar:H:RdR:bar:H}  and the fact that $\bar V$ is monotone increasing, implies
\begin{equation}
 \tfrac{|\bar V(R) - \bar v_0|}{\cxbar}  
\leq 
|R \p_R \bar K(R)|    
\leq \tfrac{|\bar v_0|(1+ \frac{\gamma d}{2})}{\cxbar+\bar v_0}  .
\label{eq:RdR:barH:at:infinity}
\end{equation}
for all $R>0$. After applying $R \p_R$ to~\eqref{eq:bar:H:RdR:bar:H}, and using~\eqref{eq:barV:barQ:at:infinity}--\eqref{eq:RdR:barH:at:infinity}, we obtain
\begin{equation}
|(R \p_R)^2 \bar K(R)|    
\leq \tfrac{|\bar v_0|(1+ \frac{\gamma d}{2})+   M_0}{\cxbar+\bar v_0}    .
\label{eq:RdR:twice:barH:at:infinity}
\end{equation}
\end{subequations}
Combining~\eqref{eq:ringW:ringZ:ringA:bar} with~\eqref{eq:Taylor:R=0:all} (recall that $\NNN\geq 1$), \eqref{eq:barV:barQ:at:infinity}--\eqref{eq:RdR:twice:barH:at:infinity} we  
we deduce that there exists a sufficiently large constant $M_1>0$, which only depends on $d,\alpha$, and $\NNN$, such that 
\begin{equation}
|\Ybar(R)| + |R \p_R \Ybar(R)| 
\leq M_1 R^2 \brak{R}^{-2-\frac{1}{\cxbar}},  
\label{eq:Ybar:at:infinity}
\end{equation}
for all $R>0$. 
\end{remark}

\begin{remark}[\bf Estimates for $\JJJ_{\MMM}$ applied to $\bar V$, $\bar Q$, $\bar H$ and $\Ybar$ for $R$ small]
From~\eqref{eq:Taylor:poly:1}, \eqref{eq:J:M:def}, the Taylor reminder theorem, the fact that $R_{\sf in} \leq \frac 14$, and recalling $\theta = \frac 12 \min\{1,\frac{1}{\cxbar} \}$, we have that for any $f \in C^{\MMM+1}$:
\begin{equation}
\|\brak{R}^\theta \JJJ_{\MMM} f\|_{L^\infty} \leq \tfrac{2}{(\MMM+1)!} \| \p_R^{\MMM+1} f\|_{L^\infty(0,2R_{\sf in})} + R_{\sf in}^{-\MMM-1} \| \brak{R}^\theta f\|_{L^\infty(R_{\sf in},\infty)}
. 
\end{equation}
Using the fact that $\bar V, \bar Q, \bar H$ are real analytic (hence so is $\Ybar$), and appealing to the bounds~\eqref{eq:barV:barQ:at:infinity}, \eqref{eq:Ybar:at:infinity}, for $\MMM=17\NNN$ we deduce that there exists a constant $M_2>0$, which only depends on $d,\alpha$, and $\NNN$, such that 
\begin{subequations}
\begin{align}
\|\brak{R}^\theta (\JJJ_{\MMM} \bar V,\JJJ_{\MMM} \bar Q, \JJJ_{\MMM} \Ybar,\JJJ_{\MMM} (R\p_R \Ybar))\|_{L^\infty} 
&\leq   M_2 R_{\sf in}^{-\MMM-1},
\label{eq:JM:barV:barQ:barY}
\\
\|(\III_{\MMM} \bar V,\III_{\MMM} \bar Q, \III_{\MMM} \Ybar,\III_{\MMM} (R\p_R \Ybar))\|_{L^\infty} 
&\leq   M_2 .
\label{eq:IM:barV:barQ:barY}
\end{align}
\end{subequations}
\end{remark}

\subsubsection{Bounds for $\Ytilde$}
From~\eqref{eq:Taylor:poly:1}, \eqref{eq:ringW:ringZ:ringA:def}, and \eqref{eq:ringW:ringZ:ringA:tilde}, we obtain that 
\begin{align*}
\Ytilde_1 &= R \p_R \tilde V - R \p_R \tilde Q + \tfrac{1}{\gamma} \bigl(\tilde Q  R \p_R \bar K + (\bar Q + \tilde Q)  R \p_R \tilde K \bigr)
,\\
\Ytilde_2 &=  \tfrac{\alpha}{\gamma} \bigl(\tilde Q  R \p_R \bar K + (\bar Q + \tilde Q)  R \p_R \tilde K \bigr)
,\\
\Ytilde_3 &= R \p_R \tilde V + R \p_R \tilde Q - \tfrac{1}{\gamma} \bigl(\tilde Q  R \p_R \bar K + (\bar Q + \tilde Q)  R \p_R \tilde K \bigr)
,
\end{align*}
which may be combined with~\eqref{eq:Taylor:R=0:all}--\eqref{eq:Taylor:R=0:tilde}, \eqref{eq:Taylor:poly:1}, \eqref{eq:boot:0}--\eqref{eq:boot:2},  and the fact that $ R_{\sf in} \leq \tfrac 14$ to deduce that 
\begin{align}
\label{eq:I:M:Ytilde:bound}
|\III_{\MMM} \Ytilde(R,\tau) |  + |\III_{\MMM} (R \p_R \Ytilde)(R,\tau) | 
&\leq C_{\eqref{eq:I:M:Ytilde:bound}} R^2 \underline{\eps}^{\frac 25}  {\bf 1}_{R\leq 2 R_{\sf in}} e^{-\underline{\lambda} \tau} 
,
\end{align}
where $C_{\eqref{eq:I:M:Ytilde:bound}} >0$ is a constant that only depends on $\gamma, d$, and $\NNN$ (through the Taylor series coefficients of $\bar Q$, $\bar V$, and $\bar K$ at $R=0$, of order $\leq \MMM+1 = 17 \NNN +1$). 
From~\eqref{eq:J:M:def}, \eqref{eq:boot:4}, and~\eqref{eq:I:M:Ytilde:bound} we immediately obtain
\begin{align}
|\Ytilde(R,\tau)| \leq  R^2 \brak{R}^{-2-\theta}
\bigl( \underline{\eps}^{\frac 15}  R_{\sf in}^{-2} \psi(R_{\sf out})^{\MMM+1}  + 
2 C_{\eqref{eq:I:M:Ytilde:bound}}  \underline{\eps}^{\frac 25}    \bigr) 
e^{-\underline{\lambda} \tau} 
\label{eq:Ytilde:bound}
\end{align}
for all $R > 0$ and $\tau\geq 0$.

\subsubsection{Bounds for $\tilde V$ and $\tilde Q$}
Appealing to~\eqref{eq:ringW:ringZ:ringA:bar}--\eqref{eq:ringW:ringZ:ringA:tilde}, we note that 
\begin{equation}
R\p_R \tilde V = \tfrac 12 (\Ytilde_1 + \Ytilde_3) ,
\quad
R\p_R \tilde Q = \tfrac 12 (\Ytilde_3 - \Ytilde_1 + \tfrac{2}{\alpha} \Ytilde_2) ,
\label{eq:tilde:V:tilde:Q:tilde:Y}
\end{equation}
and thus the estimates in~\eqref{eq:J:M:R:pR:bound} and the bootstrap~\eqref{eq:boot:4} imply
\begin{subequations}
\begin{align}
&\brak{R}^{\theta} |\JJJ_{\MMM} \tilde V(R,\tau)| 
+ \alpha \brak{R}^{\theta}|\JJJ_{\MMM} \tilde Q(R,\tau)|
\notag\\
&\quad 
\leq (2+\alpha) \Bigl( {\bf 1}_{0\leq R \leq 2 R_{\sf in}} \tfrac{2}{\MMM+1} 
+ {\bf 1}_{R\geq R_{\sf out}}  \tfrac{2^\theta}{\theta} \Bigr) 
\| \brak{R}^\theta \JJJ_{\MMM} \Ytilde\|_{L^\infty} 
\notag\\
&\qquad
+ (2+\alpha)  {\bf 1}_{2 R_{\sf in} \leq R \leq R_{\sf out}} \min \left\{ \brak{R}^\theta (1+ \log\tfrac{R}{2R_{\sf in}}), \bigl(\tfrac{\psi(R_{\sf out})}{\psi(R)}\bigr)^{\MMM+1} \right\} \| \brak{R}^\theta \JJJ_{\MMM} \Ytilde\|_{L^\infty} 
\notag\\
&\qquad 
+16 (1+\alpha) {\bf 1}_{ R_{\sf in} \leq R \leq R_{\sf out}}  
 \tfrac{\brak{R}^\theta}{\psi(R)^{\MMM+1}}  \underline{\eps}^{\frac 25} e^{-\underline{\lambda} \tau}
 \label{eq:J:M:tilde:V:Q:bound:alt}
 \\
 &\quad 
\leq (2+\alpha) \Bigl( {\bf 1}_{0\leq R \leq 2 R_{\sf in}} \tfrac{2}{\MMM+1} 
+ {\bf 1}_{R\geq R_{\sf out}}  \tfrac{2^\theta}{\theta} \Bigr) 
\underline{\eps}^{\frac 15}  e^{-\underline{\lambda} \tau}
\notag\\
&\qquad
+ (2+\alpha)  {\bf 1}_{2 R_{\sf in} \leq R \leq R_{\sf out}} \min \left\{ \brak{R}^\theta (1+ \log\tfrac{R}{2R_{\sf in}}), \bigl(\tfrac{\psi(R_{\sf out})}{\psi(R)}\bigr)^{\MMM+1} \right\} \underline{\eps}^{\frac 15}  e^{-\underline{\lambda} \tau} 
\notag\\
&\qquad 
+16 (1+\alpha) {\bf 1}_{ R_{\sf in} \leq R \leq R_{\sf out}}  
 \tfrac{\brak{R}^\theta}{\psi(R)^{\MMM+1}}  \underline{\eps}^{\frac 25} e^{-\underline{\lambda}\tau}
 .
 \label{eq:J:M:tilde:V:Q:bound}
\end{align}
\end{subequations}
Together with the bounds
\begin{align}
 \label{eq:I:M:tilde:V:Q:bound}
|\III_{\MMM} \tilde V(R,\tau)| 
\leq 3 \, {\bf 1}_{R\leq 2 R_{\sf in}}  \underline{\eps}^{\frac 25}  e^{-\underline{\lambda} \tau}
,\qquad 
|\III_{\MMM} \tilde Q(R,\tau)|
\leq 3 \, {\bf 1}_{R\leq 2 R_{\sf in}}  \underline{\eps}^{\frac 25} e^{-\underline{\lambda} \tau}
,
\end{align}
which are a direct consequence of the bootstraps~\eqref{eq:boot:0}--\eqref{eq:boot:2}  and definition~\eqref{eq:Taylor:poly:1}, the bound~\eqref{eq:J:M:tilde:V:Q:bound} implies
\begin{align}
& \brak{R}^{\theta} |\tilde V(R,\tau)| + \alpha \brak{R}^{\theta}|\tilde Q(R,\tau)|
\notag\\
&\leq 
(2+\alpha){\bf 1}_{0\leq R \leq 2 R_{\sf in}} 
\Bigl(  (2 R_{\sf in})^{\MMM+1}   \underline{\eps}^{\frac 15}  
+ 19  \underline{\eps}^{\frac 25} 
\Bigr) e^{-\underline{\lambda} \tau}
\notag\\
&\quad
+ (2+\alpha)  {\bf 1}_{2 R_{\sf in} \leq R \leq R_{\sf out}}
\left( \min \left\{\psi(R)^{\MMM+1} \brak{R}^\theta (1+ \log\tfrac{R}{2R_{\sf in}}),  \psi(R_{\sf out})^{\MMM+1} \right\} \underline{\eps}^{\frac 15} + 16  \brak{R_{\sf out}}^\theta   \underline{\eps}^{\frac 25} \right) e^{-\underline{\lambda} \tau} 
\notag\\
&\quad
+ 
(2+\alpha) \tfrac{2^\theta  }{\theta} {\bf 1}_{R\geq R_{\sf out}} \psi(R_{\sf out})^{\MMM+1}
\underline{\eps}^{\frac 15}  e^{-\underline{\lambda} \tau}
,
\label{eq:tilde:V:Q:bound} 
\end{align}
for all $R>0$ and all $\tau \geq 0$.

\subsection{Closure of the bootstraps for the Taylor coefficients and modulation functions}

\begin{proposition}[\bf Bootstraps for fundamental Taylor coefficients and modulations functions]
\label{prop:Taylor:boot:closure}
There exists a sufficiently small $\underline{\eps}^*=\underline{\eps}^*(\alpha,d,\NNN) \in (0, \frac 12]$ such that the following holds. 
Assume $0<\underline{\eps} \leq \underline{\eps}^*$, and assume that the Taylor coefficients at $\tau=0$ satisfy
\begin{subequations}
\label{eq:Taylor:coeff:assumption:time:zero}
\begin{align}   
| (\tilde v_{0}(0),\tilde q_{0}(0),\tilde k_{0}(0)  )|  &\leq \tfrac 14 \underline{\eps} , 
\label{eq:Taylor:coeff:assumption:time:zero:a}\\
| (\tilde v_{\NNN}(0),\tilde q_{\NNN}(0),\tilde k_{\NNN}(0)  )|  &\leq  \underline{\eps} 
\label{eq:Taylor:coeff:assumption:time:zero:b}
.
\end{align}
\end{subequations}
Moreover, assume that the bootstrap bounds~\eqref{eq:boot:0}, \eqref{eq:boot:1}, and \eqref{eq:boot:3} hold.
Then, the bootstrap bounds~\eqref{eq:boot:0}, \eqref{eq:boot:1}, and \eqref{eq:boot:3} hold with right-side constants $\frac 12 \underline{\eps}$, $\frac 12 \underline{\eps}^{\frac 35}$, respectively $\frac 12 \underline{\eps}^{\frac 45}$.
\end{proposition}

\begin{proof}[Proof of Proposition~\ref{prop:Taylor:boot:closure}]
We start the proof by closing the bootstrap assumption~\eqref{eq:boot:3}.
In order to estimate the right sides of~\eqref{eq:modulation:final}, we need to obtain a bound for the function $\varphi^\dagger(\tau)$, which appears as the third component of~\eqref{eq:tilde:vqh:2N:alt:alt}. In turn, to achieve this we need to bound the term $\mathsf{Nonlinear}_{\eqref{eq:euler:tilde:vqh:2N}}$ defined in~\eqref{eq:euler:tilde:vqh:2N:d}, and which may be written explicitly by appealing to~\eqref{eq:N:Q:def}--\eqref{eq:N:V:def}, \eqref{eq:Taylor:R=0:all}--\eqref{eq:Taylor:R=0:tilde}, and~\eqref{eq:euler:tilde:vqh:2N:d}. Using the bootstrap assumptions~\eqref{eq:boot:1}, we obtain
\begin{equation*}
|\mathsf{Nonlinear}_{\eqref{eq:tilde:vqh:2N}}(\tau)|
\les_{\alpha,d,\NNN} 
\begin{pmatrix}
|\tilde q_{\NNN}\tilde v_0| +  |\tilde v_{\NNN}\tilde q_0| 
+ \underline{\eps}^{\frac 45} e^{- \underline{\lambda} \tau} \bigl(|\tilde q_{\NNN}| + |\tilde v_{\NNN}| + |\tilde v_0|   + |\tilde q_0| \bigr) 
\\
|\tilde v_{\NNN}\tilde v_0| +  |\tilde q_{\NNN}\tilde q_0| +  |\tilde k_{\NNN}\tilde q_0| + |\tilde q_0|^2
+ \underline{\eps}^{\frac 45} e^{- \underline{\lambda} \tau} \bigl(|\tilde q_{\NNN}| + |\tilde v_{\NNN}| + |\tilde v_0|   + |\tilde q_0| \bigr) 
\\
|\tilde{v}_{\NNN} \tilde k_0| + |\tilde{k}_{\NNN}  \tilde{v}_0| + |\tilde v_0 \tilde k_0|
+ \underline{\eps}^{\frac 45} e^{- \underline{\lambda} \tau } (|\tilde v_0| + |\tilde{v}_{\NNN}| + |\tilde k_0| + |\tilde{k}_{\NNN}| )      
\end{pmatrix}
\end{equation*}
which together with~\eqref{eq:boot:0}--\eqref{eq:boot:1} yields
\begin{align}
\label{eq:modulation:bootstrap:nonlinear:1}
|\mathsf{Nonlinear}_{\eqref{eq:tilde:vqh:2N}}(\tau)|
&\leq C_{\eqref{eq:modulation:bootstrap:nonlinear:1}}
\underline{\eps}^{\frac 75}
e^{-  2 \underline{\lambda} \tau }
.
\end{align}
In a similar fashion, we deduce from~\eqref{eq:boot:0} and~\eqref{eq:boot:1} that the matrix $\tilde{A}_2^{(\NNN)}(\tau)$ defined in~\eqref{eq:euler:tilde:vqh:2N:c} satisfies
\begin{equation}
|\tilde{A}_2^{(\NNN)}(\tau)|
\les_{\alpha,d,\NNN}  
\underline{\eps}  e^{- \underline{\lambda} \tau}
,
\label{eq:modulation:bootstrap:nonlinear:2}
\end{equation}
and the matrix 
$P$ defined in~\eqref{eq:A1:eigensystem} satiesfies
\begin{equation}
 |P| + |P^{-1}| \les_{\alpha,d,\NNN} 1 .
 \label{eq:modulation:bootstrap:nonlinear:3}
\end{equation}
Returning to~\eqref{eq:tilde:vqh:2N:alt:alt}, we deduce from~\eqref{eq:modulation:bootstrap:nonlinear:1}--\eqref{eq:modulation:bootstrap:nonlinear:3}, \eqref{eq:varphi:sharp:flat:def}, and the Gr\"onwall inequality, and the initial data assumption~\eqref{eq:Taylor:coeff:assumption:time:zero:b}, that 
\begin{align}
|\varphi^\dagger(\tau)|
&\leq  
C_{\eqref{eq:modulation:bootstrap:nonlinear:4}}
\bigl|\bigl(\tilde q_{\NNN}(0),\tilde v_{\NNN}(0),\tilde k_{\NNN}(0) \bigr)\bigr| e^{-\lambda_{\NNN}^\sharp \tau}
+
C_{\eqref{eq:modulation:bootstrap:nonlinear:4}}
\int_0^\tau e^{-\lambda_{\NNN}^\sharp (\tau-s)}  
\bigl( \underline{\eps} e^{- \underline{\lambda} s}  +
 \underline{\eps}^{\frac 75}
e^{- 2 \underline{\lambda} s }
\bigr)
ds
\notag\\
& \leq 
C_{\eqref{eq:modulation:bootstrap:nonlinear:4}}
\underline{\eps} e^{-\lambda_{\NNN}^\sharp \tau}
+
C_{\eqref{eq:modulation:bootstrap:nonlinear:4}} \tfrac{2}{ |\lambda_{\NNN}^\sharp - \underline{\lambda}|} \underline{\eps} \bigl(e^{-\lambda_{\NNN}^\sharp \tau} +  e^{-\underline{\lambda} \tau} \bigr)
,
\label{eq:modulation:bootstrap:nonlinear:4:ancient}
\end{align}
for some constant $C_{\eqref{eq:modulation:bootstrap:nonlinear:4}}>0$, which only depends on $\alpha$, $d$, and $\NNN$. 

At this stage, we recall cf.~\eqref{eq:A1:eigensystem} that $\lambda_{\NNN}^\sharp = 2\NNN (\cxbar + \bar v_0)$ and $\lambda_{\NNN}^\flat = 4\NNN (\cxbar + \bar v_0) + 1 + 2 \bar v_0$, which together with the definition of $\cxbar = \cxstar(d,\gamma,\NNN)$ in~\eqref{eq:cx:admissible},  the second bullet in Lemma~\ref{lem:properties:exponents}, and the standing assumption that $\alpha \leq d$, yields 
\begin{equation}
 \lambda_{\NNN}^\sharp
 \geq \tfrac{2}{1+\alpha d} \sqrt{\tfrac{\alpha \gamma d}{2}}
 \geq \tfrac{\sqrt{2}\alpha}{1+\alpha d}
 ,
 \qquad
  \lambda_{\NNN}^\flat
  \geq \tfrac{4}{1+\alpha d} \sqrt{\tfrac{\alpha \gamma d}{2}} 
  \geq \tfrac{2 \sqrt{2}  \alpha}{1+\alpha d}
  \label{eq:eigs:lower:bounds}
\end{equation}
for all $\NNN\geq 1$.  
Thus, from the constraint~\eqref{eq:underline:lambda:1}, namely $\underline{\lambda} \leq \frac{\alpha}{1+\alpha d}$, we obtain $\lambda_{\NNN}^\sharp, \lambda_{\NNN}^\flat > \underline{\lambda}$, which justifies the estimate of the integral in \eqref{eq:modulation:bootstrap:nonlinear:4:ancient}.
Combining~\eqref{eq:eigs:lower:bounds} with~\eqref{eq:modulation:bootstrap:nonlinear:4:ancient}, and the fact that $\underline{\lambda} \leq \frac{\alpha}{1+\alpha d}$ (cf.~\eqref{eq:underline:lambda:1}), we thus obtain
\begin{align}
|\varphi^\dagger(\tau)|
& \leq 
C_{\eqref{eq:modulation:bootstrap:nonlinear:4}}
 \tfrac{41(1+\alpha d)}{\alpha }   \underline{\eps}
e^{- \underline{\lambda} \tau}
.
\label{eq:modulation:bootstrap:nonlinear:4}
\end{align}
Similarly to~\eqref{eq:modulation:bootstrap:nonlinear:4}, we have that the first two components of the solution of~\eqref{eq:tilde:vqh:2N:alt} satisfy
\begin{align}
|\varphi^\sharp(\tau)|
&\leq 
C_{\eqref{eq:modulation:bootstrap:nonlinear:4}}
\underline{\eps}
e^{-\lambda_{\NNN}^\sharp \tau}
+
\tfrac{2 C_{\eqref{eq:modulation:bootstrap:nonlinear:4}} }{|\lambda_{\NNN}^\sharp - \underline{\lambda} |} \underline{\eps}
\bigl(e^{-\lambda_{\NNN}^\sharp \tau} +  e^{- \underline{\lambda} \tau } \bigr)
\leq 
C_{\eqref{eq:modulation:bootstrap:nonlinear:4}}
\tfrac{41(1+\alpha d)}{\alpha }  \underline{\eps}
e^{- \underline{\lambda} \tau}
,
\label{eq:modulation:bootstrap:nonlinear:4b}
\\
|\varphi^\flat(\tau)|
&\leq 
C_{\eqref{eq:modulation:bootstrap:nonlinear:4}}
\underline{\eps}
e^{-\lambda_{\NNN}^\flat \tau}
+
\tfrac{2C_{\eqref{eq:modulation:bootstrap:nonlinear:4}}}{|\lambda_{\NNN}^\flat - \underline{\lambda}|} \underline{\eps}
\bigl(e^{-\lambda_{\NNN}^\flat \tau} +  e^{- \underline{\lambda} \tau } \bigr)
\leq 
C_{\eqref{eq:modulation:bootstrap:nonlinear:4}}
\tfrac{41(1+\alpha d)}{\alpha } \underline{\eps}
e^{- \underline{\lambda} \tau}
.
\label{eq:modulation:bootstrap:nonlinear:4c}
\end{align}
With~\eqref{eq:boot:0} and~\eqref{eq:modulation:bootstrap:nonlinear:4} we return to~\eqref{eq:modulation:2:final} and deduce that
\begin{align}
\label{eq:cxtilde:close}
|\cxtilde (\tau)|
&\leq  \bigl( 1
+ C_{\eqref{eq:modulation:bootstrap:nonlinear:4}} |\cxbar + \bar v_0|  
\tfrac{41(1+\alpha d)}{\alpha }  \bigr) \underline{\eps} e^{ - \underline{\lambda} \tau } 
.
\end{align}
In a similar fashion,~\eqref{eq:modulation:1:final} yields
\begin{align}
\label{eq:cutilde:close}
|\cutilde (\tau)|
&\leq \bigl( 2 + 3 \alpha d 
+ C_{\eqref{eq:modulation:bootstrap:nonlinear:4}} |\cxbar + \bar v_0| \tfrac{41(1+\alpha d)}{\alpha } \bigr) \underline{\eps} e^{- \underline{\lambda} \tau } 
.
\end{align}
Lastly, from~\eqref{eq:modulation:3:final} we obtain that 
\begin{align}
\label{eq:cbtilde:close}
| \cbtilde(\tau) | 
&\leq  \bigl(
\tfrac{2\alpha d}{1+\alpha d}  
+ C_{\eqref{eq:modulation:bootstrap:nonlinear:4}} |\cxbar + \bar v_0|  
\tfrac{41(1+\alpha d)}{\alpha} \bigr) \underline{\eps}  e^{- \underline{\lambda} \tau} .
\end{align}
Putting together the bounds~\eqref{eq:cxtilde:close}, \eqref{eq:cutilde:close}, \eqref{eq:cbtilde:close}, we deduce that if $\underline{\eps}$ is chosen to be sufficiently small to ensure that 
\begin{equation*}
3 + 3 \alpha d 
+ \tfrac{2\alpha d}{1+\alpha d}  
+ 3 C_{\eqref{eq:modulation:bootstrap:nonlinear:4}} |\cxbar + \bar v_0|  
\tfrac{41(1+\alpha d)}{\alpha }  
 \leq \tfrac 12 \underline{\eps}^{-\frac 15}
 ,
\end{equation*}
 then the bootstrap~\eqref{eq:boot:3} is closed.

Next, we turn to the bootstrap assumption~\eqref{eq:boot:0}, concerning the coefficients of $R^0$ in~\eqref{eq:Taylor:R=0:tilde}. We note that~\eqref{eq:euler:tilde:vq:hq:zero:total} implies
\begin{equation*}
 \bigl| \bigl(\tilde v_{0},\tilde q_{0} ,\tilde k_{0} \bigr)(\tau)\bigr|
 \leq
 \bigl| \bigl(\tilde v_{0},\tilde q_{0} ,\tilde k_{0} \bigr)(0)\bigr| 
 e^{- \frac{2\alpha d \tau}{1+\alpha d}} 
 + \int_0^\tau e^{-\frac{2\alpha d (\tau-s) }{1+\alpha d}} \bigl| \bigl(\mathcal{N}_{\tilde V},\mathcal{N}_{\tilde Q} ,0 \bigr)(0,s)\bigr| ds
 .
\end{equation*}
Restricting~\eqref{eq:N:Q:def}--\eqref{eq:N:V:def} to $R=0$, appealing to the bootstrap in~\eqref{eq:boot:0}, the definitions~\eqref{eq:modulation:2:final}--\eqref{eq:modulation:3:final}, and assuming that $\underline{\eps}$ is taken to be sufficiently small in terms of $\alpha$ and $d$, we deduce
\begin{align*}
\bigl| \bigl(\mathcal{N}_{\tilde V},\mathcal{N}_{\tilde Q} , 0 \bigr)(0,s)\bigr| 
\leq 
\underline{\eps}^{\frac 85}  e^{- 2 \underline{\lambda} s}
.
\end{align*}
Combining the two above estimates,   appealing to the initial data assumption in~\eqref{eq:Taylor:coeff:assumption:time:zero:a}, and to the $\underline{\lambda}$ bound in~\eqref{eq:underline:lambda:1}, we arrive at
\begin{align*}
 \bigl| \bigl(\tilde v_{0},\tilde q_{0} ,\tilde k_{0} \bigr)(\tau)\bigr|
 &\leq 
 e^{- \frac{2\alpha d \tau}{1+\alpha d}} 
  \bigl| \bigl(\tilde v_{0},\tilde q_{0} ,\tilde k_{0} \bigr)(0)\bigr| 
 + 
 \tfrac{  1  }{\frac{2\alpha d}{1+\alpha d} - 2 \underline{\lambda}} \underline{\eps}^{\frac 85} e^{-2 \underline{\lambda} \tau}
\leq 
 \bigl( \tfrac{1}{4}  + \tfrac{1+\alpha d}{ \alpha d } \underline{\eps}^{\frac 35} \bigr)
 \underline{\eps}
 e^{- \underline{\lambda} \tau}
 .
\end{align*}
Upon ensuring that  $\underline{\eps}^{\frac 35}\tfrac{1+\alpha d}{\alpha d}  \leq \frac 14 $, the above estimate closes the bootstrap~\eqref{eq:boot:0}.

At last, we turn to the bootstrap assumption~\eqref{eq:boot:1},
which concerns the coefficients of $R^{2\NNN}$. We recall from~\eqref{eq:varphi:sharp:flat:def} that the vector $(\tilde v_{\NNN},\tilde q_{\NNN},\tilde k_{\NNN})$ is obtained from the vector $(\varphi^\sharp,\varphi^\flat,\varphi^\dagger)$ upon multiplication by the matrix $P$, which satisfies the bound~\eqref{eq:modulation:bootstrap:nonlinear:3}. Thus, from~\eqref{eq:modulation:bootstrap:nonlinear:4}--\eqref{eq:modulation:bootstrap:nonlinear:4c} and the initial data assumption~\eqref{eq:Taylor:coeff:assumption:time:zero:b}, 
we deduce that 
\begin{align}
\bigl| \bigl(\tilde v_{\NNN}(\tau),\tilde q_{\NNN}(\tau),\tilde k_{\NNN}(\tau) \bigr)\bigr|
&\leq
C_{\eqref{eq:modulation:bootstrap:nonlinear:3}}
C_{\eqref{eq:modulation:bootstrap:nonlinear:4}}
\bigl( 1 + \tfrac{40(1+\alpha d)}{ \alpha} \bigr)  
\underline{\eps} e^{- \underline{\lambda} \tau}
.
\end{align}
We deduce that if $\underline{\eps}$  is chosen to be sufficiently small to ensure that 
$
C_{\eqref{eq:modulation:bootstrap:nonlinear:3}}
C_{\eqref{eq:modulation:bootstrap:nonlinear:4}}
\bigl(1 + \tfrac{40(1+\alpha d)}{ \alpha} \bigr)
 \leq \tfrac 12 \underline{\eps}^{-\frac 25},
$
 then the bootstrap~\eqref{eq:boot:1} is closed.
\end{proof}

\begin{proposition}[\bf Bootstraps for higher order Taylor coefficients]
\label{prop:Taylor:boot:closure:higher}
Let $\underline{\eps}^*$ be as in Proposition~\ref{prop:Taylor:boot:closure}, and let $0<\underline{\eps} \leq \underline{\eps}^*$.
Define
\begin{equation}
\underline{\lambda} = \underline{\lambda}(\gamma,d,\NNN)
:=   \min\Bigl\{ \tfrac{\alpha}{2(1+\alpha d)},  (\NNN+1) \bigl( \cxstar(d,\gamma,\NNN) - \tfrac{1}{1+ \alpha d}\bigr) , (\NNN+1) \bigl( \cxstar(d,\gamma,\NNN) -  \cxstar(d,\gamma,\NNN+1)  \bigr)
\Bigr\}
.
\label{eq:underline:lambda:def}
\end{equation}
Upon potentially reducing the value of $\underline{\eps}^*$, the following holds. Assume that the initial data at $\tau=0$ is such that for any even integer $\ell \in \{2 \NNN+1, \ldots, \MMM \}$ the Taylor coefficients at $R=0$ satisfy
\begin{align}   
\tfrac{1}{\ell!} |(\p_R^\ell \tilde Q)(0,0)| 
+
\tfrac{1}{\ell!} |(\p_R^\ell \tilde V)(0,0)| 
+
\tfrac{1}{\ell!} |(\p_R^\ell \tilde K)(0,0)| 
\leq 
\underline{\eps}^{\frac 35}
\label{eq:higher:Taylor:coeff:assumption:initial}
,
\end{align}
and assume that the bootstrap bounds~\eqref{eq:boot:0}--\eqref{eq:boot:1} and~\eqref{eq:boot:3} hold. 
Then, the bootstrap bound~\eqref{eq:boot:2} holds for all $\tau \geq 0$.
\end{proposition}
\begin{proof}[Proof of Proposition~\ref{prop:Taylor:boot:closure:higher}]
Let $\MMM = 17 \NNN$, and fix an integer $\ell \in \{ 2 \NNN + 1, \ldots, \MMM\}$. 
When $\ell$ \emph{is odd}, due to the smoothness of the initial data in original variables, we have that $(\p_R^\ell \tilde Q,\p_R^\ell \tilde V,\p_R^\ell \tilde K)(0,0)= 0$. Thus, by Lemma~\ref{lem:vanishing:at:R=0}, this vanishing is propagated forward in time for $\tau>0$. 

We are left to consider \emph{even integers} $\ell \in \{ 2 \NNN + 1, \ldots, \MMM\}$. For convenience of notation, we let $\ell = 2 \LLL$, where $\LLL \in \{ \NNN+ 1, \ldots, \lfloor \frac{\MMM}{2} \rfloor \}$, and for any sufficiently smooth function $F$ we denote
\begin{equation}
f_{\LLL}(\tau) := \tfrac{1}{(2 \LLL)!} (\p_R^{2 \LLL} F)(0,\tau)
.
\end{equation}
The above notation is  consistent with the notation used in~\eqref{eq:Taylor:R=0:all} and~\eqref{eq:Taylor:R=0:tilde}.
Our goal is to prove that there exists a constant $C_\MMM^\sharp > 0$, which only depends on $\alpha, d$, and $\NNN$, such that 
\begin{equation}
|\tilde q_{\LLL}(\tau)| + |\tilde v_{\LLL}(\tau)| + |\tilde k_{\LLL}(\tau)|
\leq \underline{\eps}^{\frac 35} (C_\MMM^\sharp)^{\LLL-\NNN} e^{-\underline{\lambda} \tau}
,
\label{eq:drac}
\end{equation}
for all $\LLL \in \{ \NNN , \ldots, \lfloor \frac{\MMM}{2} \rfloor \}$ and all $\tau \geq 0$. 
We shall prove~\eqref{eq:drac} inductively on $\LLL$. The bootstrap~\eqref{eq:boot:1} implies the bound~\eqref{eq:drac} for $\LLL = \NNN$, and this serves as the induction base. For the induction step, fix $\LLL \geq \NNN+1$, and assume that \eqref{eq:drac} holds with $\LLL$ replaced by $i$, for all $i$ such that $\NNN \leq i \leq \LLL-1$.

Next, we recall the notation for the matrices $E_0$ and $E_1$ from~\eqref{eq:E0:E1:def}. We denote by $\bar E_0$ and $\bar E_1$ the same matrices, with $(\cx,\cu,\cb)$ replaced by $(\cxbar,\cubar,\cbbar)$, and $(v_0,q_0)$ replaced by $(\bar v_0,\bar q_0)$. We also denote $\tilde E_0 = E_0 - \bar E_0$ and $\tilde E_1 = E_1 - \bar E_1$. With this notation, it follows from~\eqref{eq:coeff:evo} that for all $\LLL \in \{ \NNN+ 1, \ldots, \lfloor \frac{\MMM}{2} \rfloor \}$ we have
\begin{subequations}
\label{eq:euler:tilde:vqh:2L}
\begin{align}
&\tfrac{d}{d\tau}
(\tilde q_{\LLL} , \tilde v_{\LLL} , \tilde k_{\LLL}  )^\intercal
+ 
(E_0  + 2 \LLL  E_1  )
(\tilde q_{\LLL} , \tilde v_{\LLL} , \tilde k_{\LLL}  )^\intercal
\notag\\
&\qquad 
=
-
(\tilde E_0  + 2 \LLL  \tilde E_1  )
(\bar q_{\LLL} , \bar v_{\LLL} , \bar k_{\LLL}  )^\intercal
-   {\bf 1}_{\LLL \geq 2 \NNN}
(\tilde Z^{(\LLL)}_q,\tilde Z^{(\LLL)}_v, \tilde Z^{(\LLL)}_k)^\intercal
,
\label{eq:euler:tilde:vqh:2L:a}
\end{align}
where 
\begin{align}
\tilde Z^{(\LLL)}_q
&:=
\sum_{i=\NNN}^{\LLL -\NNN}  \bigl(\bar v_i \tilde q_{\LLL-i} + \tilde v_i \bar q_{\LLL-i} + \tilde v_i \tilde q_{\LLL-i} \bigr)  \bigl( 2 \LLL - 2 (1-\alpha) i   + (1+\alpha d)  \bigr)
,
\label{eq:euler:tilde:vqh:2L:b}
\\
\tilde Z^{(\LLL)}_v
&:= \sum_{i=\NNN}^{\LLL-\NNN} 
 (2i + 1)  \bigl(\bar v_i \tilde v_{\LLL-i} + \tilde v_i \bar v_{\LLL-i} + \tilde v_i \tilde v_{\LLL-i} \bigr)  
+ \bigl(2\alpha i + \tfrac{2\alpha^2}{\gamma}\bigr)  \bigl(\bar q_i \tilde q_{\LLL-i} + \tilde q_i \bar q_{\LLL-i} + \tilde q_i \tilde q_{\LLL-i} \bigr) 
\notag\\
&\quad 
- \tfrac{\alpha}{\gamma} \sum_{i=\NNN}^{\LLL-\NNN}  \sum_{j=0}^{\LLL-i} 
2 i \bigl(\bar k_i  \bar q_j \tilde q_{\LLL-i-j} + \bar k_i \tilde q_j \bar q_{\LLL-i-j} + \tilde k_i  \bar q_j \bar q_{\LLL-i-j} \notag\\
&\qquad \qquad \qquad \qquad\qquad\qquad
+ \bar k_i  \tilde q_j \tilde q_{\LLL-i-j} + \tilde k_i \bar q_j \tilde q_{\LLL-i-j} + \tilde k_i  \tilde q_j \bar q_{\LLL-i-j} +  \tilde k_i  \tilde q_j \tilde q_{\LLL-i-j} \bigr) 
,
\label{eq:euler:tilde:vqh:2L:c}
\\
\tilde Z^{(\LLL)}_k
&:= \sum_{i=\NNN}^{\LLL-\NNN} 2 (\LLL-i)  (\bar v_i \tilde k_{\LLL-i} + \tilde v_i \bar k_{\LLL-i} + \tilde v_i \tilde k_{\LLL-i}  ) 
.
\label{eq:euler:tilde:vqh:2L:d}
\end{align}
In deriving~\eqref{eq:euler:tilde:vqh:2L:b}--\eqref{eq:euler:tilde:vqh:2L:d} we have used that $\bar q_i =  \bar v_i =  \bar k_i = 0$ if $\frac{i}{2\NNN} \not \in \Naturals$, and the fact that Lemma~\ref{lem:vanishing:at:R=0} implies $\tilde q_i =  \tilde v_i =  \tilde k_i = 0$ for $1 \leq i \leq \NNN-1$.
\end{subequations}

Using the inductive bound~\eqref{eq:drac} and the fact that $(\bar Q,\bar V,\bar K)$ are real-analytic at $R=0$, we deduce that there exists a constant $C_{\eqref{eq:drac:1} }>0$, which only depends on $\alpha, d,$ and $\NNN$, such that for all $2\NNN  \leq \LLL \leq \lfloor \frac{\MMM}{2} \rfloor$ we have
\begin{align}
\bigl| (\tilde Z^{(\LLL)}_q,\tilde Z^{(\LLL)}_v, \tilde Z^{(\LLL)}_k) (\tau) \bigr|
\leq 
C_{\eqref{eq:drac:1} } \underline{\eps}^{\frac 35}  (C_\MMM^\sharp)^{\LLL - 2\NNN} e^{-\underline{\lambda} \tau} 
,
\label{eq:drac:1} 
\end{align}
for all $\tau\geq0$. We emphasize that $C_{\eqref{eq:drac:1} }$ is independent of $C_\MMM^\sharp$.

Next, we appeal to the bootstrap assumptions~\eqref{eq:boot:0} and~\eqref{eq:boot:3}, and to the definitions of the matrices $\tilde E_0$ and $\tilde E_1$ (cf.~\eqref{eq:E0:E1:def}), to deduce that 
\begin{equation}
| \tilde E_0 + 2 \LLL \tilde E_1 |
\leq 
C_{\eqref{eq:drac:2}}
\underline{\eps}^{\frac 45} e^{-\underline{\lambda} \tau}
 \label{eq:drac:2} 
\end{equation}
where $C_{\eqref{eq:drac:2}}>0$ is a constant that only depends on $\alpha$, $d$, and $\NNN$.

The last, and most important, ingredient of the proof concerns the eigenvalues of the matrix 
\begin{equation*}
\bar E_0 + 2 \LLL \bar E_1
:=
\begin{pmatrix}
0
&  (1+\alpha d)  \bar q_0  
& 0\\
\frac{4\alpha^2}{\gamma} \bar q_0  
& 1 + 2 \bar v_0  
& 0 \\
0
& 1
& 0
\end{pmatrix}
+2 \LLL
\begin{pmatrix}
\cxbar  +  \bar v_0
&\alpha  \bar  q_0 
& 0\\
\alpha  \bar  q_0 
&\cxbar +  \bar  v_0
& -\frac{\alpha}{\gamma}  \bar q_0^2 \\
0
&0
& \cxbar  +  \bar  v_0
\end{pmatrix}
.
\end{equation*}
Observe that $\bar E_0 + 2 \LLL \bar E_1 = A_1^{(\LLL)}$, cf.~\eqref{eq:euler:tilde:vqh:2N:b}; thus in analogy to~\eqref{eq:A1:eigensystem}, and recalling the notation in~\eqref{eq:cx:admissible},  
we may show that the eigenvalues of $\bar E_0 + 2 \LLL \bar E_1$ are given by 
\begin{equation}
\label{eq:dummy:equation:0}
2 \LLL(\cxbar + \bar v_0) ,
\quad
2\LLL\bigl( \cxbar -  \cxstar(d,\gamma,\LLL) \bigr) \, ,
\quad
2 \LLL \bigl( \cxbar -  \cxstar(d,\gamma,\LLL) \bigr) + \tfrac{4 \LLL}{1+\alpha d}  \sqrt{\tfrac{\alpha \gamma d}{2} +  \mathsf{E}_{\LLL}+ \tfrac{(1-\alpha d)^2}{16 \LLL^2 }} 
.
\end{equation}
By \emph{definition} we have $\cxbar = \cxstar(d,\gamma,\NNN)$, and recalling the second item in Lemma~\ref{lem:properties:exponents}, 
we have that 
\[
\cxbar -  \cxstar(d,\gamma,\LLL) =  \cxstar(d,\gamma,\NNN) -  \cxstar(d,\gamma,\LLL) \geq \cxstar(d,\gamma,\NNN) -  \cxstar(d,\gamma,\NNN+1) > 0,
\] 
whenever $\LLL \geq \NNN+1$. Thus, for all $\LLL \geq \NNN+1$ we have shown that all three eigenvalues of the matrix $\bar E_0 + 2 \LLL \bar E_1$ are larger than or equal to
\begin{equation*}
2 (\NNN+1) \min\bigl\{ \cxstar(d,\gamma,\NNN) - \tfrac{1}{1+ \alpha d} , \cxstar(d,\gamma,\NNN) -  \cxstar(d,\gamma,\NNN+1) 
\bigr\}
\geq  2 \underline{\lambda}
.
\end{equation*}
Here we have appealed to the definition of $\underline{\lambda}$ in~\eqref{eq:underline:lambda:def}. Since the eigenvalues of a matrix are continuous with respect to small perturbations in the matrix (the amplification factor is the related to the condition number of the original matrix, via the Bauer-Fike Theorem),  the bound~\eqref{eq:drac:2} implies that if $  \underline{\eps}^{\frac 45}$ is sufficiently small with respect to $\alpha$, $d$, and $\NNN$, all the eigenvalues of the matrix $  E_0 + 2 \LLL   E_1$ are larger than or equal to $\frac 32 \underline{\lambda}$.

Combining the above spectral information for the matrix $E_0 + 2\LLL E_1$ with the bounds~\eqref{eq:drac:1}, \eqref{eq:drac:2}, the assumption~\eqref{eq:higher:Taylor:coeff:assumption:initial}, and the evolution equation~\eqref{eq:euler:tilde:vqh:2L:a}, we deduce that 
\begin{align}
\bigl| (\tilde q_{\LLL} , \tilde v_{\LLL} , \tilde k_{\LLL}  )(\tau) \bigr|
&\leq 
\bigl| (\tilde q_{\LLL} , \tilde v_{\LLL} , \tilde k_{\LLL}  )(0) \bigr| 
e^{- \frac{3}{2} \underline{\lambda} \tau} 
\notag\\
&\quad + 
\int_0^\tau
e^{- \frac{3}{2} \underline{\lambda} (\tau-s)}
\Bigl( C_{\eqref{eq:drac:2}} |(\bar q_{\LLL}, \bar v_{\LLL}, \bar k_{\LLL})|
\underline{\eps}^{\frac 45}   e^{-\underline{\lambda} s} 
+   {\bf 1}_{\LLL \geq 2 \NNN}
C_{\eqref{eq:drac:1} } \underline{\eps}^{\frac 35}  (C_\MMM^\sharp)^{\LLL - 2\NNN} e^{-\underline{\lambda} s} 
\Bigr) 
ds
\notag\\
&\leq 
\eps_{\sf \NNN}
e^{- \frac{3}{2} \underline{\lambda} \tau} 
+ 
\Bigl( C_{\eqref{eq:drac:2}} |(\bar q_{\LLL}, \bar v_{\LLL}, \bar k_{\LLL})|
 \underline{\eps}^{\frac 45}  
+   {\bf 1}_{\LLL \geq 2 \NNN}
C_{\eqref{eq:drac:1} } \underline{\eps}^{\frac 35}  (C_\MMM^\sharp)^{\LLL - 2\NNN} 
\Bigr) \tfrac{2}{\underline{\lambda}} e^{-\underline{\lambda} \tau} 
.
\end{align}
Thus, the induction bound~\eqref{eq:drac} is established if we ensure that 
\begin{equation}
1
+ 
 \tfrac{2}{\underline{\lambda}}   C_{\eqref{eq:drac:2}} |(\bar q_{\LLL}, \bar v_{\LLL}, \bar k_{\LLL})|  \underline{\eps}^{\frac 15} 
+  \tfrac{2}{\underline{\lambda}}   {\bf 1}_{\LLL \geq 2 \NNN}
C_{\eqref{eq:drac:1} }    (C_\MMM^\sharp)^{\LLL - 2\NNN} 
\leq 
(C_\MMM^\sharp)^{\LLL-\NNN} 
\label{eq:drac:3}
.
\end{equation}
Using the fact that $ \max_{\NNN \leq \LLL \leq \MMM/2} |(\bar q_{\LLL}, \bar v_{\LLL}, \bar k_{\LLL})|$ is bounded solely in terms of $d$, $\alpha$, and $\NNN$, and recalling definition~\eqref{eq:underline:lambda:def}, we obtain that there exists a constant $C_{\eqref{eq:drac:4}}\geq 1$, which only depends on $d, \alpha$, and $\NNN$, such that 
\begin{equation}
{\sf LHS}_{\eqref{eq:drac:3}} \leq  C_{\eqref{eq:drac:4}} + {\bf 1}_{\LLL \geq 2 \NNN} C_{\eqref{eq:drac:4}} (C_{\MMM}^\sharp)^{\LLL - 2 \NNN} 
.
\label{eq:drac:4}
\end{equation}
Since $\LLL \geq \NNN+1$, it follows from~\eqref{eq:drac:4} that the inequality~\eqref{eq:drac:3} holds as soon as we let $C_{\MMM}^\sharp = 2 C_{\eqref{eq:drac:4}} \geq 2$. This concludes the proof of the inductive bound~\eqref{eq:drac}.
 
The proof  is completed by noting that if $\underline{\eps}$ is sufficiently small, then $(C_{\MMM}^\sharp)^{\lfloor \MMM/2 \rfloor - \NNN} < \underline{\eps}^{-\frac 15}$.
\end{proof}

\subsection{Closure of the derivative bootstrap}
In order to close the bootstrap~\eqref{eq:boot:4} we first need to deduce an evolution equation for $\JJJ_{\MMM} \Ytilde$. We apply the operator $\JJJ_{\MMM}$ to the equations~\eqref{eq:Ytilde:evo} and then multiply the resulting PDE by $\brak{R}^\theta$; using~\eqref{eq:J:M:properties:b}--\eqref{eq:J:M:properties:d} we deduce  
\begin{align}
& \p_\tau ( \brak{R}^{\theta} \JJJ_{\MMM} \Ytilde )
+ \TT \, R \p_R  (\brak{R}^{\theta}  \JJJ_{\MMM}\Ytilde) 
+ \bigl( (\MMM+1) \tfrac{R \p_R \psi}{\psi} - \theta \tfrac{R^2}{\brak{R}^2} \bigr) \TT \, (\brak{R}^{\theta}  \JJJ_{\MMM}\Ytilde) 
+ \DD \, (\brak{R}^{\theta} \JJJ_{\MMM} \Ytilde)
\notag\\
&\quad 
+ \brak{R}^{\theta} \JJJ_{\MMM} \Ttilde \,  (R\p_R \Ybar)
+ \brak{R}^{\theta} \JJJ_{\MMM}\Dtilde \,   \Ybar 
+   \NN(\Ybar,\brak{R}^{\theta} \JJJ_{\MMM} \Ytilde)
+   \NN(\brak{R}^{\theta}\JJJ_{\MMM} \Ytilde,\Ybar)
\notag\\
& = \brak{R}^\theta \Bigl( 
- \JJJ_{\MMM} \TT \, \III_{\MMM} (R \p_R \Ytilde)
- \JJJ_{\MMM} \DD \, \III_{\MMM} \Ytilde
- \III_{\MMM} \Ttilde \, \JJJ_{\MMM} (R\p_R \Ybar)
- \III_{\MMM} \Dtilde \, \JJJ_{\MMM} \Ybar
\notag\\
&\qquad \qquad
- \TT \III_{\MMM}^\flat \Ytilde
- \III_{\MMM}^\sharp(\TT, R\p_R \Ytilde) 
- \III_{\MMM}^\sharp(\DD,\Ytilde)
- \III_{\MMM}^\sharp(\Ttilde,R\p_R \Ybar)
- \III_{\MMM}^\sharp(\Dtilde,\Ybar)
\notag\\
&\qquad \qquad
-   \NN(\JJJ_\MMM \YY,\III_{\MMM} \Ytilde)
-   \NN(\III_{\MMM} \Ytilde,\JJJ_{\MMM} \Ybar)
-   \NN \III_{\MMM}^{\sharp}(\YY,\Ytilde)
-   \NN\III_{\MMM}^\sharp(\Ytilde,\Ybar)
-   \NN(\JJJ_{\MMM} \Ytilde,\Ytilde)  \Bigr)
.
\label{eq:JM:tildeY:evo}
\end{align}
The remainder of this section is devoted to obtaining global-in-time stability estimates for the solution $\JJJ_{\MMM}\Ytilde$ of \eqref{eq:JM:tildeY:evo}; the terms on the right side of~\eqref{eq:JM:tildeY:evo} will be treated as error terms which decay exponentially in time, while for the terms on the left side of~\eqref{eq:JM:tildeY:evo} we need to quantify the amount of helpful damping, by carefully designing the weight $\psi$, and choosing $\MMM$ suitably. The main ingredient for treating the left side of~\eqref{eq:JM:tildeY:evo} is the following result:\footnote{The motivation for considering the quantity $\DD_{ii} - \sum_{j\neq i} |\DD_{ij}|$ in \eqref{eq:diagonally:dominated} is provided by the classical Gershgorin circle theorem.}

\begin{proposition}[\bf 
Global outgoing condition 
\eqref{eq:Omega:invariant:c} and $\MMM$ sufficiently large give damping for $\Ytilde$]
\label{prop:choice:of:M:psi}
Let $\MMM = 17 \NNN$ and let $\theta = \frac 12 \min\{1, \frac{1}{\cxbar} \}$. Assume that $\underline{\eps}$ is sufficiently small, with respect to $\gamma,d$, and $\NNN$. There exists $0 < R_{\sf in} \ll 1 \ll R_{\sf out}$, and a smooth non-decreasing weight $\psi$ satisfying~\eqref{eq:psi:basic:properties}, with $R_{\sf in}, R_{\sf out}$, and $\psi$  only depending on $\gamma,d$, and $\NNN$ (also through the profiles $\bar V$ and $\bar Q$), such that  for each fixed $i \in \{1,2,3\}$, we have
\begin{align}
 &\left(  (\MMM+1) \frac{R \p_R \psi(R)}{\psi(R)}  
 - \theta \frac{R^2}{\brak{R}^2} \right) \TT_{ii}(R,\cdot)
 + \DD_{ii}(R,\cdot) 
- \sum_{j\in\{1,2,3\}, j \neq i} |\DD_{ij}(R,\cdot)|  
 \notag\\
 &\quad 
 -   \frac{ | \NN_i(\Ybar, \brak{R}^\theta \JJJ_{\MMM} \Ytilde)) + \NN_i(\brak{R}^\theta \JJJ_{\MMM}\Ytilde, \Ybar)|(R,\cdot)}{\|\brak{R}^\theta \JJJ_{\MMM} \Ytilde(R,\cdot)\|_{L^\infty}}   
 \notag\\
 &\quad 
- \sum_{j=1}^3 | R\p_R \Ybar_j(R)| \, \frac{\brak{R}^\theta|\JJJ_{\MMM} \Ttilde_{ij}(R,\cdot)|}{\|\brak{R}^\theta\JJJ_{\MMM} \Ytilde(R,\cdot)\|_{L^\infty}}  +  |\Ybar_j(R)| \frac{\brak{R}^\theta |\JJJ_{\MMM}\Dtilde_{ij}(R,\cdot)|}{\|\brak{R}^\theta\JJJ_{\MMM} \Ytilde(R,\cdot)\|_{L^\infty}}   
 \notag\\
 &\geq \frac{\alpha + \frac 14}{1+\alpha d} 
 - \frac{C_{\eqref{eq:diagonally:dominated}} }{\|\brak{R}^\theta\JJJ_{\MMM} \Ytilde(R,\cdot)\|_{L^\infty}}   \underline{\eps}^{\frac 25}  
 e^{-\underline{\lambda} \tau }
  - C_{\eqref{eq:diagonally:dominated}}  \underline{\eps}^{\frac 15}  e^{- \underline{\lambda} \tau}
, 
 \label{eq:diagonally:dominated}
\end{align} 
for all $R\geq 0$ and all $\tau \geq0$, where $C_{\eqref{eq:diagonally:dominated}}>0$ is a constant that only depends on $\alpha, d$, and $\NNN$.
\end{proposition}

\begin{remark}[\bf Global outgoing condition \eqref{eq:Omega:invariant:c} and $\MMM$ sufficiently large]
The key stability mechanism is the global outgoing condition \eqref{eq:Omega:invariant:c}, which ensures that the transport matrix satisfies $\TT_{ii}(R) \geq \frac{1}{6 \NNN}$ for any $R \geq 0$.  By taking $\MMM$ sufficiently large and subtracting the $\MMM$-th order Taylor polynomials, we obtain stability for the perturbation $  \JJJ_{\MMM}\Ytilde$ near the origin. The outgoing condition allows us to propagate the local stability near the origin to the \emph{entire domain}. The non-decreasing weight $\psi$ plays a crucial role in capturing this stability mechanism.
\end{remark}

\begin{proof}[Proof of Proposition~\ref{prop:choice:of:M:psi}]
The proof consists of four steps, which consider the cases of $R=0$, $R \leq 2 R_{\sf in}$, $R \geq  R_{\sf out}$, and lastly $2 R_{\sf in} \leq R \leq R_{\sf out}$.

\noindent\emph{ \underline{Step 1}: $R = 0$ and the choice of $\MMM$}.  
We start the proof by first considering the left side of~\eqref{eq:diagonally:dominated} at $R=0$, which enforces the only condition on the value of $\MMM$. By definition (recall~\eqref{eq:ringW:ringZ:ringA:def}--\eqref{eq:YY:def} and~\eqref{eq:ringW:ringZ:ringA:bar}), we have that $\YY(0,\tau) = \Ybar(0,\tau) = (R\p_R \Ybar)(0,\tau) = 0$. Thus, recalling the definition of the bilinear form $\NN_i$ from~\eqref{eq:ringW:ringZ:ringA:bar}, and recalling that $\psi(R) = R$ in a vicinity of $R=0$, we obtain that 
\begin{equation*}
{\sf LHS}_{\eqref{eq:diagonally:dominated}}(0,\tau)
- \tfrac{\alpha + \frac 14}{1+\alpha d}
= (\MMM+1) \TT_{ii}(0,\tau) + 2\alpha \bar v_0 + \DD_{ii}(0,\tau) - \sum_{j \in \{1,2,3\}, j\neq i} |\DD_{ij}(0,\tau)|
.
\end{equation*}
Further, appealing to~\eqref{eq:values:at:zero:and:cb}, \eqref{eq:YY:evo:aux}, and the bootstrap bounds~\eqref{eq:boot:0}, and \eqref{eq:boot:3}, we arrive at
\begin{align}
{\sf LHS}_{\eqref{eq:diagonally:dominated}}(0,\tau)
- \tfrac{\alpha + \frac 14}{1+\alpha d}
&\geq (\MMM+1) \Bigl(\cxbar + \bar v_0 + (i-2) \alpha \bar q_0
 \Bigr) 
 + 1  
 -  \bigl( 2 d  \MMM  + 10 d^2 \bigr) (\underline{\eps}^{\frac 45} +\underline{\eps})
\notag\\
&\quad
+ \bar v_0  \bigl(\alpha + \tfrac 14 + {\rm D}_{ii}^{(v)} + \tfrac{|1-\alpha d|}{2} \bigr)
- \alpha  \bar q_0 \bigl( \tfrac{10 + d + 2\alpha(d-2)}{2\gamma} - {\rm D}_{ii}^{(q)}  \bigr)
.
\label{eq:damping:diag:dominated:R=0}
\end{align}

For $i=1$, using the last bullet of Lemma~\ref{lem:properties:exponents}, and using that $\MMM = 17 \NNN$, for any $\alpha \in (0,d]$  we have   
\begin{align}
\tfrac{1}{|\bar v_0|} {\sf RHS}_{\eqref{eq:damping:diag:dominated:R=0}}
&\geq  \tfrac{\MMM+1}{4 \NNN}(1+\tfrac 23 \alpha d)  + 1 + \alpha d  -  \tfrac{2\alpha + \frac 72 +\alpha d + |1-\alpha d|}{2}  
-  \sqrt{\tfrac{\alpha \gamma d}{2}} \tfrac{8 + (d+2)(\gamma+1) + 2\alpha(d-2) }{2\gamma}  
\notag\\
&\quad
-  (1+\alpha d) \bigl( 2 d  \MMM  + 10 d^2 \bigr) (\underline{\eps}^{\frac 45} +\underline{\eps})
\notag\\
&\geq  \tfrac{17 \NNN+ 1}{4 \NNN}(1+\tfrac 23 \alpha d)  - \tfrac{17}{4}(1+\tfrac 23 \alpha d)
-  (1+\alpha d) \bigl(34 d \NNN   + 10 d^2 \bigr) (\underline{\eps}^{\frac 45} +\underline{\eps})
\notag\\
&\geq
\tfrac{1}{|\bar v_0|} \Bigr( \tfrac{1}{6 \NNN} - \bigl(34 d \NNN   + 10 d^2  \bigr) (\underline{\eps}^{\frac 45} +\underline{\eps}) \Bigr)
.
\label{eq:vomitum:1}
\end{align}
Hence, ${\sf RHS}_{\eqref{eq:damping:diag:dominated:R=0}} \geq \tfrac{1}{7 \NNN}$ once $\underline{\eps}^{\frac 45} + \underline{\eps}$ is taken to be sufficiently small, solely with respect to $d$ and $\NNN$.
For $i \in \{2,3\}$, using the last bullet of Lemma~\ref{lem:properties:exponents}, and using that $\MMM = 17 \NNN$, for any $\alpha \in (0,d]$  we have
\begin{align}
\tfrac{1}{|\bar v_0|} {\sf RHS}_{\eqref{eq:damping:diag:dominated:R=0}}
&\geq  \tfrac{\MMM+1}{4 \NNN}(1+\tfrac 23 \alpha d) + (\MMM+1) \sqrt{\tfrac{\alpha \gamma d}{2}}  + 1 + \alpha d  -  \tfrac{2\alpha + \frac 52 + 3 \alpha d + |1-\alpha d|}{2}  
-  \sqrt{\tfrac{\alpha \gamma d}{2}} \tfrac{10 +  d + 2\alpha(d-2) }{2\gamma}  
\notag\\
&\quad
-  (1+\alpha d) \bigl( 2 d  \MMM  + 10 d^2 \bigr) (\underline{\eps}^{\frac 45} +\underline{\eps})
\notag\\
&\geq  \tfrac{17 \NNN+1}{4 \NNN}(1+\tfrac 23 \alpha d) - (1+\tfrac 23 \alpha d) -  (1+\alpha d) \bigl(34 d \NNN   + 10 d^2 \bigr) (\underline{\eps}^{\frac 45} +\underline{\eps}) 
\notag\\
&\geq \tfrac{1}{|\bar v_0|} \Bigr(2 + \tfrac{1}{6 \NNN} - \bigl(34 d \NNN  + 10 d^2 \bigr) (\underline{\eps}^{\frac 45} +\underline{\eps}) \Bigr) ,
\label{eq:vomitum:2}
\end{align}
which is a lower bound strictly better than the one obtained for $i=1$. The conclusion of Step 1 of the proof is that if $\underline{\eps}^{\frac 45} + \underline{\eps}$ is sufficiently small (with respect to $d$ and $\NNN$), then for any $\alpha \in (0,d]$ and $i \in \{1,2,3\}$ we have 
\begin{equation}
{\sf LHS}_{\eqref{eq:diagonally:dominated}}(0,\tau)
\geq \tfrac{\alpha + \frac 14}{1+\alpha d}
+ \tfrac{1}{7 \NNN} ,
\label{eq:damping:diag:dominated:R=0:final}
\end{equation}
an estimate which is strictly better than the desired~\eqref{eq:diagonally:dominated}. 
Obtaining the lower bound in~\eqref{eq:damping:diag:dominated:R=0:final} is the only part of the proof that imposes the requirement $\MMM = 17 \NNN$.

\noindent \emph{\underline{Step 2}: $R \leq 2 R_{\sf in}$ and the choice of $R_{\sf in}$}.  
For the terms on the second line of the left side of~\eqref{eq:diagonally:dominated}, we use~\eqref{eq:YY:evo:NN} and~\eqref{eq:ringW:ringZ:ringA:bar} to conclude that 
\begin{align}
\frac{| \NN_i(\Ybar, \brak{R}^\theta  \JJJ_{\MMM} \Ytilde)) + \NN_i(\brak{R}^\theta  \JJJ_{\MMM}\Ytilde, \Ybar)|(R,\tau)}{\|\brak{R}^\theta\JJJ_{\MMM} \Ytilde(\cdot,\tau)\|_{L^\infty}} 
\leq  6(1+\alpha) |\Ybar(R)| \leq C_{\alpha,d,\NNN} R_{\sf in}^2
,
\label{eq:vomitus:3}
\end{align}
where $C_{\alpha,d,\NNN} >0$ is a constant that only depends on the profiles $\bar V$, $\bar Q$, and $\bar H$.

For the terms on the third line of the left side of~\eqref{eq:diagonally:dominated}, we use definitions~\eqref{eq:Ytilde:evo:Ttilde} and~\eqref{eq:Ytilde:evo:Dtilde} to rewrite 
\begin{subequations}
\label{eq:vomitus:3:new}
\begin{align}
|\JJJ_{\MMM} \Ttilde_{ij}| 
&\leq {\bf 1}_{i=j} \bigl( \tfrac{1}{R_{\sf in}^{\MMM+1}} {\bf 1}_{R\geq R_{\sf in}} |\cxtilde|  +   |\JJJ_{\MMM} \tilde V|   + \alpha   |\JJJ_{\MMM} \tilde Q|  \bigr) 
,\\
|\JJJ_{\MMM} \Dtilde_{ij}| 
&\leq \tfrac{1}{R_{\sf in}^{\MMM+1}} {\bf 1}_{R\geq R_{\sf in}} |\cxtilde-\cutilde|  +  |\JJJ_{\MMM} \tilde V|   |{\rm D}^{(v)}_{ij}| 
+ \alpha  |\JJJ_{\MMM} \tilde Q|   |{\rm D}^{(q)}_{ij}|
.
\end{align}
\end{subequations}
Appealing to~\eqref{eq:boot:3}, \eqref{eq:J:M:tilde:V:Q:bound:alt}, \eqref{eq:J:M:tilde:V:Q:bound}, and the fact that $R_{\sf in} \leq \frac 14$, we deduce from~\eqref{eq:vomitus:3:new}
that 
\begin{align}
&\sum_{j=1}^3 | R\p_R \Ybar_j(R)|  \frac{\brak{R}^\theta|\JJJ_{\MMM} \Ttilde_{ij}(R,\cdot)|}{\|\brak{R}^\theta\JJJ_{\MMM} \Ytilde(R,\cdot)\|_{L^\infty}} 
+
 |\Ybar_j(R)|   \frac{\brak{R}^\theta|\JJJ_{\MMM} \Dtilde_{ij}(R,\cdot)|}{\|\brak{R}^\theta\JJJ_{\MMM} \Ytilde(R,\cdot)\|_{L^\infty}}
 \notag\\
 &
 \leq C_{\alpha,d,\NNN} R_{\sf in}^2  
 \Bigl( 1 +  
 \frac{1}{R_{\sf in}^{\MMM+1}\|\brak{R}^\theta\JJJ_{\MMM} \Ytilde(R,\cdot)\|_{L^\infty}}   \underline{\eps}^{\frac 25} e^{- \underline{\lambda} \tau }
\Bigr)
,
\label{eq:vomitus:4}
\end{align}
where $C_{\alpha,d,\NNN}>0$ is a constant that only depends on $\alpha$, $d$, and $\NNN$.

For the terms on the first line of the left side of~\eqref{eq:diagonally:dominated}, we recall that in \emph{Step 1} we have evaluated them at $R=0$, and have shown (cf.~\eqref{eq:damping:diag:dominated:R=0:final}) that the resulting expression is $\geq \tfrac{\alpha + \frac 14}{1+\alpha d}
+ \tfrac{1}{7 \NNN}$. We also recall that for $R\leq 2 R_{\sf in}$ we have that $\frac{R\p_R \psi(R)}{\psi(R)} = 1$, and we note cf.~\eqref{eq:Omega:invariant:c} that the lower bound of $\frac{1+ \frac 23 \alpha d}{4\NNN(1+\alpha d)}$ for $\Tbar_{ii}$ holds not just at $R=0$, but uniformly in $R$. Lastly, we note that using~\eqref{eq:boot:4} and~\eqref{eq:I:M:Ytilde:bound}, in analogy to~\eqref{eq:tilde:V:Q:bound} we may show that 
\begin{align}
& |V(R,\cdot) - V(0,\cdot)|  + \alpha |Q(R,\cdot) - Q(0,\cdot)| 
\notag\\
& \leq 
|\bar V(R) - \bar v_0|  + \alpha |\bar Q(R) - \bar q_0| 
+ (2+\alpha)
\int_0^R \bigl| \III_{\MMM} \Ytilde (R^\prime,\tau) + \psi(R^\prime)^{\MMM+1} \JJJ_{\MMM} \Ytilde(R^\prime,\tau) \bigr| \tfrac{dR^\prime}{R^\prime}
\notag\\
& \leq R^2 \bigl( \|\p_R^2 \bar V\|_{L^\infty(0,\frac 12)} + \alpha \|\p_R^2 \bar Q\|_{L^\infty(0,\frac 12)} 
+ \tfrac{3(2+\alpha)}{2}  C_{\eqref{eq:I:M:Ytilde:bound}} \underline{\eps}^{\frac 25} 
+ \tfrac{2+\alpha}{\MMM+1} (2 R_{\sf in})^{\MMM-1} \underline{\eps}^{\frac 15} \bigr)
\end{align}
for all $R \leq 2 R_{\sf in}$ and $\tau \geq 0$.

Using these observations, definitions~\eqref{eq:YY:evo:TT}--\eqref{eq:YY:evo:DD}, bounds~\eqref{eq:boot:3} and~\eqref{eq:barV:barQ:at:infinity}, in analogy with~\eqref{eq:vomitum:1}--\eqref{eq:vomitum:2} we may show that for any $i \in \{1,2,3\}$ and $R \leq 2 R_{\sf in}$, we have
\begin{align}
 &  (\MMM+1) \tfrac{R \p_R \psi(R)}{\psi(R)} \TT_{ii}(R,\cdot) 
 - \theta \tfrac{R^2}{\brak{R}^2} \TT_{ii}(R,\cdot)
 + \DD_{ii}(R,\cdot) 
- \sum_{j\in\{1,2,3\}, j \neq i} |\DD_{ij}(R,\cdot)|  
- \tfrac{\alpha + \frac 14}{1+\alpha d}
\notag\\
&
\geq  \tfrac{1}{7 \NNN} 
- \theta \tfrac{R^2}{\brak{R}^2} \TT_{ii}(R,\cdot)
-  \sum_{j\in\{1,2,3\}} |\DD_{ij}(R,\cdot) - \DD_{ij}(0,\cdot)|  
\geq  \tfrac{1}{7 \NNN} 
-   C_{\alpha ,d, \NNN} R_{\sf in}^2
,
\label{eq:vomitus:5}
\end{align}
where $C_{\alpha,d,\NNN}>0$ is a constant that only depends on $\alpha$, $d$, and $\NNN$. 

By combining~\eqref{eq:vomitus:3}, \eqref{eq:vomitus:4}, and~\eqref{eq:vomitus:5}, we deduce that for all $R\leq 2 R_{\sf in}$ we have the bound
\begin{equation}
{\sf LHS}_{\eqref{eq:diagonally:dominated}}
\geq  
\bigl( \tfrac{1}{7 \NNN} - C_{\eqref{eq:vomitus:6}} R_{\sf in}^2  \bigr)
+ \tfrac{\alpha + \frac 14}{1+\alpha d}
- \tfrac{C_{\eqref{eq:vomitus:6}}}{R_{\sf in}^{\MMM+1}\|\brak{R}^\theta\JJJ_{\MMM} \Ytilde(R,\cdot)\|_{L^\infty}}     \underline{\eps}^{\frac 25}  e^{- \underline{\lambda} \tau}
,
\label{eq:vomitus:6}
\end{equation}
where
$C_{\eqref{eq:vomitus:6}}>0$ is a constant that depends only on $\alpha ,d$, and $\NNN$.
Upon choosing $R_{\sf in} \ll 1$ such that 
\begin{equation}
R_{\sf in}:= (14 \NNN C_{\eqref{eq:vomitus:6}})^{- \frac 12}
\label{eq:R:in:def}
,
\end{equation}
so that $R_{\sf in}$ only depends on $\alpha, d$, and $\NNN$, we deduce that the bound 
\eqref{eq:vomitus:6} implies \eqref{eq:diagonally:dominated}, as long as 
$ C_{\eqref{eq:diagonally:dominated}} \geq C_{\eqref{eq:vomitus:6}} (14 \NNN C_{\eqref{eq:vomitus:6}})^{\frac{\MMM+1}{2}}$; in fact, we have left ourselves a $\frac{1}{14 \NNN}$ amount of ``room to spare''.

\noindent \emph{\underline{Step 3}: $R \geq  R_{\sf out}$ and the choice of $R_{\sf out}$}.  
Recall from~\eqref{eq:psi:basic:properties} that by construction we have $\p_R \psi = 0$ for $R\geq R_{\sf out}$.
Also, we recall from~\eqref{eq:YY:evo:DD} that the damping term $\DD$ contains $(\cx-\cb) {\rm Id}$, while~\eqref{eq:YY:evo:TT} shows that the transport matrix $\TT$ contains $\cx {\rm Id}$; these terms are independent of $R$, and are hence dominant as $R\to \infty$. Using this information, the fact that $\theta \leq 1$, together with the first inequality in~\eqref{eq:vomitus:3}, and the bounds~\eqref{eq:vomitus:3:new}, \eqref{eq:JM:barV:barQ:barY}, \eqref{eq:J:M:tilde:V:Q:bound:alt}, \eqref{eq:barV:barQ:at:infinity}, \eqref{eq:tilde:V:Q:bound}, and~\eqref{eq:Ybar:at:infinity}, we obtain that for any $R\geq R_{\sf out}$ and $\tau \geq 0$ it holds that 
\begin{align}
{\sf LHS}_{\eqref{eq:diagonally:dominated}} 
&\geq
 1 - \cxbar \theta - 3 \underline{\eps}^{\frac 45}
 - 2 (1+\alpha) 
C_{\rm D}^{(v,q)}
M_0 \brak{R_{\sf out}}^{-\frac{1}{\cxbar}}
 \notag\\
& 
- (1+\alpha) \Bigl(12   + \tfrac{2^{\theta+2}}{\theta} + \tfrac{M_2}{R_{\sf in}^{\MMM+1} \|\brak{R}^\theta\JJJ_{\MMM} \Ytilde(R,\cdot)\|_{L^\infty}} \underline{\eps}^{\frac 45}  e^{-\underline{\lambda} \tau} \Bigr) 
C_{\rm D}^{(v,q)} M_1 \brak{R_{\sf out}}^{-\frac{1}{\cxbar}}
 \notag\\
 &  
 -  \tfrac{(1+\alpha)  2^{\theta+2}  }{\theta} 
C_{\rm D}^{(v,q)}
\psi(R_{\sf out})^{\MMM+1} \brak{R_{\sf out}}^{-\theta}  \underline{\eps}^{\frac 15}  e^{-\underline{\lambda} \tau} 
\label{eq:vomitus:7}
, 
\end{align} 
where $M_0$, $M_1$, and $M_2$ are constants that only depend on $\alpha, d, \NNN$, and we have denoted
\begin{equation}
C_{\rm D}^{(v,q)} = C_{\rm D}^{(v,q)} (\alpha,d) := 1 + \sum_{j\in\{1,2,3\}} |{\rm D}_{ij}^{(v)}| + \sum_{j\in\{1,2,3\}}  |{\rm D}_{ij}^{(q)}|
. \label{eq:vomitus:7:b}
\end{equation}
Taking into account~\eqref{eq:YY:evo:DD:vq}, and the choice $\theta = \frac 12 \min\{ 1, \frac{1}{\cxbar} \}$, which implies $1 - \cxbar \theta \geq \frac 12$, we deduce that upon letting $R_{\sf out} = R_{\sf out}(\alpha, d,\NNN)\gg 1$ be defined by (note in particular that $R_{\sf out}$ is defined independently of the value of $\psi(R_{\sf out})$)
\begin{equation}
8(1+\alpha d) \Bigl(14 (1+\alpha) + \tfrac{(2+\alpha)2^{\theta+1}}{\theta}  \Bigr) 
C_{\rm D}^{(v,q)} \max\{M_0,M_1\} =:  \brak{R_{\sf out}}^{\frac{1}{\cxbar}}
,
\label{eq:R:out:def}
\end{equation}
we have that for all $R\geq R_{\sf out} $  and $\tau \geq 0$: 
\begin{align}
{\sf LHS}_{\eqref{eq:diagonally:dominated}} - \tfrac{\alpha + \frac 14}{1+\alpha d}
&\geq
\tfrac{1}{8(1+\alpha d)} - 3 \underline{\eps}^{\frac 45}
- \tfrac{M_2}{112 (1+\alpha d) R_{\sf in}^{\MMM+1} \|\brak{R}^\theta\JJJ_{\MMM} \Ytilde(R,\cdot)\|_{L^\infty}} \underline{\eps}^{\frac 45}  e^{-\underline{\lambda} \tau }  
 \notag\\
 &\; \;  
 -  \tfrac{1}{8(1+\alpha d)  \max\{M_0,M_1\} } \psi(R_{\sf out})^{\MMM+1}  \underline{\eps}^{\frac 15}  e^{- \underline{\lambda} \tau} 
. 
\end{align} 
Thus, if $ \underline{\eps}^{\frac 45} \leq \frac{1}{28(1+\alpha d)}$ and $C_{\eqref{eq:diagonally:dominated}} \geq \max \{ \tfrac{M_2}{112 (1+\alpha d) R_{\sf in}^{\MMM+1}}  , \tfrac{1}{8(1+\alpha d)  \max\{M_0,M_1\} } \psi(R_{\sf out})^{\MMM+1} \}$, an expression which only depends on $\alpha, d$, and $\NNN$, we have proven that \eqref{eq:diagonally:dominated} holds for all $R\geq R_{\sf out}$ and $\tau \geq 0$;  in fact, we have left ourselves a $\frac{1}{16(1+\alpha d)} $ amount of ``room to spare''.

\noindent \emph{\underline{Step 4}: $2 R_{\sf in} \leq R \leq  R_{\sf out}$ and the choice of $\psi$}.  At this stage of the proof the values $R_{\sf in}$ and $R_{\sf out}$ have been chosen (cf.~\eqref{eq:R:in:def} and~\eqref{eq:R:out:def}), depending solely on $\alpha, d$, and $\NNN$. Only the function $\psi$ remains to be chosen. 
For this purpose, we crucially use the outgoing condition~\eqref{eq:Omega:invariant:c} 
in Corollary~\ref{cor:Omega:invariant} for the transport matrix  
$\TT$ \eqref{eq:YY:evo:TT}:
\[
 \cxbar + \bar V(R) + \alpha \bar Q(R) > \cxbar + \bar V(R) > \cxbar + \bar V(R) - \alpha \bar Q(R) \geq \tfrac{1+ \frac 23 \alpha d}{4 \NNN (1+\alpha d)} \geq \tfrac{1}{6 \NNN} \quad  \Rightarrow \quad  \bar \TT_{ii}(R)  \geq \tfrac{1}{6 \NNN}
\]
for any $R \geq 0$. In particular, for a non-decreasing weight $\psi$, the coefficient 
associated with the term $\frac{ \pa_R \psi(R)}{\psi(R)}$ in \eqref{eq:diagonally:dominated} has a favorable sign for stability.

Analogously to \eqref{eq:vomitus:7}, we may use the above condition together with \eqref{eq:vomitus:3}, \eqref{eq:vomitus:3:new}, \eqref{eq:YY:evo:aux}, \eqref{eq:Ytilde:evo:aux}, \eqref{eq:JM:barV:barQ:barY}, \eqref{eq:J:M:tilde:V:Q:bound:alt}, \eqref{eq:barV:barQ:at:infinity}, \eqref{eq:tilde:V:Q:bound}, \eqref{eq:Ybar:at:infinity}, and choose $\underline{\eps}^{\frac 45}$ to be sufficiently small with respect to $\alpha, d, \NNN$, to obtain that for any $2 R_{\sf in} \leq R \leq R_{\sf out}$ and $\tau \geq 0$ 
\begin{align}
&{\sf LHS}_{\eqref{eq:diagonally:dominated}} 
-
\tfrac{\MMM+1}{7\NNN}   \tfrac{R \p_R \psi(R)}{\psi(R)}  
\notag\\
&\geq 
1  
- \cxbar \theta - 3 \underline{\eps}^{\frac 45}
- 2(1+\alpha) C_{\rm D}^{(v,q)} M_0 \brak{R}^{-\frac{1}{\cxbar}} 
 \notag\\
 &  
\quad - (1+\alpha)
\Bigl(12  + 4 \min \left\{ \brak{R}^\theta (1+ \log\tfrac{R}{2R_{\sf in}}), \bigl(\tfrac{\psi(R_{\sf out})}{\psi(R)}\bigr)^{\MMM+1} \right\}   \Bigr) 
C_{\rm D}^{(v,q)} M_1 \brak{R}^{-\frac{1}{\cxbar}} 
\notag\\
&  
\quad  -  
 \tfrac{(1+\alpha) M_2}{R_{\sf in}^{\MMM+1} \|\brak{R}^\theta\JJJ_{\MMM} \Ytilde(R,\cdot)\|_{L^\infty}}  \underline{\eps}^{\frac 45} e^{-\underline{\lambda} \tau} 
C_{\rm D}^{(v,q)} M_1 \brak{R}^{-\frac{1}{\cxbar}}
\notag\\
&  
\quad - 4 (1+\alpha)  \Bigl((\MMM+1) \tfrac{R \p_R \psi(R)}{\psi(R)} +  C_{\rm D}^{(v,q)}  \Bigr) \psi(R_{\sf out})^{\MMM+1}   \brak{R}^{-\theta}  \underline{\eps}^{\frac 15}  e^{-\underline{\lambda} \tau}
\label{eq:vomitus:8}
.
\end{align} 
Note that ${\sf RHS}_{\eqref{eq:vomitus:8}} = {\sf RHS}_{\eqref{eq:vomitus:8}}(R)$  is a continuous function of $R$. Since $\frac{1}{\theta} 2^\theta \geq 1$, and $\p_R \psi(R_{\sf out}) =0$, by inspection we may verify that as $R \to R_{\sf out}$ we have that ${\sf RHS}_{\eqref{eq:vomitus:8}}(R) \to {\sf RHS}_{\eqref{eq:vomitus:8}}(R_{\sf out}) > {\sf RHS}_{\eqref{eq:vomitus:7}}$. Recall that at the end of \emph{Step 3} we have noted that at $R = R_{\sf out}$ the bound~\eqref{eq:diagonally:dominated} holds even if we add $\frac{1}{16(1+\alpha d)}$ to ${\sf RHS}_{\eqref{eq:diagonally:dominated}}$. By continuity in $R$, this ``room to spare'' allows us to show the existence of $\delta R \in (0,1]$,  which only depends on $\alpha, d$, and $\NNN$,  such that ${\sf RHS}_{\eqref{eq:vomitus:8}}\geq {\sf RHS}_{\eqref{eq:diagonally:dominated}}$, for all $[R_{\sf out} - \delta R, R_{\sf out}]$.  Moreover, since $R\p_R \psi(R) \geq 0$, we obtain that $\eqref{eq:diagonally:dominated}$ holds for all $[R_{\sf out} - \delta R, R_{\sf out}]$.

Recall also that at the end of \emph{Step 2} we have shown that at $R = 2 R_{\sf in}$ the bound~\eqref{eq:diagonally:dominated} holds even if we add $\frac{1}{14 \NNN}$ to ${\sf RHS}_{\eqref{eq:diagonally:dominated}}$. By continuity in $R$, and using an argument similar to the one above, we may show that  this ``room to spare'' gives the existence of an $\delta R^\prime \in (0, 2 R_{\sf in})$,  which only depends on $\alpha, d$, and $\NNN$, such that \eqref{eq:diagonally:dominated} holds true on $[2 R_{\sf in},2 R_{\sf in} + \delta R^\prime]$. To avoid redundancy, we omit these details.

We are left to establish the validity of~\eqref{eq:diagonally:dominated} for $ R \in [2 R_{\sf in} + \delta R^\prime, R_{\sf out} + \delta R]$. It is in this step that the function $\psi$ is designed more carefully, based on the lower bound established in~\eqref{eq:vomitus:8}. 
We remark that all $\psi$-dependent terms in ${\sf RHS}_{\eqref{eq:vomitus:8}}$ are either multiplied by the small parameter $\underline{\epsilon}$, or they are bounded by constants independent of $\psi$. 
Hence, the weight $\psi(R)$ can be chosen almost freely for $R \in [2 R_{\sf in} + \delta R^\prime, R_{\sf out} + \delta R]$.
More precisely, on the interval $(2 R_{\sf in} + \delta R^\prime, R_{\sf out} + \delta R)$ we let $\psi$ satisfy the ODE
\begin{subequations}
\label{eq:vomitus:9}
\begin{align}
\tfrac{\MMM + 1}{7 \NNN} \tfrac{R \p_R \psi(R)}{\psi(R)} 
= C_{\alpha,d,\NNN}^\sharp  R^{-\frac{1}{\cxbar}} 
 ,
\label{eq:vomitus:9:a}
\end{align} 
where
\begin{equation}
C_{\alpha,d,\NNN}^\sharp:= (1+\alpha) 
\bigl(14 + 4 \brak{R_{\sf out}}^\theta (1 + \log \tfrac{R_{\sf out}}{2 R_{\sf in}}) \bigr)
C_{\rm D}^{(v,q)} \max\{M_0, M_1\} 
 .
 \label{eq:vomitus:9:b}
\end{equation}
\end{subequations}

Assuming that we succeed in solving~\eqref{eq:vomitus:9} for $\psi$, from~\eqref{eq:vomitus:8}, recalling that $\theta= \frac 12 \min\{1,\frac{1}{\cxbar}\}$, taking $\underline{\eps}^{\frac 45}$ to be sufficiently small, noting that $R^{-\frac{1}{\cxbar}} \leq (2 R_{\sf in})^{-\frac{1}{\cxbar}}$ for $R \geq 2 R_{\sf in}$, we deduce that 
\begin{align}
{\sf LHS}_{\eqref{eq:diagonally:dominated}} 
&\geq 
\tfrac{\alpha + \frac 14}{1+\alpha d}
-  
  \tfrac{(1+\alpha)M_2 C_{\rm D}^{(v,q)} M_1}{R_{\sf in}^{\MMM+1} \|\brak{R}^\theta\JJJ_{\MMM} \Ytilde(R,\cdot)\|_{L^\infty}}  \underline{\eps}^{\frac 45} e^{-\underline{\lambda} \tau}    
\notag\\
&  
- 4 (1+\alpha)\Bigl( 7 \NNN  C_{\alpha,d,\NNN}^\sharp    (2 R_{\sf in})^{-\frac{1}{\cxbar}} +  C_{\rm D}^{(v,q)}  \Bigr) \psi(R_{\sf out})^{\MMM+1}   \underline{\eps}^{\frac 15}  e^{-\underline{\lambda} \tau}
. 
\label{eq:vomitus:10}
\end{align} 
Hence, if the constant  $C_{\eqref{eq:diagonally:dominated}}$ is chosen to be large enough in terms of $\alpha,d$, and $\NNN$, in order to ensure that $C_{\eqref{eq:diagonally:dominated}} \geq \frac{(1+\alpha)M_2 C_{\rm D}^{(v,q)} M_1 }{R_{\sf in}^{\MMM+1} }   $ and  $C_{\eqref{eq:diagonally:dominated}} \geq  4 (1+\alpha) ( 7 \NNN  C_{\alpha,d,\NNN}^\sharp    (2 R_{\sf in})^{-\frac{1}{\cxbar}} +  C_{\rm D}^{(v,q)}   ) \psi(R_{\sf out})^{\MMM+1}$, then 
\eqref{eq:vomitus:10} shows that  \eqref{eq:diagonally:dominated} holds also on the interval $(2 R_{\sf in} + \delta R^\prime, R_{\sf out} + \delta R)$, thereby completing the proof of the Proposition.

It thus only remains to solve~\eqref{eq:vomitus:9:a} on $(2 R_{\sf in} + \delta R^\prime, R_{\sf out} + \delta R)$, and to extend the resulting solution to match $\psi(R) = R$ for $R\leq 2 R_{\sf in}$ and $\p_R \psi (R) = 0$ for $R\geq R_{\sf out}$. The solution of~\eqref{eq:vomitus:9:a} with $\psi(2 R_{\sf in} + \delta R^\prime) = 4 R_{\sf in}$ is given by 
\begin{equation}
\psi(R) 
= 4 R_{\sf in} \exp \left(  \tfrac{7 \NNN C_{\alpha,d,\NNN}^\sharp \cxbar }{\MMM + 1} \Bigl((2 R_{\sf in} + \delta R^\prime)^{-\frac{1}{\cxbar}} -  R^{- \frac{1}{\cxbar}} \Bigr) \right)
.
\label{eq:psi:nasty:def}
\end{equation}
The choice $\psi(2 R_{\sf in} + \delta R^\prime) = 4 R_{\sf in} > 2 R_{\sf in}$ allows for a smooth and monotone increasing choice $\psi \colon [2 R_{\sf in}, 2 R_{\sf in} + \delta R^\prime] \to \Reals_+$, with $\psi(2 R_{\sf in} ) = 2 R_{\sf in}$. It is clear that the function defined in~\eqref{eq:psi:nasty:def} is monotone increasing in $R$. 
Lastly, upon declaring 
\begin{equation*}
 \psi(R_{\sf out}) 
= 4 R_{\sf in} \exp \left(  \tfrac{7 \NNN C_{\alpha,d,\NNN}^\sharp \cxbar }{17 \NNN + 1} \Bigl((3 R_{\sf in})^{-\frac{1}{\cxbar}} -  (2 R_{\sf out})^{- \frac{1}{\cxbar}} \Bigr) \right)
,
\end{equation*}
we may choose $\psi \colon [R_{\sf out} - \delta R, R_{\sf out}] \to \Reals_+$ to be smooth and monotone increasing, such that $\p_R \psi(R_{\sf out}) = 0$. For $R\geq R_{\sf out}$, we declare $\psi(R) = \psi(R_{\sf out})$. It is now clear that the above constructed function $\psi$ satisfies all the properties postulated in~\eqref{eq:psi:basic:properties}.
\end{proof}

With Proposition~\ref{prop:choice:of:M:psi} in hand, we return to the evolution equation for $\brak{R}^\theta \JJJ_{\MMM} \Ytilde$ (cf.~\eqref{eq:JM:tildeY:evo}), and establish:
\begin{proposition}
\label{prop:derivative:boot:closure}
Assume that the initial data satisfies 
\begin{align}
\sup_{R>0} \brak{R}^\theta |\JJJ_{\MMM} \Ytilde(R,0)| \leq \underline{\eps}^{\frac 25} 
\label{eq:derivative:boot:closure}
\end{align} 
and that $\underline{\eps}$ is sufficiently small, with respect to $\alpha$, $d$, and $\NNN$. Then, the bootstrap~\eqref{eq:boot:4} is closed.
\end{proposition}

\begin{proof}[Proof of Proposition~\ref{prop:derivative:boot:closure}]
Appealing to~\eqref{eq:YY:evo:NN}, \eqref{eq:Ytilde:evo:aux}, \eqref{eq:J:M:properties}, \eqref{eq:boot:0}--\eqref{eq:boot:3}, \eqref{eq:boot:4}, \eqref{eq:barV:barQ:at:infinity}, \eqref{eq:JM:barV:barQ:barY}, \eqref{eq:IM:barV:barQ:barY}, \eqref{eq:I:M:Ytilde:bound}, \eqref{eq:Ytilde:bound}, \eqref{eq:I:M:tilde:V:Q:bound}, \eqref{eq:vomitus:3:new}, we deduce that 
\begin{align}
\|{\sf RHS}_{\eqref{eq:JM:tildeY:evo}}(\cdot,\tau)\|_{L^\infty}
&\leq 
 C_{\eqref{eq:JM:tildeY:evo:RHS:bound}}
\underline{\eps}^{\frac 25}  e^{-  \underline{\lambda} \tau } 
,
\label{eq:JM:tildeY:evo:RHS:bound}
\end{align}
for all $\tau \geq 0$, 
where $C_{\eqref{eq:JM:tildeY:evo:RHS:bound}} >0$ is a constant that depends only on $\alpha, d, \NNN$, also through the previously defined parameters $\MMM$, $\theta$, $M_0$, $M_1$, $M_2$, $C_{\eqref{eq:I:M:Ytilde:bound}}$, $R_{\sf in}$, $C_{\rm D}^{(v,q)}$, and $\psi(R_{\sf out})$.

From~\eqref{eq:JM:tildeY:evo}, \eqref{eq:diagonally:dominated}, \eqref{eq:JM:tildeY:evo:RHS:bound}, and recalling from~\eqref{eq:YY:evo:TT} that the matrix $\TT$ is diagonal, we obtain that 
for each $i \in \{1, 2, 3 \}$ (we emphasize that here we do not use a summation convention on repeated indices), we have 
\begin{align}
&(\p_\tau + \TT_{ii}(R,\tau) R \p_R ) \bigl | \brak{R}^\theta \JJJ_{\MMM} \Ytilde_i(R,\tau) \bigr|
\notag\\
&\leq 
 - \Bigl( \tfrac{\alpha + \frac 14}{1+\alpha d} + \sum_{j \in \{1,2,3\}, j \neq i}  |\DD_{ij}(R,\tau)| \Bigr) \bigl | \brak{R}^\theta \JJJ_{\MMM} \Ytilde_i(R,\tau) \bigr|
 + \sum_{j \in \{1,2,3\}, j \neq i} |\DD_{ij}(R,\tau)|  \bigl | \brak{R}^\theta \JJJ_{\MMM} \Ytilde_j(R,\tau) \bigr|
\notag\\
&\quad
+ \bigl( C_{\eqref{eq:diagonally:dominated}} + C_{\eqref{eq:JM:tildeY:evo:RHS:bound}} \bigr) 
  \underline{\eps}^{\frac 25}   e^{-\underline{\lambda} \tau} 
  \label{eq:bazooka}
  .
\end{align}
Denote by $\Phi_i = \Phi_i(R,\tau)$ the ODE flow associated with the transport operator  $\p_\tau + \TT_{ii}(R,\tau) R \p_R$; that is, recalling from~\eqref{eq:YY:evo:TT} that $\TT_i =  \cx +   V + (i-2)  \alpha Q$,  the flow $\Phi_i$ is defined as the solution of
\begin{subequations}
\label{eq:Phi:i:def}
\begin{align}
\tfrac{d}{d\tau} \Phi_i(R,\tau) 
&=  \Phi_i(R,\tau)  \Bigl(\cx(\tau) + V\bigl(\Phi_i(R,\tau),\tau \bigr) + (i-2) \alpha Q\bigl(\Phi_i(R,\tau),\tau \bigr) \Bigr) 
,
\\
\Phi_i(R,0) &= R 
.
\end{align}
\end{subequations}
We also denote
\begin{equation}
\mathcal{G}_i(R,\tau) := \sum_{j \in \{1,2,3\}, j \neq i}    |\DD_{ij}(\Phi_i(R,\tau),\tau)| 
,
\end{equation}
and note that~\eqref{eq:YY:evo:DD}, \eqref{eq:YY:evo:DD:vq}, \eqref{eq:barV:barQ:at:infinity}, and~\eqref{eq:J:M:tilde:V:Q:bound} imply the existence of a constant $C_{\eqref{eq:Gi:bound}}>0$, which only depends on $\alpha, d$, and $\NNN$, such that 
\begin{equation}
\label{eq:Gi:bound}
\sup_{\tau\geq0} \| \mathcal{G}_i(\cdot,\tau) \|_{L^\infty} \leq C_{\eqref{eq:Gi:bound}}
.
\end{equation}
Upon composing~\eqref{eq:bazooka} with $\Phi_i$ in space, 
appealing to the bootstrap~\eqref{eq:boot:4}, to the bound $\underline{\lambda} \leq \frac{\alpha}{1+\alpha d}$, to the Gr\"onwall inequality, and to the fact that $\mathcal{G}_i \geq 0$, we deduce that 
\begin{align}
e^{\underline{\lambda} \tau} \bigl | \brak{\Phi_i(R,\tau)}^\theta \JJJ_{\MMM} \Ytilde_i(\Phi_i(R,\tau),\tau) \bigr|
&\leq  | \brak{R}^\theta \JJJ_{\MMM} \Ytilde_i(R,0) \bigr| 
+ 4(1+\alpha d) \bigl( C_{\eqref{eq:diagonally:dominated}} + C_{\eqref{eq:JM:tildeY:evo:RHS:bound}} \bigr) 
 \underline{\eps}^{\frac 25}    
\notag\\
&\quad 
 +  \underline{\eps}^{\frac 15}  
 \int_0^\tau e^{- \frac{\tau-s}{4(1+\alpha d)} -  \int_s^\tau  \mathcal{G}_i(R,s^\prime) ds^\prime} \mathcal{G}_i(R,s)  
 ds
  .
  \label{eq:bazooka:2}
\end{align}
To conclude, we note that the bound~\eqref{eq:Gi:bound} and a simple barrier argument imply
\begin{equation}
 \int_0^\tau e^{- \frac{\tau-s}{4(1+\alpha d)} -  \int_s^\tau  \mathcal{G}_i(R,s^\prime) ds^\prime} \mathcal{G}_i(R,s)  
 ds
 \leq  1 - \tfrac{1}{1 + 8(1+\alpha d) C_{\eqref{eq:Gi:bound}}}
\end{equation}
for all $\tau \geq 0$. By combining the above bound with~\eqref{eq:bazooka:2} we conclude
\begin{align}
&e^{\underline{\lambda} \tau} \bigl | \brak{\Phi_i(R,\tau)}^\theta \JJJ_{\MMM} \Ytilde_i(\Phi_i(R,\tau),\tau) \bigr|
\notag\\
&\leq \| \brak{R}^\theta \JJJ_{\MMM} \Ytilde_i(R,0)\|_{L^\infty} 
+ 4(1+\alpha d) \bigl( C_{\eqref{eq:diagonally:dominated}} + C_{\eqref{eq:JM:tildeY:evo:RHS:bound}} \bigr) 
\underline{\eps}^{\frac 25}  
+  \underline{\eps}^{\frac 15}  \tfrac{8(1+\alpha d) C_{\eqref{eq:Gi:bound}}}{1 + 8(1+\alpha d) C_{\eqref{eq:Gi:bound}}} 
  ,
  \label{eq:bazooka:3}
\end{align}
for all $R>0$ and $\tau \geq 0$. Since the map $R \mapsto \Phi_i(R,\cdot)$  is a bijection from $\Reals_+ \to \Reals_+$, 
upon taking the supremum over $R>0$ we deduce from~\eqref{eq:bazooka:3} that
\begin{align}
e^{\underline{\lambda} \tau} \| \brak{\cdot}^\theta \JJJ_{\MMM} \Ytilde_i(\cdot,\tau) \|_{L^\infty}
&\leq \| \brak{\cdot}^\theta \JJJ_{\MMM} \Ytilde_i(0,\cdot) \|_{L^\infty}
\notag\\
&\quad 
+ 4(1+\alpha d) \bigl( C_{\eqref{eq:diagonally:dominated}} + C_{\eqref{eq:JM:tildeY:evo:RHS:bound}} \bigr) 
\underline{\eps}^{\frac 25}    
 +  \underline{\eps}^{\frac 15} \tfrac{8(1+\alpha d) C_{\eqref{eq:Gi:bound}}}{1 + 8(1+\alpha d) C_{\eqref{eq:Gi:bound}}} 
  .
    \label{eq:bazooka:4}
\end{align}
Thus, if we assume that 
\begin{equation}
\| \brak{\cdot}^\theta \JJJ_{\MMM} \Ytilde_i(0,\cdot) \|_{L^\infty} \leq \underline{\eps}^{\frac 25} \leq \tfrac{1}{4(1 + 8(1+\alpha d) C_{\eqref{eq:Gi:bound}})}  \underline{\eps}^{\frac 15},
\end{equation}
and  $\underline{\eps}$ is taken to be  sufficiently small with respect to $\alpha, d, \NNN$ to ensure that
\begin{equation}
 4(1+\alpha d) \bigl( C_{\eqref{eq:diagonally:dominated}} + C_{\eqref{eq:JM:tildeY:evo:RHS:bound}} \bigr) 
\underline{\eps}^{\frac 25}    
  \leq \tfrac{1}{4(1 + 8(1+\alpha d) C_{\eqref{eq:Gi:bound}})}  \underline{\eps}^{\frac 15}
  ,
\end{equation}
then~\eqref{eq:bazooka:4} implies
\begin{equation}
 e^{\underline{\lambda} \tau} \| \brak{\cdot}^\theta \JJJ_{\MMM} \Ytilde_i(\cdot,\tau) \|_{L^\infty}
 \leq \underline{\eps}^{\frac 15} \tfrac{\frac 12 + 8(1+\alpha d) C_{\eqref{eq:Gi:bound}}}{1 + 8(1+\alpha d) C_{\eqref{eq:Gi:bound}}} 
 < \underline{\eps}^{\frac 15} .
\end{equation}
This closes the bootstrap~\eqref{eq:boot:4}.
\end{proof}

\subsection{Sharp bounds with respect to \texorpdfstring{$R$}{R}}
\label{sec:sharp:bounds:R}
\subsubsection{Bounds for \texorpdfstring{$\Ytilde$}{Y tilde}}
We return to the $\Ytilde$ evolution in~\eqref{eq:Ytilde:evo}; upon multiplying this equation by the weight $\brak{R}^{\frac{1}{\cxbar}}$, we obtain
\begin{align}
&\p_\tau \bigl( \brak{R}^{\frac{1}{\cxbar}} \Ytilde \bigr)
+ \TT \, R \p_R \bigl(\brak{R}^{\frac{1}{\cxbar}} \Ytilde \bigr)
+  \mathcal{D}_{\sf sharp} \, \bigl( \brak{R}^{ \frac{1}{\cxbar}}  \Ytilde  \bigr)
+ \Ttilde \, \brak{R}^{\frac{1}{\cxbar}}  R\p_R \Ybar 
+ \Dtilde \, \brak{R}^{\frac{1}{\cxbar}}  \Ybar 
\notag\\
&\quad 
+   \NN(\Ybar, \brak{R}^{\frac{1}{\cxbar}} \Ytilde)
+   \NN(\brak{R}^{\frac{1}{\cxbar}} \Ytilde,\Ybar)
+   \NN(\brak{R}^{\frac{1}{\cxbar}}\Ytilde,\Ytilde) 
= 0 , 
\label{eq:airplane:0}
\end{align}
where
\begin{align*}
\mathcal{D}_{\sf sharp}
:=
\brak{R}^{-2} 
&+ \bar V(R)  \bigl( {\rm D}^{(v)} -\tfrac{1}{\cxbar}  {\rm Id} R^2 \brak{R}^{-2} \bigr) 
\notag\\
& 
+ \alpha \bar Q(R) \bigl( {\rm D}^{(q)} 
- \tfrac{1}{\cxbar} {\rm diag}(-1,0,1) R^2 \brak{R}^{-2}     \bigr) 
+ \Dtilde(R,\tau) 
- \tfrac{1}{\cxbar} \Ttilde(R,\tau)     R^2 \brak{R}^{-2} 
.
\end{align*}
Using the bounds established earlier in~\eqref{eq:barV:barQ:at:infinity}, \eqref{eq:Ybar:at:infinity}, \eqref{eq:Ytilde:bound}, and~\eqref{eq:tilde:V:Q:bound} we deduce that there exists an implicit constant  which only depends on $\alpha,d,\NNN$ (via $\MMM$ and the choices for $R_{\sf in}$ and $R_{\sf out}$), such that 
\begin{align}
|\mathcal{D}_{\sf sharp}|
\les
\brak{R}^{-\frac{1}{\cxbar}} 
+
\underline{\eps}^{\frac 25} e^{-\underline{\lambda} \tau}
,
\qquad
|\Ttilde|
+
|\Dtilde|
+ 
|\Ytilde|
\les
\brak{R}^{-\theta} \underline{\eps}^{\frac 15} e^{-\underline{\lambda} \tau}
.
\label{eq:airplane:1} 
\end{align}
In particular, using the previously established bound~\eqref{eq:Ytilde:bound}, we may estimate
\begin{align}
\bigl| \mathcal{D}_{\sf sharp} \, \bigl( \brak{R}^{ \frac{1}{\cxbar}}  \Ytilde  \bigr)\bigr|
\les  \underline{\eps}^{\frac 15} e^{-\underline{\lambda} \tau} + \underline{\eps}^{\frac 25} e^{-\underline{\lambda} \tau}\bigl| \brak{R}^{ \frac{1}{\cxbar}}  \Ytilde  \bigr|.
\label{eq:airplane:2}  
\end{align}
By combining~\eqref{eq:airplane:0}--\eqref{eq:airplane:2},  with the bound~\eqref{eq:Ybar:at:infinity} and the definition~\eqref{eq:YY:evo:NN}, and composing with the ODE flow $\Phi_i$ defined in~\eqref{eq:Phi:i:def}, we obtain
\begin{align}
\p_\tau \Bigl| \bigl( \brak{R}^{\frac{1}{\cxbar}} \Ytilde_i\bigr) \circ \Phi_i (R,\tau) \Bigr|
\leq 
C_{\eqref{eq:airplane:3}} 
\underline{\eps}^{\frac 15} e^{-\underline{\lambda} \tau}
+
C_{\eqref{eq:airplane:3}} 
 \underline{\eps}^{\frac 25} e^{-\underline{\lambda} \tau}
\max_{i \in \{1,2,3\}}
\Bigl| ( \brak{R}^{\frac{1}{\cxbar}} \Ytilde_i) \circ \Phi_i (R,\tau)\Bigr| 
,
\label{eq:airplane:3}
\end{align}
for each $i \in \{1,2,3\}$,
where $C_{\eqref{eq:airplane:3}} \geq 1$ is a constant which only depends on $\alpha,d,\NNN$ (via $\MMM$ and the choices for $R_{\sf in}$ and $R_{\sf out}$). Here we have also used that $\theta \leq \frac{1}{2\cxbar}$. It follows from~\eqref{eq:airplane:3} that for any $R>0$ and $\tau\geq0$, we have 
\begin{align}
\max_{i \in \{1,2,3\}}
\Bigl| ( \brak{R}^{\frac{1}{\cxbar}} \Ytilde_i) \circ \Phi_i (R,\tau)\Bigr| 
&\leq 
e^{C_{\eqref{eq:airplane:3}} \underline{\lambda}^{-1} \underline{\eps}^{\frac 25}}
\left( \max_{i \in \{1,2,3\}}
\Bigl| ( \brak{R}^{\frac{1}{\cxbar}} \Ytilde_i) (R,0)\Bigr| 
+ C_{\eqref{eq:airplane:3}} \underline{\lambda}^{-1} \underline{\eps}^{\frac 15}
\right).
\label{eq:airplane:4}
\end{align}
Taking a supremum over all $R>0$, using that for each $\tau\geq 0$ and all $i \in \{1,2,3\}$ the map $\Phi_i(\cdot,\tau) \colon (0,\infty) \to (0,\infty)$ is a bijection, and taking $\underline{\eps}$ to be sufficiently small, estimate~\eqref{eq:airplane:3} implies that
\begin{equation}
 \brak{R}^{\frac{1}{\cxbar}} \bigl| \Ytilde (R,\tau)\bigr| 
\leq 
2 C_{\eqref{eq:airplane:3}} \underline{\lambda}^{-1} \underline{\eps}^{\frac 15}
+ 
2  \sup_{R>0} \brak{R}^{\frac{1}{\cxbar}} \bigl| \Ytilde (R,0)\bigr| 
.
\label{eq:airplane:5}
\end{equation}
We note  that ${\sf RHS}_{\eqref{eq:airplane:5}}$ does not vanish as $\tau\to \infty$.

\subsubsection{Bounds for $(\tilde V, \tilde Q)$ and their derivatives}

From~\eqref{eq:tilde:V:tilde:Q:tilde:Y}, \eqref{eq:Ytilde:bound}, \eqref{eq:tilde:V:Q:bound}, we obtain the bound
\begin{subequations}
\label{eq:sharp:V:Q:derivative}
\begin{align}
 | \p_R \tilde V(R,\tau)  | + \alpha  | \p_R \tilde Q(R,\tau)  |
\leq C_{\eqref{eq:sharp:V:Q:derivative}} R \brak{R}^{-2-\theta} \underline{\eps}^{\frac 15} e^{-\underline{\lambda} \tau}
,
\qquad \forall R>0,
\end{align}
where we recall that $\theta = \tfrac 12 \min\{1,\tfrac{1}{\cxbar}\}$, while~\eqref{eq:airplane:5} further implies
\begin{align}
 | \p_R \tilde V(R,\tau)  | + \alpha  | \p_R \tilde Q(R,\tau)  |
\leq R^{-1 - \frac{1}{\cxbar}} \Bigr( C_{\eqref{eq:sharp:V:Q:derivative}} \underline{\eps}^{\frac 15}  + (2+\alpha) \sup_{R>0} \brak{R}^{\frac{1}{\cxbar}} \bigl| \Ytilde (R,0)\bigr|  \Bigr)
,
\qquad \forall R\geq 1,
\end{align}
for some constant $C_{\eqref{eq:sharp:V:Q:derivative}}\geq 1$ which only depends on $\alpha,d,\NNN$. 
Using~\eqref{eq:tilde:V:Q:bound} and the fundamental theorem of calculus, we may further obtain
\begin{align}
 | \tilde V(R,\tau)  | + \alpha  |  \tilde Q(R,\tau)  |
&\leq C_{\eqref{eq:sharp:V:Q:derivative}}  \brak{R}^{-\theta} \underline{\eps}^{\frac 15} e^{-\underline{\lambda} \tau}
,
\qquad \forall R>0,
 \\
 | \tilde V(R,\tau) | + \alpha  |  \tilde Q(R,\tau)  |
&\leq \cxbar R^{- \frac{1}{\cxbar}} \Bigr( C_{\eqref{eq:sharp:V:Q:derivative}} \underline{\eps}^{\frac 15}  + (2+\alpha) \sup_{R>0} \brak{R}^{\frac{1}{\cxbar}} \bigl| \Ytilde (R,0)\bigr|  \Bigr)
,
\qquad \forall R\geq 1.
\end{align}
\end{subequations}
The bounds in~\eqref{eq:sharp:V:Q:derivative} are consistent with the estimates for $(\bar V,\bar Q, \p_R \bar V,\p_R \bar Q)$ established earlier in~\eqref{eq:barV:barQ:at:infinity}. 
In particular, we have proven:
\begin{proposition}[\bf Lower bound for $Q$]
Assume that the initial data for $\Ytilde$ is taken to be small enough to satisfy
\begin{equation}
 \sup_{R>0} \brak{R}^{\frac{1}{\cxbar}} |\Ytilde(R,0)|
\leq \tfrac{\alpha }{4 \cxbar (2+\alpha) M_0} 
\underline{\eps}^{\frac 16},
\label{eq:extra:ass:infinity}
\end{equation}
where the constant $M_0$ is as in~\eqref{eq:barQ:global:lower}. 
Then, upon taking $\underline{\eps}$ to be sufficiently small with respect to $\alpha, d ,\NNN$, we have
\begin{equation}
| \tilde Q(R,\tau) | \leq \underline{\eps}^{\frac 16} \bar Q(R)
\qquad 
\Rightarrow
\qquad
(1- \underline{\eps}^{\frac 16}) \bar Q(R) \leq    Q(R,\tau) \leq (1+ \underline{\eps}^{\frac 16}) \bar Q(R)  
\label{eq:barQ:dominates:at:infinity}
\end{equation}
for all $R\geq 0$ and all $\tau\geq 0$.
\end{proposition}

\subsubsection{Bounds for $ \tilde K$ and its derivatives}
By definition, we have that $\alpha Q R \p_R \tilde K = \gamma \Ytilde_2 - \alpha \tilde Q R \p_R \bar K$. Thus, from~\eqref{eq:RdR:barH:at:infinity}, \eqref{eq:Ytilde:bound}, \eqref{eq:sharp:V:Q:derivative}, and~\eqref{eq:barQ:dominates:at:infinity}, we deduce 
\begin{subequations}
\label{eq:sharp:K:derivative}
\begin{align}
|\p_R \tilde K(R,\tau)|
\leq C_{\eqref{eq:sharp:K:derivative}} R \brak{R}^{-2-\theta+\frac{1}{\cxbar}} \underline{\eps}^{\frac 15} e^{-\underline{\lambda} \tau}
\end{align}
and also
\begin{align}
 | \p_R \tilde K(R,\tau)  |  
& \leq C_{\eqref{eq:sharp:K:derivative}} R  \underline{\eps}^{\frac 15} e^{-\underline{\lambda} \tau} {\bf 1}_{R<1} 
+ \tfrac{1}{2 \cxbar(2+\alpha)} R^{-1} \underline{\eps}^{\frac 16} {\bf 1}_{R\geq 1} 
+ \underline{\eps}^{\frac 16} 
|\p_R \bar K(R)|
\leq 
\underline{\eps}^{\frac 17}
|\p_R \bar K(R)|,
\end{align}
\end{subequations}
for all $R>0$ and $\tau\geq 0$; upon taking $\underline{\eps}$ to be sufficiently small.

The first bound in~\eqref{eq:sharp:K:derivative}, together with \eqref{eq:boot:0}, the fact that  $\theta < \frac{1}{\cxbar} \leq \theta +1$, and the fundamental theorem of calculus yield the local-in-$R$ decay of the field $\tilde K$, according to
\begin{equation}
|\tilde K(R,\tau)| \leq C_{\eqref{eq:rough:K:no:derivative}} \brak{R}^{\frac{1}{\cxbar}-\theta} \underline{\eps}^{\frac 15} e^{-\underline{\lambda}\tau} 
\label{eq:rough:K:no:derivative}
\end{equation}
In order to obtain sharp estimates for $\tilde K$ with respect to $R$ (as $R\to \infty$), we appeal directly to its evolution in~\eqref{eq:euler:tilde:c}, which implies together with~\eqref{eq:RdR:barH:at:infinity}, \eqref{eq:sharp:V:Q:derivative},  and~\eqref{eq:boot:3} that 
\begin{align*}
\bigr| \partial_\tau \bigl( \tilde K(\Phi_2(R,\tau),\tau) \bigr) \bigl|
\leq (1+\cxbar) C_{\eqref{eq:sharp:V:Q:derivative}}   \underline{\eps}^{\frac 15} e^{-\underline{\lambda} \tau}  + 2 (1+\cxbar) \underline{\eps}^{\frac 45} e^{-\underline{\lambda}\tau}.
\end{align*}
Upon integrating the above expression in time and then taking a supremum with respect to $R$, using that $\Phi_2 (\cdot,\tau) \colon (0,\infty) \to (0,\infty)$ is a bijection, we deduce that 
\begin{align}
|\tilde K(R,\tau)| 
\leq \sup_{R>0} |\tilde K(R,0)| + 2 \underline{\lambda}^{-1} (1+\cxbar) C_{\eqref{eq:sharp:V:Q:derivative}}   \underline{\eps}^{\frac 15}
\leq \sup_{R>0} |\tilde K(R,0)| +  \underline{\eps}^{\frac 16}. 
\label{eq:sharp:K:no:derivative}
\end{align}

\subsection{Propagation of higher regularity}
Recall that the propagation of the order of vanishing claimed in Lemma~\ref{lem:vanishing:at:R=0} (which is to be applied with $\LLL = \MMM = 17 \NNN$ in light of the bootstrap~\eqref{eq:boot:2}) required that the solution $(\tilde V,\tilde Q,\tilde K)$ of \eqref{eq:euler:tilde:a}--\eqref{eq:euler:tilde:c} belongs to $L^\infty([0,T]; W^{17 \NNN+1,\infty}_{\rm loc}(\Reals_+))$, for any $T>0$. On the other hand, the bootstraps from Section~\ref{sec:boot} only provide bounds for $(\tilde V,\tilde Q,\tilde K)$ and their derivatives $(\p_R \tilde V,\p_R \tilde Q,\p_R \tilde K)$. In order to bridge this apparent mismatch, we establish the following \emph{qualitative results}, which are analogous to classical propagation of higher regularity results, using only a control on lower-regularity norms.

%%%%

%%%%

\begin{proposition}[\bf Qualitative propagation of higher regularity]
\label{prop:reg:propagation}
Let $\underline{\eps}$ be sufficiently small with respect to $\alpha, d, \NNN$. 
Assume that the initial datum $(\tilde V_0,\tilde Q_0,\tilde K_0)$ of \eqref{eq:euler:tilde:a}--\eqref{eq:euler:tilde:c} is such that:
\begin{itemize}[leftmargin=1em]
 \item the initial data is smooth, with $( \p_R^{\ell+1} \tilde V_0,   \p_R^{\ell+1} \tilde Q_0,  \p_R^{\ell} (Q_0 \p_R \tilde K_0) )   \in L^\infty(\Reals_+)$  for all $0\leq \ell \leq 17 \NNN$;
 \item let $\mathcal{L}:= \{ \ell \in \Naturals \colon 1 \leq \ell \leq 2 \NNN - 1 \} \cup \{\ell \in 2 \Naturals+ 1\colon 2\NNN+1 \leq \ell \leq 17 \NNN \}$; we assume that 
$(\p_R^\ell \tilde V_0, \p_R^\ell \tilde Q_0, \p_R^\ell \tilde K_0)(0) = 0$ for all integers $\ell \in \mathcal{L}$;
 \item the leading order Taylor coefficients satisfy the bounds~\eqref{eq:Taylor:coeff:assumption:time:zero}, while the higher order even Taylor coefficients satisfy the bound~\eqref{eq:higher:Taylor:coeff:assumption:initial};
 \item the first order derivatives are such that $\Ytilde(R,0)$ satisfies the bounds~\eqref{eq:derivative:boot:closure} and~\eqref{eq:extra:ass:infinity}.
 \end{itemize}
Furthermore, assume\footnote{This assumption automatically holds true on the time interval $[0,T]$, by the closure of the bootstrap argument carried out in Sections~\ref{sec:boot}--\ref{sec:sharp:bounds:R}. Proposition~\ref{prop:reg:propagation} is therefore to be read as a regularity-propagation statement which is to be used in tandem with the bootstrap closure, not in isolation.} that on the time interval $[0,T]$ on which the proposition is to be applied, the bulk perturbation $(\tilde V,\tilde Q,\tilde K)$ enjoys the sharp pointwise bounds~\eqref{eq:sharp:V:Q:derivative} and~\eqref{eq:Ytilde:bound}.
Then, for any $T>0$ and $0 \leq \ell \leq 17\NNN$, we have $(\p_R^{\ell+1} \tilde V, \p_R^{\ell+1} \tilde Q, \p_R^\ell (Q \p_R \tilde K) ) \in L^\infty([0,T]; L^{\infty}(\Reals_+))$. In particular, $\p_R^{\ell+1} \tilde K \in L^\infty([0,T]; L_{\rm loc}^{\infty}(\Reals_+))$.
\end{proposition}

Note that for simplicity we did not include any $R$ weights (neither at $R=0$ nor as $R\to \infty$) in the statement of Proposition~\ref{prop:reg:propagation}; sharp weighted bounds for these higher order derivatives may be obtained by following the proof below and adding the natural weights at $R=0$ and as $R\to \infty$; we do not pursue these sharper bounds here.

\begin{proof}[Proof of Proposition~\ref{prop:reg:propagation}]
Let $\MMM = 17 \NNN$.
Note from~\eqref{eq:QVH:via:ring} that
$\p_R V = \frac{1}{2R} (\YY_1 + \YY_3)$, $\p_R Q= \frac{1}{2R} (\YY_3 - \YY_1 + \frac{2}{\alpha} \YY_2)$, and $Q \p_R K = \tfrac{\gamma}{\alpha R} \YY_3$. Therefore, for any $1\leq \ell \leq \MMM$, we have the identities:
\begin{equation}
\p_R^{\ell+1} V = \tfrac{1}{2}  \p_R^{\ell} (\tfrac{1}{R} \YY_1 + \tfrac{1}{R}\YY_3),
 \quad
 \p_R^{\ell+1} Q = \tfrac{1}{2}   \p_R^{\ell}(\tfrac{1}{R} \YY_3 - \tfrac{1}{R}  \YY_1 + \tfrac{2}{\alpha} \tfrac{1}{R} \YY_2),
 \quad
 \p_R^\ell( Q \p_R K) = \tfrac{\gamma}{\alpha}   \p_{R}^{\ell} (\tfrac{1}{R} \YY_3) .
 \label{eq:V:Q:YY:relationship}
\end{equation}
Thus, it appears to be useful to estimate $\p_R^\ell (\frac{1}{R} \YY)$ for all $1\leq \ell \leq \MMM$.

Dividing~\eqref{eq:YY:evo} by $R$ and then applying $\p_R^\ell$ to the resulting PDE, we obtain
\begin{align}
&\p_\tau \bigl( \p_R^\ell (\tfrac 1R \YY)\bigr) 
+ \TT \, R \p_R \bigl(\p_R^\ell (\tfrac 1R \YY)\bigr) 
+ \bigl(\DD + (\ell-1) \TT+ \ell R \p_R \TT \bigr) \bigl(\p_R^\ell (\tfrac 1R \YY)\bigr)   
+ \p_R^\ell \Bigl( R \NN\bigl(\tfrac{1}{R} \YY, \tfrac{1}{R} \YY\bigr) \Bigr)
\notag\\
&
= - \sum_{j=1}^{\ell-1} {\ell \choose j+1} 
\bigl((j+1) \p_R^{j} \TT + R \p_R^{j+1} \TT\bigr) \p_R^{\ell-j} (\tfrac 1R \YY)
- \sum_{j=1}^{\ell} {\ell \choose j} \p_R^{j} ( \DD- \TT) \p_R^{\ell-j}(\tfrac 1R \YY).
\label{eq:Leibniz:extra:weird}
\end{align}
Since the nonlinear term $\NN$ is quadratic with constant coefficients (see~\eqref{eq:YY:evo:NN}), we may directly apply estimate~\eqref{eq:appendix:interpolate} to bound
\begin{align*}
&\|\p_R^\ell \bigl( R \NN\bigl(\tfrac{1}{R} \YY, \tfrac{1}{R} \YY\bigr) \bigr)\|_{L^\infty([0,\infty))}
\notag\\
&\qquad
\les_{\ell,\alpha}
\| \tfrac{\brak{R}}{R} \YY \|_{L^\infty([0,\infty))}
\|\p_R^{\ell} (\tfrac{1}{R} \YY) \|_{L^\infty([0,\infty))}
+
\|  \tfrac{\brak{R}}{R} \YY \|_{L^\infty([0,\infty))}^{1 + \frac{1}{\ell}}
\|\p_R^{\ell} (\tfrac{1}{R} \YY)\|_{L^\infty([0,\infty))}^{1 - \frac{1}{\ell}}
.
\end{align*}
Similarly, using definitions~\eqref{eq:YY:evo:TT} and~\eqref{eq:YY:evo:DD}, estimate~\eqref{eq:LK:ineq} and identity~\eqref{eq:V:Q:YY:relationship} we may estimate
\begin{align*}
& \| \p_R^{j} \TT\|_{L^\infty([0,\infty))}
\les_{\ell,j,\alpha}
\| \tfrac{1}{R} \YY\|_{L^\infty([0,\infty))}^{1-\frac{j-1}{\ell}}
\| \p_R^{\ell} (\tfrac{1}{R} \YY) \|_{L^\infty([0,\infty))}^{\frac{j-1}{\ell}},
\\
& \| \p_R^{j} ( \DD- \TT) \|_{L^\infty([0,\infty))}
\les_{\ell,j,\alpha}
\| \tfrac{1}{R} \YY\|_{L^\infty([0,\infty))}^{1-\frac{j-1}{\ell}}
\| \p_R^{\ell} (\tfrac{1}{R} \YY) \|_{L^\infty([0,\infty))}^{\frac{j-1}{\ell}},
\\
& \| \p_R^{\ell-j}(\tfrac 1R \YY)\|_{L^\infty([0,\infty))}
\les_{\ell,j}
\| \tfrac{1}{R} \YY\|_{L^\infty([0,\infty))}^{\frac{j}{\ell}}
\| \p_R^{\ell} (\tfrac{1}{R} \YY) \|_{L^\infty([0,\infty))}^{1-\frac{j}{\ell}}.
\end{align*}
By appealing to Lemma~\ref{lem:pointwise} we have that 
\[
\| R \, \p_R^{j+1} \TT \, \p_R^{\ell-j} (\tfrac 1R \YY) \|_{L^\infty([0,\infty))}
\les_{\ell,j,\alpha}
\| \tfrac{\brak{R}}{R} \YY \|_{L^\infty([0,\infty))}
\|\p_R^{\ell} (\tfrac{1}{R} \YY) \|_{L^\infty([0,\infty))}
.
\]
Lastly, we note that because $\underline{\eps}$ is sufficiently small, using that $\cxbar - \cubar =1$, $\cxbar + \bar V(R)>0$, and appealing to the bound~\eqref{eq:sharp:V:Q:derivative}, we deduce via~\eqref{eq:YY:evo:TT}, \eqref{eq:YY:evo:DD} and~\eqref{eq:QVH:via:ring}, that 
\[
\|
\bigl(\DD + (\ell-1) \TT+ \ell R \p_R \TT \bigr)_-
\|_{L^\infty([0,\infty))}
\les_{\alpha,d}
\|V \|_{L^\infty([0,\infty))} + \|\alpha Q\|_{L^\infty([0,\infty))} 
+ \|\YY \|_{L^\infty([0,\infty))}.
\]
Here $(x)_-$ denotes the negative part of $x$.

With these bounds, we return to\eqref{eq:Leibniz:extra:weird}, 
compose with the ODE flow $\Phi_i$ defined in~\eqref{eq:Phi:i:def}, and appealing to the six above-established inequalities, we obtain that 
\begin{align*}
&\tfrac{d}{d\tau} \bigl| \p_R^{\ell} (\tfrac{1}{R} \YY_i) \circ \Phi_i(R,\tau) \bigr| 
\notag\\
&
\les_{\alpha,d,\ell}
\bigl( \|V \|_{L^\infty([0,\infty))} + \|\alpha Q\|_{L^\infty([0,\infty))}  + \| \tfrac{\brak{R}}{R} \YY \|_{L^\infty([0,\infty))}\bigr)
\|\p_R^{\ell} (\tfrac{1}{R} \YY) \|_{L^\infty([0,\infty))}
+ \| \tfrac{\brak{R}}{R} \YY \|_{L^\infty([0,\infty))}^2.
\end{align*}
Integrating this inequality, taking the supremum over all $R>0$ and the maximum over $i \in \{1,2,3\}$, we obtain that 
\begin{align}
&\|\p_R^{\ell} (\tfrac{1}{R} \YY)(\cdot,\tau) \|_{L^\infty([0,\infty))} 
\notag\\
&\leq \left( \|\p_R^{\ell} (\tfrac{1}{R} \YY)(0,\cdot) \|_{L^\infty([0,\infty))} + C_{\eqref{eq:YY:Gronwall}} \int_0^\tau \| \tfrac{\brak{R}}{R} \YY(\cdot,\tau') \|_{L^\infty([0,\infty))}^2 d\tau' \right)
\notag\\
&\times \exp\left( C_{\eqref{eq:YY:Gronwall}} \int_0^\tau \|V(\cdot,\tau') \|_{L^\infty([0,\infty))} + \|\alpha Q(\cdot,\tau')\|_{L^\infty([0,\infty))}  + \| \tfrac{\brak{R}}{R} \YY(\cdot,\tau') \|_{L^\infty([0,\infty))} d\tau' \right)
\label{eq:YY:Gronwall}
\end{align}
where $C_{\eqref{eq:YY:Gronwall}}>0$ only depends on $\alpha$, $d$, and $\ell$.
By further appealing to~\eqref{eq:barV:barQ:at:infinity}, ~\eqref{eq:Ybar:at:infinity}, \eqref{eq:Ytilde:bound}, and the third estimate in~\eqref{eq:sharp:V:Q:derivative}, we deduce from~\eqref{eq:YY:Gronwall} that there exists a constant $C_{\eqref{eq:YY:Gronwall:2}}>0$, which only depends on $\alpha, d, \NN$, and $\ell$, such that 
\begin{align}
\|\p_R^{\ell} (\tfrac{1}{R} \YY)(\cdot,\tau) \|_{L^\infty([0,\infty))} 
\leq \left( \|\p_R^{\ell} (\tfrac{1}{R} \YY)(0,\cdot) \|_{L^\infty([0,\infty))} +\tau \, C_{\eqref{eq:YY:Gronwall:2}} \right) e^{ \tau\, C_{\eqref{eq:YY:Gronwall:2}}}
\label{eq:YY:Gronwall:2}
\end{align}
for all $\tau \geq 0$. In light of~\eqref{eq:V:Q:YY:relationship} the above estimate provides bounds for $\|\p_R^{\ell+1} V\|_{L^\infty}$, $\| \p_R^{\ell+1} Q\|_{L^\infty}$, and $\| \p_R^{\ell} (Q\p_R K)\|_{L^\infty}$ for all $\tau \geq 0$.

To conclude the proof, we combine the above conclusion with the available regularity for $\bar V, \bar Q, \bar K$, and the triangle inequality. The $L^\infty_{\rm loc}$ regularity for $\p_R^{\ell+1} \tilde K$ now follows from the available upper and lower bounds for $Q$ in~\eqref{eq:barQ:dominates:at:infinity}.
\end{proof}
 
\subsection{Main result: stability of the implosion profiles in radial symmetry}

The outcome of the analysis in this section is that we have estalished:

\begin{theorem}[\bf Stability in Radial Symmetry]
\label{thm:stability:radial}
Fix $d\in \{1,2,3\}$, $1 < \gamma \leq 2d+1$, and let $\alpha = \frac{\gamma-1}{2}$; also, fix $\NNN\geq 1$. Let $(\bar V,\bar Q,\bar H)$ denote the smooth self-similar profile constructed in Theorem~\ref{thm:main:profiles}, with associated similarity exponents $(\cxbar, \cubar, \cbbar) = (\cxstar(d,\gamma,\NNN),\cxstar(d,\gamma,\NNN)-1,\cbstar(d,\gamma,\NNN))$ given by~\eqref{eq:cx:admissible}--\eqref{eq:cb:admissible}. Let $\bar K = \log \bar H$. Define $\MMM := 17 \NNN$, $\theta := \frac{1}{2}\min\{1, \frac{1}{\cxbar}\}$, and let $\underline{\lambda}=\underline{\lambda}(d,\gamma,\NNN)>0$ be defined explicitly by \eqref{eq:underline:lambda:def}. 
There exist a sufficiently small constant $\underline{\eps}^* = \underline{\eps}^*(\alpha,d,\NNN) >0$, such that for any $\underline{\eps} \leq \underline{\eps}^*$, the following holds.

Denote the initial data for~\eqref{eq:euler:final} by $(V,Q,K)|_{\tau=0}$. 
With~\eqref{eq:self-similar:ansatz} and~\eqref{eq:class:of:perturbations}, this  corresponds to initial data for density, velocity and pressure in the original formulation of the Euler equations~\eqref{eq:euler0} given by 
\[
\uu(x,-1) = x V(|x|,0), 
\quad
\rho(x,-1) = (\tfrac{\alpha Q}{H})^{\frac{1}{\alpha}}(|x|,0),
\quad
p(x,-1) = \tfrac{1}{\gamma} |x|^2 ( \tfrac{(\alpha Q)^\gamma}{H})^{\frac{1}{\alpha}}(|x|,0) .
\] 
In particular,  $p(0,-1)=0$.
In self-similar radial coordinates, we denote the \emph{initial perturbation} by 
\[
(\tilde V_0,\tilde Q_0,\tilde K_0) = (V|_{\tau=0} - \bar V, Q|_{\tau=0}-\bar Q, K|_{\tau=0} - \bar K);
\]
this serves as a set of initial conditions for the PDE system~\eqref{eq:euler:tilde}.  We assume this initial perturbation is chosen such that:

\begin{enumerate}[label=(\roman*),leftmargin=2em]
\item \textsl{(Regularity).} We have $(u,\rho,p)|_{t=-1} \in C^{\MMM+1}(\mathbb{R}^d)$ and $(\partial_R \tilde V_0, \partial_R \tilde Q_0, Q_0 \partial_R \tilde K_0) \in W^{\MMM,\infty}([0,\infty))$.

\item \textsl{(Vanishing Taylor coefficients at $R=0$).} For all \emph{even} integers $\ell \in \{2, \ldots, 2\NNN-2\}$, we have $(\partial_R^\ell \tilde V_0, \partial_R^\ell \tilde Q_0, \partial_R^\ell \tilde K_0)(0) = (0,0,0)$.

\item \textsl{(Smallness of leading Taylor coefficients).} 
The first two sets of nonzero Taylor coefficients satisfy
\begin{subequations}
\label{eq:thm:init:Taylor}
\begin{align}
|\tilde V_0(0)| + |\tilde Q_0(0)| + |\tilde K_0(0)| &\leq \tfrac{1}{4} \underline{\eps}, \label{eq:thm:init:0}\\
\tfrac{1}{(2\NNN)!}|\partial_R^{2\NNN} \tilde V_0(0)| + \tfrac{1}{(2\NNN)!}|\partial_R^{2\NNN} \tilde Q_0(0)| + \tfrac{1}{(2\NNN)!}|\partial_R^{2\NNN} \tilde K_0(0)| &\leq \underline{\eps}. \label{eq:thm:init:N}
\end{align}

\item \textsl{(Smallness of higher order Taylor coefficients).} For all \emph{even} integers $\ell \in \{2\NNN+2,  \ldots, \MMM\}$ we have
\begin{equation}
\tfrac{1}{\ell!}|\partial_R^{\ell} \tilde V_0(0)| + \tfrac{1}{\ell!}|\partial_R^{\ell} \tilde Q_0(0)| + \tfrac{1}{\ell!}|\partial_R^{\ell} \tilde K_0(0)| \leq \underline{\eps}^{\frac 35}.
\end{equation}
\end{subequations}

\item \textsl{(Smallness of first order derivatives).} 
Using the initial conditions $(\tilde V_0,\tilde Q_0,\tilde K_0)$, define  $\Ytilde_0(R) := \YY(R,0) - \Ybar(R)$, where $\YY$ is given by via~\eqref{eq:ringW:ringZ:ringA:def}--\eqref{eq:YY:def}. We assume that 
\begin{subequations}
\label{eq:thm:init:JM}
\begin{equation}
\sup_{R > 0} \langle R \rangle^{\frac{1}{\cxbar}} |\Ytilde_0(R)| 
\leq  \underline{\eps}^{\frac 15}. \label{eq:thm:init:decay}
\end{equation}
Moreover, let $\psi$ be the non-decreasing weight constructed in Proposition~\ref{prop:choice:of:M:psi}. In terms of this weight $\psi$, define   
the operator $\JJJ_{\MMM}$ via~\eqref{eq:Taylor:poly:1} and~\eqref{eq:J:M:def}. Our last assumption is that
\begin{equation}
\sup_{R > 0} \langle R \rangle^\theta |\JJJ_{\MMM} \Ytilde_0(R)| 
\leq \underline{\eps}^{\frac 25}. \label{eq:thm:init:deriv}
\end{equation}
\end{subequations}
\end{enumerate}

Under these assumptions, there exists a unique global-in-$\tau$ smooth solution $(V,Q,K)$ of~\eqref{eq:euler:final} with the prescribed initial data, with modulation functions $(\cx,\cu,\cb)$ given by~\eqref{eq:L:infty:perturbation:b} and~\eqref{eq:modulation:final}. Moreover:
\begin{enumerate}[label=(\alph*),leftmargin=2em]
\item \textsl{(Exponential decay of perturbation).} The perturbation $(\tilde V, \tilde Q, \tilde K) := (V - \bar V, Q - \bar Q, K - \bar K)$ satisfies
\begin{subequations}
\label{eq:thm:conc:Taylor}
\begin{align}
|\tilde V(0,\tau)| + |\tilde Q(0,\tau)| + |\tilde K(0,\tau)| 
&\leq  \underline{\eps} \, e^{-\underline{\lambda} \tau}, \label{eq:thm:conc:0}\\
\tfrac{1}{(2\NNN)!}|\partial_R^{2\NNN} \tilde V(0,\tau)| + \tfrac{1}{(2\NNN)!}|\partial_R^{2\NNN} \tilde Q(0,\tau)| + \tfrac{1}{(2\NNN)!}|\partial_R^{2\NNN} \tilde K(0,\tau)| &\leq \underline{\eps}^{\frac 35} e^{-\underline{\lambda} \tau}
\end{align}
\end{subequations}
and the higher order Taylor coefficients at $R=0$ decay exponentially according to~\eqref{eq:boot:2}. The fields $\tilde V$, $\tilde Q$, $\tilde K$ and their derivatives decay according to
\begin{subequations}
\label{eq:thm:conc:JM}
\begin{align}
\brak{R}^\theta   |\tilde V(R,\tau)| + \brak{R}^\theta |\tilde Q(R,\tau)|  + \brak{R}^{\theta -\frac{1}{\cxbar}} |\tilde K(R,\tau)| 
&\leq C^* \underline{\eps}^{\frac 15} e^{-\underline{\lambda} \tau}, \label{eq:thm:conc:VQ}\\
\brak{R}^\theta   |\p_R \tilde V(R,\tau)| + \brak{R}^\theta |\p_R \tilde Q(R,\tau)| + \brak{R}^{\theta -\frac{1}{\cxbar}} |\p_R \tilde K(R,\tau)|  &\leq R \brak{R}^{-2} C^* \underline{\eps}^{\frac 15} e^{-\underline{\lambda} \tau}
, \label{eq:thm:conc:deriv}
\end{align}
\end{subequations}
for all $R>0$ and all $\tau \geq 0$, where $C^* = C^*(\alpha,d,\NNN)>0$ is a constant.   

\item \textsl{(Modulation functions).} The modulation functions $(\cxtilde, \cutilde, \cbtilde) := (\cx - \cxbar, \cu - \cubar, \cb - \cbbar)$ defined in~\eqref{eq:modulation:final} decay exponentially fast according to
\[
|\cxtilde(\tau)| + |\cutilde(\tau)| + |\cbtilde(\tau)| \leq  \underline{\eps}^{\frac 45} e^{-\underline{\lambda} \tau},
\qquad
\mbox{for all } \tau \geq 0.
\]

\item \textsl{(Sharp decay at infinity).} For all $R \geq 1$ and $\tau \geq 0$ we have
\begin{subequations}
\label{eq:thm:conc:infty}
\begin{align}
R^{ \frac{1}{\cxbar} }|\tilde V(R,\tau)| + R^{ \frac{1}{\cxbar} } |\tilde Q(R,\tau)| &\leq C^* \underline{\eps}^{\frac 16} , \label{eq:thm:conc:sharp:VQ}\\
R^{1+\frac{1}{\cxbar}} |\partial_R \tilde V(R,\tau)| + R^{1+\frac{1}{\cxbar}} |\partial_R \tilde Q(R,\tau)| + R |\partial_R \tilde K(R,\tau)| &\leq C^* \underline{\eps}^{\frac 16} . \label{eq:thm:conc:sharp:deriv}
\end{align}
\end{subequations}
Additionally, we have that $\tilde K$ is globally bounded as $\sup_{R>0}|\tilde K(R,\tau)| \leq \sup_{R>0} |\tilde K(R,0)| +  \underline{\eps}^{\frac 16}$ for all $\tau \geq 0$, and $Q$ obeys  global upper and lower bounds: $(1- \underline{\eps}^{\frac 16}) \bar Q(R) \leq    Q(R,\tau) \leq (1+ \underline{\eps}^{\frac 16}) \bar Q(R)$.  

\item \textsl{(Propagation of higher-order regularity).} For any $T > 0$ and $1 \leq \ell \leq \MMM$:
\[
(\partial_R^{\ell+1} \tilde V, \partial_R^{\ell+1} \tilde Q, \partial_R^\ell(Q \partial_R \tilde K)) \in L^\infty([0,T]; L^\infty([0,\infty))).
\]
\end{enumerate}
\end{theorem}

An immediate consequence of Theorem~\ref{thm:stability:radial} is:
\begin{corollary}[\bf Blowup in physical variables]
\label{cor:radial:blowup:physical}
The solution $(V,Q,K)$ constructed in Theorem~\ref{thm:stability:radial} corresponds, via the transformations~\eqref{eq:self-similar:ansatz} and~\eqref{eq:class:of:perturbations}, to a smooth radially symmetric solution $(u^r, \sigma, b)$ of the compressible Euler equations~\eqref{eq:euler2} on $[0,\infty) \times [-1, t_*)$, where the blowup time is defined by
\[
t_* = -1 + \int_0^\infty \mathsf{C_r}(\tau') \mathsf{C_u}^{-1}(\tau') \, d\tau',
\]
and satisfies the bound
\[
|t_*| \leq 2 \underline{\lambda}^{-1} \underline{\eps}^{\frac 45}.
\] 
As $t \to t_*^-$, the solution develops a singularity at the origin with self-similar blowup rates determined by the profile $(\bar V, \bar Q, \bar K)$.
\end{corollary}

\begin{remark}[\bf Stability holds modulo finitely many, explicit, compatibility conditions]
\label{rem:ground:state:stable}
Inspecting the assumptions on the initial data in Theorem~\ref{thm:stability:radial}, we notice that items (i), (iii), (iv), and (v) are open conditions. These assumptions allow us to perturb the stationary self-similar profiles $(\bar V,\bar Q,\bar K)$ in an open set, with respect to a suitable topology. The \emph{finite co-dimension stability} aspect of Theorem~\ref{thm:stability:radial} is entirely due to assumption (ii). Assumption (ii)  poses the compatibility conditions $(\partial_R^\ell \tilde V_0, \partial_R^\ell \tilde Q_0, \partial_R^\ell \tilde K_0)(0) = (0,0,0)$ for all \emph{even} integers $\ell$, with $2\leq \ell \leq 2\NNN-2$. There are $(\NNN-1)$ many such integers. Thus, for smooth perturbations we are only imposing $3(\NNN-1)$ many compatibility conditions on the initial perturbations at $R=0$ (these compatibility conditions are automatically propagated forward in time due to Lemma~\ref{lem:vanishing:at:R=0}). 
The important aspect of this precise characterization of the compatibility conditions for the initial data is that for the ``ground state'' corresponding to $\NNN=1$, the number of compatibility conditions vanishes, and thus \emph{we may perturb the globally self-similar ground state profiles in an open set}.
\end{remark}

\begin{proof}[Proof of Theorem~\ref{thm:stability:radial}]
We provide the roadmap for the proof of Theorem~\ref{thm:stability:radial}; all estimates were already proven earlier in this section, so here we only indicate where/how each component is proven.

\noindent\emph{\underline{Setup and local existence.}}
The system~\eqref{eq:euler:final} for $(V,Q,K)$ is derived in Section~\ref{sec:true:ss:ansatz}. Local-in-$\tau$ existence and uniqueness of smooth solutions follows from standard theory for symmetric hyperbolic systems, given the smooth initial data satisfying assumptions (i)--(v); in our setting, this local existence also follows from Proposition~\ref{prop:reg:propagation}. The perturbation equations~\eqref{eq:euler:tilde} are derived in Section~\ref{sec:stability}.

\noindent\emph{\underline{Taylor structure at $R=0$.}}
By the assumed regularity of $(u,\rho,p)|_{t=-1}$, for all \emph{odd} integers $1\leq \ell\leq \MMM$ we have that $(\p_R^\ell V, \p_R^\ell Q,\p_R^\ell K)|_{\tau=0}(0) = (0,0,0)$. Since the profiles $\bar V,\bar Q, \bar H$ also enjoy this property, we deduce that $(\partial_R^\ell \tilde V_0, \partial_R^\ell \tilde Q_0, \partial_R^\ell \tilde K_0)(0) = (0,0,0)$ for all \emph{odd} integers $1\leq \ell \leq \MMM$. Assumption (ii) also implies that  $(\partial_R^\ell \tilde V_0, \partial_R^\ell \tilde Q_0, \partial_R^\ell \tilde K_0)(0) = (0,0,0)$ for all \emph{even} integers $2\leq \ell \leq 2 (\NNN-1)$. The Taylor structure at $R=0$ provided by the class of such functions is preserved by the flow due to Lemma~\ref{lem:vanishing:at:R=0} and the propagation of  $W^{\MMM+1,\infty}_{\rm loc}$ regularity established in Proposition~\ref{prop:reg:propagation}.

\noindent\emph{\underline{Choice of modulation functions.}}
The modulation functions $(\cxtilde, \cutilde, \cbtilde)$ are defined in~\eqref{eq:modulation:final}. Note that the constraint~\eqref{eq:modulation:1} removes the unstable eigenvalue of the matrix $A_0$ (cf.~\eqref{eq:euler:tilde:vq:zero}), the constraint~\eqref{eq:modulation:2} ensures damping for $\tilde k_0$, and the constraint~\eqref{eq:modulation:3} removes the neutral mode at order $R^{2\NNN}$.

\noindent\emph{\underline{Bootstrap assumptions.}}
We introduce the bootstrap assumptions~\eqref{eq:boot:0}--\eqref{eq:boot:4}. These control:
\begin{itemize}[leftmargin=1em]
\item The Taylor coefficients $(\tilde v_0, \tilde q_0, \tilde k_0)$ at $R=0$ via~\eqref{eq:boot:0};
\item The Taylor coefficients $(\tilde v_\NNN, \tilde q_\NNN, \tilde k_\NNN)$ at order $R^{2\NNN}$ via~\eqref{eq:boot:1};
\item Higher Taylor coefficients via~\eqref{eq:boot:2};
\item The modulation functions via~\eqref{eq:boot:3};
\item The weighted derivative $\langle R \rangle^\theta \JJJ_\MMM \Ytilde$ via~\eqref{eq:boot:4}.
\end{itemize}

\noindent\emph{\underline{Closure of  bootstraps.}}
Propositions~\ref{prop:Taylor:boot:closure}, \ref{prop:Taylor:boot:closure:higher}, and~\ref{prop:derivative:boot:closure} establishes the ``closure'' (meaning, improvement by a factor $<1$) of the bootstrap assumptions. 

\noindent\emph{\underline{Sharp bounds.}}
The stated decay estimates for the perturbations are either direct consequences of the bootstraps, see Section~\ref{sec:boot:consequence}, or they are established in Section~\ref{sec:sharp:bounds:R}.
 
\noindent\emph{\underline{Global existence and blowup.}}
Since all bootstrap assumptions close with strictly improved constants, a standard continuity argument shows that the bootstraps hold for all $\tau \geq 0$. This establishes global-in-$\tau$ existence. The transformation back to physical variables $(u^r, \sigma, b)$ via~\eqref{eq:self-similar:ansatz} yields a solution on $[-1, t_*)$ that blows up at time $t_*$. The bound on $t_*$ follows from the exponential decay of $(\cxtilde, \cutilde)$ and definition~\eqref{eq:self-similar:ansatz:tau}.
\end{proof}

%%%%%%%%%%%%%%%%%%%

\section{Stability outside of radial symmetry}
\label{sec:nonradial}

In Section~\ref{sec:stability} we have analyzed the stability of the globally self-similar solutions of the Euler equations
$(\bar \uu(x,t),\bar \sigma(x,t),\bar b(x,t)) = (\vec{e}_r \bar u^r(|x|,t), \bar \sigma(|x|,t), \bar b(|x|,t))$, 
defined via~\eqref{eq:globally:SS:exact} in terms of  smooth implosion profiles $(\bar U,\bar \Sigma,\bar B)$, within the class of smooth radially symmetric solutions of Euler; that is, as solutions to~\eqref{eq:euler2}. In this section, we analyze the stability of these solutions with respect to smooth \emph{non-radial} perturbations to Euler; that is, as solutions to~\eqref{eq:euler1}, or equivalently,~\eqref{eq:euler:primary}.

The main result of this section is Theorem~\ref{thm:nonrad:stab}, below. This result proves that the global-in-time stability of the stationary profiles is controlled by the behavior of finitely many Taylor coefficients at $x=0$. In particular, this implies finite co-dimension stability; see Remark~\ref{rem:data:for:nonrad:stability}. Theorem~\ref{thm:nonrad:stab} does not seek a sharp bound on the number of potentially unstable directions (see assumption~\eqref{eq:ass:EOM}); instead, the  sharp characterization of non-radial instabilities is analyzed in Section~\ref{sec:ODE}, see Theorem~\ref{thm:ODE_stab}.

\subsection{A convenient change of unknowns}
\label{sec:convenient:unknowns}
In the radial setting of Section~\ref{sec:stability}, the stability analysis is carried out via the modulated self-similar transformation~\eqref{eq:self-similar:ansatz}, renormalization at $|x|=0$ through~\eqref{eq:class:of:perturbations}, and the evolution of the smooth functions $(V,Q,K)$ satisfying~\eqref{eq:euler:final}. A key underlying ingredient is the classical fact that radial scalar functions which are smooth in $x$ near the origin admit a representation as smooth functions of $|x|^2$. 

In the non-radial setting considered in this section we lose this equivalence between regularity in $x$ and smoothness with respect to $|x|^2$. As such, if we were to work with unknowns $(V,Q,K)$ defined analogously to~\eqref{eq:class:of:perturbations}, we would be dealing with functions that are a priori not smooth. To avoid this issue we recall that for the Euler system in terms of the primary flow variables $(\rho,u,p)$, namely~\eqref{eq:euler:primary}, the initial data is smooth with respect to $x$. Therefore, it is natural to perform our non-radial stability analysis directly at the level of~\eqref{eq:euler:primary} (for $(\rho,u,p)$), not at the level of~\eqref{eq:euler1} (for $(\sigma,u,b)$). For technical reasons, we do not do this. Instead of the pressure equation~\eqref{eq:euler:p}, which has a quadratic non-linear term, we use the transport equation for entropy which has no (non-transport) non-linear term. Then, the nonlinear term in the velocity equation~\eqref{eq:euler:u} must then be changed from $\rho^{-1} \nabla p$ to another  nonlinear term which involves entropy gradients; in order to keep this nonlinear term quadratic and without fractions (which do not behave nicely with respect to the Leibniz rule), and taking advantage of the fact that the density does not vanish, we choose to work with the evolution equation for $\rho^{\gamma-1}$ instead of~\eqref{eq:euler:rho}.  For convenience, let us denote
\begin{equation} \label{eq:nonrad_var}
  \vvrho := \rho^{\gamma-1} = \rho^{2\alpha}, \qquad \mbb := e^s =  \tfrac{\gamma p }{\rho^{\gamma}}.
\end{equation}
Note that $\mbb = b^2$, where $b$ is as defined in~\eqref{eq:b:def}.
With this notation, the sound speed is
\[
c = \sqrt{ \tfrac{\gamma p}{\rho}} = \sqrt{ \vvrho \mbb }
,
\]
and we may  rewrite the nonlinear term $ \rho^{-1} \na p$  in terms of $\vvrho, \mbb$ as
\begin{equation*}
\rho^{-1} \na p
= \rho^{-1} \nabla \bigl(\tfrac{1}{\gamma} \rho^{\gamma} \mbb \bigr)
= \tfrac{1}{\gamma} \rho^{\gamma-1} \nabla \mbb
    + \tfrac{1}{\gamma-1} \mbb \na \bigl(\rho^{\gamma-1}\bigr) 
= \tfrac{1}{\gamma} \vvrho \nabla \mbb + \tfrac{1}{ 2\alpha}  \mbb \nabla \vvrho 
.
\end{equation*}
We may thus rewrite the Euler system~\eqref{eq:euler:primary} in an equivalent form as an evolution equation for $(\eta,\uu,\mbb)$:
\begin{subequations}
\label{eq:sys_UB}
\begin{align} 
\pa_t \vvrho +  \uu  \cdot \nabla \vvrho  + 2 \alpha \vvrho \, \div \uu &=  0, \\
\pa_t \uu +  \uu  \cdot \nabla  \uu  + \tfrac{1}{2\alpha} \mbb \nabla \vvrho + \tfrac{1}{\gamma} \vvrho \nabla \mbb &= 0 , \\
 \pa_t \mbb +  \uu  \cdot \nabla \mbb &= 0.
 \end{align}
\end{subequations}

\begin{remark}[\bf Velocity and pressure vanish at the origin] 
The initial data for~\eqref{eq:sys_UB} is taken to be smooth; this initial data is such that the non-negative function $\mbb$ attains a global minimum at a unique point, and that $\min_{x\in\mathbb{R}^3} \mbb(\cdot,-1) = 0$. By the Gallilean symmetry of the Euler equations, we may assume without loss of generality that $\mbb(0,-1) =0$; being a global mimimum of a smooth function, it automatically follows that $\nabla \mbb(0,-1)=0$. 
Also by Gallilean symmetry, we may assume without loss of generality that $\uu(0,-1)=0$.
It follows from~\eqref{eq:sys_UB} that 
\begin{align*}
& \pa_\tau \uu(0,\cdot) + \uu(0,\cdot) \cdot \nabla  \uu(0,\cdot)  + \tfrac{1}{2\alpha} \mbb(0,\cdot) \nabla \vvrho (0,\cdot) + \tfrac{1}{\gamma} \vvrho(0,\cdot) \nabla \mbb(0,\cdot)  = 0 \, , \\
& \pa_\tau \mbb (0,\cdot)  +  \uu(0,\cdot)  \cdot \nabla \mbb(0,\cdot)  =  0 \, , \\
& \pa_\tau \nabla \mbb(0,\cdot)  + \nabla \uu(0,\cdot)  \cdot \nabla \mbb(0,\cdot)  + \uu(0,\cdot) \cdot \nabla \nabla \mbb(0,\cdot)  = 0 \, .
\end{align*}
It follows that as long as $\uu$ remains $C^1$ smooth near $x=0$ and $\mbb$ remains $C^2$ smooth near $x=0$ for $t\in [-1,T)$, then we have
\begin{equation}
\label{eq:vanishing:order:at:origin:default}
 \uu(0,t) = 0 \,, \qquad 
 \mbb(0,t) = 0 \,, \qquad 
 \nabla \mbb(0,t) = 0 \,, 
\end{equation}
for all $t\in [-1,T)$. Throughout this section we work within the class of solutions which satisfies~\eqref{eq:vanishing:order:at:origin:default}. 
\end{remark}

\subsection{Modulated self-similar ansatz}
\label{sec:modulate:SS:ansatz:nonradial}
In analogy to~\eqref{eq:cr:cu:cb:first} and~\eqref{eq:Cr:Cu:Cb:def}, we define
\begin{subequations}
\label{eq:normal-1}
\begin{align}
C_r(\tau)  =  \exp\Bigl( - \int_0^{\tau}  \cx (\tau') d \tau' \Bigr), \qquad
C_u(\tau)  =  \exp\Bigl(  - \int_0^{\tau} \cu (\tau') d \tau' \Bigr), \\ 
C_\mb(\tau) =  \exp\Bigl(  - \int_0^{\tau} \cbb (\tau') d \tau' \Bigr), \qquad
C_{\vrho}(\tau) =  \exp\Bigl(  - \int_0^{\tau} \crho(\tau') d \tau' \Bigr), 
\end{align}
where the modulation functions 
\[
\cx, \cu,\cbb,\crho \colon [0,\infty) \to \Reals.
\]
are connected via
\begin{equation}
\label{eq:normal0}
\crho (\tau) + \cbb (\tau) = 2 \cu(\tau)
\end{equation}
\end{subequations}
which ensures that $c = \sqrt{ \vvrho \mbb }$ has dynamically the same scaling as $u$. Then we have three degrees of freedom to rescale the system \eqref{eq:sys_UB} dynamically: either $(\cx, \cu , \cbb)$, or equivalently $(\cx, \cu, \crho)$. 
As in~\eqref{eq:self-similar:ansatz} we   let 
\begin{subequations}
\label{eq:SS_ansatz}
\begin{equation}
t = t(\tau) := -1 + \int_0^{\tau}   C_r(\tau') C_u^{-1}(\tau')  d \tau' , 
\end{equation}
rescale the space coordinate $x$ as 
\begin{equation}
y = x \, C_r^{-1}(\tau), \qquad |y| = R = r \,  C_r^{-1}(\tau),
\end{equation}
and rescale the unknowns from~\eqref{eq:sys_UB} as
\begin{align}
\label{eq:self-similar:new:variables}
\vvrho(x, t ) & =  C_{\vrho}(\tau) \vrho( y , \tau),  \\
\uu(x, t) & =  C_u(\tau) \UU( y, \tau),   \\
\mbb(x,t) & =  C_{\mb}(\tau) \mb( y, \tau).  
\end{align}
\end{subequations}
Using the self-similar transform \eqref{eq:SS_ansatz} and the chain rule, we rewrite the Euler system 
\eqref{eq:sys_UB} equivalently as 
\begin{subequations}
\label{eq:sys}
\begin{align} 
\pa_\tau \vrho +  (\cx y + \UU) \cdot \nabla \vrho  + 2 \alpha \vrho \, \div \UU &=  \crho \vrho, 
\label{eq:sys:a} \\
\pa_\tau \UU + (\cx y + \UU) \cdot \nabla  \UU  + \tfrac{1}{2\alpha} \mb \nabla \vrho + \tfrac{1}{\gamma} \vrho \nabla \mb &= \cu  \UU, \\
 \pa_\tau \mb + (\cx y + \UU) \cdot \nabla \mb &= \mcb \mb 
 \label{eq:sys:c}.
 \end{align}
\end{subequations}

\begin{remark}[Speed of sound]
In original $(x,t)$ coordinates the speed of sound $c = \alpha \sigma$ relates to $\vvrho$ and $\mbb$ via $c^2 = \vvrho \mbb$. Using the self-similar transformation in~\eqref{eq:SS_ansatz} we may thus define the non-negative self-similar sound speed $\cc = \cc(y,\tau)$ via
\begin{equation}
\label{eq:cc:def}
\cc^2 := \vrho \mb,
\end{equation}
and note that $c(x,t) = \alpha \sigma(x,t) = C_u(\tau) \cc(y,\tau)$.
\end{remark}

\subsection{The globally self-similar solution}
Fix $d\in \{1,2,3\}$, $\alpha \in (0,d]$, and $\NNN \geq 1$.
The radially symmetric globally self-similar solution $(\bar u^r,\bar \sigma,\bar b)(r,t)$ of~\eqref{eq:euler2} constructed in Section~\ref{sec:profiles}, see Remark~\ref{eq:globally:SS:exact}, also defines a radially symmetric globally self-similar solution of~\eqref{eq:sys_UB} via
\begin{subequations}
\label{eq:globally:SS:exact:nonradial}
\begin{equation}
\cx = \cxbar , \qquad
\cu = \cubar = \cxbar -1, \qquad
\cbb = \cbbbar = 2 \cbbar, \qquad
\crho = \crhobar = 2 (\cubar  - \cbbar) = - \tfrac{2\alpha d}{1+\alpha d},
\label{eq:globally:SS:exact:nonradial:a}
\end{equation}
and
\begin{align}
\bar \vvrho (x,t) &:=  (-t)^{\crhobar} \bar \vrho\Bigl(\tfrac{x}{(-t)^{\cxbar}}\Bigr) ,  
&&\bar \vrho(y) := \bigl(\tfrac{\alpha \bar \Sigma}{\bar B}\bigr)^2(|y|) = \bigl(\tfrac{\alpha \bar Q}{\bar H}\bigr)^2 (|y|)
\label{eq:globally:SS:exact:nonradial:b}
\\
\bar \uu (x,t) &:= (-t)^{\cxbar -1}  \bar \UU\Bigl(\tfrac{x}{(-t)^{\cxbar}}\Bigr) ,  
&&\bar \UU(y) := \vec{e}_r \bar U(|y|) = y \bar V(|y|), 
\label{eq:globally:SS:exact:nonradial:c}
\\
\bar \mbb(x,t) &:= (-t)^{\cbbbar} \bar \mb \Bigl(\tfrac{x}{(-t)^{\cxbar}} \Bigr),  
&&\bar \mb(y) := \bar B^2 (|y|) = |y|^2 \bar H^2(|y|),
\label{eq:globally:SS:exact:nonradial:d} \\
\bar c(x,t) &:= (-t)^{\cxbar -1} \bar \cc \Bigl(\tfrac{x}{(-t)^{\cxbar}} \Bigr), 
&&\bar \cc(y) := \alpha \bar \Sigma (|y|) = \alpha |y|  \bar Q(|y|).
\label{eq:globally:SS:exact:nonradial:e}
\end{align}
\end{subequations} 
The solution $(\bar \vvrho,\bar \uu,\bar \mbb)$ defined above is smooth for times $t\in [-1,0)$, implodes at time $t(\infty) = 0$, and is radially symmetric.

\begin{remark}[\bf Properties of the globally self-similar profiles]\label{rem:decay}
We record a few useful properties of the stationary solution to~\eqref{eq:sys} defined by the profiles $(\bar \vrho, \bar \UU, \bar \mb)$ and the   modulation parameters $(\cxbar,\cubar,\cbbbar,\crhobar)$ given in \eqref{eq:globally:SS:exact:nonradial}. 
The profiles $(\bar \vrho, \bar \UU, \bar \mb) \in C^{\infty}(\Reals^d)$ satisfy the lower bounds
\begin{subequations}\label{eq:profile_van}
\begin{align}
  &\bvr(y) \approx 1, 
  \quad \bc(y), |\bu(y)| \approx |y|,  
  \quad \barb(y) \approx |y|^2,  &&\forall |y| \leq 1, \\
 &\bvr(y) > 0, \quad \bc(y)>0 \quad \barb(y) > 0,  &&\forall y \neq 0.
\end{align}
the upper bounds
\begin{equation}\label{eq:profile_decay}
   |\nabla^k \bvr (y) | \les_k \la y \ra^{\frac{\crhobar}{\cxbar} - k},
   \qquad 
   |\nabla^k \bu (y) |
   +
   |\nabla^k \bc (y) | \les_k \la y \ra^{\frac{\cubar}{\cxbar} - k},
   \qquad 
   |\nabla^k \barb (y) | \les_k \la y \ra^{\frac{\cbbbar}{\cxbar}  - k},
\end{equation}
for any integer $k \geq 0$. The implicit constants in the $\gtr$ and $\les$ symbols in~\eqref{eq:profile_van} are assumed to depend on the fixed values of $d$, $\alpha$, and $\NNN$. Moreover, cf.~\eqref{eq:Omega:invariant:c}, we have the outgoing condition  
\begin{equation}
     \bcr |y|+  \bar U(y)  - \bc(y)  \geq  C_{\eqref{eq:profile_outgoing}}   |y|, \qquad
      C_{\eqref{eq:profile_outgoing}} :=  \tfrac{1+ \frac 23 \alpha d}{4 \NNN (1+\alpha d)} 
      > \tf{1}{6 \NNN }.
    \label{eq:profile_outgoing}
\end{equation}
\end{subequations}
The scaling parameters satisfy 
\begin{subequations}\label{eq:profile_scal}
 \begin{equation} 
  \cubar = \cxbar - 1, 
  \qquad \bcb = 2 (\bcr + \uua ) = 2 \cbbar,
  \qquad \bcvr  = 2 \bcu - \bcb = 2\alpha d \uua, 
  \qquad \uua  = - \tfrac{1}{1 + \alpha d}, 
 \end{equation}
 so that by the second bullet in Lemma~\ref{lem:properties:exponents}, 
 \begin{equation}
\bcr > 0, \qquad
0 <   \tfrac{ \bcb}{2 \bcr}  < 1,  \qquad 
\bcu < \tfrac{1}{2} \bcb .
 \end{equation} 
\end{subequations}
For compactness of notation, it is convenient to denote 
\begin{equation}
\label{eq:ODE_main:kp}
  \kp = \bcr + \uua = \cbbar. 
\end{equation}
Lastly, from~\eqref{eq:Taylor:R=0:all} we deduce that  the $\NNN$-th profiles $(\bar \vrho, \bar \UU, \bar \mb)$ in this family satisfy
\begin{subequations}
\label{eq:ODE:main:profi}
\begin{align} 
&\bvr = \aaa + \aan  |y|^{2\NNN} + \OO_{|y|\to0}(|y|^{4\NNN}),    
 \qquad 
 \aaa = (\alpha \bar q_0)^2 = \tfrac{\alpha\gamma d}{2} \bar v_0^2, 
 \qquad  
 \aan = 2 \alpha^2 \bar q_0 (\bar q_{\NNN} -\bar q_0 \bar h_{\NNN}),  \\
  &\bu = \uua y + \uun   y |y|^{2\NNN} + \OO_{|y|\to0}(|y|^{4\NNN + 1}),   
  \qquad
  \uua = \bar v_0 = - \tfrac{1}{1+\alpha d}, 
  \qquad 
  \uun = \bar v_{\NNN}
 ,  \\
  &\barb =   \bb |y|^2 + 
  \bbn |y|^{ 2 \NNN + 2} + \OO_{|y|\to0}(|y|^{ 4 \NNN + 2}),  
\qquad \bb = 1,  
\qquad \bbn = 2 \bar h_{\NNN}, 
\end{align}
so that 
\begin{equation} 
 \na \bvr  = \OO_{|y|\to0}(|y|^{2\NNN-1}), \quad 
 \na \bu  = \uua \Id + \OO_{|y|\to0}(|y|^{2\NNN}), \quad
 \na \barb  =  2 \bb y + 
\OO_{|y|\to0}(|y|^{2 \NNN + 1}).
\end{equation}
\end{subequations}
\end{remark}

\subsection{Evolution of the perturbation}
Our goal is to analyze solutions $(\vrho,\UU,\mb)$ of~\eqref{eq:sys} which are close the the stationary profile $(\bvr,\bu,\barb)$ discussed in Remark~\ref{rem:decay}. For this purpose, we denote 
\[
(\tvr,  \tu, \tb)
:=
(\vrho-\bvr,\UU-\bu,\mb-\barb),
\]
and write the modulation functions as
\[
\cx = \cxbar + \cxtilde,\qquad
\cu = \cubar + \cutilde,\qquad
\cbb=\cbbbar + \cbbtilde,\qquad
\crho = \crhobar + \crhotilde.
\]
Note that due to~\eqref{eq:vanishing:order:at:origin:default} and~\eqref{eq:profile_van} we have
\begin{equation}
\label{eq:vanishing:order:at:origin:SS}
 \UU(0,\tau) = \tu(0,\tau) = 0 \,, \qquad 
 \mb(0,\tau) = \tb(0,\tau) = 0 \,, \qquad 
 \nabla \mb(0,\tau) = \nabla \tb(0,\tau) = 0 \,, 
\end{equation}
for all $\tau \geq 0$.

\subsubsection{Full linear-nonlinear decomposition}
It is useful to introduce the following compact notation.  Denote the full solution, the profile, and the perturbation as the five-dimensional vectors
\[
\WW := (\vrho, \UU, \mb) ,
\qquad
\bar \WW := (\bvr,\bu,\barb),
\qquad
\tw := \WW - \bar \WW
= (\tvr,  \tu, \tb).
\]  
In order to rewrite quadratic nonlinear terms in~\eqref{eq:sys}, for $\WW^{(i)}=(\vrho^{(i)},\UU^{(i)},\mb^{(i)})$, $i \in \{1,2\}$, we  introduce the following bilinear operators, which are associated to the nonlinearities in \eqref{eq:sys} 
\begin{subequations}\label{eq:bilin}
\begin{align}
  \cB_1(\WW^{(1)}, \WW^{(2)}) & = - \UU^{(1)} \cdot \nabla \vrho^{(2)} - 2 \alpha \vrho^{(1)} \div  \UU^{(2)}  ,  \\
  \cB_2(\WW^{(1)}, \WW^{(2)}) & = -  \UU^{(1)} \cdot \nabla  \UU^{(2)} - \tfrac{1}{2\alpha} \mb^{(1)} \nabla \vrho^{(2)} - \tfrac{1}{\gamma} \vrho^{(1)} \nabla \mb^{(2)}  , \\
  \cB_3(\WW^{(1)}, \WW^{(2)})& =  - \UU^{(1)} \cdot \nabla \mb^{(2)} . 
\end{align}
\end{subequations}
In particular, with the notation in~\eqref{eq:bilin} the derivative falls on the second input variable $\WW^{(2)}$. We let $\cB = (\cB_1,\cB_2,\cB_3)$.

Similarly, in order to write the linear terms arising from the linearization of~\eqref{eq:sys} in a more compact fashion, we denote $\cL = (\cL_1,\cL_2,\cL_3)$, where 
\begin{subequations}\label{eq:lin_ODE:cL}
\begin{align}
 \cL_1 \tw &:= - ( \bcr y + \bu) \cdot \nabla \tvr - 2 \alpha \bvr \, \div \tu  + \bcvr \tvr 
 -  \tu  \cdot \nabla \bvr - 2 \alpha \tvr \, \div \bu
  -   \tcr y  \cdot \nabla \bvr  + \tcvr  \bar \vrho , \label{eq:lin_ODE:a} \\
 \cL_2 \tw &:=  -  ( \bcr y + \bu) \cdot \nabla \tu  - \tfrac{1}{2\alpha} \barb \nabla \tvr
  - \tfrac{1}{\gamma} \bvr \nabla \tb    +  \bcu  \tu  -  \tu  \cdot \nabla \bu -  \tfrac{1}{2\alpha} \tb \nabla \bvr
  - \tfrac{1}{\gamma} \tvr \nabla \barb \notag \\
& \quad -  \tcr y  \cdot \nabla \bu + \tcu \bu ,  \label{eq:lin_ODE:b} \\
 \cL_3 \tw &:= - ( \bcr y + \bu) \cdot \nabla \tb   +  \bcb \tb 
 -  \tu \cdot \nabla \barb    - \tcr y  \cdot \nabla \barb  +  \tcb \barb . \label{eq:lin_ODE:c} 
\end{align}
\end{subequations}

With this notation, we obtain from~\eqref{eq:sys} that the equation satisfied by the perturbation $\tw = (\tvr, \tu, \tb)$ is 
\begin{equation}
\label{eq:lin_ODE}
\p_\tau \tw  
=  \cL \tw  
- \tcr y \cdot \na \tw  
+ (\tcvr \tvr, \tcu \tu, \tcb \tb ) 
+ \cB(\tw, \tw).
\end{equation}

\subsubsection{A partially-nonlinear decomposition}
The full linear-to-nonlinear decomposition in~\eqref{eq:lin_ODE} is not always the most practical way to write down the evolution equations for $\tw = (\tvr,\tu,\tb)$. An alternative form that can be derived from~\eqref{eq:sys} is 
\begin{subequations}\label{eq:lin:full}
\begin{align}
\pa_\tau \tvr &= -( \cx y + \UU) \! \cdot  \! \nabla \tvr - 2 \alpha \vrho \, \div  \tu   + \cvr \tvr 
- (  \tcr y + \tu)  \! \cdot  \!  \nabla \bvr - 2 \alpha \tvr \, \div \bu
+ \tcvr  \bar \vrho , 
\label{eq:lin:full:a} \\
\pa_\tau \tu &=  -  ( \cx y + \UU) \!  \cdot \! \nabla   \tu  - \tfrac{1}{2\alpha} \mb \nabla \tvr - \tfrac{1}{\gamma} \vrho \nabla \tb    
+ \cu  \tu  
- ( \tcr y + \tu ) \!  \cdot  \!  \nabla \bu 
-  \tfrac{1}{2\al} \tb \nabla \bvr - \tfrac{1}{\gamma} \tvr \nabla \barb 
+ \tcu \bu ,  
\label{eq:lin:full:b} \\
\pa_\tau \tb &=  - (\cx y + \UU)  \! \cdot  \! \nabla \tb   +  \cbb \tb 
 - (\tcr y + \tu) \! \cdot \! \nabla \barb  + \tcb \barb  . 
 \label{eq:lin:full:c}
 \end{align}
\end{subequations}
The equation set~\eqref{eq:lin:full} is different from 
\eqref{eq:lin_ODE} as it combines pieces of the linear and nonlinear parts in the transport term. 

\subsubsection{Constant density in the far field}

Recall from Theorem~\ref{thm:main:profiles}, item (iii), and~\eqref{eq:globally:SS:exact:nonradial:a}--\eqref{eq:globally:SS:exact:nonradial:b} that the density profile $\bvr(y)= (\frac{\alpha \bar Q}{\bar H})^2 (|y|)$ behaves asymptotically as $|y|^{\crhobar/\cxbar}$ as $|y|\to \infty$, with $\cxbar>0$ and $\crhobar = - \frac{2\alpha d}{1+\alpha d}<0$; in particular, $\bvr(y) \to 0$ as $|y|\to \infty$. 

We wish however to allow for \emph{non-decaying densities}, a physically more realistic setting; for the initial data, we should consider densities which converge to a nonzero-constant in the far field. For this purpose, instead of working with the spatially decaying self-similar profile $\bvr$, we introduce a time-dependent, spatially cutoff version of it, defined by
\begin{equation}
\label{eq:vrho_s}
\vrs(y,\tau) :=  \chi\Bigl( \tfrac{|y|}{\rs} \Bigr) \bvr(y) + \Bigl(1 - \chi\Bigl( \tfrac{|y|}{\rs} \Bigr) \Bigr) \bvr( \rs ), 
  \quad \mbox{where} \quad 
\rs = e^{ \bcr \tau } R_0,
\end{equation}
with $R_0\geq 1$ a sufficiently large constant to be chosen later (see~\eqref{eq:bar_eps_cond:4} and \eqref{eq:bar_eps_cond:6} below),
where $\chi : \Reals \to \Reals_+$ is a smooth non-negative cutoff function satisfying 
\begin{equation}\label{eq:cutoff}
  \chi( z ) = 1, \quad \forall |z| \leq 1, \qquad  \chi(z) = 0,\ \quad \forall |z| \geq 2,
  \qquad 
  |\chi'(z)|\leq 2 \ \quad \forall 1 < |z| < 2.
\end{equation}
The function $\vrs(y,\tau)$ serves as an approximation of $\vrho(y,\tau)$, which solves~\eqref{eq:sys}.
The initial data $\vrho_{\mathsf{in}} = \vrho(0,\cdot)$ for \eqref{eq:sys}  will be constructed (see Theorem~\ref{thm:nonrad:stab}) by perturbing $\vrs$, which does not decay for large $y$.  

\subsection{Main result: finite co-dimension non-radial stability}
The remainder of this section is dedicated to establishing the global existence and decay to equilibrium for the perturbation $(\tvr,\tu,\tb)$ which solves~\eqref{eq:lin:full}, or equivalently,~\eqref{eq:lin_ODE}.  A global existence and decay result for non-radial solutions of~\eqref{eq:lin:full} implies---by definition---the stability of the radially symmetric globally self-similar solution defined in~\eqref{eq:globally:SS:exact:nonradial}, to perturbations outside this symmetry class. The main result is Theorem~\ref{thm:nonrad:stab} below.

\subsubsection{Measuring Taylor coefficients at the origin}
As we have already seen in Section~\ref{sec:stability} for the radially-symmetric case, proving such a global existence result requires an a-priori control (globally in time $\tau$) of \emph{finitely many Taylor coefficients at $y=0$} for the functions $(\tvr,\tu,\tb)$. In Theorem~\ref{thm:stability:radial} we have shown that Taylor coefficients of sufficiently high index are under control~\eqref{eq:boot:2}, and that the modulation functions may be used to control the fundamental Taylor coefficients, see~\eqref{eq:boot:0} and~\eqref{eq:boot:1}. How about the ``in between'' Taylor coefficients, whose order is not high? In Theorem~\ref{thm:stability:radial} we have chosen to turn these coefficients ``off'' at the initial time, see assumption $(ii)$, and Lemma~\ref{lem:vanishing:at:R=0} then ensured these coefficients remain ``off'' for all $\tau>0$. We should thus expect that a similar difficulty must be faced when dealing with non-radial perturbations of the radial profiles, i.e., for solutions of~\eqref{eq:lin:full}. Controlling Taylor coefficients for $(\tvr,\tu,\tb)$ at $y=0$ whose index is not large, is quite complicated;  Section~\ref{sec:ODE} is dedicated solely to this global-in-time ODE analysis. 

In this section, the control of all Taylor coefficients of sufficiently high index is taken as an assumption; we do not need to assume that these coefficients vanish identically, we only need them to decay as $\tau \to \infty$. For this purpose, for any $\MMM \in \Naturals_0$, we define a norm that controls the modulation functions and the Taylor coefficients up to order $\MMM$:
\begin{align}
\label{eq:E:O:L}
E_{O, \MMM}(\tau) 
&:=  
|\tcr (\tau)| + |\tcu (\tau)| + |\tcb (\tau)| + |\tcvr (\tau)| 
\notag\\
&  \qquad
+
{\sum}_{0 \leq |\bfa| \leq \MMM}  |(\p^\bfa \tvr)(0,\tau)| 
+ 
 |(\p^\bfa \nabla \tu)(0,\tau)| 
+ 
 |(\p^\bfa \nabla^2 \tb)(0,\tau)| 
.   
\end{align}
A comment regarding~\eqref{eq:E:O:L} is in order: we 
recall from~\eqref{eq:vanishing:order:at:origin:SS} that $\tu(0,\tau) =0$, $\tb(0,\tau)=0$, and $\nabla \tb(0,\tau)=0$ for all $\tau \geq 0$; it is because of this fact that we have chosen $E_{O,\MMM}$ (cf.~\eqref{eq:E:O:L}) to measure the Taylor coefficients of order $\leq \MMM$ for $\nabla \tu$ and $\nabla^2\tb$ (instead of $\tu$ and $\tb$).

\subsubsection{Measuring the solution on the bulk}
Assuming that $E_{O, \MMM}$ is under control for a sufficiently large integer $\MMM$---see assumption~\eqref{eq:ass:EOM} below---our aim is show that the solution $(\tvr,\tu,\tb)$ of~\eqref{eq:lin:full} is global in time. In order to achieve this we need to estimate the solution in a carefully designed norm, which takes into account not just regularity (to take advantage of the energy structure of the system), but also (nearly) optimal asymptotic behavior as $|y|\to0$ and $|y|\to \infty$. For this purpose we introduce:
\begin{itemize}[leftmargin=1em]
\item Let $\phi: \Reals^d \to \Reals_+$ be a smooth non-negative cutoff function with 
\[
\phi(z) = 1, \quad \forall|z| \leq 1, \qquad 
\phi(z) = 0, \quad \forall |z| \geq 2.
\]
For an integer $k\geq 0$ and a function $f  \in C^k(\Reals^d)$, we define (compare with~\eqref{eq:Taylor:poly:1} and~\eqref{eq:J:M:def}) the cutoff Taylor polynomial $\tl_k f$ and the remainder $\rl_k f$ by\footnote{Note that as opposed to the definition of $\JJJ_k$ in~\eqref{eq:J:M:def} we do not divide $\rl_k f$ by any weight in $y$.}
\begin{subequations}
\label{eq:nonrad:Taylor}
\begin{align}
  \tl_k f(y) & := \phi(y) \Bigl( {\sum}_{ |\al | \leq k } \tfrac{1}{\al!} \pa^{\al} f(0) y^{\al} \Bigr), 
 \label{eq:nonrad:Taylor:a}\\
  \rl_k f(y) & := f(y) - \tl_k f(y).
  \label{eq:nonrad:Taylor:b}
\end{align}
Using these operators, we decompose functions $f$ into three parts: the profile $\bar f$, the cutoff Taylor polynomials $\tl_k \tilde f$, and the remainder $\rl_k \tilde f$, meaning that 
\begin{equation}
\label{eq:nonrad:Taylor:c}
f = \bar f + \tilde f = \bar f + \tl_k \tilde f + \rl_k \tilde f.
\end{equation}
\end{subequations}
\item The values of $k$ (in the above decomposition) that we need to choose for the various components of the vector $\tw = (\tvr,\tu,\tb)$ are different. Indeed, we recall from~\eqref{eq:vanishing:order:at:origin:SS} that $\tu(0,\tau) =0$, $\tb(0,\tau)=0$, and $\nabla \tb(0,\tau)=0$ for all $\tau \geq 0$, and using the fact that the profiles $\bar \WW$ satisfy~\eqref{eq:ODE:main:profi} we have 
that~\eqref{eq:sys} and~\eqref{eq:lin:full} preserve the vanishing orders:
 \begin{equation}\label{eq:vanish1}
 \vrho, \tvr = \OO_{|x|\to 0}(1), 
 \quad 
 \UU, \tu = \OO_{|x|\to 0}(|x|), 
 \quad 
 \mb , \tb= \OO_{|x|\to 0}(|x|^2).
 \end{equation}
Therefore, if $\tw \in C^{\mm + 3}$ for some integer $\MMM \geq 0$, we may define 
\begin{subequations}
\label{eq:pertb:main}
\begin{align}
\vrho &= \bvr + \tl_{\MMM} \tvr + \rl_{\MMM} \tvr,
&&\tvr_{\mm} := \rl_{\mm} \tvr,\\
\UU &= \bar \UU + \tl_{\MMM+1} \tu + \rl_{\MMM+1} \tu, 
&&\tu_{\mm} : = \rl_{\mm+1}  \tu, \\
\mb &= \barb + \tl_{\MMM+2} \tb + \rl_{\MMM+2} \tb, &&\tb_{\mm}: = \rl_{\mm+2} \tb,
\end{align}
and
\begin{equation}
\label{eq:pertb:main:tw}
\tw_{\mm} := ( \tvr_{\mm},  \tu_{\mm}, \tb_{\mm}). 
\end{equation}
\end{subequations}
Since $E_{O,\MMM}$ controls $(\tl_{\mm} \tvr, \tl_{\mm+1} \tu, \tl_{\mm+2}\tb)$, see~\eqref{eq:E:O:L} and~\eqref{eq:nonrad:Taylor}, the entire analysis in this section is dedicated to bounding $\tw_{\mm}$, for a carefully chosen $\mm$.

\item In order to bound $\tw_{\mm}$ we define a (weighted Sobolev) energy density by
\begin{subequations}
\label{energy:EE0}
\begin{equation}
E_k(\twm) = {\sum}_{|\bfa| = k} E_{\bfa}(\twm),
\end{equation}
with
\begin{equation}
\label{energy:EE0:alpha}
  E_{\bfa}(\twm) := \frac{ | \pa^{\bfa} \tvrm|^2}{ (2\al \vrho)^2}
  + \frac{ |\pa^{\bfa} \tum|^2}{\cc^2} 
  + \kp_{\mb} \frac{ | \pa^{\bfa} \tbm|^2}{ \mb^2 }
  + \frac{2}{\gamma}   \frac{ \pa^{\bfa} \tvrm}{ 2\al \vrho} \cdot \frac{ \pa^{\bfa} \tbm}{ \mb } , 
\end{equation}
where $\cc$ is the self-similar speed of sound (recall~\eqref{eq:cc:def}), and $\kp_{\mb} \geq 1$ is a sufficiently large constant, to be chosen later (see~\eqref{eq:def:kp_mb}).
Upon integrating the energy density $E_k(\twm)$ against a suitable weight in $y$, we define the $k^{\rm th}$ order  energy functional $E_k$ as
\begin{equation}
\label{energy:EE0:k}
\EE_k(\tau) := \int E_k(\twm)(y,\tau) \vp_0(y) |y|^{2 k} d y , \quad  \forall \, k \geq 0, 
\quad 
\end{equation}
for a  weight $\vp_0 >0 $ which satisfies~\eqref{eq:wgk_1} and is defined precisely in Lemma~\ref{lem:nonrad_wg}.

\item The total energy $\EE_{\msf{tot}}$ in which we measure $\twm$ is given by
\begin{equation}
\label{energy:EE0:tot}
\EE_{\msf{tot}} = \EE_\kk(\tau) + \nu \EE_0(\tau),  
\end{equation}
\end{subequations}
where the coupling parameter $\nu>0$ is to be determined later~\eqref{energy:tot}, 
and the regularity parameter $\kk$ is 
\begin{equation}
\label{eq:k_star}
      \kk = 2 d + 10.
\end{equation}
The above value of $\kk$ is large enough to close the nonlinear energy estimates, but it is not too large; in particular, it is \emph{explicit}.\footnote{This differs from the weighted $H^k$ stability analysis for implosion in the isentropic compressible Euler equations 
\cite{MRRS2022b,BCG2025,CCSV2024}, where $k$ is required to be \emph{sufficiently large} and is 
\emph{implicit}.}
\end{itemize}

\subsubsection{The stability result} 
We recall from \eqref{eq:E:O:L} that $E_{O, \mm}$ measures the perturbations in the modulation functions $(\cxtilde, \cutilde,\cbbtilde, \crhotilde)$, and the Taylor coefficients of order $\leq \MMM$ for $\tvr$, $\nabla \tu$, and $\nabla^2 \tb$, at $y=0$. Our main stability in this section is:
 
\begin{theorem}[\bf Finite co-dimension non-radial stability]
\label{thm:nonrad:stab}
There exists a sufficiently large integer $\mm$ such that the following holds.
Assume that there exists a decay rate $\lam > 0$ and a constant $C_{\eqref{eq:ass:EOM}} >0$,  such that 
\begin{equation}
\label{eq:ass:EOM} 
E_{O,\mm}(\tau) \leq C_{\eqref{eq:ass:EOM}} e^{-\lam \tau} E_{O, \mm}(0) , \qquad \forall \tau \geq 0.
\end{equation}
There exist 
\begin{itemize}[leftmargin=2em]
\item a sufficiently large radius $ \bar R_0 \geq 1$,  
\item sufficiently small parameters $\bar \eps, \nu, \delta_0 >0$,
\item a weight $\vp_0$ function  satisfying the bounds 
\[
C_1(\vp_0) \hat \vp_0(y)  \leq \vp_0(y) \leq C_2(\vp_0) \hat \vp_0(y),
\qquad \mbox{where}\qquad  \hat \vp_0(y) :=  |y|^{-2\mm -d} + \la y \ra^{- d - 4 \ddd} , 
\]
for some positive constants $C_2(\vp_0) > C_1(\vp_0) > 0$,
\end{itemize}
all depending on $d,\alpha,\NNN$, on the profiles $(\bar \vrho,\bar \UU,\bar \mb)$, and on $\mm,\lam, C_{\eqref{eq:ass:EOM}}$, such that the following holds. 

For any $R_0 > \bar R_0$, consider  initial data $(\vrho_{\iin}, \UU_{\iin}, \mb_{\iin})$ for \eqref{eq:sys} which satisfies
\begin{subequations}\label{ass:nonrad_init}
\begin{align} 
    \tfrac{1}{2} \vrs(y,0)   \la y \ra^{-\ddd}  & \leq \vrho_{\iin}( y)  \leq  2 \vrs(y,0)
  \la y \ra^{\ddd}      , \label{ass:nonrad_init:a} \\ 
 \quad  \tfrac{1}{2}  \barb(y)   \la y \ra^{-\ddd}  & \leq \mb_{\iin}(y) 
 \leq 2   \barb(y)   \la y \ra^{\ddd},  \label{ass:nonrad_init:b}
 \end{align} 
where $\vrs$ is defined in \eqref{eq:vrho_s} in terms of $R_0$, 
and such that the initial data $(\tilde \varrho_{\mathsf{in}}, \tilde \UU_{\mathsf{in}}, \tilde \mb_{\mathsf{in}}) := (\vrho_{\iin}, \UU_{\iin}, \mb_{\iin}) - (\bar \vrho,\bar \UU,\bar \mb)$ for~\eqref{eq:lin:full} satisfies
\begin{equation}
\label{eq:nonrad_init:small}
  E_{O, \mm}(0)  < \bar \eps,  \qquad \EEt(0) < \bar \eps^2,
\end{equation}
\end{subequations}
where the total energy $\EEt$ is  defined in \eqref{energy:EE0}, and depends on $\nu$ and $\vp_0$.

For any such initial data, the solution $(\tvr,\tu,\tb)$ of~\eqref{eq:lin:full} is global in time and satisfies the bounds
\begin{align*}
 \tfrac{1}{4} \vrs(y,\tau)  \la y \ra^{-\ddd} 
 &\leq \vrho(y,\tau)  
 \leq  4 \vrs(y,\tau)  \la y \ra^{\ddd}, \\
 \qquad
 \tfrac{1}{4} \barb(y)  \la y \ra^{-\ddd}  
 &\leq \mb(y,\tau)
 \leq 4  \barb(y)  \la y \ra^{\ddd}  , 
\end{align*}
For all $y\in\Reals^d$ and $\tau\geq 0$, 
together with the stability estimate 
\begin{equation}\label{eq:nonrad_decay}
   \EEt(\tau) \leq C \bar \eps^2 e^{- \lam_2 \tau},
\end{equation}
for some constants $\lam_2 , C >0$ depending on $\vp_0, \mm$, and on the profile. 
\end{theorem}

\begin{remark}[\bf Stability assumption of $E_{O,\mm}$]
The Taylor coefficients and the scaling parameters appearing in $E_{O,\mm}$ (recall~\eqref{eq:E:O:L})
form a \emph{finite dimensional and closed} system ODEs. In Section~\ref{sec:ODE}, we analyze this  ODE system and prove that~\eqref{eq:ass:EOM} indeed holds for initial data selected from a finite co-dimension set.  See Theorem~\ref{thm:ODE_stab}. 
We also refer to Remark \ref{rem:relation} for a discussion of the relation among the stability of the ODEs associated with \(E_{O,\mm}\), the stability estimates for the PDE away from the origin, and the local-in-time regularity of the solution.

\end{remark}

\begin{remark}
By requiring the parameter $\bar \eps$ in~\eqref{eq:nonrad_init:small} to be small enough, assumption \eqref{eq:ass:EOM} implies
\begin{equation}\label{eq:EOM}
 E_{O,\mm}(\tau) \leq C E_{O,\mm}(0) < C \bar \eps  < 1
\end{equation}
for all $\tau\geq 0$. The bound~\eqref{eq:EOM} is useful to simplify energy estimates.
\end{remark}

\begin{remark}[\bf The set for initial data is non-empty]
\label{rem:data:for:nonrad:stability}
We give examples of initial data $(\tilde \varrho_{\mathsf{in}}, \tilde \UU_{\mathsf{in}}, \tilde \mb_{\mathsf{in}})$ which satisfy the assumptions of Theorem~\ref{thm:nonrad:stab}. 
First, we may turn ``off'' the modulation functions, i.e., set $(\cxtilde, \cutilde,\cbbtilde, \crhotilde)=(0,0,0,0)$. Second, we note that for any non-negative integer $\MMM \geq 2$, if the initial perturbation is smooth enough at the origin and satisfies\footnote{We use the notation $\na^{\leq k} f(0)=0$ to indicate that $\partial^{\bfa} f(0) =0$ for all $|\bfa|\leq k$.} 
\begin{equation}\label{eq:initial:co-dim}
\na^{\leq \MMM} \tvr_{\mathsf{in}}(0) = 0 , \qquad  
\na^{\leq \MMM+1} \tu_{\mathsf{in}}(0) = 0,  \qquad 
\na^{\leq \MMM + 2} \tb_{\mathsf{in}}(0) = 0, 
\end{equation}
then we have 
$\na^{\leq \MMM} \tvr(0,\tau) = 0, \na^{\leq \MMM+1} \tu(0,\tau) =0,  \na^{\leq \MMM + 2} \tb(0,\tau)=0$
for all $\tau >0$. This fact is a direct consequence of the ODE analysis in Section \ref{sec:ODE}; that is to say, the vanishing order at the origin is preserved forward in time.  As a result, for initial data satisfying \eqref{eq:initial:co-dim}, cf.~\eqref{norm:ODE} we have that $E_{O, \mm}(\tau) \equiv 0$. Hence, the estimate \eqref{eq:ass:EOM} and assumption \eqref{eq:nonrad_init:small} on $E_{O, \mm}(0)$ hold trivially. 

Then, we note that the energy $\EEt(0)$  defined in \eqref{energy:EE0}  
\emph{only} depends on  $\tw_{\MMM}(0,\cdot)$, defined via~\eqref{eq:pertb:main} as $(\rl_{\mm} \tvr_{\mathsf{in}},\rl_{\mm+1} \tu_{\mathsf{in}},\rl_{\mm+2} \tb_{\mathsf{in}}) = (\tilde \varrho_{\mathsf{in}}, \tilde \UU_{\mathsf{in}}, \tilde \mb_{\mathsf{in}})$; the last equality holds due to~\eqref{eq:initial:co-dim}. The condition \eqref{eq:nonrad_init:small} on $\EEt(0)$ defines an open set $\mathsf{O}_1$ of initial data for $\tw_{\MMM}(0,\cdot)$ in an appropriate weighted Sobolev space. 
By the embedding estimates of Lemma~\ref{lem:embed}, condition \eqref{eq:nonrad_init:small} 
implies the smallness of $\na^2 \tw_{\MMM}$. Thus, for initial data satisfying \eqref{eq:initial:co-dim} and \eqref{eq:nonrad_init:small} with $\bar \eps$ small enough, the inequalities in \eqref{ass:nonrad_init:a}--\eqref{ass:nonrad_init:b} hold for $|y| \leq 2$. For $|y| >2$, since the cutoff Taylor polynomial $\tl_k f$ vanishes by definition \eqref{eq:nonrad:Taylor}, inequalities \eqref{ass:nonrad_init:a}--\eqref{ass:nonrad_init:b} define an open set $\mathsf{O}_2$ of initial data for $\tw_{\MMM}$ in a weighted $L^{\infty}$ space. 
For fixed $\bar \eps$, by choosing $R_0$ large enough, we obtain that  
the initial perturbation $( \vrs - \bvr, 0 ,0)$ 
associated with the initial data $(\vrs, \bu, \barb)$ belongs to these two open sets $\mathsf{O}_1$ and $\mathsf{O}_2$.
Therefore, for initial data of $\tw_{\MMM}$ in $\mathsf{O}_1 \cap \mathsf{O}_2$, and assuming that the first $\MMM$-th order Taylor coefficients satisfy~\eqref{eq:initial:co-dim}, the assumptions in Theorem \ref{thm:nonrad:stab} are satisfied. 
In particular, the set for initial data is non-empty. 

In Theorem \ref{thm:ODE_stab}, we show that a more general class of initial data leads to 
\eqref{eq:ass:EOM} and we relax \eqref{eq:initial:co-dim}.
\end{remark}

%%%%

%%%%
\begin{remark}[\bf Initial data may be constant outside a compact set]
\label{rem:compact:support}
At $\tau=0$, the cutoff self-similar profile $\vrs(\,\cdot\,,0)$ defined in~\eqref{eq:vrho_s}, satisfies $\vrs(y,0) = \bar \vrho(R_0)$ (a strictly positive constant) for all $|y| \geq 2 R_0$. Consequently, both the upper and lower bounds in~\eqref{ass:nonrad_init:a} on the initial density are themselves constant on $\{|y| \geq 2 R_0\}$, and the analogous upper/lower bounds in~\eqref{ass:nonrad_init:b} on the initial pseudo-entropy reduce, on the same set, to upper/lower bounds for $\barb(y) = |y|^2 \bar H^2(|y|)$. Moreover, the smallness condition $\EEt(0) < \bar \eps^2$ in~\eqref{eq:nonrad_init:small} only constrains the bulk remainder $\tw_{\MMM} = (\rl_{\mm}\tvr_{\sf in}, \rl_{\mm+1} \tu_{\sf in}, \rl_{\mm+2} \tb_{\sf in})$ defined in~\eqref{eq:pertb:main}, and does not constrain the values of $(\vrho_{\sf in},\UU_{\sf in},\mb_{\sf in})$ at any individual point $|y| \geq 2 R_0 \geq 2 $.

Therefore, the assumptions of Theorem~\ref{thm:nonrad:stab} are compatible with initial data that are \emph{identically equal to a constant non-vacuous state outside of a fixed compact set in $\Reals^d$}: density and pressure are constant, and velocity vanishes, on the complement of a ball. An explicit class of such initial data is obtained by smoothly cutting the asymptotic profile $(\bar \vrho,\bar \UU,\barb)$ of Theorem~\ref{thm:main:profiles} against the constant state $(\bar \vrho(R_0),0,\barb(R_0))$ at a radius $\sim R_0$, and then verifying~\eqref{ass:nonrad_init} and~\eqref{eq:nonrad_init:small} for $R_0$ sufficiently large; see also Remark~\ref{rem:data:for:nonrad:stability}. 
\end{remark}

The remainder of this section is dedicated to the proof of Theorem~\ref{thm:nonrad:stab}, whose summary is given in~\S~\ref{sec:proof:thm:nonrad:stab}.

\subsection{Estimates of the Taylor expansion}
We recall the definitions of the operator $\tl_k$ and $\rl_k = {\rm Id} - \tl_k$ from~\eqref{eq:nonrad:Taylor}. In this section we record a few properties of these operators, and discuss their relationship to the vanishing order of a given function at $y=0$.

\subsubsection{Functional spaces}

For any $k$ times continuously differentiable function $f$, we denote 
\begin{equation}
   | f(y)|_{\cC^k} = {\sum}_{ |\bfa| \leq k } | \pa_y^{\bfa} f(y)|.
 \label{def:cC_norm}
\end{equation}
Note: $\cC^k$ is not the $C^k$-H\"older norm; no supremum over $y$ is taken.

In order to measure the vanishing order of a function, for $y\neq 0$ we define
\begin{subequations}
\label{norm:gam}
\begin{equation}
 |f(y)|_{\Gam_{l}^k }   :=  {\sum}_{0\leq |\bfa|\leq k }  |y|^{|\bfa|-l} |\partial^\bfa  f (y) | ,   
 \qquad \mbox{and} \qquad
  \| f \|_{\Gam_{l}^k }  := {\sup}_{0<|y| \leq 2 }  | f(y)|_{\Gam_{l}^k} .
\end{equation}
The supremum in the above definition is taken over $|y|\leq 2$ to match the support of the function $\phi$ appearing in~\eqref{eq:nonrad:Taylor:a}, and hence, to match the support of $\tl_l f $.
Note that the above definition implies
\begin{equation}
 |\partial^{\bfa} f (y) | \leq |y|^{l - k} |f(y)|_{\Gam_{l}^{k}},
\qquad
\mbox{for all}
\qquad
|\bfa| = k.
\label{norm:gam:b}
 \end{equation}
Formally, a function $f \in \Gam_l^k$ has a vanishing order similar to $|x|^l$ as $|x|\to 0$, with bounded derivatives up to $k$-th order.

In order to measure the vanishing order of the triple $(f, g, h)$, e.g.~$\WW$, $\tw$, or $\twm$, we define 
\begin{align} 
  |(f, g, h)(y)|_{\Gam^k_{ (a, b, c)}} &:= | f(y)| _{\Gam^k_{ a}} +
  | g(y)| _{\Gam^k_{ b}} + |h(y)|_{\Gam^k_{ c}}, 
  \\
   \|(f, g, h) \|_{\Gam^k_{ (a, b, c)}}  &:= \sup_{|y| \leq 2 }  | (f,g,h)(y)|_{\Gam^k_{ (a,b,c)}} .
\end{align}
 \end{subequations}
When $\norm{f}_{\Gamma^k_l}<\infty$ or $\norm{(f,g,h)}_{\Gamma^k_{(a,b,c)}}<\infty$ we sometimes write $f \in \Gamma^k_l$, respectively $(f,g,h) \in \Gamma^k_{(a,b,c)}$.

\begin{lemma}\label{lem:TL_prod}
With the notation in~\eqref{norm:gam} and~\eqref{eq:nonrad:Taylor}, we have:
\begin{itemize}[leftmargin=2em]
\item[(a)] For any $i, p, q \geq 0$, $l \leq p+ q$, and $0< |y| \leq 2$,  we have 
\begin{subequations}\label{eq:TL_prod1}
\begin{equation}\label{eq:TL_prod1:a}
\bigl| (f g)(y) \bigr|_{ \Gam^i_{ l} }\les_{i,p,q,l} | f(y) |_{\Gam^i_{ p}}  \; | g(y) |_{\Gam^i_{ q}}.
\end{equation}
Moreover, for any $i, p, q \geq 0$ with $p \leq q$, and $0< |y| \leq 2$, we have 
\begin{equation}\label{eq:TL_prod1:b}
 |\na f(y)|_{ \Gam^i_{ p} } \les_{i,  p} |f(y)|_{\Gam^{i+1}_{ p+1}}, 
 \qquad
 |f(y)|_{\Gam^i_p}   \les |f(y)|_{\Gam^i_q} .
\end{equation}
\end{subequations}

 \item[(b)] For $f \in C^{k}(B_2(0))$ with $\|f\|_{\Gamma_{l}^{k}} < \infty$ and $k \geq l > j$, we have $\rl_j f = f$. In particular, for $f, g \in C^{k,\delta}(B_2(0))$ for some $\delta>0$, we have
\begin{subequations}\label{eq:TL_prod2}
\begin{equation}\label{eq:TL_prod2:a}
\rl_l ( f \cdot \rl_k g)  =  f \cdot \rl_k g ,
\quad \mbox{for any}\quad l \leq k ,  \\
\end{equation}
For any $l \leq k-1$, or for $l = k$ and $f(0) = 0$, or for $l = k+1$ and $f(0) = \nabla f(0) = 0$, we have 
\begin{equation}\label{eq:TL_prod2:b}
  \rl_l (f \cdot \na \rl_k g)  = f \cdot  \na \rl_k g ,
\end{equation}
for all $|y|\leq 2$.
\end{subequations}

\item[(c)]
For integers $i,k \geq 0$, any $f \in C^n (B_2(0))$ with $n=\max\{ k+1, i\}$, and any $0<|y|\leq 2$ we have 
\begin{subequations}\label{eq:TL_prod3}
  \begin{equation}\label{eq:TL_prod3:a}
  |\rl_k f(y)|_{\Gam^i_{k+1}} \les_{i, k} \| f \|_{ C^{n}(\overline{ B_{|y|}(0)  } ) },
  \quad n = \max \{ k+1, i\}.
  \end{equation}  
 Moreover, if $l \geq 0$ and $g$ is $C^l$ smooth near the origin, we have $\supp(f \cdot \tl_l g ) \subset B_2(0)$, and 
\begin{equation}\label{eq:TL_prod3:b}
   | \rl_k (  f \cdot \tl_l g ) |_{\Gam^i_{k+1}} 
   \les_{i, k}  \one_{|y|\leq 2}   \| f\|_{C^{n}(\overline{ B_2(0)} ) } |g(0)|_{\cC^l} ,
   \quad n = \max\{k+1,  i\}.
\end{equation}
\end{subequations}
\end{itemize}
\end{lemma}

\begin{proof}[Proof of Lemma~\ref{lem:TL_prod}]
Item~(a). Estimate \eqref{eq:TL_prod1} follows from the definition \eqref{norm:gam} and the  Leibniz  rule.
 
Item~(b). If $\| f \|_{\Gam_l^k} < \infty$ and $k \geq l > j$, by definition \eqref{norm:gam}, we have $\partial^{\bfa} f(0) =0$ for $|\alpha| \leq l-1$. Since $l -1 \geq j$, by definition \eqref{eq:nonrad:Taylor}, we obtain $ \tl_j f = 0$; hence $f= \rl_j f$.

For \eqref{eq:TL_prod2},  the regularity assumption and Taylor's theorem implies $\rl_k g(y) = \OO_{|y|\to 0}(|y|^{k+\delta})$. 
Thus, it is easy to verify $\partial^\bfa( f \cdot \rl_k g )(0) = 0$ for $|\bfa| \leq k$. When $l\leq k$ we thus have $\tl_l(f \cdot \rl_k g) = 0$, which implies~\eqref{eq:TL_prod2:a}. Analogously, when $l\leq k-1$ we have $\tl_l(f \cdot \nabla \rl_k g) = 0$. When $f(0)=0$, then $f \cdot \na \rl_k g = \OO_{|y|\to 0}(|y|^{k+\delta})$, and hence $\tl_k( f \cdot \na \rl_k g )=0$. The third case is similar, and \eqref{eq:TL_prod2} follows.
 
Item~(c). 
Estimate \eqref{eq:TL_prod3:a} follows from the integral form of the remainder in the Taylor expansion; 
when $j\leq k$ it implies that $|y|^{j-(k+1)} |\nabla^j \rl_k f(y)| \les \norm{f}_{C^{k+1}(\overline{B_{|y|}(0)})}$, while for $j\geq k+1$, we have $|y|^{j-(k+1)} |\nabla^j \rl_k f(y)| \les 2^{j-(k+1)}  \norm{f}_{C^{j}(\overline{B_{|y|}(0)})}$. Taking the sum over $j\leq i$, estimate~\eqref{eq:TL_prod3:a} now follows.
 
For \eqref{eq:TL_prod3:b}, since $ f \cdot \tl_l g$ is supported in  
$\overline{B_2(0)}$, we use \eqref{eq:TL_prod3:a}  
to prove
\[
     | \rl_k (  f \cdot \tl_l g ) |_{\Gam^i_{k+1}} 
    \les_{i,k} \| f \cdot \tl_l g \|_{ C^n(\overline{B_2(0)})}
  \les_{i,k} \| f \|_{C^n(\overline{B_2(0)})} \| \tl_l g \|_{ C^n(\overline{B_2(0)})}
  \les_{i, k} \| f \|_{C^n(\overline{B_2(0)})} |  g(0) |_{ \cC^l},
\]
completing the proof.
\end{proof}

\subsection{Derivations of the remainder equation}
Here we use \eqref{eq:lin:full} to derive the evolution equation for $\tw_{\mm}$ (which is defined in~\eqref{eq:pertb:main:tw}), up to error terms which we shall denote by 
$\cE_{\mm, 0}$ in $\rl_{\mm} \tvr$-equation, $\cE_{\mm, 1}$ in $\rl_{\mm+1} \tu$-equation, 
and $\cE_{\mm, 2}$ in $\rl_{\mm+2}\tb$-equation. See~\eqref{eq:lin}--\eqref{eq:error} below.  

We first make some preliminary remarks concerning vanishing orders at $y=0$. From \eqref{eq:vanish1} and the definition \eqref{eq:pertb:main}, we obtain 
\begin{subequations}
\label{eq:vanish2}
\begin{equation}  
 \bw  = (\bvr, \bu, \barb) \in \Gam^{\infty}_{(0,1,2)}, \qquad 
  \tw \in \Gam^{\mm + 3}_{( 0,1,2)} , \qquad 
  \tw_{\mm} \in \Gam^{\mm+3}_{(\mm+1,  \mm+ 2, \mm + 3)} , 
\end{equation}
with pointwise-in-$y$ bounds which follow from~\eqref{norm:gam:b}
\begin{align}
 |\na^i \tvr_{\mm}| &\les_{\mm, i} |y|^{\mm + 1 - i}   | \tvr_{\mm} |_{\Gam^i_{ \mm+1}},\\ 
    |\na^{i} \tu_{\mm}| &\les_{\mm, i} |y|^{\mm + 2 - i }  | \tu_{\mm}|_{\Gam^i_{ \mm+2}},\\ 
    |\na^{i} \tb_{\mm}| &\les_{\mm, i} |y|^{\mm + 3 - i }   | \tb_{\mm}|_{\Gam^i_{ \mm+3}} . 
\end{align}  
\end{subequations}
Note: in the remainder of the section we only use $\tw_{\mm} \in \Gam^{\mm+3}_{(\mm+1, \mm+2, \mm+3)}$ \textit{qualitatively} to derive the PDE that this vector satisfies (e.g., by using identities \eqref{eq:TL_prod2}); when aiming for \textit{quantitative} estimates, we use $\tw_{\mm} \in \Gam^{k}_{(\mm+1, \mm+2, \mm+3)}$ with $k < \mm$.

To measure the error terms in various Taylor expansions, we introduce a function class
which contains functions $f$ are similar to $|x|^{\mm+j+1}$ for $0<|x|\leq 2$; to make this precise, for an integer $j\geq 0$ we write
\begin{subequations}
\label{eq:error_func}
\begin{equation}
 f= \OO(\cE_{\mm, j})
\end{equation}
if $f\in C^{\MMM+j+1}(\overline{B}_{2}(0))$ and we have the pointwise bound
\begin{equation} 
  |f(y)|_{\Gam^i_{\mm+j +1}} 
  \les_{i, j, \mm}  \cR_{\mm, i}(y) , 
\end{equation}
where
\begin{equation} 
\cR_{\mm, i}(y)
:=E_{O , \mm}  \Bigl( 1+        |\tw_{\mm}(y)|_{ \Gam^i_{ (\mm+1, \mm+2, \mm+3)} } \Bigr).
\end{equation}
for all $y \in \overline{B}_2(0)$ and 
for any $ i \leq \mm+1$. 
\end{subequations}
Here we recall cf.~\eqref{eq:E:O:L} that
\begin{equation}\label{norm:ODE}
E_{O, \mm}  =    |\tvr(0)|_{\cC^{\mm}}
 + |\tu(0)|_{\cC^{\mm+1}}  + |\tb(0)|_{\cC^{\mm+2}} 
+ |\tcr | + |\tcu| + |\tcb| + |\tcvr| .   
\end{equation}

\subsubsection{Transport terms}
Consider the transport terms $f \cdot \na \td g$ in \eqref{eq:lin:full} with 
\[
f = \cx y + \UU.
\]
Using linearity of $\rlm$ and applying \eqref{eq:TL_prod2} in Lemma  \ref{lem:TL_prod},  we obtain
\begin{align} 
\label{eq:TL_van1}
\rlk (f \cdot \na \td g) 
& = \rlk (f \cdot \na \rlk \td g +  ( \bar f + \tl_{\mm+1} \td f + \rl_{\mm+1} \td f ) \cdot \na \tlk \td g   ) 
\notag \\
& =  f \cdot \na \rlk \td g + 
 \rl_k(\rl_{\mm+1} \td f  \cdot \na \tlk \td g)
+ \rlk(   ( \bar f + \tl_{\mm+1} \td f  ) \cdot \na \tlk \td g ) 
\notag\\
&=: I + II + III.
\end{align}
For $(k, \td g ) \in \{ (\mm, \tvr), (\mm+1, \tu), (\mm+2, \tb) \}$, 
applying the product rule \eqref{eq:TL_prod1}, using the vanishing conditions \eqref{eq:vanish2} for $\tb \in \Gam^i_2$, and recalling the definition of $\tu_{\mm}$ in~\eqref{eq:pertb:main}, we obtain 
\begin{subequations}\label{eq:TL_van2}
\begin{equation}
  |\rl_{\mm+1} \td f  \cdot \na \tlk \td g|_{\Gam^i_{k+1}}
\les_{i, \mm} |\rl_{\mm+1} \td f|_{\Gam^i_{\mm+2}} \cdot 
|\na \tl_k \td g|_{ \Gam^i_{\max(k- \mm-1, 0 ) } }
\les_{i, \mm} ( |\tcr| + | \tu_{\mm} |_{\Gam^i_{\mm+2}}  ) E_{O , \mm}
\end{equation}
for any $0 \leq i \leq \mm+3$. Thus, using \eqref{eq:TL_prod2}, we obtain 
\begin{equation}
II = \rl_k (\rl_{\mm+1} \td f  \cdot \na \tlk \td g) = \rl_{\mm+1} \td f  \cdot \na \tlk \td g, 
\qquad \mbox{and} \qquad 
| II |_{\Gam^i_{k+1}} \les ( |\tcr| + | \tu_{\mm} |_{\Gam^i_{\mm+2}}  ) E_{O , \mm}.
\end{equation}
\end{subequations}
For the term $III$ appearing in~\eqref{eq:TL_van1}, applying the product rule \eqref{eq:TL_prod3}, and using $\bar f , \tl_{k} g , \rl_{\mm+1} \td f \in C^{\mm+1}$, together with~\eqref{eq:EOM}, we obtain
\begin{equation}\label{eq:TL_van3}
  |III|_{\Gam^i_{\mm+1}} \les_{\mm, i}
  (1 + | \td f(0) |_{\cC^{\mm+1}})  | \td g(0)|_{\cC^k}
\les_{\mm, i}  (1 + E_{O,\mm}) E_{O , \mm}
\les_{\mm, i} E_{O , \mm}
\end{equation}
for $i \leq \mm + 1$.   Thus, using the notation for the error \eqref{eq:error_func} 
and  \eqref{eq:pertb:main}, for $f = \cx y + \UU$ we obtain 
\begin{subequations}\label{eq:TL_tran}
\begin{align} 
 \rlm ( f \cdot \na \tvr ) 
 & = f \cdot \na \rl_{\mm} \tvr  + \OO(\cE_{\mm, 0})
 =  f \cdot \na \tvr_{\mm} + \OO(\cE_{\mm, 0}) \\
 \rl_{\mm+1} ( f \cdot \na \tu )  
 &= f \cdot \na \rl_{\mm+1} \tu  + \OO(\cE_{\mm, 1})
 = f \cdot \na \tu_{\mm}  + \OO(\cE_{\mm, 1}),  \\
 \rl_{\mm+2} ( f \cdot \na \tb ) 
 & = f \cdot \na \rl_{\mm+2} \tb  + \OO(\cE_{\mm, 2})
 =  f \cdot \na  \tb_{\mm}  + \OO(\cE_{\mm, 2}). 
  \end{align}
\end{subequations}

\subsubsection{Other derivative terms}
Next, we consider other terms involving $\na (\tvr, \tu, \tb)$ in \eqref{eq:lin:full}. Using 
\eqref{eq:TL_prod2}, Lemma~\ref{lem:TL_prod}, and the vanishing order in \eqref{eq:vanish2}, 
$(\tvr, \tu, \tb), (\vrho, \UU, \mb) \in \Gam^{\mm+3}_{(0,1,2)}$, we obtain 
\begin{subequations}\label{eq:TL_other_deri}
\begin{align} 
  \rl_{\mm}(  \vrho \, \div \tu )
  & = \rl_{\mm} (  \vrho \, \div  \rl_{\mm+1} \tu 
  + (\bvr + \rl_{\mm} \tvr + \tl_{\mm} \tvr) \, \div  \tl_{\mm+1} \tu  )  
  \notag \\
  & =  \vrho \, \div \tu_{\mm}
  +   \tvr_{\mm}  \, \div \tl_{\mm+1} \tu
    + \rl_{\mm}( (\bvr + \tl_{\mm} \tvr) \, \div \tl_{\mm+1} \tu  ) 
  := I_1 + I_2 + I_3,  \\
\rl_{\mm+1}( \mb \na \tvr )
& = \rl_{\mm+1}( \mb \na \rl_{\mm} \tvr
+ ( \barb + \rl_{\mm+2} \tb + \tl_{\mm+2} \tb ) \na \tl_{\mm}\tvr ) 
\notag \\
& = \mb \na  \tvr_{\mm}  +  \tb_{\mm}    \na \tl_{\mm}\tvr
+ 
\rl_{\mm+1} (  ( \barb + \tl_{\mm+2} \tb ) \na \tl_{\mm}\tvr ) 
:=  II_1 + II_2 + II_3\\
\rl_{\mm+1}( \vrho \na \tb )
& = \rl_{\mm+1}\Big( \vrho \na  \rl_{\mm+2} \tb 
+ ( \bvr + \rl_{\mm} \tvr + \tl_{\mm} \tvr) \na  \tl_{\mm+2} \tb  \Big), 
\notag \\
& = \vrho \na  \tb_{\mm}  
+ \tvr_{\mm}  \cdot  \na  \tl_{\mm+2} \tb  
+  \rl_{\mm+1} ( ( \bvr +  \tl_{\mm} \tvr) \na  \tl_{\mm+2} \tb  ) 
:= III_1 + III_2 + III_3.
\end{align}
Note that from the vanishing order of $\tw, \tw_{\mm}$ \eqref{eq:vanish2}, it is easy to check the vanishing order
\[
I_1, I_2 = \OO(|x|^{\mm+1}), \quad II_1 , II_2, III_1, III_2 = \OO(|x|^{\mm+2}), 
\]
which along with the regularity of $\tw$ and \eqref{eq:TL_prod2}, implies the above identities. 

Applying the product rule \eqref{eq:TL_prod1} in Lemma \ref{lem:TL_prod}, 
for any $i\leq \mm+2$, we estimate $I_2, II_2, III_2$
\begin{align} 
   |I_2|_{\Gam_{\mm+1}^i} & \les_{\mm, i} 
  |  \tvr_{\mm} |_{\Gam_{\mm+1}^i} | \tu(0)|_{\cC^{\mm+1}}
  \les   |\tvr_{\mm}|_{\Gam_{\mm+1}^i} E_{O , \mm} \les \cR_{\mm, i}, \\
 |II_2 |_{\Gam_{\mm+2}^i} & \les_{\mm, i} 
  | \tb_{\mm} |_{\Gam_{\mm+3}^i} |\tvr(0)|_{\cC^{\mm}}
\les |\tb_{\mm}|_{\Gam_{\mm+3}^i} E_{O , \mm} \les \cR_{\mm, i} , \\
|III_2|_{\Gam_{\mm+ 2}^i} & \les_{\mm, i} 
  | \tvr_{\mm} |_{\Gam_{\mm+1}^i} 
  |  \na  \tl_{\mm+2} \tb  |_{ \Gam_1^i } 
\les_{\mm, i}
|\tvr_{\mm}|_{\Gam_{\mm+1}^i}
  |   \tb(0)  |_{ \cC^{\mm+2} } 
  \les |\tvr_{\mm}|_{\Gam_{\mm+1}^i} E_{O , \mm} \les \cR_{\mm, i}.
\end{align}
The terms $I_3, II_3, III_3$ only depend on the profile $\bw$ and $\na^i \tw(0)$ and have compact support. Applying the product rule \eqref{eq:TL_prod3} in Lemma \ref{lem:TL_prod}, using $\bw, \tl_{l} \tw  \in C^{\infty}$, and the bound~\eqref{eq:EOM}, we estimate 
\begin{align} 
  |I_3|_{\Gam_{\mm+1}^i} & \les_{\mm, i} (1 + |\tvr(0)|_{\cC^{\mm}}) |\tu(0)|_{\cC^{\mm+1}} 
  \les (1 + E_{O , \mm})  E_{O , \mm} \les \cR_{\mm, i},  \\
    |II_3|_{\Gam_{\mm+2}^i} & \les_{\mm, i} (1 + |\tb(0)|_{\cC^{\mm+2}}) |\tvr(0)|_{\cC^{\mm}} 
  \les (1 + E_{O , \mm})  E_{O , \mm} \les \cR_{\mm, i},  \\
  | III_3 |_{\Gam_{\mm+2}^i} & \les_{\mm, i}  (1 + |\tvr(0)|_{\cC^{\mm}}) |\tb(0)|_{\cC^{\mm+2}} 
    \les (1 + E_{O , \mm})  E_{O , \mm} \les \cR_{\mm, i} .
\end{align}
The above terms $I_2, I_3, II_2, II_3,  III_2, III_3$ have compact support. From \eqref{eq:error_func}, we obtain
\begin{equation}
  I_2, I_3 = \OO(\cE_{\mm, 0}),
  \quad II_2, II_3 = \OO(\cE_{\mm, 1}), \quad III_2, III_3 = \OO(\cE_{\mm, 1}).
\end{equation}
\end{subequations}

\subsubsection{Estimate of remaining terms}
The estimates of remaining terms in \eqref{eq:lin:full} are similar, and are obtained by applying Lemma \ref{lem:TL_prod}. We sketch the estimates. 
Using the notation \eqref{eq:pertb:main} and \eqref{eq:error_func}, we have 
\begin{subequations}\label{eq:TL_rem}
\begin{equation}
  \rl_{\mm} (\crho \tvr ) = \crho \tvr_{\mm},
  \quad   \rl_{\mm+1} (\cu \tu) = \cu  \tu_{\mm},
  \quad 
  \rl_{\mm+2} (\cbb \tb) = \cbb \tb_{\mm}.
\end{equation}
Since $\bw \in C^{\infty}$ and $\tl_l \tw \in C_c^{\infty}$ for any $l \geq 0$ such that $\tw \in C^l$ at the origin,  using \eqref{eq:vanish2} and Lemma \ref{lem:TL_prod}, we obtain
\begin{align}
 \rl_{\mm}(  -  \tu \cdot \na \bvr - 2 \al \tvr \, \div \bu )
& =   - \rl_{\mm} \Big( ( \rl_{\mm+1} \tu + \tl_{\mm+1} \tu ) \cdot \na \bvr
+ 2 \al ( \rl_{\mm}  \tvr + \tl_{\mm} \tvr ) \, \div \bu \Big) \notag \\
& = -  \tu_{\mm}  \cdot \na \bvr  
- 2 \al   \tvr_{\mm}  \, \div \bu  + \OO(\cE_{\mm, 0}),   \\
\rl_{\mm+1} \big( -  \tu  \cdot \na \bu -  \tfrac{1}{2\alpha} \tb \na \bvr 
  - \tfrac{1}{\gamma} \tvr \na \barb \big)  
  & = \rl_{\mm+1} \Big( -  ( \rl_{\mm+1} \tu + \tl_{\mm+1} \tu)  \cdot \na \bu  \notag \\
  & \qquad \qquad -  \tfrac{1}{2\alpha} 
  (\rl_{\mm+2}  \tb + \tl_{\mm+2} \tb) \na \bvr  
  - \tfrac{1}{\gamma} (\rl_{\mm}  \tvr + \tl_{\mm} \tvr) \na \barb \Big)  \notag \\ 
  & =   -   \tu_{\mm}   \cdot \na \bu  -  \tfrac{1}{2\alpha} \tb_{\mm}  \na \bvr  
  - \tfrac{1}{\gamma}  \tvr_{\mm}  \na \barb + \OO(\cE_{\mm, 1}) ,   \\
 - \rl_{\mm+2} ( \tu \cdot \na \barb )  
 & =  - \rl_{\mm+2} ( (\rl_{\mm+1} \tu + \tl_{\mm+1} \tu) \cdot \na \barb )  \notag \\
  & = - \tu_{\mm} \cdot \na \barb + \OO(\cE_{\mm, 2}). 
\end{align}

It remains to bound the terms on ${\sf RHS}_{\eqref{eq:lin:full}}$ which are a product of the perturbed modulation functions $(\cxtilde, \tcvr, \cutilde, \cbtilde)$ and the stationary profiles $(\bvr,\bu,\barb)$. 
For this purpose, we note that for any smooth function $f$ satisfying $f(y) >0$ for $y \neq 0$, $f(y) \asymp |y|^a$ for $|y| \leq 1$, $f(y) \asymp |y|^{-c}$ for $|y| \leq 1$, and $ |\na^i f(y)| \les_{i, c} \la y \ra^{-c -i}$ for any $ i\geq 0$, using \eqref{eq:TL_prod3} in Lemma \ref{lem:TL_prod} and the fact that $\rl_l(f)(y) = f(y)$ for any $|y| \geq 2$, we have 
\[
\begin{aligned} 
  |\na^i \rl_k f(y) | & \les_{i,k} |y|^{k+1 - i} 
  \les_{i,k} |y|^{k+1 - i - a} |f(y)|,  && |y| \leq 2, \\
    |\na^i \rl_k f (y) | & = |\na^i f (y)| \les_{i,c} \la y \ra^{-c - i}
  \les_{i, c} \la y \ra^{-i} |f(y)| , && |y| \geq 2.
  \end{aligned}
\]
Applying the above estimates to pairs $(f,k) \in\{ (\bvr,\mm), (\bu,\mm+1), (\barb,\mm+2)\}$, using the decay estimates of $\bw$ in Lemma \ref{rem:decay}, the vanishing conditions in \eqref{eq:vanish2}, and the asymptotic $|\bc| \asymp |x|$ near $x=0$, we obtain
\begin{align} 
  {\bvr}^{-1} \cdot  ( |  \na^i \rl_{\mm}( \tcvr \bvr )| 
+ | \tcr \na^i \rl_{\mm}( y \cdot \na \bvr )| ) \les  |y|^{\mm + 1 -i} \brak{y}^{-\mm-1} E_{O,\mm}, \\
      {\bc}^{-1} \cdot ( | \na^i \rl_{\mm+1}(  \tcu \bu )| +
 | \na^i \rl_{\mm+1}( \tcr y \cdot \na \bu )| ) \les  |y|^{\mm + 1 -i} \brak{y}^{-\mm-1} E_{O,\mm} , \\
        {\barb}^{-1} \cdot ( | \na^i \rl_{\mm+2}( \tcb \barb )| 
        + 
  | \na^i \rl_{\mm+2}( \tcr y \cdot \na \barb )| ) \les  |y|^{\mm + 1 -i} \brak{y}^{-\mm-1} E_{O,\mm} , 
  \end{align}
\end{subequations}
where we recall that $E_{O,\mm} \leq  \cR_{\mm, i}$ (for any $0\leq i \leq \mm+1$) is defined in \eqref{norm:ODE}.

\subsubsection{Summary of the derivations}
 
Combining \eqref{eq:TL_tran}, \eqref{eq:TL_other_deri}, \eqref{eq:TL_rem}, 
and using the notation in~\eqref{eq:pertb:main}, we derive the equations for $( \tvr_{\mm}, \tu_{\mm}, \tb_{\mm})$ 
by applying the linear operator $\rl_{\mm}$ to \eqref{eq:lin:full:a}, $\rl_{\mm+1}$ to 
\eqref{eq:lin:full:b}, and $\rl_{\mm+2}$ to \eqref{eq:lin:full:c}, resulting in
\begin{subequations}
\label{eq:lin}
\begin{align}
\pa_\tau \tvrm & = \cN_{\vrho} + \cE_{\mm, \vrho}, &&\cE_{\mm, \vrho} = \OO(\cE_{\mm, 0}), \\
\pa_\tau \tum &  = \cN_U +   \cE_{\mm, U} , &&\cE_{\mm, U} = \OO(\cE_{\mm, 1}),\\
 \pa_\tau \tbm & = \cN_{\mb} +  \cE_{\mm, \mb},  &&\cE_{\mm, \mb} = \OO(\cE_{\mm, 2}),
 \end{align}
 \end{subequations}
where the nonlinear terms $\cN_{\bullet}$ are defined as 
\begin{subequations}
\label{eq:non}
\begin{align}
 \cN_{\vrho} := & -( \cx y + \UU) \! \cdot \! \na \tvr_{\mm} - 2 \al \vrho \, \div \tu_{\mm}   + \cvr \tvr_{\mm}  -  \tu_{\mm} \! \cdot \! \na \bvr - 2 \al \tvr_{\mm} \, \div \bu , 
 \\
  \cN_U := &  -  ( \cx y + \UU) \! \cdot \! \na \tu_{\mm}  - \tfrac{1}{2\alpha} \mb \na \tvr_{\mm}
  - \tfrac{1}{\gamma} \vrho \na \tb_{\mm}   +  \cu  \tu_{\mm}   -  \tu_{\mm} \! \cdot \! \na \bu -  \tfrac{1}{2\alpha} \tb_{\mm} \na \bvr
  - \tfrac{1}{\gamma} \tvr_{\mm} \na \barb,
    \\
\cN_{\mb} := &  - ( \cx y + \UU) \! \cdot \! \na \tb_{\mm}  +  \cbb \tb_{\mm} 
 -  \tu_{\mm} \! \cdot \! \na \barb   ,
\end{align}
\end{subequations}
and the error terms $\cE_{\mm,\bullet}$ satisfy the estimates 
\begin{subequations}\label{eq:error}
\begin{align}
 |\na^i \cE_{\mm,\vrho}| 
& \les_{i, \mm} 
\bvr(y)  |y|^{\mm + 1 -i} \brak{y}^{-\mm-1} \cR_{\mm, i}(y)
, \\
 |\na^i \cE_{\mm,U}| 
& \les_{i, \mm} 
\bc(y) |y|^{\mm + 1 -i} \brak{y}^{-\mm-1} \cR_{\mm, i}(y)
, \\
 |\na^i \cE_{\mm,\mb}| 
& \les_{i, \mm} 
\barb(y) |y|^{\mm + 1 -i} \brak{y}^{-\mm-1}\cR_{\mm, i}(y)
,
\end{align}
\end{subequations}
for any $ i \leq \mm+3$. Here, $\cR_{\mm,i}$ is defined in \eqref{eq:error_func} and we used that $\brho(y)  \asymp 1, \bc(y) \asymp |y|, \barb(y) \asymp |y|^2$ for $|y| \leq 1$, 
and $\bvr, \bc, \barb > 0$  for all $|y| > 0$.

We note that the nonlinear terms in~\eqref{eq:non} may also be written in a more compact form by appealing to the bilinear operators $\cB = (\cB_1,\cB_2,\cB_3)$ defined earlier in~\eqref{eq:bilin}; indeed, we can rewrite  \eqref{eq:non} as 
\begin{equation}\label{eq:bilin:d}
(\cN_{\vrho}, \cN_U, \cN_{\mb})
 = - \cx y \cdot \na \twm 
 + ( \crho \tvrm,  \cu \tum, \cbb \tbm )
 + \cB(\WW, \twm) + \cB(\twm, \bar \WW),
\end{equation}
where $\WW = (\vrho, \UU, \mb)$, and $\twm = ( \tvr_{\mm}, \tu_{\mm}, \tb_{\mm})$.

\subsection{Bootstrap assumptions}
We first define the $\del_{\bullet}>0$ parameters which appear at various exponents 
\begin{subequations}
\begin{equation}\label{eq:para_del}
   \ddb := \tfrac{\bcb}{2 \bcr},
   \quad 
   \ddd := \tfrac{1}{100} \min \bigl( 1 - \tfrac{\cubar}{\cxbar}, 1 - \ddb  \bigr) , 
   \quad 
   \dda := \ddb + 8 \ddd  .
\end{equation}

Using~\eqref{eq:profile_scal},  we obtain
\begin{equation}
   \ddd > 0, \quad  0<  \ddb <  \dda < 1.
\end{equation}
\end{subequations}

Recall the modified profile $\vrs = \vrs(y,\tau)$ from \eqref{eq:vrho_s} and the profile $\barb=\barb(y)$ from \eqref{eq:globally:SS:exact:nonradial}. We introduce the upper (with lower index ${\sf u}$) and lower (with lower index ${\sf l}$) barrier functions $\vrl, \vru, \bbl, \bbu, \ccl, \ccu$ via
\begin{subequations}
\label{eq:boot1:old}
\begin{align}
& \vrl :=  \vrs \cdot \la y \ra^{-\ddd} , 
&& \vru :=  \vrs \cdot \la y \ra^{\ddd} ,  \\
& \bbl := \barb \cdot \la y \ra^{-\ddd}  ,   
&& \bbu := \barb \cdot \la y \ra^{ \ddd} , \\
&\ccl := ( \vrs \barb)^{1/2} \la y \ra^{-\ddd} , 
&&\ccu :=   ( \vrs \barb)^{1/2} \la y \ra^{\ddd}. 
\end{align}
\end{subequations}
Then, we impose the bootstrap assumptions
\begin{equation}
\label{eq:boot1} 
 \tfrac{1}{4} \vrl \leq \vrho  \leq  4 \vru, 
 \qquad
 \tfrac{1}{4} \bbl   \leq \mb   \leq 4  \bbu, 
\end{equation}
and 
\begin{subequations}
\label{eq:boot2}
\begin{align}
\vrho^{-1} |\na \vrho| 
& \leq \bvr^{-1} | \na \bvr| + | y |^{-1}  \la y \ra^{ 4 \ddd},\\
|\na (\UU - \bu)| 
& \leq  \la y \ra^{- 1 + \ddb + 4 \ddd}, \\
\mb^{-1} |\na \mb| 
& \leq \barb^{-1} |\na \barb|  +  | y |^{-1} \la y \ra^{ 4 \ddd},    
\end{align}
\end{subequations}
for all $y \neq 0$ and all $\tau \geq 0$.

\subsection{Consequence of bootstrap assumptions}
\subsubsection{Properties of barrier functions}
From Remark~\ref{rem:decay} and the definition \eqref{eq:vrho_s},  we obtain 
\begin{subequations}\label{eq:vrho_s_est}
\begin{equation}
  \vrs  \asymp  \la y \ra^{ \frac{\bcvr}{\bcr} } + \rs^{ \frac{\bcvr}{\bcr}}, 
  \qquad   
  \bvr \les \vrs,
  \qquad 
  |\na^k \vrs | \les_k  \la y \ra^{-k} \vrs , \\
\end{equation}
with implicit constants independent of $R_0, \rs$. Applying the Leibniz  rule and  Remark~\ref{rem:decay},
we obtain
\begin{equation}
|\na^k \barb | \les_k \la y \ra^{-k}(\barb + \chi), 
\qquad 
|\na^k (\vrs \barb )| \les_k \la y \ra^{-k} (\vrs \barb + \chi)
\end{equation}
\end{subequations}
for all $|y|>0$; here $\chi$ is the cutoff function from~\eqref{eq:cutoff}. 
For  $g \in \{ \vrs, \barb, (\vrs \barb)^{1/2} \}$ we have $g^2 + \chi >0 $ for any $y$. Using \eqref{eq:vrho_s_est}, and the  Leibniz  rule, for  any $a \in \Reals$, we obtain 
\[
   |\na^k  (g^2 \la y \ra^{2 a} + \chi)^{1/2}    |
\les_{k, a}  \la y \ra^{-k} (g^2 \la y \ra^{2 a} + \chi)^{1/2} ,
\qquad
g \in \{ \vrs, \barb, (\vrs \barb)^{1/2} \}.
\]
As a result, for the barrier functions $g \in \{ \vrl, \vru, \bbl, \bbu,\ccl, \ccu\}$ defined in \eqref{eq:boot1:old}, 
taking $a = \pm \ddd$ in the above estimates, we obtain 
\begin{subequations}\label{eq:weight_est}
\begin{equation}
  |\na^k ( g^2 + \chi )^{1/2}| \les_k   \ang y^{-k}  (g^2 + \chi)^{1/2},
  \quad g  \in \{   \vrl, \vru, \bbl, \bbu, \ccl, \ccu \}.
\end{equation}
Moreover, for $k\geq 0$, we have 
\begin{equation}\label{eq:weight_est:b}
 |\na^k g | \les_k |y|^{-k} g, \qquad g  \in \{ \vrl, \vru, \bbl, \bbu, \ccl, \ccu \}.
\end{equation}
\end{subequations}

\subsubsection{Bounds related to the sound speed}
Since $\cc = (\vrho \mb)^{1/2}$,  using assumptions \eqref{eq:boot1}, \eqref{eq:boot2}, we obtain 
\begin{subequations}\label{eq:boot1_res1}
\begin{align} 
&\cc^{-1} |\na \cc| 
\leq \tfrac 12 \bvr^{-1} | \na \bvr| + \tfrac 12 \barb^{-1} |\na \barb| + |y|^{-1} \la y \ra^{4 \ddd } 
\les |y|^{-1} \la y \ra^{4 \ddd } , 
\label{eq:boot1_res1:a} \\
\quad 
&\tfrac{1}{4} \ccl  \leq \cc \leq  4 \ccu
 ,
\label{eq:boot1_res1:b}
\end{align}for all $y\neq 0$. The implicit constant appearing on ${\sf RHS}_{\eqref{eq:boot1_res1:a}}$ depends only on the profiles.

Since $\vrs \les 1, \barb \les \la y \ra^{\frac{\bcb}{\bcr}} = \la y \ra^{2 \ddb}$ (see \eqref{eq:para_del}), $\ccu(0) = 0$, and $\dda < 1$, we further obtain 
\begin{equation}\label{eq:boot1_res1:c}
  \cc \leq 4 \ccu \les   |y| \la y \ra^{\ddb + \ddd - 1} 
  \les  |y| \la y \ra^{\dda - 4 \ddd- 1}
  \les \ang y^{\dda- 4 \ddd}.
\end{equation}
\end{subequations}
Using that $\cc = (\vrho \mb)^{1/2}$, the above estimate, \eqref{eq:boot1}, and \eqref{eq:vrho_s_est}, imply that
\begin{subequations}\label{eq:Wbar_upper}
\begin{align}
  &\bvr \les \vrs \les \vrho \la y \ra^{\ddd}, 
   && \bc = (\bvr \barb)^{1/2}  \les (\vrs \barb )^{1/2}  \les \cc \la y \ra^{\ddd},
   && \barb \les \mb \la y \ra^{\ddd} , \\
& \vru \les \vrho \la y \ra^{2 \ddd},
  && \ccu \les  \cc  \la y \ra^{2 \ddd},   
  && \bbu \les \mb \la y \ra^{2 \ddd}.
  \end{align}
\end{subequations}

\subsubsection{Relative decay of quadratic nonlinearities}
Under the bootstrap assumptions \eqref{eq:boot2}, using bounds on $\bu$ from Remark~\ref{rem:decay}, the fact that 
$\UU(0) = 0$ and $\bcu \leq \frac{1}{2}\bcb$,  we deduce
\begin{subequations}\label{eq:boot_decay}
\begin{equation}\label{eq:boot_decay:a}
|y|^{-1} |\UU| +  |\na \UU|  
\les \ang y^{-1 + \frac{\bcu}{\bcr}} + \ang y^{-1 + \ddb + 4\ddd} 
\les \ang y^{-1 +  \ddb + 4\ddd}  = \ang y^{-1 + \dda - 4 \ddd} .
\end{equation}
Under the bootstrap assumptions \eqref{eq:boot1}, \eqref{eq:boot2}, Remark~\ref{rem:decay}, \eqref{eq:boot1_res1:c}, \eqref{eq:boot_decay:a}, and $\cc = (\vrho \mb)^{1/2}$, we obtain 
\begin{align}
\cc^{-1} | \UU \cdot \na \cc |
+ \mb^{-1} | \UU \cdot \na \mb | 
+ \vrho^{-1} | \UU \cdot \na \vrho |
& \les |\UU | ( |y|^{-1} + |y|^{-1} \la  y \ra^{ 4 \ddd}  )
\les \la y \ra^{-1 + \dda} , \\
\vrho^{-1}  |\cc \na \vrho|
+ \mb^{-1} |\cc \na \mb|
& \les  \cc  |y|^{-1} \la y \ra^{ 4 \ddd}
\les \la  y \ra^{-1 + \dda  },
\label{eq:boot_decay:c}
\end{align}
\end{subequations}
where the implicit constants only depend on the profiles.

\subsubsection{Weights}
We recall that $k_* = 2d+10$ was fixed earlier in~\eqref{eq:k_star}.
Henceforth, we require the parameter $\mm$ (appearing for instance in the definition of $\twm$) to satisfy
\begin{equation}\label{eq:mm_low}
 \mm \geq k_*  + 2 .
\end{equation}
The precise value of $\mm$ is determined later in \eqref{def:nonrad_m}.

We recall that the weighted $H^k$ energy $\EE_k$ defined in~\eqref{energy:EE0:k} contains the weight function
\begin{equation}\label{eq:wgk_1:old}
\vp_0(y) |y|^{2k}  =: \vp_k(y).
\end{equation}
We impose the following constraints on the weight $\vp_0$, which is defined precisely in Section \ref{sec:choose_wg}: 
\begin{subequations}\label{eq:wgk_1}
\begin{align}
 C_1(\vp_0) \hat \vp_0
& \leq \vp_0 \leq C_2(\vp_0) \hat \vp_0,
\qquad \hat \vp_0 :=  |y|^{ - 2 \mm - d} + \la y \ra^{-d - 4 \ddd} ,  \label{eq:wgk_1:a} \\
  |\nabla \vp_0| & \leq C_3(\vp_0) \cdot |\vp_0| \cdot |y|^{-1},   \label{eq:wgk_1:b} 
\end{align}
\end{subequations}
for any $|y|>0$, where $C_i(\vp_0)>0$ are three constants which depend on $\vp_0$.

In terms of the weight $\hat \vp_0$, it is convenient to define an auxiliary norm for the vector $\tw = (\tvr, \tu, \tb)$:
\begin{equation}\label{norm:X}
  \| \tw \|_{\cX^l}^2 
  := \sum_{0\leq i\leq  l } \bigl \|  \hat \vp_0^{1/2}   |y|^i ( {\vru}^{-1} \na^i \tvr, \ccu^{-1} \na^i \tu, \bbu^{-1} \na^i \tb ) \bigr \|_{L^2}^2 .
\end{equation}

\subsubsection{Pointwise estimates}
Under the bootstrap assumption \eqref{eq:boot1}, and its consequence~\eqref{eq:boot1_res1:b}, since $\vrs \les 1$ and $\barb \les 1$ for $|y| \leq 1$, 
using the definition of $\hat \vp_0$ in \eqref{eq:wgk_1}, for any $i \leq k_* \leq \MMM$, we obtain 
\[
\hat \vp_0 |y|^{2 i} g^{-2} 
\gtr \hat \vp_0 |y|^{2 i} g_{\msf{u}}^{-2} 
\gtr \la y \ra^{- d - 4 \ddd + 2 i } (g_{\msf{u}}^2 + \chi)^{-1}, 
\qquad  g \in \{ \vrho, \cc, \mb \}.
\]

From \eqref{eq:weight_est}, we obtain that the weight $h =(g_{\msf{u}}^2 + \chi)^{ \frac{1}{2} }$, with $g_{\msf{u}} \in \{ \vru, \bbu, \ccu\}$, satisfies $|\nabla^i h|\les_i h \langle y \rangle^{-1}$, which is precisely assumption \eqref{ass:lem:embed} of Lemma~\ref{lem:embed}. Thus, 
for $k \leq (\kk - 1) - d $, applying Lemma~\ref{lem:embed} and recalling~\eqref{eq:mm_low} and~\eqref{eq:wgk_1}--\eqref{norm:X}, we obtain
\begin{subequations}
\label{eq:boot_point}
\begin{align} 
 (\bbu^2 + \chi)^{-1/2}  |\na^k \tbm| 
 &\les_{\kk} 
 \la y \ra^{ 2 \ddd  - k} 
  \| \tbm \cdot   \la y \ra^{  \frac{ - d -  4\ddd }{2}}   ( \bbu^2 + \chi)^{-1/2} \|_{L^2} 
  \notag\\
  &\qquad
 +\la y \ra^{ 2 \ddd  - k} 
  \| \na^{\kk-1} \tbm \cdot  \la y \ra^{  \frac{ - d -  4\ddd }{2} + \kk-1}    ( \bbu^2 + \chi)^{-1/2} \|_{L^2} 
  \notag \\
 & \les_{\kk} 
  \la y \ra^{ 2 \ddd  - k}    
 \Bigl( \|  \tbm \cdot  \hat \vp_0^{1/2}    \bbu^{-1} \|_{L^2}  
  + \| \na^{\kk-1} \tbm \cdot   \hat \vp_0^{1/2} |y|^{\kk-1}   \bbu^{-1} \|_{L^2}  \Bigr)
  \notag \\
  & \les_{\kk} \la y \ra^{ 2 \ddd  - k}  \| \twm\|_{\cX^{\kk -1}} .
     \label{eq:boot_point:a}
  \end{align}
Similarly, we obtain 
\begin{equation}
   ({\vru}^2 + \chi)^{-1/2}  |\na^k \tvrm| 
   +    (\ccu^2 + \chi)^{-1/2}  |\na^k \tum|  \les_{\kk} \la y \ra^{ 2 \ddd  - k}  \| \twm\|_{\cX^{\kk-1}} .
   \label{eq:boot_point:b}
\end{equation}
In particular, the bounds in~\eqref{eq:boot_point:a}--\eqref{eq:boot_point:b} show that $\twm(y) \in C^{\kk -1-d}$ near $y=0$, with $C^{\kk-1-d}$-norm (in the unit ball around the origin) bounded by $\| \twm\|_{\cX^{\kk -1}}$.
Since $\mm > \kk - 1 -d $ by \eqref{eq:mm_low}, it follows that  $\na^i \twm(0) =0$ for all $ i \leq \kk - d - 2$; thus,  we obtain
\begin{equation}
|\na^k \twm (y)| \les_{\kk} |y|^2 \max_{|y|\leq 1 } | \twm|_{\cC^{k+2}} 
\les_{\kk} |y|^2 \| \twm \|_{\cX^{\kk -1}} , \quad \forall |y| \leq 1.
\end{equation}
for all $|y|\leq 1$ and all $k \leq \kk-d-4$.
\end{subequations}
We deduce that for any $k$ with $0\leq k \leq  \kk/2 + 1$ (thus, $k \leq \kk -d-4$ cf.~\eqref{eq:k_star}), by combining estimates \eqref{eq:boot_point},  \eqref{eq:boot1}, \eqref{eq:boot1_res1:b}, \eqref{eq:Wbar_upper}, and using that $\vru \gtr 1, \ccu \gtr |y|, \bbu \gtr |y|^2 $ for $|y| \leq 2$, we arrive at  
\begin{subequations}\label{eq:boot_res2}
\begin{align}
|\na^k \tvrm (y)| &\les_{k_*} \vru(y) \la y \ra^{-k + 2\ddd}  \| \twm \|_{\cX^{k_*-1}}
\les_{k_*} \vrho(y) \la y \ra^{-k +  4 \ddd  }  \| \twm \|_{\cX^{k_*-1}}  , \\
 |\na^k \tum (y)| & \les_{k_*} \ccu(y) \la y \ra^{-k + 2 \ddd}  \| \twm \|_{\cX^{k_*-1}}
  \les_{k_*} \cc(y) \la y \ra^{-k +  4 \ddd  }  \| \twm \|_{\cX^{k_*-1}}
  , \\
 |\na^k \tbm (y)| & \les_{k_*} \bbu(y) \la y \ra^{-k + 2 \ddd}  \| \twm \|_{\cX^{k_*-1}} 
  \les_{k_*} \mb(y) \la y \ra^{-k + 4 \ddd  }  \| \twm \|_{\cX^{k_*-1}},
\end{align}  
\end{subequations}
for all $|y|>0$. In the above displays we have suppressed the time dependence relevant functions.
 
We recall from~\eqref{eq:nonrad:Taylor:c} that any smooth function $g$ may be decomposed as $g = \bar g + \tl_l \td g + \rl_l \td g$; see also~\eqref{eq:pertb:main}. 
For any $ k \geq 0$,  using Remark~\ref{rem:decay},  
and the upper bounds  for $(\bvr, \bc, \barb)$ in~\eqref{eq:Wbar_upper}, we obtain 
\begin{subequations}
\label{eq:boot_res3}
\begin{align} 
|\na^k \bvr|   
& \les_{ k} 
\one_{|y| \leq 1} +  \bvr \la y \ra^{-k} 
\les_{k } 
\bvr | y |^{- k} \les_k \vrho  | y |^{- k} \la y \ra^{ \ddd}, 
\label{eq:boot_res3:a}
\\
|\na^k \bu |
& \les_{ k}   |y|^{1 -k}  \one_{|y| \leq 1}+ \one_{|y| \geq 1}  \bc \la y \ra^{-k} 
\les_{ k }  \bc | y |^{- k} 
\les  \cc | y |^{- k} \la y \ra^{ \ddd} ,
\label{eq:boot_res3:b}
 \\ 
|\na^k \barb |
&    \les_{  k}  |y|^{2 -k}  \one_{|y| \leq 1} + \one_{|y| \geq 1}  \barb \la y \ra^{-k} 
\les_{ k } \barb | y |^{- k}   
\les_k 
 \mb  | y |^{- k} \la y \ra^{ \ddd}. 
 \label{eq:boot_res3:c}
\end{align}
Similarly, since $\tl_\mm$ applied to any function is smooth and has support in $B_2(0)$, we obtain 
\begin{align} 
 |\na^k  \tl_\mm \tvr|
& \les_{ k, \mm} 
  \one_{|y| \leq 2}  E_{O, \mm}
\les_{k, \mm }  \bvr | y |^{- k}   E_{O, \mm}  
\les_{k, \mm }  \vrho | y |^{- k} \la y \ra^{  \ddd}   E_{O, \mm} , 
\label{eq:boot_res3:d}
\\
 |\na^k  \tl_{\mm+1} \tu |
 & \les_{ k, \mm}   |y|^{1 -k}  \one_{|y| \leq 2} 
   E_{O,\mm }
\les_{ k, \mm }  \bc | y |^{-k}   E_{O, \mm} 
\les_{ k, \mm }  \cc | y |^{-k}  \la y \ra^{  \ddd}    E_{O, \mm}  ,
\label{eq:boot_res3:e}
 \\ 
  |\na^k   \tl_{\mm+2} \tb  | 
&    \les_{  k, \mm}   |y|^{2 -k}  \one_{|y| \leq 2}   E_{O, \mm} 
\les_{ k, \mm } \barb | y |^{- k}       E_{O, \mm}  
\les_{k, \mm} \mb | y |^{- k}  \la y \ra^{  \ddd}      E_{O, \mm}  . 
\label{eq:boot_res3:f}
\end{align}
\end{subequations}

The final pointwise estimate is for $\td \cc := \cc  - \bc$. Using  the estimates \eqref{eq:boot_res2}, \eqref{eq:boot_res3} with $k = 0$, and the definition $\cc^2 = \vrho \mb$, we bound 
\begin{align*}
|\td \cc | 
= \frac{| \cc^2 - \bc^2|}{ \cc + \bc  }
= \frac{| \mb \vrho - \barb \bvr | }{\cc + \bc} 
&\leq   \frac{ |\mb - \barb| \vrho + \barb | \vrho - \bvr|  }{\cc + \bc} 
\notag\\
&\les ( C_{\mm} E_{O, \mm} + \| \twm \|_{\cX^{\kk-1}} )  \frac{ \mb \vrho \la y \ra^{4 \ddd} 
+ \barb \vrho \la y \ra^{4 \ddd}  }{\cc + \bc}  .
\end{align*} 
The implicit constant in the above bound is independent of $\mm$; to emphasize that the implicit constant in~\eqref{eq:boot_res3:d}--\eqref{eq:boot_res3:f} does depend on $\mm$, we have included the constant $C_\mm$ in the above estimate.
Using \eqref{eq:boot1}, \eqref{eq:boot1_res1}, and $\cc^2 = \vrho \mb$, we obtain 
\[
  \mb \vrho \la y \ra^{4 \ddd} 
+ \barb \vrho \la y \ra^{4 \ddd}  
\les \mb \vrho \la y \ra^{5 \ddd}
=  \cc^2 \la y \ra^{5 \ddd}
\les \cc |y| \la y \ra^{ \ddb + 6 \ddd - 1 }
\les \cc |y| \ang y^{\dda - 1}.
\]
Combining the above two estimates, we derive the pointwise estimate
\begin{equation}
\label{eq:boot_est_tc}
  |\td \cc(y)| 
  \les  \bigl( C_M E_{O, \mm} + \| \twm \|_{\cX^{\kk-1}} \bigr) |y|  \ang y^{\dda - 1} ,
\end{equation}
with implicit constant that depends on $\kk$, but not on $\mm$.

\subsection{Weighted \texorpdfstring{$H^k$}{H k} estimates}
Our goal is to  perform weighted $H^k$ estimates on the evolution equation \eqref{eq:lin}, for $0\leq k \leq \kk$.  While the system \eqref{eq:lin} may be diagonalized, we do not do so to simplify the derivations. For convenience of the reader, we recall from~\eqref{energy:EE0} the definitions of the relevant energy norms. 
For a multi-index $\al$ with $|\al | =k$ and a large coupling parameter $\kp_{\mb} > \gamma^{-2}$, to be determined, we consider
\begin{equation}\label{energy:Hk}
  E_{\al}(\twm) 
  := \frac{ | \pa^{\al} \tvrm|^2}{ (2\al \vrho)^2}
  +  \frac{ |\pa^{\al} \tum|^2}{\cc^2} 
  +\kp_{\mb} \frac{ | \pa^{\al} \tbm|^2}{ \mb^2 }
  +  \frac{2}{\gamma}    \frac{ \pa^{\al} \tvrm}{ 2\al \vrho} \cdot \frac{ \pa^{\al} \tbm}{ \mb } , 
  \quad 
  E_k(\twm) = \sum_{|\al| = k} E_{\al}(\twm) .
\end{equation}
The associated energy norm $\EE_k$ is defined by integrating the above against a time-independent weight $\vp_k$, which is to be defined precisely later:
\begin{equation}\label{energy:Hk_norm}
   \EE_k = \int E_k(\twm) \vp_k d y.
\end{equation}

\subsubsection{Choosing the coupling parameter $\kp_{\mb}$ }
Note that although the fourth term appearing in the definition of~$E_{\al}(\twm)$ is not signed, the total energy density $E_{\al}(\twm)$ is signed (positive). To see this, note that for $\gamma > 1$ we have
\[
 \frac{2}{\gamma}   \left|   \frac{\pa^{\al} \tvrm}{ 2\al \vrho}    \frac{\pa^{\al} \tbm}{\mb} \right|
  \leq  2 \left| \frac{ \pa^{\al} \tvrm}{ 2\al \vrho}   \frac{\pa^{\al} \tbm}{\mb}\right|
\leq \kp_{\mb}^{-1/2}
\left(\kp_{\mb}  \frac{|\pa^{\al} \tbm|^2}{ \mb^2}  +   \frac{ | \pa^{\al} \tvrm|^2}{ (2\al \vrho)^2} \right).
\]
Upon choosing $\kp_{\mb}\geq 100 $, we obtain
\begin{equation}\label{eq:kp_large}
   \kp_{\mb}^{-1/2} < \tfrac{1}{2}
   \quad 
   \Longrightarrow
   \quad
   (1 + \kp_{\mb}^{-1} ) (1 - \kp_{\mb}^{-1/2})^{-1} < 1 + 3 \kp_{\mb}^{-1/2}
,
\end{equation}
and therefore the first inequality in~\eqref{eq:kp_large} implies 
\begin{align}
&(\underbrace{1-\kp_{\mb}^{-1/2}}_{\geq 1/2}) \Big(
   \frac{| \pa^{\al} \tvrm|^2}{ (2\al \vrho)^2}+  \frac{ | \pa^{\al} \tum|^2}{\cc^2} 
  +\kp_{\mb} \frac{ | \pa^{\al} \tbm|^2}{ \mb^2 } \Big)  
\notag\\
&  \leq  E_{\al}(\twm)    \leq
(\underbrace{1+\kp_{\mb}^{-1/2}}_{\leq 3/2}) \Big(
   \frac{| \pa^{\al} \tvrm|^2}{ (2\al \vrho)^2}+  \frac{ | \pa^{\al} \tum|^2}{\cc^2} 
  +\kp_{\mb} \frac{ | \pa^{\al} \tbm|^2}{ \mb^2 } \Big)  .
  \label{energy:L20_equiv}
\end{align}
The other constraint on $\kp_{\mb}$ relates to the outgoing property of the stationary profiles (cf.~\eqref{eq:profile_outgoing}). Since $\bc(0) = 0$, 
$\bc(y) > 0$ for $y \neq 0$, using the decay estimates \eqref{eq:profile_decay}, we obtain that there exists a finite constant $C  = C_{\eqref{eq:C:sublinear}}>0$ such that 
\begin{equation}
|\bc(y) | \leq |y| C_{\eqref{eq:C:sublinear}}  
\label{eq:C:sublinear}
\end{equation}
for any $y\neq 0$. Thus, if we define 
\begin{equation}
  \kp_{\mb} = \max \Bigl\{ 
  \frac{1}{100} ,  
  \frac{36 C_{\eqref{eq:C:sublinear}}^2}{C_{\eqref{eq:profile_outgoing}}^2}  \Bigr\}, 
  \label{eq:def:kp_mb}
\end{equation}
then the outgoing condition \eqref{eq:profile_outgoing} implies
\begin{equation}
    \bcr |y| +  \bar U(y)  - ( 1 + 3 \kp_{\mb}^{-1/2} ) \bc(y)  \geq
    \bigl(  C_{ \eqref{eq:profile_outgoing} } -  3 C_{\eqref{eq:C:sublinear}} \kp_{\mb}^{-1/2} \bigr) |y|
    \geq 
    \tfrac{1}{2} C_{ \eqref{eq:profile_outgoing}} |y| .
    \label{eq:outgoing_strong}
\end{equation}

\begin{remark}[\bf Treating $\kp_{\mb}$ as an absolute constant]\label{rem:kp_mb}
Since $\kp_{\mb}$, defined in \eqref{eq:def:kp_mb}, depends only on the profiles $\bar U$ and $\bc$, and on the parameters $\al$, $d$, and $\NNN$, we treat $\kp_{\mb}$ as an absolute constant
throughout the rest of the section. As such, we do not explicitly track the dependence of
most constants on  $\kp_{\mb}$.
\end{remark}

\subsubsection{Energy identities}
The time derivative of $\EE_{k}$ is calculated via~\eqref{eq:lin}--\eqref{eq:non}, and may be decomposed as
\begin{subequations}\label{eq:Hk_est}
\begin{equation}
 \label{eq:Hk_est:a}
 \frac{1}{2} \frac{d \EE_k}{d \tau} 
 =
  \frac{1}{2} \frac{d}{d \tau} \int E_{k}(\twm) \vp_k 
  = \sum_{|\al| =k}  \frac{1}{2} \frac{d}{d \tau} \int E_{\al}(\twm) \vp_k 
  =  I_{\cE, k} +  \sum_{|\al| = k} (\cI_{\cN, \al} + \cI_{s, \al} ),
\end{equation}
where $\cI_{s, \al}$ denotes the terms with time derivative $\pa_\tau$ acts on the denominators in~\eqref{energy:Hk}
\begin{equation}
 \label{eq:Hk_est:b}
 \cI_{s, \al} \! = \! - \int \Big(   \frac{\pa_\tau \vrho}{\vrho} \cdot \frac{ | \pa^{\al} \tvrm|^2}{ (2\al \vrho)^2}+
  \frac{\pa_\tau \cc}{\cc} \cdot  \frac{ | \pa^{\al} \tum|^2}{\cc^2} 
  +  \frac{\pa_\tau \mb}{\mb} \cdot \kp_{\mb} \frac{ | \pa^{\al} \tbm|^2}{ \mb^2 }
  +  (   \frac{\pa_\tau \vrho}{\vrho} +  \frac{\pa_\tau \mb}{\mb}   )  \frac{1}{\gamma}   \frac{ \pa^{\al} \tvrm}{ 2\al \vrho} \cdot \frac{ \pa^{\al} \tbm}{ \mb } 
  \Big) \vp_k ,
\end{equation}
the terms $\cI_{\cN, \al}$ are defined in terms of the nonlinear terms $\cN_{\bullet}$ appearing in~\eqref{eq:lin} as 
 \begin{align}
  \label{eq:Hk_est:c}
\cI_{\cN, \al} :=&  
\int \Big( \pa^{\al} \cN_{\vrho} \cdot  \frac{\pa^{\al} \tvrm}{ (2 \al \vrho)^2 }+ \pa^{\al} \cN_U \cdot   \frac{\pa^{\al} \tum}{\cc^2}
  + \kp_{\mb} \pa^{\al} \cN_{\mb} \cdot   \frac{\pa^{\al} \tbm}{\mb^2} 
  \notag \\
  &\qquad \qquad \qquad \qquad 
  +  \frac{1}{2\alpha \gamma \cc^2} ( \pa^{\al} \cN_{\vrho} \cdot  \pa^{\al} \tbm + \pa^{\al} \cN_{\mb} \cdot \pa^{\al} \tvrm  )
   \Big)   \vp_k , 
 \end{align}
 and $\cI_{\cE,k}$ captures the error terms $\cE_{\mm,\bullet}$ appearing in~\eqref{eq:lin}
 \begin{align} 
 \label{eq:Hk_est:d}
\cI_{\cE, k} &:=  
\sum_{|\al| = k} 
\int \Big( \pa^{\al} \cE_{\mm, \vrho} \cdot  \frac{\pa^{\al} \tvrm}{ (2 \al \vrho)^2 }+ \pa^{\al} \cE_{\mm, U} \cdot   \frac{\pa^{\al} \tum}{\cc^2}
  + \kp_{\mb} \pa^{\al} \cE_{\mm, \mb} \cdot   \frac{\pa^{\al} \tbm}{\mb^2}  
  \notag \\
  &\qquad \qquad \qquad \qquad \qquad \qquad 
  +  \frac{1}{2\alpha \gamma \cc^2}  ( \pa^{\al} \cE_{\mm, \vrho} \cdot  \pa^{\al} \tbm + \pa^{\al} \cE_{\mm, \mb} \cdot \pa^{\al} \tvrm  )
   \Big)   \vp_k .
   \end{align}
 \end{subequations}
The remainder of this subsection is dedicated to bounding the terms $ \cI_{s, \al}$, $\cI_{\cN, \al}$, and $\cI_{\cE, k}$. 

To simplify the notation, we shall denote the \textit{advection} operator $\cA$ and its action on $\log \vp_k$ by
\begin{equation}\label{nota:tran}
  \cA f := ( \cx y + \UU ) \cdot \na f ,
  \qquad 
    d_{\vp_k} :=  \frac{(\cx y + \UU) \cdot \na \vp_k}{\vp_k} = \cA (\log \vp_k).
\end{equation}
In the following derivations, we use ``mathcal" notation $\cI_{\bullet}$ to denote integrals, and  $I_{\bullet}$ to denote the corresponding integrands.  

\subsubsection{Rewriting $\cI_{s, \al}$}
We rewrite the $\cI_{s, \al}$ term defined in \eqref{eq:Hk_est:b} by using the equations \eqref{eq:sys}, the definition $\cc = (\vrho \mb)^{1/2}$, and the identity $\cu = \frac{1}{2}( \crho + \cbb )$, to arrive at
\begin{align*} 
  \frac{\pa_\tau \vrho}{\vrho} 
  &=  - \frac{ \cA \vrho }{\vrho} + \crho - 2 \alpha \div \UU 
  \\
  \frac{\pa_\tau \mb}{\mb} 
  &= - \frac{\cA \mb}{\mb} + \cbb ,  
  \\
  \frac{\pa_\tau \cc}{\cc} 
  &= \frac{1}{2} (   \frac{\pa_\tau \vrho}{\vrho} +   \frac{\pa_\tau \mb}{\mb}  )
  = - \frac{1}{2} ( \frac{ \cA \vrho }{\vrho}  +   \frac{\cA \mb}{\mb})
  + \cu - \alpha \div \UU 
    = -   \frac{ \cA \cc }{\cc} + \cu - \alpha \div \UU  .
\end{align*}
Recalling~\eqref{eq:Hk_est:b}, it follows that 
\begin{align}
  \cI_{s, \al} & =  \int 
  \Big(  (   \frac{ \cA \vrho }{\vrho}  - \crho + 2 \alpha \div \UU  ) \cdot \frac{ | \pa^{\al} \tvrm|^2}{ (2\al \vrho)^2}
  + ( \frac{ \cA \cc }{\cc} - \cu + \alpha \div \UU) \cdot  \frac{ | \pa^{\al} \tum|^2}{\cc^2} \notag \\
& \qquad  \qquad
+ (\frac{\cA \mb}{\mb} - \cbb ) \cdot \kp_{\mb} \frac{ | \pa^{\al} \tbm|^2}{ \mb^2 } 
+ (\frac{ \cA \cc }{\cc} - \cu + \alpha \div \UU) \frac{2}{\gamma}    \frac{ \pa^{\al} \tvrm}{ 2\al \vrho} \cdot \frac{ \pa^{\al} \tbm}{ \mb } 
  \Big) \vp_k
  \notag\\
  & =: \int I_{s, \al} \vp_k
  .
\label{eq:lin_ds}
\end{align}
Note that the bound~\eqref{eq:profile_decay} and the bootstrap \eqref{eq:boot2}, together with $\bcu < \tfrac{1}{2} \bcb$, imply the pointwise bound
\begin{align}
|\div \UU|(y)
\leq 
 \la y \ra^{- 1 + \ddb + 4 \ddd}
+  |\div \bu|(y)
& \leq
 \la y \ra^{-1 + \ddb + 4 \ddd}
+
 C_{\eqref{eq:profile_decay}}\la y \ra^{-1 + \ddb}
 \notag\\
&\leq 
(1 + C_{\eqref{eq:profile_decay}}) \la y \ra^{-1 + \dda - 4 \ddd}.
\label{eq:div:U:useful}
\end{align}

\subsubsection{Decomposition of $\cI_{\cN, \al}$ into transport, dissipation, and remainders}
For any multi-index $\al$ with $|\al| =k$. Applying $\pa^{\al}$ to \eqref{eq:non} and recalling the notation in~\eqref{eq:bilin:d} and \eqref{eq:bilin}, we obtain 
\begin{subequations}
\label{eq:Hk_main}
\begin{align} 
  \pa^{\al} \cN_{\vrho} := 
  & \underbrace{ -( \cx y + \UU) \cdot \na \pa^{\al} \tvr_{\mm} - 2 \al \vrho \div \pa^{\al} \tu_{\mm}}_{:= \cT_{\vrho, \al} } + \underbrace{  ( \cvr - k \cx) \pa^{\al}\tvr_{\mm} }_{:= \cD_{\vrho, \al} } 
  + \cR_{\al, 1},  
 \\
  \pa^{\al} \cN_U := 
  & \underbrace{-  ( \cx y + \UU) \cdot \na \pa^{\al} \tu_{\mm}  - \tfrac{1}{2\alpha} \mb \na \pa^{\al} \tvr_{\mm}
  - \tfrac{1}{\gamma} \vrho \na \pa^{\al} \tb_{\mm}  }_{:=\cT_{U, \al}} +
\underbrace{   (\cu - k \cx) \pa^{\al} \tu_{\mm}  }_{:=\cD_{U,\al} }
  + \cR_{\al, 2},  
    \\
\pa^{\al} \cN_{\mb} := 
& \underbrace{ - (\cx y + \UU) \cdot \na \pa^{\al} \tb_{\mm} }_{:= \cT_{\mb, \al}} +\underbrace{ ( \cbb - k \cx) \pa^{\al} \tb_{\mm} }_{:= \cD_{ \mb, \al}}
+   \cR_{\al, 3},  
\end{align} 
where the $\cR_{\al, i}$ terms only contain derivatives of order $\leq k$ acting on $\twm$, and are given by 
\begin{equation}
\cR_{\al, i} =    \pa^{\al} \cB_i( \WW,  \twm ) - \cB_i(\WW, \pa^{\al} \twm)
  + \pa^{\al} \cB_i(\twm, \bw)  . 
  \label{eq:Hk_main:d}
\end{equation}
\end{subequations}
In \eqref{eq:Hk_main}, we have used the notation $\cT, \cD, \cR$ to single out \textit{transport, dissipative},  and  \textit{remainder} terms.

\subsubsection{Estimates of the transport terms  $\cT_{\bullet, \al}$}
We rewrite the contribution of the transport terms from~\eqref{eq:Hk_main} to $I_{\cN,\al}$  
in \eqref{eq:Hk_est:c} as
\begin{align*} 
\cI_{\cT, \al} 
&=  - \int \Big\{
\cA  \pa^{\al} \tvrm  \cdot   \frac{\pa^{\al} \tvrm }{(2\al \vrho)^2}
+ \cA \pa^{\al} \tum  \cdot  \frac{\pa^{\al} \tum }{\cc^2}
+ \kp_{\mb}  \cA  \pa^{\al} \tbm  \cdot  \frac{\pa^{\al} \tbm}{\mb^2}   \\
&\qquad \quad
+ \frac{1}{2 \al \gamma \cc^2} (  \cA \pa^{\al} \tvrm \cdot  \pa^{\al} \tbm  + \cA \pa^{\al} \tbm \cdot \pa^{\al} \tvrm )   
+ 2 \al \vrho \div \pa^{\al} \tum   \cdot \frac{\pa^{\al} \tvrm}{ (2\al \vrho)^2}     \\
&\qquad \quad
+ \frac{\pa^{\al} \tum}{\cc^2 } \Big( \frac{1}{2\alpha} \mb \na \pa^{\al} \tvrm + \frac{1}{\gamma} \vrho \na \pa^{\al} \tbm \Big)  
+ \frac{1}{2 \al \gamma \cc^2}  \cdot  2 \al \vrho \div \pa^{\al} \tum  \cdot \pa^{\al} \tbm
\Big\}  \vp_k \\
&=  - \int \Big\{
\frac{1}{2 (2\al \vrho)^2} \cA \bigl| \pa^{\al} \tvrm\bigr|^2  
+ \frac{1}{2 \cc^2}  \cA \bigl| \pa^{\al} \tum\bigr|^2
+  \frac{\kp_{\mb} }{2 \mb^2} \cA  \bigl| \pa^{\al} \tbm\bigr|^2    \\
&\qquad \quad
+ \frac{1}{2 \al \gamma \cc^2}   \cA \bigl( \pa^{\al} \tvrm \cdot  \pa^{\al} \tbm \bigr )   
+ \frac{1}{2\al \vrho}  \div \bigl(\pa^{\al} \tvrm \, \pa^{\al} \tum  \bigr) 
+ \frac{1}{  \gamma \mb}       \div \bigl( \pa^{\al} \tbm \, \pa^{\al} \tum\bigr)
\Big\}  \vp_k.
\end{align*}
In the second equality above, we have used that $\cc^2 = \vrho \mb$.
Therefore, using integration by parts the divergences and of the advection operator $\cA  = ( \cx y + \UU ) \cdot \na$,  we obtain
\begin{align} 
\cI_{\cT, \al}
&=  \int  \Big\{ 
\frac{ \na \cdot ( ( \cx y + \UU) \vrho^{-2} \vp_k) }{ 
2 \vrho^{-2} \vp_k }  \,  \frac{ | \pa^{\al} \tvrm|^2}{(2\al \vrho)^2}
+  \frac{ \na \cdot ( (  \cx y + \UU) {\cc}^{-2}\vp_k) }{ 
2 {\cc}^{-2}\vp_k }  \,  \frac{| \pa^{\al} \tum|^2 }{\cc^2}  
\notag \\
&\qquad \quad  
+ \kp_{\mb}  \frac{ \na \cdot ( (  \cx y + \UU) \mb^{-2}\vp_k) }{ 
2 \mb^{-2}\vp_k } \,  \frac{|\pa^{\al} \tbm|^2 }{\mb^2} 
+ \frac{\na \cdot ( (  \cx y + \UU) \cc^{-2} \vp_k)}{\cc^{-2} \vp_k} \,  \frac{\pa^{\al} \tvrm \pa^{\al} \tbm }{\gamma (2 \al \vrho) \mb} 
\notag \\
&\qquad \quad 
+ \frac{ \na (\vrho^{-1} \vp_k) }{\vrho^{-1} \vp_k} \cdot \frac{\pa^{\al} \tum \pa^{\al} \tvrm}{ 2\al \vrho }
+ \frac{\na (\mb^{-1} \vp_k)}{ \mb^{-1} \vp_k} \cdot \frac{\pa^{\al} \tbm \cdot \pa^{\al} \tum }{\gamma \mb} \Big\} \vp_k
\notag\\
&:= \int \bigl( I_{\cT,1} + I_{\cT,2}  + I_{\cT,3} + I_{\cT,4} + I_{\cT,5} + I_{\cT,6} \bigr)  \vp_k,
\label{eq:lin_cT0}
\end{align}
where $\{I_{\cT,i}\}_{i=1}^6$ are defined in the natural way, in order of appearance on the right side of \eqref{eq:lin_cT0}.

Next, we estimate the first-order-derivative-coefficients in the integrands of $\{ I_{\cT,i} \}_{i=1}^6$.  
Under the bootstrap assumptions \eqref{eq:boot1}, \eqref{eq:boot2}, using the decay estimate in~\eqref{eq:div:U:useful}, and using the notation in~\eqref{nota:tran} for $\cA$  and $d_{\vp_k}$, we extract the main terms, which do not posses sufficient decay  in $y$ 
\begin{subequations}\label{eq:Hk_est_coe}
\begin{align} 
\frac{ \na \cdot ( ( \cx y + \UU) \vrho^{-2} \vp_k) }{ 
\vrho^{-2}  \vp_k }
& = d_{\vp_k} - 2 \frac{\cA \vrho}{\vrho}  + d \cx + \div \UU
= (d_{\vp_k} + d \cx)  - 2 \frac{ \cA \vrho}{\vrho} + \OO(  \la y \ra^{-1 + \dda} ),  \\
\frac{ \na \cdot ( ( \cx y + \UU) \cc^{-2} \vp_k) }{\cc^{-2} \vp_k }
&= d_{\vp_k} - 2 \frac{\cA \cc}{\cc} + d \cx + \div \UU
= (d_{\vp_k} + d \cx) - 2 \frac{\cA \cc }{\cc}  + \OO(  \la y \ra^{-1 + \dda} ) , \\
\frac{ \na \cdot ( ( \cx y + \UU) \mb^{-2} \vp_k) }{\mb^{-2} \vp_k}
&= d_{\vp_k} -2  \frac{\cA \mb  }{\mb}
+ d \cx + \div \UU 
=  (d_{\vp_k} +  d \cx)  - 2 \frac{\cA  \mb}{\mb} +   \OO(  \la y \ra^{-1 + \dda} ).
\end{align}
Applying \eqref{eq:boot_decay}, we also obtain
\begin{equation}
  \frac{ \cc \na ( \vrho^{-1} \vp_k) }{ \vrho^{-1} \vp_k} 
    =   \frac{ \cc \na \vp_k}{\vp_k} + \OO (\la y \ra^{ -1 + \dda }) ,
   \qquad 
     \frac{ \cc \na ( \mb^{-1} \vp_k) }{ \mb^{-1} \vp_k} 
    = 
    \frac{ \cc \na \vp_k}{\vp_k} + \OO(  \la y \ra^{-1 + \dda} ).
\end{equation}
\end{subequations}
It is important to remark that the implicit constants in the $\OO$ symbols appearing in~\eqref{eq:Hk_est_coe} are \emph{independent of $\vp_k$ and $\mm$}. 

Using the above estimates, and recalling the definition of the energy density $E_{\al}$ from~\eqref{energy:Hk} (or equivalently, using the bound~\eqref{energy:L20_equiv}), we may bound the terms $I_{\cT, i}$ defined in~\eqref{eq:lin_cT0}, for $1\leq i \leq 4$, as  
\begin{subequations}\label{eq:lin_cT1}
\begin{align}
 &  I_{\cT,1} + I_{\cT,2} + I_{\cT,3} + I_{\cT,4}  \notag \\
  & \quad \leq  
 \Bigl( \frac{1}{2} ( d_{\vp_k} +  d \cx) + C_{\eqref{eq:lin_cT1:a}}  \la y \ra^{-1 + \dda} \Bigr) E_{\al}(\twm) 
 \notag  \\
  &  \qquad -   \frac{\cA \vrho}{\vrho} \cdot \frac{ |\pa^{\al} \tvrm|^2}{(2\al \vrho)^2}
 -   \frac{\cA \cc}{\cc}  \cdot \frac{ | \pa^{\al} \tum|^2}{\cc^2}
 -  \frac{\cA \mb}{\mb}  \cdot  \kp_{\mb} \frac{ | \pa^{\al} \tbm |^2}{\mb^2}
  -  \frac{\cA \cc}{\cc}   \cdot \frac{2}{\gamma} \frac{ \pa^{\al} \tvrm}{2\al \vrho} \cdot \frac{ \pa^{\al} \tbm}{\mb}    , \label{eq:lin_cT1:a}
\end{align}
for a constant $C_{\eqref{eq:lin_cT1:a}} >0 $ which is independent of $k$, $\mm$, $\kp_{\mb}$, and the weight $\vp_k$.
Similarly, for $I_{\cT, 5}, I_{\cT,6}$, we apply the Cauchy-Schwarz inequality and the bound \eqref{eq:Hk_est_coe} to deduce 
\begin{align} 
  |I_{\cT, 5}
  + I_{\cT, 6}| 
 &  =  \Big|  \frac{ \cc \na ( \vrho^{-1} \vp_k) }{ \vrho^{-1} \vp_k} 
\cdot \frac{ \pa^{\al} \tum }{\cc} \frac{ \pa^{\al} \tvrm}{2 \al \vrho} 
+        \frac{ \cc \na ( \mb^{-1} \vp_k) }{ \mb^{-1} \vp_k } 
 \cdot \frac{1}{\gamma} \frac{\pa^{\al} \tbm}{\mb} \frac{ \pa^{\al} \tum}{\cc} \Big|   
 \notag \\
 & \leq C_{\eqref{eq:lin_cT1:ab}}
  \la y \ra^{-1 + \dda} E_{\al}(\twm)
+ \frac{\cc |\na \vp_k|}{\vp_k} \Big(  \frac{ |\pa^{\al} \tum| }{\cc} \frac{| \pa^{\al} \tvrm|}{2 \al \vrho} + 
\frac{1}{\gamma} \frac{|\pa^{\al} \tbm|}{\mb} \frac{|\pa^{\al} \tum|}{\cc} \Big), 
\label{eq:lin_cT1:ab}
\end{align} 
for a constant $C_{\eqref{eq:lin_cT1:ab}} >0 $ which is independent of $k$, $\mm$, and the weight $\vp_k$.
Since $\gamma > 1$, applying the Cauchy-Schwarz inequality and the bound \eqref{energy:L20_equiv}, we note that 
\begin{align*} 
 \frac{ |\pa^{\al} \tum| }{\cc} \frac{| \pa^{\al} \tvrm|}{2 \al \vrho} + 
\frac{1}{\gamma} \frac{|\pa^{\al} \tbm|}{\mb} \frac{|\pa^{\al} \tum|}{\cc}  
& \leq \frac{1}{2} \Bigl(  \frac{|\pa^{\al} \tum|^2}{\cc^2} +
\frac{|\pa^{\al} \tvrm|^2}{(2 \al \vrho)^2} 
+\kp_{\mb}^{-1}  \frac{|\pa^{\al} \tum|^2}{\cc^2}
+ \kp_{\mb} \frac{|\pa^{\al}\tbm|^2}{\mb^2} \Bigr)  \\
& \leq \frac{1 + \kp_{\mb}^{-1}}{2}
(1 -  \kp_{\mb}^{-1/2} )^{-1} E_{\al}(\twm) . 
\end{align*} 
Combining the above two estimates and using the second bound in \eqref{eq:kp_large}, we derive 
\begin{equation}\label{eq:lin_cT1:b}
  |I_{\cT, 5}
  + I_{\cT, 6}| 
  \leq C_{\eqref{eq:lin_cT1:ab}} \la y \ra^{-1 + \dda} E_{\al}(\twm)  
  + \frac{1 + 3 \kp_{\mb}^{-1/2}}{2} \cdot
\frac{ |\cc \na \vp_k|  }{\vp_k} E_{\al}(\twm).
\end{equation}
\end{subequations}

\subsubsection{Rewriting the dissipative terms  $\cD_{\bullet, \al}$ }
Recall from \eqref{eq:Hk_main} the definitions of the dissipative terms $\cD_{\bullet, \al}$; their contribution to the energy estimate is via~\eqref{eq:Hk_est:c}, and is given by
\begin{align}
 \cI_{\cD, \al} 
&=  \int \Big(  \cD_{\vrho,\al}    \frac{\pa^{\al} \tvrm}{ (2 \al \vrho)^2 }+  \cD_{U,\al}    \frac{\pa^{\al}  \tum }{\cc^2}  + \kp_{\mb}   \cD_{\mb,\al}     \frac{\pa^{\al} \tbm}{\mb^2} + \frac{1}{2\al \gamma \cc^2}  ( \cD_{\vrho,\al}   \pa^{\al} \tbm +  \cD_{\mb,\al}   \pa^{\al} \tvrm  ) \Big)   \vp_k 
\notag\\
& = \int \Bigl( ( \crho - k \cx) \frac{| \pa^{\al} \tvrm|^2 }{ (2 \al \vrho)^2 }
+ (\cu - k \cx) \frac{ |\pa^{\al}\tum|^2}{ \cc^2} 
  +  \kp_{\mb} (\cbb - k \cx) \frac{|\pa^{\al} \tbm|^2}{\mb^2}
  \notag\\
  &\qquad \qquad \qquad \qquad \qquad \qquad \qquad \qquad 
  +( \cu -  k \cx ) \frac{2   }{ \gamma }  \frac{ \pa^{\al} \tbm }{\mb}  \frac{\pa^{\al} \tvrm}{2\al \vrho}\Bigr) \vp_k
 =: \int I_{\cD, \al}   \vp_k, 
   \label{eq:lin_cD}
\end{align}
Note that the main terms in $I_{\cD, \al}$ (cf.~\eqref{eq:lin_cD}), $I_{s, \al}$ (cf.~\eqref{eq:lin_ds}), and the terms on the last row of \eqref{eq:lin_cT1:a} cancel each other.

\subsubsection{Combining the bounds for the principal terms}
Combining the bounds in~\eqref{eq:lin_cT1:a} and \eqref{eq:lin_cT1:b}, with the identities in~\eqref{eq:lin_ds} and~\eqref{eq:lin_cD},
and then using $ 2 \cu = \cbb + \crho$ (cf.~\eqref{eq:normal0}), the bound in~\eqref{eq:div:U:useful}, and the energy $E_{\al}$ (cf.~\eqref{energy:Hk}), for $|\alpha|=k$ we obtain 
\begin{align}\label{eq:Hk_est_top}
  &\sum_{1\leq i \leq 6} I_{\cT, i} + I_{\cD, \al} + I_{s, \al} \notag \\
  &\leq \Bigl( \frac{1}{2} ( d_{\vp_k} +  d \cx) - k \cx \Bigr) E_{\al}(\twm)
  \notag\\
  &\qquad +( 1 + C_{\eqref{eq:lin_cT1:a}}  + C_{\eqref{eq:profile_decay}} + C_{\eqref{eq:lin_cT1:ab}}) \la y \ra^{-1 + \dda} E_{\al}(\twm)  
  + \frac{1 + 3 \kp_{\mb}^{-1/2}}{2} \frac{ |\cc \na \vp_k|  }{\vp_k} E_{\al}(\twm)
  \notag\\
& =  \frac{1}{2} \left( d_{\vp_k} + (d - 2 k) \cx + (1 + 3 \kp_{\mb}^{-1/2})     \frac{ |\cc \na \vp_k|  }{\vp_k} \right) E_{\al}(\twm)
+ C   \la y \ra^{-1 + \dda}  E_{\al}(\twm) ,
\end{align}
with implicit constant $C = C_{\eqref{eq:Hk_est_top}}>0$ independent of the choice of $\vp_k$ and $ \MMM$ .

\subsubsection{Estimates for the lower order terms $\cR_{\al,\bullet}$}
We bound the contribution of the remainders $\cR_{\al,1}$, $\cR_{\al,2}$, and $\cR_{\al,3}$ appearing in \eqref{eq:Hk_main}.

We begin with the pointwise estimate for $\cR_{\al,1}$. From the definition of $\cB_i$ in~\eqref{eq:bilin}, applying $\p^\alpha$ to \eqref{eq:Hk_main:d}, using the decomposition in~\eqref{eq:pertb:main}, and rewriting functions as $g = \bar g + \tl_l \td g + \rl_l \td g $, we estimate 
\begin{align*} 
|\cR_{\al, 1} |  
& \leq    |\pa^{\al} \cB_1( \WW,  \twm ) - \cB_1(\WW, \pa^{\al} \twm)| 
+ |\pa^{\al} \cB_1(\twm, \bw)| \\
& \les_k  
 \sum_{ 1\leq i \leq k}  |\na^i  \UU| \cdot |\na^{k + 1-i} \tvrm| 
+ |\na^i \vrho |  \cdot |\na^{k+1 - i} \tum| 
\notag\\
&\qquad
+ \sum_{0 \leq i \leq k}
  |\na^i  \tum  | \cdot |\na^{k + 1-i} \bvr| 
+ |\na^i \tvrm |  \cdot |\na^{k + 1-i} \bu|     \\
& \les_k \sum_{1 \leq i \leq k} |\na^i \tum \cdot \na^{k+1 - i} \tvrm |
\notag\\
&\qquad
+ \sum_{0 \leq i \leq k}
 |\na^i  \tum  | \cdot |\na^{k + 1-i} (\bvr, \tl_{\mm} \tvr) | 
+ |\na^i \tvrm |  \cdot |\na^{k + 1-i} (\bu, \tl_{\mm+1} \tu)|    \\
& := \cR_{\al, 1, I} + \cR_{\al, 1, II} .
\end{align*} 
For $\cR_{\al, 1, I}$, since $k \leq k_*$ and $\min(i, k+1 -i) \leq k_* /2 + 1 $ we may apply estimate \eqref{eq:boot_res2} to obtain
\[
  \cR_{\al, 1, I} 
  \les_k  \sum_{1\leq i\leq k}
    \la y \ra^{-(k+1 - i) + 4 \ddd} ( \vrho |\na^i \tum| 
   + \cc |\na^i \tvrm|
   ) \| \twm \|_{\cX^{k_*-1}}.
\]
For the second summation, we apply estimates in \eqref{eq:boot_res3} to $\na^{k+1-i}(\bar g, \tl_l g)$ with $l\in \{\mm, \mm+1\}$, and deduce
\[
|\cR_{\al, 1,II}|
\les_{k } \sum_{0 \leq i \leq k}
\bigl( 1 + C_{\mm} E_{O, \mm}  ) |y|^{-(k + 1 - i)} \la y \ra^{ \ddd}  ( \cc  |\na^{ i } \tvrm|
  + \vrho   | \na^{ i } \tum| \bigr).
\]
Combining the two bounds above it follows that
\begin{subequations}\label{eq:Hk_cR_est1}
\begin{equation}
  |\cR_{\al, 1}| \les_k ( \| \twm \|_{\cX^{k_*-1}} + 
      1 + C_{\mm} E_{O, \mm}  ) \sum_{0 \leq i\leq k}
    |y|^{-(k + 1 - i)} \la y \ra^{ 4 \ddd} 
      \cdot ( \cc  |\na^{ i } \tvrm|
  + \vrho   | \na^{ i } \tum| ).
\end{equation}
Similarly, for the terms $\cR_{\al,2}, \cR_{\al,3}$ defined in \eqref{eq:Hk_main}, using the Leibniz rule and 
the pointwise estimates \eqref{eq:boot_res2}, \eqref{eq:boot_res3},  we obtain 
\begin{align}
  | \cR_{\al, 2} |  
  &\les_k  \sum_{1 \leq i \leq k} 
  |\na^i \UU| \cdot |\na^{k+1 - i} \tum|
  + |\na^i \mb| \cdot |\na^{k+1 - i} \tvrm|
  + |\na^i \vrho| \cdot |\na^{k+1 - i} \tbm| \\
  &\qquad +  \sum_{0 \leq i \leq k} 
  |\na^i \tum| \cdot |\na^{k+1 - i} \bu|
  + |\na^i \tbm| \cdot |\na^{k+1 - i} \bvr|
  + |\na^i \tvrm| \cdot |\na^{k+1 - i} \barb|   \notag \\
  &\les_k 
    ( \| \twm\|_{\cX^{k_*- 1}} +  1 + C_{\mm} E_{O, \mm}  ) 
   \sum_{0 \leq i \leq k}      |y|^{-(k + 1 - i)}  
   \la y \ra^{4 \ddd } 
 ( \cc |\na^{ i } \tum|
  + \mb |\na^i \tvrm| 
  + \vrho |\na^i \tbm|), \notag 
\end{align}
and also 
\begin{align}
|\cR_{\al, 3}| 
& \les_k   \sum_{1 \leq i \leq k} 
|\na^i \UU| \cdot |\na^{k+1 - i} \tbm|
+   \sum_{0 \leq i \leq k}  |\na^i \tum| \cdot |\na^{k+1 - i} \barb|    \\
 & \les_k
( \| \twm\|_{\cX^{k_*- 1}} +  1 + C_{\mm} E_{O, \mm}  ) 
\sum_{0\leq i\leq k}     |y|^{-(k + 1 - i)} 
\la y \ra^{ 4 \ddd } 
 ( \mb |\na^{ i } \tum| + \cc |\na^i \tbm| ) . \notag 
\end{align}
\end{subequations}

For the bounds in~\eqref{eq:Hk_cR_est1} to be useful, we still need to estimate the $\WW$--to--$\twm$ products present in the summations over $0\leq i \leq k$ in~\eqref{eq:Hk_cR_est1}.  Using that $\cc \les |y|\la y \ra^{\ddb + \ddd - 1}$ (cf.~\eqref{eq:boot1_res1:c}), the bootstrap in~\eqref{eq:boot1} and its consequence~\eqref{eq:boot1_res1:b}, the $ g_{\msf u} \les g   \la y \ra^{2 \ddd}$ for $g \in \{ \cc, \vrho, \mb\}$ (cf.~\eqref{eq:Wbar_upper}), we  obtain 
\begin{align*} 
 \vrho^{-1} ( \cc |\na^{ i } \tvrm|
   + \vrho  | \na^{ i } \tum| )
& = \cc ( \vrho^{-1} |\na^{ i } \tvrm| 
+ \cc^{-1}  | \na^{ i } \tum |) \\
& \les |y| \la y \ra^{\ddb + 3 \ddd - 1} \bigl( {\vru}^{-1} |\na^{ i } \tvrm| 
+ \ccu^{-1} | \na^{ i } \tum | \bigr).
\end{align*}
Similarly, using $\cc^2 = \vrho \mb$, we obtain 
\begin{align*} 
  \cc^{-1} 
( \cc |\na^{ i } \tum|
  + \mb |\na^i \tvrm| 
  + \vrho |\na^i \tbm|)  
& = \cc (  \cc^{-1} |\na^{ i } \tum| + \vrho^{-1} |\na^i \tvrm| 
+  \mb^{-1}  |\na^i \tbm|) \\
& \les |y| \la y \ra^{\ddb + 3 \ddd - 1} \bigl( \ccu^{-1} |\na^{ i } \tum| + {\vru}^{-1} |\na^i \tvrm| 
+  \bbu^{-1}  |\na^i \tbm| \bigr), 
\end{align*}
and
\begin{align*}
  \mb^{-1} ( \mb |\na^{ i } \tum|
  + \cc |\na^i \tbm| )
  & = \cc ( \cc^{-1 }  |\na^{ i } \tum| + \mb^{-1} |\na^i \tbm| ) \\
  & \les |y| \la y \ra^{\ddb + 3 \ddd - 1} ( \ccu^{-1 }  |\na^{ i } \tum| + \bbu^{-1} |\na^i \tbm| ) .
\end{align*}

Since $\dda < 1$ and $\ddb + 8 \ddd - 1 = \dda -1  < \frac{1}{2} (\dda - 1)$ (see definition \eqref{eq:para_del}), plugging the above estimates in \eqref{eq:Hk_cR_est1}, 
we obtain 
\begin{align}\label{eq:Hk_low0}
    & \vrho^{-1} |\cR_{\al,1} | 
    + \cc^{-1} |\cR_{\al,2}| 
    + \mb^{-1} |\cR_{\al,3} | \\
  & \les_k 
( \| \twm \|_{\cX^{k_* -1}}
    +  1 + C_{\mm} E_{O, \mm}  )
   \sum_{0 \leq i\leq k}
   \la y \ra^{\ddb + 7 \ddd -1} |y|^{-(k - i)} \bigl( {\vru}^{-1} |\na^{ i } \tvrm| 
+ \ccu^{-1} |  \na^{ i } \tum| + \bbu^{-1} |\na^i \tbm| \bigr)   \notag \\
& \les_k 
( \| \twm \|_{\cX^{k_* -1}}
    +  1 + C_{\mm} E_{O, \mm}  )
   \sum_{0 \leq i\leq k}
   \la y \ra^{ \frac{1}{2} ( \dda -1) } |y|^{-(k - i)} \bigl( {\vru}^{-1} |\na^{ i } \tvrm| 
+ \ccu^{-1} |  \na^{ i } \tum| + \bbu^{-1} |\na^i \tbm| \bigr) . \notag
\end{align}

With~\eqref{eq:Hk_low0} in hand, we turn to estimating the  $\cR_{\al,\bullet}$ contribution to ${\cI}_{\cN, \al}$, which is derived from~\eqref{eq:Hk_est:c} and~\eqref{eq:Hk_main}, and is bounded using Cauchy-Schwartz as (we also recall the definitions of the energies in~\eqref{energy:Hk}--\eqref{energy:Hk_norm} and their estimates \eqref{energy:L20_equiv}
)
\begin{subequations}
\label{eq:Hk_low1}
\begin{align} 
&\int \Big| \cR_{\al,1}    \frac{\pa^{\al} \tvrm}{ (2 \al \vrho)^2 }+ \cR_{\al,2}  \frac{\pa^{\al} \tum}{\cc^2}
  + \kp_{\mb} \cR_{\al,3} \frac{\pa^{\al} \tbm}{\mb^2} 
  +  \frac{1}{2\alpha \gamma \vrho \mb} ( \cR_{\al,1}  \pa^{\al} \tbm + \cR_{\al,3} \pa^{\al} \tvrm  )
   \Big|   \vp_k 
\les_{k}{\cI}_{\cR, \al}^{1/2} \EE_k^{1/2} ,
 \label{eq:Hk_low0a}
\end{align}
where\footnote{
We recall from Remark \ref{rem:kp_mb} that $\kp_{\mb}$, as defined in \eqref{eq:def:kp_mb}, 
is treated as an absolute constant.
}
\begin{align}
{\cI}_{\cR, \al}:= \int  \vp_k 
 \bigl(  \vrho^{-2} \cR_{\al, 1}^2  + \cc^{-2} \cR_{\al, 2}^2  +  \mb^{-2} \cR_{\al, 3}^2 \bigr).
\end{align}
To estimate ${\cI}_{\cR, \al}$, we use the fact that  $\vp_k = |y|^{2k} \vp_0$ (cf.~\eqref{eq:wgk_1:old}), the estimate~\eqref{eq:Hk_low0}, and recall the definition of the $\cX^l$-norm in~\eqref{norm:X}, to conclude
\begin{align}
{\cI}_{\cR, \al}
&\les_k   ( \| \twm \|_{\cX^{k_* -1}}  +  1 + C_{\mm} E_{O, \mm}  )^{2}
 \int  \la y \ra^{  \dda - 1 }  \vp_k 
( {\vru}^{-2} |\na^k \tvrm|^2 + \ccu^{-2} |\na^k \tum |^2  + \bbu^{-2} |\na^k \tbm|^2 ) 
 \notag \\
 & \quad +
 ( \| \twm \|_{\cX^{k_* -1}}
    +  1 + C_{\mm} E_{O, \mm}  )^{2}
    C_2(\vp_0) \| \twm\|_{\cX^{k-1}}^2 . \label{eq:Hk_low0b}  
\end{align}
For the second term in the above bound we have used~\eqref{eq:wgk_1}, have recalled the definition~\eqref{norm:X}, and the fact that $\dda<1$. Using that $g_{\msf{u}}^{-1}\les g^{-1} $ for $g \in \{ \mb, \cc, \vrho\}$ (cf.~\eqref{eq:boot1}, \eqref{eq:boot1_res1}),   and recalling the notation for $E_k$ from \eqref{energy:Hk}, we arrive at 
\begin{equation}
{\cI}_{\cR, \al}
 \leq C_{\eqref{eq:Hk_low0c}}  ( \| \twm \|_{\cX^{k_* -1}}
    +  1 +  C_{\mm} E_{O, \mm}  )^{2} 
    \Big( 
    \int   \la y \ra^{ \dda -1 } E_k (\twm) \vp_k 
    +  C(\vp_0)   \| \twm\|_{\cX^{k-1}}^2 \Big) ,
    \label{eq:Hk_low0c}
\end{equation}
where the constant $C_{\eqref{eq:Hk_low0c}}>0$ depends on $k$, but is \emph{independent of the weight} $\vp_0$.
\end{subequations}

\subsubsection{Estimate for the error term $\cI_{\cE, k}$}
Returning to~\eqref{eq:Hk_est:a}, we are only missing a bound for the error term $\cI_{\cE, k}$ defined in~\eqref{eq:Hk_est:d}. 
First, analogously to~\eqref{eq:Hk_low0a} we have
\begin{subequations}
\label{eq:IEk:JEk}
\begin{equation}
{\cI}_{\cE,k} 
\les_{k}
{\cJ}_{\cE, k}^{1/2} \EE_k^{1/2}\label{eq:IEk:JEk:a}
\end{equation}
where
\begin{equation}
 {\cJ}_{\cE, k}
 := \int \vp_k \bigl(  \vrho^{-2} |\na^k \cE_{\mm, \vrho}|^2 
 + \cc^{-2} |\na^k \cE_{\mm, U}|^2 
+ \mb^{-2} |\na^k \cE_{\mm,\mb}|^2 \bigr).
\end{equation}
We recall that the error terms $\cE_{\mm,\bullet}$ satisfy the bounds~\eqref{eq:error}; by also appealing to~\eqref{eq:Wbar_upper}, we obtain
\begin{align}  
 \cJ_{\cE, k} 
&\les_{k,\mm} 
\int  \vp_0 \la y \ra^{2 \ddd} \min\bigl( |y|^{\mm + 1}, 1 \bigr)^2  \; \cR_{\mm, k}^2 d y 
\notag\\
&\les_{k,\mm}  E_{O, \mm}^2 
\underbrace{\int  \vp_0 \la y \ra^{2 \ddd} \min\bigl( |y|^{\mm + 1}, 1 \bigr)^2  \Bigl(1 +  | \twm |_{\Gam^k_{\mm+1 , \mm+2, \mm+3}}^2 
\Bigr)d y}_{=:
 \cJ_{\cE, k,1}  +  \cJ_{\cE, k,2}} ,
\label{eq:IEk:JEk:c}
\end{align}
\end{subequations}
where $\cJ_{\cE, k, 1}, \cJ_{\cE, k, 2}$ denote the integrals associated with $1$ and the function $ | \twm |_{\Gam^k_{\mm+1 , \mm+2, \mm+3}}^2$, respectively.
In the second inequality above we have used the definition of $\cR_{\mm,k}$ in~\eqref{eq:error_func}, namely
\[
\cR_{\mm, k} =  E_{O, \mm}  + E_{O, \mm}  | \twm |_{\Gam^k_{\mm+1 , \mm+2, \mm+3}} .
\]

Using the definition of  $\Gam^k_l$ in~\eqref{norm:gam} and the bound $\min( |y|^{\mm + 1}, 1 ) |y|^{-\mm-1} 
\les_{\mm} \la y \ra^{-\mm-1}$, we obtain 
\begin{align*} 
 & \min(|y|^{\mm+1} , 1) | \twm(y) |_{\Gam^k_{\mm+1 , \mm+2, \mm+3} }  \\
 & \les   \min(|y|^{\mm+1} , 1)  \sum_{0 \leq i\leq k}  \bigl( |y|^{ i - \mm - 1 } |\na^i \tvrm (y)| +  |y|^{ i - \mm - 2 } |\na^i \tum (y)|
 + |y|^{ i - \mm - 3 } |\na^i \tbm (y)| \bigr) \\
 & \les_{\mm}  \la y \ra^{-\mm-1}  \sum_{0 \leq i\leq k} |y|^i \bigl(  |\na^i \tvrm (y)|
 +  |y|^{  - 1 } |\na^i \tum (y)| + |y|^{  - 2 } |\na^i \tbm (y)| \bigr)  .
\end{align*} 
Using the above estimate and the interpolation Lemma \ref{lem:interp},\footnote{Here Lemma~\ref{lem:interp} is applied with the weights $\vp_0 \brak{y}^{2(\ddd-\mm-1)} |y|^{2(n-\ell)}$ for $\ell \in \{0,1,2\}$. One may check that these weights satisfy~\eqref{ass:interp} with $C(A_1,n) = 1$ and $C(A_2,n) = C(\vp_0, \mm, n) $. Here we are using that the weight $\vp_0$ is constructed to satisfy the bound $|\nabla \vp_0| \leq C(\vp_0) |y|^{-1} \vp_0$ for some constant $C(\vp_0)>1$; see \eqref{eq:wgk_1:b}.}
we obtain 
\[
\begin{aligned} 
\cJ_{\cE,k, 2} & := \int \vp_0 \la y \ra^{2 \ddd} 
 \min( |y|^{\mm + 1 }, 1 )^2    | \twm(y) |_{\Gam^k_{\mm+1 , \mm+2, \mm+3} }^2 \\
& \leq C(\vp_0, k, \mm) \int \vp_0 \la y \ra^{2 \ddd - 2 \mm - 2}  
\sum_{i \in \{0, k\}}  |y|^{2i} ( |\na^i \tvrm|^2
 +  |y|^{  - 2  } |\na^i \tum|^2
 + |y|^{ -4 } |\na^i \tbm|^2 ) .
 \end{aligned}
\]
Note that Remark~\ref{rem:decay}, the bootstrap assumptions \eqref{eq:boot1}, \eqref{eq:boot1_res1}, 
and the bound $\ddb + 2 \ddd < 1$ (cf.~\eqref{eq:para_del}), imply that $\vrho \les 1, \cc \les |y| ,  \mb \les |y|^2$.
Further using the energy \eqref{energy:Hk}, the bound in the previous display, and also  $ 2 \ddd - 2 \mm - 2 < 0$  (cf.~\eqref{eq:para_del}), we thus obtain 
\begin{align*} 
  \cJ_{\cE, k, 2} 
  & \leq 
  C(\vp_0, k, \mm) \sum_{i \in \{0, k\}}
  \int \vp_0 \la y \ra^{2 \ddd - 2 \mm - 2}   |y|^{ 2 i }  \bigl( \vrho^{-2}  |\na^i \tvrm|^2 + \cc^{-2}   |\na^i \tum|^2 + \mb^{-2}  |\na^i \tbm|^2 \bigr)  \\
 & \leq   C(\vp_0, k, \mm) 
 \int \la y \ra^{2 \ddd - 2 \mm - 2}  \bigl( \vp_0 E_0(\twm)  + \vp_k E_k(\twm) \bigr)
 \\
 &\leq  C(\vp_0, k, \mm)     ( \EE_0 + \EE_k ) .
\end{align*}

Returning to~\eqref{eq:IEk:JEk:c}, for the first integral we appeal to the bound \eqref{eq:wgk_1} for  $\vp_0$ and obtain 
\[
\cJ_{\cE,k, 1} :=
  \int \vp_0 \la y \ra^{2 \ddd}  \min( |y|^{\mm + 1 }, 1 )^2
  \leq C_2(\vp_0) \int ( |y|^{-d + 2} \one_{|y| \leq 1} + \la y \ra^{-2 \ddd - d } \one_{|y| \geq 1} )
  \leq C(\vp_0) .
\]
Combining the estimates in the above two displays, we thus have 
\begin{subequations}\label{eq:Hk_est_error}
\begin{equation}
\label{eq:Hk_est_error:a}
\cJ_{\cE, k} \leq C(\vp_0, k, \mm)    E_{O, \mm}^2 \big(1 + \EE_0 + \EE_k \big) .
\end{equation}
Together with~\eqref{eq:IEk:JEk:a}, the bound~\eqref{eq:Hk_est_error:a} provides an estimate for $\cI_{\cE, k}$ (defined in~\eqref{eq:Hk_est:d}):
\begin{equation}
\label{eq:Hk_est_error:b}
 |\cI_{\cE, k}| 
 \leq
 C(\vp_0, k, \mm)    E_{O, \mm}  \big(1 + \EE_0 + \EE_k \big)^{1/2} \EE_k^{1/2}.
 \end{equation}
\end{subequations}

\subsection{Summary of the estimates} 
We now have bounded all terms on the right side of~\eqref{eq:Hk_est:a}. 

By summing~\eqref{eq:Hk_est_error:b} with the sum over $|\alpha| = k \geq 1$ in~\eqref{eq:Hk_est_top}, \eqref{eq:Hk_low1},  we arrive at 
\begin{align}
\frac{1}{2} \frac{d \EE_k}{d \tau} 
& \leq  \frac{1}{2} \int   \Big( d_{\vp_k} + (d- 2 k) \cx + ( 1 + 3 \kp_{\mb}^{-1/2} )   \frac{ |\cc \na \vp_k|  }{\vp_k} + C_{\eqref{eq:EE_Hk}}^{(0)}\la y \ra^{ \dda -1  } \Big) E_{k }(\twm) \vp_k 
 \notag \\
&\qquad + C_{\eqref{eq:EE_Hk}}^{(k)}
\bigl(1 + \| \twm \|_{\cX^{k_* -1}}^2 +  C_{\eqref{eq:EE_Hk}}^{(\mm)} E_{O, \mm}^2  \bigr)  
\Bigl( 
    \int   \la y \ra^{ \dda -1 } E_k (\twm) \vp_k 
    +  C_{\eqref{eq:EE_Hk}}^{(\vp_0)}   \| \twm\|_{\cX^{k-1}}^2  \Bigr)
  \notag\\
  &\qquad 
  + C_{\eqref{eq:EE_Hk}}^{(\vp_0,k,\mm)}    E_{O, \mm}^2  \big(1 + \EE_0 + \EE_k \big) 
 + \frac{\bcr \ddd}{4} \EE_k ,
 \label{eq:EE_Hk}
  \end{align}
where the various constants appearing in~\eqref{eq:EE_Hk} have the following dependences:
\begin{itemize}[leftmargin=2em]
\item $C_{\eqref{eq:EE_Hk}}^{(0)} >0$ is a constant that is \emph{independent} of $\vp_0,k,  \MMM$; 
\item $C_{\eqref{eq:EE_Hk}}^{(\mm)}>0$ is a constant that only depends on $\mm$; 
\item $C_{\eqref{eq:EE_Hk}}^{(k)}>0$ is a constant that only depends on $k$; 
\item $C_{\eqref{eq:EE_Hk}}^{(\vp_0)}>0$ is a constant that only depends on $\vp_0$; 
\item $C_{\eqref{eq:EE_Hk}}^{(\vp_0,k,\mm)}>0$ is a constant that  only depends on $\vp_0, k, \mm$. 
\end{itemize}
Note that since $\ddd, \kp_{\mb}$ depend only on $d, \al, \NNN$, and the self-similar exponents (cf.~\eqref{eq:para_del}, \eqref{eq:def:kp_mb}), we do not track the dependence of various constants  on $\ddd, \kp_{\mb}$.

In the case of $k = 0$, we do not have terms with $i\leq k-1$ derivatives in the estimates of $\cR_{\al,i}$ (cf.~\eqref{eq:Hk_low0} and~\eqref{eq:Hk_low1}), which contributes to the terms $\| \twm \|_{\cX^{k-1}}$ in~\eqref{eq:EE_Hk}.  Thus, for $k=0$ the bound we obtain is
\begin{align} 
\frac{1}{2}  \frac{d \EE_0}{d \tau}  
&\leq
\frac{1}{2} \int   \Big( d_{\vp_0} + d \cx + ( 1 + 3 \kp_{\mb}^{-1/2} )   \frac{ |\cc \na \vp_0|  }{\vp_0} +  C_{\eqref{eq:EE_L2}}^{(0)} \la y \ra^{\dda  -1  }    \Big) E_{0 }(\twm) \vp_0 
\notag  \\
&\qquad + C_{\eqref{eq:EE_L2}}^{(0)}
\bigl( \| \twm \|_{\cX^{k_* -1}}^2 +  C_{\eqref{eq:EE_L2}}^{(\mm)} E_{O, \mm}^2   \bigr) 
    \int   \la y \ra^{ \dda -1 } E_0 (\twm) \vp_0
  \notag\\
  &\qquad
  + C_{\eqref{eq:EE_L2}}^{(\vp_0,\mm)}   E_{O, \mm}^2  \big(1 + \EE_0   \big) 
 + \frac{\bcr \ddd}{4} \EE_0 \, , 
    \label{eq:EE_L2}
\end{align}
where the various constants appearing in~\eqref{eq:EE_L2} have the following dependences:
\begin{itemize}[leftmargin=2em]
\item $C_{\eqref{eq:EE_L2}}^{(0)} >0$ is a constant that is \emph{independent} of $\vp_0, \MMM$; 
\item $C_{\eqref{eq:EE_L2}}^{(\mm)}>0$ is a constant that only depends on $\mm$; 
\item $C_{\eqref{eq:EE_L2}}^{(\vp_0 ,\mm)}>0$ is a constant that  only depends on $\vp_0,  \mm$. 
\end{itemize}

\subsection{Choosing \texorpdfstring{$\mm$ and $\vp_0$}{M and varphi 0}}
\label{sec:choose_wg}
We aim to choose $\mm$ sufficiently large and suitably construct the weight function $\vp_0$, in order to obtain linear stability for $\EE_0$ and nonlinear stability for $\EE_{\kk} + \nu \EE_0$, for a suitable coupling parameter $\nu>0$. In this direction, we  have:

\begin{lemma}[\bf Choice of $\mm$ and $\vp_0$]
\label{lem:nonrad_wg}
Suppose that the bootstrap assumptions \eqref{eq:boot1} and \eqref{eq:boot2} hold. 
Let $\vp_\kk = |y|^{2 \kk} \vp_0$, with $\kk$ defined in \eqref{eq:k_star}. Let $\kp_{\mb}$ be as in \eqref{eq:def:kp_mb}. Fix the constants $C_{\eqref{eq:EE_L2}}^{(0)} >0$ appearing in \eqref{eq:EE_L2}, and the constants $C_{\eqref{eq:EE_Hk}}^{(0)}, C_{\eqref{eq:EE_Hk}}^{(\kk)}> 0 $ appearing  in \eqref{eq:EE_Hk} when $k = \kk$.

There exists a sufficiently large integer $\mm$ satisfying $\mm \geq \kk + 2$  (cf.~\eqref{eq:mm_low}), and a weight function $\vp_0$ satisfying \eqref{eq:wgk_1}, such that the following estimates hold for any $y \in \Reals^d:$
\begin{subequations}\label{eq:wg_stab}
\begin{align} 
d_{\vp_0} + d \cx + ( 1 + 3 \kp_{\mb}^{-\frac12} )  \frac{ |\cc \na \vp_0|  }{\vp_0} +   C_{\eqref{eq:EE_L2}}^{(0)} \la y \ra^{ -1 + \dda }   
 & \leq  - 2 \bcr \ddd + C_{\eqref{eq:wg_stab}}  \hat E_{\mm} , 
 \label{eq:wg_L2_stab} \\
   d_{\vp_\kk } + (d- 2 \kk) \cx + ( 1 + 3 \kp_{\mb}^{-\frac12} )  \frac{ |\cc \na \vp_\kk|  }{\vp_\kk} + 
 ( C_{\eqref{eq:EE_Hk}}^{(0)} + 2   C_{\eqref{eq:EE_Hk}}^{( \kk)}  )  \la y \ra^{ -1 + \dda }  
 &   \leq - 2 \bcr \ddd 
+ C_{\eqref{eq:wg_stab}} \hat E_{\mm} ,
 \label{eq:wg_Hk_stab}
\end{align}
\end{subequations}
where $ \hat E_{\mm} := \| \twm \|_{\cX^{k_* - 1 }} + E_{O, \mm}$, $d_{\vp_k} =  \frac{(\cx y + \UU) \cdot \na \vp_k}{\vp_k}$ (cf.~\eqref{nota:tran} for $k \in \{0,\kk\}$), and $
C_{\eqref{eq:wg_stab}}>0$ is a constant that depends on the choice of $\vp_0$ and $\mm$.
\end{lemma}

\begin{proof}[Proof of Lemma~\ref{lem:nonrad_wg}]
Note that $\kk$, $C_{\eqref{eq:EE_L2}}^{(0)}$, $C_{\eqref{eq:EE_Hk}}^{(0)}$, $C_{\eqref{eq:EE_Hk}}^{(\kk)}$, and $\kp_{\mb}$ only depend on $\al, d, \NNN$ and the profile; we treat these as absolute constants. 

Let $0 < R_2 < R_3$ and an integer $\ell\geq 1$, to be chosen later, and define
\begin{subequations}
\label{def:vp0_1}
\begin{equation}
  \vp_0 := \hat \vp_0 \; \psi^\ell,\qquad 
  \hat \vp_0 :=   |y|^{-2 \mm -d } + \la  y \ra^{-d - 4 \ddd} ,
    \label{def:vp0_1:a}
\end{equation}
where $\psi \in C^{\infty}(\Reals^d)$ is radially symmetric, to be chosen later, and satisfies
\begin{align}
  &\psi(y) = 1 , \quad  |y| \leq \tfrac{1}{2} R_2, 
  \qquad \quad
  \psi(y) = \tfrac{1}{2}, \quad |y| \geq 2 R_3, 
  \qquad \quad 
  \tfrac{1}{2} \leq  \psi(y) \leq 1, \quad y \in \Reals^d, 
  \label{def:vp0_1:b}
  \\
  &\pa_r \psi(y) \leq 0, \quad y  \in \Reals^d,
  \qquad \quad
  \pa_r \psi(y) < 0 , \quad |y| \in [R_2, R_3].
  \label{def:vp0_1:c}
\end{align}
\end{subequations}
Here $\partial_r \psi$ denotes the radial derivative $\frac{y}{|y|} \cdot \nabla \psi$.
We determine the aforementioned parameters and the weight $\psi$ in the following order 
\[
   \mm \leadsto R_2 \leadsto  R_3 \leadsto \psi \leadsto \ell ,
\]
with each parameter in this sequence being allowed to depend on those chosen earlier.

We verify that~\eqref{def:vp0_1} implies~\eqref{eq:wgk_1}. First, we note that by definition we have the useful property
\begin{equation}
\label{eq:wg_decrease}
   \pa_r \vp_0(y) \leq 0,  \qquad  \forall \, y \in \Reals^d.
\end{equation}
Then, using the definition of $\vp_0$ in~\eqref{def:vp0_1:c},  we obtain 
\begin{subequations}
\label{eq:wg_pf1}
\begin{equation}
2^{-\ell} \hat \vp_0 \leq \vp_0 \leq   \hat \vp_0, 
\qquad 
\frac{|y \cdot \na \vp_0|}{\vp_0} 
\leq \frac{\ell |y| |\nabla \psi|}{\psi} + \frac{|y| |\nabla \hat \vp_0|}{\hat \vp_0}
\leq C_3(\vp_0) < \infty,
\label{eq:wg_pf1:a}
\end{equation}
where $C_3(\vp_0) = 2 \mm + 2 d + 1 + 2 \ell \max_{R_2/2 \leq r \leq 2 R_3} |r \p_r \psi|$.
Thus, $\vp_0$ satisfies properties \eqref{eq:wgk_1} with constants $C_1(\vp_0) = 2^{-\ell}$, $C_2(\vp_0) =1$, and $C_3(\vp_0)$ as defined before.

Further, using \eqref{eq:boot_decay}, \eqref{eq:boot1_res1:c}, \eqref{eq:boot_res3}, \eqref{eq:boot_res2}, \eqref{eq:boot_est_tc}, and \eqref{eq:wg_pf1:a}, we obtain
\begin{equation}
  \quad (|\UU| + |\cc|) \frac{ |\na \vp_0| }{ \vp_0 }
  \les |y| \frac{ |\na \vp_0| }{ \vp_0 } \les_{\vp_0} 1,
  \qquad  
  ( |\tu| + |\td \cc| ) \frac{ |\na \vp_0| }{ \vp_0 } \les_{\vp_0, \mm}  \hat E_{\mm},
\end{equation}
where $\hat E_{\mm} =     E_{O, \mm} + \| \twm \|_{\cX^{\kk-1}} $.
\end{subequations}

In order to conclude the proof of Lemma~\ref{lem:nonrad_wg}, it remains to establish the two bounds in~\eqref{eq:wg_stab}.

\noindent
\underline{\emph{Derivation of the main terms}.} Recall the definition of $d_{\vp_k}$ from \eqref{nota:tran}; in particular,  since $\vp_{\kk } = |y|^{2 \kk} \vp_0$, we have $ \frac{\na \vp_{\kk}}{\vp_{\kk}}
= \frac{\na \vp_0}{\vp_0} + 2 \kk \frac{y}{|y|^2}$. Thus, we obtain
\[
\begin{aligned} 
  \msf{LHS}_{\eqref{eq:wg_Hk_stab}}
  & \leq \frac{ (\cx y + \UU) \cdot \na \vp_0 }{\vp_0} 
  +  (1 + 3 \kp_{\mb}^{-1/2}) \frac{ |\cc \na \vp_0| }{\vp_0}   \\
 & \quad  + \Big( (d - 2 \kk) \cx  
  + 2 \kk   \frac{(\cx y + \UU) \cdot y}{|y|^2} + 2 \kk  (1 + 3 \kp_{\mb}^{-1/2})   \frac{\cc}{ |y| }  \Big)  
  + ( C_{\eqref{eq:EE_Hk}}^{(0)} +   2 C_{\eqref{eq:EE_Hk}}^{(\kk)}  )  \la y \ra^{-1 + \dda} \\
& =: I_1 + I_2 + I_3   +
( C_{\eqref{eq:EE_Hk}}^{(0)} +  2  C_{\eqref{eq:EE_Hk}}^{( \kk )})   \la y \ra^{-1 + \dda} .
\end{aligned} 
\]
For $I_3$, since $\kp_{\mb}$ defined in \eqref{eq:def:kp_mb} is an absolute constant, and since $\kk = 2d+10$ (cf.~\eqref{eq:k_star}) is an absolute constant, using \eqref{eq:boot1_res1:c} and \eqref{eq:boot_decay:a}, we obtain 
\begin{equation}
  I_3 =  d \cx + 2\kk  \frac{ \UU(y)  \cdot y}{|y|^2} + 2\kk (1 + 3 \kp_{\mb}^{-1/2}) \frac{ \cc }{ |y| }  
\leq d \cx + C_{\eqref{eq:I3_est}} \la y \ra^{-1 + \dda}.
\label{eq:I3_est}
\end{equation}
for some absolute constant $ C_{\eqref{eq:I3_est}}>0$. Combining the above estimates and the formula of $\msf{LHS}_{\eqref{eq:wg_L2_stab}}$, we obtain 
\begin{align} 
\max \bigl\{ \msf{LHS}_{\eqref{eq:wg_L2_stab}},   \msf{LHS}_{\eqref{eq:wg_Hk_stab}} \bigr\}  
& \leq d  \cx + \frac{ (\cx y + \UU) \cdot \na \vp_0 }{\vp_0} 
  +  (1 + 3 \kp_{\mb}^{-1/2}) \frac{ |\cc \na \vp_0| }{\vp_0}  +  \mu \ang y^{-1 + \dda} 
  \notag \\
& =: d  \cx + I_1 + I_2  + \mu \ang y^{-1 + \dda}.
\label{eq:wg_pf2}
\end{align}
where  
\begin{equation}
\mu := \max \bigl\{   C^{(0)}_{ \eqref{eq:EE_L2}} , \ C^{(0)}_{ \eqref{eq:EE_Hk}}  
+  2 C^{(\kk)}_{ \eqref{eq:EE_Hk}}  +  C_{\eqref{eq:I3_est}} \bigr\}.
\label{def:mu}
\end{equation}
is  an absolute constant. 

Decomposing $\cx = \bcr + \tcr , \UU = \bu + \tu, \cc = \bc + \td \cc$, using \eqref{eq:wg_pf1},
and $|\tcr| \leq E_{O, \mm}$ \eqref{norm:ODE}, we derive the main terms 
\begin{align*} 
I_1 &=  \frac{ (\bcr y + \bu) \cdot \na \vp_0 }{\vp_0}  
+  \frac{ (\tcr y + \tu) \cdot \na \vp_0 }{\vp_0}  
 \leq \frac{ (\bcr y + \bu) \cdot \na \vp_0 }{\vp_0}  + C(\vp_0, \mm) \hat E_{\mm}, \\
I_2  &\leq   (1 + 3 \kp_{\mb}^{-1/2} )   \Bigl( \frac{ |\bc \na \vp_0| }{ \vp_0}
+  \frac{ |\td \cc \na \vp_0| }{ \vp_0 }  \Bigr)
\leq (1 + 3 \kp_{\mb}^{-1/2} )  \frac{ |\bc \na \vp_0| }{ \vp_0}  + C(\vp_0, \mm)  \hat E_{ \mm} ,
 \\
 d  \cx &\leq   d \cxbar + d \cxtilde   \leq  d \bcr + d \hat E_{\mm}.
\end{align*}
Note that $\bu = \bar U(|y|) \frac{y}{|y|}$, and $\vp_0$ is radially symmetric and decreasing (cf.~\eqref{eq:wg_decrease}), so that $|\na \vp_0| = - \pa_r \vp_0$.  Combining this information with the  above three estimates, we obtain 
\begin{subequations}\label{eq:wg_pf3}
\begin{equation}
\label{eq:wg_pf3:a}
d \cx +  I_1 + I_2  +  \mu \ang y^{-1 + \dda}
  \leq \cI(\vp_0) +d \bcr + C(\vp_0, \mm) \hat E_{\mm} 
  +  \mu \ang y^{-1 + \dda},
\end{equation}
where $\cI(\vp_0)$ is defined by setting $g= \vp_0$ in 
\begin{equation}
\label{eq:wg_pf3:b}
\cI(g)  := \bigl( \bcr r  + \bar U  -  (1 + 3 \kp_{\mb}^{-1/2} ) \bc \bigr) \frac{ \pa_r g}{g} ,
\end{equation}
with $r = |y|$ and $\p_r = \frac{y}{|y|} \cdot \nabla$.
Note that the definition \eqref{def:vp0_1} implies $\frac{\pa_r \vp_0}{\vp_0} =\frac{\pa_r \hat \vp_0}{\hat \vp_0} + \ell \frac{\pa_r \psi}{\psi}$; thus, with the notation in~\eqref{eq:wg_pf3:b} we may decompose 
\begin{equation}
\label{eq:wg_pf3:c}
    \cI(\vp_0) =  \cI(\hat \vp_0) + \ell    \cI(\psi).
\end{equation}
\end{subequations}

At this stage we record a crucial lower bound, which is a consequence of the outgoing condition \eqref{eq:profile_outgoing} and our choice of $\kp_{\mb}$ in~\eqref{eq:def:kp_mb}: from  \eqref{eq:outgoing_strong} and $\bc(y) \geq 0$  we obtain  
\begin{equation}
   \bcr r + \bar U(y) -  (1 + 3 \kp_{\mb}^{-1/2} )  \bc(y) \geq a_* r > 0 , 
   \qquad 
   a_* = \tfrac{1}{2} C_{\eqref{eq:profile_outgoing}} > 0,
  \label{eq:wg_pf4}
\end{equation}
for all $r = |y| >0$.

\noindent
\underline{\emph{The choice of $\mm$ and $R_2: |y|\ll 1$}.}
By definition \eqref{def:vp0_1:a} we have that $\hat \vp_0 = r^{-2\mm-d} + \brak{r}^{-d - 4 \ddd}$, and using the fact that $| r \la r \ra^{-1} - 1 | \leq  \la r \ra^{-1}$, we obtain 
\begin{align} 
  \frac{ r \pa_r \hat \vp_0}{\hat \vp_0} 
  & = \frac{ - (2\mm + d) r^{- 2\mm - d} - (d + 4\ddd)  r \la r \ra^{-d -4 \ddd-1} }{  r^{- 2\mm - d} + \la r \ra^{-d -4 \ddd} } 
  \notag\\
  &\leq  \frac{ - (2\mm + d) r^{- 2\mm - d} - (d + 4\ddd)   \la r \ra^{-d -4 \ddd} }{  r^{- 2\mm - d} + \la r \ra^{-d -4 \ddd} }  + (d + 4\ddd) \la r \ra^{-1} 
  \notag \\
  & = - (d + 4 \ddd) -  \frac{ (2 \mm - 4 \ddd) r^{-2 \mm-d}}{   r^{- 2\mm - d} + \la r \ra^{-d -4 \ddd} } 
  +(d + 4\ddd) \la r \ra^{-1} 
  \notag\\
  &\leq  - (d + 4 \ddd) - \frac{ 2 \mm - 4 \ddd   }{ 1  + r^{2 \mm  -4 \ddd} } 
  + (d + 4\ddd) \la r \ra^{-1}  .
\label{eq:wg_pf5}
\end{align}
Using \eqref{eq:wg_pf3}, the fact that $\pa_r \hat \vp_0 < 0$ (this follows from~\eqref{def:vp0_1}), and the above estimates, we obtain 
\begin{align*}
\cI(\hat \vp_0) + \mu \ang y^{-1 + \dda} 
&\leq \frac{a_* r \pa_r \hat \vp_0}{\hat \vp_0} +  \mu \ang y^{-1 + \dda} 
\notag\\
&\leq a_* \Big( - (d + 4 \ddd) -  \frac{ 2 \mm - 4 \ddd  }{   1  + r^{2 \mm  -4 \ddd} }  \Big) 
+ \bigl( a_* (d + 4\ddd) + \mu \bigr) \la r \ra^{-1 + \dda},
\end{align*}
where we recall that $\mu$ is the absolute constant defined in~\eqref{def:mu}.

In light of the above estimate, we choose $\mm$ as 
\begin{equation}\label{def:nonrad_m}
  \mm  := \max\left( \left\lceil \frac{ d \bcr + \mu + 1 +  a_* (4\ddd+1) }{2 a_* } \right\rceil, \kk + 2 \right).
\end{equation}
so that the constraint \eqref{eq:mm_low} is satisfied. 
Using continuity in $r$, the bounds $\dda\leq 1$ and  $2 \mm - 1 - 4 \ddd > 0$, we deduce that there exists $0< R_2 \ll 1 $, depending on $\mm$, such that 
\begin{align} 
  \cI(\hat \vp_0) + \mu \ang y^{-1 + \dda}  & \leq 
  a_* ( - (d + 4 \ddd) - (2 \mm - 4 \ddd - 1) ) + \bigl( a_* (d + 4\ddd) + \mu \bigr)   
  \notag \\
  & \leq - a_* (2 \mm - 4 \ddd - 1) +    \mu     \notag\\
  & \leq -1 - d \bcr ,\label{eq:damp_near}
\end{align}
for any $r = |y| \leq R_2$. Specifically, we may choose $R_2 := (2 \mm - \ddd -1)^{1/(2\mm - 4\ddd)}$.

\noindent{\emph{The choice of $R_3: |y|\gg 1$}.}
We combine definition~\eqref{eq:wg_pf3:b} with the bounds  \eqref{eq:profile_decay}, $\cubar/\cxbar < \ddb < \dda$ (cf.~\eqref{eq:para_del}), 
$3 \kp_{\mb}^{-1/2}  < 1$, and \eqref{eq:def:kp_mb}, to obtain 
\begin{align*}
  \cI(\hat \vp_0) 
  &\leq \bcr \frac{r \pa_r \hat \vp_0}{\hat \vp_0}
  + \frac{|\bu| + 2 |\bc|}{|y|} \cdot \frac{ r |\pa_r \hat \vp_0| }{\hat \vp_0} 
  \\
  &\leq \bcr \frac{r \pa_r \hat \vp_0}{\hat \vp_0}
  + (2\mm + 2 d+1) \frac{|\bu| + 2 |\bc|}{|y|}  
  \\
  &\leq \bcr \frac{r \pa_r \hat \vp_0}{\hat \vp_0}  + (2\mm + 2 d+1) C_{\eqref{eq:profile_decay}} \ang y^{-1 + \dda},
\end{align*}
where $C_{\eqref{eq:profile_decay}}>0$ is the absolute constant in the second inequality of \eqref{eq:profile_decay}, for $k \in \{0,1\}$. Since $\mm > 2\ddd$, using \eqref{eq:wg_pf5} we obtain 
\begin{subequations}\label{eq:damp_far}
\begin{equation}\label{eq:damp_far:a}
    \cI(\hat \vp_0) + \mu \ang y^{-1 + \dda} 
    \leq - (d + 4 \ddd) \bcr  + \mu_2 \la y \ra^{-1 + \dda},
\end{equation}
for all $y\in \Reals^d$, where $\mu_2:= d + 4 \ddd + (2\mm + 2 d+1) C_{\eqref{eq:profile_decay}}$  depends $\mm$. Since $\dda<1$, there exists a sufficiently large $R_3$, such that 
\begin{equation}
      \cI(\hat \vp_0) + \mu \ang y^{-1 + \dda}   \leq - (d +  3 \ddd) \bcr,  \qquad \forall |y| \geq R_3.
\end{equation}
\end{subequations}
Specifically, we may choose $R_3:= 1 + (\mu_2/ (\ddd \cxbar))^{1/(1+\dda)}$.

\noindent 
\underline{\emph{Choosing $\psi$ and $\ell: |y| = \OO(1)$}.}
With $\mm, R_2, R_3$ already chosen, we let $\psi$ be any smooth function  satisfying \eqref{def:vp0_1:b}--\eqref{def:vp0_1:c}.
Since $\pa_r \psi \leq 0$, using \eqref{eq:wg_pf4} and~\eqref{eq:wg_pf3:b}, we obtain 
\begin{equation*}
 \cI(\psi ) \leq 0,
\end{equation*}
which along with \eqref{eq:damp_near} and \eqref{eq:damp_far} imply 
\begin{subequations}
\label{eq:damp_psi_pf1}
\begin{align}
       \cI(\hat \vp_0) + \ell \cI(\psi) + d \bcr + \mu \ang y^{-1 + \dda} 
       & \leq - 1 - d \bcr + d \bcr \leq -1 ,
       \qquad \forall |y| \leq R_2,  \\
       \cI(\hat \vp_0) + \ell \cI(\psi) + d \bcr  + \mu \ang y^{-1 + \dda} 
       & \leq  - (d +  3 \ddd) \bcr 
       + d \bcr \leq - 3 \ddd \bcr, 
       \qquad \forall |y| \geq R_3 .
\end{align}
\end{subequations}

Since $\frac{ \pa_r \psi }{\psi}<0$ for $|y| \in [R_2, R_3]$ (cf.~\eqref{def:vp0_1:c}), there exists a strictly positive constanat $C(R_2, R_3)>0$ such that 
\[
  \frac{ r \pa_r \psi(y) }{\psi(y)} \leq - C(R_2, R_3), \quad \forall \  R_2 \leq |y| \leq R_3.
\]
Recall the parameters $\mu_2$ from \eqref{eq:damp_far:a} and $a_*$ from~\eqref{eq:wg_pf4}; in terms of these, we define $\ell$ as 
\begin{equation}\label{def:vp0_N}
 \ell =  \left \lceil \frac{ \mu_2 }{ a_* C(R_2, R_3) } \right\rceil.
\end{equation}
Using \eqref{eq:wg_pf4}, $\pa_r \psi\leq 0$ and the above estimates, for $|y| \in [R_2, R_3]$, we derive 
\begin{equation*}
  \ell \cdot \cI(\psi) \leq \ell a_*  \, \frac{r \pa_r \psi}{ \psi}
  \leq - \ell a_* C(R_2,R_3)
  \leq  - \mu_2,
\end{equation*}
which along with \eqref{eq:damp_far:a} implies 
\begin{align}
\label{eq:damp_mid}
       \cI(\hat \vp_0) + \ell \cI(\psi) + d \bcr + \mu \ang y^{-1 + \dda}
       &\leq - (d + 4 \ddd) \bcr 
       - \mu_2 + d \bcr  + \mu_2 \ang y^{-1 + \dda}
       \notag\\
       &\leq -4 \ddd \bcr , \qquad \forall |y| \in [R_2, R_3].
\end{align}

Combining \eqref{eq:wg_pf2}, \eqref{eq:wg_pf3}, \eqref{eq:damp_psi_pf1}, and \eqref{eq:damp_mid}, we prove 
\[
\max\{  \msf{LHS}_{\eqref{eq:wg_L2_stab}},   \msf{LHS}_{\eqref{eq:wg_Hk_stab}} \}  
\leq - 3 \ddd \bcr + C(\vp_0 , \mm)  \hat E_{\mm} .
\]
This  concludes the proof of  \eqref{eq:wg_stab}, and hence of the Lemma.
\end{proof}

\subsection{Closing the stability estimates}

We are in a position to complete the nonlinear stability estimate stated in~\eqref{eq:nonrad_decay}. As we have determined $\vp_0$ and $\mm$ in Section \ref{sec:choose_wg}, in order to simplify notation we do not track their dependence in various constants in this section. The implicit constants $C$ appearing below may depend on $\kp_{\mb}$, $\vp_0$, $\mm$, $d$, and on the profile.

Assuming the bootstrap assumptions \eqref{eq:boot1}--\eqref{eq:boot2}, upon combining \eqref{eq:wg_L2_stab} and \eqref{eq:EE_L2}, and using that $\dda < 1$ and $E_{O, \mm} < 1$ (cf.~\eqref{eq:EOM}), we obtain 
\begin{align}
  \frac{1}{2} \frac{d}{d\tau} \EE_0 
  &\leq
  \big( -  \bcr \ddd + C \| \twm \|_{\cX^{k_* - 1 }} + C  E_{O, \mm}  \big) \EE_{0}
\notag  \\
&\qquad + C
\bigl(  \| \twm \|_{\cX^{k_* -1}}^2 +   E_{O, \mm}^2   \bigr) 
    \int   \la y \ra^{ \dda -1 } E_0 (\twm) \vp_0
  + C   E_{O, \mm}^2  \big(1 + \EE_0   \big) 
 + \frac{\bcr \ddd}{4} \EE_0
 \notag
 \\
& \leq - \frac{3\bcr \ddd}{4}  \EE_0 + C_{\eqref{eq:EE_L2_simp1}}  ( \| \twm \|_{\cX^{\kk-1}} + \| \twm \|_{\cX^{\kk-1}}^2+ E_{O, \mm})  \EE_0
  + C_{\eqref{eq:EE_L2_simp1}} E_{O, \mm}^2 , \label{eq:EE_L2_simp1}
\end{align}
for a constant $C_{\eqref{eq:EE_L2_simp1}}>0$.
Similarly, upon combining \eqref{eq:wg_Hk_stab} with \eqref{eq:EE_Hk} when $k = \kk$, and using 
that $\dda < 1$ and $E_{O, \mm} < 1$ (cf.~\eqref{eq:EOM}), we derive 
\begin{align}
\frac{1}{2} \frac{d}{d\tau} \EE_{\kk}
&\leq \big( -  \bcr \ddd + C \| \twm \|_{\cX^{k_* - 1 }} + C  E_{O, \mm}  \big)  \EE_{\kk} 
+ C  \| \twm\|_{\cX^{\kk-1}}^2 
 \notag \\
&\qquad + C
\bigl(  \| \twm \|_{\cX^{k_* -1}}^2 +   E_{O, \mm}^2  \bigr)  
\Bigl( \int   \la y \ra^{ \dda -1 } E_{\kk} (\twm) \vp_{\kk} +    \| \twm\|_{\cX^{\kk-1}}^2  \Bigr)
  \notag\\
  &\qquad 
  + C  E_{O, \mm}^2  \big(1 + \EE_0 + \EE_{\kk} \big) 
 + \frac{\bcr \ddd}{4} \EE_{\kk}
 \notag\\
  & \leq  - \frac{3\bcr \ddd}{4}  \EE_{\kk}
    +  C_{\eqref{eq:EE_Hk_simp1} }  \| \twm\|_{\cX^{\kk-1}}^2 
    +   C_{\eqref{eq:EE_Hk_simp1}} E_{O,\mm}^2
      \notag \\
   & \quad  + C_{\eqref{eq:EE_Hk_simp1}}  ( \| \twm \|_{\cX^{k_* -1}} + \| \twm \|_{\cX^{k_* -1}}^2 +  E_{O, \mm}) (\EE_0 +  \EE_{\kk} + \| \twm\|_{\cX^{\kk-1}}^2 ) ,
  \label{eq:EE_Hk_simp1}
\end{align}
for a constant $C_{\eqref{eq:EE_Hk_simp1}}>0$.
Note that since we have determined $\vp_0$ and $\mm$ in Section \ref{sec:choose_wg},  
and $\kk$ has been chosen in \eqref{eq:k_star}, we do not track their dependence in constants in the above estimates. 
 
\subsubsection{Interpolation}
We recall that the norm $\|\cdot\|_{\cX^{\kk-1}}$ is defined by setting $l=\kk-1$ in~\eqref{norm:X}. Our aim is to interpolate $\| \twm\|_{\cX^{\kk-1}}$ between $\EE_0$ and $\EE_{\kk}$.

To do this, we consider the interpolation Lemma~\ref{lem:interp} with the sequence of weights $w_n(y) = |y|^{2n} \vp_0 g_{\msf{u}}^{-2}$ for $g \in \{ \vrs , \cc, \mb \}$. In order to verify that these weights $w_n$ satisfy assumption \eqref{ass:interp}, we first note that by definition it holds that
\[
w_n = (w_{n-1} w_{n+1})^{1/2}.
\]
Second, we note that~\eqref{eq:wgk_1:b} and \eqref{eq:weight_est:b} imply
\[
  \frac{ |\na w_n|}{w_n} \leq \frac{2n}{|y|}
  + \frac{C_3(\vp_0)}{|y|} + \frac{ |\na g_{\msf{u}}| }{g_{\msf{u}}}
\leq \frac{2n + C_3(\vp_0) + C_{\eqref{eq:weight_est:b}}}{|y|}.
\]
Since $|y|^{-1} w_n = (w_n w_{n-1})^{1/2}$, this shows that assumption \eqref{ass:interp} is satisfied. Consequently, for any $\upsilon>0$, by appealing to the bounds $ \hat \vp_0 \les \vp_0 $ (cf.~\eqref{eq:wgk_1}), $g_{\msf u}^{-1} \les g^{-1}$ for $g \in \{ \vrs,\cc,\mb\}$ (cf.~\eqref{eq:boot1} and~\eqref{eq:boot1_res1:b}), and to the norm equivalence in~\eqref{energy:L20_equiv}, 
we obtain from Lemma~\ref{lem:interp} that
\begin{align*} 
        \|\tw_{\mm}\|_{\cX^{\kk-1} }^2 & \leq 
        \upsilon
        \|  \vp_0^{1/2} |y|^{\kk} ( {\vru}^{-1} \na^{\kk} \tvr, \ccu^{-1} \na^{\kk} \tu, \bbu^{-1} \na^{\kk} \tb ) \|_{L^2}^2 
        + C(\upsilon)   \|  \vp_0^{1/2} ( {\vru}^{-1}  \tvr, \ccu^{-1}  \tu, \bbu^{-1}  \tb ) \|_{L^2}^2  \\
        & \les   \upsilon \EE_{\kk} +    C(\upsilon) \EE_0,
  \end{align*} 
where the implicit constant in the last inequality is independent of $\upsilon$. 
Upon changing $\upsilon$ to absorb the implicit constant, we obtain that for any $\upsilon>0$, there exists a constant $\bar C(\upsilon)>0$ such that 
\begin{equation}\label{eq:Xm_interp}
          \|\tw_{\mm}\|_{\cX^{\kk-1} }^2 \leq  \upsilon \EE_{\kk} +    \bar C(\upsilon) \EE_0.
\end{equation}

\subsubsection{Conclusion of the energy estimates}
Using \eqref{eq:Xm_interp} with 
\[
\upsilon = \frac{\bcr \ddd}{4  C_{\eqref{eq:EE_Hk_simp1}  } }
\]
we deduce from~\eqref{eq:EE_Hk_simp1} that
\begin{align}
    \frac{1}{2} \frac{d}{d\tau} \EE_{\kk} 
    &\leq 
    - \frac{\bcr \ddd}{2}  \EE_{\kk}
    +  C_{ \eqref{eq:EE_Hk_simp2}} \EE_0 +  C_{ \eqref{eq:EE_Hk_simp1}} E_{O,\mm}^2
    \notag\\
    &\quad + C_{ \eqref{eq:EE_Hk_simp1}}  ( \| \twm \|_{\cX^{k_* -1}} +  \| \twm \|_{\cX^{k_* -1}}^2 + E_{O, \mm}) (\EE_0 + \EE_\kk + \| \twm\|_{\cX^{\kk-1}}^2  )
   ,
  \label{eq:EE_Hk_simp2}
  \end{align}
for some absolute constant $ C_{ \eqref{eq:EE_Hk_simp2}}> 0$. 

Now, we define the total energy 
\begin{equation}
  \EE_{\msf{tot}} = \EE_{\kk} + \nu \EE_0,
  \qquad \nu := \frac{ 4  C_{ \eqref{eq:EE_Hk_simp2}} }{ \bcr \ddd }.
  \label{energy:tot}
\end{equation}
From  \eqref{eq:EE_Hk_simp2} + $ \nu \times$ \eqref{eq:EE_L2_simp1}, and appealing to \eqref{energy:tot} and \eqref{eq:Xm_interp} with $\upsilon=1$ for the second inequality, we obtain 
\begin{align}
  \frac{1}{2} \frac{d}{d\tau} \EEt
& \leq -\frac{\bcr \ddd}{2} \EEt  + \Bigl( C_{ \eqref{eq:EE_Hk_simp2}}  - \frac{\bcr \ddd}{4} \nu  \Bigr) \EE_0  
\notag  \\
& \quad + C (   \| \twm \|_{\cX^{k_* -1}} +   \| \twm \|_{\cX^{k_* -1}}^2 + E_{O, \mm})
 (\EE_0 + \EE_\kk + \| \twm \|_{\cX^{\kk-1}}^2   )
 + C E_{O, \mm}^2
\notag \\
 & \leq -\frac{\bcr \ddd}{2} \EEt   + C ( \EEt^{1/2}  + \EEt + E_{O, \mm}) \EEt + C E_{O, \mm}^2
\notag \\
 & \leq -\frac{\bcr \ddd}{4} \EEt   + C_{\eqref{eq:EE_tot1a}} \EEt^2 + C_{\eqref{eq:EE_tot1a}} E_{O, \mm}^2,
 \label{eq:EE_tot1a}
\end{align}
for some computable absolute constant $C_{\eqref{eq:EE_tot1a}}>0$. 
 Using assumption \eqref{eq:ass:EOM}, we further obtain 
\begin{equation}
  \frac{1}{2} \frac{d}{d\tau} \EEt  \leq 
   -\lam_2  \EEt +  C_{\eqref{eq:EE_tot1}}  e^{-2 \lam \tau} E_{O,\mm}(0)^2 + C_{\eqref{eq:EE_tot1a}}  \EEt^2,
   \qquad \lam_2 := \min \bigl\{ \tfrac{1}{4} \bcr \ddd, \lam \bigr\},
   \label{eq:EE_tot1}
\end{equation}
where $\lam$ is the exponent in assumption \eqref{eq:ass:EOM}, and $C_{\eqref{eq:EE_tot1}} = C_{\eqref{eq:ass:EOM}}^2 C_{\eqref{eq:EE_tot1a}} >0$. 

Upon defining 
\begin{subequations}
\label{eq:EE_tot:ass}
\begin{equation}
 \eps_* :=  \frac{ \lam_2}{1 + 4 C_{ \eqref{eq:EE_tot1a}} +  C_{\eqref{eq:EE_tot1}}  }, 
\end{equation}
if at the initial time $\tau=0$ we have
\begin{equation}
   E_{O, \mm}(0)  < \eps_*, 
   \qquad   \EEt(0)  < \eps_*, 
   \label{eq:EE_tot:ass:b}
\end{equation}
\end{subequations}
since
\[
   - \lam_2 \cdot 2 \eps_* +   C_{\eqref{eq:EE_tot1}}    \eps_*^2 + C_{\eqref{eq:EE_tot1a}}   \cdot (2 \eps_*)^2 
  <  0, 
\]
as long as the bootstrap assumptions \eqref{eq:boot1}--\eqref{eq:boot2} holds, the ODE \eqref{eq:EE_tot1} implies that
\begin{equation}\label{eq:EE_tot2}
  \EEt(\tau) < 2 \eps_*.
\end{equation}

Thus if the bootstrap assumptions \eqref{eq:boot1}--\eqref{eq:boot2} hold for all $\tau \in [0,\tau_{\sf max})$, and if \eqref{eq:EE_tot:ass:b} holds, then from \eqref{eq:EE_tot1}--\eqref{eq:EE_tot2}, we derive 
\[
\frac{1}{2} \frac{d}{d\tau} \EEt  
\leq 
-\lam_2 \EEt + C_{ \eqref{eq:EE_tot1} }  e^{-2 \lam \tau} E_{O,\mm}(0)^2 + 2 C_{ \eqref{eq:EE_tot1a} }  \eps_*  \EEt
= 
- \frac{1}{2} \lam_2\EEt  + C_{ \eqref{eq:EE_tot1} }  e^{-2 \lam \tau} E_{O,\mm}(0)^2.
\]
Since $\lam_2 \leq \lam < 2 \lam$, integrating the above inequality we arrive at 
\begin{align}
\EEt(\tau) 
&\leq e^{- \lam_2 \tau } \EEt(0) 
+ 2 C_{ \eqref{eq:EE_tot1} }  \int_0^{\tau} e^{- \lam_2(\tau -\tau')} 
e^{- 2 \lam \tau'} E_{O, \mm}^2(0) d \tau' 
\notag\\
&\leq e^{-\lam_2 \tau}    \EEt(0) +  C_{ \eqref{eq:EE_tot3} } e^{-\lam_2 \tau} E_{O, \mm}^2(0)  ,
\label{eq:EE_tot3}
\end{align}
where $C_{ \eqref{eq:EE_tot3} }  = 2 C_{ \eqref{eq:EE_tot1} } \lambda^{-1} >0$, 
for all $\tau \in [0,\tau_{\sf max})$.

\subsection{Improving the bootstrap assumptions}
\label{sec:JC:improve:bootstrap}
In this section we show that by requiring  
\begin{equation}\label{eq:bar_eps_cond}
E_{O, \mm}(0) < \bar \epsilon,  \qquad  
\EEt(0) < \bar \epsilon^2, 
\qquad 
\bar \epsilon < \min(1, \eps_*)
\end{equation}
with $\bar \epsilon $ sufficiently small, and $R_0$ 
(appearing in the definition of $\vrs$ in \eqref{eq:vrho_s}) being sufficiently large, the  bootstrap assumptions \eqref{eq:boot1}--\eqref{eq:boot2} can be improved

Assume that the bootstrap assumptions \eqref{eq:boot1}--\eqref{eq:boot2} hold for all $\tau \in [0,\tau_{\sf max})$, for some $\tau_{\sf max}>0$. Due to~\eqref{eq:bar_eps_cond}, assumption~\eqref{eq:EE_tot:ass:b} holds, and thus $\EEt$ satisfies the bounds in~\eqref{eq:EE_tot2} and~\eqref{eq:EE_tot3} for all $\tau \in [0,\tau_{\sf max})$. Using \eqref{eq:Xm_interp}, \eqref{eq:EE_tot2}, and 
assumption \eqref{eq:ass:EOM},  on this time interval we also have 
\begin{equation}\label{eq:boot_impr1}
  \| \twm(\tau)\|_{\cX^{\kk-1}} \les \EEt(\tau)^{1/2} \les  1,
  \qquad 
  E_{O, \mm}(\tau) \les E_{O, \mm}(0) \les 1.
\end{equation}
Moreover, by combining~\eqref{eq:ass:EOM} and~\eqref{eq:EE_tot3} with~\eqref{eq:bar_eps_cond} we deduce that 
\begin{equation}
\label{eq:E:tot:E:0M:sharp}
\EEt^{1/2}(\tau) + E_{O, \mm}(\tau) 
\leq C_{\eqref{eq:E:tot:E:0M:sharp}} \bar \eps e^{- \frac 12 \lambda_2 \tau },
\end{equation}
for all $\tau \in [0,\tau_{\sf max})$, 
where $C_{\eqref{eq:E:tot:E:0M:sharp}} = 1 + C_{\eqref{eq:ass:EOM}} +  C_{ \eqref{eq:EE_tot3} } >0$.

\subsubsection{Improving bootstrap~\eqref{eq:boot2}}
Appealing to \eqref{eq:boot_res2} and \eqref{eq:boot_res3} with $k = 1$,  to the bound  \eqref{eq:boot_impr1}, to the second inequality in~\eqref{eq:boot1_res1:c}, recalling the $\delta_{\bullet}$ definitions in~\eqref{eq:para_del}, and the fact that $\tilde g = g - \bar g = \tl_l \td g + \rl_l \td g$, we obtain 
\begin{subequations}\label{eq:boot_impr2}
\begin{align}
|\na  \tu| 
&\leq C \ccu(y)  \la y \ra^{-1+ 2 \ddd}   \| \twm(\tau)\|_{\cX^{\kk-1}}
+ C E_{O, \mm} \one_{|y| \leq 2}
\notag\\
&\leq  C_{ \eqref{eq:boot_impr2:a} } \la y \ra^{-1 + 3 \ddd + \ddb} (  E_{O, \mm} + \EEt^{1/2} ),
\label{eq:boot_impr2:a}
\end{align}
for some constant $C_{ \eqref{eq:boot_impr2:a} } >0$.
Similarly, for $g \in \{ \vrho , \mb\}$, we may estimate 
\[
   \Big|\frac{\na g}{g} - \frac{\na \bar g}{\bar g}\Big|
  \leq \Big|\frac{\na (g - \bar g)}{g} - \frac{\na \bar g  \cdot (g - \bar g)}{g \bar g} \Big|
  \leq   \frac{|\na (g - \bar g)|}{g} +  \frac{ |\na \bar g| \cdot |g-\bar g| }{g \bar g} . 
\]
Using Remark~\ref{rem:decay}, the bound~\eqref{eq:boot_impr1},  and   estimates \eqref{eq:boot_res2}, \eqref{eq:boot_res3:d}--\eqref{eq:boot_res3:f} with $k \in \{0, 1\}$, for $g \in \{ \vrho, \mb\}$, we obtain 
\begin{align}
     \Big|\frac{\na g}{g} - \frac{\na \bar g}{\bar g}\Big| 
    & \leq C( g^{-1} |\na (g - \bar g)|   +  g^{-1} |y|^{-1}  |g-\bar g|   ) 
    \notag\\
     &\leq  C( |y|^{-1} \ang y ^{4 \ddd}  \| \twm(s)\|_{\cX^{\kk-1}}
     +  |y|^{-1} \ang y ^{ \ddd} E_{O, \mm} ) \notag \\
     & \leq C_{ \eqref{eq:boot_impr2:b} } |y|^{-1} \ang y ^{ 4 \ddd}  (\EEt^{1/2} + E_{O, \mm})
    \label{eq:boot_impr2:b}
\end{align}
\end{subequations}
for some constant $C_{ \eqref{eq:boot_impr2:b} }>0$.

Combining~\eqref{eq:E:tot:E:0M:sharp} with the bounds in~\eqref{eq:boot_impr2}, it follows that if $\bar \eps$ is chosen small enough to ensure
\begin{align}
\bar \eps <  \min \left\{ \frac{1}{C_{\eqref{eq:E:tot:E:0M:sharp}} C_{\eqref{eq:boot_impr2:a}}},  \frac{1}{C_{\eqref{eq:E:tot:E:0M:sharp}} C_{\eqref{eq:boot_impr2:b}}} \right\}
\label{eq:bar_eps_cond:2}
\end{align}
then the bootstrap \eqref{eq:boot2} is \emph{strictly improved} on $[0,\tau_{\sf max})$.

\subsubsection{Improving the $\vrho$ bootstrap in~\eqref{eq:boot1}}
Define the transport operator $\cJ_\tau$ as
\begin{equation}\label{eq:tran_J}
  \cJ_\tau f  := \pa_\tau f + (\bcr y + \bu) \cdot \na f.
\end{equation}
Using the $\vrho$-evolution in~\eqref{eq:sys:a} and the definition of $\vrs$ in~\eqref{eq:vrho_s}, we have
\begin{align}
  \pa_\tau \frac{\vrho}{\; \vrs \; } 
  &+ ( \cx y + \UU) \cdot \na  \frac{\vrho}{\; \vrs \;} 
  = \Big( \crho - 2\alpha \div \UU - \frac{ ( \cx y + \UU )\cdot \na \vrs}{\vrs} - \frac{\pa_\tau \vrs}{\vrs} \Big)   \, \frac{\vrho}{\; \vrs \;} 
  \notag\\
  &=- \Big(\underbrace{\cJ_\tau  \vrs  + 2\alpha \div \bu\vrs  - \crhobar \vrs}_{=: I_1}   \Big)   \frac{\vrho}{{\vrs}^2}
  + \Big( \underbrace{\crhotilde \vrs - 2\alpha \div \tu \vrs-  ( \cxtilde y + \tu )\cdot \na \vrs }_{=: I_2} \Big)   \frac{\vrho}{{\vrs}^2}
  .
  \label{eq:boot_vrho1}
\end{align}

Upon recalling the definition of $\vrs$ from \eqref{eq:vrho_s}, we have that
\begin{subequations}
  \label{eq:boot_vrho1:new}
\begin{align}
\nabla \vrs(y,\tau)
=  
\chi_{\msf{in}}(y,\tau) \nabla \bvr(y) 
+
\nabla \chi_{\msf{in}}(y,\tau) \bigl( \bvr(y) - \bvr(R(\tau)) \bigr),
\end{align}
where
\[
\chi_{\msf{in}}(y,\tau) := \chi \Bigl(  \tfrac{|y|}{\rs} \Bigr),
\qquad\mbox{and}\qquad
\rs = e^{ \bcr \tau } R_0.
\]
Moreover, since $\bvr$ (with $\crhobar$ and $\bu$) is the steady state to~\eqref{eq:sys:a}, we have
\begin{align}
 \cJ_\tau  \vrs(y,\tau)
&= - \frac{\cxbar |y|}{R(\tau)} \chi^\prime\Bigl(  \tfrac{|y|}{\rs} \Bigr) \bigl( \bvr(y) - \bvr(R(\tau)) \bigr)
+ \bigl( 1 - \chi_{\msf{in}}(y,\tau)\bigr) \cxbar R(\tau) \bigl( \tfrac{y}{|y|}\cdot\nabla \bvr\bigr)(R(\tau))
\notag\\
&\quad + \chi_{\msf{in}}(y,\tau) (\bcr y + \bu) \cdot  \nabla \bvr(y) 
+
 \frac{(\bcr y + \bu) \cdot y}{|y| R(\tau)} \chi^\prime\Bigl(  \tfrac{|y|}{\rs} \Bigr) \bigl( \bvr(y) - \bvr(R(\tau)) \bigr)
 \notag\\
 &=  \bigl( 1 - \chi_{\msf{in}}(y,\tau)\bigr)    \bigl( \crhobar \bvr - 2\alpha \div \bu \bvr - \bu \cdot \nabla \bvr \bigr)(R(\tau))
\notag\\
&\quad + \chi_{\msf{in}}(y,\tau) \bigl(\crhobar - 2\alpha \div \bu \bigr) \bvr(y) 
+
 \frac{  \bu  \cdot y}{|y| R(\tau)} \chi^\prime\Bigl(  \tfrac{|y|}{\rs} \Bigr) \bigl( \bvr(y) - \bvr(R(\tau)) \bigr).
\end{align}
\end{subequations}
From~\eqref{eq:vrho_s} and \eqref{eq:boot_vrho1:new}, we deduce that the terms $I_1$ and $I_2$ defined in~\eqref{eq:boot_vrho1} may be rewritten as
\begin{subequations}
\label{eq:boot_res2:new}
\begin{align}
I_1 &=
- \bigl( 1 - \chi_{\msf{in}}(y,\tau)\bigr) \bigl(  \bu \cdot \nabla \bvr \bigr)(R(\tau)) 
+ \bu(y) \cdot \nabla  \chi_{\msf{in}}(y,\tau) \bigl( \bvr(y) - \bvr(R(\tau)) \bigr)
\\
I_2 &=
\bigl( \crhotilde - 2\alpha \div \tu \bigr) \vrs(y,\tau)
-  \chi_{\msf{in}}(y,\tau) ( \cxtilde y + \tu )\cdot  \nabla \bvr(y) 
-  ( \cxtilde y + \tu )\cdot \nabla \chi_{\msf{in}}(y,\tau) \bigl( \bvr(y) - \bvr(R(\tau)) \bigr)
\end{align} 
\end{subequations}
The bound for $I_2$ is slightly easier to obtain; using \eqref{eq:boot_res3:a}, \eqref{eq:cutoff}, \eqref{eq:vrho_s_est}, \eqref{norm:ODE},  \eqref{eq:boot_res2}, 
\eqref{eq:boot1_res1:c}, \eqref{eq:boot_impr1}, and~\eqref{eq:para_del}, we obtain 
\begin{align}
|I_2| 
&\leq
|\crhotilde - 2\alpha \div \tu  | \vrs 
+  \chi_{\msf{in}}  \bigl| ( \cxtilde y + \tu )\cdot  \nabla \bvr  \bigr| 
+ \bigl| ( \cxtilde y + \tu )\cdot \nabla \chi_{\msf{in}} \bigr| \bigl| \bvr  - \bvr(R(\tau)) \bigr|
\notag\\
&\les
\bigl( |\crhotilde | + |\cxtilde| + |\div \tu |  + |y|^{-1} |\tu| \bigr)\vrs
+  \bigl( |\crhotilde |  + |y|^{-1} |\tu| \bigr) \bvr(R(\tau)) \one_{ |y| \in [\rs, 2 \rs]}  
\notag\\
&\les 
\bigl(E_{O, \mm} + |y|^{-1} \cc(y) \brak{y}^{4\ddd} \| \twm\|_{ \cX^{\kk-1} } \bigr) \vrs
\notag\\
&\les 
\bigl(1 +  \brak{y}^{\ddb + 5\ddd - 1} \bigr) \bigl(E_{O, \mm} +   \EEt^{1/2} \bigr) \vrs
\notag\\
&\leq C_{\eqref{eq:boot_vrho2}} \bigl(E_{O, \mm} +   \EEt^{1/2} \bigr) \vrs, 
\label{eq:boot_vrho2}
\end{align}
for a computable constant $C_{\eqref{eq:boot_vrho2}}>0$. 
Here we have implicitly used the fact that  $\bvr(y), \bvr(\rs) \les \vrs(y,\tau)$ for all $|y| \in [\rs,2\rs]$.
Similarly, for the term $I_1$  in~\eqref{eq:boot_res2:new}, we use that $\rs \geq R_0 \geq 1$, and then appeal to~\eqref{eq:vrho_s}, \eqref{eq:boot_res3:a}, \eqref{eq:boot_res3:b}, \eqref{eq:profile_decay},  
in order to obtain
\begin{align}
|I_1|
&\leq
\bigl( 1 - \chi_{\msf{in}}(y,\tau)\bigr) \bigl|  (\bu \cdot \nabla \bvr )(R(\tau)) \bigr|
+\bigl| \bu(y) \cdot \nabla  \chi_{\msf{in}}(y,\tau) \bigr| \, \bigl| \bvr(y) - \bvr(R(\tau)) \bigr|
\notag\\
&\les \bigl( \brak{\rs}^{\frac{\cubar}{\cxbar} - 1} 
+    \brak{y}^{\frac{\cubar}{\cxbar} - 1}  \one_{ |y| \in [\rs, 2 \rs]}  \bigr) \vrs 
\notag\\
&\leq C_{\eqref{eq:boot_vrho3}}  \rs^{\frac{\cubar}{\cxbar} - 1} 
\vrs,
\label{eq:boot_vrho3}
\end{align}
for a computable constant $C_{\eqref{eq:boot_vrho3}}>0$. 

Thus, by combining~\eqref{eq:boot_vrho2} and \eqref{eq:boot_vrho3}, \eqref{eq:E:tot:E:0M:sharp},  we obtain that 
\begin{equation}\label{eq:boot_vrho_D}
  | {\sf RHS}_{\eqref{eq:boot_vrho1}} | 
  \leq \frac{\vrho}{\, \vrs \,}
  \Bigl( C_{\eqref{eq:boot_vrho2}} \bigl(E_{O, \mm} +   \EEt^{1/2} \bigr) +  C_{\eqref{eq:boot_vrho3}}  \rs^{\frac{\cubar}{\cxbar} - 1} \Bigr).
\end{equation}
Note however that closing the $\vrho$-bootstrap~\eqref{eq:boot1}, amounts to controlling not just the ratio $ \vrho/ \vrs $, but this ratio multiplied with $\brak{y}^{\pm \ddd}$ (cf.~\eqref{eq:boot1:old}).

For this purpose, we multiply \eqref{eq:boot_vrho1} by $\la y \ra^{\theta}$, with $\theta  = \pm \ddd$, and obtain 
\begin{equation}\label{eq:boot_vrho_eqn2}
  \pa_\tau \Bigl(\frac{\vrho}{\, \vrs \, } \la y \ra^{\theta} \Bigr)
  + (\cx y + \UU) \cdot \na  \Bigl(\frac{\vrho}{\, \vrs \, } \la y \ra^{\theta}\Bigr) 
  = \Bigl( \frac{\, \vrs \,} {\vrho}{\sf RHS}_{\eqref{eq:boot_vrho1}}  
  +  \theta \frac{ ( \cx y + \UU) \cdot y}{ \la y \ra^2 } \Bigr) \Bigl(\frac{\vrho}{\, \vrs \, } \la y \ra^{\theta}\Bigr).
\end{equation}

When $\theta >0$, using estimates similar to those in \eqref{eq:boot_vrho2}, \eqref{eq:boot_vrho3}, 
and the bound $\bcr |y| + \bar U(|y|) >0 $ (cf.~\eqref{eq:outgoing_strong}), we derive 
\begin{align}
\theta \frac{ ( \cx y + \UU) \cdot y}{ \la y \ra^2 } 
&= \theta \frac{ ( \cxbar |y| + \bar U(|y|)) |y|}{ \la y \ra^2 }   
+  \theta \frac{ ( \cxtilde y + \td \UU) \cdot y}{ \la y \ra^2 } 
\notag\\
&\geq  \theta \frac{ ( \cxtilde y + \td \UU) \cdot y}{ \la y \ra^2 }
\geq - C_{\eqref{eq:boot_vrho_D2}} \theta ( \EEt^{1/2} + E_{O, \mm} ),
\label{eq:boot_vrho_D2}
\end{align}
for some computable constant $C_{\eqref{eq:boot_vrho_D2}}>0$.
Similarly, for $\theta < 0$ we have
\begin{align}
\theta \frac{ ( \cx y + \UU) \cdot y}{ \la y \ra^2 } 
&= \theta \frac{ ( \cxbar |y| + \bar U(|y|)) |y|}{ \la y \ra^2 }   
+  \theta \frac{ ( \cxtilde y + \td \UU) \cdot y}{ \la y \ra^2 } 
\notag\\
&\leq  \theta \frac{ ( \cxtilde y + \td \UU) \cdot y}{ \la y \ra^2 }
\leq -  C_{\eqref{eq:boot_vrho_D3}}  \theta ( \EEt^{1/2} + E_{O, \mm} ),
\label{eq:boot_vrho_D3}
\end{align}
for some computable constant $C_{\eqref{eq:boot_vrho_D3}}>0$.

By combining~\eqref{eq:boot_vrho_D}, \eqref{eq:boot_vrho_eqn2}, \eqref{eq:boot_vrho_D2}, and \eqref{eq:boot_vrho_D3}, we thus deduce
\begin{subequations}
\label{eq:boot_vrho_almost_final}
\begin{align}
\bigl(\p_\tau + (\cx y + \UU)\cdot \nabla \bigr) \log \Bigl(\frac{\vrho}{\, \vrs \, } \la y \ra^{\ddd}\Bigr)
&\geq 
- C_{\eqref{eq:boot_vrho_almost_final}} \bigl(E_{O, \mm} +   \EEt^{1/2} + \rs^{\frac{\cubar}{\cxbar} - 1} \bigr) 
\\
\bigl(\p_\tau + (\cx y + \UU)\cdot \nabla \bigr) \log \Bigl(\frac{\vrho}{\, \vrs \, } \la y \ra^{-\ddd}\Bigr)
&\leq 
C_{\eqref{eq:boot_vrho_almost_final}} \bigl(E_{O, \mm} +   \EEt^{1/2} + \rs^{\frac{\cubar}{\cxbar} - 1} \bigr) 
\end{align} 
where
$C_{\eqref{eq:boot_vrho_almost_final}} = C_{\eqref{eq:boot_vrho2}}  + C_{\eqref{eq:boot_vrho3}} + C_{\eqref{eq:boot_vrho_D2}} \ddd   + C_{\eqref{eq:boot_vrho_D3}} \ddd >0$.
\end{subequations}
Note that under the bootstrap assumption \eqref{eq:boot1}--\eqref{eq:boot2}, and the smallness of initial perturbation \eqref{eq:bar_eps_cond}  we have \eqref{eq:E:tot:E:0M:sharp}, and the bound \eqref{eq:ass:EOM} was assumed. Recalling that $\rs = e^{\bcr \tau} R_0 $, with $R_0\geq 1$, and that $\cxbar - \cubar = 1$, it follows that 
\begin{align}
\int_0^\tau \Bigl( E_{O, \mm}(\tau') +   \EEt^{1/2}(\tau') + R(\tau')^{\frac{\cubar}{\cxbar} - 1}\Bigr) d\tau'
&\leq 
\int_0^\infty \Bigl( C_{\eqref{eq:E:tot:E:0M:sharp}} \bar \eps e^{- \frac 12 \lambda_2 \tau' } + R_0^{-\frac{1}{\cxbar}} e^{- \tau'}\Bigr)
d\tau'
\notag\\
&\leq 
\frac{2C_{\eqref{eq:E:tot:E:0M:sharp}}}{\lambda_2} \bar \eps
+ R_0^{- \frac{1}{\cxbar}}
.
\label{eq:boot_vrho_D3:new}
\end{align}
Thus, if we choose $\bar \eps$ to be small enough to ensure
\begin{equation}
\frac{2 C_{\eqref{eq:boot_vrho_almost_final}}  C_{\eqref{eq:E:tot:E:0M:sharp}}}{\lambda_2} \bar \eps
< \frac{\log(3/2)}{2}
\label{eq:bar_eps_cond:3}
\end{equation}
and $R_0\geq 1$ to be large enough to ensure
\begin{equation}
C_{\eqref{eq:boot_vrho_almost_final}} 
R_0^{- \frac{1}{\cxbar}}
< \frac{\log(3/2)}{2},
\label{eq:bar_eps_cond:4}
\end{equation}
we obtain from~\eqref{eq:boot_vrho_D3:new}
that
\[
C_{\eqref{eq:boot_vrho_almost_final}} \int_0^\tau \Bigl( E_{O, \mm}(\tau') +   \EEt^{1/2}(\tau') + R(\tau')^{\frac{\cubar}{\cxbar} - 1}\Bigr) d\tau'
< \log (3/2).
\]
Then, integrating \eqref{eq:boot_vrho_almost_final} along the characteristics of $\cx y + \UU$, and using the assumption  \eqref{ass:nonrad_init:a} on the initial data, we obtain the pointwise bounds 
\begin{subequations}
\label{eq:boot_vrho_final}
\begin{align}
 \frac{\vrho}{\, \vrs\, }(y,\tau) \la y \ra^{\ddd}
&\geq e^{-\log (3/2)} \min_{y}\Big( \frac{\vrho}{\, \vrs\, }(y,0) \la y \ra^{\ddd} \Big) 
\geq \frac 13, \\
 \frac{\vrho}{\, \vrs\, }(y,\tau) \la y \ra^{-\ddd}
&\leq e^{\log (3/2)} \max_{y}\Big( \frac{\vrho}{\, \vrs\, }(y,0) \la y \ra^{-\ddd} \Big) 
\leq 3,
\end{align}
\end{subequations}
thereby improving the $\vrho$ bootstrap in~\eqref{eq:boot1}.

\subsubsection{Improving the $\mb$ bootstrap in~\eqref{eq:boot1}}
In order to estimate $\mb/\barb$ we appeal to \eqref{eq:sys:c} 
and the stationary form of~\eqref{eq:sys:c}; this argument is very similar, but much simpler. First, we have
\begin{equation}
\bigl(\p_\tau + (\cx y + \UU)\cdot \nabla \bigr)
\log \Bigl(\frac{\mb}{\barb} \brak{y}^\theta \Bigr)
=  \tcb - \frac{(\cxtilde y + \tu)\cdot \nabla\barb}{\barb} +  \theta \frac{ ( \cx y + \UU) \cdot y}{ \la y \ra^2 }.  
\label{eq:boot_b_med}
\end{equation}
Second, analogously to~\eqref{eq:boot_vrho2}--\eqref{eq:boot_vrho_D3}, we may show that 
\begin{subequations}
\label{eq:boot_b_med2}
\begin{align}
{\sf RHS}_{\eqref{eq:boot_b_med}}
\geq  
- C_{\eqref{eq:boot_b_med2}} \bigl(E_{O, \mm} +   \EEt^{1/2} + \rs^{\frac{\cubar}{\cxbar} - 1} \bigr), \qquad \theta = \ddd, \\
{\sf RHS}_{\eqref{eq:boot_b_med}}
\leq  
C_{\eqref{eq:boot_b_med2}} \bigl(E_{O, \mm} +   \EEt^{1/2} + \rs^{\frac{\cubar}{\cxbar} - 1} \bigr), \qquad \theta = - \ddd,  
\end{align}
\end{subequations}
for some computable constant $C_{\eqref{eq:boot_b_med2}}>0$.
Thus, if we choose $\bar \eps$ to be small enough to ensure
\begin{equation}
\frac{2 C_{\eqref{eq:boot_b_med2}}  C_{\eqref{eq:E:tot:E:0M:sharp}}}{\lambda_2} \bar \eps
< \frac{\log(3/2)}{2}
\label{eq:bar_eps_cond:5}
\end{equation}
and $R_0\geq 1$ to be large enough to ensure
\begin{equation}
C_{\eqref{eq:boot_b_med2}}
R_0^{- \frac{1}{\cxbar}}
< \frac{\log(3/2)}{2},
\label{eq:bar_eps_cond:6}
\end{equation}
we conclude from~\eqref{eq:boot_b_med}--\eqref{eq:boot_b_med2}, and using the assumptions \eqref{ass:nonrad_init:b}, that
\begin{subequations}
\label{eq:boot_b_final} 
\begin{align}
 \frac{\mb}{\, \barb \, }(y,\tau) \la y \ra^{\ddd}
&\geq e^{-\log (3/2)} \min_{y}\Big( \frac{\mb}{\, \barb\, }(y,0) \la y \ra^{\ddd} \Big) 
\geq \frac 13, \\
 \frac{\mb}{\, \barb \, }(y,\tau) \la y \ra^{-\ddd}
&\leq e^{\log (3/2)} \max_{y}\Big( \frac{\mb}{\, \barb \, }(y,0) \la y \ra^{-\ddd} \Big) 
\leq 3,
\end{align}
\end{subequations}
thereby improving the $\mb$ bootstrap in~\eqref{eq:boot1}.

\subsection{The proof of Theorem~\ref{thm:nonrad:stab}}
\label{sec:proof:thm:nonrad:stab}
Let $\mm$ be as in Lemma~\ref{lem:nonrad_wg}. Assume that the modulation functions and the Taylor coefficients satisfy the bound~\eqref{eq:ass:EOM}, with $E_{O,\mm}$ as defined in~\eqref{eq:E:O:L}.
Let $R_0\geq 1$ be large enough to ensure that~\eqref{eq:bar_eps_cond:4} and \eqref{eq:bar_eps_cond:6} hold, and let $\bar \eps>0$ to be small enough to ensure that~\eqref{eq:bar_eps_cond}, \eqref{eq:bar_eps_cond:2}, \eqref{eq:bar_eps_cond:3}, and \eqref{eq:bar_eps_cond:5} hold. Define $\nu>0$ as in~\eqref{energy:tot}, and define $\ddd>0$ as in~\eqref{eq:para_del}.
Lastly, let $\vp_0$ be as defined in Lemma~\ref{lem:nonrad_wg}.

Assuming that $(\vrho_{\sf in},\mb_{\sf in})$ obeys the bounds in~\eqref{ass:nonrad_init}, and take  $\EEt(0) < \bar \epsilon^2$ and $E_{O,\mm}(0) < \bar \epsilon$. Under these assumptions,  
we have proven in Section~\ref{sec:JC:improve:bootstrap} that the bootstraps~\eqref{eq:boot1} and~\eqref{eq:boot2} are strictly improved. Thus, these bootstraps hold for all $\tau \in [0,\infty)$. In particular, the bounds in~\eqref{eq:boot1} are exactly the pointwise estimates claimed in Theorem~\ref{thm:nonrad:stab}. Finally, estimate~\eqref{eq:nonrad_decay} follows from \eqref{eq:bar_eps_cond} and \eqref{eq:EE_tot3}. This completes the proof of Theorem~\ref{thm:nonrad:stab}.

\section{Sharp stability analysis of the Taylor coefficients outside radial symmetry}
\label{sec:ODE}

In Section~\ref{sec:nonradial} we have established a global-in-$\tau$ exponential decay estimate for the bulk part of a non-radial perturbation $(\tvr,\tu,\tb)$ of the globally self-similar solution $(\bar \vrho, \bar \UU,\barb)$, see Theorem~\ref{thm:nonrad:stab}, conditional on assumption~\eqref{eq:ass:EOM}: the functional $E_{O,\mm}(\tau)$, which collects the modulation functions $(\cxtilde,\cutilde,\cbbtilde,\crhotilde)$ and the Taylor coefficients of $(\tvr,\nabla \tu, \nabla^2 \tb)$ at $y=0$ up to order $\mm$, decays exponentially in $\tau$. The purpose of the present section is to characterize \emph{sharply} the set of initial data for which the assumption~\eqref{eq:ass:EOM} holds.

Our analysis  proceeds in four steps. First, we show that the Taylor coefficients $V_{\leq n}(\tau)$ of the perturbation at $y=0$, together with the modulation functions, satisfy a \emph{closed and finite-dimensional system of ODEs} (cf.~Theorem~\ref{thm:ODE_stab}, item~(ii)). In particular, controlling these coefficients does not require any information from the bulk PDE. Second, we provide an \emph{explicit and computable upper bound} for the dimension of the unstable subspace $\Sigma_{\mathsf{uns}}$ of this ODE system, of the form $\dim(\Sigma_{\mathsf{uns}}) \leq \bar C(d) \NNN^d$, valid for every $d\in\{1,2,3\}$, every adiabatic exponent $1<\gamma\leq 2d+1$, and every $\NNN\geq 1$. Third, we prove that the unstable subspace \emph{stabilizes at a finite order}: for all orders $n \geq n_1 := (18d)^2 \NNN$, the unstable subspace $\Sigma_{\mathsf{uns},\leq n}$ is a trivial extension of $\Sigma_{\mathsf{uns},\leq n_1}$, an $n$-independent subspace whose definition does not require knowledge of arbitrarily many Taylor coefficients (cf.~Theorem~\ref{thm:ODE_stab}, item~(iii)). Fourth, in five distinguished cases that include the physically most relevant monatomic and diatomic gases for the ground state $\NNN=1$, we determine $\dim(\Sigma_{\mathsf{uns}})$ \emph{explicitly}, and we show that the unstable directions are confined to the Taylor coefficients of order $\leq 2$ at the origin (cf.~Theorem~\ref{thm:ODE_stab}, items~(iv)--(v)).
Combining these four steps with Theorem~\ref{thm:nonrad:stab}, we obtain the full nonlinear PDE stability picture for the implosions constructed in Section~\ref{sec:profiles}, outside of radial symmetry. 

\subsection{Main results of stability estimates for the ODEs}

The main result of this section is Theorem~\ref{thm:ODE_stab} below. Before stating this result, we first need to introduce the following notation.

\vspace{0.1in}
\paragraph{\bf Taylor coefficients at the origin}
We recall from~\eqref{def:cC_norm} that
\[
| f(y)|_{\cC^n} = {\sum}_{ |\bfa| \leq n } | \pa^{\bfa} f(y)|
,
\]
controls the Taylor coefficients up to order $n$ at a given point. We also recall from~\eqref{norm:ODE} that the functional 
\[
  E_{O, n}(\tau) =    |\tvr(0,\tau)|_{\cC^{n}}
 + |\tu(0,\tau)|_{\cC^{n +1}}  + |\tb(0,\tau)|_{\cC^{n+2}} 
+ |\tcr(\tau) | + |\tcu(\tau)| + |\tcb(\tau)| + |\tcvr(\tau)|
\]
controls the Taylor coefficients of order $\leq n$ at the origin, and also the modulation parameters from~\eqref{eq:lin:full}.

For any $n\geq 0$, we denote by $\nabla^n f$ the tensor of $n$-th order mixed partial derivatives of $f$.\footnote{For $n\geq0$, we let $\nabla^n f$ denote a row-vector which enumerates the elements of the set $\{ \p^\bfa f \colon \bfa \in \Naturals_0^d, |\bfa|=n \}$; this is instead of the usual meaning, of a symmetric $n$-tensor. The canonical enumeration of the aforementioned set is that following graded lexicographic order on $d$-dimensional multi-indices of length $n$. At no point in the proof will we use a specific ordering for the elements of $\nabla^n f$, so we do not insist on this detail.}
We denote 
\begin{subequations}\label{def:TL_nth_order} 
\begin{align}
V_{=n}(\tau) & := ( \nabla^{n} \tvr(0,\tau), \nabla^{  n+1} \tu(0,\tau), \nabla^{n+2} \tb(0,\tau) )  , \\
V_{\leq n}(\tau) & := (V_{=0}(\tau), V_{=1}(\tau),.., V_{=n}(\tau) ) 
 \in \Reals^{d_{\leq n} }
\end{align}
\end{subequations}
where  $d_{\leq n}$ denotes the dimension of the vector $ V_{\leq n}$.\footnote{
\label{foot:d:leq:n} If $f$ is a scalar function, then the length of the vector $\nabla^n f$ is ${n + d - 1 \choose d-1}$. Consequently, given that $\tu$ is a $d$-vector, we obtain that $V_{=n}$ has length ${n + d - 1 \choose d-1} + d {n + d \choose d-1} + {n + d +1 \choose d-1}$. Summing the lengths of $V_{=j}$ for $0\leq j \leq n$, we obtain that $d_{\leq n} = {n+d \choose d} + d {n+d+1 \choose d} + {n+d+2 \choose d} - (2d +1)$. The specific value of $d_{\leq n}$ is not used in the proof, so we do not insist on this detail.} We view $V_{=n}, V_{\leq n}$ as row vectors. Recall that by  \eqref{eq:vanishing:order:at:origin:SS}, we have 
\[
   \tu(0, \tau ) = 0,  \quad \tb(0, \tau ) = 0, \quad \na \tb(0, \tau ) = 0.
\] 
Thus, we do not include $  \tu(0, \tau )$ and $\na^{\leq 1} \tb(0, \tau )$ in the definition and analysis of $V_{\leq n}(\tau)$. 

For compactness of notation, in this section, we sometimes write $f(0)$ instead of $f(0,\tau)$, or $f$ instead of $f(\tau)$, whenever confusion cannot arise.

\vspace{0.1in}
\paragraph{\bf Multiplicity of eigenvalues}
Given an eigenvalue $\lam$ of a matrix $M$, we denote by 
\bseq\label{def:Algm}
\beq
  \Algm(\lam, M)
\eeq
the \textit{algebraic multiplicity} of the eigenvalue $\lam$. We also introduce the notation $ \Algm_{\geq 0}$ to denote the sum of the algebraic 
multiplicities of all eigenvalues of $M$ with non-negative real part:
\beq
 \Algm_{\geq 0}(M) := \sum_{\substack{\lam \textrm{ \ is an eigenvalue of } M\\  \Re(\lam) \geq 0}} \Algm( \lam, M).
\eeq
\eseq

\vspace{0.1in}
\paragraph{\bf The Main Result}
We have the following finite co-dimension stability of the modulation functions and the ODE system for the high-order Taylor coefficients at $y=0$. 

\begin{theorem}[\bf Finite co-dimension stability for ODEs at the origin]
\label{thm:ODE_stab}
Fix $d \in \{ 1, 2, 3 \}, 1 < \gamma \leq 2d+1$, $\NNN \geq 1$, and let $\alpha = \frac{\gamma-1}{2}$. 

\begin{enumerate}[label=(\roman*),leftmargin=2em]
\item \textsl{(Modulation functions).} The modulation functions $\tcr, \tcvr, \tcu, \tcb $ appearing in~\eqref{eq:lin:full} 
are determined as linear combinations of $\tvr(0, \tau), \div\tu(0, \tau), \Delta \tb(0, \tau), \Delta^{\NNN } \tvr(0,\tau), \Delta^{\NNN+1} (y \cdot \tu)(0,\tau),
\Delta^{\NNN+1} \tb(0, \tau) $. The precise definition is given in \eqref{eq:modulated:nosym:summary}. 

\item \textsl{(ODE system).}
For any $ n  \geq 2 \NNN$,\footnote{Our proof shows that the ODE system for $V_{\leq n}$ is closed for any $n\geq 1$, \emph{were it not for the modulation functions}. The fact that the modulation functions require information about $V_{=2\NNN}$ means that the ODE system is closed only when $n\geq 2\NNN$.} the Taylor coefficients at $y=0$ of the solution $(\tvr,\tu,\tb)$ to \eqref{eq:lin:full} form a closed ODE system
 \begin{equation}\label{eq:ODE_stab}
   \tfrac{d}{d \tau} V_{\leq n}(\tau)^\intercal = M_{\leq n} V_{\leq n}(\tau)^\intercal  + Q_n(V_{\leq n}(\tau), V_{\leq n}(\tau) ) ,
 \end{equation}
for some constant (in $\tau$) matrix $M_{\leq n} \in \Reals^{d_{\leq n} \times d_{\leq n} }$, and some bilinear form $Q_n: \Reals^{d_{\leq n}} \times \Reals^{d_{\leq n} } \to \Reals^{d_{\leq n}} $ with constant (in $\tau$) coefficients); here $d_{\leq n}$ is the dimension of the vector $V_{\leq n}$.

\item \textsl{(Nonlinear stability estimates of the ODE system).}  Suppose that the matrix $M_{\leq n}$ has $ k_{\leq n}$ different eigenvalues with non-negative real parts, which are given by $\{ \lam_{ M_{\leq n}, i}\}_{ i=1}^{k_{\leq n}} $, and denote their algebraic multiplicity by $ \Algm( \lam_{ M_{\leq n}, i}, M_{\leq n}   )$, for $1\leq i \leq k_{\leq n}$. We recall $M_{\leq n} \in \Reals^{d_{\leq n}, d_{\leq n}}$. We define the following real linear subspace of $\Reals^{d_{\leq n}}$:
\beq\label{def:Sigma_uns_n}
 \Sigma_{\mathsf{uns}, \leq n} := \Re \bigoplus_{1\leq i \leq k_{\leq n} } \ker\Bigl(  (\lam_{ M_{\leq n} , i} \Id - M_{\leq n} )^{ 
  \Algm( \lam_{ M_{\leq n}, i}, M_{\leq n}   ) } \Bigr)  \subset  \Reals^{d_{\leq n}}.
\eeq
For any $\NNN \geq 1$, let 
\[
n_1:= (18 d)^2 \NNN.
\]
There exists an explicitly computable constant $\bar C(d)\geq 1$ such that the following statement holds. For any $n \geq n_1 $, a column vector 
$ v_{\leq n} $ belongs to $ \Sigma_{\mathsf{uns}, \leq n}$ if and only if 
\[
 v_{\leq n} = ( v_{\leq n_1}^\intercal, 0,.., 0)^\intercal , 
\]
for some column vector $v_{\leq n_1} \in  \Sigma_{\mathsf{uns}, \leq n_1}$. As a result, for any $n \geq n_1$, $\Sigma_{\mathsf{uns}, \leq n}$ is a trivial lifting of 
the $n$-independent space $\Sigma_{\mathsf{uns}, \leq n_1}$ from $\Reals^{ d_{\leq n_1}}$ to $\Reals^{ d_{\leq n}}$, and we have
\[
\dim ( \Sigma_{\mathsf{uns}, \leq n}) = \dim ( \Sigma_{\mathsf{uns}, \leq n_1}),
\qquad\mbox{with} \qquad
\dim ( \Sigma_{\mathsf{uns}, \leq n_1}) 
\leq \bar C(d) \NNN^{d}.
\]

Moreover, there exists $\delta = \delta(n, M_{\leq n} ) > 0$ sufficiently small such that given any 
initial data $V_{1, 0} \in \Reals^{d_{\leq n}}$ with $\| V_{ 1, 0}  \|_2 < \delta$, there exists 
an initial data in the $n$-independent subspace $V_{\mathsf{uns}, 0} \in \Sigma_{\mathsf{uns}, \leq n_1}$ and a global-in-$\tau$ solution $V_{\leq n}$ to the ODE system \eqref{eq:ODE_stab}, with the initial condition
\[
   V_{\leq n}(0) = V_{ 1 , 0} + ( V_{\mathsf{uns}, 0}, 0,..,0) , 
   \qquad 
   \|  V_{\mathsf{uns}, 0} \|_2 \les_n  \| V_{ 1, 0} \|_2,  
\]
and exponential decay estimate
\begin{equation}\label{eq:ODE_decay}
    E_{O, n }(\tau) \leq C_n  e^{-\lam \tau} E_{O, n }(0) , \quad  \forall \ \tau \geq 0,
\end{equation}
where $\lam > 0$ only depends on 
$n, M_{\leq n}$, and the constant $C_n > 0$ depends on the profile and on $n\geq n_1$.

\item \textsl{(Quantitative characterization of  $\Sigma_{\mathsf{uns}}$).}  Consider the ground state $\NNN = 1$ and five special cases: 
$(\gamma, d) \in \bigl\{ (\tf53, 3) , \ (\tf75, 3), \ ( 2, 2 ), ( \tf53, 2)\bigr\}$, or  $\gamma \in (1,3], d = 1$. The same structural result for  $\Sigma_{\mathsf{uns} , \leq n}$ 
and nonlinear stability results in (iii) hold with $n_1 = 2$. In particular, 
 for any $n\geq 2$, $\Sigma_{\mathsf{uns} , \leq n}$ is a trivial lifting of 
$\Sigma_{\mathsf{uns} , \leq 2}$, and  $  \dim( \Sigma_{\mathsf{uns}, \leq n} )= \dim( \Sigma_{\mathsf{uns}, \leq 2} ) $. 
Moreover,  we have\footnote{Recall cf.~\eqref{def:Sigma_uns_n} that $\Sigma_{\mathsf{uns, \leq 2} } \subset \Reals^{d_{\leq 2}}$. Also, note that $d_{\leq 2} = {d+2 \choose d} + d {d+3 \choose d} + {d+4 \choose d} - (2d+1)$. In particular, for $d=1$ we have $d_{\leq 2} = 9$, for $d=2$ we have $d_{\leq 2} = 36$, and for $d=3$ we have $d_{\leq 2} = 98$. This highlights the fact that the dimension of $\Sigma_{\mathsf{uns,\leq 2}}$  computed in item (iv) is \emph{much smaller} than the full dimension $d_{\leq 2}$.} 
 \begin{align*}
  \dim( \Sigma_{\mathsf{uns, \leq 2} } ) & = d^2 - 1 + d  , && \mw{for \ } (\gamma, d) \in \bigl\{ (\tf53, 3) ,  \ ( 2, 2 )\bigr\}, \ \mw{or}  \ \gamma \in (1,3], d = 1 ,  \label{eq:dim_Sigma_uns:a}  \\ 
\dim( \Sigma_{\mathsf{uns}, \leq 2 } ) & = 18  , \quad && \mw{for \ } (\gamma, d) = (\tf75, 3), \\
\dim( \Sigma_{\mathsf{uns, \leq 2}} )  & = 7 , && \mw{ for \ }    (\gamma, d) = ( \tfrac53, 2 ) .
\end{align*}
In other words, the  unstable or neutrally stable modes for the ODE system of $V_{\leq n}$ 
\eqref{eq:ODE_stab} only arise from a linear subspace of 
$\{ V_{\leq 2} \in \Reals^{d_{\leq 2 }} \}$, for any $n\geq 2$.

\item \textsl{(Full stability for initial data with higher vanishing order)}
Let $\NNN=1$. Consider five special cases: 
$(\gamma, d) \in \bigl\{ (\tf53, 3) , \ (\tf75, 3), \ ( 2, 2 ), ( \tf53, 2)\bigr\}$, or  $\gamma \in (1,3], d = 1$.
For any $n \geq 2$, there exists $\delta(n, M_{\leq n}) > 0$ such that the following statement holds. 
If the initial data for the mixed-derivatives satisfies $V_{\leq 2}(0) = 0$ 
and $ \| V_{\leq n}(0) \|_{2} \leq \delta$, then
$V_{\leq 2}(\tau) = 0 $ for any $\tau \geq 0$ and  
\[
    E_{O, n }(\tau ) \leq C_n  e^{-\lam \tau } E_{O, n  }(0) , \quad  \forall \ \tau \geq 0 ,
\]
for some $\lam >0$ depending on the matrix $M_{ \leq n}$ and $n\geq 2$.
\end{enumerate}
\end{theorem}

\begin{remark}[\bf Complete characterization of  $\Sigma_{\mathsf{uns}, \leq n}$]
In Section \ref{sec:lower_block} we derive a relation between the matrix $M_{\leq n}$ in \eqref{eq:ODE_stab} and the lifted matrices $\{\mathbf{M}_i\}_{1\leq i\leq n}$ defined in \eqref{eq:upper_block:a}. 
Each matrix $\mathbf{M}_i$ is a block lower
triangular matrix with diagonal blocks given by the matrices $ \HHH_{i/2} $ in \eqref{eq:ODE_Nth:spec:a},  $\HH_{|\bet|, k}$ in \eqref{eq:ODE_Nth}, $\HHT_i$ in \eqref{eq:ODE_Nth_ZB}, or the scalar $ (- i \kp - 2\uua -1 )$ in \eqref{eq:ODE_Nth_U}. 
 Since these matrices on the diagonal blocks are given \emph{explicitly} 
and have size $m \times m$ with $m\leq 4$, we can obtain explicit closed-form formulas for the eigenvalues of each $\mathbf{M}_i$.\footnote{
Polynomials of degree at most four have explicit closed-form formulas for their roots.}
Moreover, by Lemma~\ref{lem:lift}, the eigenvalues of $M_{\leq n}$ form a subset of those of $\{\mathbf{M}_i\}_{1\leq i\leq n}$. Hence, we also obtain the  explicit closed-form formulas for the eigenvalues of $M_{\leq n}$.
For the five special cases in item (iv) of Theorem~\ref{thm:ODE_stab}, in the proof of Theorem \ref{thm:uns_specific} we construct an invertible map between certain sub-matrices of $\{\mathbf{M}_i\}_{1\leq i\leq n}$ and $M_{\leq n}$.  
Using these eigenvalues together with the definition~\eqref{def:Sigma_uns_n}, we obtain a
complete characterization of the linear subspace $ \Sigma_{\mathsf{uns}, \leq n}$, which determines the unstable + center manifold.
\end{remark}

\begin{remark}[\bf Relation between the ODEs, the bulk stability estimates, and LWP]
\label{rem:relation}
Theorem~\ref{thm:ODE_stab}, which provides the existence of a global, exponentially decaying solution to the closed ODE system~\eqref{eq:ODE_stab}, and Theorem~\ref{thm:nonrad:stab}, which proves global exponential decay of the bulk perturbation, conditional on assumption~\eqref{eq:ass:EOM}, interact through the local well-posedness (LWP) of the full PDE system~\eqref{eq:sys}. We summarize this interaction below; this yields a self-consistent global-in-$\tau$ stability theorem for both the perturbation bulk and the Taylor coefficients at the origin.

Given small initial data $V_{1,0} \in \Reals^{d_{\leq n}}$ for the ODE system~\eqref{eq:ODE_stab}, item~(iii) of Theorem~\ref{thm:ODE_stab} produces an exponentially decaying solution $V_{\mathsf{ODE},\leq n}(\tau) $ from 
initial data $V_{\mathsf{ODE},\leq n}(0) = V_{1,0} + (V_{\mathsf{uns},0},0,\ldots,0)$ defined for all $\tau \geq 0$. The construction of $V_{\mathsf{ODE},\leq n}$ and the decay estimate~\eqref{eq:ODE_decay} are valid \emph{regardless} of the lifespan of the solution $\tw = (\tvr,\tu,\tb)$ to the bulk PDE~\eqref{eq:lin_ODE}; that is, $V_{\mathsf{ODE},\leq n}$ is constructed without any input from the bulk analysis of Section~\ref{sec:nonradial}.

Now consider PDE initial data for~\eqref{eq:lin_ODE} whose Taylor coefficients at the origin agree with $V_{\mathsf{ODE},\leq n}(0)$, which satisfies the smallness conditions~\eqref{eq:nonrad_init:small} and the upper/lower bounds in~\eqref{ass:nonrad_init:a}--\eqref{ass:nonrad_init:b}.
The interaction between the ODE analysis (Theorem~\ref{thm:ODE_stab}), the bulk PDE estimates (Theorem~\ref{thm:nonrad:stab}), and the local well-posedness of~\eqref{eq:sys} proceeds via the following three implications, each of which is valid on any interval $[0,T]$ on which the high-regularity LWP of~\eqref{eq:sys} holds:
\begin{enumerate}[label=(\roman*),leftmargin=2em]
\item \textsl{(LWP $\Rightarrow$ ODE system).} Local well-posedness of~\eqref{eq:sys} in a sufficiently high-regularity class on $[0,T]$ (e.g.~$H^k$ with $k$ large) ensures that the Taylor expansion of $\tw=(\tvr,\tu,\tb)$ at the origin is well-defined for all $\tau \in [0,T]$. By the uniqueness of the ODE flow~\eqref{eq:ODE_stab}, the Taylor coefficients $V_{\mathsf{PDE}, \leq n}(\tau)$ of the PDE solution coincide on $[0,T]$ with $V_{\mathsf{ODE},\leq n}(\tau)$, and therefore inherit the decay estimate~\eqref{eq:ODE_decay}. In particular, assumption~\eqref{eq:ass:EOM} of Theorem~\ref{thm:nonrad:stab} holds on~$[0,T]$.

\item \textsl{(LWP $\Rightarrow$ bulk stability estimate).} On the same interval $[0,T]$, the validity of~\eqref{eq:ass:EOM} allows us to perform the nonlinear weighted energy estimates of Section~\ref{sec:nonradial}, yielding the bulk decay estimate~\eqref{eq:nonrad_decay} of Theorem~\ref{thm:nonrad:stab} for the perturbation $\tw_\mm$.

\item \textsl{(Stability estimates $\Rightarrow$ extension of LWP).} The ODE decay~\eqref{eq:ODE_decay} together with the bulk decay~\eqref{eq:nonrad_decay} ensure that $\tw=(\tvr,\tu,\tb)$ remains uniformly small in the high-regularity norm on $[0,T]$. The standard continuation criterion for the symmetric hyperbolic system~\eqref{eq:sys} then implies that the local solution may be extended beyond~$T$, to an interval $[0,T+\eps]$ for some $\eps>0$.
\end{enumerate}

A standard continuity argument in $\tau$ closes the cycle: the set of times $T \geq 0$ on which all of (i)--(iii) above hold is non-empty, open, and closed in $[0,\infty)$, hence equal to $[0,\infty)$. We conclude that the PDE solution $\tw$ exists globally in $\tau$, that its Taylor coefficients at the origin coincide with the prescribed ODE solution $V_{\mathsf{ODE},\leq n}$ for all $\tau \geq 0$, and that both the ODE decay estimate~\eqref{eq:ODE_decay} and the bulk decay estimate~\eqref{eq:nonrad_decay} hold.
\end{remark}

In the rest of this section, we first analyze the lower order ODE system $V_{\leq 2}$ 
in Sections \ref{sec:ODE_0th}, \ref{sec:ODE_N1}, and determine the modulation functions in Section \ref{sec:modulate_nosym}. In Section \ref{sec:ODE_high}, we present the key structural 
result Proposition \ref{prop:ODE_Nth} for the high order ODE system of $V_{\leq n}$, whose proof is defered to Appendix \ref{app:derive_ODE}. Using Proposition \ref{prop:ODE_Nth}, 
we analyze the number of unstable modes for the specific case and $\NNN=1$ in Section 
\ref{sec:unstable_modes_specific}, and for the general case in Section \ref{sec:unstable_modes_general}. In Section \ref{sec:global_ODE}, we construct the stable manifold of the ODE system \eqref{eq:ODE_stab} and prove the decay estimates in \eqref{eq:ODE_decay}. 
In Section \ref{sec:ODE_stab_proof}, we summarize the estimates in this section and prove 
Theorem \ref{thm:ODE_stab}.

\vspace{0.1in}
\paragraph{\bf Notation for error terms}
To track the lower order terms in the high order ODEs,  for any $l \geq 0$, we introduce the long vectors
\begin{subequations}\label{norm:ODE_low}
\begin{align}
  \FFs_{k, l} & :=  ( \na^k \Del^l \tvr(0) , \,  \na^k \Del^l (\div \tu)(0) , \,
      \na^k \Del^{l+1} ( y \cdot \tu)(0) , \,
  \na^k \Del^{l+1} \tb(0) ),  \label{norm:ODE_low:a} \\
   \GGs_k & := 
    \Big( \one_{k\geq 0} \na^{\leq k} \tvr(0) , \,
  \one_{k\geq -1} |\na^{\leq k+1} \tu(0),  \, \one_{k\geq -2} \na^{\leq k+2} \tb(0) ,
  \  \tcr , \tcu,  \tcvr, \tcb \Big) .
    \label{norm:ODE_low:b}
\end{align}

To track the \emph{linear} lower order terms, 
for a vector $g$, we introduce notation $ \OOL{h}( g )$ to denote a scalar or vector $f$, which depends linearly on $g$
\begin{equation}\label{norm:ODE_low:OOL}
f \in  \OOL{h}( g ) \ \Longrightarrow \  f = M_{h} \,  g
\end{equation}
\end{subequations}
for some constant matrix $M_{ h}$ that depends on the parameter $h$.  
We mainly apply the notation $ \OOL{}$ to the vectors $\FFs_{k, l}, \GGs_k$ defined above. Using the standard big-$\OO$ notation: $f = \OO_h(g)$, which denotes $|f| \leq C_h g$ for constant $C_h>0$, we obtain 
the relation 
\begin{equation}\label{eq:iden:OOL_O}
 \OOL{h} (\FFs_{k, l}) 
 = \OO_h( |\FFs_{k, l} | ),
 \quad  
 \OOL{h} (\GGs_{k}) 
 = \OO_h( |\GGs_{k} | ) 
\end{equation}
with $C_h = |M_h|$. 
Thus, \eqref{norm:ODE_low:OOL} generalizes the standard big-$\OO$ notation to additionally indicate linear structure.

\vspace{0.1in}
\paragraph{\bf Principal Taylor coefficients of the stationary self-similar profile}
Throughout this Section, we refer to the leading Taylor series coefficients $(\aaa,\aan,\uua,\uun,\bb,\bbn)$ of the radial profiles $(\bvr,\bu,\barb)$, see~\eqref{eq:ODE:main:profi}. For the convenience of the reader, we recall here that~\eqref{eq:ODE:main:profi} gives
\begin{equation}
\aaa  = (\alpha \bar q_0)^2 = \tfrac{\alpha \gamma d}{2}\uua^2 
\,, \qquad
\uua = \bar v_0 = -\tfrac{1}{1+\alpha d}\,, \qquad
\bb = 1 \,.
\label{eq:rho0:U1:B2}
\end{equation}
We also recall from~\eqref{eq:ODE:main:profi} and~\eqref{eq:Taylor:R=0:all} that
\begin{subequations}
\label{eq:rho2:U3:B4}
\begin{align}
\aan &= 2 \alpha^2 \bar q_0 (\bar q_{\NNN} - \bar q_0 \bar h_{\NNN})
= -  (\alpha \bar q_0)^2 \tfrac{\alpha(d+2\NNN)}{\NNN(\cxbar + \bar v_0)}
= -  \aaa \tfrac{\alpha(d+2\NNN)}{\NNN \kappa}, 
\\
\uun &=\bar v_{\NNN} = 1, 
\\
\bbn &= 2 \bar h_{\NNN} = - \tfrac{1}{\NNN(\bcr + \bar v_0)} = - \tfrac{1}{\NNN \kp},
\end{align}
\end{subequations}
where we recall from~\eqref{eq:cb:def} and~\eqref{eq:ODE_main:kp} that we have denoted 
\[
\kp = \cxbar + \bar v_0 = \bcr + \uua = \cbbar.
\]
We choose to work with $\kappa$ instead of the already-defined $\cbbar$, so as to not confuse with the $\tcb$ or $\bcb$ modulation parameters.

\subsection{The \texorpdfstring{$0$-th}{zeroth} order ODE system}\label{sec:ODE_0th} 

In this section, we derive and analyze the $0$-th order ODE , i.e., the evolution equations for the components of the vector $V_{\leq 0} = ( \tvr(0), \na^{\leq 1} \tu(0) , \na^{\leq 2} \tb(0) )$. In light of the vanishing condition \eqref{eq:vanishing:order:at:origin:SS}, we start by analyzing the behavior of $\tvr(0)$, $\na \tu(0)$, and $\na^2 \tb(0)$, which are the leading order terms for these perturbations at $y=0$. We introduce the strain matrix for $\td \UU$ and the Hessian of $\td B$, defined according to
\begin{align*}
 \Std := \tfrac 12 (\na \tu)(0,\cdot) + \tfrac 12 (\na \tu)^\intercal(0,\cdot) , 
\quad 
\Htd :=  (\na^2 \tb)(0,\cdot) .
\end{align*}
Here and throughout this section, we use the convention $(\nabla \tu)_{ij} := \partial_j \tu_i$.

The analysis of $(\tvr(0), \na \tu(0),\na^2 \tb(0))$ proceeds in three steps:
\begin{itemize}[leftmargin=2em]
\item we first analyze the behavior of the scalar fields (which have direct radially symmetric analogues)
\[
\tvr(0,\cdot) \,, \qquad \div \tu  (0,\cdot)\,, \qquad (\Delta \tb)(0,\cdot) 
\,,
\]
\item then we consider the evolution of the traceless part of the strain matrix for $\tu$ and the traceless part of the Hessian of $\tb$, namely the traceless symmetric matrices
\begin{align*}
\Sringbb &:= \tfrac 12 (\na \tu)(0,\cdot) + \tfrac 12 (\na \tu)^\intercal(0,\cdot) - \tfrac{1}{d} (\nabla \cdot \tu) (0,\cdot) {\rm Id} \,, \\
\Hringbb &:=  (\na^2 \tb)(0,\cdot) - \tfrac{1}{d}(\Delta \tb)(0,\cdot) {\rm Id} \,,
\end{align*}
\item we conclude with the analysis of the anti-symmetric part of $\na \tu$, which contains the components of vorticity, namely the traceless matrix
\[
\Aringbb := (\na \tu)(0,\cdot) - (\na \tu)^\intercal(0,\cdot)
\,.
\]
\end{itemize}

\subsubsection{ODEs for $\tvr, \div \tu $ and $\Delta \tb$}
From \eqref{eq:lin_ODE}-- \eqref{eq:lin_ODE:cL} (equivalently, from~\eqref{eq:lin:full}), and using \eqref{eq:vanishing:order:at:origin:SS}, we derive the ODE system for $(\tvr, \nabla \tu,\nabla^2 \mb)(0,\cdot)$
\begin{subequations}\label{eq:ODE_0th_full}
\begin{align}
\tf{d}{d \tau} \tvr  +     2 \al \div \tu \aaa - \tcvr   \aaa  & = \NNs_{\vrho} , \\ 
\tf{d}{d \tau} \na \tu 
+ \tfrac{2}{\gamma} \tvr \bb \Id  + \tfrac{1}{\gamma} \aaa \Htd  
+ (1 + 2 \uua) (\na \tu) + (\tcr - \tcu)  \uua \Id  & = \NNs_U , \label{eq:ODE_0th_full:b} \\
\tf{d}{d \tau} \Htd  + 4 \bb \Std
+    (2 \tcr - \tcb ) 2 \bb \Id  & = \NNs_{\mb}, \label{eq:ODE_0th_full:c}
\end{align}
\end{subequations}
where  $\NNs_{\vrho},  \NNs_U, \NNs_{\mb}$ denote the nonlinear terms
\begin{subequations}\label{eq:ODE_0th_full:non}
\begin{align}
  \NNs_{\vrho} &:= -  (   2 \al \div \tu - \tcvr ) \tvr,  \\
  \NNs_U  & := 
   - (\na \tu)^2 - (\tcr - \tcu) (\na \tu)- \tfrac{1}{\gamma} \tvr \Htd \\
  \NNs_{\mb} & :=  - \bigl(  (\na \tu)^{\intercal} \Htd + \Htd (\na \tu) \bigr)
  - ( 2 \tcr - \tcb ) \Htd. 
  \label{eq:ODE_0th_full:non:c}
\end{align}
\end{subequations}
Taking the trace in ~\eqref{eq:ODE_0th_full} and \eqref{eq:ODE_0th_full:non}, we obtain 
\begin{align}
&\tfrac{d}{d\tau} (\tvr,  \div \tu ,  \Delta \tb  )(0)^\intercal 
+
\mathsf{G}_0 (\tvr ,  \div \tu ,  \Delta \tb )(0)^\intercal 
\notag\\
&\qquad +
\bigl( (\tcb - 2 \tcu) \aaa, (\tcr-\tcu) d \uua, (2 \tcr - \tcb) 2 d \bb \bigr)^\intercal
= ( \NNs_{\vrho}, \tr( \NNs_U ), \tr( \NNs_{\mb} ) )^{\itl} 
\,,
\label{eq:mahomes:is:done:0}
\end{align}
where in the exact expression for the modulation term we have used that $\tcvr = 2 \tcu - \tcb$ 
derived from \eqref{eq:normal0} and \eqref{eq:profile_scal},  and we have denoted
\begin{align}
\mathsf{G}_0
:= 
\begin{pmatrix}
0 & 2 \alpha \aaa & 0 \\
\tfrac{2d}{\gamma} \bb & 1 + 2 \uua & \tfrac{1}{\gamma} \aaa \\
0 & 4 \bb & 0 
\end{pmatrix}
\,.
\end{align}

We may explicitly compute the eigen-system associated to the matrix $\mathsf{G}_0$ as:
\begin{itemize}[leftmargin=2em]
 \item Neutral eigenvalue $0$, with eigenvector $(\aaa, 0, -2 d \bb)^\intercal = ( (\alpha \bar q_0)^2 , 0 , - 2d)^\intercal$.
 \item Unstable eigenvalue $\frac 12 + \uua - \frac 12 \sqrt{(1+ 2 \uua)^2 + \frac{16}{\gamma} \aaa \bb (1+\alpha d)} = - 1$, with eigenvector \\ $(2\alpha \aaa, -1,  4 \bb)^\intercal = (2\alpha (\alpha \bar q_0)^2 , -1 , 4)^\intercal$.
 \item Stable eigenvalue $\frac 12 + \uua + \frac 12 \sqrt{(1+ 2 \uua)^2 + \frac{16}{\gamma} \aaa \bb (1+\alpha d)} = \frac{2\alpha d}{1+\alpha d}$, with explicit eigenvector $(2 \alpha  \aaa, \frac{2\alpha d}{1+\alpha d}, 4 \bb   )^\intercal = (2\alpha (\alpha \bar q_0)^2, \frac{2\alpha d}{1+\alpha d},   4 )^\intercal$.
\end{itemize}
As in the radially symmetric case, we use two (relations between the) modulation functions to remove the unstable direction corresponding to the eigenvalue $-1$ and the neutral direction corresponding to the eigenvalue $0$, for the matrix $\mathsf{G}_0$. To achieve this, we consider the matrix $\mathsf{P}$ whose columns are the above-mentioned eigenvectors, that is: 
\begin{equation}\label{def:P_sf}
\mathsf{P}:=
\begin{pmatrix}
2 \alpha \aaa &  \aaa & 2 \alpha \aaa \\
\frac{2\alpha d}{1+\alpha d} & 0 & -1 \\
4 \bb & -2 d \bb & 4 \bb
\end{pmatrix},
\qquad
\mathsf{P}^{-1}
=
\tfrac{1}{1+3\alpha d}
\begin{pmatrix}
\frac{d}{2\aaa}
& 1+\alpha d
& \frac{1}{4\bb}
\\ 
\frac{1+3\alpha d}{\aaa(1+\alpha d)}
& 0
& -\frac{\alpha (1+3\alpha d)}{2\bb(1+\alpha d)}
\\
\frac{\alpha d^{2}}{\aaa(1+\alpha d)}
& - (1+\alpha d)
& \frac{\alpha d}{2\bb(1+\alpha d)}
\end{pmatrix}
.
\end{equation}
Then, upon denoting 
\[
\mathscr{V} := \mathsf{P}^{-1} (\tvr,  \div \tu ,  \Delta \tb  )(0)^\intercal,
\]
we obtain from~\eqref{eq:mahomes:is:done:0} that
\begin{align*}
& \tfrac{d}{dt} \mathscr{V} + {\rm diag}(\tfrac{2\alpha d}{1+\alpha d}, 0 , -1) \mathscr{V} 
+ \mathsf{P}^{-1} \bigl( (\tcb - 2 \tcu) \aaa, (\tcr-\tcu) d \uua, (2 \tcr - \tcb) 2 d \bb \bigr)^\intercal
\\
&\qquad =  \mathsf{P}^{-1}  ( \NNs_{\vrho}, \tr( \NNs_U ), \tr( \NNs_{\mb} ) )^{\intercal} .
\end{align*}
Remarkably, the third term on the left side of the above equation does not have a component in the first entry; using that  $\bb = 1$ (cf.~\eqref{eq:rho0:U1:B2}), we may explicitly calculate
\[
\mathsf{P}^{-1} \bigl( (\tcb - 2 \tcu) \aaa, (\tcr-\tcu) d \uua, (2 \tcr - \tcb) 2 d \bb \bigr)^\intercal
=
\bigl(0, \tcb + 2 \uua (\tcu + \alpha d \tcr), -d \uua (\tcr - \tcu) \bigr)^\itl
.
\]
Therefore, upon letting 
\begin{subequations}
\label{eq:modulated:no:symmetry}
\begin{align}
\tcb(\tau) + 2 \uua \bigl(\tcu (\tau) + \alpha d \tcr(\tau) \bigr) 
&:= \tfrac{2\alpha d}{1+\alpha d} \mathscr{V}_2 (\tau)
\notag\\
&=\tfrac{2\alpha d}{(1+\alpha d)^2} \bigl(\tfrac{1}{\aaa} \tvr - \tfrac{\alpha}{2} \Delta \tb \bigr)(0,\tau)
\label{eq:modulated:no:symmetry:1}
\\
-d \uua \bigl(\tcr(\tau) - \tcu(\tau) \bigr) 
&:= \bigl( 1 + \tfrac{2\alpha d}{1+\alpha d} \bigr) \mathscr{V}_3 (\tau)
\notag\\
&= \tfrac{1 }{(1+\alpha d)^2}
\bigl( \tfrac{\alpha d^2}{\aaa} \tvr - (1+\alpha d)^2 \div\tu + \tfrac{\alpha d}{2} \Delta \tb \bigr)(0,\tau)  , 
\label{eq:modulated:no:symmetry:2}
\end{align} 
\end{subequations}
we obtain 
\[
\tfrac{d}{dt} \mathscr{V} +  \tfrac{2\alpha d}{1+\alpha d} \mathscr{V} 
=   \mathsf{P}^{-1}  ( \NNs_{\vrho}, \tr( \NNs_U ), \tr( \NNs_{\mb} ) )^{\itl}.
\]
In turn, multiplying the above identity from the left by  $\mathsf{P}$, we obtain 
\begin{align}
&\tfrac{d}{d\tau} (\tvr,  \div \tu ,  \Delta \tb  )(0)^\intercal 
+
\tfrac{2\alpha d}{1+\alpha d} (\tvr ,  \div \tu ,  \Delta \tb )(0)^\intercal 
= 
  ( \NNs_{\vrho}, \tr( \NNs_U ), \tr( \NNs_{\mb} ) )^{\itl} 
\,.
\label{eq:mahomes:is:done:1}
\end{align}
Thus, under the choice of modulation functions in~\eqref{eq:modulated:no:symmetry}, the system of ODEs~\eqref{eq:mahomes:is:done:0} becomes the fully stable ODE system~\eqref{eq:mahomes:is:done:1}.

\subsubsection{ODEs for $\Sringbb$, and $\Hringbb$.}

Taking the symmetric part of \eqref{eq:ODE_0th_full:b}, \eqref{eq:ODE_0th_full:c},  and subtracting the trace part, we derive the ODE system satisfied by the traceless symmetric matrices $\Sringbb$ and $\Hringbb$:
\begin{subequations}\label{eq:ODE_0th_SH}
\begin{align}
  \tf{d}{d \tau} \Sringbb + \tf{1}{\gamma} \aaa \Htd  
+ (1 + 2 \uua)  \Sringbb   & = \mr{\NNs_S} , \\
\tf{d}{d \tau} \Hringbb  + 4 \bb \Sringbb   & = \mr{\NNs_{\mb}},  
\end{align}
where the nonlinear terms $\NNs_{S}$ and the traceless nonlinear terms $\mr{\NNs_S}, \mr{\NNs_{\mb}}$ are given by 
 \begin{align} 
 \NNs_S  & = \tf{1}{2} (\NNs_U  + \NNs_U^{\itl} )
  = - (\Std^2 + \Aringbb^2) + (\tcu - \tcr) \Std - \tf{1}{\gamma} \tvr \Htd, \\
\mr{\NNs_S} &=  \NNs_S - \tf{1}{d} \tr{\NNs_S} \cdot \Id,   \\
\mr{\NNs_{\mb}} &= \NNs_{\mb} - \tf{1}{d} \tr \NNs_{\mb} \cdot \Id, 
 \end{align} 
\end{subequations}
and $\NNs_{\mb}$ is defined in \eqref{eq:ODE_0th_full:non:c}. To derive $\NNs_S$, we have used  $\na \tu = \Aringbb + \Std , \ (\na \tu)^{\itl} = \Std - \Aringbb$,  
and \eqref{eq:ODE_0th_full:non}.

Fix any $(i,j) \in \{1,\ldots,d\}^2$. Since the matrices $\Sringbb$ and $\Hringbb$ are traceless and symmetric, we are only interested in $\frac{d(d+1)}{2}-1$ many pairs $(i,j)$. 
 Then, using the notation $\GGs_0$ in \eqref{norm:ODE_low}, we have that 
 \begin{align}
&\tfrac{d}{d\tau} ( \Sringbb^{ij}, \Hringbb^{ij} )(0)^\intercal 
 +
 \mathsf{G}_0^{\sf sym} ( \Sringbb^{ij}, \Hringbb^{ij} )(0)^\intercal
 = (\mr{\NNs_S}^{ij}, \mr{\NNs_{\mb}}^{ij})^{\itl} = \OO( |\GGs_0|^2)
\,,
\label{eq:deviatoric:strain:evolution:tensor}
\end{align}
where we have denoted
\begin{align}
\mathsf{G}_0^{\sf sym}
:= 
\begin{pmatrix}
1 + 2 \uua  & \tfrac{1}{\gamma} \aaa \\
4 \bb & 0 
\end{pmatrix}
\,.
\end{align}
Note that here we do not have any modulation functions to help us.
Using the information $\aaa \bb = (\alpha \bar q_0)^2 = \frac{\alpha \gamma d}{2}   \uua^2$ and $\uua = -\frac{1}{1+\alpha d}$ from \eqref{eq:rho0:U1:B2}, we may explicitly compute the eigen-system associated to the matrix $\mathsf{G}_0^{\sf sym}$ as:
\begin{itemize} 
 \item Stable eigenvalue $\frac 12 + \uua + \frac 12 \sqrt{(1+ 2 \uua)^2 + \frac{16}{\gamma} \aaa \bb } = \frac{\alpha d -1 + \sqrt{1+6\alpha d +(\alpha d)^2}}{2(1+\alpha d)}> 0 $
 \item Unstable eigenvalue $\frac 12 + \uua - \frac 12 \sqrt{(1+ 2 \uua)^2 + \frac{16}{\gamma} \aaa \bb } = \frac{\alpha d -1 - \sqrt{1+6\alpha d +(\alpha d)^2}}{2(1+\alpha d)}< 0 $.
\end{itemize}
Note that the signs of the above-computed eigenvalues are  \emph{independent} of the value of $\alpha$.

\subsubsection{ODEs for $\Aringbb$.}
Fix any of the $\frac{d(d-1)}{2}$ independent pairs $(i,j) \in \{1,\ldots,d\}^2$ for which $\Aringbb^{ij}$ is not equal to $0$. Taking the anti-symmetric part of \eqref{eq:ODE_0th_full:b},
for each such pair $(i,j)$ we have
\begin{subequations}\label{eq:vorticity:evolution:tensor}
\begin{equation}
\tfrac{d}{d\tau} \Aringbb^{ij} + (1+ 2 \uua)  \Aringbb^{ij} 
= \NNs_A^{ij}, 
\end{equation}
where $\NNs_A$ denotes the nonlinear term
\begin{equation}
\NNs_{A}= 
 \bigl(\tcu - \tcr - \tfrac{2}{d} \div \tilde \UU(0) \bigr) \Aringbb
- \Sringbb \Aringbb
- \Aringbb \Sringbb 
= \OO( |\GGs_0|^2) .
\end{equation}
\end{subequations}

Here, $\Sringbb \Aringbb$ and $\Aringbb \Sringbb$ denotes matrix multiplication. Note that $1+ 2 \uua = \frac{\alpha d - 1}{1+\alpha d}$ changes sign as $\alpha$ crosses the ``monatomic gas'' value $\alpha = \frac 1d$.

\begin{remark}[\bf Stabilization of $0$-th order ODE via a higher co-dimension constraint]
From ~\eqref{eq:ODE_0th_SH} and~\eqref{eq:vorticity:evolution:tensor} it is clear that if at $\tau=0$ we specify initial data such that $\Sringbb(0) = \Hringbb(0)= \Aringbb(0) = 0$, which are $d(d+1) -2 + \frac{d(d-1)}{2}$ many constraints, then for all $\tau\geq 0$ we have  $\Sringbb(\tau) = \Hringbb(\tau)= \Aringbb(\tau) = 0$. 
\end{remark}

\begin{remark}[\bf Vorticity for the monatomic case $\alpha = 1/d$]
The coefficient $\tcu - \tcr - \tfrac{2}{d} \div \tu(0) $ appearing in~\eqref{eq:vorticity:evolution:tensor} may be rewritten using~\eqref{eq:modulated:no:symmetry:2} as 
\[
\tcu - \tcr - \tfrac{2}{d} \div \tu(0) 
=
- \tfrac{2(1+\alpha d)}{\gamma} \tvr(0) 
+ \tfrac{\alpha d-1}{d} \div \tu(0) 
- \tfrac{\alpha}{2(1+\alpha d)}   (\Delta \tb)(0)
.
\]
As such, when $\alpha = \frac 1d$, the evolution equation for the vorticity in~\eqref{eq:vorticity:evolution:tensor} becomes
\[
\tfrac{d}{d\tau} \Aringbb^{ij} 
=
- \bigl(\tfrac{2(1+\alpha d)}{\gamma} \tvr(0) 
+ \tfrac{\alpha}{2(1+\alpha d)}   (\Delta \tb)(0)  \bigr) \Aringbb^{ij}
- (\Sringbb \Aringbb)^{ij}
- (\Aringbb \Sringbb)^{ij}.
\]
While this is not pursued in the present paper, it would be very interesting to further analyze the (neutral) stability of the above ODE system for the vorticity (the components of $\Aringbb$), in the monatomic case $\alpha = 1/d$.
\end{remark}

\subsection{The ODE system for \texorpdfstring{$\Delta^{\NNN} \tvr, \Del^{\NNN+1} \xu, \Del^{\NNN+1} \tb $}{Laplacian to a power applied to rho, Z, B}}\label{sec:ODE_N1}
We introduce the scalar field $\xu$, defined as
\begin{equation}\label{eq:XU_V}
  \xu(y,\tau)  := y \cdot \tu(y,\tau).
\end{equation}
This is the analog of the radial velocity $R^2 \td V(R,\tau)$, which played an important role in the stability analysis of Section~\ref{sec:stability}. 

The evolution equations for the scalar quantities $\Delta^{\NNN} \tvr$, $\Del^{\NNN+1} \xu$, and $\Del^{\NNN+1} \tb$ are coupled (see~\eqref{eq:mahomes:is:done:del0} below); they are obtained as a particular case of the evolution equations derived in Proposition~\ref{prop:ODE_Nth},
and \eqref{eq:ODE_Nth:spec:c}--\eqref{eq:ODE_Nth:spec:d} below. To avoid redundancy, we do not repeat here the derivation of this system; we simply note that \eqref{eq:ODE_Nth:spec:c}--\eqref{eq:ODE_Nth:spec:f} yield 
\begin{align}
&\tfrac{d}{d\tau} (\Delta^{\NNN} \tvr,  \Delta^{\NNN+1} \xu ,  \Delta^{\NNN+1} \tb  )(0)^\intercal 
+ \mathsf{G}_1 (\Delta^{\NNN} \tvr,  \Delta^{\NNN+1} \xu ,  \Delta^{\NNN+1} \tb  )(0)^\intercal
+ \mathsf{G}_2 (\tvr,  \div \tu ,  \Delta \tb  )(0)^\intercal 
\notag\\
&= 
\bigl( (2 \tcu -\tcb  - 2\NNN \tcr) \ccs_{\NNN, 1}   \aan , (\tcu - ( 2\NNN+1) \tcr) \ccs_{\NNN+1, 1} \uun  , (\tcb- (2 \NNN+2) \tcr) \ccs_{\NNN+1, 1} \bbn \bigr)^\intercal
\notag\\
& \qquad 
+ \OO_{\NNN} (  |\GGs_{2\NNN}|^2)
\,,
\label{eq:mahomes:is:done:del0}
\end{align}
where using the convention $\kappa =  \bcr + \uua $ (see~\eqref{eq:ODE_main:kp}), and we have denoted
\begin{equation} 
  \mathsf{G}_1 :=  
  - \HHH_{\NNN} = 
  \begin{pmatrix}
    2\NNN \kp    &  \frac{\al}{\NNN+1} \aaa  & 0  \\ 
    2 ( \frac{\NNN}{\al} + \frac{2}{\gamma}) 
 (d  + 2 \NNN ) (\NNN+1) \bb  &  2\NNN \kp + 1 + 2\uua    &      \frac{2 (\NNN+1) }{\gamma} \aaa   \\
   0  &  2 \bb &  2 \NNN \kp    
   \end{pmatrix} , 
\end{equation}
we have used the formulas \eqref{eq:ODE_Nth:spec:g}, namely 
\[
\msf{c}_{ \NNN , \Del} := ( \NNN +1) \prod_{1 \leq i \leq \NNN} 2i (2 i + d),   \quad 
\msf c_{ \NNN, 1} = \prod_{1\leq i \leq \NNN} 2i (2i + d - 2),
\]
and we have denoted the coupling matrix to the terms with lower order derivatives by
\begin{align}
\mathsf{G}_2
:= 
\begin{pmatrix}
2 \al  (2 \NNN + d ) \msf \msf c_{\NNN, 1} \uun 
&  2 (\al \msf c_{\NNN, 1} +  2 \NNN \msf c_{\NNN-1, \Del} ) \aan
& 0 \\
\tfrac{1}{\gamma} (2\NNN+2 )  \ccs_{\NNN+1, 1}  \bbn 
& 4\NNN  \ccs_{\NNN, \Del} \uun 
&  ( \tfrac{2}{\gamma} + \tfrac{ \NNN }{\al} ) \ccs_{\NNN,\Del} \aan  \\
0 
& (4 \NNN+4)\ccs_{\NNN, \Del}   \bbn 
& - 2 \ccs_{\NNN, \Del}  \uun
\end{pmatrix},
\label{eq:gambling:G2}
\end{align}
where $(\aaa,\uua,\bb)$ are given by~\eqref{eq:rho0:U1:B2}, and $(\aan,\uun,\bbn)$ are given by~\eqref{eq:rho2:U3:B4}. 

We recall from~\eqref{eq:cx:admissible:1} that $\kp = \cxbar + \uua$ is the larger (of the two roots) of the quadratic equation 
\[
\kp^2 
+\kp \tfrac{1 + 2 \uua}{2 \NNN}  
- \uua^2 \bigl( \tfrac{\alpha \gamma d}{2} + \mathsf{E}_\NNN \bigr) 
= 0,
\]
with $\mathsf{E}_\NNN$ as given by \eqref{eq:En:def}.
Using this fact and~\eqref{eq:rho2:U3:B4}, we may explicitly compute the eigen-system associated to the matrix $\mathsf{G}_1$ as:
\begin{itemize} 
 \item Neutral eigenvalue $0$, with eigenvector $(-\alpha  \aaa, 2 \NNN (\NNN+1) \kappa , - 2 (\NNN+1)\bb   )^\intercal$.
 \item Stable eigenvalue $2 \NNN \kp >0 $, with eigenvector $(\alpha \aaa , 0 , -(d+2\NNN) (\NNN \gamma + 2 \alpha) \bb )^\intercal$.
 \item Stable eigenvalue $4\NNN \kp + 1+ 2 \uua >0$, with   eigenvector $(\alpha  \aaa , (\NNN+1)(2 \NNN \kappa + 1 + \uua), 2 (\NNN+1) \bb  )^\intercal $.
\end{itemize}
As in the spherically symmetric case, the final constraint on the modulation functions is used to modulate the neutral eigenvalue of $\mathsf{G}_1$. For this purpose, we need to rewrite the first term on ${\sf RHS}_{\eqref{eq:mahomes:is:done:del0}}$, taking into account the definitions in~\eqref{eq:modulated:no:symmetry}. Indeed, using~\eqref{eq:modulated:no:symmetry}, after a tedious computation we arrive at
\begin{align}
& \bigl( (2 \tcu -\tcb  - 2\NNN \tcr) \ccs_{\NNN, 1}   \aan , (\tcu - ( 2\NNN+1) \tcr) \ccs_{\NNN+1, 1} \uun  , (\tcb- (2 \NNN+2) \tcr) \ccs_{\NNN+1, 1} \bbn \bigr)^\intercal    
\notag \\
&= 
-2\NNN\,\tcr 
\bigl(
\ccs_{\NNN,1}\,\aan, 
\ccs_{\NNN+1,1}\,\uun, 
\ccs_{\NNN+1,1}\,\bbn
\bigr)^\intercal
+ \mathsf{G}_3
(\tvr,  \div \tu ,  \Delta \tb  )(0)^\intercal 
\label{eq:mahomes:is:done:2}
\end{align}
where
\begin{align}
\mathsf{G}_3 :=
\begin{pmatrix}
-\frac{2\alpha d}{(1+\alpha d)\aaa}\,\ccs_{\NNN,1}\aan 
& 2\alpha\,\ccs_{\NNN,1}\aan 
& 0 
\\ 
-\frac{\alpha d}{(1+\alpha d)\aaa}\,\ccs_{\NNN+1,1}\uun 
& \frac{1+\alpha d}{d}\,\ccs_{\NNN+1,1}\uun 
& -\frac{\alpha}{2(1+\alpha d)}\,\ccs_{\NNN+1,1}\uun 
\\ 
0 
& \frac{2}{d}\,\ccs_{\NNN+1,1}\bbn 
& -\frac{\alpha}{1+\alpha d}\,\ccs_{\NNN+1,1}\bbn
\end{pmatrix}
\label{eq:gambling:G3}
\end{align}
Combining~\eqref{eq:mahomes:is:done:del0} with~\eqref{eq:mahomes:is:done:2}--\eqref{eq:gambling:G3} we thus arrive at
\begin{align}
&\tfrac{d}{d\tau} (\Delta^{\NNN} \tvr,  \Delta^{\NNN+1} \xu ,  \Delta^{\NNN+1} \tb  )(0)^\intercal 
+ \mathsf{G}_1 (\Delta^{\NNN} \tvr,  \Delta^{\NNN+1} \xu ,  \Delta^{\NNN+1} \tb  )(0)^\intercal
\notag\\
&\qquad 
+ 2\NNN\,\tcr 
\bigl(
\ccs_{\NNN,1}\,\aan, 
\ccs_{\NNN+1,1}\,\uun, 
\ccs_{\NNN+1,1}\,\bbn
\bigr)^\intercal
\notag\\
&= ( \mathsf{G}_3 - \mathsf{G}_2) (\tvr,  \div \tu ,  \Delta \tb  )(0)^\intercal 
+ \OO_{\NNN}( |\GGs_{2\NNN}|^2)
\,.
\label{eq:mahomes:is:done:3}
\end{align}
As before, we next introduce the matrix whose columns are the eigenvectors of the matrix $\mathsf{G}_1$; that is, we define
\begin{equation*}
\mathsf{P} 
:= 
\begin{pmatrix}
\alpha \aaa 
& \alpha \aaa 
& - \alpha \aaa
\\
(\NNN+1)(2 \NNN \kp + 1 + 2\uua) 
& 0 
& 2\NNN (\NNN+1) \kp \\
2(\NNN+1) \bb 
& - (d+2\NNN)(\NNN\gamma+2\alpha) \bb  
& - 2 (\NNN+1)\bb 
\end{pmatrix},
\end{equation*}
and we define 
\begin{equation}\label{def:mathscr_W}
 \mathscr{W}:= \mathsf{P}^{-1}(\Delta^{\NNN} \tvr,  \Delta^{\NNN+1} \xu ,  \Delta^{\NNN+1} \tb  )(0)^\intercal.  
\end{equation}
Then, \eqref{eq:mahomes:is:done:3} may be recast as
\begin{align}
&\tfrac{d}{dt} \mathscr{W} 
+ {\rm diag}(4 \NNN \kp + 1 + 2 \uua, 2 \NNN \kp , 0) \mathscr{W} 
+ 2\NNN \tcr \mathsf{P}^{-1} \bigl(
\ccs_{\NNN,1}\,\aan, 
\ccs_{\NNN+1,1}\,\uun, 
\ccs_{\NNN+1,1}\,\bbn
\bigr)^\intercal
\notag\\
&
= \mathsf{P}^{-1} \bigl(\mathsf{G}_3 - \mathsf{G}_2\bigr) (\tvr,  \div \tu ,  \Delta \tb  )(0)^\intercal
+ \OO_{\NNN}(| \GGs_{2\NNN}|^2).
\label{eq:mahomes:is:done:4}
\end{align}
Using~\eqref{eq:rho0:U1:B2}, \eqref{eq:rho2:U3:B4}, the definition of $\mathsf{P}$ above, and definition~\eqref{eq:ODE_Nth:spec:g} (which gives $2(\NNN+1) (2\NNN+d) \ccs_{\NNN,1}= \ccs_{\NNN+1,1}$), we may verify the remarkable identity 
\begin{align}
\mathsf{P}^{-1} \bigl(
\ccs_{\NNN,1}\,\aan, 
\ccs_{\NNN+1,1}\,\uun, 
\ccs_{\NNN+1,1}\,\bbn
\bigr)^\intercal
&=- \tfrac{1}{\NNN \kp}
\mathsf{P}^{-1} \bigl(
  \alpha\ccs_{\NNN,1} (d+2\NNN) \aaa , 
- \ccs_{\NNN+1,1}  \NNN \kp , 
\ccs_{\NNN+1,1} \bb
\bigr)^\intercal
\notag\\
&= 
\tfrac{\ccs_{\NNN+1,1}}{2\NNN (\NNN+1) \kp}
\mathsf{P}^{-1} \mathsf{P} 
\vec{e}_3
= \bigl(0, 0, \tfrac{\ccs_{\NNN+1,1}}{2\NNN (\NNN+1) \kp} \bigr)^\intercal,
\label{eq:modulated:no:symmetry:3:semi:new}
\end{align}
which allows us the drastically simplify the third term on ${\sf LHS}_{\eqref{def:mathscr_W}}$.
With~\eqref{eq:modulated:no:symmetry:3:semi:new} and~\eqref{def:mathscr_W}, it is thus natural to define the modulation function $\tcr$ as 
\begin{align}
\tcr(\tau) 
:= 
\tfrac{2 \NNN (\NNN+1) }\kp^2{\ccs_{\NNN+1,1}} \mathscr{W}_3(\tau).
\label{eq:modulated:no:symmetry:3}
\end{align}
With the choice in~\eqref{eq:modulated:no:symmetry:3}, 
the evolution~\eqref{eq:mahomes:is:done:4}
may be recast as
\begin{align}
&\tfrac{d}{dt} \mathscr{W} + {\rm diag} \bigl(4 \NNN \kp + 1 + 2 \uua, 2 \NNN \kp , 2 \NNN \kp \bigr) \mathscr{W} 
\notag\\
&\qquad 
= \mathsf{P}^{-1} \bigl(\mathsf{G}_3 - \mathsf{G}_2\bigr) (\tvr,  \div\tu ,  \Delta \tb  )(0)^\intercal
+ \OO_{\NNN}( |\GGs_{2\NNN}|^2).
\label{eq:mahomes:is:done:5}
\end{align}
In particular, \eqref{eq:mahomes:is:done:5} implies that the full vector $\mathscr{W}$ inherits the temporal decay which is the worst between  $(\tvr,  \div \tu ,  \Delta \tb  )(0)$ and $|\GGs_{2\NNN}|^2$. The decay in time of $\mathscr{W}(\tau)$ then directly implies the decay in time of $\mathsf{P}  \mathscr{W}(\tau) =  (\Delta^{\NNN} \tvr,  \Delta^{\NNN+1} \xu ,  \Delta^{\NNN+1} \tb  )(0,\tau)^\intercal$.

\subsection{Definition of modulation functions}\label{sec:modulate_nosym}

Recall from \eqref{eq:rho0:U1:B2} that $\uua = - \frac{1}{1 + \al d }$, recall that $\kp = \cxbar + \uua$, and recall from \eqref{eq:normal0} the constraint: 
\begin{equation}\label{eq:modulated:tcvr}
\tcvr(\tau) + \tcb(\tau) = 2 \tcu(\tau) ,\quad \forall  \tau \geq 0.
\end{equation}
Combining the formulas for the modulation functions \eqref{eq:modulated:no:symmetry:1}, \eqref{eq:modulated:no:symmetry:2}, \eqref{eq:modulated:no:symmetry:3}, 
we determine the modulation functions as follows.
First, solve~\eqref{eq:mahomes:is:done:0} for the evolution of the vector $(\tvr,\div\tu,\Delta\tb)(0,\tau)$. Second, solve \eqref{eq:mahomes:is:done:5} for the evolution of the vector $\mathscr{W}(\tau)$. Finally, define: \begin{subequations}\label{eq:modulated:nosym:summary}
\begin{align}
\tcr(\tau)  & = \tfrac{2\NNN(\NNN+1)}{\ccs_{\NNN+1,1}} \kp^2 \mathscr{W}_3(\tau), \\
\tcu (\tau) & = \tcr(\tau) - 
 \tfrac{1 }{(1+\alpha d) d } \bigl( \tfrac{\alpha d^2}{\aaa} \tvr(0,\tau) - (1+\alpha d)^2 \div \tu(0,\tau) + \tfrac{\alpha d}{2} (\Delta \tb)(0,\tau) \bigr) , \\
 \tcb(\tau) & = \tfrac{2}{1 + \al d}  \bigl(\tcu (\tau) + \alpha d \tcr(\tau) \bigr) + \tfrac{2\alpha d}{(1+\alpha d)^2} \bigl(\tfrac{1}{\aaa} \tvr(0,\tau) - \tfrac{\alpha}{2} (\Delta \tb)(0,\tau) \bigr) ,  \\
\tcvr(\tau) & = 2 \tcu(\tau) - \tcb(\tau). 
\end{align}
\end{subequations}
In the above definitions we do not substitute the formula \eqref{def:mathscr_W} for $\mathscr{W}_3$, nor do we expand the formula for $\tcb$,  since these \emph{explicit} formulas are not used in later derivations; the formulas in~\eqref{eq:modulated:nosym:summary} are only used to show that the modulation functions decay exponentially fast in $\tau$, under the standing bootstrap assumptions.

\subsection{The fundamental high order ODE systems}\label{sec:ODE_high}
In order to derive the evolution equations for $V_{=n}(\tau) =(\na^{n} \tvr,\na^{n+1} \tu, \na^{n+2} \tb)(0,\tau)$, we apply $\na^n$ to \eqref{eq:lin:full:a}, $\na^{n+1}$ to \eqref{eq:lin:full:b}, $\na^{n+2}$ to \eqref{eq:lin:full:c}, and then restrict the resulting PDE to $y=0$. When $n\geq 1$ is general, the resulting system of ODEs is quite complicated, and the ``accounting'' becomes delicate. In order to deal with this ``accounting problem'', for an integer $k \geq 0$ and a multi-index $\beta \in \Naturals_0^d$, we introduce the mixed derivative vector 
\begin{equation}\label{def:VVs}
 \VVs_{\bet, k}(\tau)  := (\pa^{\bet} \Del^{k} \tvr(0,\tau), \  \pa^{\bet} \Del^{k} \div \tu (0,\tau),  
 \ \pa^{\bet} \Del^{k+1} \xu (0,\tau),  \ \pa^{\bet} \Del^{k+1} \tb(0,\tau)   )^\intercal \in \Reals^{4 \times 1}. 
\end{equation}
The reason for singling out powers of the Laplacian $k$, as opposed to just keeping track of all multi-indices $\beta$ at once, is as follows. Heuristically, since the profile $\brho, \bu, \bar \mb$ is radially symmetric, 
by taking a more symmetric derivative combination
of the perturbation $(\na^{k} \tvr(0,\tau), \na^{k+1} \tu(0,\tau)$, $\na^{k+2} \tb(0,\tau) ) $, 
one expects that the ODE system associated with this combination is simpler; one of the most symmetric combinations of mixed derivatives is the Laplacian. 
Guided by the above heuristic,  we first analyze the ODE system for $ \VVs_{\bet, k}$ with high powers of $\Delta$ (larger $k$), which has a simpler form, and then we analyze other mixed derivatives.

Below, we fix the parameter $\NNN$ for the profile and treat it as an absolute constant.  
We recall the notations $\FFs ,\GGs, \OOL{}$ from \eqref{norm:ODE_low}. We have the following structural results for the ODE system. 

\begin{proposition}[\bf Structure of the ODE system]
\label{prop:ODE_Nth}
Let $\aaa, \uua, \bb, \aan, \uun,\bbn$  be as in \eqref{eq:rho0:U1:B2}--\eqref{eq:rho2:U3:B4}, and let $\kp = \cxbar+\uua$ be as defined in~\eqref{eq:ODE_main:kp}. Denote $a_+ = \max(a, 0)$. 
Suppose that the system \eqref{eq:sys} admits a solution with regularity $C^{\mm + 1}( B_1(0))$ 
on a given time interval; then the following conclusions hold on the same time interval.

Let $n,k \in\Naturals_0$ be such that $ 1 \leq n \leq \mm$ and $ 0\leq k \leq \frac{n}{2} $.  Let $\bet$ be a multi-index with $|\bet| = n - 2k$. Then:

\begin{itemize}[leftmargin=2em]

\item[(i)] {\bf ODEs for powers of the Laplacian.}
Let $n$ be even, $|\bet| =0$, and $k = \frac{n}{2}$; we have 
\begin{subequations}\label{eq:ODE_Nth:spec}
\begin{equation}\label{eq:ODE_Nth:spec:a}
 \tf{d}{d\tau} 
(\Del^k \tvr, \Del^{k+1} \xu, \Del^{k+1} \tb)(0)^\intercal
= \HHH_{ k} (\Del^k \tvr, \Del^{k+1} \xu, \Del^{k+1} \tb)(0)^\intercal
  + \OOL{n} (\GGs_{ (n -2\NNN )_+} )  + \OO_n( |\GGs_{ n  }|^2) .
\end{equation}
where the matrix $\HHH_{ k}$ associated with the top order linear terms is given by
\begin{equation}\label{eq:ODE_Nth:spec:b}
\begin{aligned} 
  \HHH_{ k} =  
  \begin{pmatrix}
   - n \kp    & - \frac{ 2 \al}{n+2} \aaa  & 0  \\ 
   - ( \frac{n}{2 \al} + \frac{2}{\gamma}) 
 (2 d  + 4k ) (k+1) \bb  &   - n \kp  - 2\uua -1    &     - \frac{ (n+2) }{\gamma} \aaa   \\
   0  &    - 2 \bb &  - n \kp    
   \end{pmatrix} , \quad  n = 2k.
   \end{aligned}
\end{equation}
In the special case $|\bet| =0, n  = 2\NNN, k =\NNN $, we have 
\begin{equation}\label{eq:ODE_Nth:spec:c}
 \frac{d}{d\tau} 
 \begin{pmatrix}
\Del^{\NNN} \tvr \\ 
\Del^{\NNN+1} \xu \\ 
\Del^{\NNN+1} \tb 
\end{pmatrix}
= \HHH_{\NNN}   \begin{pmatrix}
\Del^{\NNN} \tvr \\ 
\Del^{\NNN+1} \xu \\ 
\Del^{\NNN+1} \tb 
\end{pmatrix} 
+ 
\begin{pmatrix}
\LLs_1\\ 
\LLs_2  \\ 
\LLs_3 \\
\end{pmatrix} 
  + \OO_{\NNN} (  |\GGs_{2\NNN}|^2) ,
\end{equation}
where vector $(\LLs_1,\LLs_2,\LLs_3)$ depends on the modulation functions and on the vector $(\tvr,\div\tu,\Delta\bu)(0,\tau)$, and is defined as
\begin{align}
\LLs_1  
&= ( 2 \tcu -\tcb  - 2\NNN \tcr ) \ccs_{\NNN, 1}   \aan 
\notag\\
&\qquad 
- 2 \al  (2 \NNN + d )\ccs_{\NNN, 1} \uun   \tvr 
- ( 2\al \msf c_{\NNN, 1} +  4 \NNN \msf c_{\NNN-1, \Del} ) \aan \, \div \tu
 , 
\label{eq:ODE_Nth:spec:d}
\\
\LLs_2 
&= ( \tcu - ( 2\NNN+1) \tcr  ) \ccs_{\NNN+1, 1} \uun 
\notag\\
&\qquad 
- \tfrac{1}{\gamma} (2\NNN+2 )  \ccs_{\NNN+1, 1}  \bbn \tvr 
-  4\NNN  \ccs_{\NNN, \Del} \uun \, \div \tu
- ( \tfrac{2}{\gamma} + \tfrac{ \NNN }{\al} )\ccs_{\NNN,\Del}  \aan \Del \tb 
, 
\label{eq:ODE_Nth:spec:e}\\
\LLs_3 
&= (\tcb- (2 \NNN+2) \tcr)\ccs_{\NNN+1, 1} \bbn  
\notag\\
&\qquad 
- (4 \NNN+4) \ccs_{\NNN, \Del}  \bbn 
+ 2 \ccs_{\NNN, \Del}  \uun   \Del \tb \, \div \tu ,
\label{eq:ODE_Nth:spec:f}
\end{align}
where $\ccs_{\NNN,\Del}, \ccs_{\NNN,1}$ are defined as 
\begin{align} 
\label{eq:ODE_Nth:spec:g}
\msf{c}_{ \NNN , \Del} := ( \NNN +1) \prod_{1 \leq i \leq \NNN} 2i (2 i + d),   \quad 
\msf c_{ \NNN, 1} = \prod_{1\leq i \leq \NNN} 2i (2i + d - 2).
\end{align}

\end{subequations}

\item[(ii)] {\bf ODEs for mixed derivatives.}
\begin{subequations}\label{eq:ODE_Nth}
Let $0\leq k \leq \frac{n}{2}$ and $|\bet| + 2k = n$; we have
\begin{equation}
\label{eq:ODE_Nth:a}
 \tf{d}{d \tau}  \VVs_{\bet, k} := \HH_{|\bet|, k} \VVs_{\bet, k} +  
  \one_{ |\bet| \geq 2 } \OOL{ \bet, k}(\FFs_{ |\bet|-2, k+1 }) 
  + \OOL{\bet, k} ( \GGs_{ ( n - 2 \NNN )_+} )  + \OO_n( | \GGs_{ n  }|^2) ,
\end{equation}
where the matrix $\HH_{|\bet|, k}$ associated with the top order linear terms is given by
\begin{align}
  & \HH_{|\bet|, k} :=  
  \label{eq:ODE_Nth:b}\\
  &\begin{pmatrix}
   - n \kp  & - 2 \al \aaa &  0 & 0  \\ 
   -   \frac{ ( d + 2 |\bet| + 2 (k-1) ) k    \bb }{\al}  -  ( \frac{n}{\al} + \frac{2( d+ n)}{\gamma} ) \bb 
     & - n \kp - 2 \uua -1  & 0  &  - \frac{1}{\gamma} \aaa  \\
   -  ( \frac{n}{2 \al} + \frac{2}{\gamma}) 
 (2 d + 4 |\bet| + 4k ) (k+1) \bb  &  0 & - n \kp - 2\uua -1   &     - \frac{ (n+2) }{\gamma} \aaa   \\
   0  &  0 &  -2 \bb &  - n \kp   \notag 
  \end{pmatrix} .
  \end{align}
\end{subequations}

\item[(iii)] {\bf ODEs for derivatives of $\xu$ and $\tb$.} After estimating the ODEs~\eqref{eq:ODE_Nth:a} for all multi-indices $(\bet, k)$ with $|\bet | + 2k = n$, we may treat $\na^n \tvr(0)$ as a given lower order term. 

Then, for any multi-index $\bet$ with $|\bet| = n+2$, we have
 \begin{subequations}\label{eq:ODE_Nth_ZB}
 \begin{align} 
    \tfrac{d}{d\tau} \begin{pmatrix}
   \pa^{\bet} \xu \\
    \pa^{\bet} \tb  \\
    \end{pmatrix}
    &= \HHT_{n}  \begin{pmatrix}
   \pa^{\bet} \xu  \\
    \pa^{\bet} \tb   \\
    \end{pmatrix} 
    +  \OOL{\beta}\bigl(\na^n \tvr(0)\bigr)
    + \OOL{n}(\GGs_{ (n- 2\NNN )_+ }) 
    + \OO_{n} ( |\GGs_n|^2  ) ,  
  \label{eq:ODE_Nth_ZB:a}  \\
 \HHT_{n} & = \begin{pmatrix}
  - n \kp - 2 \uua -1    &  - \frac{n+2}{\gamma} \aaa  \\
     - 2\bb   & - n \kp \\  
    \end{pmatrix}, \quad  |\bet| = n+2.
  \end{align}
 \end{subequations}

\item[(iv)] {\bf ODEs for $\tu$.}
After estimating the ODEs for all mixed derivatives $\na^n \tvr, \na^{n+2} \tb$, cf.~\eqref{eq:ODE_Nth:a} and~\eqref{eq:ODE_Nth_ZB:a}, we treat these terms as given lower order terms. Then, for any multi-index $\bet$ with $|\bet| = n+1$, we have 
\begin{equation}\label{eq:ODE_Nth_U}
  \tfrac{d}{d\tau} \pa^{\bet} \tu = 
    (- n \kp - 2 \uua -1 ) \pa^{\bet} \tu
    + \OOL{\beta}\bigl( \na^n \tvr(0) , \na^{n+2} \tb(0) \bigr)
    + \OOL{n}(\GGs_{ ( n- 2 \NNN )_+}) 
    + \OO_{n}(|\GGs_n|^2).
\end{equation}

\item[(v)] {\bf ODEs in 1D.}
In the special case of $d = 1$, we do not have mixed-derivatives and~\eqref{eq:ODE_Nth}
simplifies to
\begin{subequations}\label{eq:ODE_Nth:1D}
\begin{equation}
\tf{d}{d\tau} (\pa^{n} \tvr, \pa^{n+1} \tu, \pa^{n+2} \tb)^\intercal
= \HH^{(1)}_n (\pa^{n} \tvr, \pa^{n+1} \tu, \pa^{n+2} \tb)^\intercal
+ \OOL{n} ( \GGs_{(n - 2\NNN)_+} )  
+ \OO_{n}( |\GGs_{n}|^2) ,
\end{equation}
where
\begin{equation}
\HH^{(1)}_n = 
\begin{pmatrix}
- n \kp  
& - 2 \al \aaa  
& 0 \\
- (  \frac{  n (n+1) }{2 \al} + \frac{2 (n+1)}{\gamma} ) \bb  
& - n  \kp - 2\uua -1  
& - \frac{1}{\gamma} \aaa \\
0 
& - 2(n+2)\bb   
& - n \kp  
\end{pmatrix} 
.
\end{equation}
\end{subequations}

\end{itemize}
\end{proposition}

The proof of Proposition~\ref{prop:ODE_Nth} follows from a careful Taylor expansion near $y=0$.
Since the argument is standard but technically involved and lengthy, we defer the proof to Appendix~\ref{app:derive_ODE}.

\begin{remark}[\bf ODEs for $\div \tu$ and  $y \cdot \tu$]
When $\bet = 0$ and $ k = n / 2$, the ODEs \eqref{eq:ODE_Nth} for the quantities $\Del^k ( \div \tu )(0)$ 
and $\Del^{k+1}( y \cdot \tu )(0)$ are the same since 
\begin{equation}
\label{eq:iden_Delta_xu_divu}
 \Delta^{k+1} (y \cdot \tu)(0) 
 = 
 2(k+1) \Del^k (\div \tu)(0), \quad \forall k \geq 0,
\end{equation}
due to the identity \eqref{eq:lap_id:div}. 
Note that in the non-radial case, and for $\bet \neq 0$, in general $\partial^{\beta} \Del^k ( \div \tu )(0)$ is not a constant multiple of  $\partial^{\beta} \Del^{k+1}( y \cdot \tu )(0)$ with a constant independent of $\tu$.
\end{remark}

\vspace{0.1in}
\paragraph{\bf Schematic form of the ODEs for $V_{=n}$ and $V_{\leq n}$}
Recall the notation for the $n$-th order Taylor coefficients from \eqref{def:TL_nth_order}, namely the vector $V_{=n} = (\na^{n} \tvr(0),  \na^{n+1} \tu(0), \na^{n+2} \tb(0) )$.  Since the ODE system satisfied by $ (\Del^{\NNN} \tvr ,  \Del^{\NNN+1} \xu , \Del^{\NNN+1} \tb)$ (cf.~\eqref{eq:ODE_Nth:spec:c}, equivalently~\eqref{eq:mahomes:is:done:del0}) is quite special, we decompose the vector of $2 \NNN$-th order mixed-derivatives  $V_{=2\NNN}$ as 
\begin{equation}\label{eq:Vn_decompose}
V_{=2\NNN}
=
 V_{ \Del^{\NNN } }
 \oplus
 V_{ 2 \NNN \backslash \Del^{\NNN } }
 ,
 \qquad 
 V_{ \Del^{\NNN } } := 
 (\Del^{\NNN} \tvr ,  \Del^{\NNN+1} \xu , \Del^{\NNN+1} \tb)(0, \tau).
\end{equation}
The vector $V_{2 \NNN \backslash \Del^{\NNN } } $ includes all mixed derivatives of order $2\NNN$
except for three mixed derivatives, 
chosen so that the combinations of mixed derivatives 
in $V_{2 \NNN \backslash \Del^{\NNN } } $ 
and in $\VVN$ are linearly independent 
and together span the space of all mixed derivatives of order $2\NNN$. 
For example, we can choose $ V_{2 \NNN \backslash  \Del^{\NNN } }$ as all mixed derivatives of order $2\NNN$ except for $\pa_{d}^{2 \NNN} 
\tvr, \pa_d^{2 \NNN+2} \tb,  \pa_d^{2\NNN+1} \tu_d$.
Indeed, using identity \eqref{eq:iden_Delta_xu_divu}, we note that 
$(\Del^{\NNN} \tvr ,  \Del^{\NNN+1} \xu ,\Del^{\NNN+1} \tb) = (\Del^{\NNN} \tvr ,  2 (\NNN+1) \Del^{\NNN} \div \tu ,\Del^{\NNN+1} \tb)$, and thus
the linear space of 
$V_{= 2 \NNN }$ is that same as that of 
$ V_{ \Del^{\NNN } } \oplus V_{ 2 \NNN \backslash \Del^{\NNN } }$. Moreover, the mixed derivatives 
in $ V_{ \Del^{\NNN } }$ and $V_{ 2 \NNN \backslash \Del^{\NNN } }$ are linearity independent.

\vspace{0.1in}
\paragraph{\bf  Schematic form of the ODEs for $V_{=n}$ }

The ODE system for $V_{=0}$, and for $V_{\Del^{\NNN}}$ has been derived in \eqref{eq:ODE_0th_full}, \eqref{eq:ODE_Nth:spec:c} (equivalent to \eqref{eq:mahomes:is:done:del0}), respectively. 
Since the modulation functions chosen in \eqref{eq:modulated:nosym:summary} only depend on $V_{=0}$ and $V_{ \Del^\NNN}$, the ODE system for  $V_{=0}, V_{\Del^{\NNN}}$ is \emph{closed}, modulo nonlinear terms. Using the notation $\GGs$ from \eqref{norm:ODE_low:b},  we may schematically rewrite the closed system \eqref{eq:ODE_0th_full}, \eqref{eq:ODE_Nth:spec:c} as 
\begin{equation}\label{eq:scheme_ODE_V0}
 \tfrac{d}{ d \tau} (V_{=0}, V_{\Del^{\NNN}})^\intercal
 = \mathbf{m}_0 (V_{=0}, V_{\Del^{\NNN}})^\intercal 
 + \OO_{\NNN}( |\GGs_{2\NNN}|^2). 
\end{equation}
for some matrix $\mathbf{m}_0$. Since $\GGs$ from \eqref{norm:ODE_low:b} also depends on
the modulation functions $\tcvr , \tcr, \tcb, \tcu$, which in turn depend linearly 
on  $V_{=0}, V_{\Del^{\NNN}}$ by \eqref{eq:modulated:nosym:summary}, using the notation 
$\OOL{}$ in \eqref{norm:ODE_low:OOL}, we obtain
\begin{equation}\label{eq:GGs_V}
\GGs_k = \OOL{k}( V_{=0}, V_{\Delta^\NNN } ,   V_{\leq k} ), \quad \forall k \geq 0.
\end{equation}

For $ n \geq 1$, we Taylor expand equations
to \eqref{eq:lin_ODE:cL}--\eqref{eq:lin_ODE} around $y=0$, and similarly to Appendix~\ref{app:derive_ODE} we obtain the ODEs\footnote{
The derivation of~\eqref{eq:scheme_ODE_Vn} follows from Taylor expanding \eqref{eq:lin_ODE:cL}--\eqref{eq:lin_ODE} around $y=0$, and is much simpler than the proof of Proposition~\ref{prop:ODE_Nth}. This is because here in~\eqref{eq:scheme_ODE_Vn} we do not explicitly track terms other than the top order linear terms; to avoid redundancy with the arguments in Appendix~\ref{app:derive_ODE}, this proof is omitted. 
}
\begin{subequations}\label{eq:scheme_ODE_Vn}
\begin{align}
 \tfrac{d}{d \tau} 
V_{=n}^\intercal  
& = \mathbf{m}_n   V_{=n}^\intercal + \OOL{ n}( \GGs_{(n-2 \NNN)_+}) +   \OO_n ( |\GGs_n|^2), \notag \\
 & = \mathbf{m}_n   V_{=n}^\intercal + \OOL{ n}(V_0,   V_{\Delta^\NNN } ,  V_{ \leq (n-2 \NNN)_+}, V_0 ) +   \OO_n ( |\GGs_n|^2) .
 \label{eq:scheme_ODE_Vn:a}
\end{align}
for some matrix $\mathbf{m}_n$ that only depends on the profile and parameters $n,\al , d$. 
By definition in \eqref{eq:Vn_decompose},  $\VVNother$ is a sub-vector of $V_{= 2\NNN}$. 
Thus, restricting the ODEs \eqref{eq:scheme_ODE_Vn:a} with $n = 2\NNN$ to the ODEs 
 for $\VVNother $, we obtain
\begin{equation}\label{eq:scheme_ODE_Vn:b}
  \tfrac{d}{d \tau} \VVNother^\intercal  = \mathbf{m}_{2 \NNN}^{\prime} 
  \VVNother^\intercal  + \OOL{\NNN}(V_0,   V_{\Delta^\NNN } ,  V_{ \leq (n-2 \NNN)_+}, V_0 ) +   \OO_{2\NNN} ( |\GGs_{2\NNN}|^2) .
\end{equation}
\end{subequations}
Combining the ODEs \eqref{eq:scheme_ODE_V0} and \eqref{eq:scheme_ODE_Vn}, for any $n\geq 2\NNN+1$, 
we that the ODEs system for $V_{\leq n}$ is given by 
\begin{subequations}\label{eq:scheme_ODE_V_leqn}
\begin{align}
  & \tfrac{d}{d \tau} ( (V_{= 0}, \VVN ) , \, V_{=1}, .., V_{= 2 \NNN-1}, \VVNother, V_{=2\NNN+1},.., V_{=n} )^\intercal 
  \notag \\ 
  & \qquad = M_{\leq n} ( (V_{= 0}, \VVN ) , \, V_{=1}, .., V_{= 2 \NNN-1}, \VVNother, V_{=2\NNN+1},.., V_{=n} )^\intercal
  +  \OO_n ( |\GGs_n|^2),
 \end{align}
where $M_{\leq n}$ is a block lower triangular matrix (see \eqref{eq:upper_block} for an illustration)
with matrices on the diagonal blocks of $M_{\leq n}$ given by $\bmm_i$ and $\bmm_{2\NNN}^{\prime}$:
\begin{align}
& M_{\le n}  = (\td A_{ij})_{1 \le i,j \le n+1}, 
\qquad \td A_{ij}=0 \quad \text{for } i < j, 
\notag \\
& (\td A_{11}, \td A_{22}, ..., \td A_{n+1, n+1} )  = ( \bmm_0,\bmm_1,.., \bmm_{2\NNN-1}, \bmm_{2\NNN}^{\pr},\bmm_{2\NNN+1}, .., \bmm_{n}) .
\end{align}
\end{subequations}
We do not derive the matrix $M_{\leq n}$ explicitly, since its exact entries are not used in the stability analysis. At this stage, we emphasize the significance of the fact that the matrix $M_{\leq n}$ is lower block triangular: this indicates that the ODE system for the mixed derivatives 
$V_{=i}$ is allowed to depend on $V_{=j}$ for all $j\leq i-1$ linearly, but is is \emph{independent of} $V_{=j}$ for $j > i$. This fact is in accordance with~\eqref{eq:ODE_stab} above.

\subsubsection{Lower block triangular matrices} \label{sec:lower_block}
The key consequence of Proposition \ref{prop:ODE_Nth} is the following structure. 
From~$V_{=n}$ we define an \emph{overdetermined} vector of mixed derivatives of order $n$ for $\tvr$, $n+1$ for $\tu$, and $n+2$ for $\tb$; this vector is denoted as $\mathbf{P}_n V_{=n}^\intercal$, and the constant-coefficient matrix $\mathbf{P}_n$ is defined as: 
\begin{itemize}[leftmargin=2em]
\item If $n$ is even, we choose a constant matrix $\mathbf{P}_n$ such that 
\begin{subequations}\label{eq:Vn_lift}
\begin{equation}
\bal
\mathbf{P}_n V_{=n}^\intercal & = 
\Bigl( (\Del^{n/2} \tvr, \Del^{n/2+1} \xu, \Del^{n/2+1} \tb), \quad
 \underbrace{ \VVs_{ \beta_2, n/2 - 1}  }_{ \forall \beta_2 : |\beta_2| = 2  }
, \  \underbrace{  \VVs_{ \beta_4, n/2 - 2} }_{ \forall \beta_4 : |\beta_4| = 4  }
, \, ...,  \   \underbrace{  \VVs_{ \beta_n,  0} }_{ \forall \beta_n : |\beta_n| = n  }
 ,  \\
& \qquad  \underbrace{  ( \partial^{\beta_{n+2}} \xu, \, \partial^{\beta_{n+2}} \tb) }_{ \forall \beta_{n+2} : |\beta_{n+2}| = n + 2  } ,  \quad \underbrace{ \pa^{\beta_{n+1}} \tu  }_{ \forall \beta_{n+1} : |\beta_{n+1}| = n + 1  } \Bigr)^\intercal,
\eal
\end{equation}
where each $\beta_i \in \Naturals_0^d$ denotes a multi-index with $|\beta_i| = i$, 
and we order these vectors by first listing all the vectors $\VVs_{ \beta_2, n/2 - 1} $ with different $\beta_2$, and then the vectors $\VVs_{ \beta_4, n/2 - 1} $  with different $\beta_4$, and so on.  

\item If $n$ is odd, similarly, we choose a constant matrix $\PP_n$ such that 
\begin{equation}\label{eq:Vn_lift:b}
\mathbf{P}_n V_{=n}^\intercal  = 
\Big( \underbrace{ \VVs_{ \beta_1, (n-1)/2 }  }_{ \forall \beta_1 : |\beta_1| = 1  }
, \  \underbrace{  \VVs_{ \beta_3, (n-3)/2 } }_{ \forall \beta_3 : |\beta_3| = 3  }
, \, ...,     \underbrace{  \VVs_{ \beta_n,  0} }_{ \forall \beta_n : |\beta_n| = n  }
 ,  \, \underbrace{  ( \partial^{\beta_{n+2}} \xu, \, \partial^{\beta_{n+2}} \tb) }_{ \forall \beta_{n+2} : |\beta_{n+2}| = n + 2  } ,  \  \underbrace{ \pa^{\beta_{n+1}} \tu  }_{ \forall \beta_{n+1} : |\beta_{n+1}| = n + 1  }  \Big)^\intercal .
\end{equation}
\end{subequations} 
\end{itemize}

Proposition~\ref{prop:ODE_Nth} implies that we can lift the  ODEs \eqref{eq:scheme_ODE_Vn:a}  to the \emph{overdetermined} ODE systems
\begin{equation}\label{eq:ODE_Vn_lift}
    \tfrac{d}{d \tau} ( \mathbf{P}_n V_{=n}^\intercal)
    =  \mathbf{M}_n (\mathbf{P}_n  V_{=n}^\intercal) + \OOL{n} ( \GGs_{ ( n - 2 \NNN )_+} ) + \OO_n( \GGs_{ n  }^2) ,
\end{equation}
where $\mathbf{M}_n $ is a block lower triangular matrix
\bseq\label{eq:upper_block}
\begin{equation}
 \MMn =  \begin{pmatrix}
A_{11} & 0 & \cdots & 0 \\
A_{21} & A_{22} & \cdots & 0 \\
\vdots & \vdots & \ddots & \vdots \\
A_{\ell_n 1} & A_{\ell_n 2} & \cdots & A_{\ell_n \ell_n}
\end{pmatrix},
\label{eq:upper_block:a}
\end{equation}
with diagonal blocks $\{A_{ii}\}_{i=1}^{\ell_n}$ given by the matrices
$ \HHH_{n/2} $ in \eqref{eq:ODE_Nth:spec:a} (when $n$ is even),  $\HH_{|\bet|, k}$ in \eqref{eq:ODE_Nth}, $\HHT_n$ in \eqref{eq:ODE_Nth_ZB}, or the scalar $ (- n \kp - 2\uua -1 )$ in \eqref{eq:ODE_Nth_U}, arranged in the following manner as $i=1,\ldots,\ell_n$ increases:
\beq
\underbrace{    \HHH_{ k}, \mbox{with \ }  k = \tfrac{n}{2}}_{ \mathrm{one \ matrix}, n \, \mathrm{ even}},   \quad 
 \underbrace{  \HH_{ |\bet| , k},  \mbox{with \ } |\beta| + 2 k = n}_{ 
 \mathrm{for \ each \ } 0\leq k \leq \frac{n-1}{2}, \ 
\mathrm{repeat \ p_{d, n-2k} \ times }  },  \quad \underbrace{ \HHT_n }_{ \mathrm{  repeat \ p_{d, n+2} \ times}
 }, 
\quad  \underbrace{ - n \kp - 2 \uua -1}_{ \mathrm{  repeat \ p_{d, n+1} \ times} } ,
\label{eq:upper_block:b}
\eeq
where $p_{d,k}$ denotes 
\beq\label{eq:upper_block_index}
p_{d,k} = | \{    \beta \in \Reals^d : |\beta | = k  \}  | 
 = \binom{d+k-1}{d-1}.
 \eeq
 \eseq
Each matrix $A_{ii}$ on the diagonal of $\MMn$ corresponds to a matrix for the main linear terms in the ODE systems \eqref{eq:ODE_Nth:spec:a} (one $3\times3$ ODE, when $n$ is even), \eqref{eq:ODE_Nth} ($p_{d,n-2k}$ different $4\times 4$ ODEs, for each $0\leq k\leq \frac{n-1}{2}$), \eqref{eq:ODE_Nth_ZB} ($p_{d,n+2}$ different $2\times 2 $ ODEs), and \eqref{eq:ODE_Nth_U} ($p_{d,n+1}$ different scalar ODEs).

\vspace{0.1in}
\paragraph{\bf Relation between the original and the lifted ODE systems}
We view the ODE system \eqref{eq:ODE_Vn_lift} as a ``lift'' of the ODE system \eqref{eq:scheme_ODE_Vn} because the 
linear space of the mixed derivatives in $\mathbf{P}_n V_{=n}^\intercal$ \eqref{eq:Vn_lift} 
is the same as that of $V_{=n}$. Moreover, the ODE system \eqref{eq:ODE_Vn_lift} is overdetermined since the dimension of $\mathbf{P}_n V_{=n}^\intercal$ in \eqref{eq:Vn_lift} is larger than that of $V_{=n}$. 
As a result, we obtain $\mathbf{P}_n \in \Reals^{ a_n \times b_n }$ for some  $a_n \geq b_n$ and $\mathbf{P}_n$ has full column rank $b_n$.  By comparing \eqref{eq:scheme_ODE_Vn} and \eqref{eq:ODE_Vn_lift}, we obtain
\begin{equation}\label{eq:ODE_lift_relation}
\mathbf{P}_n \mathbf{m}_n = \mathbf{M}_n   \mathbf{P}_n,
\quad    \mathbf{P}_n \in \Reals^{ a_n \times b_n },
\quad   \mathbf{m}_n  \in \Reals^{b_n \times b_n },
\quad  \mathbf{M}_n  \in \Reals^{ a_n \times a_n },  
\quad  a_n \geq b_n.
\end{equation}
Since we do not use the exact values of $a_n, b_n$ in the proof, although these integers are explicitly computable, we do not specify their values.

Recall the notation for multiplicity $\Algm(\lam, M)$ from \eqref{def:Algm}.  Given any real-valued matrix $M$ and a real eigenvalue $\lam$ of $M$, using Jordan norm form, we obtain 
\footnote{
Note that the algebraic multiplicity $\Algm(\lam, M)$ equals the total number of times that $\lam$ appears on the diagonal of the Jordan normal form. 
While the identity \eqref{eq:iden_Algm_dim} can be extended to complex eigenvalues $\lam$, we do not pursue this generalization, as we only apply \eqref{eq:iden_Algm_dim} for $\lam \in \Reals$.
}
\beq\label{eq:iden_Algm_dim}
 \dim_{\Reals}( \ker ( \lam \Id - M)^{ \Algm( \lam, M)} ) = \Algm(\lam, M).  
\eeq
Since $\lam \Id - M$ is real, in the above identity, we view $ \ker ( \lam \Id - M)^{ \Algm( \lam, M)}$ 
as a real linear subspace.

The lifting relation \eqref{eq:ODE_lift_relation} between 
$ \mathbf{m}_n$ and $\mathbf{M}_n$ allows us to relate the eigenspaces of $\mathbf{m}_n$ 
and $\MMn$. In particular, we have the following abstract linear algebra lemma.

\begin{lemma}\label{lem:lift}
Suppose that $\bmm \in \Reals^{b \times b}, \PP \in \Reals^{a \times b},  
\mathbf{M} \in \Reals^{ a \times a}$, $a \geq b$, $\PP$ has full column rank, and 
\beq\label{eq:ass_lift}
\PP \bmm = \mathbf{M}  \PP .
\eeq
Suppose that $(\lambda, v \neq 0)$ is a generalized eigenpair of $\bmm$;  namely, $( \bmm -\lam \Id)^k v = 0$ for some integer $k\geq 1$. Then $(\lambda, \PP v \neq 0)$ is a generalized eigenpair of $\mathbf{M}$ with 
$ (\mathbf{M} -\lambda \Id)^k \PP v = 0$. 
\end{lemma}

\begin{proof}[Proof of Lemma~\ref{lem:lift}]
Using \eqref{eq:ass_lift} repeatedly, we obtain
\[
 \PP (\bmm - \lam \Id)^k = (\mathbf{M} -\lam \Id) \PP  (\bmm - \lam \Id)^{k-1}
 =  \ldots =  ( \mathbf{M} -\lam \Id)^k \PP.
 \]
Using $(\bmm - \lam \Id)^k v = 0$, we obtain $ ( \mathbf{M}  -\lam \Id)^k \PP v = 0$. Since $\PP$ 
has a full column rank and $v\neq 0$, we obtain $\PP v \neq 0$. 
We conclude the proof. 
\end{proof}
 
Due to the block lower triangular structure of $\MMn$, using that $ A_{ii} \in \Reals^{ n_i \times n_i}$ with $ n_i \leq 4$, we can \emph{derive} all the eigenvalues of $\MMn$ explicitly. 
If $\mathbf{M}_n$ has no eigenvalues with non-negative real parts, then, by Lemma~\ref{lem:lift}, $\mathbf{m}_n$ 
has no eigenvalues with non-negative real parts.

\subsubsection{Relation between $\HHH_k$ and $\HH_{0, k}$}

Since $\Del^{k+1} (y \cdot \tu )(0) = 2(k+1)  \Del^k (\div \tu) (0)$ due to the identity \eqref{eq:iden_Delta_xu_divu},  the matrices $\HHH_k \in \Reals^{3 \times 3}$ and $\HH_{ |\bet|, k} \in \Reals^{4 \times 4}$ are not completely unrelated. We have the following relation.

\begin{lemma}\label{lem:HH_connect}
For any $z_1, z_2, z_3$, denote $(v_1, v_2, v_3)^\intercal= \HHH_k (z_1, z_2, z_3)^\intercal $. Then 
\begin{equation}\label{eq:HH_connect}
(v_1, \tfrac{1}{2 k +2} v_2, v_2, v_3 )^\intercal = \HH_{0, k}(z_1, \tfrac{1}{2 k +2}z_2 , z_2, z_3)^\intercal .
\end{equation}
As a result, if $\lam$ is an eigenvalue of $\HHH_k$, then it is an eigenvalue of $\HH_{0, k}$.
\end{lemma}

\begin{proof}[Proof of Lemma~\ref{lem:HH_connect}]
Identity \eqref{eq:HH_connect} follows from the definitions of $\HHH_k$ \eqref{eq:ODE_Nth:spec}
and $\HH_{0, k}$ \eqref{eq:ODE_Nth} with $\bet = 0$ and $n = 2k+ |\bet| = 2k$.   Suppose that $(\lam, z)$ is an eigen-pair of $\HHH_k$. Then $v = \HHH_k z = \lam z$. Using \eqref{eq:HH_connect}, we obtain 
\[
 \HH_{0, k}(z_1, \tfrac{1}{2 k +2}z_2 , z_2, z_3)^T =
 (v_1, \tfrac{1}{2 k +2} v_2, v_2, v_3 )^T = \lam (z_1, \tfrac{1}{2 k +2} z_2, z_2, z_3)^T.
\]
Thus, $\lam$ is an eigenvalue of $\HH_{0, k}$.
\end{proof}

\subsection{Complete characterization of unstable directions for the ground state in five important cases}
\label{sec:unstable_modes_specific}
For the ground-state profile at $\NNN = 1$, we analyze five important special cases and compute the \emph{exact dimension} of 
the linear space of unstable 
or neutrally stable directions for the matrices $\mathbf{M}_n$ appearing in~\eqref{eq:upper_block},
and $M_{\leq n}$ appearing in
~\eqref{eq:scheme_ODE_V_leqn}. 

For this purpose, we introduce the following notation.
Let $A$ be a matrix with real coefficients, and let $\{ \lam_i\}$ be the eigenvalues of $A$. 
We introduce $  \Lam_{s}(A),  s \in \{ +, 0, -\}$ to denote the number of unstable, neutrally stable, and unstable eigenvalues of $A$: 
\begin{equation}\label{eq:eigen_count}
  \Lam_{\bullet}(A) := \bigl|\{  \lam_i(A): \text{sign}( \Re(\lam_i(A) ) ) = \bullet \} \bigr|,  \qquad \bullet \in \{ +, 0, -\},
\end{equation}
where $|S|$ denotes the cardinality of the set $S$. With this notation, we have:

\begin{theorem}[\bf Exact number of unstable directions for specific cases]\label{thm:uns_specific}
Let $\NNN=1$. Consider the following five pairs of adiabatic exponents and dimensions:   $(\gamma, d) \in \{ (\frac{5}{3}, 3), (2, 2) \}$ monotonic gas in three \& two dimensions; $(\gamma, d) \in \{ ( \frac{7}{5}, 3 ), (\frac{5}{3}, 2) \}$ diatomic gas in three \& two dimensions;  $d=1$ one space dimension.

\paragraph{\bf Part I} We have the following results concerning the exact number of stable, neutral, and unstable eigenvalues of the matrices appearing in~\eqref{eq:upper_block:b}, in the aforementioned five special cases.

\begin{enumerate}[label=(\roman*),leftmargin=2em]

\item \textsl{(Eigenvalues of $\HHH_{k}$).}
For $n=2k$, we have 
\begin{align*} 
 (\Lam_-, \Lam_0, \Lam_+)(\HHH_k)  & = (2, 1, 0), \qquad  k = 1 , \\
  (\Lam_-, \Lam_0, \Lam_+)(\HHH_k) & = (3, 0, 0),  \qquad  \forall \, k \geq 2 .
\end{align*}

\item \textsl{(Eigenvalues of $\HH_{ |\bet|, k}$).}
For $n = |\bet| + 2k \geq 1$, we have 
\begin{align*} 
   (\Lam_-, \Lam_0, \Lam_+)(\HH_{ |\bet|, k})  & = (4, 0, 0), \qquad \forall |\bet| + 2k \geq 2, \ (\bet,k) \neq (0,1) , \\
    (\Lam_-, \Lam_0, \Lam_+)(\HH_{ |\bet|, k})  & = (3, 0, 1),  \qquad  (\bet,k) = (1, 0). 
   \end{align*}
Note that when $n = 1$, we have $|\bet| = 1$ and $k = 0$.

\item \textsl{(Eigenvalues of $\HHT_n$).}
We have 
\begin{align*} 
(\Lam_-, \Lam_0, \Lam_+)(\HHT_{n}) & = (2, 0,0) , \qquad  \forall  \ n \geq 1, \ (\gamma, d) \in \{ (\tf53, 3) , \ ( 2, 2 ) \},  \\
(\Lam_-, \Lam_0, \Lam_+)(\HHT_{n}) & = (2, 0,0) , \qquad \forall  \ n \geq 2, \ (\gamma, d) \in \{ (\tf75, 3),  (\tf 53, 2) \},  
\\
(\Lam_-, \Lam_0, \Lam_+)(\HHT_{n}) & = (1, 0, 1) , \qquad  \quad   n = 1,   \ (\gamma, d) \in \{   (\tf75, 3),   (\tf 53, 2)\} .
\end{align*}

\item \textsl{(Sign of the ODE for $\tu$).}
Consider  $(\gamma, d) \in \{ (\tf53, 3) ,   (\tf75, 3),   ( 2, 2 ),  (\tf53, 2)\} $.  For any $n \geq 1$, we have
\begin{equation}\label{eq:sign_U_ODE}
  - n \kp - 2 \uua -1  <0 .
\end{equation}

\item \textsl{(Eigenvalues in the one-dimensional case).}
For $d = 1$, we have 
\begin{align*} 
(\Lam_-, \Lam_0, \Lam_+)(\HH^{(1)}_{n}) & = (3,0,0), \qquad \forall  \ n \geq 3, \\
(\Lam_-, \Lam_0, \Lam_+)(\HH^{(1)}_{n}) & = (2,1,0), \qquad n = 2.  \\
(\Lam_-, \Lam_0, \Lam_+)(\HH^{(1)}_{n}) & = (2,0,1), \qquad n = 1. 
\end{align*}
\end{enumerate}
When $n=1$, $\HH_{1, 0}$ has exactly one positive eigenvalue. When $n=1, 
(\gamma, d)  \in \{ (\tf 75, 3), (\tf 53,  2 )  \}$, $\HHT_1$ has exactly one positive eigenvalue.

\vspace{0.1in}

\paragraph{\bf Part II} 
Consider $n \geq 2$. Recall from~\eqref{eq:ODE_stab} the matrix $M_{\leq n}$ appearing in the closed ODE system for $V_{\leq n}$, where $V_{\leq n} = (\nabla^{\leq n} \tvr, \nabla^{\leq n+1} \tu, \nabla^{\leq n+2} \tb)$. Denote 
the eigenvalues of $M_{\leq n}$ with non-negative real parts 
by $ \{ \lam_{ M_{\leq n}, i} \}_{i=1}^{k_{\leq n}}$, and define the linear subspace $\Sigma_{\mathsf{uns}, \leq n}$ as in~\eqref{def:Sigma_uns_n}.
We have the following results:
\begin{enumerate}[label=(\roman*),leftmargin=2em]

\item \textsl{(Real eigenvalues).} The eigenvalues $\lam_{ M_{\leq n}, i}$ are real, for all $1\leq i\leq k_{ \leq n}$. 

\item \textsl{(Structure of 
$\Sigma_{\mathsf{uns}, \leq n}$)}
For any $n \geq 2 $, a column vector 
$ v_{\leq n} $ belongs to $ \Sigma_{\mathsf{uns}, \leq n}$ if and only if 
\[
 v_{\leq n} = ( v_{\leq 2}^\intercal, 0,.., 0)^\intercal , 
\]
for some column vector $v_{\leq 2} \in  \Sigma_{\mathsf{uns}, \leq 2}$. As a result, for any $n \geq 2$, $\Sigma_{\mathsf{uns}, \leq n}$ is a trivial lifting of $\Sigma_{\mathsf{uns}, \leq 2}$ from $\Reals^{ d_{\leq 2}}$ to $\Reals^{ d_{\leq n}}$, and $ \dim ( \Sigma_{\mathsf{uns}, \leq n}) = \dim ( \Sigma_{\mathsf{uns}, \leq 2})$.

\item \textsl{(Dimension of $\Sigma_{\mathsf{uns, \leq 2} }$).}
For the following five particular cases we have 
\bseq\label{eq:dim_Sigma_uns}
\begin{align}
\dim( \Sigma_{\mathsf{uns, \leq 2} } ) & = d^2 - 1 + d  , && \mw{for \ } (\gamma, d) \in \bigl\{ (\tf53, 3) ,  \ ( 2, 2 )\bigr\}, \ \mw{or}  \ \gamma \in (1,3], d = 1 ,  \label{eq:dim_Sigma_uns:a}  \\ 
\dim( \Sigma_{\mathsf{uns}, \leq 2 } ) & = 18  , \quad && \mw{for \ } (\gamma, d) = (\tf75, 3), \\
\dim( \Sigma_{\mathsf{uns, \leq 2}} )  & = 7 , && \mw{ for \ }    (\gamma, d) = ( \tfrac53, 2 ) .
\label{eq:dim_Sigma_uns:b}
\end{align}
\eseq

\end{enumerate}
\end{theorem}

\begin{proof}[Proof of Theorem~\ref{thm:uns_specific}]

\textbf{Proof of Part I}
The matrices in Theorem~\ref{thm:uns_specific} are given explicitly in Proposition~\ref{prop:ODE_Nth}. We prove the results in Step I  using the Routh–Hurwitz criterion, which allows us to analyze the number of unstable eigenvalues; this computation involves proving a number of inequalities connecting the parameters $\al, d, k,\beta$. We defer the proof of Part I to Appendix~\ref{sec:eigen}.

\vspace{0.1in}
\paragraph{\bf Proof of Part II}
We consider $n \geq 2$. We recall the matrix $M_{\leq n}$ and the ODEs for $\nabla^{\leq n} \tvr$, $\nabla^{\leq n+1} \tu$, and $\nabla^{\leq n+2} \tb$ from \eqref{eq:scheme_ODE_V_leqn}.

\vspace{0.1in}
\paragraph{\bf Proof of Part II, item (i)}
From the results in Part I of the Theorem, all the eigenvalues of $ \HHH_{ k} $ with $k\geq 2$, $\HH_{|\bet|, k}$ with $|\beta| + 2 k \geq 2$, $\HHT_n$ with $n\geq 2$, and $ (-  n \kp - 2 \uua -1 ) \Id$ with $n\geq 1$ have negative real parts. In fact, any eigenvalue $\lambda$ of the matrices $\HHH_k$ with $k \geq 1$, $\HH_{|\beta|, k}$ with $|\beta| + 2k \geq 1$, $\HHT_n$ with $n \geq 1$, and $(-n \kp - 2 \uua -1 ) \Id$ with $n \geq 1$ that has nonnegative real part, is simple. In particular, such an eigenvalue $\lambda$ must be real.

Using this information, using the lower triangular block structure of $\mathbf{M}_\ell$ for $\ell \geq 3 $ (see~\eqref{eq:upper_block:a}), using the fact that the diagonal blocks are given explicitly in~\eqref{eq:upper_block:b}, and using Lemma~\ref{lem:lift}, we see that all the eigenvalues of $\mathbf{M}_\ell$ and $\bmm_\ell$ with $\ell\geq 3$ have negative real parts. Moreover, the eigenvalues of $\bmm_\ell$ with non-negative real part are real, for any $\ell \geq 1$. 
This proves item (i).

\vspace{0.1in}
\paragraph{\bf Proof of Part II, item (ii)}
Recall that the linear subspace $\Sigma_{ \mw{uns}, \leq n}$ is defined in~\eqref{def:Sigma_uns_n}. Let 
$n \geq 2$ and suppose that  $v_{\leq n}^\intercal \in  \ker(  (\lam \Id -  M_{\leq n} )^{  \Algm( \lam, M_{\leq n}   ) } ) $ for some 
eigenvalue $\lam$ of $M_{\leq n}$ with $\Re(\lam) \geq 0$; here $v_{\leq n}$ is a row vector. By definition, we have
\beq\label{eq:eigenvector_Mleqn}
   (\lam \Id - M_{\leq n} )^{  \Algm( \lam, M_{\leq n}   )} v_{\leq n}^\intercal = 0.
\eeq 
Since $M_{ \leq n}$ is the matrix associated with the ODE for $V_{\leq n}$ in \eqref{eq:scheme_ODE_V_leqn}, and since $\NNN=1$, we can assume that $v_{\leq n}$ is a row vector with the following partition:
\[
v_{\leq n} = ( (v_{= 0}, v_{ \Del}) , \, v_{=1},  v_{2 \backslash \Del}, v_{= 3},.., v_{=n} ) .
\]
In the proof of item (i) above, we have shown that for each $i\geq 3$ the matrix $\bmm_i$ does not have eigenvalues with non-negative real parts. By assumption, we have $\Re(\lam) \geq 0$;  
we obtain that $ \lam \Id -\bmm_i$ is an invertible matrix for all $i\geq 3$. 
Since $M_{\leq n} - \lam \Id$ is a block lower triangular matrix  with matrices $\bmm_0 - \lam \Id, \bmm_1 - \lam \Id,.., \bmm_{2}^{\pr} -\lam \Id, .. \bmm_n - \lam \Id$ on the diagonal (see~\eqref{eq:scheme_ODE_V_leqn}),
from~\eqref{eq:eigenvector_Mleqn}
we obtain that 
\[
 v_{=i} =0 , \  \forall \  i \geq 3 ,
 \qquad 
 \Algm( \lam, M_{\leq n} ) = \Algm( \lam, M_{\leq 2}).
 \]
Here we recall from \eqref{def:Algm} that $\Algm(\cdot)$ denotes the algebraic multiplicity of an eigenvalue.  Thus, equation \eqref{eq:eigenvector_Mleqn} is equivalent to 
\[
      (\lambda \Id - M_{\leq 2} )^{ \Algm( \lam, M_{\leq 2}) }   ( (v_{= 0}, v_{ \Del}) , \, v_{=1}, v_{2 \backslash \Del}  )^\intercal  = 0, \qquad 
      v_{\leq n } = (  ( (v_{= 0}, v_{ \Del}) , v_{=1}, v_{2 \backslash \Del},  \,  0, .. , 0).
\]
The above equation is further equivalent to $ ( (v_{= 0}, v_{ \Del}) , \, v_{=1}, v_{2 \backslash \Del}  )^\intercal \in \ker(  ( M_{\leq 2 } -\lam)^{ \Algm( \lam, M_{\leq 2 }) } )$.
Using the definition of $ \Sigma_{ \mw{uns}, \leq 2 }$ from \eqref{def:Sigma_uns_n} (with $n=2$), we conclude the proof of item (ii).

\vspace{0.1in}
\paragraph{\bf Proof of Part II, item (iii)}

Since  $M_{\leq n}$ is block lower triangular due to \eqref{eq:scheme_ODE_V_leqn}, for any $\mu \in \Compl$, we have the following factorization of the characteristic polynomial
\beq\label{eq:det_factor}
\det( \mu \Id - M_{\leq n}) =   \det( \mu \Id - \bmm^{\pr}_{2}) \cdot  \prod_{0\leq i \leq n, i \neq 2} 
\det( \mu \Id - \bmm_{i} ).
\eeq
As discussed above, when for all $i\geq 3$ and any $\lam \geq 0$, we have $\Algm(\lam, \bmm_i) = 0$. Therefore, using 
\eqref{eq:det_factor}, we obtain
\beq\label{eq:Algm_factor}
 \Algm( \lam, M_{\leq n} ) 
 = \Algm( \lam, \bmm_0) 
 + \Algm( \lam, \bmm_1)
 + \Algm( \lam, \bmm_2^{\pr}) 
 = \Algm( \lam, M_{\leq 2} ), \quad \forall \lam \geq 0. 
\eeq
Recall from \eqref{def:Sigma_uns_n} the definition of $\Sigma_{\mathsf{uns}, \leq n}$. 
Since $\lam_{M_{\leq n}, i} $ is real (cf. Part II, item (i)), using identity \eqref{eq:iden_Algm_dim} and the notation $\Algmpos$ from \eqref{def:Algm}, we obtain
\[
 \dim(\Sigma_{\mathsf{uns}, \leq n} )
 = \sum_{ i\leq k_{\leq n}}
   \dim(  \ker ( \lam_{M_{\leq n}, i} \Id - M_{\leq n} )^{ \Algm( \lam_{M_{\leq n}, i}, M_{\leq n} )} )
   = \sum_{ i\leq k_{\leq n}} 
   \Algm( \lam_{M_{\leq n}, i}, M_{\leq n} ) .
\]
Summing the identity \eqref{eq:Algm_factor} over all eigenvalues of $M_{\leq n}$ with non-negative real part, we obtain
\beq\label{eq:dim_uns_factor}
 \dim(\Sigma_{\mathsf{uns}, \leq n} ) = \Algmpos( M_{\leq n}) 
 =  \Algmpos( M_{\leq 2}) =
 \Algmpos(  \bmm_0) +  \Algmpos(  \bmm_1) +   \Algmpos( \bmm_2^{\pr}).
\eeq
Thus, to prove the results claimed in item (iii) of Part II, we only need to count the multiplicity 
of eigenvalues of $\bmm_0, \bmm_1, \bmm_2^{\pr}$ with non-negative real part.

Note that we have computed $ \Algmpos(  \bmm_0)$ in Sections \ref{sec:ODE_0th}--\ref{sec:ODE_N1}:
\beq\label{eq:Algmpos_m0}
   \Algmpos( \bmm_0 ) = d^2 - 1, \   \mbox{ for \ } (d, \gamma) \in \bigl\{ (3, \tfrac53), \, (3, \tf75 ), \, (2, 2), \, (2, \tf53)\bigr\}, \quad \mbox{or for} \quad d = 1 .
\eeq
It thus remains to calculate $\Algmpos(  \bmm_1)$ and $\Algmpos( \bmm_2^{\pr})$.

We recall the system \eqref{eq:scheme_ODE_V_leqn} with $n=2, \NNN=1$ and the matrix $M_{\leq 2}$:
\begin{subequations}\label{eq:ODE2_before_lift}
\begin{align}
  &  \tfrac{d}{d \tau} ( (V_{= 0}, V_{\Delta} ) , \, V_{=1} , V_{2 \backslash \Delta} )^\intercal  = M_{\leq 2}  ( (V_{= 0} , V_{\Delta} ) , \, V_{=1} , V_{2 \backslash \Delta}  )^\intercal
  +  O ( |\GGs_2|^2), \\
&  M_{\le 2}  = (\td A_{ij})_{1 \le i,j \le 3}, \ \td A_{ij}=0, \quad \text{for } i>j, 
\qquad
 (\td A_{11}, \td A_{22},  \td A_{33} )  = (\bmm_0, \bmm_1, \bmm_2^{\pr}).
\end{align}
\end{subequations}
where the matrix $\bmm_2^{\pr}$ is the square matrix in the ODEs for $V_{2 \backslash \Delta}$ in \eqref{eq:scheme_ODE_Vn:b} with $\NNN=1$. 
In the following three steps, we show that all the eigenvalues of the matrix $\bmm_2^{\pr}$ are negative (so that $\Algmpos( \bmm_2^{\pr}) =0$), and we compute $\Algmpos(\bmm_1)$.

\vspace{0.05in}
\paragraph{\bf Step 1: Case $d=1$ }

When $d=1$, since there are no mixed derivatives and the Laplacian $\Del$ reduces to $\pa_1^2$, from the definition of  $V_{2 \backslash \Delta}$ in \eqref{eq:Vn_decompose}, we obtain $V_{2 \backslash \Delta} = \emptyset $. 
Thus $ \bmm_2^{\pr} = \emptyset$. 

By the definition of $V_{=1}$ (see~\eqref{def:TL_nth_order}), we obtain $V_{=1} 
= (\pa_1 \tvr, \pa_1^2 \tu, \pa_1^3 \tb)$. Thus, the ODE for $V_{=1}$ is given exactly by \eqref{eq:ODE_Nth:1D} (with $n=1$); appealing also to~\eqref{eq:scheme_ODE_Vn} with $n=1$, we conclude that $\bmm_1 =  \HH^{(1)}_1$. 
From Part I item (v) in this Theorem, we know that $\bmm_1$ has exactly one unstable eigenvalue, which is simple. Thus, 
for $d=1$, we have 
\beq\label{eq:Algmpos_1D}
 \Algmpos( \bmm_2^{\pr}) = 0, \qquad  \Algmpos( \bmm_1 ) = 1.
\eeq

We are left to consider the four interesting cases $(\gamma, d) \in \{ (\tf53, 3) ,  (\tf75, 3),  ( 2, 2 ), (\tf53, 2)\}$; specifically, we analyze $\Algmpos(\bmm_2^{\pr})$ and $\Algmpos( \bmm_1)$.

\vspace{0.1in}
\paragraph{\bf Step 2: No unstable eigenvalues from $\bmm_2^{\pr}$ }

In this step, we show that $\bmm_2^{\pr}$ does not have unstable eigenvalues by relating it 
to the matrices $\HH_{|\bet|, k}$ 
in \eqref{eq:ODE_Nth}, $\HHT_n$ in \eqref{eq:ODE_Nth_ZB}, and to $ (- n \kp - 2 \uua -1 ) \Id$ in \eqref{eq:ODE_Nth_U}.

The ODE system in \eqref{eq:ODE2_before_lift} is equivalent to the overdetermined system for 
$( (V_{= 0} , V_{\Delta} ) , \, V_{=1} , \td V_{2 \backslash \Delta})^\intercal  $, where $\td V_{2 \backslash \Delta} $ denotes the mixed derivatives 
\[
  \VVs_{ \beta_2,  0}, 
\quad  ( \partial^{\beta_4 } \xu, \, \partial^{\beta_4} \tb) ,  \quad \pa^{\beta_3} \tu .
\]
Here, $\VVs_{ \beta_2, 0 }$ denotes all the 
vectors $\VVs_{ \beta, 0 } \in \Reals^4$ with multi-index $\beta \in \Naturals_0^d$ such that 
$|\beta| = 2$;
$ (\partial^{\beta_4} \xu, \, \partial^{\beta_{4}} \tb) $ denotes 
all the vectors $(\partial^{\beta} \xu, \, \partial^{\beta } \tb) \in \Reals^2$ with multi-index $\beta \in \Naturals_0^d$ such that  $|\beta| = 4$; and 
$\pa^{\beta_3} \tu$ denotes all the scalars $ \pa^{\beta} \tu_i$ with multi-index $\beta \in \Naturals_0^d$ such that $|\beta|=3$.
Using the ODEs~\eqref{eq:ODE_Nth} with $k=0, |\beta|=2$, $n=2$, the ODEs~\eqref{eq:ODE_Nth_ZB} with $|\beta|=4$, $n=2$, and the ODEs~\eqref{eq:ODE_Nth_U} with $|\beta| = 3$, $n=2$, we derive the system of ODEs satisfied by $\td V_{2 \backslash \Delta} $:
\beq\label{eq:ODE_td_V2other}
 \tfrac{d}{d \tau} \td V_{2 \backslash \Delta}^\intercal
 = \td {\mathbf{M}}_2 \td V_{2 \backslash \Delta}^\intercal
 + \OOL{} ( (\FFs_{0, 1}, \GGs_{0})) + O( |\GGs_2|^2). 
\eeq
Here recall the notation $\FFs$ and $\GGs$ from \eqref{norm:ODE_low}.
Using \eqref{eq:GGs_V}, \eqref{norm:ODE_low:OOL}, 
and the identity \eqref{eq:iden_Delta_xu_divu}, we obtain
\[
 \FFs_{0, 1} = \OOL{}( V_{\Del}),
 \quad \GGs_0 = \OOL{}(V_{=0}, V_{\Delta}).
\]
Similarly to~\eqref{eq:upper_block}, we may show that $\td {\mathbf{M}}_2$ (appearing in~\eqref{eq:ODE_td_V2other}) is a block lower triangular matrix (see \eqref{eq:upper_block} for an illustration) with matrices on the diagonal 
$\td {\mathbf{M}}_2$ given by $\bmm_i$ and $\bmm_{2}^{\prime}$:
\[
 \td {\mathbf{M}}_2 = ( C_{ij})_{1 \le i,j \le c},  
 \quad
 C_{ij} = 0, \
 \forall i < j , 
 \quad  C_{ii} \in \bigl\{ \HH_{2, 0},   
 \HHT_2, (- 2 \kp - 2 \uua -1) \Id \}.
\]
From the results in Part I of this theorem, the 
eigenvalues of the matrices $\HH_{2, 0}$,   $\HHT_2$, and $ (- 2 \kp - 2 \uua -1) \Id$ all have negative 
real parts. Due to the block lower triangular structure, we obtain
\begin{equation}\label{eq:sign_lift_M2}
  \Re( \lam) < 0,
  \qquad\mbox{for all eigenvalues }
  \; \lambda \;\mbox{ of } \; \td {\mathbf{M}}_2.
\end{equation}

Combining~\eqref{eq:ODE_td_V2other} and \eqref{eq:scheme_ODE_V0}, we obtain that 
\bseq\label{eq:ODE2_after_lift}
\begin{equation}
     \tfrac{d}{d \tau} ( (V_{= 0}, V_{\Delta}  ) , \, V_{=1} , \td V_{2 \backslash \Delta} )^\intercal  = \td M_{\leq 2}  ( (V_{= 0} , V_{\Delta}  ) , \, V_{=1} ,\td  V_{2 \backslash \Delta}  )^\intercal
  +  O ( |\GGs_2|^2) ,
\end{equation}
where $\td M_{\leq 2}$ is a block lower triangular  matrix with diagonal blocks given by the matrices
\begin{equation}
  \td M_{\le 2}  = (\td A_{ij})_{1 \le i,j \le 3}, 
  \quad
  \td A_{ij}=0, \  \text{for } i < j, 
  \quad  (\td A_{11}, \td A_{22},  \td A_{2, 2} )  = (\bmm_0, \bmm_1, \td {\mathbf{M}}_2 ).
\end{equation}
\eseq

Next, we note that by construction, we have the following relation 
\beq\label{eq:lift_P2}
( (V_{= 0} , V_{\Delta} ) , \, V_{=1} , \td V_{2 \backslash \Delta}  )^\intercal 
= \td {\mathbf{P}}_2  ( (V_{= 0} , V_{\Delta} ) , \, V_{=1} ,  V_{2 \backslash \Delta} )^\intercal ,
\quad 
\td {\mathbf{P}}_2 = \begin{pmatrix}
 \Id &  0 & 0  \\
 0 & \Id & 0 \\
 A_1 & A_2 & A_3 \\  
\end{pmatrix}
\eeq
for some constant matrices $A_i$.
Moreover,  $A_3 $, and hence also $\td {\mathbf{P}}_2$, have full column rank. 
Here, one should interpret $\Id$ as an identity matrix of an appropriate size.
Comparing the ODEs 
$ \td \PP_2 \times $ \eqref{eq:ODE2_before_lift} and \eqref{eq:ODE2_after_lift}, we obtain
\[
 \td {\mathbf{P}}_2  M_{\leq 2} = \td {\mathbf{M}}_{\leq 2}  \td {\mathbf{P}}_2, 
 \quad 
 \Longrightarrow \quad   \td {\mathbf{P}}_2 ( M_{\leq 2} - c \cdot \Id) =  ( \td {\mathbf{M}}_{\leq 2} - c \cdot  \Id) 
  \td {\mathbf{P}}_2 , \qquad \forall \, c \in \Reals.
\]
Since $M_{\leq 2}$ is a block lower triangular matrix (see~\eqref{eq:ODE2_before_lift}), if $\bmm_2^{\pr}$ has an eigenpair
$(\lam, v)$ with $\Re(\lam) \geq 0$ and a right column eigenvector $v \neq 0$, we obtain
that $ (0, 0, v^\intercal)^\intercal$ is a right eigenvector for $M_{\leq 2}$. Thus, using the above identity with $c =\lam$ and using the precise structure of $\td{ \mathbf{P}}_2$, we obtain
\[
0 = \td{ \mathbf{P}}_2 ( M_{\leq 2} - \lam \Id)  (0, 0, v^\intercal)^\intercal
= ( \td {\mathbf{M}}_{\leq 2} - \lam \Id) \td{ \mathbf{P}}_2 (0, 0, v^\intercal)^\intercal
= ( \td {\mathbf{M}}_{\leq 2} - \lam \Id) (0, 0, (A_3 v)^\intercal )^\intercal.
\]
Due to the block lower triangular structure of $\td {\mathbf{M}}_{\leq 2} $ in \eqref{eq:ODE2_after_lift}, 
we obtain
\[
0 = ( \td {\mathbf{M}}_{\leq 2} - \lam \Id) (0, 0, (A_3 v)^\intercal )^\intercal
 = (0, 0, (  ( \td {\mathbf{M}}_2 -\lam \Id) A_3 v  )^\intercal )^\intercal .
\]
Since $v\neq 0$ and $A_3$ has a full column rank, we obtain that $(\lam, A_3 v \neq 0)$ is an eigenpair for $ \td {\mathbf{M}}_2$. However, the assumption $\Re(\lam) \geq 0$ contradicts 
\eqref{eq:sign_lift_M2}. Thus, we prove that \emph{all the eigenvalues} of $\bmm_2^{\pr}$ are strictly negative.

To complete the proof, it remains to compute $\Algmpos( \bmm_1 )$ 
for $(\gamma, d) \in \{ (\tf53, 3) , \ (\tf75, 3), \ ( 2, 2 ), (\tf53, 2) \}$.

\vspace{0.05in}
\paragraph{\bf Step 3: 
Invertible map for the case $d=2, n=1$}
For $d=2$ and $n=1$, we can find a one-to-one linear map from the vector of mixed partial derivatives $V_{=1}^\intercal$ to a subset of the mixed derivatives in $\PP_1 V_{=1}^\intercal$ (as defined in~\eqref{eq:Vn_lift:b} with $n=1$). 

Denote $\ee_1 = (1,0), \ee_2 = (0,1)$. We define a matrix $\PP_s$ so that 
\beq\label{eq:lift_V1_spec}
\PP_s V_{=1}^\intercal  = 
\bigl(  \VVs_{ \ee_1, 0 } ,
\VVs_{\ee_2, 0}
 ,  ( \partial_1^3 \xu, \, \partial_1^3 \tb), ( \partial_1^2 \pa_2 \xu, \, \partial_1^2 \pa_2 \tb)  \bigr)^\intercal .
\eeq

It is not difficult to obtain that for $d=2$, the vector $V_{=1} = (\na \tvr, \na^2 \tu, \na^3 \tb)$ of all 
first-order mixed derivatives has length $ 2 + 2 \cdot 3 + 4 = 12$, and the above vector $\PP_s V_{=1}$
has length $ 4 \times 2 + 4 = 12$. Thus, $\PP_s \in \Reals^{12 \times 12}$ is a square matrix
and $V_{=1}, \PP_s V_{=1}^\top \in \Reals^{12}$.

Next, we show that $\PP_s$ is invertible. To simplify the derivation, we use the notation $\OOL{}$ from \eqref{norm:ODE_low}. It suffices to show that we can obtain all mixed derivatives in the vector 
$V_{=1} = (\na \tvr, \na^2 \tu, \na^3 \tb)$ from linear combinations of mixed derivatives in  $ \msf{RHS}_{ \eqref{eq:lift_V1_spec}}$:
\beq\label{eq:mix_inv}
\na \tvr, \na^2 \tu, \na^3 \tb 
= \OOL{}\bigl(  \VVs_{ \ee_1, 0 } ,
\VVs_{\ee_2, 0}
 ,  ( \partial_1^3 \xu, \, \partial_1^3 \tb), ( \partial_1^2 \pa_2 \xu, \, \partial_1^2 \pa_2 \tb)  \bigr) 
=\OOL{}(\msf{RHS}_{ \eqref{eq:lift_V1_spec}} ).
\eeq
Firstly, using the definition of $\VVs_{\ee_i,0}$ in~\eqref{def:VVs}, and using the definition of $\PP_s V_{=1}^\intercal$ in~\eqref{eq:lift_V1_spec}, we obtain
\beq
  \na \tvr(0),  \ \na \Del \xu ,  \ \na \Del \tb, \ \na ( \div \tu), \ 
     \pa_1^2 \na \xu , \pa_1^2 \na \tb
  = \OOL{}(\msf{RHS}_{ \eqref{eq:lift_V1_spec}} ).
  \label{eq:mix_inv_pf1}
\eeq
Next, we show that $\na^3 \xu,  \na^3 \tb  =\OOL{}(\msf{RHS}_{ \eqref{eq:lift_V1_spec}} )$. A direct calculation yields 
\[
 \na \pa_2^2 \xu = \na \Del \xu - \na \pa_1^2 \xu   =  \OOL{}(\msf{RHS}_{ \eqref{eq:lift_V1_spec}} ),
 \qquad  \na \pa_2^2 \tb = \na \Del \tb - \na \pa_1^2 \tb    =  \OOL{}(\msf{RHS}_{ \eqref{eq:lift_V1_spec}} ),
\]
which along with~\eqref{eq:mix_inv_pf1} implies 
\beq
\na^3 \tb, \na^3 \xu =  \OOL{}(\msf{RHS}_{ \eqref{eq:lift_V1_spec}} ).
  \label{eq:mix_inv_pf2}
\eeq 
Recall that $\xu = y \cdot \tu$. Using \eqref{eq:mix_inv_pf2} , for $i=1,2$, we obtain 
\[
\pa_i^2 \tu_i(0) = \tf13 \pa_i^3 ( y \cdot \tu)(0) =  \OOL{}(\msf{RHS}_{ \eqref{eq:lift_V1_spec}} ).
\]
Using the above identity and \eqref{eq:mix_inv_pf1}, we obtain
\[
 \pa_{12} \tu_1 = \pa_2 ( \div \tu - \pa_2 \tu_2) =  \OOL{}(\msf{RHS}_{ \eqref{eq:lift_V1_spec}} ), 
 \qquad 
 \pa_{12} \tu_2 = \pa_1( \div \tu - \pa_1 \tu_1) =  \OOL{}(\msf{RHS}_{ \eqref{eq:lift_V1_spec}} ).
\]
Finally, using \eqref{eq:mix_inv_pf1} and the above identity, we derive 
\[
\pa_1^2 \tu_2(0) =  \pa_2 \pa_1^2( y_1 \tu_1 + y_2 \tu_2)(0)  - 2 \pa_{12} \tu_1(0)  =   \OOL{}(\msf{RHS}_{ \eqref{eq:lift_V1_spec}} ).
\]
Similarly, we obtain $\pa_2^2 \tu_1 =  \OOL{}(\msf{RHS}_{ \eqref{eq:lift_V1_spec}} )$. Combining the above identities, \eqref{eq:mix_inv} follows; as a result, $\PP_s$, as defined in~\eqref{eq:lift_V1_spec}, is invertible. 

Following the derivation in Section \ref{sec:lower_block} with $n=1$, we obtain that $\PP_s V_{=1}^\intercal$ solves the ODE
\beq\label{eq:lift_ODE_spec}
  \tfrac{d}{d \tau} \PP_s V_{=1}^\intercal = \mathbf{M}_s \PP_s V_{=1}^\intercal
  + \OOL{} ( \GGs_{ 0} ) + \OO( \GGs_{ 1  }^2) ,
\eeq
where $ \mathbf{M}_{s}$ is a block lower triangular matrix and has a similar structure as 
\eqref{eq:upper_block} with diagonal blocks given by 
\[
 \underbrace{  \HH_{ 1 , 0} }_{ 
\mathrm{repeat \ twice }  },  \quad \underbrace{ \HHT_1 }_{ \mathrm{  repeat \ twice }
 } .
\]

The matrices $\HH_{1,0}, \HHT_1$ 
repeat twice since we have two different $\VVs_{\ee_i, 0}$ 
and two different $( \partial^{\beta_{3}} \xu, \, \partial^{\beta_{3}} \tb)$ in \eqref{eq:lift_V1_spec}. 
Moreover, by comparing \eqref{eq:lift_ODE_spec} and \eqref{eq:scheme_ODE_Vn} with $n=1$, we get 
 \[
  \PP_s \bmm_1 = \mathbf{M}_s \PP_s .
\]
Since $\PP_s$ is invertible, $\bmm_1$ and $\mathbf{M}_s$ are conjugate, which implies 
\[
  \Algmpos (\bmm_1) =  \Algmpos ( \mathbf{M}_s).
\]
From the results proven in Part I of this Theorem, when $n=1, \gamma = 2, d=2$, $\HH_{1,0}$ has one positive eigenvalue, 
while the eigenvalues of $ \HHT_1$ all have negative real part. When $n=1, \gamma = \tfrac53, d=2$, both $\HH_{1,0}$ and $ \HHT_1$ have one positive eigenvalue. By counting these eigenvalues, we obtain 
\beq\label{eq:Algmpos_m1_d2}
   \Algmpos (\bmm_1)  = 2, \quad \mbox{for \ }  \, \gamma = 2,  \, d=2 ,
   \quad   \Algmpos (\bmm_1)  = 4,  \quad \mbox{for \ }   \, \gamma = \tfrac53, \, d=2 .
\eeq

\vspace{0.1in}
\paragraph{\bf Step 4: Invertible map for the case $d=3,n=1$ }

For $d=3$ and $n=1$, we find a one-to-one linear map from the vector of mixed partial derivatives $V_{=1}^\intercal$ to a subset of the mixed derivatives in $\PP_1 V_{=1}^\intercal$ (as defined in~\eqref{eq:Vn_lift:b} with $n=1$). Below, all the mixed derivatives evaluate at $y=0$.

Denote $\ee_1 = (1,0,0), \ee_2 = (0,1,0), \ee_3 = (0,0,1)$. We define a matrix $\PP_s$ so that 
\beq\label{eq:lift_V1_spec_d3}
\PP_{s,3} V_{=1}^\intercal  = 
\bigl(  \VVs_{ \ee_1, 0 } ,
\VVs_{\ee_2, 0}, \VVs_{ \ee_3, 0 } 
 ,  \underbrace{  ( \partial^{\bet} \xu, \, \partial^{\bet} \tb) }_{ |\bet| = 3, \bet_3 \leq 1 ,
 \text{7 different } \bet },  
 \pa_{12} \tu_1, \pa_{13} \tu_1, \pa_{23} \tu_1, 
  \pa_{33} \tu_1, \pa_{13} \tu_2 
  \bigr)^\intercal .
\eeq

It is not difficult to obtain that for $d= 3$, the vector $V_{=1} = (\na \tvr, \na^2 \tu, \na^3 \tb)$ of all 
first-order mixed derivatives has length $ 3 + 3 \cdot 6 + 10 = 31$, and the above vector $\PP_s V_{=1}$
has length $ 4 \times 3 + 7 \times 2 + 5 = 31$. Thus, $\PP_s \in \Reals^{31 \times 31}$ is a square matrix
and $V_{=1}, \PP_{s, 3} V_{=1}^\top \in \Reals^{31}$.

Following Step 3, we show that $\PP_{s,3}$ is invertible by proving 
\beq\label{eq:mix_inv_d3}
\na \tvr, \na^2 \tu, \na^3 \tb 
=\OOL{}(\msf{RHS}_{ \eqref{eq:lift_V1_spec_d3}} ).
\eeq

Clearly, from the definition of $\VVs_{\ee_i, 0}$ in \eqref{def:VVs}, 
for any $i \in \{ 0,1,2,3\}, j \in \{ 0, 1, 2\}$, we have 
\beq
  \na \tvr(0),  \ \na \Del \xu ,  \ \na \Del \tb, \ \na ( \div \tu), \ 
     \pa_1^i \pa_2^{3-i} ( \xu ,  \tb),  \
          \pa_1^j \pa_2^{2- j} \pa_3 ( \xu ,  \tb) = \OOL{}(\msf{RHS}_{ \eqref{eq:lift_V1_spec_d3}} ).
\label{eq:mix_inv_d3_pf1}        
 \eeq

Next, we consider $\na^3 \tb(0)$. 
For $i \in \{1,2,3\}$, since $\pa_i (\pa_1^2 + \pa_2^2)\tb$ contains at most one derivative with respect to $y_3$, namely $\pa_3$, using \eqref{eq:mix_inv_d3_pf1}, we obtain
\[
 \pa_i \pa_3^2 \tb 
 = \pa_i \Del \tb - \pa_i ( \pa_1^2 + \pa_2^2 ) \tb  = \OOL{}(\msf{RHS}_{ \eqref{eq:lift_V1_spec_d3}} ).
\]

Using the above estimates and \eqref{eq:mix_inv_d3_pf1}, we obtain 
$\na^3 \tb =  \OOL{}(\msf{RHS}_{ \eqref{eq:lift_V1_spec_d3}} )$. 
Applying the same argument to $\xu$, we establish
\beq\label{eq:mix_inv_d3_pf2}   
  \na^3 \xu = \na^3 ( y \cdot \tu ), \ \na^3 \tb = \OOL{}(\msf{RHS}_{ \eqref{eq:lift_V1_spec_d3}} ).
\eeq

\paragraph{\bf Derivation for $\na^2 \tu$}

 For $i\in \{1,2,3\}$, using the identities 
\[
  3  \pa_i^2 \tu_i = \pa_i^3 ( y \cdot \tu)  ,  \quad 
 \pa_i  \Delta( y \cdot \tu) 
 = \pa_i ( 2 \na \cdot \tu + y \cdot \Delta \tu )
 =  \Delta \tu_i + 2 \pa_i ( \na \cdot \tu ) ,
\]
and the estimates \eqref{eq:mix_inv_d3_pf1} and \eqref{eq:mix_inv_d3_pf2}, we obtain
\beq
 \pa_i^2 \tu_i , \ \Del \tu_i =  \OOL{}(\msf{RHS}_{ \eqref{eq:lift_V1_spec_d3}} ) , \quad i \in \{1,2,3\}.
 \label{eq:mix_inv_d3_pf2b} 
\eeq

Using $\pa_1^2 \tu_1 , \Del \tu_1 =  \OOL{}(\msf{RHS}_{ \eqref{eq:lift_V1_spec_d3}} )$ from the above estimate and $\pa_i \pa_j \tu_1 =  \OOL{}(\msf{RHS}_{ \eqref{eq:lift_V1_spec_d3}} )$ for $(i,j) 
\in \{ (1,2), (1,3), (2,3), (3,3) \}$ from \eqref{eq:lift_V1_spec_d3}, 
 we obtain
\beq
 \pa_2^2 \tu_1 = (\Del - \pa_1^2 - \pa_3^3) \tu_1 \in 
  \OOL{}(\msf{RHS}_{ \eqref{eq:lift_V1_spec_d3}}) \, \Rightarrow \, \na^2 \tu_1 
  =   \OOL{}(\msf{RHS}_{ \eqref{eq:lift_V1_spec_d3}}).
  \label{eq:mix_inv_U1} 
\eeq

Using \eqref{eq:mix_inv_U1} and  \eqref{eq:mix_inv_d3_pf2}, we obtain
\[
 \pa_1^i \pa_2^{3-i} ( y_2 \tu_2)
 =  \pa_1^i \pa_2^{3-i} ( y_2 \tu_2 + y_3 \tu_3)
 =  \pa_1^i \pa_2^{3-i} ( y \cdot \tu - y_1 \cdot \tu_1 )  = \OOL{}(\msf{RHS}_{ \eqref{eq:lift_V1_spec_d3}} ).
\]

By choosing $i=2, 1, 0$ in order, we derive 
\beq\label{eq:mix_inv_U2_pf1}
\pa_1^i \pa_2^{2-i} \tu_2 
=\OOL{}(\msf{RHS}_{ \eqref{eq:lift_V1_spec_d3}} )
\eeq

Using  \eqref{eq:mix_inv_U2_pf1}, 
 \eqref{eq:mix_inv_d3_pf1}, and \eqref{eq:mix_inv_d3_pf2b} , we derive
\[
 \pa_3^2 \tu_2 = (\Del - \pa_1^2 - \pa_2^2 ) \tu_2 
 = \OOL{}(\msf{RHS}_{ \eqref{eq:lift_V1_spec_d3}}),
 \quad \pa_{23} \tu_2 = \pa_3 ( \na \cdot \tu - \pa_1 \tu_1 - \pa_3 \tu_3 ) 
 = \OOL{}(\msf{RHS}_{ \eqref{eq:lift_V1_spec_d3}}).
\]

Using $\pa_{13} \tu_2 =\OOL{}(\msf{RHS}_{ \eqref{eq:lift_V1_spec_d3}})$ and the above two estimates, we obtain
\beq\label{eq:mix_inv_U2}
\na^2  \tu_2  =\OOL{}(\msf{RHS}_{ \eqref{eq:lift_V1_spec_d3}} )
\eeq

Finally, using \eqref{eq:mix_inv_d3_pf2}, 
\eqref{eq:mix_inv_U1}, and \eqref{eq:mix_inv_U2}, 
for any $\bet$ with $|\bet|=3$, we obtain 
\[
 \pa^{\bet} ( y_3 \cdot \tu_3)
 = \pa^{\bet} (y \cdot \tu - y_1 \cdot \tu_1 - y_2 \cdot \tu_2)  =\OOL{}(\msf{RHS}_{ \eqref{eq:lift_V1_spec_d3}} ) , \  \Rightarrow \  \na^2 \tu_3 =
 \OOL{}(\msf{RHS}_{ \eqref{eq:lift_V1_spec_d3}} ).
\]

Combining the above estimates, we prove \eqref{eq:mix_inv_d3}. Thus, $\PP_{s,3}$ defined via \eqref{eq:lift_V1_spec_d3} is an invertible square matrix. 

Following the argument in Step 3, we derive the following ODE for 
\beq
  \tfrac{d}{d \tau} \PP_{s,3} V_{=1}^\intercal = \mathbf{M}_{s,3} \PP_{s,3} V_{=1}^\intercal
  + \OOL{} ( \GGs_{ 0} ) + \OO( \GGs_{ 1  }^2) ,
\eeq
where $ \mathbf{M}_{s,3}$ is a lower triangular matrix and has a similar structure as 
\eqref{eq:upper_block} with diagonal blocks given by 
\[
 \underbrace{  \HH_{ 1 , 0} }_{ 
\mathrm{repeat \ 3\ times }  },  \quad \underbrace{ \HHT_1 }_{ \mathrm{  repeat \ 7 \ times}
 }, 
\quad  \underbrace{ - 2 \kp - 2 \uua -1}_{ \mathrm{  repeat \ 5 \ times} } .
\]

We have the above structure of the diagonal blocks since we have $3$ different $\VVs_{\ee_i, 0}$, 
$7$ different $( \partial^{\beta_{3}} \xu, \, \partial^{\beta_{3}} \tb)$, 
and $5$ different $\pa^{\bet_2} \tu_i$ in \eqref{eq:lift_V1_spec_d3}. Moreover, $ \mathbf{M}_{s,3}$ and $\bmm_1$ are conjugate 
\[
   \PP_{s,3} \bmm_1 = \mathbf{M}_{s, 3}\PP_{s, 3} , \quad \PP_{s, 3}  \quad \mw{ is  \ invertible}.
\]

From the results proven in Part I of this Theorem, when $n=1, (\gamma ,d)= (\tf53, 3)$, $\HH_{1,0}$ has one positive eigenvalue, 
while the eigenvalues of $ \HHT_1$ all have negative real part. When $n=1, (\gamma ,d)= (\tf75, 3) $, both $\HH_{1,0}$ and $ \HHT_1$ have one positive eigenvalue. 
In both cases, we have 
$- 2 \kp - 2 \uua -1 < 0$ from \eqref{eq:sign_U_ODE} of this Theorem.  By counting these eigenvalues
and following the argument in Step 3, we obtain 
\beq\label{eq:Algmpos_m1_d3}
   \Algmpos (\bmm_1)  = 3, \quad \mbox{for \ }  \, \gamma = \tf53,  \, d=3 ,
   \quad   \Algmpos (\bmm_1)  = 3 + 7=10,  \quad \mbox{for \ }   \, \gamma = \tfrac75, \, d= 3.
\eeq

Combining identities \eqref{eq:Algmpos_m0}, \eqref{eq:Algmpos_1D}, \eqref{eq:Algmpos_m1_d3}, \eqref{eq:Algmpos_m1_d2} 
to \eqref{eq:dim_uns_factor}, we have thus proven:
\[
\bal
\dim( \Sigma_{\mathsf{uns}, \leq 2 } ) & = d^2 - 1 + d  , \quad && \mw{for \ } (\gamma, d) \in \{ (\tf53, 3) ,  ( 2, 2 ) \}, \ \mw{or}  \ \gamma > 1, d = 1 ,  \\ 
\dim( \Sigma_{\mathsf{uns}, \leq 2 } ) & = 18  , \quad && \mw{for \ } (\gamma, d) = (\tf75, 3), \\
\dim( \Sigma_{\mathsf{uns}, \leq 2} )  & = 7 , \quad  && \mw{ for \ }    (\gamma, d) = ( \tfrac53, 2 ) ,
\eal 
\]
which completes the proof of item (iii) in Part II of Theorem~\ref{thm:uns_specific}.  
\end{proof}

\begin{corollary}\label{cor:stable}
Consider the ground state profile at $\NNN=1$ and the five special cases:  
\[
(\gamma, d) \in \bigl\{ (\tf53, 3) , \ (\tf75, 3), \ ( 2, 2 ), ( \tf53, 2)\bigr\}, \quad \mbox{or} \quad \gamma \in (1,3], d = 1.
\]
For any $n \geq 2$, there exists $\delta = \delta(n, M_{\leq n}) > 0$ such that the following statement holds. 
If the initial data for the mixed-derivatives satisfies $V_{\leq 2}(0) = 0$ 
and $ \| V_{\leq n}(0) \|_{2} \leq \delta$, then 
$V_{\leq 2}(\tau) = 0 $ for any $\tau \geq 0$ and the exponential decay estimate
\[
    E_{O, n }(\tau ) \leq C_n  e^{-\lam \tau } E_{O, n  }(0) , \quad  \forall \ \tau \geq 0 , 
\]
holds for some $\lam >0$ depending on the matrix $M_{ \leq n}$ and $n$.
 \end{corollary}

\begin{proof}[Proof of Corollary~\ref{cor:stable}]
Since $V_{\leq 2}|_{\tau=0} = ( \na^{\leq 2} \tvr, \na^{\leq 3} \tu, \na^{\leq 4} \tb)|_{\tau = 0, y=0} = 0$, since
and $V_{\leq 2}(\tau)$ obeys a closed ODE system when $\NNN=1$ (see~\eqref{eq:scheme_ODE_V_leqn}), we obtain 
$V_{\leq 2}(\tau) \equiv 0$ for any $\tau \geq 0$. 
Since the modulation function $\tcr, \tcvr, \tcu, \tcb $ are linear combinations of $V_{=0}, V_{\Delta}$ (see~\eqref{eq:modulated:nosym:summary}), 
we also obtain that $\tcr(\tau), \tcvr(\tau), \tcu(\tau), \tcb(\tau) = 0$ for all $\tau\geq 0$. Thus, using the definition of 
$E_{O,n}$ in \eqref{norm:ODE}, we have 
\beq\label{eq:EOM_reduce}
E_{O,n}  =    |\tvr(0)|_{\cC^{n}}
 + |\tu(0)|_{\cC^{n+1}} + |\tb(0)|_{\cC^{n+2}}  .
\eeq

For $n \geq 3$, the ODE system \eqref{eq:scheme_ODE_V_leqn} reduces an ODE system for $
V_{3 \leq \ast \leq n}=  (  V_{=3}, .., 
V_{=n} )$, which may be written as 
\beq\label{eq:reduce_ODE}
   \tfrac{d}{d \tau} 
  V_{3 \leq \ast \leq n}^\intercal   = M_{3 \leq \ast \leq n}   V_{3 \leq \ast \leq n}^\intercal
  +  \OO_n ( | V_{3 \leq \ast \leq n} |^2),
\eeq
where $M_{3 \leq \ast \leq n}$ is a submatrix of $M_{\leq n}$ in \eqref{eq:scheme_ODE_V_leqn} and 
is block lower triangular:
\begin{align*}
 & M_{3 \leq \ast \leq n}  = (\td A_{ij})_{3 \le i,j \le n+1}, 
\qquad \td A_{ij}=0 \quad \text{for } i < j, 
\notag \\
& (\td A_{33}, \td A_{44}, \ldots , \td A_{n+1, n+1} )  =  (\bmm_{3}, \bmm_4, \ldots , \bmm_{n}) .
\end{align*}
Since $M_{3 \leq \ast \leq n}$  is block lower triangular,  
and since all the eigenvalues of the diagonal blocks $\bmm_\ell$ (with $l \geq 3$) have negative real  (see the proof of item (i) in Part II of Theorem~\ref{thm:uns_specific}), we obtain that all the eigenvalues of
$M_{3 \leq \ast \leq n}$ have negative real part. Thus, using the classical nonlinear stability result
of an equilibrium in ODE theory, see e.g.~\cite[Chapter 4.3, Theorem 2]{braun1983differential},
we obtain that there exists a small $\delta = \delta(n, M_{\leq n}) >0$ such that 
\[
\mbox{if} \quad 
 \| V_{3 \leq \ast \leq n}(0) \|_2 \leq \delta, 
\quad \mbox{then} \quad
\| V_{3 \leq \ast \leq n}( \tau ) \|_2 \leq C_n e^{- \lam \tau } \| V_{3 \leq \ast \leq n}( 0 ) \|_2
\]
for some $\lam
= \lam( M_{\geq 3, \leq n}, n) >0$ and all $\tau \geq 0$. Using identity \eqref{eq:EOM_reduce}, and recalling definition~\eqref{def:TL_nth_order}, we complete the proof. 
 \end{proof}

\subsection{Upper bound of unstable directions in general case}
\label{sec:unstable_modes_general}

For a general radially symmetric, globally self-similar profile with $\NNN\in \Naturals$, dimension $d\in \{1,2,3\}$, and for $1<\gamma \leq 2d+ 1$, we have the following estimates on the number of unstable directions for \eqref{eq:ODE_stab}.

\begin{theorem}\label{thm:uns_general}
Fix $d \in \{ 1, 2, 3 \}, 1 < \gamma \leq 2d+1$, $\NNN \geq 1$, and let $\alpha = \frac{\gamma-1}{2}$.
\begin{itemize}[leftmargin=2em]
\item[(i)] For any 
 \[
 n \geq n_1:= (18 d)^2  \NNN,
 \]
 and any $\bet, k$ with $|\bet| + 2 k = n$, the real part of eigenvalues $\lam$ of the matrices $\HH_{|\bet|, k}, \HHH_k, \HHT_{n}$ in \eqref{eq:ODE_Nth} are strictly negative and satisfy $\Re(\lam) \leq -1$. As a result, 
 the real part of the eigenvalues $\lam$ of the matrix $\mathbf{m}_n$ for the ODE system \eqref{eq:scheme_ODE_Vn:a} are strictly negative and satisfy $\Re(\lam) \leq -1$.

\item[(ii)] Recall the notation $\Sigma_{\mw{uns},\leq n}$ from \eqref{def:Sigma_uns_n}.  For any $n \geq n_1 $, the column vector 
$ v_{\leq n} $ belongs to $ \Sigma_{\mathsf{uns}, \leq n}$ if and only if 
\[
 v_{\leq n} = ( v_{\leq n_1}^\intercal, 0,.., 0)^\intercal , 
\]
for some column vector $v_{\leq n_1} \in  \Sigma_{\mathsf{uns}, \leq n_1}$. As a result, for any $n \geq n_1$, $\Sigma_{\mathsf{uns}, \leq n}$ is a trivial lifting of $\Sigma_{\mathsf{uns}, \leq n_1}$ from $\Reals^{ d_{\leq n_1}}$ to $\Reals^{ d_{\leq n}}$, and $ \dim ( \Sigma_{\mathsf{uns}, \leq n}) = \dim ( \Sigma_{\mathsf{uns}, \leq n_1})$.

\item[(iii)] We have 
\[
 \dim ( \Sigma_{\mathsf{uns}, \leq n_1}) \leq   C( d)  \NNN^d.
\]
for some explicitly computable constant $C(d)>0$ which depends only on $d$.
\end{itemize}
\end{theorem}

To prove Theorem \ref{thm:uns_general}, we use a classical result on the eigenvalues of a weighted diagonally-dominated matrix, a version of the Gershgorin circle theorem. 

\begin{lemma}[\bf Gershgorin-type bound]
\label{lem:diag_eig}
Let $A \in \Reals^{n \times n}$. Suppose that for every $1\leq i\leq n$ we have $A_{ii} < 0$, and assume that there exists $\mu_i >0$ and $\theta> 0$ such that
\begin{equation}\label{eq:diag_test}
   (|A_{ii}| - \theta) \mu_i \geq \sum_{j \neq i} |A_{i j}\mu_j |.
\end{equation}
Then all eigenvalues $\lambda$ of $A$ have negative real and $\Re( \lam_i ) \leq - \theta$.
\end{lemma}

For the sake of completeness, we provide the proof.

\begin{proof}[Proof of Lemma~\ref{lem:diag_eig}]
Let $(\lam, v )$ be an arbitrary eigen-pair of $A$: $A v = \lam v$ with $v \neq 0$. Let $i_1 = \arg\max_i \frac{|v_i|}{\mu_i}$. We have $v_{i_1}\neq 0$ and
\[
  \sum_j A_{i_1 j}  v_j 
  = (A v)_{i_1}
  = \lam v_{i_1} ,
  \quad (\lam - A_{i_1 i_1}) v_{i_1} =  \sum_{j \neq i_1} A_{i_1 j} v_j. 
\]
Using \eqref{eq:diag_test}, the definition of $i_1$, and the fact that $|v_{i_1}| > 0$, we obtain 
\[
  |(\lam - A_{i_1 i_1})| \, |v_{i_1}| =
  |(\lam - A_{i_1 i_1}) v_{i_1}| \leq \max_j\frac{|v_j|}{ \mu_j  } \cdot \sum_{j \neq i_1} |A_{i_1, j} \mu_j| 
\leq \frac{|v_{i_1}|}{\mu_{i_1}} ( |A_{i_1, i_1}| - \theta) \mu_{i_1} 
=  ( |A_{i_1, i_1}| - \theta) | v_{i_1}|. 
\]
Canceling $|v_{i_1}|\neq 0$, we obtain that  $|\lam - A_{i_1 i_1} | \leq |A_{i_1 i_1}| - \theta$. Since $A_{i_1, i_1} < 0$, it thus follows that $\lambda$ lies in the disk of radius $|A_{i_1 i_1}| - \theta$ centered at $-|A_{i_1 i_1}|$, which then implies $\Re(\lam ) \leq - \theta < 0$.
\end{proof}

Next, we prove Theorem \ref{thm:uns_general}.

\begin{proof}[Proof of Theorem \ref{thm:uns_general}]

Our first goal is to show that if $n\geq 1$ is chosen large enough, then the matrices $\HH_{|\beta|,k}$ (as defined in~\eqref{eq:ODE_Nth:b}) with $|\beta|+2k=n$, are diagonally dominated, meaning they satisfy the assumptions of Lemma~\ref{lem:diag_eig}.
For this purpose, we choose the following weights 
\begin{subequations}\label{eq:diag_est0}
 \begin{equation} 
  (\mu_1, \mu_2, \mu_3, \mu_4) 
  := \Bigl( 1, \frac{n}{2\al \cca} , \frac{n^2}{ 2 \al \cca }, \frac{n}{\al \cca^2} \Bigr),
\end{equation}
and use the notation\footnote{We will prove later on that $( \msf{H}_{|\bet|, k} )_{ i i } < 0$, as required by Lemma~\ref{lem:diag_eig}; thus, in the definition of $d_i$, the term $-( \msf{H}_{|\bet|, k} )_{ i i }$ is the same as $| ( \msf{H}_{|\bet|, k} )_{ i i }|$, the form used in~\eqref{eq:diag_test}.}
\begin{equation}
d_i := -( \msf{H}_{|\bet|, k} )_{ i i } \mu_i - \sum_{ j \neq i } |( \msf{H}_{|\bet|, k} )_{ij}| \mu_j,
    \qquad \mbox{for} \; 1\leq i \leq 4.
\end{equation}
\end{subequations}
Here $\cca$ (see~\eqref{eq:globally:SS:exact:nonradial:e}) denotes leading order coefficient of the sound speed at $y=0$
\begin{equation}\label{def:C1}
\cca := \lim_{|y| \to 0} |y|^{-1} \bc(y)
= \lim_{|y| \to 0} |y|^{-1} (\bvr \barb(y))^{1/2}  
= \aaa^{1/2} =  \alpha \bar q_0 = \tfrac{1}{ 1  + \al d} \sqrt{ \tfrac{\al \gamma d}{2} } ,
\end{equation}
where we recall $\aaa = (\alpha \bar q_0)^2 $ and $\uua $ from \eqref{eq:rho0:U1:B2}. 
Since $\gamma = 1 + 2 \al < (1 + \al d)^2$, we obtain 
\beq\label{eq:bound_C1}
 \cca < \sqrt{ \al d}. 
\eeq

Next, we show that $d_i >0$ for $n$ sufficiently large.

Recall from \eqref{eq:ODE_main:kp} that $\kp = \bcr + \uua$. Dividing the outgoing inequality \eqref{eq:profile_outgoing} by $|y|$, evaluating the resulting inequality at $y=0$, and then using the definitions of $\cca$ in \eqref{def:C1} and $\uua$ in \eqref{eq:ODE:main:profi}, we derive
\begin{equation}\label{eq:diag_est1}
 \kp + \cca \geq \kp \geq \kp - \cca \geq a , \quad 
 a :=\frac{1}{6\NNN}.
 \end{equation}

In the following estimates, we derive the leading order terms for $d_i$ (as defined in~\eqref{eq:diag_est0}) in terms of $n$. For the matrices $\HH_{|\beta|,k}$, inspecting their definition in~\eqref{eq:ODE_Nth:b}, we see that only the $(2,1)$ and $(3,1)$ components have a complicated form. 
For $(\HH_{|\bet|, k})_{2,1}$ and $(\HH_{|\bet|, k})_{3,1}$, since 
$\bb = 1$ \eqref{eq:rho0:U1:B2}, $|\bet| + 2k = n$ and $(n-k) k \leq \frac{n^2}{4}$, we may bound 
\begin{align*} 
\bigl| (\HH_{|\bet|, k})_{2,1} \bigr| &= \Big| 
   -   \frac{ ( d + 2 |\bet| + 2 (k-1) ) k    \bb }{\al}  -  \bigl( \frac{n}{\al} + \frac{2( d+ n)}{\gamma} \bigr) \bb 
\Big| 
\\
& \leq  
\frac{(d+2n -2 k -2)k + n }{\alpha} + \frac{2(d+n)}{\gamma}
\leq \frac{ n^2 }{2\al}  + \frac{ (d-2)_+ k + n + d + n}{\al}
<  \frac{n^2}{2\al} + 
\frac{4 d n}{\alpha}, 
\\
\bigl| (\HH_{|\bet|, k})_{3,1} \bigr|  
& = \Big|    -  \bigl( \frac{n}{2 \al} + \frac{2}{\gamma} \bigr) 
 (2 d + 4 |\bet| + 4k ) (k+1) \bb  \Big|
\notag\\
&\leq 
\bigl( \frac{n}{2 \al} + \frac{2}{\gamma} \bigr) \bigl( n^2 + 2d (k+1) + 4 (n-k) \bigr) 
\leq \frac{ n^3 }{2\al} + \frac{12 d n^2}{\alpha} .
\end{align*}
Here we have denoted $x_+ = \max(x, 0)$ and we have used the bounds $k+1 \leq n$, $2 d(k+1) + 4(n-k) \leq 6 d n$, and $\tf{2}{\gamma} < \frac{1}{\al}$ in the last inequality. Using \eqref{def:C1}, $\uua = - \frac{1}{1 + \al d}$ \eqref{eq:rho0:U1:B2},
$\al \leq d$, and  $\gamma > 1$, we obtain
\[
\frac{|1 + 2 \uua|}{2\alpha \cca} 
= \frac{|1 + 2 \uua|}{2\alpha \sqrt{ \frac{ \al \gamma d}{2} } \cdot \frac{1}{1 + \al d} } 
= \frac{ | 1 - \al d| }{ \al^{3/2} \sqrt{2 \gamma d}} 
< \frac{ \max(1 , \al d) }{ \al^{3/2} d^{1/2}} 
\leq \frac{d^2}{ \al^{3/2} d^{1/2}}  < \frac{d^{3/2}}{\al^{3/2}},
\quad \frac{1}{\al \gamma} < \frac{d^{1/2}}{\al^{3/2}} .
\]

Using the above estimates, $\mu_i$ from \eqref{eq:diag_est0}, $\aaa = \cca^2$ \eqref{def:C1}, 
$\al \leq d$, the bounds in~\eqref{eq:diag_est1}, 
and $n \kp +  2 \uua + 1 \geq  n \kp - |1 + 2 \uua|$, we obtain 
\begin{align*} 
d_1 
&= \mu_1 \cdot n \kp  - \mu_2 \cdot 2\al \aaa 
= n \kp - n \cca
\geq n a, \\
d_2 
&= \mu_2 \cdot ( n \kp + 2 \uua + 1) 
- \mu_1 | (\HH_{|\bet|, k})_{2,1} | 
- \mu_4 \frac{\aaa}{\gamma}\notag\\
&
\geq \frac{n^2 \kp}{2\al \cca} 
- \frac{n^2}{2\al} 
- \Bigl(\frac{|1 + 2 \uua|}{2\alpha \cca} + \frac{4 d}{\alpha} + \frac{1}{\alpha \gamma} \Bigr) n
\geq \frac{n^2 a }{2 \al \cca} 
- \frac{6 d^{3/2}}{\alpha^{3/2}} n, \\
d_3 &= \mu_3 \cdot ( n \kp + 2 \uua +1 ) 
- \mu_1  | (\HH_{|\bet|, k})_{3,1} | 
- \mu_4  \frac{(n+2)\aaa}{\gamma} \notag \\
&
\geq \frac{n^3 \kp }{2\al \cca} 
- \frac{n^3}{2\al} 
- \Bigl(\frac{|1+2\uua|}{2\alpha \cca} + \frac{12 d}{\alpha} + \frac{3}{\alpha \gamma} \Bigr) n^2
\geq \frac{n^3 a}{2\al \cca} -  \frac{16 d^{3/2}}{\alpha^{3/2}}  n^2
, \\
d_4 & = \mu_4 \cdot n \kp - 2 \mu_3 
=  
\frac{n^2 \kp}{\al \cca^2} - \frac{n^2}{\al \cca}
\geq \frac{n^2 a}{\al \cca^2}.
\end{align*}
From the above estimates, definition~\eqref{eq:bar:Q:zero}, 
the standing assumption $\alpha \in (0,d]$, 
\eqref{eq:bound_C1} for $\cca$, and  \eqref{eq:diag_est1}, we obtain that if we choose 
\begin{equation}\label{eq:diag_est2}
n\geq   (18 d)^2\; \NNN
\end{equation}
then
\begin{equation*}
  n a \geq  1 + 2 \al (\al d)^{1/2}  \cdot 18 d^{3/2} \al^{-3/2} >   1 + 2\alpha \cca \cdot 18 d^{3/2} \cdot \alpha^{-3/2},
\end{equation*}
and hence $d_i\geq \mu_i $ for all $1\leq i \leq 4$. 
Note that the above inequality and~\eqref{eq:diag_est1} also implies 
\begin{equation}
n\kappa + 2 \uua + 1
\geq 
n a
+ 
\tfrac{\alpha d - 1}{1+\alpha d}
\geq 
n a - 1
> 1.
\label{eq:diag_est3}
\end{equation}
Thus, for $n$ as in \eqref{eq:diag_est2}, condition~\eqref{eq:diag_test} is satisfied with $\theta=1$, the matrix $\HH_{|\bet|, k}$ is weighted diagonally dominated, and using Lemma \ref{lem:diag_eig}, we obtain that all eigenvalues $\lam$ of $\HH_{| \bet |, k}$ 
have negative real part and satisfy $\Re(\lam) \leq -1$.

For $n= 2k$ satisfying \eqref{eq:diag_est2}, using the relation between the eigenvalues of 
$\HHH_k$ and $\HH_{0, k}$ from Lemma \ref{lem:HH_connect}, we obtain that all eigenvalues 
$\lam$ of $\HHH_k$ have negative real part and satisfy $\Re(\lam) \leq -1$.

For the matrix $\HHT_n$ appearing in the ODE~\eqref{eq:ODE_Nth_ZB}, we note that $\HH_{|\bet| }^{(2)}$ is the lower-right $2\times 2$ diagonal block of the matrix $\HH_{|\bet|, k}$ (rows $3$ \& $4$  and columns $3$ \& $4$). 
Since the  weighted diagonal dominated condition~\eqref{eq:diag_test} 
holds with $\theta=1$ for $A = \HH_{ |\bet|, k}$, it follows that $\HH_{|\bet|}^{(2)}$ 
also satisfies this condition with $\theta = 1$; 
using~\eqref{eq:diag_est3} and Lemma~\ref{lem:diag_eig}, we obtain that all eigenvalues $\lam$ of $\HH_{|\bet|}^{(2)}$ satisfy $\Re(\lam) \leq -1$. 

Thus, for $n$ satisfying \eqref{eq:diag_est2},  all eigenvalues $\lam$ of the matrices $ \HHH_{ k} $ in \eqref{eq:ODE_Nth:spec:a}, 
$\HH_{ |\bet|, k}$ in \eqref{eq:ODE_Nth}, $\HHT_n$ in \eqref{eq:ODE_Nth_ZB}, and
$ (- n \kp -2 \uua -1 ) \Id$ in \eqref{eq:ODE_Nth_U}  have  negative real parts 
and satisfy $\Re(\lam ) \leq -1$. In view of the block lower triangular structure for $\MMn$ in~\eqref{eq:upper_block} and appealing to Lemma~\ref{lem:lift}, it follows that all eigenvalues $\lam$ of the matrices $\bmm_n$ and $\MMn$ have negative real parts with $\Re(\lam) \leq -1$ (for $n$ satisfying \eqref{eq:diag_est2}). 

Let $n_1 = (18 d)^2 \NNN$. The proof of the results in item (ii) of Theorem~\ref{thm:uns_general} is the same as the proof of the results in Part II, item (ii), of~Theorem \ref{thm:uns_specific}. 
To avoid redundancy, we omit this proof here.

For the proof of item (iii) of Theorem~\ref{thm:uns_general}, we note that the number of mix-derivatives with order less than $n_1 = (18 d)^2  \NNN$ are bounded by $C(d) \NNN^d $ for some explicitly computable $C(d)>0$. We obtain that $M_{\leq n_1} \in \Reals^{
d_{\leq n_1},  d_{\leq n_1}}$ has a size $d_{\leq n_1}$ at most $C(d) \NNN^d $. By definition of $ \Sigma_{\mathsf{uns}, \leq n_1} $ in \eqref{def:Sigma_uns_n}, it thus follows that
\[
 \dim( \Sigma_{\mathsf{uns}, \leq n_1}) \leq  d_{\leq n_1} \leq C(d) \NNN^d,
\]
which completes the proof.
\end{proof}

\subsection{Global solution to the ODEs}\label{sec:global_ODE}

In this section, we construct global solutions to the ODEs. Given a matrix $M \in \mathbb{R}^{n \times n}$, we denote by $\lambda_{M, i}, i\leq n_M \leq n$ the eigenvalues of $M$. We can decompose 
the complex space $\mathbb{C}^n$ as 
\begin{equation}\label{eq:spec_decomp}
   \mathbb{C}^n = \bigoplus_{1 \leq i \leq n_M } \ker(  (\lam_{M, i} - M )^{ \Algm( \lam_{M, i}, M) } )
\end{equation}
where $\Algm( \lam_{M, i}, M) $ is 
the multiplicity of the eigenvalues 
$\lam_{M, i}$. Given a set $\sigma$,  the spectral projection $\Pi_{\sigma}$ associated with eigenvalues in $\sigma$ is defined as the linear map
\beq\label{def:Sigma_sigma}
  \Pi_{\sigma}: \mathbb{C}^n  \to \Sigma_{\sigma} := \bigoplus_{ \lam_{M, i} \in \sigma }\ker(  (\lam_{M, i} - M )^{ \Algm( \lam_{M, i}, M) } )
\eeq
with range $\Sigma_{\sigma}$ and kernel  $\bigoplus_{ \lam_{M, i} \notin \sigma }\ker(  (\lam_{M, i} - M )^{ \Algm( \lam_{M, i}, M) } )$.

We have the following abstract theorem for solving the ODE system.
\begin{theorem}\label{thm:ODE_solve}
Let $ V(\tau) \in \Reals^l $ and $\msf{H} \in \Reals^{l \times l}$ be a $\tau$-independent matrix and 
$\cN(V, V) \in \Reals^l$ be a bilinear operator in $V$. Consider the ODE system
\begin{equation}\label{eq:ODE_global}
  \tfrac{d}{d \tau} V(\tau) = \HH V + \cN(V, V).
\end{equation}
Let $\lam_{\HH,i}$ be the eigenvalues of $\HH$ and denote\footnote{Note that $\sig_u$ can contain eigenvalues with $\Re(\lam_{\HH,i}) = 0$.}
\begin{equation}\label{eq:ODE_spectral}
\sig_s := \{\lam_{\HH,i}: \Re(\lam_{\HH, i} ) < 0 \},
\quad  \sig_u := \{\lam_{\HH,i}: \Re(\lam_{\HH, i} ) \geq 0 \} ,
\quad 
  \mu_{\HH} := -  \max_{ \lam \in \sig_s  } \lam .
\end{equation}
If $\sig_s = \emptyset$, we denote $\mu_{\HH} = 0$. Let $\Pi_{\sigma_s}, \Pi_{\sigma_u}$ be the spectral projections associated with $\sig_s, \sig_u$, respectively.
There exists $\del_* = \del_*(\HH, \cN) >0$, such that for any initial data $V_{1, 0} \in \Reals^l$ with 
$\|V_{1, 0} \|_2 \leq \del_*$,  there exists $V_{u, 0} \in \Re \Pi_u \Compl^l
=\Re \Sigma_{\sigma_u}$ and a global solution $V$ to \eqref{eq:ODE_global} such that 
\begin{subequations}\label{eq:ODE_solve_est}
\begin{align}
  V(0) = V_{1, 0} + V_{u, 0} , 
  \quad  \| V_{u, 0} \|_2 \les  \| V_{1, 0} \|_2,  \\
  \| V(\tau) \|_2 \les e^{- \frac{3}{4} \mu_{\HH} \tau} \| V_{1, 0} \|_2,
\end{align}
\end{subequations}
where the implicit constants depend only on $\HH$ and the form of $\cN$.
\end{theorem}

For completeness, we provide a simple proof, which is a variant of 
the classical method for constructing stable manifolds of nonlinear ODEs~\cite[Section 9.2]{teschl2012ordinary}. Using the splitting method in \cite{ChenHou2023a,CCSV2024}, we avoid the analysis of central manifold when $\HH$ has eigenvalues $\lam$ 
with $\Re(\lam) = 0$.

\begin{proof}[Proof of Theorem~\ref{thm:ODE_solve}]
 If $\sig_s = \emptyset$, 
 from the definition \eqref{def:Sigma_sigma}, we obtain
 $\Sigma_{\sigma_u} = \Compl^l$ 
 and $\Re \Sigma_{\sigma_u} = \Reals^l$. 
 Thus, 
 given any $V_{1,0}$, we can choose 
 $V_{u,0} = -V_{1,0}$ to obtain 
 trivial initial condition $V(0)=0$
 and a trivial solution $V(\tau) \equiv 0$. Otherwise, we have $\mu_{\HH} > 0$.

We consider the following decomposition of the ODEs \eqref{eq:ODE_global} $V = V_1 + V_2$ with $V_i$ satisfying 
\footnote{
The decomposition \eqref{eq:ODE_split} is similar to the compact perturbation in \cite{ChenHou2023a,CCSV2024}. Here, since the problem is finite-dimensional, we can change $\HH V_1$ to 
$ A V_1$ with any matrix $A$ by a compact perturbation.
}
\begin{subequations}\label{eq:ODE_split}
\begin{equation}\label{eq:ODE_split:a}
\begin{aligned} 
\tfrac{d}{d \tau} V_1 &= - \mu_{\HH} V_1 + \cN(V, V),  \quad V_1 |_{\tau =0} = V_{1, 0} , \\
\tfrac{d}{d \tau} V_2 &= \HH V_2 + K_{\HH} V_1, \quad K_{\HH} := \mu_{\HH} + \HH.
\end{aligned} 
\end{equation}

Recall the projection $\Pi_s, \Pi_u$. We solve $V_2$ using the Duhamel's formula 
\begin{equation}\label{eq:ODE_split:b}
 \begin{gathered} 
  V_2(\tau) := \cA_2(V_1)(\tau),   \\
  \quad \cA_2 (V_1)(\tau) := \Re \int_0^{\tau} e^{\HH( \tau - t )} \Pi_s ( K_{\HH} V_1)(t)  d t
 - \Re \int_{\tau}^{ \infty } e^{\HH( \tau - t )} \Pi_u ( K_{\HH} V_1)( t )  d t ,
  \end{gathered} 
\end{equation}
\end{subequations}
which determines the initial data for $V_2$  implicitly. By definition, we obtain
\[
V_2(0) \in \Re( \Pi_u \Reals^l)
\subset \Re( \Sigma_{ \sigma_u} ).
\]

Since $\Pi_u$ is the projection onto the $\HH$-invariant subspace associated with eigenvalues with $\Re(\lam_{\HH,i}) \geq 0$, by definition, for any vector $\vv \in \Reals^l$, we have 
\begin{equation}\label{eq:ODE_semi}
\begin{aligned} 
  \| e^{\HH t} \Pi_s \vv \|_2 & \les e^{- \frac78 \mu_{\HH} t} \| \vv\|_2 , \quad  \forall t \geq 0 , \\
  \| e^{\HH t} \Pi_u \vv \|_2 & \les e^{ - \frac{1}{4}  \mu_{\HH} t} \| \vv\|_2 , \quad  \forall t < 0 . \\
  \end{aligned}
\end{equation}

\paragraph{\bf Step 1: Fixed point map}

Let $V_{1, 0}$ be the initial data in Lemma \ref{thm:ODE_solve}. 
We solve \eqref{eq:ODE_split} using a fixed point method. 
We define
\begin{equation}\label{eq:ODE_Ynorm}
   \| V \|_{Y} := \| e^{  \frac{3}{4} \mu_{\HH}  \tau } V(\tau) \|_{L^{\infty}},
   \quad B_{Y}(r) := \{ V: \| V\|_Y < r,  V(0) = V_{1, 0}  \}.
\end{equation}

The set $B_Y(r)$ is nonempty if $r > \| V_{1,0}\|_2$ as $ V_{1,0} e^{- \frac{3}{4}\mu_{\HH} \tau} \in B_Y(r)$.
For any input $\hat V_1 \in Y$, we first construct $V_2 = \cA_2(\hat V_1)$ via \eqref{eq:ODE_split:b} 
and then construct $V_1$ by solving \eqref{eq:ODE_split:a}. This defines a map
\begin{equation}\label{eq:ODE_fix_pt}
  V_1 = \cA(\hat V_1).
\end{equation}

The fixed point of \eqref{eq:ODE_fix_pt} corresponds to a global solution to \eqref{eq:ODE_global}.  
Below, we show that for  
\begin{equation}\label{eq:ODE_size_init}
   \del = \| V_{1, 0}\|_2 
\end{equation}
sufficient small,  $\cA$ is a contraction mapping from $B_Y( 2 \del) \to B_Y( 2\del)$.

Using \eqref{eq:ODE_Ynorm} and \eqref{eq:ODE_semi}, for any $\vv \in Y$, we obtain 
\begin{equation}\label{eq:ODE_semi2}
\begin{aligned} 
  \| \cA_2(  \vv)(\tau)\|_2 &\les  \int_0^{\tau} e^{-\frac78 \mu_{\HH}(\tau -t)} \| (\mu_{\HH} + \HH) \vv(t)\|_2 d t 
  + \int_{\tau}^{\infty} e^{ -\frac{1}{4} \mu_{\HH}(\tau -t) } \| (\mu_{\HH} + \HH) \vv(t)\|_2 d t  \\
  & \les \| \vv\|_Y \Big( \int_0^{\tau}  e^{- \frac78 \mu_{\HH}(\tau -t)} e^{- \frac{3}{4} \mu_{\HH} t}
  + \int_{\tau}^{\infty} e^{ -\frac{1}{4} \mu_{\HH}(\tau -t) } e^{- \frac{3}{4} \mu_{\HH} t}  d t  \Big)
\les e^{- \frac{3}{4}\mu_{\HH} \tau} \| \vv\|_Y .
\end{aligned}
\end{equation}

\paragraph{\bf Step 2: Estimate of $V_2, V_1$}

Suppose that $\| \hat V_1\|_Y < 2 \del$. Using \eqref{eq:ODE_semi2} and $V_2 = \cA_2(\hat V_1)$, we obtain 
\beq\label{eq:solve_ODE_V2}
  \| V_2(\tau) \|_2 \les e^{- \frac{3}{4}\mu_{\HH} \tau} | \hat V_1\|_Y . 
\eeq

Since $\cN(\cdot, \cdot)$ in \eqref{eq:ODE_global} is bilinear and $V = V_1 + V_2$, 
using the above estimate, we obtain 
\[
  \| \cN( V, V )(\tau) \|_2 \les  \| V(\tau)\|_2^2 
  \les \|V_1(\tau) \|_2^2  + \| V_2(\tau)\|_2^2 
\les \|V_1(\tau) \|_2^2  + e^{-3/2 \mu_{\HH}\tau } \| \hat V_1\|_Y^2 .
\]

Since $\cN(V, V)$ is Lipschitz in $V_1$, the $V_1$-ODE \eqref{eq:ODE_split} admits local-in-time solution. Multiplying $V_1$ to the $V_1$-equation and using the above estimate, we obtain 
\begin{equation}\label{eq:ODE_fix_EE1}
  \tfrac{d}{d \tau} \|V_1\|_2^2 \leq - 2 \mu_{\HH} \|V_1\|_2^2  + C
  \|V_1\|_2 ( \|V_1 \|_2^2  + e^{- \frac{3}{2} \mu_{\HH} \tau}  \| \hat V_1\|_Y^2).
\end{equation}

Using $\| \hat V_1\|_Y < 2 \del$ and multiplying both sides by $e^{3/2 \mu_{\HH} \tau}$,  we estimate
$v(\tau) = e^{ \frac{3 \tau }{4} \mu_H}  \|V_1(\tau)\|_2 $:
\[
  \tfrac{d}{d \tau}   v(\tau)^2 
  \leq - \frac{1}{2}\mu_{\HH} v(\tau)^2 
  + C_1 e^{ - \frac{3}{4} \mu_{\HH} \tau}  v(\tau) ( v(\tau)^2  + \del^2   ) ,
\]
for some absolute constant $C_1$. Since $\|V_1(0)\|_2 = \del$ \eqref{eq:ODE_size_init}, by choosing $\del < \frac{\mu_{\HH}}{5 C_1}$ so that 
\[
  - \tfrac{1}{2} \mu_{\HH} (2 \del)^2 + C_1 ( 2 \del) \cdot ( 4 \del^2 + \del^2) 
  = 
    - 2  \mu_{\HH}  \del^2 + 10 C_1  \del^3
   < 0,
\]
and using a bootstrap argument, we prove  $v(\tau) < 2 \del$ for all $\tau \geq 0$. 
Thus, we obtain a global in time solution $V_1$ with $\|V_1\|_Y < 2 \del$, and the map 
$\cA$ maps $B_Y(2 \del)$ into itself.

\paragraph{\bf Step 3: Contraction}

Given any two different inputs $\hat V_{1, a}, \hat V_{2, a} \in Y$. Denote 
\begin{equation}\label{eq:ODE_solve_nota}
 \begin{gathered} 
V_{2, \al} = \cA_2(\hat V_{1, \al}),
\quad V_{1,\al} = \cA( \hat V_{1,\al} ), 
\quad V_{\al} = V_{1, \al } + V_{2,\al}, 
\ \al \in \{a, b\},  \\
 V_{1, \Del} = V_{1, a} - V_{1, b},
\quad V_{2, \Del} = V_{2, a} - V_{2, b},
 \end{gathered} 
\end{equation}

From Step 2, for  $ \al \in \{a, b\}$,  we obtain 
\begin{equation}\label{eq:ODE_fix_est1}
   \| V_{1,\al}(\tau) \|_2 <  2 e^{-\frac34 \mu_{\HH} \tau} \del,
   \quad  \| V_{2,\al}(\tau) \|_2 \les e^{- \frac34 \mu_{\HH} \tau } \del,
   \quad \| V_{\al} (\tau) \|_2 \les e^{-\frac34 \mu_{\HH} \tau } \del .
\end{equation}

Since $\cA_2$ is linear, applying \eqref{eq:ODE_semi2}, we obtain 
\begin{equation}\label{eq:ODE_fix_est2}
\| V_{2, \Del}(\tau) \|_2 = \| \cA_2( \hat V_{1, \Del} )(\tau) \|_2 
   \les e^{- \frac34 \mu_{\HH} \tau} \| \hat V_{1, \Del} \|_Y. 
\end{equation}

Since $\cN$ is bilinear, using  \eqref{eq:ODE_fix_est1} and \eqref{eq:ODE_fix_est2}, we bound 
\[
\begin{aligned} 
  \|\cN(V_a, V_a )(\tau) -   \cN(V_b, V_b )(\tau) \|_2
& \les ( \| V_a(\tau)\|_2 + \|V_b(\tau)\|_2 ) \| (V_a - V_b)(\tau) \|_2 \\
&\les  e^{- \frac34 \mu_{\HH} \tau} \del ( \| V_{1, \Del}(\tau)\|_2 + \| V_{2, \Del}(\tau) \|_2) \\
& \les  
e^{-\frac34 \mu_{\HH} \tau} \del ( \| V_{1, \Del}(\tau) \|_2 + e^{-\frac34 \mu_{\HH} \tau} \| \hat V_{1, \Del} \|_Y  ).
\end{aligned}
\]

Using the $V_1$-ODE in \eqref{eq:ODE_split}, we estimate 
\[
\begin{aligned} 
  \tfrac{d}{d \tau} \| V_{1,\Del}\|_2^2 
  &  = - 2 \mu_{\HH} \| V_{1,\Del}\|_2^2  
   + 2 V_{1, \Del} \cdot ( \cN(V_a, V_a ) -   \cN(V_b, V_b ) ) \\
& \leq - 2 \mu_{\HH} \| V_{1,\Del}\|_2^2   
+ C e^{-\frac34 \mu_{\HH} \tau} \del  \| V_{1, \Del}  \|_2
  ( \| V_{1, \Del}\|_2 + e^{- \frac34 \mu_{\HH} \tau} \| \hat V_{1, \Del} \|_Y  ).
\end{aligned} 
 \]

 Since $V_{1,\Del}(0) =0$, applying estimates similar to those in \eqref{eq:ODE_fix_EE1} and by further requiring $\del$ small enough, and solving the above ODE, we prove 
 \[
   \| V_{1,\Del}(\tau) \|_2 \leq \tfrac{1}{2} e^{-\frac34 \mu_{\HH} \tau} \| \hat V_{1, \Del} \|_Y .
 \]
 Form the definitions \eqref{eq:ODE_solve_nota}, \eqref{eq:ODE_fix_pt}, \eqref{eq:ODE_Ynorm}, the above estimate implies that $\cA$ is a contraction. 
Using the Banach fixed point argument, we prove the existence of fixed point solution 
 $V_1$ to  \eqref{eq:ODE_split} 
 with $\| V_1 \|_Y < 2 \delta$. Moreover, from \eqref{eq:ODE_fix_pt}, \eqref{eq:solve_ODE_V2}, 
 $ \|V_1 \|_Y < 2 \delta$,  the solution $V_1$ and $V_2 = \cA_2(V_1)$ satisfy 
 \[
   \| V_2(0) \|_2 \les \del = \| V_{1, 0}\|_2,
   \quad \| V(\tau)\|_2 \les  e^{-\frac34 \mu_{\HH} \tau} \del 
   = e^{-\frac34 \mu_{\HH} \tau}  \| V_{1, 0}\|_2.
 \]
We complete the proof of \eqref{eq:ODE_solve_est} and Theorem \ref{thm:ODE_solve}.
\end{proof}

\subsection{Proof of  Theorem \ref{thm:ODE_stab} }\label{sec:ODE_stab_proof}

We combine the results in the previous subsection to prove Theorem~\ref{thm:ODE_stab}.

Result (i) in Theorem \ref{thm:ODE_stab} follows from~\eqref{eq:modulated:nosym:summary}. 

We have derived the $n$-th order ODE system in Sections~\ref{sec:ODE_0th}--\ref{sec:ODE_high}; these ODEs have the schematic form given in \eqref{eq:ODE_Vn_lift}. 
By extracting the linearly independent mixed derivatives from $\PP_i V_{=i}^\intercal$ for  $i\leq n$, 
we derive the ODEs for $V_{\leq n}$; see also \eqref{eq:scheme_ODE_V_leqn}. Since the nonlinear terms in these ODE systems are quadratic, result (ii) in Theorem \ref{thm:ODE_stab} follows.

For any $n \geq 2\NNN$, by the definitions of $\Sigma_{\mw{uns}, \leq n}$ in \eqref{def:Sigma_uns_n},
and that of $\Sigma_{\sigma_u}$ in \eqref{def:Sigma_sigma} and  \eqref{eq:ODE_spectral}, 
we obtain that the linear space $\Sigma_{\mathsf{uns},\leq n} $ is given by $\Re \Sigma_{\sigma_u}$ in Theorem \ref{thm:ODE_solve}, with $\HH = M_{\leq n}$. 
Thus, the first part of result (iii) in Theorem \ref{thm:ODE_stab} 
is established in Theorem \ref{thm:uns_general}. 
Applying Theorem \ref{thm:ODE_solve} 
with $\HH = M_{\leq n}, V = V_{\leq n}$ to the ODE system \eqref{eq:ODE_stab}, 
and using $\Re \Sigma_{\sigma_u} =\Sigma_{\mathsf{uns},\leq n }$, we establish the nonlinear stability results in the second part of result (iii) in Theorem \ref{thm:ODE_stab}.

Result (iv) in Theorem \ref{thm:ODE_stab} follows from Theorem \ref{thm:uns_specific}
and the nonlinear stability results in Theorem \ref{thm:ODE_solve}.

Result (v) in Theorem \ref{thm:ODE_stab} is established in Corollary \ref{cor:stable}, and this completes the proof of Theorem~\ref{thm:ODE_stab}.

%%%%%%%%%%%%%%%%%%%

\appendix

\section{Sign of the leading coefficient for velocity at infinity}
\label{sec:appendix:sign}
Recall from~\eqref{eq:V:Q:R=infty} that the leading order asymptotic of the stationary profiles $(\bar V,\bar Q)$ as $R\to \infty$ is given by $\bar V(R) = \underline{v_1} R^{-\frac{1}{\cxbar}} + \OO_{R\to \infty}(R^{-\frac{2}{\cxbar}})$, and $\bar Q(R) = \underline{q_1} R^{-\frac{1}{\cxbar}} + \OO_{R\to \infty}(R^{-\frac{2}{\cxbar}})$. We have already shown in  Proposition~\ref{prop:power:series:infinity} that $0<\underline{q_1}<\infty$; however, we were not able to determine the sign of $\underline{v}_1$. While for the purpose of this manuscript this sign of $\underline{v}_1$ is irrelevant, in a future work we will address the problem of continuing the solutions constructed in this manuscript \emph{past the time of implosion}, as a \emph{self-similar explosion}. In this future work, the fact that $\underline{v}_1 < 0$ plays an important role. The purpose of this Appendix is to show that for the ground state (corresponding to $\NNN=1$), we  have $\underline{v}_1 < 0$ (see~Proposition~\ref{prop:sign:v1}). In fact, $\underline v_1<0$ holds not just for the ground state. Our main result is:

\begin{proposition}[\bf Lower bound for $\underline v_1 / \underline q_1$]
\label{prop:sign:v1}
Let $d \in \{1,2,3\}$, $1 < \gamma \leq 2d+1$, and $\NNN \geq 1$. Let $(\bar V, \bar Q)$ be the global-in-$R$ solution from Proposition~\ref{prop:global:profile}, and let $\underline v_1, \underline q_1$ be the leading coefficients of~\eqref{eq:V:Q:R=infty}. Then 
\begin{equation}
\label{eq:v1:q1:enclosure}
- \sqrt{\tfrac{2\alpha}{d\gamma}} \leq \tfrac{\underline v_1}{\underline q_1} .
\end{equation}
Moreover, if  either $\frac{d}{2(1+\alpha d)} + 1 <   \cxbar$, or if $\frac{d}{2(1+\alpha d)} + 1 \geq   \cxbar$ and $\NNN=1$,  
then we have $\underline{v}_1 < 0$.
\end{proposition}

\subsection{Setup}
Since $\bar Q(R)$ is strictly monotone decreasing on $(0,\infty)$ (cf.~Corollary~\ref{cor:Omega:invariant}), the map $\bar Q \colon [0,\infty) \to (0,\bar q_0]$ is a bijection. We may thus define a map $\cV \colon (0,\bar q_0] \to \Reals$ by $\cV(q) = \bar V(\bar Q^{-1}(q))$, or equivalently, via the implicit definition  
\[
\bar V(R) = \cV(\bar Q(R)).
\] 
By construction, $\cV(\bar q_0) = \bar v_0$, and taking limits as $R\to \infty$ we may define $\cV(0) = 0$, so that $\mathcal{V}$ extends as a continuous function to all of $[0,\bar q_0]$. By definition and~\eqref{eq:euler7} we have
\begin{equation}
\frac{d \cV}{dq}
= \frac{d \bar V}{dR} \Bigr|_{R = \bar Q^{-1}(q)} \frac{d \bar Q^{-1}}{dq}
= \frac{\frac{d \bar V}{dR}}{\frac{d \bar Q}{dR}} \Bigr|_{R = \bar Q^{-1}(q)}
= \frac{\Delta_{\bar V}[\cV(q),q]}{\Delta_{\bar Q}[\cV(q),q]} 
\label{eq:mathcal:V:q:evo}
\end{equation}
for all $q\in(0,q_0)$. 

It is convenient to denote 
\[
\lambda_0 := \tfrac{\bar v_0}{\bar q_0}  = - \sqrt{\tfrac{2\alpha}{d\gamma}} < 0, 
\]
and recall the identities $(1+\alpha d) \bar v_0  = - 1$ and  $\bar v_0 + \bar v_0^2 + \frac{2\alpha^2}{\gamma} \bar q_0^2 = 0$.

\subsection{Lower bound}
\begin{lemma}[\bf Lower barrier]
\label{lem:lower:barrier}
We have 
$\cV(q) > \lambda_0 q$ for all $q \in (0,\bar q_0)$.
\end{lemma}
\begin{proof}[Proof of Lemma~\ref{lem:lower:barrier}]
Define $\phi(R) := \bar V(R) - \lambda_0 \bar Q(R)$. Then by~\eqref{eq:euler7} we have $R \p_R \phi = (F/\Delta)[\bar V ,\bar Q ]$, where $F:= \Delta_{\bar V} - \lambda_0  \Delta_{\bar Q}$. 

From~\eqref{cor:power:series:V:Q} we deduce that $\phi(0)  = 0$, and more precisely, that $\phi(R) = (\bar v_\NNN - \lambda_0 \bar q_\NNN) R^{2\NNN} + \OO(R^{4\NNN})$ as $R\to 0^+$. 
From~\eqref{eq:first:nontrivial:Taylor:coefficient}, \eqref{eq:bar:V:zero},   \eqref{eq:cb:def},~\eqref{eq:cb:admissible} we have
\begin{equation}
\bar v_\NNN - \lambda_0 \bar q_\NNN
= 1 - \tfrac{1+\alpha d + 2 \alpha \NNN }{2\NNN (1+\alpha d) \cbbar} 
= \tfrac{1}{ (1+\alpha d)\cbbar} \bigl( \Sigma_\NNN - \tfrac{1+3\alpha d}{4\NNN}   -   \alpha   \bigr) ,
\label{eq:fancy:app:1}
\end{equation}
where
\[
\Sigma_\NNN := \sqrt{ \tfrac{\alpha \gamma d}{2} + \mathsf{E}_\NNN  + \tfrac{(1-\alpha d)^2}{16 \NNN^2}  }
.
\]
Thus, appealing also to~\eqref{eq:En:def}, we have 
\begin{align*}
{\sf RHS}_{\eqref{eq:fancy:app:1}} >0
&\Leftrightarrow
 \Sigma_\NNN> \tfrac{1+3\alpha d}{4\NNN}   +   \alpha 
  \\
 &\Leftrightarrow
\tfrac{\alpha \gamma d}{2} +  \tfrac{\alpha \gamma d (d+2)}{4\NNN} 
+ \tfrac{\alpha d (1+ \alpha d  ) }{2 \NNN^2}     + \tfrac{(1-\alpha d)^2}{16 \NNN^2} > \tfrac{(1+3\alpha d)^2}{16\NNN^2} + \tfrac{\alpha(1+3\alpha d)}{2\NNN}  +  \alpha^2 
\\
 &\Leftrightarrow
 \alpha^2 (d-1) + \tfrac{\alpha d}{2} +  \tfrac{\alpha (  d (d+2)-  2 )}{4\NNN} + \tfrac{\alpha^2 d (d-1)}{2\NNN}   >  0,
\end{align*}
which is true.  We have thus proven that $\bar v_\NNN - \lambda_0 \bar q_\NNN > 0$, and hence $\phi(R) > 0$ for all sufficiently small $R>0$.

Assume by contradiction that there exists a smallest/first $R_*>0$ such that $\phi(R_*)=0$. In particular, $\phi > 0$ on $(0,R_*)$, and  
\begin{equation}
0 \geq (R \p_R \phi)|_{R_*} = (F/\Delta)[\lambda_0 \bar Q(R_*),\bar Q(R_*)].
\label{eq:fancy:app:2}
\end{equation}
If we are able to show that $(F/\Delta)[\lambda_0 q,q]>0 \Leftrightarrow F[\lambda_0 q,q] >0$ for all $q\in (0,\bar q_0)$ (cf.~\eqref{eq:Omega:invariant:c}), then we would contradict~\eqref{eq:fancy:app:2}, and hence no such $R_*$ exists, concluding the proof. 

Let $s := 1- q/\bar q_0 \in (0,1)$ and $A := \cxbar + \lambda_0 q = \cbbar + s |\bar v_0|  > 0$. Evaluating~\eqref{eq:Delta:V:alt}--\eqref{eq:Delta:Q:alt} at $\bar V = \lambda_0 q$, $\bar Q = q$,   using that $(1+\alpha d)|\bar v_0| = 1$ and  $\bar v_0 + \bar v_0^2 + \frac{2\alpha^2}{\gamma} \bar q_0^2 = 0$, after a tedious computation we  obtain
\begin{align}
\label{eq:F:factored}
F[\lambda_0 q, q] 
&=
(\Delta_{\bar V} - \lambda_0  \Delta_{\bar Q})[\lambda_0 q, q] 
\notag\\
&= 
 \alpha  |\bar v_0| s(1-s)^2 \Bigl(
 \tfrac{ d^2  + 2(d-1)  (1+ \alpha d)   }{2(1+\alpha d)^2} A 
+ \tfrac{d}{2(1+\alpha d)^3 }   (1-s)   
\Bigr),
\end{align}
which is clearly strictly positive. 
\end{proof}

\begin{corollary}
\label{cor:lower:barrier}
The lower bound in~\eqref{eq:v1:q1:enclosure} holds.
\end{corollary}
\begin{proof}[Proof of Corollary~\ref{cor:lower:barrier}]
Dividing $\cV(q) > \lambda_0 q$ by $q$ and sending $q\to 0^+$, we obtain with the help of~\eqref{eq:V:Q:R=infty} that  $\underline v_1/\underline q_1 \geq \lambda_0$.
The proof now follows by the definition of $\lambda_0$.
\end{proof}

\subsection{Upper bound}
We first note that from the power series expansion~\eqref{eq:V:Q:R=infty} as $R\to \infty$, we have 
\[
\mathcal{V}(q) 
= \tfrac{\underline{v}_1}{\underline{q}_1} q  + \tfrac{\underline{v}_2 \underline{q}_1 - \underline{v}_1 \underline{q}_2}{\underline{q}_1^3} q^2 + \OO_{q\to 0^+} (q^3). 
\]
By further appealing to the recursion relations in~\eqref{eq:recursions:at:infinity} for $n=2$, we may compute $\underline{v}_2 = \frac{\cxbar-1}{\cxbar} \underline{v}_1^2 - \beta \underline{q}_1^2$, where
\begin{equation}
\beta:= \tfrac{\alpha^2}{\gamma \cxbar} \bigl( \tfrac{d}{1+\alpha d} + 2 - 2 \cxbar \bigr),
\label{eq:app:beta:def}
\end{equation}
and $\underline{q}_2 = \frac{\cxbar(1+\alpha d)-1-\alpha}{\cxbar} \underline{q}_1\underline{v}_1$. Together with the previous expression for $\mathcal{V}$, we obtain 
\begin{equation}
\mathcal{V}(q) 
= \tfrac{\underline{v}_1}{\underline{q}_1} q  - \Bigl( \tfrac{\alpha(d \cxbar -1)}{\cxbar} \tfrac{\underline{v}_1^2}{\underline{q}_1^2} + \beta\Bigr) q^2 + \OO_{q\to 0^+} (q^3). 
\label{eq:fancy:app:3}
\end{equation}

The parameter $\beta$ plays one more distinguished role; to see this, note that 
\[
\Delta_{\bar V}[0,\bar Q] = 
\tfrac{\alpha^2(2+\gamma d)}{\gamma(1+\alpha d)}  \cxbar  \bar Q^2  - \cxbar ^2    \tfrac{2\alpha^2}{\gamma} \bar Q^2  
= \cxbar \bar Q^2 \tfrac{\alpha^2}{\gamma}  (\tfrac{d}{ 1+\alpha d }  +2   - 2 \cxbar  )
=\cxbar^2 \bar Q^2  \beta.
\]
Since \eqref{eq:Omega:invariant:a} implies 
\[
-\cxbar^2\leq \tfrac{1}{\bar Q}\Delta_{\bar Q}[0,\bar Q]
\leq - \tfrac{\alpha^2 d}{(1+\alpha d)^3}  ,
\]
it follows that when $\beta < 0$ we have
\[
(-\beta) \bar Q 
\leq 
\tfrac{\Delta_{\bar V}[0,\bar Q]}{\Delta_{\bar Q}[0,\bar Q]} 
\leq 
(-\beta) \tfrac{\cxbar^2 (1+\alpha d)^3 }{\alpha^2 d} \bar Q.
\]
Thus, in this $\beta<0$, if there exists a largest value $q_* \in (0,\bar q_0]$ such that $\cV(q_*)= 0$, then due to \eqref{eq:mathcal:V:q:evo} and the above inequalities we have that 
\[
0 < (-\beta) q_*
\leq  \tfrac{d \cV}{dq}\bigl|_{q=q_*}
\leq (-\beta) \tfrac{\cxbar^2 (1+\alpha d)^3 }{\alpha^2 d} q_* 
\]
contradicting the maximality of $q_*$. We have thus proven:
\begin{lemma}
\label{lem:upper:barrier:easy}
Assume that the parameter $\beta$ defined in~\eqref{eq:app:beta:def} is strictly negative. Then, we have $\cV(q)<0$ for all $q\in (0,\bar q_0]$. In particular, $\bar v_1 < 0 $.
\end{lemma}
\begin{proof}[Proof of Lemma~\ref{lem:upper:barrier:easy}]
We have proven $\cV(q)<0$, and thus $\cV(q)/q <0$ for all $q\in (0,\bar q_0]$. Passing $q\to 0^+$, we deduce from~\eqref{eq:fancy:app:3} that $\bar v_1\leq 0$. If the extreme case $\bar v_1=0$ were to occur, then~\eqref{eq:fancy:app:3} implies $\cV(q) = - \beta q^2 +   \OO_{q\to 0^+} (q^3)$; hence, if $\beta<0$ we must have $\cV(q)>0$ for $q\ll 1$, a contradiction. 
\end{proof}

Therefore, we are only left to analyze the situation in which 
\begin{equation}
\beta \geq 0 \quad \Leftrightarrow \quad \cxbar \leq \tfrac{d}{2(1+\alpha d)} + 1 .
\label{eq:fancy:app:4} 
\end{equation}
Recall the notation 
\[
s = 1 -\tfrac{q}{\bar q_0}
\]
used in the proof of Lemma~\ref{lem:lower:barrier}.
Our goal is find a function $F \colon [0,1] \to [0,1]$ such that 
\begin{equation}
\cW(q) := \bar v_0 (1 - F(s))
\label{eq:upper:barrier:definition}
\end{equation}
serves as an upper barrier for the function $\cV$; that is, we must ensure: 
\begin{enumerate}
\item Endpoint matching: $F(0) = 0$ and $F(1)=1$. This condition ensures $\cW(0) = 0 = \cV(0)$ and $\cW(\bar q_0) = \bar v_0 = \cV(\bar q_0)$.
\item Negativity: $0 < F(s) < 1$ for all $s\in (0,1)$. This condition ensures $\cW(q) < 0$ for all $q\in(0,\bar q_0)$.
\item Barrier condition: we must ensure that for all $s\in (0,1)$ we have
\begin{equation}
\Bigl(\tfrac{\Delta_{\bar V}}{\Delta_{\bar Q}}\Bigr)[\bar v_0(1-F(s)),\bar q_0(1-s)] \geq \tfrac{\bar v_0}{\bar q_0} F'(s) .
\label{eq:upper:barrier:condition}
\end{equation}
Equation~\eqref{eq:upper:barrier:condition} is a restatement of the inequality $\tfrac{\Delta_{\bar V}}{\Delta_{\bar Q}}[\cW(q),q] \geq \cW'(q)$. In light of~\eqref{eq:mathcal:V:q:evo}, which gives $\tfrac{\Delta_{\bar V}}{\Delta_{\bar Q}}[\cV(q),q] = \cV'(q)$, if we are able to design $F$ such that~\eqref{eq:upper:barrier:condition} holds, it directly follows that $\cW(q) \geq \cV(q)$ for all $q\in [0,\bar q_0]$, which is our goal.
\end{enumerate}

\begin{remark}
We note that the function $F$ must satisfies an priori  linear lower bound. Since we already know $\cV(q) > \bar v_0 q/\bar q_0$, the desired inequality $\cW(q) \geq \cV(q)$, which is equivalent to  $F(s)  > 1- \frac{1}{\bar v_0} \cV\bigl(\bar q_0 (1-s)\bigr)$, implies that $F(s) > s$.
\end{remark}

Our next result identifies a function $F$ which satisfies the above-listed assumptions (i)--(iii), for the \emph{ground state profiles} corresponding to $\NNN=1$, in the case that is complementary to Lemma~\ref{lem:upper:barrier:easy}.
\begin{proposition}
\label{prop:ground:state:good}
Let $\NNN=1$, $d\in \{1,2,3\}$, and $\alpha \in (0,d]$ be such that $\cxbar = \cxstar(d,2\alpha+1,1)$ satisfies inequality~\eqref{eq:fancy:app:4}. Then, the function
\[
F(s) = 2 s - s^2
\]
satisfies $F(0)=0$, $F(1) =1$, $0<F<1$ on $(0,1)$, and inequality~\eqref{eq:upper:barrier:condition} holds.
\end{proposition}

\begin{remark}
From \eqref{eq:cxstar:ground:state}, it follows that condition~\eqref{eq:fancy:app:4} is satisfied if and only if
\begin{align}
&\cxstar(d,2\alpha+1,1) \leq \tfrac{d}{2(1+\alpha d)} + 1 
\notag\\
\Leftrightarrow \quad &
 \tfrac{1}{1+\alpha d} \Bigl(1  +  \sqrt{   \tfrac{4 (\alpha d)^2 +  (2 +  d + 4 \gamma) \alpha d  }{4} + \tfrac{(1-\alpha d)^2}{16} } + \tfrac{1-\alpha d}{4}   \Bigr)  
 \leq \tfrac{d}{2(1+\alpha d)} + 1 
 \notag\\
\Leftrightarrow \quad &
\sqrt{   \tfrac{4 (\alpha d)^2 +  (2 +  d + 4 \gamma) \alpha d  }{4} + \tfrac{(1-\alpha d)^2}{16} }  
 \leq \tfrac{d}{2} +\alpha d -  \tfrac{1-\alpha d}{4} 
  \notag\\
\Leftrightarrow \quad &
16 (\alpha d)^2 +  4 (2 +  d + 4 \gamma) \alpha d    +  (1-\alpha d)^2  
 \leq (2 d -   1 + 5 \alpha d)^2
   \notag\\
\Leftrightarrow \quad &
 (2 +  d + 4 \gamma) \alpha d    
 \leq   d^2 -   d  +   (5 d - 2) \alpha d + 2 (\alpha d)^2 
    \notag\\
\Leftrightarrow \quad &  
 2(4-d) \alpha^2 + 4 (2-d) \alpha - (d- 1)    \leq 0.
 \end{align}
This last inequality is never satisfied when $d=1$, and hence we must only consider Proposition~\ref{prop:ground:state:good} for dimensions $2$ and $3$. 
\end{remark}

By combining Proposition~\ref{prop:ground:state:good} with \eqref{eq:fancy:app:3}, we obtain that not just $\underline{v}_1\leq 0$; the strict inequality holds. 

\begin{corollary}
\label{cor:ground:state:good}
When $\NNN=1$, $d\in \{1,2,3\}$, and $\alpha \in (0,d]$, we have that  $\cV(q) \leq \cW(q) = \bar v_0 (1-s)^2$, for all $q\in [0,\bar q_0]$ (or equivalently, $s\in[0,1]$). In particular, $\underline v_1 < 0$.
\end{corollary}
\begin{proof}[Proof of Corollary~\ref{cor:ground:state:good}]
The fact that $\cV(q) \leq \cW(q) = \bar v_0 (1-s)^2$ implies that $\underline{v}_1\leq0$. If $\underline{v}_1 =0$, then \eqref{eq:fancy:app:3} would imply that $-\beta q^2 \leq \cW(q)$ as $q\to 0^+$, which is equivalent to $-\beta \bar q_0^2(1-s)^2 \leq \bar v_0 (1-s)^2$, which is equivalent to $\beta \geq \frac{|\bar v_0|}{\bar q_0^2}$. Recalling~\eqref{eq:app:beta:def}, this condition is equivalent to $\tfrac{\alpha d}{1+\alpha d} \bigl( \tfrac{d}{1+\alpha d} + 2 - 2 \cxbar \bigr) \geq  2   \cxbar$, which is equivalent to $\tfrac{(2 +d + 2\alpha d)\alpha d}{2(1+2\alpha d)}   \geq(1+\alpha d)   \cxbar  $. To arrive at a contradiction,  we appeal to \eqref{eq:cxstar:ground:state}, and obtain that 
\begin{align*}
\beta \geq \tfrac{|\bar v_0|}{\bar q_0^2}
\quad \Leftrightarrow \quad 
&\tfrac{(2 +d + 2\alpha d)\alpha d}{2(1+2\alpha d)}   \geq 1  +  \sqrt{   \tfrac{4 (\alpha d)^2 +  (2 +  d + 4 \gamma) \alpha d  }{4} + \tfrac{(1-\alpha d)^2}{16} } + \tfrac{1-\alpha d}{4} 
\\
\quad \Leftrightarrow \quad 
&\tfrac{-5 - 5 \alpha d +2 \alpha d^2 + 6 (\alpha d)^2 }{ 1+2\alpha d }     \geq    \sqrt{ 16 (\alpha d)^2 +  4 (6 +  d + 8 \alpha) \alpha d  +  (1-\alpha d)^2 }  .
\end{align*}
The last inequality fails because the right side (involving the square root) is strictly larger than $4 \alpha d$, and we may readily verify that $-5 - 5 \alpha d +2 \alpha d^2 + 6 (\alpha d)^2 < 4\alpha d(1+2\alpha d)$. This contradicts the assumption that $\underline{v}_1 =0$; thus, $\underline{v}_1 < 0$, as desired.
\end{proof}

We conclude this appendix with:
\begin{proof}[Proof of Proposition~\ref{prop:ground:state:good}]
We use \eqref{eq:Delta:alt} to express ${\sf LHS}_{\eqref{eq:upper:barrier:condition}}$ in terms of $F$.
First, we note that upon replacing $[\bar V,\bar Q]\mapsto [\bar v_0(1-F(s)),\bar q_0(1-s)]$, we have
\begin{align*}
 \cxbar +\bar V &= \cxbar + \bar v_0 (1-F(s)) = \cbbar  + |\bar v_0| F(s) =: A(s) >0, \\
  \bar V - \bar v_0 &= |\bar v_0|  F(s) = A(s) - \cbbar , \\
   \bar V + \bar V^2 + \tfrac{2\alpha^2}{\gamma} \bar Q^2 
   &= \bar v_0(1-F(s)) + \bar v_0^2 (1-F(s))^2+ \tfrac{2\alpha^2}{\gamma} \bar q_0^2 (1-s)^2 \\
   &= \bar v_0(2s - s^2 - F(s) ) + \bar v_0^2 (s- F(s)) (2 - s- F(s)) \\
   &= - \bar v_0^2 s (1-s)^2 (2-s)
\end{align*}
where in the last equality we have used the explicit representation $F(s) = 2s - s^2$.
Inserting these expressions into \eqref{eq:Delta:alt}, we obtain
\begin{align*}
\Delta_{\bar V}[\bar v_0(1-F(s)),\bar q_0(1-s)]
=   \bar v_0^2 s (1-s)^2 (2-s) A(s) \Bigl( \tfrac{\alpha (2+\gamma d)d}{2}  |\bar v_0|
+A(s)   \Bigr),
\end{align*}
and
\begin{align*}
\Delta_{\bar Q}[\bar v_0(1-F(s)),\bar q_0(1-s)]
=
\bar q_0 \bar v_0^2 s  (1-s)^3 (2-s) \Bigl( 
\tfrac{\alpha  d}{2}    |\bar v_0| 
 - \alpha  A(s)   \Bigr)
- \bar q_0  s (1-s) (2-s) A(s)^2  
.
\end{align*} 
Since ${\sf RHS}_{\eqref{eq:upper:barrier:condition}} = \tfrac{\bar v_0}{\bar q_0} 2 (1-s) $, $0<s<1$, $\bar v_0 < 0$, and  cf.~\eqref{eq:Omega:invariant:a} we have that $\Delta_{\bar Q} \leq 0$, we obtain
\begin{align}
\eqref{eq:upper:barrier:condition}
\quad
\Leftrightarrow
\quad
& - |\bar v_0| A(s) \Bigl( \tfrac{\alpha (2+\gamma d)d}{2}  |\bar v_0| + A(s)   \Bigr) 
+ 2 \alpha \bar v_0^2    (1-s)^2  \Bigl(A(s)  - \tfrac{d}{2}    |\bar v_0|  \Bigr) 
+ 2  A(s)^2
\geq 0
\notag\\
\Leftrightarrow
\quad
&- \tfrac{1}{2(1+\alpha d)}  \hat{A}(s) \Bigl( \tfrac{\alpha (2+\gamma d)d}{2}  + \hat{A}(s)   \Bigr) 
+ \tfrac{\alpha}{1+\alpha d}  (1-s)^2  \Bigl(\hat{A}(s)  - \tfrac{d}{2} \Bigr) 
+  \hat{A}(s)^2
\geq 0
\label{eq:upper:barrier:condition:2},
\end{align}
where we have denoted 
\[
\hat{A}(s) = \tfrac{A(s)}{|\bar v_0|} = \tfrac{\cbbar}{|\bar v_0|} + s(2-s)
= \bar \mu - (1-s)^2,
\qquad 
\bar \mu:= \tfrac{\cbbar}{|\bar v_0|} + 1 = \tfrac{\cxbar}{|\bar v_0|} > 1.
\]
We expand $\hat{A}$ in the last inequality in~\eqref{eq:upper:barrier:condition:2}
and deduce (upon multiplying by $1+\alpha d>0$) that 
\begin{align}
\eqref{eq:upper:barrier:condition}
\quad
\Leftrightarrow
\quad
&
\mu_0
- \mu_1 (1-s)^2   
+ \mu_2  (1-s)^4 
\geq 0
\label{eq:upper:barrier:condition:3},
\end{align}
where
\begin{align*}
 \mu_0  =  \tfrac{1 + 2\alpha d}{2}  \bar\mu^2
- \tfrac{ (2+\gamma d) \alpha d}{4}  \bar \mu, 
\quad
 \mu_1  = (1 -\alpha  + 2\alpha d) \bar\mu  
 - \tfrac{\alpha  \gamma d^2}{4},
 \quad 
 \mu_2  =\tfrac{1 + 2\alpha (d-1)}{2} .
\end{align*}
Next, we note that 
\begin{align*}
\bar \mu -1 = \tfrac{\cxbar}{|\bar v_0|} - 1
&= \sqrt{   \tfrac{4 (\alpha d)^2 +  (2 +  d + 4 \gamma) \alpha d  }{4} + \tfrac{(1-\alpha d)^2}{16} } + \tfrac{1-\alpha d}{4}
> \tfrac{4+ 7\alpha d}{8}
\notag\\
\Leftrightarrow
&\quad  \tfrac{4 (\alpha d)^2 +  (2 +  d + 4 \gamma) \alpha d  }{4} + \tfrac{(1-\alpha d)^2}{16}   
> \tfrac{(1+\frac 92\alpha d)^2}{16}  
\notag\\
\Leftrightarrow
&\quad  16 (\alpha d)^2 + 4 (6 + d + 8 \alpha) \alpha d    
> 11\alpha d + \tfrac{77}{4}(\alpha d)^2
\notag\\
\Leftrightarrow
&\quad   13 + 4  d + 32 \alpha
>   \tfrac{13}{4} \alpha d 
, \qquad \mbox{which is true since } d \in \{2,3\}, \alpha \in (0,d].
\end{align*}
Thus,
\begin{equation}
\label{eq:app:bar:mu:lower:bound}
\bar \mu > \tfrac{12 + 7\alpha d}{8} > \tfrac{3+\alpha d}{2}.
\end{equation}
Using this information, we may bound
\[
\mu_0
= \tfrac{\bar\mu (2(1 + 2\alpha d) \bar\mu -  (2+\gamma d) \alpha d   ) }{4}  
> \tfrac{\bar\mu ( (1 + 2\alpha d) (3+\alpha d) -  (2+\gamma d) \alpha d   )
}{4}  
= \tfrac{\bar\mu (3 +  \alpha d (5-d))}{4}  
\geq \tfrac{\bar\mu (3 + 2  \alpha d)}{4}  
> \tfrac{(3+\alpha d)(3+2\alpha d)}{8}.
\]
Thus, \eqref{eq:upper:barrier:condition:3} automatically holds true when $s$ is close to $1$. 
Next, we show that $\mu_1 > 2 \mu_2$, so that the minimum of the expression \eqref{eq:upper:barrier:condition:3} is attained at $s=0$; indeed, we have
\[
\mu_1 - 2 \mu_2
=
(1 -\alpha  + 2\alpha d) \bar\mu  
 - \tfrac{\alpha  \gamma d^2}{4}  -(1 + 2\alpha (d-1))
 >
\tfrac{(1 -\alpha  + 2\alpha d) (3+\alpha d)}{2} 
- \tfrac{4 + 8\alpha (d-1) + \alpha  \gamma d^2}{4} 
>
\tfrac{
1+   \alpha  +  \alpha d +  \alpha^2 d }{2}.
\]
Thus, we have reduced the problem of verifying~\eqref{eq:upper:barrier:condition:3}, meaning that $\mu_0 - \mu_1 (1-s)^2 + \mu_2  (1-s)^4  \geq 0$ for all $s\in[0,1]$, to the problem of verifying that $\mu_0 - \mu_1 + \mu_2 \geq 0$, which is the value at $s=0$. We verify this fact explicitly; first,
\begin{align*}
\mu_0 - \mu_1 + \mu_2
&= \tfrac{1 + 2\alpha d}{2}  \bar\mu^2
- \tfrac{   \alpha (d^2-4) + 2 (1+\alpha d)^2 + 2 (1+3 \alpha d)}{4}  \bar \mu
+ \tfrac{\alpha (d^2-4) + 2 (1+\alpha d)^2 }{4}
=: g(\bar \mu).
\end{align*}
Since 
\begin{align*}
\tfrac{   \alpha (d^2-4) + 2 (1+\alpha d)^2 + 2 (1+3 \alpha d)}{4(1 + 2\alpha d)}
< \tfrac{3+\alpha d}{2}
\Leftrightarrow
 \alpha d^2
< 2 + 4\alpha +  4\alpha d + 2(\alpha d)^2,
\quad \mbox{which is true},
\end{align*}
and since $\bar \mu > \frac{3+\alpha d}{2}$, it follows that $g(\bar \mu)$ is a strictly increasing function of $\bar \mu$, in the range of $\bar \mu$ we are interested in.
Finally, in light of~\eqref{eq:app:bar:mu:lower:bound}, we verify
\begin{align*}
 g(\tfrac{12 + 7\alpha d}{8})
&=     \tfrac{1 + 2\alpha d}{2} \tfrac{(12 + 7\alpha d)^2}{64}
- \tfrac{   \alpha (d^2-4) + 2 (1+\alpha d)^2 + 2 (1+3 \alpha d)}{4} \tfrac{(12 + 7\alpha d)}{8}
+ \tfrac{\alpha (d^2-4) + 2 (1+\alpha d)^2 }{4}
\notag\\
&
=
\begin{cases}
 \tfrac{4(1-\alpha) + \alpha^2 (73 + 84\alpha)}{32} ,&d=2, \\
 \tfrac{ (4-13\alpha)^2 +2 \alpha^2 (34 + 567\alpha)}{128},&d=3.
\end{cases}
\end{align*}
Since $\alpha \in(0,d]$, we deduce that $g(\bar \mu) > g(\tfrac{12 + 7\alpha d}{8}) > 0$; therefore \eqref{eq:upper:barrier:condition} is true, which in turn shows that \eqref{eq:upper:barrier:condition:3} is true. This concludes the proof of the Proposition. 
\end{proof}

\section{Functional inequalities}

In this Appendix, we establish an $L^{\infty}$-based interpolation inequality for the stability analysis in radial symmetry in Section \ref{sec:stability}, and $L^2$-based functional inequalities  for the stability analysis outside of radial symmetry in  Section  \ref{sec:nonradial}.

\subsection{An \texorpdfstring{$L^\infty$-based}{L infinity based} interpolation inequality}

We prove an inequality for the $N$th derivative of the product $x f^2(x)$, where $x\in \Omega$ and $f\colon \Omega \to \Reals$, and $\Omega = [0,\infty)$. The main result is:

\begin{theorem}[Weighted interpolation inequality]\label{thm:appendix:interpolate}
Let $N \geq 1$. There exists a constant $C_N > 0$ such that for all $f \in W^{N,\infty}([0,\infty))$ with $\norm{\brak{x}  f}_{L^\infty} < \infty$:
\begin{equation}
\norm{x f^2}_{\dot{W}^{N,\infty}([0,\infty))} \leq C_N \left( \norm{\brak{x}  f}_{L^\infty} \norm{f}_{\dot{W}^{N,\infty}} + \norm{\brak{x}  f}_{L^\infty}^{1+1/N} \norm{f}_{\dot{W}^{N,\infty}}^{1-1/N} \right).
\label{eq:appendix:interpolate}
\end{equation}
A mutatis-mutandi estimate holds when $x f(x)^2$ is replaced by $x \, f(x) \, g(x)$ with $f\neq g$.
\end{theorem}

Throughout this Appendix, we denote the $n$th derivative of a function $f\colon (0,\infty) \to \Reals$ by $f^{(n)}$. Before proving Theorem~\ref{thm:appendix:interpolate}, we  recall a classical interpolation inequality on half-lines.

\begin{lemma}\label{lem:LK}
Let $N \geq 1$, $0 \leq j \leq N$, and $x_0 \geq 0$. For any $f \in W^{N,\infty}([x_0,\infty))$, we have
\begin{equation}
\| f^{(j)}\|_{L^\infty([x_0,\infty))} \leq C_{N,j} 
\|f\|_{L^\infty([x_0,\infty))}^{1-j/N} 
\| f^{(N)}\|_{L^\infty([x_0,\infty))}^{j/N},
\label{eq:LK:ineq}
\end{equation}
where $C_{N,j}$ is a constant depending only on $N$ and $j$; in particular, the constant is independent of $x_0$.
\end{lemma}

\begin{proof}[Proof of Lemma~\ref{lem:LK}]
We proceed by induction on $N$.  
For the base case $N = 1$, we note that the only  cases are $j = 0$ or $j = 1$, for which~\eqref{eq:LK:ineq} holds immediately, with $C_{1,0} = C_{1,1} = 1$.

For the inductive step, assume that~\eqref{eq:LK:ineq} holds at level $N-1$, for all $0\leq j \leq N-1$. We aim to prove~\eqref{eq:LK:ineq} at level $N$. The bound trivially holds when $j=0$ and $j=N$, with $C_{N,0} = C_{N,N} = 1$. Fix $x > x_0$ and $h > 0$. By Taylor's theorem with integral remainder:
\[
f(x+h) = \sum_{k=0}^{N-1} \tfrac{1}{k!} h^k f^{(k)}(x) + \tfrac{1}{(N-1)!} \int_x^{x+h} (x+h-t)^{N-1} f^{(N)}(t) \, dt.
\]
Setting $M_0 := \norm{f}_{L^\infty([x_0,\infty))}$ and $M_N := \| f^{(N)}\|_{L^\infty([x_0,\infty))}$, we have
\[
\abs{f(x+h)} \leq M_0, \qquad \abs{\int_x^{x+h} (x+h-t)^{N-1} f^{(N)}(t) \, dt} \leq M_N \tfrac{1}{N!} h^N.
\]
Rearranging for $f^{(N-1)}(x)$ we obtain
\[
f^{(N-1)}(x) =  (N-1)! h^{1-N} \Bigl( f(x+h) - \sum_{k=0}^{N-2} \tfrac{1}{k!} h^kf^{(k)}(x) - R_N \Bigr),
\]
where $R_N$ quantifies the integral remainder term, and is bounded as $\abs{R_N} \leq M_N h^N / N!$. Using the inductive hypothesis to bound $\abs{f^{(k)}(x)}$ for $k \leq N-2$, we obtain the pointwise bound
\begin{align*}
| f^{(N-1)}(x) | 
&\leq  h^{1-N}  \left( 2 (N-1)!M_0 +  \sum_{k=1}^{N-2} \tfrac{C_{N-1,k}(N-1)!}{k!}  h^k M_0^{1- \frac{k}{N-1}}  \|f^{(N-1)}\|_{L^\infty( [x_0,\infty) )}^{\frac{k}{N-1}}  \right) +  \tfrac{1}{N}h M_N \notag\\
&\leq h^{1-N} \left( C_N^\prime M_0 +  \tfrac{1}{2} h^{N-1} \|f^{(N-1)}\|_{L^\infty( [x_0,\infty) )} \right) +  \tfrac{1}{N} h M_N 
\end{align*}
for some constant $C_N^\prime>0$. Taking the supremum over $x\in [x_0,\infty)$, absorbing the $\frac{1}{2} \|f^{(N-1)}\|_{L^\infty( [x_0,\infty) )}$ term into the left side, and then optimizing in $h$, we obtain the bound~\eqref{eq:LK:ineq} for $j=N-1$, namely
\[
\|f^{(N-1)}\|_{L^\infty( [x_0,\infty) )} \leq C_{N,N-1} M_0^{1/N} M_N^{(N-1)/N}.
\]
For the intermediate derivatives $1\leq j \leq N-2$, we apply the inductive hypothesis, interpolating $f^{(j)}$ between $f$ and $f^{(N-1)}$, and deduce 
\[
\| f^{(j)}\|_{L^\infty([x_0,\infty))} \leq C_{N,j} M_0^{1-j/N} M_N^{j/N}.
\]
This concludes the proof of~\eqref{eq:LK:ineq}.
\end{proof}

Next, we apply Lemma~\ref{lem:LK} locally, to obtain pointwise decay estimates.

\begin{lemma}[Pointwise decay of derivatives]\label{lem:pointwise}
Let $N \geq 1$ and $0 \leq j \leq N$. Suppose $f \in W^{N,\infty}([0,\infty))$ satisfies
\[
M_0 := \| \brak{x} f\|_{L^\infty([0,\infty))} < \infty, \qquad M_N := \|f^{(N)}\|_{L^\infty([0,\infty))} < \infty.
\]
Then, for all $x > 0$ we have
\[
|f^{(j)}(x)| \leq C_{N,j} ( \brak{x}^{-1} M_0)^{1-j/N} M_N^{j/N}.
\]
\end{lemma}

\begin{proof}[Proof of Lemma~\ref{lem:pointwise}]
Fix $x > 0$. We aim to apply Lemma~\ref{lem:LK} on the half-line $[x,\infty)$. First, we note the trivial bounds
\[
\|f\|_{L^\infty([x,\infty))} \leq \brak{x}^{-1} M_0,
\qquad
\| f^{(N)}\|_{L^\infty([x,\infty))} \leq M_N.
\]
Then,~\eqref{eq:LK:ineq} implies
\[
\| f^{(j)}\|_{L^\infty([x,\infty))} 
\leq C_{N,j} 
\|f\|_{L^\infty([x,\infty))}^{1-j/N} 
\| f^{(N)}\|_{L^\infty([x,\infty))}^{j/N}
\leq C_{N,j} (\brak{x}^{-1} M_0)^{1-j/N} M_N^{j/N}.
\]
\end{proof}

\begin{proof}[Proof of Theorem~\ref{thm:appendix:interpolate}]
As in Lemma~\ref{lem:pointwise}, set $M_0 := \norm{\brak{x}  f}_{L^\infty([0,\infty))}$ and $M_N := \norm{f^{(N)}}_{L^\infty([0,\infty))}$.
Using the Leibniz rule we deduce
\begin{equation}\label{eq:expansion}
(xf^2)^{(N)} = \sum_{j=0}^{N} \binom{N}{j} x f^{(j)} f^{(N-j)} + N \sum_{j=0}^{N-1} \binom{N-1}{j} f^{(j)} f^{(N-1-j)}.
\end{equation}
Consider a term from the first sum in~\eqref{eq:expansion}, namely $x f^{(j)} f^{(N-j)}$ with $0 \leq j \leq N$. Using Lemma~\ref{lem:pointwise} we may estimate
\[
|x f^{(j)}(x) f^{(N-j)}(x)| \leq C_{N,j} C_{N,N-j} \bigr( x \brak{x}^{-1} \bigl) M_0 M_N \leq C_{N,j} C_{N,N-j}  M_0 M_N .
\]
Next, consider a term from the second sum in~\eqref{eq:expansion}, namely $ f^{(j)} f^{(N-1-j)}$ with $0 \leq j \leq N-1$. By Lemma~\ref{lem:pointwise}, 
\[
|x f^{(j)}(x) f^{(N-1-j)}(x)| 
\leq C_{N,j} C_{N,N-1-j} x (\brak{x}^{-1} M_0)^{1 + 1/N} M_N^{1-1/N}
\leq C_{N,j} C_{N,N-1-j}  M_0^{1 + 1/N} M_N^{1-1/N}
.
\]
Note that the powers of $M_0$ and $M_N$ are independent of $j$ in both of the above bounds. Inserting these estimates into~\eqref{eq:expansion} concludes the proof of~\eqref{eq:appendix:interpolate}.
\end{proof}

\subsection{\texorpdfstring{$L^2$-based}{L2 based} functional inequalities}

\begin{lemma}
\label{lem:interp}
Let $\{\vp_n\}_{n \geq 0}$ be a sequence of weights satisfying
\begin{equation}\label{ass:interp}
  \vp_n \leq C(A_1, n) (\vp_{n-1} \vp_{n+1})^{1/2}, \quad  |\na \vp_n| \leq C(A_2, n) (\vp_n \vp_{n-1})^{1/2} ,
\end{equation}
for some parameters $A_1, A_2$ independent of $n$, and some constants $C(A_1, n), C(A_2, n)$  depending on  $n$ and the parameters $A_1, A_2$. Denote 
\[
  I_n := \int  |\na^n f|^2 \vp_n d x . 
\]
If~\eqref{ass:interp} holds for all $1\leq n\leq N$, then for any $0 \leq n \leq N$ and $\del> 0$, we have 
\begin{equation}\label{eq:interp_convex}
  I_n \leq C(A_1, A_2, k, \del) I_0 + \del I_{n+1}.
\end{equation}
\end{lemma}
\begin{proof}[Proof of Lemma~\ref{lem:interp}]
The proof follows that of \cite[Lemma C.2]{CCSV2024}.
Throughout the proof, all implicit constants depend on the index $n$ of $I_n$, and on the parameters $A_1$ and $A_2$. 
Using integration by parts, we have 
\begin{align*}
 I_n 
 &= \sum_i \int  \pa_i \na^{n-1} f \cdot \pa_i \na^{n-1} f \vp_n = -\sum_i \int 
 \bigl( \Delta \na^{n-1} f \cdot \na^{n-1} f \vp_n
  +  \pa_i \na^{n-1} f \cdot \na^{n-1} f \pa_i \vp_n \bigr) .
\end{align*}

Using the assumption \eqref{ass:interp} and the Cauchy-Schwarz inequality we obtain
\begin{equation}
\label{eq:interp_convex2}
I_n \les (I_{n+1} I_{n-1})^{1/2} +  (I_n I_{n-1} )^{1/2}, \quad n \geq 1.
\end{equation} 

To prove~\eqref{eq:interp_convex}, it suffices to show that for any $ n \geq 0$ and $\del>0$, there exists $C_{\del,n}>0$ such that 
\begin{equation}\label{eq:interp_convex3}
I_n \leq \del I_{n+1} + C_{\del,n} I_0.
\end{equation}
We prove~\eqref{eq:interp_convex3} by induction on $n$. The base case $n=0$ is trivial, with $C_{\del,0}=1$. 
For the induction step, let $n\geq 1$. We use~\eqref{eq:interp_convex2}, the inductive hypothesis for $n-1$, and the Cauchy-Schwarz inequality, to conclude 
\[
I_n \leq \del I_{n+1} + \tfrac{1}{2} I_n + \tfrac 14 \bigl(2+ \delta^{-1} \bigr)  C_{\eqref{eq:interp_convex2},n}^2 I_{n-1} 
\leq \del I_{n+1} + \tfrac{3}{4}  I_n  + \underbrace{\tfrac 14 C_{\delta',n-1}}_{=:C_{\del,n}} I_0,
\]
where $\delta' := (2+ \delta^{-1})^{-1} C_{\eqref{eq:interp_convex2},n}^{-2}$.
The above estimate concludes the proof of the induction step for~\eqref{eq:interp_convex3}.
\end{proof}

\begin{lemma}\label{lem:embed}
Let $d$ be the dimension. Let $k\geq d$. Suppose that the weight $g \in C^{\infty}$ satisfies $ g >0$  and 
\begin{equation}\label{ass:lem:embed}
 |\na^i g(y) | \leq \mu g(y) \la y \ra^{-i} , \qquad 1 \leq i \leq d,
\end{equation}
for a constant $\mu\geq 1$. Then, for any function $f \colon \Reals^d \to \Reals$, any $a \in \Reals$ and $0 \leq i\leq k - d $,  we have 
\begin{equation}\label{eq:embed}
  |\na^i f(y)| g^{-1}(y)  \leq C(a,d, k, \mu) \la y \ra^{-i - \frac{a+d}{2} } 
  \left(
  \bigl\| f \; \la z \ra^{a/2} g^{-1} \bigr\|_{L^2}
  +
  \bigl\| \na^k f \; \la z \ra^{k+a/2} g^{-1} \bigr\|_{L^2}
  \right)  .
\end{equation}
\end{lemma}

\begin{proof}[Proof of Lemma~\ref{lem:embed}]
Fix $0 \leq i \leq k-d$. 
By a standard density argument we can assume that $\nabla^i f \in C_c^{\infty}(\Reals^d)$.
Consider the cone with vertex at $y$ extending towards infinity: $\Om(y) := \{ z \in \Reals^d \colon z_j \sgn(y_j) \geq |y_j|, \mbox{ for all } 1\leq j\leq d \}$. In particular, for any $z \in \Om(y)$ we have $|z| \geq |y|$. By integrating on rays extending to infinity, we have
\begin{equation}\label{eq:embed_pf1}
 \la y \ra^{ 2 i + a + d } g^{-2} |\na^i f(y) |^2 
 \les_d 
 \int_{\Om(y)} \Bigl| \pa_1 \pa_2 .. \pa_d \bigl( 
 \la z \ra^{ 2 i + a + d } g^{-2}(z) |\na^i f(z) |^2 \bigr)  \Bigr| d z 
\end{equation}
Applying Leibniz's rule, we obtain 
\begin{align*} 
I(z) 
&:= \bigl| \pa_1 \pa_2 .. \pa_d ( 
 \la z \ra^{ 2 i + a + d } g^{-2}(z) |\na^i f(z) |^2 ) \bigr| 
 \\
& \les  \sum_{ j + l + m + n = d }  |\na^j \la z \ra^{2 i + a + d} | 
\cdot |\na^{i+l} f(z) | \cdot |\na^{ i + m } f(z)| \cdot |\na^n g^{-2}(z) |.
\end{align*}
Further applying the Leibniz rule (the Fa\`a di Bruno formula) and using \eqref{ass:lem:embed}, we obtain 
\begin{align*} 
|\na^n g^{-2}(z)| 
& \les_n   \sum_{ m \leq n } g(z)^{-(m+2)} \sum_{\substack{i_1 + i_2 + .. + i_n = m\\ 1\cdot i_1 + 2 \cdot i_2 + \ldots n \cdot i_n = n}} \prod_{k=1}^{n} |\na^{k} g(z) |^{i_k}
\\ 
&\les_{n, \mu}  \sum_{ m \leq n } g(z)^{-(m+2)} 
\sum_{\substack{i_1 + i_2 + .. + i_n = m\\ 1\cdot i_1 + 2 \cdot i_2 + \ldots n \cdot i_n = n}}  \prod_{k=1}^{n} ( g(z) \la z \ra^{-k})^{i_k}
\les_{n, \mu} g(z)^{-2} \la z \ra^{-n}
.
\end{align*}
Since $ d - j - n = l + m$, it follows that
\[
|\na^j \la z \ra^{2 i + a + d} |  \cdot   |\na^n g^{-2}(z)|
\les_{a, d, k, \mu} g(z)^{-2} \la z \ra^{2 i + a + d - j -n}
=_{a, d, k, \mu} g(z)^{-2} \la z \ra^{2 i + a +  l + m}.
\]
Combining the above estimates, we obtain the bound
\[
 I(z) \les_{a, d, k, \mu}
 \sum_{ l + m \leq d} \la z \ra^{i+l + a/2}  
 |\na^{i+ l} f(z) | g(z)^{-1} \cdot 
\la z \ra^{i+m + a/2}   |\na^{ i + m } f(z)|  g(z)^{-1} .
\]

Next, we verify that the sequence of weights $\vp_n = \la x \ra^{2 n + a} g^{-2}$ satisfy the assumption \eqref{ass:interp}  for all $1\leq n \leq k-i$. Clearly, we have $\vp_n = (\vp_{n-1} \vp_{n+1})^{1/2}$ for $n \geq 1$. Using \eqref{ass:lem:embed}, we obtain 
\[
  |\na \vp_n | \les_{n, a} \la x \ra^{2 n + a}  |\na g| g^{-3}
  + \la x \ra^{2 n +a-1} g^{-2}
  \les \la x \ra^{2 n + a - 1} g^{-2}
  \les_{n, a, \mu} (\vp_n \vp_{n-1})^{1/2}.
\]
Applying Cauchy-Schwarz inequality, and Lemma~\ref{lem:interp} with weights $\{\vp_n\}_{n\geq 0}$, we establish 
\begin{equation}\label{eq:embed_pf2}
 \int_{\Om(y)} I(z) d z 
 \les_{a , d ,k ,\mu}
 \bigl\| \bigl( |f| + \la y \ra^{k} |\na^k f| \bigr) \la y \ra^{a/2} g^{-1} \bigr\|_{L^2}^2.
\end{equation}
Combining the estimate \eqref{eq:embed_pf2} to \eqref{eq:embed_pf1}, we  prove \eqref{eq:embed}.
\end{proof}

\section{Derivation of the high order ODEs and the Proof of Proposition \ref{prop:ODE_Nth}}
\label{app:derive_ODE}
In this Appendix, we use \eqref{eq:lin_ODE:cL}--\eqref{eq:lin_ODE} to derive the ODE system satisfied by the  high-order derivatives of the fundamental unknowns at the origin (namely $\na^n \tvr(0,\tau)$, $\na^{n+1} \tu(0,\tau)$, and $\na^{n+2} \tb(0,\tau)$), and prove Proposition \ref{prop:ODE_Nth}. Throughout the appendix, we suppress time dependence, and denote $g(0,\tau)$ simply as $g(0)$.

\subsection{High order identities}\label{sec:high_ODE_iden}
In order to simplify the arguments, we first establish several identities for high order derivatives at $y = 0$.

\begin{lemma}\label{lem:scal_id}
Suppose that $f \in C^{n, \nu}$ with some $\nu>0$ near $y = 0$. For any multi-indices $\bet, \del \in \Naturals_0^d$ with $|\bet| \leq n$, we have 
\begin{subequations}
\begin{align}
  \pa^{\bet}( y \cdot \na f ) |_{y = 0} &= |\bet| ( \pa^{\bet} f)(0) , 
\label{eq:scal_id:a}
   \\
  \pa^{\bet} (y^{\del} 
   (y \cdot \na f )  )|_{y = 0}  & =
(|\bet| - |\del| ) \one_{|\bet| \geq |\del|} \pa^{\bet} (y^{\del}  f  ) |_{y = 0}.
\label{eq:scal_id:b}
\end{align}
\end{subequations}
\end{lemma}

\begin{proof}[Proof of Lemma~\ref{lem:scal_id}]
Performing a Taylor expansion at the origin, we obtain 
$
  f = \sum_{ |\th|\leq k, \th \in \Naturals_0^d} f_{\th} y^{\th} + \cR$, where the error term $\cR$ satisfies $ | \na^i \cR|(0) = 0$ for any $i \leq k$,
and $\th$ is a multi-index. Since $ y \cdot \na y^{\th} = |\th| y^{\th} $, we obtain 
\[
    \pa^{\bet}( y \cdot \na f ) |_{y = 0}
    =  \sum_{ |\th|\leq k } \pa^{\bet} ( |\th| f_{\th} y^{\th} )(0)
    = |\bet| f_{ \bet } \pa^{\bet} (y^{\bet})(0)
    =  |\bet| ( \pa^{\bet} f)(0),
\]
thereby proving \eqref{eq:scal_id:a}. 

In order to prove~\eqref{eq:scal_id:b}, we note that  if $|\bet| < |\del|$, both sides equal $0$. If $|\bet| \geq |\del|$, using \eqref{eq:scal_id:a},
 we have 
\[
\pa^{\bet} (y^{\del} (y \cdot \na f - (|\bet| - |\del|) f))(0)
= \sum_{ |\bet_1| = |\del| , \bet_1 + \bet_2 = \bet }
C_{\bet_1, \bet_2} ( \pa^{\bet_1} y^{\del})(0)  \pa^{\bet_2} (y \cdot \na f - (|\bet| - |\del|) f)(0) 
=0,
\]
which completes the proof.
\end{proof}

The system of $\VVs_{\bet, k}$ to be derived in \eqref{eq:ODE_main}, \eqref{eq:ODE_main_xu}, \eqref{eq:ODE_main_divu} is not closed due to terms such as $|y|^2 \tvr$ in \eqref{eq:ODE_main_xu}.
We use the following  identities to show that it is closed up to some lower order terms. 
We recall the notation $\FFs_{k, l}, \OOL{}(\FFs_{k,l})$ from \eqref{norm:ODE_low}
\[
\begin{aligned}
  \FFs_{k, l} & =  ( \na^k \Del^l \tvr(0) , \,  \na^k \Del^l (\div \tu)(0) , \,
      \na^k \Del^{l+1} ( y \cdot \tu)(0) , \,
  \na^k \Del^{l+1} \tb(0) ), \\
 g & = C_{h} \cdot \FFs_{k, l}  , \qquad  \forall \, g \in   \OOL{h}( \FFs_{k, l})
\end{aligned}
\]

\begin{lemma}\label{lem:lap}
For any multi-index $\bet$ and $k \geq 0$, we have 
\begin{subequations}\label{eq:lap_id}
\begin{align}
  \Del^{k+1} (y \cdot \tu )(0) &= 2(k+1)  \Del^k (\div \tu) (0), 
  \label{eq:lap_id:div}  \\
  \pa^{\bet} \Del^{k+1} ( |y|^2 \tvr )(0)
& = (2 d + 4 |\bet| + 4k ) (k+1) \pa^{\bet} \Del^k \tvr(0) + \one_{ |\bet |\geq 2 } \OOL{n}( \FFs_{|\bet|-2, k+1} ) , 
\label{eq:lap_id:a} \\
  \pa^{\bet} \Del^{k} ( |y|^2 \Del \tvr )(0)
 & = (2 d + 4 |\bet| + 4 (k-1) ) k  \pa^{\bet} \Del^k \tvr(0) + \one_{ |\bet |\geq 2 } \OOL{n}( \FFs_{|\bet|-2, k+1} ) , 
\label{eq:lap_id:b} \\
\Del^{k+1} (|y|^{2k} f ) (0)
& =  \msf{c}_{k, \Del} \Del f (0) ,
\quad \msf{c}_{k, \Del} := (k+1) \prod_{1 \leq i \leq k} 2i (2 i + d) 
 \label{eq:lap_id:c} . 
\end{align}
We also introduce the following constant 
\beq
  \msf c_{k, 1} := \Del^{k} (|y|^{2k}) =  \prod_{1\leq i \leq k} 2i (2i + d - 2).
  \label{eq:lap_id:d} 
\eeq

\end{subequations}
\end{lemma}
 
\begin{proof}[Proof of Lemma~\ref{lem:lap}]
Throughout this proof, all identities are evaluated at $y = 0$.

\vspace{0.1in}
\paragraph{\bf Proof of \eqref{eq:lap_id:div}} 
Using Leibniz rule, for any $i \geq 1$ and $j \geq 0$, we obtain
\begin{equation}\label{eq:lap_id:div:pf}
  \Del^{i}( y \cdot \Del^j \tu )
  = \Del^{i-1} (2 \pa_l y  \cdot \Del^j \pa_l \tu
  + y \cdot \Del^{j+1} \tu)
  = 2 \Del^{i+j-1} \div \tu
  + \Del^{i-1}(    y \cdot \Del^{j+1} \tu) .
\end{equation}

Summing  \eqref{eq:lap_id:div:pf} over $i = k+1, k,.., 1$ with $j = k+ 1 - i$ 
and then evaluating at $y=0$, we obtain 
\[
    \Del^{k+1}( y \cdot  \tu )(0)
    = 2 (k+1) \Del^k \div \tu(0)
    + y \cdot \Del^{k+1} \tu (0)
    = 2 (k+1) \Del^k \div \tu(0).
\]

\vspace{0.1in}

\paragraph{\bf Proof of \eqref{eq:lap_id:a}, \eqref{eq:lap_id:b} }
Using Leibniz rule, for any $i \geq 1$ and $j \geq 0$, we obtain 
\[
  \Del^i (|y|^2 \Del^j \tvr)
= \Del^{i-1}( \Del |y|^2 \Del^j \tvr + 2 \na |y|^2  \cdot \na \Del^j \tvr 
+ |y|^2 \Del^{j + 1} \tvr  )
= \Del^{i-1}  ( 2 d \Del^j \tvr + 4  y  \cdot \na \Del^j \tvr 
+ |y|^2 \Del^{j + 1} \tvr 
  ) .
\]

Using Lemma \ref{lem:scal_id}, $|\bet| + 2k = n$ and $f = \Del^j \tvr$  we obtain 
\[
  \pa^{\bet} \Del^{i-1} ( y \cdot \na f )(0)
  = (|\bet| + 2(i-1)) \pa^{\bet} \Del^{i-1} f(0)
  = (|\bet| + 2(i-1))\pa^{\bet} \Del^{i-1 + j} \tvr (0).
\]

Combining the above two identities, we derive
\begin{equation}\label{eq:lap_recur}
   \pa^{\bet} \Del^i (|y|^2 \Del^j \tvr)(0)
     = (2 d +  4 |\bet| + 8(i-1)) \pa^{\bet} \Del^{i+j-1} \tvr(0)
     +  \pa^{\bet} \Del^{i-1} (|y|^2 \Del^{j+1} \tvr).
\end{equation}

Summing the above identities over $i = k+1, k,.., 1$ with $j = k+ 1 - i$ we obtain 
\[
     \pa^{\bet} \Del^{k+1} (|y|^2  \tvr)(0)
     = \bigl( (2 d + 4 |\bet|) (k+1) + 4 k (k+1)  \bigr)
     \pa^{\bet} \Del^{k} \tvr(0) + \pa^{\bet}( |y|^2 \Del^{k+1} \tvr )(0). 
\]
Summing \eqref{eq:lap_recur} over  $i =  k,.., 2, 1$ with $j = k+ 1 - i$, we obtain 
\[
       \pa^{\bet} \Del^{k} (|y|^2 \Del \tvr)(0)
     = \bigl( (2 d + 4 |\bet|) k + 4 k (k-1)  \bigr)
     \pa^{\bet} \Del^{k} \tvr(0) + \pa^{\bet}( |y|^2 \Del^{k+1} \tvr )(0) .
\]
Using Leibniz rule, we obtain the estimate of the error terms in \eqref{eq:lap_id:a}, \eqref{eq:lap_id:b}
\[
  \pa^{\bet}( |y|^2 \Del^{k+1} \tvr )(0) 
  = C_{\bet, k} \cdot  \one_{|\bet| \geq 2}  \na^{|\bet|-2} \Del^{k+1} \tvr(0) 
  = \one_{|\bet| \geq 2}  \OOL{n} ( \FFs_{ |\bet |-2, k+1 })  .
\]

\vspace{0.1in}
\paragraph{\bf Proof of \eqref{eq:lap_id:c}}
Using the Leibniz rule, for any $k \geq 1$, we have 
\[
  \Del^{k+1}( |y|^{2k} f  )
= \Del^k( \Del (|y|^{2k}) f + |y|^{2k} \Del f 
+ 2 \na |y|^{2k} \cdot \na f  ) := I_1 + I_2 + I_3 .
\]

Denote 
\beq\label{eq:lap_constant}
  \msf c_{k, 1} := \prod_{1\leq i \leq k} 2i (2i + d - 2), 
  \quad \msf c_{k, 2} =  \prod_{1\leq i \leq k} 2i (2i + d + 2),
  \quad \msf c_{0, 2} = 1.
\eeq

Using Lemma \ref{lem:scal_id} and a direct computation, we obtain 
\[
\begin{aligned} 
& I_1=\Del^k ( 2k (2 k + d -2) |y|^{2k-2} f ),  
\quad I_3 =  \Del^k( 4 k |y|^{2k-2} y \cdot \na f  )
= 8 k \Del^k (  |y|^{2k-2} f  ), \\
& I_2 = \Del^k(|y|^{2k}) \Del f
 = 2k (2 k + d -2)  \Del^{k-1} (|y|^{2k-2}) \Del f
= \dots = \msf{c}_{k, 1} \Del f.
\end{aligned}
\]
It follows that
\[
    \Del^{k+1}( |y|^{2k} f  ) = (2k )(2 k + d + 2) \Del^k( |y|^{2k-2} f )
    + \msf{c}_{k, 1} \Del f.
\]
Dividing both sides by $\msf{c}_{k, 2}$, we obtain
\[
  \frac{     \Del^{k+1}( |y|^{2k} f  ) }{ \msf{c}_{k, 2} }
  =   \frac{     \Del^{k}( |y|^{2k-2} f  ) }{ \msf{c}_{k-1, 2} }
  +  \frac{ \msf{c}_{k, 1}  }{ \msf{c}_{k, 2}  } \Del f
  =  \frac{     \Del^{k}( |y|^{2k-2} f  ) }{ \msf{c}_{k-1, 2} } + \frac{ d(d+2) }{ (2k + d + 2)( 2k + d ) }
  \Del f.
\]
Using the recursive relation, the telescoping identity
$\frac{ d(d+2) }{ (2k + d + 2)( 2k + d ) } = \frac{(d+2)d}{2} ( \frac{1}{2k+d } - \frac{1}{2k+d+2} )$, 
and $\msf c_{0, 2} = 1$ from
\eqref{eq:lap_constant}, we deduce that
\begin{align*}
\frac{ \Del^{k+1}( |y|^{2k} f  ) }{ \msf{c}_{k, 2} }
&=  \frac{\Del f}{\msf{c}_{0, 2}} + \frac{(d+2)d}{2} \Bigl( \frac{1}{2 + d} - \frac{1}{2 k + d + 2} \Bigr)  \Del f
\notag\\
&= \Bigl( \frac{d+2}{2} - \frac{(d+2)d}{2 ( 2 k + d + 2 )} \Bigr) \Del f
= \frac{ (d+2)(k+1) }{ 2 k + d + 2 } \Del f .
\end{align*}
It follows that
\[
\Del^{k+1}( |y|^{2k} f  )
= \msf{c}_{k, 2}  \frac{ (d+2)(k+1) }{ 2 k + d + 2 }   \Del f 
= (k+1) \prod_{1 \leq i \leq k} 2i (2 i + d) \Del f
= \msf c_{k, \Del} \Del f.
\]
We complete the proof.
\end{proof}

\subsection{Main terms in \texorpdfstring{$n$-th}{n th} order ODEs}
In this section, we derive the main terms in the ODEs. Below, we fix any multi-index $\bet_1, \bet_2, \bet_3$ with $|\bet_1| = n, |\bet_2| = n+1, |\bet_3| = n+2$.

\vspace{0.05in}

\paragraph{\bf Scaling terms}

Since $\GGs_0$   \eqref{norm:ODE_low} contains the scaling terms $\tcr, \tcu, \tcvr, \tcb$, using  Lemma \ref{lem:scal_id} and the notation $\OOL{}$ from \eqref{norm:ODE_low}, we treat the linear terms involving $\tcr, \tcu, \tcvr, \tcb$ in  \eqref{eq:lin_ODE:cL}  as
\begin{subequations}\label{eq:ODE_main_scal}
\begin{equation}\label{eq:ODE_main_scal:lin}
\begin{aligned} 
   \pa^{\bet_1}(  -  \tcr y  \cdot \na \bvr     + \tcvr    \bvr)(0) 
   &= ( - n \tcr +     \tcvr )   \pa^{\bet_1} \bvr(0)
= \OOL{n}( \GGs_{ 0 }  ), \\
 \pa^{\bet_2} (   - \tcr y  \cdot \na \bu  + \tcu  \bu   )(0)
 &= ( - (n+1) \tcr  + \tcu) \pa^{\bet_2} \bu(0)
= \OOL{n}( \GGs_{ 0 } ) 
  , \\
 \pa^{\bet_3} (  - \tcr y  \cdot \na \barb   +  \tcb \barb    )(0)
 &= ( - (n+2) \tcr + \tcb ) \pa^{\bet_3} \barb(0)
 = \OOL{n}( \GGs_{ 0 } ) .
\end{aligned}
\end{equation}

 We estimate  the nonlinear terms involving $\tcr, \tcu, \tcvr, \tcb$ in \eqref{eq:lin_ODE} perturbatively 
\begin{equation}\label{eq:ODE_main_scal:non}
\begin{aligned}
      |\na^n(  -  \tcr y  \cdot \na  \tvr   + \tcvr  \tvr  )(0) | 
     +  |\na^{n+1} (   - \tcr y  \cdot \na \tu + \tcu  \tu  )(0)| 
       +  |\na^{n+2}(  - \tcr y  \cdot \na  \tb   +  \tcb  \tb  )(0) |
   \les_n | \GGs_n|^2.
   \end{aligned}
\end{equation}
\end{subequations}

The estimates of the terms $\bcvr \tvr, \bcu \tcu, \bcb \tb$ are straightforward 
\begin{equation}\label{eq:ODE_main_scal2}
  \pa^{\bet_1} \bcvr \tvr = \bcvr   \pa^{\bet_1} \tvr, 
  \quad   \pa^{\bet_2}\bcu \tu = \bcu  \pa^{\bet_2} \tu,
  \quad \pa^{\bet_3} \bcb \tb = \bcb \pa^{\bet_3} \tb .
\end{equation}

\paragraph{\bf Transport terms}
Next, we simplify the transport term $(\bcr y + \bu) \cdot \na f$ for $f = \tvr, \tu, \tb$ in \eqref{eq:lin_ODE:cL}. Since $\bu = \uua y +  \OO(|y|^{2\NNN+1})$, using the vanishing condition \eqref{eq:ODE:main:profi}, Lemma \ref{lem:scal_id} and the Leibniz rule, we obtain
\begin{subequations}\label{eq:ODE_main_tran}
\begin{align} 
  \pa^{\bet_1} (\bcr y + \bu) \cdot \na \tvr |_{y=0}  & = n  (\bcr + \uua) ( \pa^{\bet_1} \tvr)(0)  + 
  C_n \cdot  \na^{\leq n- 2 \NNN} \tvr(0)  \cdot \one_{n - 2\NNN-1\geq0} ,  \\
  \pa^{\bet_2} (\bcr y + \bu) \cdot \na \tu |_{y=0}  & = (n +1) (\bcr + \uua) ( \pa^{\bet_2} \tu)(0)  + C_n \cdot \na^{\leq n + 1- 2 \NNN} \tu(0) , 
  \\
  \pa^{\bet_3} (\bcr y + \bu) \cdot \na \tb |_{y=0}  & = (n + 2) (\bcr + \uua) ( \pa^{\bet_3} \tb)(0)  + C_n  \cdot  \na^{\leq n+2 - 2 \NNN} \tb(0) . 
\end{align}
From the definitions of $\GGs, \OOL{}(\GGs)$ in \eqref{norm:ODE_low}, we obtain 
\begin{equation}
 ( \na^{\leq n- 2 \NNN} \tvr(0) , \na^{\leq n + 1 - 2 \NNN} \tu(0) ,  \na^{\leq n+2 - 2\NNN} \tb(0) ) 
 =  \OOL{}(\GGs_{ (n - 2\NNN)_+ }).
\end{equation}
\end{subequations}

\paragraph{\bf Nonlinear terms}
For the nonlinear terms $\cB_i(\tw, \tw)$ with $\cB_i$ defined in \eqref{eq:bilin}, using 
\eqref{eq:vanishing:order:at:origin:SS}:
\begin{equation}\label{eq:ODE_van}
 \tu(0) = 0, \quad \na^i \tb(0) = 0, \quad i = 0, 1, 
\end{equation}
and the Leibniz rule, we obtain 
\begin{align}
 |\na^n \cB_1(\tw, \tw)(0) | 
 & \les_n \sum_{ 0 \leq i \leq n  + 1 } |\na^i \tu | \cdot |\na^{n+1 - i} \tvr(0)| \les_n |\GGs_n|^2, \notag \\
 |\na^{n+1} \cB_2(\tw, \tw)(0) | 
 & \les_n \sum_{ 0 \leq  i  \leq n +  2 }
 |\na^i \tu(0)| \cdot |\na^{n+2 - i} \tu(0)|
 + |\na^i \tb(0)| \cdot |\na^{n+2 - i} \tvr(0)| 
 \les_n |\GGs_n|^2, \notag \\
|\na^{n + 2} \cB_3(\tw, \tw)(0) | 
& \les_n \sum_{ 0 \leq  i  \leq n +  3 }
|\na^i \tu| |\na^{n+3 - i} \tb| \les_n |\GGs_n|^2. \label{eq:ODE_main_non}
\end{align}

\paragraph{\bf Remaining terms}

Next, we estimate the remaining terms in $\cL_i \tw$ \eqref{eq:lin_ODE:cL}, which are not covered by 
\eqref{eq:ODE_main_scal}, \eqref{eq:ODE_main_tran}, and \eqref{eq:ODE_main_non}.  
We choose arbitrary multi-indices $\bet_i$  with $|\bet_1| = n, |\bet_2| = n+1, |\bet_3| = n+2$. 
In the following derivations, we always evaluate at $ y= 0$ and do not indicate this dependence.

For the remaining terms in $\cL_1 \tw$ \eqref{eq:lin_ODE:cL}, since $|\bet_1| = n$, using Leibniz rule and \eqref{eq:ODE:main:profi}, we obtain 
\begin{subequations}\label{eq:ODE_main_rem}
\begin{align} 
 & \pa^{\bet_1}  ( - 2 \al \bvr (\div \tu) 
 - \tu \cdot \na \bvr 
 - 2 \al \tvr (\na \cdot \bu)    )   \\
& = \pa^{\bet_1} (  
- 2 \al \aaa  (\div \tu) 
- 2 \al  \tvr  \cdot (\na \cdot \bu)(0)  )
 + C_n \cdot \big( \na^{\leq |\bet_1| + 1 - 2 \NNN} \tu , \na^{\leq |\bet_1|- 2\NNN }\tvr  \big)  \\
 & = \pa^{\bet_1} \bigl(  
- 2 \al \aaa  (\div \tu) 
- 2 \al d \uua  \tvr \bigr)
   + \OOL{ n}( \GGs_{ (n- 2\NNN)_+ }),
\label{eq:ODE_main_rem:a}
 \end{align}
where we have used $\div \bu(0) = d \pa_r \bar U(0) = d \uua$,

For the remaining terms in $\cL_2 \tw$ \eqref{eq:lin_ODE:cL}, 
since $|\bet_2| = n+1$, using \eqref{eq:ODE:main:profi}, we estimate
\begin{align} 
 \pa^{\bet_2}  
\Bigl(  - \tfrac{1}{2\al} & \barb \na \tvr
  - \tfrac{1}{\gamma} \bvr \na \tb     -   \tu  \cdot \na \bu -  \tfrac{1}{2\al} \tb \na \bvr
  - \tfrac{1}{\gamma} \tvr \na \barb   \Bigr) \\
&  = \pa^{\bet_2}  \Big( - \tfrac{1}{2 \al} \bb |y|^2 \na \tvr  
 - \tfrac{1}{\gamma} \aaa \na \tb - \uua \tu - \tfrac{ 2 }{\gamma } \bb y \tvr 
 \Big)  \\
 & \qquad  + C_n \cdot \big( \na^{\leq |\bet_2|+1 - 2\NNN} \tb(0), \,  \na^{\leq |\bet_2|  - 2 \NNN} \tu(0) , \, \na^{\leq |\bet_2| -1 - 2\NNN} \tvr(0)   \big) \\
  & =  \pa^{\bet_2}  \Big( - \tfrac{1}{2 \al} \bb |y|^2 \na \tvr  
 - \tfrac{1}{\gamma} \aaa \na \tb - \uua \tu - \tfrac{2 }{\gamma } \bb y \tvr 
 \Big) 
  + \OOL{n}( \GGs_{ (n-2\NNN)_+ } ).
  \label{eq:ODE_main_rem:b}
\end{align}
For the remaining terms in $ \cL_3 \tw$, since $\bet_3$ with $|\bet_3| = n+2$, 
using \eqref{eq:ODE:main:profi}, we obtain 
\begin{equation}\label{eq:ODE_main_rem:c}
  \pa^{\bet_3} ( - \tu \cdot \na \barb)
  = \pa^{\bet_3} ( -\tu \cdot 2 \bb y ) + \OO_n(|\na^{\leq |\bet_3| - 1 - 2\NNN}  \tu|)
  = \pa^{\bet_3} ( - 2 \bb  y \cdot \tu ) + 
  \OOL{n}( \GGs_{ (n-2\NNN)_+ } ).
\end{equation}
\end{subequations}

\vspace{0.1in}

\paragraph{\bf Summary of the main terms}
We recall the parameter $\kp$ from \eqref{eq:ODE_main:kp}
\[
  \kp = \bcr + \uua. 
\]

\begin{subequations}\label{eq:ODE_main}
Therefore, combining \eqref{eq:ODE_main_scal}, \eqref{eq:ODE_main_scal2}, \eqref{eq:ODE_main_tran},
\eqref{eq:ODE_main_non}, \eqref{eq:ODE_main_rem}, we rewrite \eqref{eq:lin_ODE} and \eqref{eq:lin_ODE:cL} as 
\begin{align}
  \pa_\tau \tvr &= ( - n \kp - 2 \al d \uua + \bcvr) \tvr -
 2 \al \aaa  (\div \tu)  + \ell_{\vrho}
 =  - n \kp  \tvr - 2 \al \aaa  (\div \tu)  + \ell_{\vrho}, 
 \label{eq:ODE_main:a} \\
   \pa_\tau \tu &= (- (n+1) \kp - \uua + \bcu) \tu
    - \tfrac{1}{2 \al} \bb |y|^2 \na \tvr 
    - \tfrac{1}{\gamma} \aaa \na \tb - \tfrac{2}{\gamma} \bb y \tvr + \ell_U  \notag \\
    & = (- n \kp - 2 \uua -1) \tu
    - \tfrac{1}{2 \al} \bb |y|^2 \na \tvr 
    - \tfrac{1}{\gamma} \aaa \na \tb - \tfrac{2}{\gamma} \bb y \tvr + \ell_U
        , \label{eq:ODE_main:b} \\
  \pa_\tau \tb & = ( - (n+2) \kp  + \bcb ) \tb
  - 2\bb (y \cdot \tu) + \ell_B =  - n \kp   \tb
  - 2\bb (y \cdot \tu) + \ell_B ,
   \label{eq:ODE_main:c}
\end{align}
where we have used 
$ \bcvr = 2 \al d \uua$ (see \eqref{eq:profile_scal}) 
in the last identity in \eqref{eq:ODE_main:a}, the following identity 
\[
- (n+1) \kp - \uua + \bcu = - n \kp - \bcr - \uua - \uua + \bcr -1
= - n \kp - 2 \uua -1 
\]
in the last identity in \eqref{eq:ODE_main:b}  (see \eqref{eq:ODE_main:kp} and \eqref{eq:profile_scal}),  and 
$\bcb = 2 (\bcr + \uua) = 2 \kp$  (see \eqref{eq:ODE_main:kp} and \eqref{eq:profile_scal}) in the last identity in \eqref{eq:ODE_main:c}, and the error terms satisfy
\begin{equation}\label{eq:ODE_main_error}
   ( \na^n \ell_{\vrho}(0) , \    \na^{n+1} \ell_U(0) , \
    \na^{n+2} \ell_B (0) ) = \OOL{ n}( \GGs_{ (n-2\NNN)_+ } ) + 
\OO_n ( |\GGs_n |^2 ).
\end{equation}

\end{subequations}

\vspace{0.1in}

\paragraph{\bf ODEs in 1D}

For $d = 1$, we have $\na f = \pa_x f$. Applying Leibniz's rule to \eqref{eq:ODE_main}, then evaluating at $x=0$, and using the error bounds \eqref{eq:ODE_main_error}, we obtain the ODEs for $\pa_x^n \tvr(0), \pa_x^{n+1} \tu(0), \pa_x^{n+2} \tb(0)$ 
\[
\begin{aligned} 
  \tf{d}{d\tau} \pa_x^n \tvr &=  - n \kp   \pa_x^n  \tvr -
 2 \al \aaa   \pa_x^{n+1}  \tu  
 + \OOL{ n}( \GGs_{ (n-2\NNN)_+ }) + \OO_n( |\GGs_n|^2 ) , \\
  \tf{d}{d\tau} \pa_x^{n+1} \tu &=  (- n \kp - 2 \uua -1)   \pa_x^{n+1}  \tu
    - \bigl( \tfrac{ (n+1) n }{2 \al} + \tfrac{2(n+1)}{\gamma} \bigr)   \bb \pa_x^n \tvr 
    - \tfrac{1}{\gamma} \aaa \pa_x^{n+2} \tb  \\
    & \quad  + \OOL{ n}( \GGs_{ (n-2\NNN)_+ }) + \OO_n( |\GGs_n|^2 )  , \\
  \tf{d}{d\tau} \pa_x^{n+2} \tb & = - n \kp   \pa_x^{n+2} \tb
  - 2\bb (n+2) \pa_x^{n+1}  \tu 
+ \OOL{ n}( \GGs_{ (n-2\NNN)_+ }) + \OO_n( |\GGs_n|^2 ), 
  \end{aligned}
\]
We prove \eqref{eq:ODE_Nth:1D}.

\subsubsection{ODEs for $ y \cdot \tu, \div \tu$}

In this section, we further estimate the ODEs for $ \xu = y \cdot \tu$, $\div \tu$ and prove Proposition \ref{prop:ODE_Nth}.  
Using \eqref{eq:ODE_main} and Lemma \ref{lem:scal_id}, we obtain 
\begin{subequations}\label{eq:xu_iden}
\begin{align}
  | \na^{n+2} ( y \cdot |y|^2 \na \tvr  - n |y|^2 \tvr )(0) | 
  &\les_n |\na^n ( y \cdot \na \tvr - n \tvr  ) (0) | = 0, \label{eq:xu_iden:a}  \\
   \na^{n+2} ( y \cdot \na \tb - (n+2) \tb )(0) &= 0, \label{eq:xu_iden:b} \\
   \na^n\Big( \na \cdot ( |y|^2 \na \tvr) 
   - |y|^2 \Del \tvr - 2 n \tvr \Big)(0)
   & = \na^n ( 2 y \cdot \na \tvr  - 2 n \tvr)(0) = 0, \label{eq:xu_iden:c} \\
     \na^n ( \na \cdot ( y  \tvr) - (d + n) \tvr )(0)
     & = \na^n ( y \cdot \na \tvr - n \tvr )(0) = 0. \label{eq:xu_iden:d}
\end{align}
\end{subequations}

From \eqref{eq:ODE_main_error}, we obtain 
\begin{equation}\label{eq:ODE_main_err2}
\begin{aligned} 
  ( \na^{n+2}( y \cdot \ell_U )(0) , \ \na^n (\na \cdot \ell_U)(0) )
  = C_n \cdot \na^{n+1} \ell_U(0) 
  = \OOL{n} (\GGs_{ (n-2\NNN)_+ }) + \OO_n( |\GGs_n|^2 ) .  
\end{aligned}
\end{equation}

Deriving the ODE for $\xu = y \cdot \tu$ from \eqref{eq:ODE_main:b}, using \eqref{eq:xu_iden:a}, \eqref{eq:xu_iden:b} and the error estimate \eqref{eq:ODE_main_err2}, we obtain 
\begin{subequations}\label{eq:ODE_main_xu}
\begin{equation}
  \pa_\tau \xu =  (- n \kp - 2 \uua -1)   \xu
    - \bigl( \tfrac{n}{2 \al} + \tfrac{2}{\gamma} \bigr) \bb |y|^2 \tvr 
    - \tfrac{ (n+2) }{\gamma} \aaa  \tb 
    + \ell_{U, 1},
\end{equation}
with error $\ell_{U, 1}$ satisfying 
\begin{equation}
 \na^{n+2} \ell_{U, 1}(0)  = \OOL{ n}( \GGs_{ (n-2\NNN)_+ }) + \OO_n( |\GGs_n|^2 ).
\end{equation} 

\end{subequations}

Taking divergence in \eqref{eq:ODE_main:b}, using \eqref{eq:xu_iden:c}, \eqref{eq:xu_iden:d}, and the error estimate \eqref{eq:ODE_main_err2}, we obtain 
\begin{subequations}\label{eq:ODE_main_divu}
\begin{equation}\label{eq:ODE_main_divu:a}
\begin{aligned} 
\pa_\tau \div \tu & =  (- n \kp - 2 \uua -1)  \div \tu 
    - \tfrac{1}{ 2 \al} \bb (|y|^2 \Del \tvr + 2 n \tvr)
    - \tfrac{1}{\gamma} \aaa \Del \tb - \tfrac{2}{\gamma} \bb (d + n) \tvr + \ell_{U, 2} \\
    & =  (- n \kp - 2 \uua -1)   \div \tu 
    - \tfrac{1}{2 \al} \bb  |y|^2 \Del \tvr
-  ( \tfrac{n}{\al} + \tfrac{2( d+ n)}{\gamma} ) \bb \tvr
    - \tfrac{1}{\gamma} \aaa \Del \tb + \ell_{U, 2}
\end{aligned}
\end{equation}
with error $\ell_{U, 2}$ satisfying 
\begin{equation}\label{eq:ODE_main_divu:b}
 \na^{n} \ell_{U, 2}(0)  = \OOL{ n} ( \GGs_{ (n-2 \NNN)_+ } ) + \OO_n( |\GGs_n|^2 ).
\end{equation} 
\end{subequations}

\subsection{Proof of Proposition \ref{prop:ODE_Nth}}\label{sec:proof_ODE_Nth}

For $\bet = 0$ and $k = \frac{n}{2}$, applying $\Del^k$ to \eqref{eq:ODE_main:a}, using 
\eqref{eq:lap_id:div} from Lemma \ref{lem:lap}, $2k + 2 = n + 2$, 
and taking $y=0$, we derive 
\begin{align} 
  \tfrac{d}{d\tau} \Del^k \tvr 
&= - n \kp  \Del^k \tvr -
 2 \al \aaa  \Del^k (\div \tu)
 + \OOL{ n}( \GGs_{ (n-2\NNN)_+ } )+ \OO_n( |\GGs_n|^2 ) 
 \notag \\
 &=   - n \kp \Del^k \tvr -
  \tfrac{2 \al \aaa }{ n+2 } \Del^{k+1} ( y \cdot \tu)
+ \OOL{ n}( \GGs_{ (n-2\NNN)_+ } )+ \OO_n( |\GGs_n|^2 ) , \quad k = \tfrac{n}{2}.
\label{eq:ODE_deri1}
\end{align}

Recall $|\bet| + 2k = n$.  Applying $\pa^{\bet} \Del^k$ to \eqref{eq:ODE_main:a} and 
$\pa^{\bet} \Del^{k+1}$ to \eqref{eq:ODE_main:c}, and then using the error bound \eqref{eq:ODE_main_error}, 
we obtain 
\begin{subequations}\label{eq:ODE_deri2}
\begin{align}
  \tfrac{d}{d\tau} \pa^{\bet} \Del^k \tvr &=  - n \kp  \pa^{\bet} \Del^k \tvr -
 2 \al \aaa  \pa^{\bet} \Del^k (\div \tu)  + \OOL{ n}( \GGs_{ (n-2\NNN)_+ } )+ \OO_n( |\GGs_n|^2 ), 
\label{eq:ODE_deri2:a}
  \\
\tfrac{d}{d \tau} \pa^{\bet} \Del^{k+1}  \tb & =  - n \kp  
\pa^{\bet} \Del^{k+1} \tb
  - 2 \bb \pa^{\bet} \Del^{k+1} (y \cdot \tu) 
+ \OOL{ n}( \GGs_{ (n-2\NNN)_+ } )+ \OO_n( |\GGs_n|^2 )
  \label{eq:ODE_deri2:b}.
\end{align}
\end{subequations}

Applying \eqref{eq:lap_id:b} in Lemma \ref{lem:lap} and 
applying $\pa^{\bet} \Del^k$ to \eqref{eq:ODE_main_divu}, and then taking $y = 0$, we obtain 
\begin{align}
 \tfrac{d}{d\tau} \pa^{\bet}&  \Del^k \divu    
 =   (- n \kp - 2 \uua -1)   \pa^{\bet} \Del^k \divu
 -   \frac{1}{2 \al} \cdot ( 2 d + 4|\bet| + 4(k-1) ) k    \bb \pa^{\bet} \Del^k \tvr \notag \\
& - \bigl( \tfrac{n}{\al} + \tfrac{2(d+ n)}{\gamma} \bigr) \bb  \pa^{\bet} \Del^k \tvr
    - \frac{1}{\gamma} \aaa \pa^{\bet} \Del^{k+1} \tb + \pa^{\bet} \Del^{k}  \ell_{U, 2} 
+  \one_{ |\bet |\geq 2 } \OOL{ n}( \FFs_{|\bet|-2, k+1} )  
\label{eq:ODE_deri:divu} .
\end{align}

Applying \eqref{eq:lap_id:a} in Lemma \ref{lem:lap} and 
applying $\pa^{\bet} \Del^{k+1}$ to \eqref{eq:ODE_main_xu}, and then taking $y = 0$, we obtain 
\begin{align}
   \tfrac{d}{d\tau} \pa^{\bet} \Del^{k+1} \xu  = &  (- n \kp - 2 \uua -1)   \pa^{\bet} \Del^{k+1} \xu
    - \bigl( \tfrac{n}{2 \al} + \tfrac{2}{\gamma} \bigr)  
 (2 d + 4 |\bet| + 4k ) (k+1) \bb \pa^{\bet} \Del^k \tvr(0) \notag \\
& 
    - \tfrac{ (n+2) }{\gamma} \aaa  \pa^{\bet} \Del^{k+1} \tb 
    + \pa^{\bet} \Del^{k+1} \ell_{U, 1}
    +  \one_{ |\bet |\geq 2 } \OOL{n}( \FFs_{|\bet|-2, k+1} ) .  
    \label{eq:ODE_deri:xu}
\end{align}

\vspace{0.1in}

\paragraph{\bf Proof of ODEs \eqref{eq:ODE_Nth:spec:a},  \eqref{eq:ODE_Nth:spec:b}  }

When $n$ is even, combining \eqref{eq:ODE_deri1} with $k = \frac{n}{2}$, 
\eqref{eq:ODE_deri2:b} with $\bet = 0, k = \frac{n}{2}$, 
\eqref{eq:ODE_deri:xu} with $\bet = 0, k = \frac{n}{2}$, 
we prove \eqref{eq:ODE_Nth:spec:a}, \eqref{eq:ODE_Nth:spec:b}.

\vspace{0.1in}

\paragraph{\bf Proof of ODEs \eqref{eq:ODE_Nth}}

For general $\bet$, \eqref{eq:ODE_deri2} implies the ODEs for $\pa^{\bet} \Del^k \tvr$ and $\pa^{\bet} \Del^{k+1} \tb$ in \eqref{eq:ODE_Nth}.

Combining the second and third term in \eqref{eq:ODE_deri:divu} 
and applying the error bound \eqref{eq:ODE_main_divu} to $\pa^{\bet} \Del^k \ell_{U,2}$ in \eqref{eq:ODE_deri:divu}, we prove the ODE of $\pa^{\bet} \Del^k (\div \tu)$ in \eqref{eq:ODE_Nth}.

Applying the error bound \eqref{eq:ODE_main_xu} to $\pa^{\bet} \Del^{k+1} \ell_{U, 1}$ in \eqref{eq:ODE_deri:xu}, we prove the ODEs for $\pa^{\bet} \Del^{k+1} \xu$ in \eqref{eq:ODE_Nth}.  We complete the proof of \eqref{eq:ODE_Nth} in Proposition \ref{prop:ODE_Nth}.

\vspace{0.1in}

\paragraph{\bf Proof of ODEs \eqref{eq:ODE_Nth_ZB}, \eqref{eq:ODE_Nth_U}
for $\tb, \xu, \tu$}

For any multi-indices $\bet$ with $|\bet| = n+2$, taking $\pa^{\bet}$ to \eqref{eq:ODE_main_xu} and 
\eqref{eq:ODE_main:c}, and then using
\[
   \pa^{\bet} ( |y|^2 \tvr )(0)  = C_{\beta} \cdot \na^n \tvr(0) , 
\]
we prove \eqref{eq:ODE_Nth_ZB}. For any multi-indices $\bet$ with $|\bet| = n+1$, taking $\pa^{\bet}$ to \eqref{eq:ODE_main:b}, and then using 
\[
  \pa^{\bet} (|y|^2 \na \tvr)(0)  = C_{\beta} \cdot \na^n \tvr(0),
  \quad   \pa^{\bet} ( y  \tvr)(0)  = C_{\beta} \cdot \na^n \tvr(0),
  \quad   \pa^{\bet} \na \tb(0) = C_{\beta} \cdot \na^{n+2 } \tb(0),
\]
we prove \eqref{eq:ODE_Nth_U}.

\subsubsection{Special case $\Del^{\NNN} \tvr, \Del^{\NNN+1} (y \cdot \tu),
 \Del^{\NNN+1} \tb$}
To derive \eqref{eq:ODE_Nth:spec} with $n = 2 \NNN, k = \NNN$, we derive the next order terms.
Applying the expansion \eqref{eq:ODE:main:profi} to \eqref{eq:lin_ODE:a}-\eqref{eq:lin_ODE:c}, 
for $k \leq \NNN +1$, we obtain 
\begin{subequations}\label{eq:L_form}
\begin{equation}
  \Del^k \cL_i = \Del^k \cL_{i,  L} + \Del^k \cL_{i, H},
\end{equation}
where $\cL_{i, L}$ depends on the profile $(\bvr, \bu, \barb)$ via the leading order terms in 
\eqref{eq:ODE:main:profi}, and we keep all the terms involving $\tcu, \tcr, \tcvr, \tcb$ in $\cL_{i, H}$
\begin{equation}
\begin{aligned} 
  \cL_{1,L} = &-( \bcr y + \uua y ) \cdot \na \tvr - 2 \al \aaa (\div \tu)  + \bcvr \tvr 
 -  \tu \cdot \na \aaa - 2 \al \tvr (\na \cdot (\uua y) )
 ,  \\
\cL_{1, H} =& - \uun y |y|^{2\NNN} \cdot \na \tvr + \Big( - 2 \al \aan |y|^{ 2\NNN } (\div \tu)   
 -  \tu \cdot \na (\aan |y|^{2\NNN}) - 2 \al \tvr (\na \cdot (\uun y |y|^{2\NNN})) \Big)  \\
 & + ( -\tcr y \cdot \na \bvr  + \tcvr  \bvr   ) 
 \notag\\
 &:= \cT_1 + \cR_1 + \cS_1
 ,\\
 \cL_{2, L}  =&  -  ( \bcr y + \uua y  ) \cdot \na \tu  - \tfrac{1}{2\al} \bb |y|^{2} \na \tvr
  - \tfrac{1}{\gamma} \aaa \na \tb    +  \bcu  \tu  -  \tu  \cdot \na (\uua y) -  \tfrac{1}{2\al} \tb \na \aaa
  - \tfrac{1}{\gamma} \tvr \na (\bb |y|^2) , \\
\cL_{2, H}  =& - \uun y |y|^{2 \NNN} \cdot \na \tu + \Big(- \tfrac{1}{2\al} \bbn |y|^{2 \NNN + 2} \na \tvr
  - \tfrac{1}{\gamma} \aan |y|^{2\NNN} \na \tb      -  \tu  \cdot \na (\uun y |y|^{2\NNN})
   \\
  & -  \tfrac{1}{2\al} \tb \na (\aan |y|^{2\NNN}) 
   - \tfrac{1}{\gamma} \tvr \na ( \bbn |y|^{2\NNN+2} ) \Big) + ( - \tcr y  \cdot \na \tu + \tcu \bu ) 
\notag\\
&:= \cT_2 + \cR_2 + \cS_2, \\
\cL_{3, L} = & - ( \bcr y + \uua y) \cdot \na \tb   +  \bcb \tb 
  -  \tu  \cdot \na (\bb |y|^2)  , \\
  \cL_{3, H} = &  - \uun y |y|^{2\NNN} \cdot \na \tb   
  -  \tu  \cdot \na ( \bbn |y|^{2\NNN+2} ) 
  +(-  \tcr y \cdot \na \barb + \tcb \barb) \notag\\
&:= \cT_3 + \cR_3 + \cS_3, 
\end{aligned}
\end{equation}
\end{subequations}
where we have used the notation $\cT, \cR, \cS$ in $\cL_{i, H}$ to single out \textit{transport, remaining},  and  \textit{scaling} terms. Below, we evaluate at $y=0$. We have estimated 
\[
  \Del^{\NNN} \cL_{1, L}(0),
  \quad 
    \Del^{\NNN + 1} ( y \cdot \cL_{2, L} )(0) ,
    \quad 
    \Del^{\NNN+1}  \cL_{3, L} (0),
\]
in the right hand side in \eqref{eq:ODE_deri1}, \eqref{eq:ODE_deri2:b} with $\bet =0$, and \eqref{eq:ODE_deri:xu} with $\bet =0$, without the error terms 
$ \OOL{}( \GGs_{ (n-2\NNN)_+ } )+ \OO_n( |\GGs_n|^2 )$.
In these equations, we treat $\cL_{i,H}$ as a lower-order linear term and track its contribution 
implicitly using $\mathcal{O}_n(\GGs_0)$.
Below, we further expand
\[
  \Del^{\NNN} \cL_{1, H}(0),
  \quad 
    \Del^{\NNN + 1} ( y \cdot \cL_{2, H} )(0) ,
    \quad 
    \Del^{\NNN+1} \cL_{3, H}(0) .
\]
 
  For  transport terms $\cT_i$ \eqref{eq:L_form}, 
  using \eqref{eq:scal_id:b} in Lemma \ref{lem:scal_id} 
 and then \eqref{eq:lap_id:c},\eqref{eq:lap_id:d}, 
 and $\Del( \tu \cdot y ) = 2 ( \div \tu ) $, we derive
\beq\label{eq:ODE:spec:tran}
\bal
   \Del^{\NNN} ( \uun y |y|^{2\NNN} \cdot \na \tvr )(0 ) & = 0  ,  \\
   \Del^{\NNN+1}  ( \uun |y|^{2\NNN}  y \cdot (y \cdot \na \tu ) ) 
  & = \Del^{\NNN+1}  ( \uun |y|^{2\NNN}  y \cdot \tu ) 
  = \ccs_{\NNN, \Del} \uun \Del( y \cdot \tu ) \\
 & =  2 \ccs_{\NNN, \Del} \uun \div \tu ,  
  \\
    \Del^{\NNN+1} ( \uun y |y|^{2\NNN} \cdot \na \tb )
 & = 2 \uun \Del^{\NNN+1} ( |y|^{2\NNN} \tb )
 = 2 \ccs_{\NNN, \Del}  \uun   \Del \tb.
\eal
\eeq

For the scaling terms $\cS_i$, applying \eqref{eq:ODE_main_scal} with $n = 2\NNN, |\bet_1| = n, |\bet_2| = n+1, |\bet_3| = n+2 $, 
and then using the constant in \eqref{eq:lap_id:d}, and $\tcvr = 2 \tcu - \tcb$ \eqref{eq:modulated:tcvr}, we derive 
\begin{equation}\label{eq:ODE:spec:scal}
\begin{aligned} 
   \Del^{\NNN} (   -  \tcr y  \cdot \na \bvr     + \tcvr    \bvr)(0) 
  & = (    - 2\NNN \tcr +   \tcvr   )   \aan \Del^{\NNN} ( |y|^{2 \NNN} ) 
   = (  2 \tcu - \tcb - 2\NNN \tcr  )   \aan \ccs_{\NNN, 1} 
   ,\\
 \Del^{\NNN+1} (  y \cdot (  - \tcr y  \cdot \na \bu  + \tcu  \bu )  )(0)
 & = \Del^{\NNN+1} \Big( y \cdot ( - ( 2\NNN+1) \tcr  + \tcu) \uun  y |y|^{2 \NNN}  \Big) \\
 &   =   (  \tcu - ( 2\NNN+1) \tcr  ) \uun \Del^{\NNN+1} ( |y|^{2 \NNN+2} )   \\
 & =  (   \tcu- ( 2\NNN+1) \tcr  ) \uun \ccs_{\NNN+1, 1} ,
   \\
  \Del^{\NNN+1} (  - \tcr y  \cdot \na \barb   +  \tcb \barb    )(0)
 & = ( - (2 \NNN+2) \tcr + \tcb ) \bbn  \Del^{\NNN+1} ( |y|^{2 \NNN+2} )  \\
 &  = ( \tcb - (2 \NNN+2) \tcr  ) \bbn \ccs_{\NNN+1, 1}.
\end{aligned} 
\end{equation}

Next, we estimate the remaining terms $\cR_i$ in \eqref{eq:L_form}.
For $\cR_1$, using \eqref{eq:lap_id:c}, \eqref{eq:lap_id:d}, we have
\begin{subequations}    \label{eq:ODE:spec:rem1}
\begin{align}
\Del^{\NNN} \cR_1 &= 
  \Del^{\NNN} ( - 2 \al \aan |y|^{ 2\NNN } (\div \tu)   
 -  \tu \cdot \na (\aan |y|^{2\NNN}) - 2 \al \tvr (\na \cdot (\uun y |y|^{2\NNN}))    ) \notag \\
 & =  
 - 2\al \aan \Del^{\NNN} ( |y|^{2\NNN} ) ( \div \tu )
 -  2\NNN \aan \Del^{\NNN}( |y|^{ 2\NNN-2 } \tu \cdot y  )
 - 2 \al  (2 \NNN + d ) \uun \tvr \Del^{\NNN}(|y|^{2 \NNN}) \notag  \\
 & = 
 -2\al \aan \msf c_{\NNN, 1} (\div \tu)
 - 2 \NNN \aan \msf c_{\NNN-1, \Del} \Del (\tu \cdot y)
   - 2 \al  (2 \NNN + d ) \uun \msf c_{\NNN, 1} \tvr .
 \end{align}
Since $\Del( \tu \cdot y ) = 2 ( \div \tu )$, we further obtain
\begin{equation}
  \Del^{\NNN} \cR_1 = 
   - ( 2\al \aan \msf c_{\NNN, 1} +  4 \NNN \aan \msf c_{\NNN-1, \Del} ) (\div \tu)
   - 2 \al  (2 \NNN + d ) \uun \msf c_{\NNN, 1} \tvr.
\end{equation}
\end{subequations}

For $\cR_2$ \eqref{eq:L_form}, by definition, we have
\begin{align*} 
  \Del^{\NNN+1}( y \cdot \cR_2 )
  = &  \Del^{\NNN+1} \Big(- \tfrac{1}{2\al} \bbn  |y|^{2 \NNN + 2} ( y \cdot \na \tvr)
  - \tfrac{1}{\gamma} \aan |y|^{2\NNN} (y \cdot \na \tb )      - y\cdot \Big( \tu  \cdot \na (\uun y |y|^{2\NNN}) 
  \Big)    \\
  &\qquad \qquad -  \tfrac{1}{2\al} \tb  y \cdot \na (\aan |y|^{2\NNN}) 
   - \tfrac{1}{\gamma} \tvr  y \cdot \na ( \bbn |y|^{2\NNN+2} ) \Big).
\end{align*}
The first term on the right side is $0$. Using Lemma \ref{lem:scal_id} and then \eqref{eq:lap_id:c}, \eqref{eq:lap_id:d}, 
we simplify the above terms as
\[
\begin{aligned} 
  \Del^{\NNN+1} (  |y|^{2 \NNN } ( y \cdot \na \tb) )
 & = 2  \Del^{\NNN+1} (  |y|^{2 \NNN } \tb )
  = 2 \ccs_{\NNN, \Del} \Del \tb, \\
    y\cdot ( \tu  \cdot \na ( y |y|^{2\NNN} ) ) 
 &  = y \cdot \Big( \tu  \cdot ( \Id |y|^{2\NNN})  + 2 \NNN ( \tu \cdot y)  y  |y|^{2\NNN-2} \Big)  
= (2 \NNN+1) (y \cdot \tu) |y|^{2\NNN} ,
     \\
  \Del^{\NNN+1} (   y\cdot ( \tu  \cdot \na ( y |y|^{2\NNN} ) ) )
 & =  (2\NNN+1)  \Del^{\NNN+1}( (y \cdot \tu) |y|^{2\NNN}  )
  =  (2\NNN+1)  \ccs_{\NNN, \Del} \Del(y \cdot \tu) ,  \\
\Del^{\NNN+1} ( \tb y \cdot \na |y|^{2\NNN} )
& = \Del^{\NNN+1} ( 2 \NNN \tb  |y|^{2\NNN} )
 = 2 \NNN \ccs_{\NNN, \Del} \Del \tb , \\
\Del^{\NNN+1} (\tvr   y \cdot \na (  |y|^{2\NNN+2} ) )
& =(2\NNN+2) \Del^{\NNN+1} (\tvr   |y|^{2\NNN+2}  )
= (2\NNN+2 ) \tvr \Del^{\NNN+1} (   |y|^{2\NNN+2}  )
=  (2\NNN+2 ) \ccs_{\NNN+1, 1} \tvr .
  \end{aligned}
\]

Using the above identities and $\Del( y \cdot \tu ) =2 \div \tu$ , we derive 
\begin{align} 
&\Del^{\NNN+1}( y \cdot \cR_2 )
\notag\\
& = - \tfrac{2}{\gamma}  \ccs_{\NNN, \Del} \aan \Del \tb
- (2\NNN+1)  \ccs_{\NNN, \Del} \uun \Del(y \cdot \tu)
- \tfrac{1}{\al} \NNN \cdot \ccs_{\NNN, \Del} \aan \Del \tb 
- \tfrac{1}{\gamma} (2\NNN+2 )  \ccs_{\NNN+1, 1}  \bbn \tvr  
\notag\\
& = - ( \tfrac{2}{\gamma} + \tfrac{ \NNN }{\al} ) \aan \ccs_{\NNN,\Del} \Del \tb
-  (4\NNN+2)  \ccs_{\NNN, \Del} \uun  (\na  \cdot \tu)
- \tfrac{1}{\gamma} (2\NNN+2 )  \ccs_{\NNN+1, 1}  \bbn \tvr .
\label{eq:ODE:spec:rem2}
 \end{align}

For $\cR_3$ \eqref{eq:L_form}, using \eqref{eq:lap_id:c} and $\Del(\tu \cdot y) = 2 \div \tu$, we derive 
\begin{align} 
\Del^{\NNN+1} \cR_3  
 &  = \Del^{\NNN+1} (   -  \tu  \cdot \na ( \bbn |y|^{2\NNN+2} )  )
  =- (2 \NNN+2)\bbn  \Del^{\NNN+1} ( \tu \cdot y |y|^{2 \NNN} ) 
  \notag \\
 &  = - (2 \NNN+2)\bbn \ccs_{\NNN, \Del} \Del(\tu \cdot y)
 = - (4 \NNN+4)\bbn \ccs_{\NNN, \Del} ( \na  \cdot \tu).
 \label{eq:ODE:spec:rem3}
\end{align}
Combining the estimates of $\cL_{i, L}$ in the right hand side of \eqref{eq:ODE_deri1}, \eqref{eq:ODE_deri2:b} and \eqref{eq:ODE_deri:xu} with $(n, k,\bet) =(2\NNN,\NNN, 0)$, and replacing $\OOL{0}( \GGs_{ (n-2\NNN)_+ } )$  by the above estimates of $\cL_{i, H}$, we prove 
\eqref{eq:ODE_Nth:spec:c}--\eqref{eq:ODE_Nth:spec:f}. We complete the proof of Proposition~\ref{prop:ODE_Nth}.

%%%%%%%%%%%%%%%%%%%%%%%%%%%%%%%%%%%%%%%%%%%

\section{Proof of the results in Part I of~Theorem~\ref{thm:uns_specific}}\label{sec:eigen}

In this section we prove the claims in Part I of Theorem~\ref{thm:uns_specific}; these concern on the number of unstable eigenvalues for the ODE systems at $y=0$ for ground state profile corresponding to $\NNN = 1$, for either $(\gamma , d) \in \bigl\{ (\frac53, 3), ( \frac75, 3), ( \frac53, 2 ), (2,2)\}$, or for $d = 1$ and $\gamma \in (1,3]$. The proofs follow from Lemma~\ref{lem:Hbk} and Lemma~\ref{lem:H1d} below.

To count the number of unstable eigenvalues, we use the Routh–Hurwitz Theorem, a classical stability result in control theory. Given an $l-$th degree polynomial $p_l(x) = a_l x^l + \ldots + a_0$. 
The first two column of the Routh table $A$ are given by the coefficients of $p_l$
\[
A = \left(
    \begin{array}{cccc}
   a_l  & a_{l-2} & a_{l-4}  & ... \\
   a_{l-1} & a_{l-3} & a_{l-5} & ... \\
   ... & ... & ...   
   \end{array}   
   \right),
   \qquad A_{1j} = a_{l+ 2 - 2 j } \one_{2j \leq l+2}, \qquad 
   A_{2 j} = a_{l+1 - 2j} \one_{2j \leq l+1}.
\]
For $i \geq 3$, given the $i-2$ and $i-1$-th rows of $A$, the $i-$th row of the Routh table is constructed as 
\[
  A_{i j} =  A_{i-2, j+1} - A_{i-1, j+1} \cdot \frac{ A_{ i-2, 1} }{A_{i-1, 1}} ,  \qquad  j \geq 1.
\]

\begin{theorem}[Routh–Hurwitz]\label{thm:RH}

Suppose that $A_{i 1} \neq 0$ for $ 1 \leq i \leq l+1$.  The number of roots of $p_l(x)$ with 
 positive real parts is equal to the number of changes in sign of the first column of the
Routh array, i.e. the number of negative terms in $\frac{A_{i+1,1}}{A_{i,1}}, 1 \leq i \leq l$. Moreover, there is no root with $0$ real part.
\end{theorem}

See \cite[Chapter V]{gantmacher2005applications} for detailed discussion and \cite{meinsma1995elementary,bodson2019explaining} for simpler proofs. The Routh–Hurwitz theorem can also handle the case with $A_{i1} = 0$ for some $i$. We only use the regular case with $A_{i1} \neq 0$ for simplicity.
We apply Theorem \ref{thm:RH} to analyze the characteristic polynomial of $H \in \Reals^{l \times l}$ 
\[
  p_l(\lam) = \det(\lam \Id - H) =  \lam^l + a_{l-1} \lam^{l-1} + ... + a_1 \lam + a_0,
  \quad a_l = 1.
\]
with $l = 3,4$. In particular, we use $(H, l) = (\HH_{|\bet|, k}, 4), (\HH^{(1)}_n, 3)$.

\subsection{Estimates of \texorpdfstring{$\HH_{|\bet|, k}$}{H beta k}}
We can derive the first column of the Routh table as 
\begin{subequations}
\begin{equation}
  A = \left(
  \begin{array}{ccc c}
  a_4  & a_2 & a_0 & 0  \\
  a_3 & a_1 & 0  & 0 \\
 \frac{\Del_2}{a_3} & a_0 & 0 & 0 \\
  \frac{\Del_3}{\Del_2} & 0 & 0 & 0\\
  a_0 & 0 & 0 & 0 \\
\end{array}
  \right) , 
  \quad 
A_{31}= a_2 - a_1 \frac{a_4}{a_3} = \frac{ \Del_2 }{a_3},  \quad
A_{41} = a_1 - a_0 \frac{ a_3}{ \frac{\Del_2}{a_3} } = \frac{\Del_3}{\Del_2}, 
\end{equation}
where
 $a_4, \Del_2, \Del_3$ are given by 
\begin{equation}
a_4 = 1, \quad 
   \Del_2 =  a_3 a_2 - a_1 a_4 ,
   \quad  \Del_3 = a_3 a_2 a_1 - a_3^2 a_0 - a_1^2 a_4 .
\end{equation}
\end{subequations}

To apply Theorem \ref{thm:RH} to $\HH_{|\bet|, k}$, we estimate the signs of $a_3, \Del_2, \Del_3$, and $a_0$, which 
determine the signs of $A_{j 1}$ for $1 \leq j \leq 5$.

\begin{lemma}\label{lem:Hbk}
Consider $\NNN = 1$, $(\gamma , d) \in \{ (\frac53, 3), ( \frac75, 3),  (2,2), ( \frac53, 2 )\}$, and $n = |\bet| + 2 k$.
Recall the matrices 
$\HH_{|\bet|, k}$ from \eqref{eq:ODE_Nth}, 
$\HHH_{ k}$ from \eqref{eq:ODE_Nth:spec:b}, 
$\HH_n^{(2)}$ from  \eqref{eq:ODE_Nth_ZB},
and $ (- n \kp - 2\uua - 1 ) \Id$ from \eqref{eq:ODE_Nth_U}.

\paragraph{\bf (I) Estimate of $\HH_{|\bet|, k} \in \Reals^{4 \times 4}$ }

For any $n \geq 1$, we have 
\begin{subequations}\label{eq:eig_Hbk}
\begin{align}
a_3 > 0, & \quad  \Del_2 > 0 , \quad  \Del_3 >0 .
\end{align}
For $a_0$,   we have the following estimates:
\begin{align}
a_0 >0, & \quad && \forall \   |\bet| + 2 k \geq 2, \quad  (|\bet|, k) \neq (0, 1), \\
a_0 < 0, & \quad  && (|\bet|, k) = (1, 0) .
\end{align}

\end{subequations}

Thus, for $|\bet| + 2 k \geq 2$ and $(|\bet|, k) \neq (0, 1)$, $\HH_{|\bet|, k}$ has $4$ eigenvalues with negative real part; 
for $(|\bet|, k) = (1, 0) $, $\HH_{|\bet|, k}$ has $1$ eigenvalue with positive real part, and $3$ eigenvalues with negative real part.

\paragraph{\bf (II) Estimate of $\HHH_{ k} \in \Reals^{3 \times 3}$ } 

For $ k \geq 2$,  $\HHH_{ k}$ has $3$ eigenvalues with negative real part.

\paragraph{\bf (III) Estimate of $\HH_n^{(2)} \in \Reals^{2 \times 2} $}

For $\HH_n^{(2)} \in \Reals^{2 \times 2}$, the signs of the real parts of its eigenvalues are determined by its trace and determinant, which satisfy the following estimates.

For any $n \geq 1$, $\gamma > 1$ and $d \geq 1$,   we have 
\begin{equation}\label{eq:eig_HHT_tr}
  \tr(\HHT_{n}) < 0.
\end{equation}

For $ \det( \HHT_{n} ) $,  we have 
\begin{subequations}\label{eq:eig_HHT_det}
\begin{align}
     \det( \HHT_{n} ) >0 , & \quad   && (\gamma , d) =(\tf53, 3),  \  (2,2), \ n \geq 1 , 
\label{eq:eig_HHT_det:a}
      \\
       \det( \HHT_{n} ) >0 , &  \quad && (\gamma, d) =  ( \tf75, 3), \  (\tf53, 2),     
       \ n \geq 2 , 
\label{eq:eig_HHT_det:b}
        \\
         \det( \HHT_{1} )< 0, &  \quad  &&  (\gamma, d) =
          ( \tf75, 3), \  (\tf53, 2) , \ n = 1. 
         \label{eq:eig_HHT_det:c}
\end{align}
\end{subequations}

Thus, $\HHT_{n}$ has $2$ eigenvalues with negative real part for $(\gamma, d)
= (\frac53, 3),  \  (2,2)$ with $n \geq 1$
and for $(\gamma, d) =  ( \frac75, 3),  \ (\frac53, 2)$ with $n \geq 2$; 
it has $1$ eigenvalue with negative real part and $1$ with positive real part for  $(\gamma, d) =
  ( \frac75, 3), \ (\frac53, 2)$ and $n = 1$.

\paragraph{\bf(IV)}  
For any $n \geq 1$, $\gamma > 1$, and $d \in [1,4]$ satisfying $\al d > \frac{1}{2}$, we have
\[
- n \kp - 2\uua - 1 < 0.
\]
In particular, the assumptions on $(\gamma, d)$ are satisfied for
$(\gamma , d) \in \{ (\frac53, 3), ( \frac75, 3),  (2,2), ( \frac53, 2 )\}$.

\paragraph{\bf(V)} 

When $n=1$, $\HH_{1, 0}$ has exactly one positive eigenvalue. When $n=1, 
(\gamma, d)  \in \{ (\tf 75, 3), (\tf 53,  2 )  \}$, $\HHT_1$ has exactly one positive eigenvalue.
\end{lemma}

Before we prove the above Lemma, we have a simple estimate of $\bcr$. We denote $w = \al d$. 

\begin{lemma}\label{lem:cr_low}
Consider $\NNN=1$. For any $\gamma > 1, d \leq 4$, we have $\bcr > 1 +  \frac{(1-\al d)_+}{2(1 + \al d)} $.
\end{lemma}

\begin{proof}[Proof of Lemma~\ref{lem:cr_low}]
Recall $\bcr$ from formula \eqref{eq:cx:admissible} and $\mathsf{E}_\NNN$ from \eqref{eq:En:def} 
\begin{equation}\label{eq:cr_recall}
\bcr(d,\gamma,\NNN)
:= \tfrac{1}{1+\alpha d} \Bigl(1  +  \sqrt{ \tfrac{\alpha \gamma d}{2}  + \mathsf{E}_\NNN  + \tfrac{(1-\alpha d)^2}{16 \NNN^2}  } + \tfrac{1-\alpha d}{4\NNN}   \Bigr),
\quad \mathsf{E}_\NNN =  
 \tfrac{\alpha \gamma d (d+2)}{4 \NNN} 
+ \tfrac{\alpha d (1+ \alpha d  ) }{2 \NNN^2}   
.
\end{equation}
Using $\gamma = 2 \al + 1$ and $d \leq 4$,  we obtain 
\[
    \gamma - \tf{1}{2} | \al d -1|
\geq 2 \al + 1  - \tf12 \al d - \tf 12
> \tf12 \al ( 4-d)  \geq 0.
\]
For $\NNN=1$, it follows 
\[
\bal
\tfrac{\alpha \gamma d}{2}  + \mathsf{E}_1  + \tfrac{(1-\alpha d)^2}{16 }  
& \geq \tfrac{\alpha \gamma d}{2} + 
\tfrac{\alpha \gamma d}{2}  + \tf{ \al ( 2 \al + 1) d^2 }{4}
+ \tf{ (\al d)^2}{2} + \tfrac{(1-\alpha d)^2}{16 }  
> \al \gamma d +  ( \al d)^2 + \tfrac{(1-\alpha d)^2}{16 }   \\
& \geq \al d \cdot \tf{ |1-\al d|}{2} +   ( \al d)^2 + \tfrac{(1-\alpha d)^2}{16 }
\geq \tf{ (|1-\al d| + 4\al d)^2}{16}. 
\eal 
\]
Plugging the above estimate in the formula of $\bcr$, we obtain 
\[
  \bcr  > \tfrac{1}{1+\al d} \Big(1 + \tfrac{4 \al d + |1-\al d|}{4} + \tfrac{1-\al d}{4} \Big)
  = 1 + \tfrac{(1-\al d)_+}{2 (1 + \al d)},
\]
which completes the proof.
\end{proof}

In the following subsections, we estimate the signs of the entries in the first column of the Routh array for the matrices $\HH_{|\bet|, k}$ \eqref{eq:ODE_Nth}, $\HHH_{k}$ \eqref{eq:ODE_Nth:spec:b}, and similar terms for $\HH_n^{(2)}$ from \eqref{eq:ODE_Nth_ZB}. Since these matrices are given explicitly, we first use symbolic computation to derive formulas for $\Delta_2$, $\Delta_3$, and $a_3$ in \eqref{eq:eig_Hbk}, 
and related variables. Based on these formulas, we prove the inequalities in Lemma \ref{lem:Hbk} analytically.

\subsection{Estimate of eigenvalues of \texorpdfstring{$\HH_{|\bet|, k}$}{H beta k}}
In this section, we estimate $a_3, \Del_2, \Del_3, a_0$ in Lemma \ref{lem:Hbk}.

\subsubsection{Estimate of $a_3$}\label{sec:est_a3}

Using a direct calculation, $w = \al d$, and $\bcr \geq 1$ from Lemma \ref{lem:cr_low}, we obtain 
\[
  a_3 = \frac{-2 + 4(\bcr - 1) n  + 2 \al (d + 2 \bcr d n )  }{1 + \al d}
\geq \frac{ -2 + 4(\bcr - 1)  + 2 w + 4 \bcr  w  }{1 + w} .
\]

Using Lemma \ref{lem:cr_low}, $ d \geq 1$, and $w >0$, we prove 
\[
  4(\bcr -1) - 2 + 2 w + 4 \bcr  w 
  \geq \frac{2 (1-w)}{1 + w} + 6 w - 2
  \geq 6 w - \frac{4w}{1 +w} > 0.
\]

\subsubsection{Estimate of $\Del_2$}\label{sec:est_del2}

To simplify notation, we denote $b = |\bet|$. 
Below, we prove $\Del_2 >0$ for $ n = |\bet|+ 2 k \geq 1$, and discuss four cases 
$(\gamma, d) \in \{ (\tf53, 3) ,  (\tf75, 3),  ( 2, 2 ), (\tf53, 2)\}$ separately. 

\vspace{0.1in}
\paragraph{\bf Case $( \gamma, d) = (\frac53, 3) $}
 We obtain
\[
  \Del_2 = \frac{1}{48} \sqrt{\frac{47}{3}} (b+2 k)
  \underbrace{ \left(235 b^2+20 b (46 k-1)+920 k^2-50 k-24\right) }_{:=I}
\]

If $k \geq 1$, since each coefficient of $b^i, i=0,1,2$ in $I$ is positive, we prove $I > 0$. If $k = 0$, 
we get  $I = 235 b^2 - 20 b  -24,$ which is positive since $b = n \geq 1$. We prove $\Del_2 > 0$.

\vspace{0.1in}
\paragraph{\bf Case $( \gamma, d) = ( \frac75, 3  ) $}
 We obtain
\[
\Del_2 =  \frac{(15 b+30 k-2) \left(1125 b^2+192 b (23 k-2)+4416 k^2-810 k-80\right)}{1024}
:=  \frac{(15 b+30 k-2) \times I}{1024},
\]
where $I$ denotes the second term in the numerator. Since $ b + 2k = n \geq 1$, we get $15 b+30 k-2 >0$. If $k \geq 1$, since each coefficient of $b^i, i=0,1,2$ in $I$ is positive, we prove $ I > 0$. If $k = 0$, 
we get  $I = 1125 b^2+192 b \cdot(-2) -80 >0$ since $b = n \geq 1$.
Thus, we prove $\Del_2>0$. 

\vspace{0.1in}
\paragraph{\bf Case $( \gamma, d) = (2, 2) $}
We obtain
\[
  \Del_2 = 2 (b+2 k) \udb{ \left(10 b^2+b (39 k-1)+39 k^2-2 k-1\right) }_{:=I}.
\]

If $k \geq 1$,  since each coefficient of $b^i, i=0,1,2$ in $I$ is positive, we prove $I > 0$. If $k = 0$, we get $I = 10 b^2 - b - 1 > 0$  since $b = n \geq 1$. Thus, we prove $\Del > 0$.

\vspace{0.1in}
\paragraph{\bf Case $( \gamma, d) = (\frac53, 2) $}
We obtain 
\[
\begin{aligned} 
  \Del_2 &= \tfrac{\left(\sqrt{321} b+b+2 k+2 \sqrt{321} k-2\right) \cdot \left(5 \left(\sqrt{321}+161\right) b^2+10 b \left(2 \left(\sqrt{321}+157\right) k-\sqrt{321}-9\right)+4 \left(5 \left(\sqrt{321}+157\right) k^2-5 \left(\sqrt{321}+9\right) k-18\right)\right)}{1000} \\
  & := \tfrac{I \cdot II}{1000},
\end{aligned} 
\]
where $I$ and $II$ denote the first and second terms, respectively, in the product appearing in the numerator.
Since $b + 2k \geq 1$, $I$ is positive. 
If $k \geq 1$,  since each coefficient of $b^i, i=0,1,2$ in $II$ is positive, we obtain $II > 0$. If $k = 0$, we get
\[
  II := 5 (\sqrt{321}+161 ) b^2+10 b (
  -\sqrt{321}-9 )+4\cdot (-18 ) .
\]
Since $\sqrt{321} < 20$ and $b \geq 1$, we obtain 
\[
  II \geq 800 b^2 - 300 b - 100  \geq 400 >0.
\]
Thus, we prove $\Del_2>0$.

\subsubsection{Estimate of $\Del_3$}\label{sec:est_del3}
Below, we prove $\Del_3 >0$ for $ n = b + 2 k \geq 1$, and discuss four cases $(\gamma, d) \in \{ (\tf53, 3) ,  (\tf75, 3),  ( 2, 2 ), (\tf53, 2)\}$ separately. 
\paragraph{\bf Case $( \gamma, d) = (\frac53, 3) $}
We have
\[
\begin{aligned} 
\tfrac{\Del_3 }{(b+2k)^2} = \tfrac{47 \left(2209 b^4+188 b^3 (89 k-5)+b^2 \left(48416 k^2-5910 k-992\right)+4 b \left(15842 k^3-3115 k^2-943 k+69\right)+\left(-178 k^2+25 k+12\right)^2\right)}{1728}
  \end{aligned}
\]

For $k \geq 1$, since each coefficient of $b^i, i\leq 4$ in $ \frac{\Del_3}{(b+2k)^2}$ is positive, we obtain $\Del_3 > 0$. 

If $ k =0$ and $b = n \geq 1$, we obtain 
\[
\tfrac{\Del_3 }{n^2} = \tfrac{47 \left(2209 b^4-940 b^3-992 b^2+276 b+144\right)}{1728}
> \tfrac{47 \left(1000 b^3 + 1000 b^2-940 b^3-992 b^2\right)}{1728} > 0 .
\]
We prove $\Del_3>0$.

\vspace{0.1in}
\paragraph{\bf Case $( \gamma, d) = (\frac75, 3) $}
We have
\[
  \Del_3  = \frac{3 (-2 + 15 b + 30 k)^2 \cdot II}{ 262144 }
\]
where $II$ is given by 
\[
\bal 
  II & =  16875 b^4+900 b^3 (143 k-17)+6 b^2 \left(62348 k^2-15349 k-522\right) \\
 & \qquad +4 b \left(122694 k^3-46761 k^2-2767 k+843\right)+245388 k^4-128700 k^3-9493 k^2+6760 k+640   .
\eal
\]

Clearly, for $b + 2 k = n \geq 1$, we have $(-2 + 15 b + 30 k)^2>1>0$.  For $k \geq 1$, since each coefficient of $b^i, i\leq 4$ in $II$ is positive, we obtain $\Del_3 > 0$.  For $k = 0$ and $ b = n \geq 1$, we obtain 
\[
\begin{aligned} 
  II &= 640 + 3372 b - 3132 b^2 - 15300 b^3 + 16875 b^4
  > 3132(b- b^2) + 15300(b^4 - b^3) \\
  & =  (b-1) ( 15300 b^2 - 3132 b ) \geq 0.
\end{aligned}
\]
We prove $\Del_3 > 0$.

\vspace{0.1in}
\paragraph{\bf Case $( \gamma, d) = (2, 2) $}

We obtain 
\[
  \tfrac{\Del_3}{(b+2k)^2} =64 b^4+32 b^3 (15 k-1)+b^2 \left(1380 k^2-184 k-27\right)+8 b \left(225 k^3-45 k^2-13 k+1\right)+4 \left(-15 k^2+2 k+1\right)^2.
\]

For $k \geq 1$, since each coefficient of $b^i, i\leq 4$ in $ \frac{\Del_3}{(b+2k)^2}$ is positive, we obtain $\Del_3 > 0$. 

For $ k =0$ and $ b = n \geq 1$, we obtain 
\[
    \tfrac{\Del_3}{(b+2k)^2} = 4 + 8 b - 27 b^2 - 32 b^3 + 64 b^4
  > - 27 b^2 - 32 b^3  + 27 b^2 + 32 b^3 =0.
\]
We prove $\Del_3 > 0$.

\vspace{0.1in}
\paragraph{\bf Case $( \gamma, d) = (\frac53, 2) $}

We have
\[
\Del_3  =  \tfrac{\left(\sqrt{321} b+b+2 k+2 \sqrt{321} k-2\right)^2}{125000} \cdot II,  
\]
where $II$ is given by
\[
\begin{aligned} 
II  = &  { \scriptstyle
  \left(161 \sqrt{321}+13121\right) b^4 } \\ 
 &  { \scriptstyle +4 b^3 \left(8 \left(39 \sqrt{321}+3079\right) k-91 \sqrt{321}-1851\right) } \\
&  {\scriptstyle  +b^2 \left(16 \left(229 \sqrt{321}+17719\right) k^2-8 \left(263 \sqrt{321}+5343\right) k+45 \sqrt{321}-4467\right) } \\
& {\scriptstyle  +2 b \left(16 \left(151 \sqrt{321}+11561\right) k^3-48 \left(43 \sqrt{321}+873\right) k^2+10 \left(9 \sqrt{321}-875\right) k+39 \sqrt{321}+799\right) } \\ 
& {\scriptstyle +4 \left(\left(604 \sqrt{321}+46244\right) k^4-16 \left(43 \sqrt{321}+873\right) k^3+5 \left(9 \sqrt{321}-875\right) k^2+\left(39 \sqrt{321}+799\right) k+180\right)} .
\end{aligned}
\]

We rewrite $II$ as
$II = {\scriptstyle \sum_{i\leq 4}} p_i(k) b^i$,
where each $p_i$ is a polynomial in $k$, given in the corresponding row of the above formula for $II$. Since $\sqrt{321} \in (10, 20)$, for $k\geq 1$, 
clearly, we get $p_4(k)>0,p_3(k)>0$. For $k\geq 1$, we estimate $p_2, p_1, p_0$ as 
\[
\begin{aligned} 
  p_2(k) & > 16 \cdot 10000 k^2 - 8 \cdot ( 300 \cdot 20 + 6000) k - 5000
  > 160000k^2 - 100000 k - 5000 > 0, \\
p_1(k) & > 2 \cdot ( 16 \cdot 10000 k^3 - 50 (50 \cdot 20 + 1000) k^2 - 10\cdot 1000 k  )
> 2( 160000k^3 - 10^5 k^2- 10^4 k )>0, \\
p_0(k) & > 4 ( 46000 k^4- 16 (50 \cdot 20 + 1000)k^3 - 5 \cdot 1000 k^2 ) 
> 4 (46000 k^4 - 32000 k^3 - 5000 k^2 ) > 0.
  \end{aligned}
\]
Since $ b + 2k \geq 1$, the factor of $II$ is strictly positive.
Thus, for $k \geq 1$, we prove $\Del_3 >0$.

For $k = 0$ and $ b = n\geq 1$, we get
\[
 II= 
  \left(161 \sqrt{321}+13121\right) b^4-4 \left(91 \sqrt{321}+1851\right) b^3+\left(45 \sqrt{321}-4467\right) b^2+2 \left(39 \sqrt{321}+799\right) b+720 .
\]

When $b =1$, we get
 \[
   II =  3568-80 \sqrt{321} > 3000 - 80 \cdot 20 >0.
 \]
 For $ b \geq 2$, since $\sqrt{321} \in (10, 20)$, we get
 \[
   II >  13000 b^4 - 4 \cdot (100 \cdot 20 + 2000) b^3 - 5000 b^2
   > 8000b^4 - 16000 b^3 \geq 0.
 \]
 Combining the above estimates, we prove $\Del_3 >0$ for $2 k + b \geq 1$.

\subsubsection{Estimate of $a_0$}\label{sec:est_a0}

The estimate of $a_0$ is more delicate.  
In this subsection, we prove 
\beq\label{eq:a0_goal}
\bal 
a_0 > 0,  \qquad && \forall  \ n = 2k + b \geq 2,  \ (b, k ) \neq (0, 1 ),  \\ 
a_0 <0  , \qquad  && n = 1, \  (b, k) = (1, 0). 
\eal 
\eeq 
The parameter $(b,k)$ in the first case is equivalent to  $k\geq 2$  or $k=1, b \geq 1$ or $k=0, b \geq 2$.

\vspace{0.1in}
\paragraph{\bf Case $( \gamma, d) = (\frac53, 3) $}
We obtain 
\[
\bal
a_0 \cdot 2304 & = 
2209 b^4 + 376 b^3 (-5 + 42 k) + 752 k^2 (-12 - 25 k + 37 k^2)  \\
& \qquad +   4 b^2 (-600 - 3055 k + 10904 k^2) + 
 16 b (-9 - 564 k - 1645 k^2 + 3478 k^3)  . 
 \eal 
\]

For $k \geq 1$, each coefficient of $b^i, i\leq 4$  is non-negative.
Moreover, since $b^4 >0$ for $b \geq 1$ and $ 752 k^2 (-12 - 25 k + 37 k^2) > 0$ for $k\geq 2$, 
we obtain $a_0>0$ for $ k \geq 2$ or $k=1, b \geq 1$.

For $ k =0$ and $ b \geq 2$, we obtain
\[
a_0 \cdot 2304 =  b (-144 - 2400 b - 1880 b^2 + 2209 b^3)
> b ( 940 b^2(b-2) + 600 b(b^2 - 4) + 100 b^3 - 144)>0.
\]

For $k =0$ and $b = 1$, we obtain 
\[
  a_0 \cdot 2304  =   -144 - 2400  - 1880  + 2209  < 0.
\]

\vspace{0.1in}
\paragraph{\bf Case $( \gamma, d) = (\frac75, 3) $}

We obtain 
\[
\bal 
  a_0 \cdot 65536 & =  45 \Big(1125 b^4+480 b^3 (17 k-3)+4 b^2 \left(5700 k^2-2209 k-184\right) \\ 
  & \quad +16 b \left(1830 k^3-1143 k^2-157 k+11\right)+16 k \left(915 k^3-797 k^2-150 k+32\right) \Big)  
\eal 
\]

Note that $p_0(k) := 16 k \left(915 k^3-797 k^2-150 k+32\right)$ 
satisfies $p_0(1)= 16 \cdot ( 915 - 797 - 150 + 32) = 0$. For $k \geq 1$, each coefficient of $b^i, i\leq 4$ is non-negative. 
Moreover, since $b^4 >0$ for $b \geq 1$ and $ p_0(k) > 0$ for $k\geq 2$,  we obtain $a_0>0$ for $ k \geq 2$ or $k=1, b \geq 1$.

For $ k =0$ and $ b \geq 2$, we obtain
\[
    a_0 \cdot 65536 = 45 b (176 - 736 b - 1440 b^2 + 1125 b^3)
    > 45b (720(b^3 - 2b^2) + 400 b^3 - 800 b)>0.
\]

For $ k =0$ and $ b = 1$, we obtain
\[
      a_0 \cdot 65536 =  45 \cdot (176 - 736  - 1440  + 1125 ) < 0     .
\]

\vspace{0.1in}
\paragraph{\bf Case $( \gamma, d) = (2, 2) $}

We obtain 
\[
  a_0 = b^4+b^3 (7 k-1)+b^2 \left(19 k^2-6 k-\tfrac{17}{16}\right)+4 b k \left(6 k^2-3 k-1\right)+4 k^2 \left(3 k^2-2 k-1\right)
\]
Clearly, we have $a_0 \geq 0$ for $k\geq 1$ and $a_0>0$ for $k \geq 2$ or $k=1, b \geq 1$.

For $k = 0$ and $ b \geq 2$, we have
\[
  a_0 = b^4-b^3-\tfrac{17 b^2}{16}
> b^4 - b^3 - 2b^2 = b^2(b-2)(b+1) \geq 0.
\]

For $k = 0$ and $ b = 1$, we have
\[
  a_0 = -\tfrac{17 }{16} <0.
\]

\vspace{0.1in}
\paragraph{\bf Case $( \gamma, d) = (\frac53, 2) $}

We obtain 
\[
\begin{aligned} 
  20000 a_0  = & { \scriptstyle   \left(161 \sqrt{321}+13121\right) b^4 } \\ 
  & { \scriptstyle  +4 b^3 \left(\left(302 \sqrt{321}+23022\right) k-101 \sqrt{321}-3461\right) } \\
  & { \scriptstyle +4 b^2 \left(\left(866 \sqrt{321}+62626\right) k^2-2 \left(283 \sqrt{321}+10363\right) k+21 \sqrt{321}-3307\right) } \\
  & { \scriptstyle +16 b \left(2 \left(141 \sqrt{321}+9901\right) k^3-21 \left(13 \sqrt{321}+493\right) k^2+\left(21 \sqrt{321}-3019\right) k+10 \left(\sqrt{321}+1\right)\right) } \\
  &  { \scriptstyle +16 (k-1) k \left(\left(141 \sqrt{321}+9901\right) k^2+\left(2999-41 \sqrt{321}\right) k-20 \left(\sqrt{321}+1\right)\right) } \\
     :=&  { \scriptstyle   \sum_{i\leq 4}} p_i(k) b^i  ,
  \end{aligned}
\]
where $p_i$ denotes the coefficient of $b^i$ and is a polynomial in $k$, given in the corresponding row of the above formula. For $ k \geq 1$, it is easy to obtain that $p_i(k) >0$ for $2\leq i \leq 4$. For $p_1(k)$ and $k \geq 1$, 
since $\sqrt{321} \in (10, 20)$, we have 
\[
  p_1(k) >16 \cdot ( 20000 k^3 - 21 \cdot (15 \cdot 20 + 500) k^2 - 3000 k )
 \geq 16 \cdot ( 17000 k^3 - 16800 k^2 ) >0.
\]
It is easy to obtain $p_0(0) = p_0(1) =0$ and $p_0(k) >0$ for $k \geq 2$.
Since $b^4 >0$ for $b \geq 1$, combining these estimates, we obtain $a_0>0$ for $k \geq 2$ or $k=1, b \geq 1$.

For $ k = 0$ and $ b \geq 2$, we obtain 
\[
    20000 a_0  = b \left(\left(161 \sqrt{321}+13121\right) b^3-4 \left(101 \sqrt{321}+3461\right) b^2+4 \left(21 \sqrt{321}-3307\right) b+160 \left(\sqrt{321}+1\right)\right).
\]

Since $\sqrt{321} \in (15, 20)$ and $50 \sqrt{321} < 1000$, for $b\geq 2$, we obtain 
\[
      20000 a_0 >  b \Big( (210 \sqrt{321} + 7000 ) b^3   + 4000 b^3 
      - (410 \sqrt{321} + 14000) b^2 -4 \cdot 4000 b  \Big)
> 4000 b (b^3- 4 b) \geq 0.
\]

For $k=0$ and $ b=1$, we obtain 
\[
\bal 
20000 a_0  &=
  (161 \sqrt{321}+13121 ) -4 (101 \sqrt{321}+3461 ) +4 (21 \sqrt{321}-3307 ) +160 (\sqrt{321}+1 )  \\
  & = \sqrt{321} 
  + 13121 - 4 \times 3461 - 4 \times 3307 + 160 < 0.
\eal 
\]

Combining the above cases, we prove \eqref{eq:a0_goal}: $a_0 > 0$ for $n = 2k + b \geq 2$ with $(k, b) \neq (1, 0)$
and $a_0<0$ for $n=1$.

Combining the estimates of $a_3, a_0, \Del_2, \Del_3$ in Sections \ref{sec:est_a3}-\ref{sec:est_a0}, we prove estimates \eqref{eq:eig_Hbk}.

\paragraph{\bf Summary}

Using Theorem \ref{thm:RH}, we prove the statement of eigenvalues of $\HH_{|\bet|, k}$ in Lemma \ref{lem:Hbk}. 

Since $\HH_{0, k}$ with $k\geq 2$ have $4$ eigenvalues with negative real part, using Lemma \ref{lem:HH_connect}, we prove that all eigenvalues of $\HHH_{k}$ with $k\geq 2$ have negative real parts. We prove Part (I) and (II) in  Lemma \ref{lem:Hbk}.

\subsection{Estimate of eigenvalues of \texorpdfstring{$\HHT_{n}$}{H lower n upper 2}}

For any $\gamma > 1, d \geq 1$, a direct calculation yields 
\[
f_3 := \tr(\HHT_n) = n (\frac{1}{ \al d+1}- \bcr )+\frac{- \al ( \bcr d n+d)- \bcr n+n+1}{\al d+1} .
\]

Using Lemma \ref{lem:cr_low} and $\bcr > 1$, we obtain that $f_3$ is decreasing in $n$. Since $n \geq 1$, we get
\[
  f_3 \leq f_3 |_{n=1} = \frac{ 3 - \al d - 2 \bcr (1 + \al d) }{1 + \al d} .
\]

Using the lower bound of $\bcr$ in Lemma \ref{lem:cr_low}, we obtain 
\[
   3 - \al d - 2 \bcr (1 + \al d)  < 3 - \al d - 2(1 + \al d) - (1 - \al d)_+
   = 1- \al d - (1 - \al d)_+ - 2 \al d  < 0.
\]

Thus, we prove $f_3 < 0$ and prove \eqref{eq:eig_HHT_tr}.

Next, we estimate $f_2 =   \det( \HHT_{n} )$. 
A direct calculation yields 
\[
\begin{aligned} 
f_2(b,k) & = \tfrac{1}{48} \left(47 b^2+4 b (47 k-3)+4 \left(47 k^2-6 k-6\right)\right) , 
\quad  &&  (\gamma, d) = (\tfrac53, 3), \\
f_2(b,k) & =  \tfrac{15}{256} \left(15 b^2+b (60 k-8)+60 k^2-16 k-8\right), 
 \quad && (\gamma, d) = (\tfrac75, 3), \\
 f_2(b,k) & =  b^2+b (4 k-\tfrac{1}{4} )+4 k^2-\tfrac{k}{2}-\tfrac{1}{2} , 
 \quad && (\gamma, d) = (2, 2), \\
 \end{aligned}
 \]
 and 
 \[
 f_2(b,k)  = 
 {\scriptstyle 
  \frac{1}{200} \left(\sqrt{321}+161\right) b^2+\frac{1}{200} b \left(4 \left(\sqrt{321}+161\right) k-2 \left(\sqrt{321}+25\right)\right)+\frac{1}{200} \left(4 \left(\sqrt{321}+161\right) k^2-4 \left(\sqrt{321}+25\right) k-96\right)   }, 
\]
for $    (\gamma, d) = ( \tfrac53, 2) $.

For $k \geq 1$, each coefficient of $b^i, i\leq 2$  is positive. Thus, $f_2 >0$ for $k \geq 1$.

For $ k =0$, we obtain 
\[
\begin{aligned} 
  f_2(b,0)  & = \tfrac{1}{48} \left(47 b^2-12 b-24\right) ,\quad  &&  (\gamma, d) = (\tfrac53, 3) , \\
  f_2(b, 0)  & =  \tfrac{15}{256} \left(15 b^2-8 b-8\right), 
   \quad && (\gamma, d) = (\tfrac75, 3) , \\ 
  f_2(b, 0) &  = b^2-\tfrac{b}{4}-\tfrac{1}{2} ,
\quad && (\gamma, d) = (2, 2) . \\
\end{aligned}
\]

Clearly, for $ (\gamma, d ) =  (\tfrac53, 3),  (2, 2)$ and $ b \geq 1$,  we obtain $f_2(b, 0) > 0$. 
 We prove \eqref{eq:eig_HHT_det:a}. 
 
For $ (\gamma, d )  = (\tfrac75, 3)$, we obtain $f_2(b, 0) > 0$ for $ b \geq 2$, and 
\[
f_2(1,0) = - \tfrac{15}{256}  < 0.
\]

For $(\gamma, d)=(\frac53, 2)$ and $b \geq 2$, we obtain 
\[
   f_2(b, 0)
   = \frac{1}{200} (\sqrt{321}+161 ) b^2-\frac{1}{100} (\sqrt{321}+25 ) b-\frac{12}{25}
   \geq \frac{161 }{100} b - \frac{25 }{100} b - \frac{12}{25} > 0.
\]

For $(\gamma, d)=(\frac53, 2)$ and $b =1$, we obtain 
\[
  f_2(1,0) = \tfrac{1}{200} (15-\sqrt{321} )  < 0.
\]

Combining the above estimates, we prove \eqref{eq:eig_HHT_det:b}, \eqref{eq:eig_HHT_det:c}. 
The signs of the real part of the eigenvalues follow from \eqref{eq:eig_HHT_tr} and \eqref{eq:eig_HHT_det}. This completes the proof of Part (III) in Lemma~\ref{lem:Hbk}.

\paragraph{\bf Proof of Part (IV) }

We estimate $f_4 := - n \kp - 2 \uua -1$. 
 For $n \geq 1$, a direct estimate yields 
\[
f_4  =\frac{- \al  ( \bcr d n+d)- \bcr n+n+1 }{ \al d+1}
\leq \frac{- \al  ( \bcr d +d)- \bcr + 2 }{ \al d+1}.
\]
Using Lemma \ref{lem:cr_low}, we prove 
\[
- \al  ( \bcr d +d)- \bcr + 2
\leq - 1 - 2 \al d - \frac{(1-\al d)_+}{2} + 2
\leq 1  - 2 \al d- \frac{(1-\al d)_+}{2} .
\]

For $\al d \geq \frac12$, the upper bound is strictly negative and we prove $f_4<0$.
In particular, for $(\gamma, d) \in \{ (\frac53, 3), (\frac75, 3), (2,2),(\frac53, 2)\}$ with $\al = \frac{\gamma-1}{2}$, we obtain $\al d > \frac{1}{2}$ and prove $f_4<0$. We prove Part (IV) in Lemma \ref{lem:Hbk}.

\subsection{Proof of result (V) in Lemma \ref{lem:Hbk}}

From result (I) and (III) in Lemma \ref{lem:Hbk}, $\HH_{1,0}$ has exactly one eigenvalue with positive real part;
for $(\gamma, d)  \in \{ (\frac53, 2) , ( \tf75, 3) \}$, $\HHT_1$ has exactly one eigenvalue with positive real  part.
Since these eigenvalues have multiplicity one, they are real. 
We prove result (V) in Lemma \ref{lem:Hbk}.

\subsection{Estimates of \texorpdfstring{$\HH_{n}^{(1)}$}{H n 1}}

Recall the matrix $\HH_n^{(1)}$ from \eqref{eq:ODE_Nth:1D}.  We write 
\[
    p_l(\lam) = \det(\lam \Id - \HH_n^{(1)}) = c_l  \lam^l + c_{l-1} \lam^{l-1} + ... + c_1 \lam + c_0,
  \quad c_l = 1, \quad  l = 3.
\]
The associated Routh table is given by
\begin{equation}\label{eq:RH_H1d}
  A = \left(
  \begin{array}{ccc }
  c_3  & c_1 & 0   \\
  c_2 & c_0 & 0   \\
 \frac{\Del_4}{c_2} & c_0 & 0  \\
  c_0 & 0 & 0  \\
\end{array}
  \right) , 
  \quad 
A_{31}= c_1 - c_0 \frac{c_3}{c_2} = \frac{ \Del_4 }{c_2},  \quad 
\Del_4 =  c_2 c_1 - c_0 c_3, \quad  c_3 =1.
\end{equation}

To apply Theorem \ref{thm:RH} to $\HH_n^{(1)}$,  we estimate the signs of $
c_2, \Del_4, c_0$, which determine the signs of $A_{j 1} $.

\begin{lemma}\label{lem:H1d}
Consider $d = 1$ and $n \geq 1$. For any $\gamma > 1$, we have 
\begin{subequations}\label{eq:eig_H1d}
\begin{equation}
  c_2 > 0 , \quad  \Del_4 > 0.
\end{equation}
For $c_0$ and any $\gamma>1$, we have 
\begin{equation}
  c_0 > 0,  \ \forall n \geq 3, \quad c_0 = 0, \ n = 2,  \quad c_0 < 0, \ n = 1.
\end{equation}
\end{subequations}
As a result, $\HH^{(1)}_n$ has $3$ eigenvalues with negative real part for $n \geq 3$;
$2$ eigenvalues with negative real part and a $0$-eigenvalue for $n =2$;
$2$ eigenvalues with negative real part and $1$ eigenvalue with positive real part for $n = 1$.
\end{lemma}

\begin{proof}[Proof of Lemma~\ref{lem:H1d}] 
We denote $w = \al d = \al$. 

\paragraph{\bf Estimate of $c_2$}

A direct calculation yields 
\[
  c_2 = \frac{3 n (\sqrt{49 w^2+26 w+1}+1 )+(4-3 n) w-4}{4 (w+1)}.
\]
Since $c_2$ is increasing in $n$, $n \geq 1$, and $w>0$, we obtain 
\[
  c_2 \geq c_2 |_{n=1} = \frac{3 \sqrt{49 w^2+26 w+1}+w-1}{4 (w+1)} >0.
\]

\vspace{0.1in}
\paragraph{\bf Estimate of $\Del_4$}

We show that $\Del_4$ is monotone in $n$ for $n \geq 1$. Firstly, a direct calculation yields 
\[
\begin{aligned} 
I:=4(1+w)^3 \Del_4 =  &  { \scriptstyle  n^3  \left(-72 w^3+\left(24 \sqrt{49 w^2+26 w+1}+35\right) w^2+\left(9 \sqrt{49 w^2+26 w+1}+35\right) w+2 \left(\sqrt{49 w^2+26 w+1}+1\right)\right) } \\
 & { \scriptstyle  +n^2    \left(102 w^3-\left(10 \sqrt{49 w^2+26 w+1}+53\right) w^2+5 \left(\sqrt{49 w^2+26 w+1}-9\right) w-4 \left(\sqrt{49 w^2+26 w+1}+1\right)\right) }\\
&  { \scriptstyle +n    \left(-10 w^3-2 \left(\sqrt{49 w^2+26 w+1}-6\right) w^2-4 \left(2 \sqrt{49 w^2+26 w+1}+1\right) w+2 \left(\sqrt{49 w^2+26 w+1}+1\right)\right)-8 w^3+8 w } ,
\end{aligned}
\]
which is a cubic polynomial in $n$. For the $n^3$-term, since 
$ (24 \sqrt{49 w^2+26 w+1}+35 ) w^2> 24 \cdot 7 w \cdot w^2> 72 w^3$, the coefficient of $n^3$ is positive. Thus, $\pa_n^3 I >0$  and $\pa_n^2 I$ is increasing in $n$ for any $n >0$.

Using a direct computation, for $n \geq 1$, we obtain 
\[
\bal 
  & \pa_n^2 I(n) \geq  (\pa_n^2 I)(1)  \\
  &=  -57 w^3+\left(31 \sqrt{49 w^2+26 w+1}+26\right) w^2+2 \left(8 \sqrt{49 w^2+26 w+1}+15\right) w+\sqrt{49 w^2+26 w+1}+1 .
\eal 
\]
Since $(31 \sqrt{49 w^2+26 w+1}+26 ) w^2 > 31 \cdot 7 w^3 > 57 w^3$, we get $  \pa_n^2 I(n) >0$ 
and $\pa_n I(n)$ is increasing in $n$ for $n \geq 1$.

Using a direct computation, for $n \geq 1$, we have 
\[
    \pa_n I(n) \geq  (\pa_n I)(1)  =  \tfrac{1}{4} w \left(-22 w^2+\left(50 \sqrt{49 w^2+26 w+1}+11\right) w+29 \sqrt{49 w^2+26 w+1}+11\right).
\]
Since $50 \sqrt{49 w^2} w > 22 w^2$, we prove $\pa_n I(n) > 0$ and $I(n)$ is increasing in $n$ for $n \geq 1$. 

Using a direct computation, for $n \geq 1$, we deduce 
\[
      I(n) \geq  I(1)  = \tfrac{3}{2} w (2 w+1) (\sqrt{49 w^2+26 w+1}+w-1 ) > 0.
\]
We prove $\Del_4 >0$ for $n \geq 1$.  

\vspace{0.1in}
\paragraph{\bf Estimate of $c_0$}

Recall $d=1$ and $w = \al d = \al$.  A direct computation yields
\[
\begin{aligned} 
 II & :=\frac{32(1+w)^3}{n} c_0   = 
 (n-2) \left(\sqrt{49 w^2+26 w+1}-w+1\right) \\
   & \qquad  \times \left((17 n+8) w^2+w \left(8-n \left(\sqrt{49 w^2+26 w+1}-8\right)\right)+n \left(\sqrt{49 w^2+26 w+1}+1\right)\right) \\
   & :=(n-2) J_1 J_2.
   \end{aligned}
\]
Clearly, the first factor $J_1>0$. Below, we show that $J_2 >0$. We can rewrite $J_2$ as 
\[
  J_2 = n \left(17 w^2-\left(\sqrt{49 w^2+26 w+1}-8\right) w+\sqrt{49 w^2+26 w+1}+1\right)+8 w^2+8 w.
\]
Since $\sqrt{49 w^2+26 w+1} < 7 w + 2$, we obtain 
\[
  17 w^2- (\sqrt{49 w^2+26 w+1}-8 ) w+\sqrt{49 w^2+26 w+1}+1
  > 17 w^2 + 8 w + 1 - (7w + 2) w > 0.
\]
We prove $J_2>0$. Thus, we obtain $\sgn( c_0) = \sgn(n-2)$. Combining the above estimates, we prove \eqref{eq:eig_H1d}. 

For $n \geq 3$ or $n=1$, the signs of the eigenvalues follow from \eqref{eq:eig_H1d} and Theorem \ref{thm:RH}. 
In particular, $\HH_n^{(1)}$ has three eigenvalues with negative real parts for $n \geq 3$;
 it has two eigenvalues with negative real parts and one with a positive real part for $n=1$.

For $n=2$, since $c_0 = 0$, we obtain $p_3(\lam) = \lam (c_3 \lam^2 + c_2 \lam + c_1)$. 
From \eqref{eq:RH_H1d} and \eqref{eq:eig_H1d}, 
since $c_0 =0$, we have $c_3 = 1, c_2 >0, \Del_4 = c_2 c_1 > 0 $, which implies $c_1, c_2, c_3>0$. 
Thus, the quadratic polynomial $c_3 \lam^2 + c_2 \lam + c_1$ has two roots with negative real part.
We complete the proof of Lemma \ref{lem:H1d}.
\end{proof}

%\bibliographystyle{plain}
%\bibliography{implosion_refs}

\begin{thebibliography}{10}

\bibitem{Barenblatt1996}
G.~I. Barenblatt.
\newblock {\em Scaling, Self-similarity, and Intermediate Asymptotics},
  volume~14 of {\em Cambridge Texts in Applied Mathematics}.
\newblock Cambridge University Press, 1996.

\bibitem{bedrossian2026finite}
Jacob Bedrossian, Jiajie Chen, Maria~Pia Gualdani, Sehyun Ji, Vlad Vicol, and
  Jincheng Yang.
\newblock Finite time singularities in the Landau equation with very hard
  potentials.
\newblock {\em arXiv preprint arXiv:2602.05981}, 2026.

\bibitem{Biasi2021}
Anxo Biasi.
\newblock Self-similar solutions to the compressible {E}uler equations and
  their instabilities.
\newblock {\em Commun. Nonlinear Sci. Numer. Simul.}, 103:Paper No. 106014,
  2021.

\bibitem{bodson2019explaining}
Marc Bodson.
\newblock Explaining the routh-hurwitz criterion.
\newblock {\em A tutorial presentation (Ingl{\'e}s)[Explicaci{\'o}n del
  criterio de Routh-Hurwitz. Un tutorial de presentaci{\'o}n]}, 15, 2019.

\bibitem{braun1983differential}
Martin Braun.
\newblock {\em Differential equations and their applications}.
\newblock Springer New York, NY, 3rd edition, 1983.

\bibitem{BrushlinskiiKazhdan1963}
K.~V. Brushlinskii and Ya.~M. Kazhdan.
\newblock On self-similar solutions of certain problems of gas dynamics.
\newblock {\em Uspekhi Mat. Nauk}, 18(2(110)):3--23, 1963.

\bibitem{BCG2025}
Tristan Buckmaster, Gonzalo Cao-Labora, and Javier G{\'o}mez-Serrano.
\newblock Smooth imploding solutions for {3D} compressible fluids.
\newblock {\em Forum Math. Pi}, 13:e6, 2025.

\bibitem{BuckChen2024}
Tristan Buckmaster and Jiajie Chen.
\newblock Blowup for the defocusing septic complex-valued nonlinear wave
  equation in {$\mathbb{R}^{4+1}$}.
\newblock {\em To appear in Commun.~Amer.~Math.~Soc.; arXiv:2410.15619}, 2024.

\bibitem{BDSV}
Tristan Buckmaster, Theodore~D. Drivas, Steve Shkoller, and Vlad Vicol.
\newblock Simultaneous development of shocks and cusps for 2d euler with
  azimuthal symmetry from smooth data.
\newblock {\em Annals of PDE}, 8(2):26, 2022.

\bibitem{caflisch1980fluid}
Russel~E. Caflisch.
\newblock The fluid dynamic limit of the nonlinear {B}oltzmann equation.
\newblock {\em Comm. Pure Appl. Math.}, 33(5):651--666, 1980.

\bibitem{CGSS_NLS2024}
Gonzalo Cao-Labora, Javier G{\'o}mez-Serrano, Jia Shi, and Gigliola Staffilani.
\newblock Non-radial implosion for the defocusing nonlinear {S}chr{\"o}dinger
  equation in {$\mathbb{T}^d$} and {$\mathbb{R}^d$}.
\newblock {\em arXiv preprint arXiv:2410.04532}, 2024.

\bibitem{CGSS2025}
Gonzalo Cao-Labora, Javier G{\'o}mez-Serrano, Jia Shi, and Gigliola Staffilani.
\newblock Non-radial implosion for compressible {E}uler and {N}avier-{S}tokes
  in {$\mathbb{T}^3$} and {$\mathbb{R}^3$}.
\newblock {\em Cambridge J. Math.}, 2025.
\newblock to appear.

\bibitem{ChenLiuZhu2025}
Gui-Qiang~G. Chen, Lihui Liu, and Shengguo Zhu.
\newblock Development of implosions of solutions to the three-dimensional
  degenerate compressible {N}avier--{S}tokes equations.
\newblock Preprint, \texttt{arXiv:2603.10141}, 2025.
\newblock Preprint.

\bibitem{ChenZhangZhu2025}
Gui-Qiang~G. Chen, Jiawen Zhang, and Shengguo Zhu.
\newblock Global regular solutions of the multidimensional degenerate
  compressible {N}avier--{S}tokes equations with large initial data of
  spherical symmetry.
\newblock Preprint, \texttt{arXiv:2512.18545}, 2025.
\newblock Preprint.

\bibitem{ChenZhangPanarella1995}
H.~B. Chen, L.~Zhang, and E.~Panarella.
\newblock Stability of imploding spherical shock waves.
\newblock {\em J. Fusion Energy}, 14:389--392, 1995.

\bibitem{chen2020singularity}
Jiajie Chen.
\newblock Singularity formation and global well-posedness for the generalized
  {Constantin--Lax--Majda} equation with dissipation.
\newblock {\em Nonlinearity}, 33(5):2502, 2020.

\bibitem{Chen2024}
Jiajie Chen.
\newblock Vorticity blowup in compressible euler equations in {$\mathbb{R}^d$},
  {$d \geq 3$}.
\newblock {\em Annals of PDE}, 11(2):21, 2025.

\bibitem{CCSV2024}
Jiajie Chen, Giorgio Cialdea, Steve Shkoller, and Vlad Vicol.
\newblock Vorticity blowup in {2D} compressible {E}uler equations.
\newblock {\em To appear in Duke Math. J., arXiv preprint arXiv:2407.06455},
  2024.

\bibitem{ChenHou2023a}
Jiajie Chen and Thomas~Y Hou.
\newblock Stable nearly self-similar blowup of the 2{D} {B}oussinesq and 3{D}
  {E}uler equations with smooth data {I}: {A}nalysis.
\newblock {\em arXiv preprint: arXiv:2210.07191v3}, 2022.

\bibitem{chen2025singularity}
Jiajie Chen and Thomas~Y Hou.
\newblock Singularity formation in 3d euler equations with smooth initial data
  and boundary.
\newblock {\em Proceedings of the National Academy of Sciences},
  122(27):e2500940122, 2025.

\bibitem{chen2019finite}
Jiajie Chen, Thomas~Y Hou, and De~Huang.
\newblock On the finite time blowup of the {D}e {G}regorio model for the 3{D}
  {E}uler equations.
\newblock {\em Communications on Pure and Applied Mathematics},
  74(6):1282--1350, 2021.

\bibitem{chen2021HL}
Jiajie Chen, Thomas~Y Hou, and De~Huang.
\newblock Asymptotically self-similar blowup of the {H}ou--{L}uo model for the
  3{D} {E}uler equations.
\newblock {\em Annals of PDE}, 8(2):24, 2022.

\bibitem{chen2025stability}
Jiajie Chen, Thomas~Y Hou, Van~Tien Nguyen, and Yixuan Wang.
\newblock On the stability of blowup solutions to the complex Ginzburg-Landau
  equation in.
\newblock {\em Annals of PDE}, 11(2):29, 2025.

\bibitem{CialdShkVic2025}
Giorgio Cialdea, Steve Shkoller, and Vlad Vicol.
\newblock Classical {E}uler flows generate the strong {G}uderley imploding
  shock wave.
\newblock {\em arXiv preprint arXiv:2510.19688}, 2025.

\bibitem{dafermos2005hyberbolic}
Constantine~M Dafermos.
\newblock {\em Hyberbolic conservation laws in continuum physics}.
\newblock Springer, 2005.

\bibitem{EggersFontelos2015}
Jens Eggers and Marco~A. Fontelos.
\newblock {\em Singularities: Formation, Structure, and Propagation}, volume~53
  of {\em Cambridge Texts in Applied Mathematics}.
\newblock Cambridge University Press, Cambridge, 2015.

\bibitem{gantmacher2005applications}
Feliks~Rouminovich Gantmacher and Joel~Lee Brenner.
\newblock {\em Applications of the Theory of Matrices}.
\newblock Courier Corporation, 2005.

\bibitem{GardnerBrookBernstein1982}
J.~H. Gardner, D.~L. Book, and I.~B. Bernstein.
\newblock Stability of imploding shocks in the {CCW} approximation.
\newblock {\em J. Fluid Mech.}, 114:41--58, 1982.

\bibitem{Giron2023}
Itamar Giron, Shmuel Balberg, and Menahem Krief.
\newblock Solutions of the converging and diverging shock problem in a medium
  with varying density.
\newblock {\em Physics of Fluids}, 35(6), 2023.

\bibitem{GoldingHenderson25}
William Golding and Christopher Henderson.
\newblock On hydrodynamic implosions and the {L}andau-{C}oulomb equation.
\newblock {\em arXiv preprint arXiv:2511.03033}, 2025.

\bibitem{Guderley1942}
K.~G. Guderley.
\newblock {Starke kugelige und zylindrische Verdichtungsst{\"o}sse in der
  N{\"a}he des Kugelmittelpunktes bzw. der Zylinderachse}.
\newblock {\em Luftfahrtforschung}, 19:302--312, 1942.

\bibitem{GuillenSilvestreLandau}
Nestor Guillen and Luis Silvestre.
\newblock The {L}andau equation does not blow up.
\newblock {\em Acta Mathematica}, 234(2):315--375, 2025.

\bibitem{guo2021larson}
Yan Guo, Mahir Had{\v{z}}i{\'c}, and Juhi Jang.
\newblock Larson--Penston self-similar gravitational collapse.
\newblock {\em Communications in mathematical physics}, 386(3):1551--1601,
  2021.

\bibitem{guo2025nonlinear}
Yan Guo, Mahir Hadzic, Juhi Jang, and Matthew Schrecker.
\newblock Nonlinear stability of the Larson--Penston collapse.
\newblock {\em arXiv preprint arXiv:2509.12435}, 2025.

\bibitem{hanawa1997stability}
Tomoyuki Hanawa and Kunji Nakayama.
\newblock Stability of similarity solutions for a gravitationally contracting
  isothermal sphere: Convergence to the Larson--Penston solution.
\newblock {\em The Astrophysical Journal}, 484(1):238--244, 1997.

\bibitem{imbert2026monotonicity}
Cyril Imbert, Luis Silvestre, and C{\'e}dric Villani.
\newblock On the monotonicity of the fisher information for the Boltzmann
  equation.
\newblock {\em Inventiones mathematicae}, 243(1):127--179, 2026.

\bibitem{JangLiuSchrecker2025}
Juhi Jang, Jiaqi Liu, and Matthew~Schrecker.
\newblock On self-similar converging shock waves.
\newblock {\em Arch. Ration. Mech. Anal.}, 249:Paper No. 37, 2025.

\bibitem{Jenssen2025}
Helge~Kristian Jenssen.
\newblock Amplitude blowup in compressible {E}uler flows without shock
  formation.
\newblock {\em arXiv preprint arXiv:2501.09037}, 2025.

\bibitem{JenssenTsikkou2018}
Helge~Kristian Jenssen and Charis Tsikkou.
\newblock On similarity flows for the compressible {E}uler system.
\newblock {\em J. Math. Phys.}, 59(12):121507, 2018.

\bibitem{jenssen2023radially}
Helge~Kristian Jenssen and Charis Tsikkou.
\newblock Radially symmetric non-isentropic euler flows: Continuous blowup with
  positive pressure.
\newblock {\em Physics of Fluids}, 35(1), 2023.

\bibitem{LandauLifshitz1987}
L.D.~Landau and E.M.~Lifshitz.
\newblock {\em Fluid Mechanics}, volume~6 of {\em Course of Theoretical
  Physics}.
\newblock Pergamon Press, Oxford, 2 edition, 1987.
\newblock Translated from the Russian by J. B. Sykes and W. H. Reid.

\bibitem{larson1969numerical}
Richard~B. Larson.
\newblock Numerical calculations of the dynamics of a collapsing proto-star.
\newblock {\em Monthly Notices of the Royal Astronomical Society},
  145(3):271--295, 1969.

\bibitem{Lax1954}
Peter~D. Lax.
\newblock Weak solutions of nonlinear hyperbolic equations and their numerical
  computation.
\newblock {\em Comm. Pure Appl. Math.}, 7:159--193, 1954.

\bibitem{Lazarus1977}
R.B.~Lazarus and R.D.~Richtmyer.
\newblock Similarity solutions for converging shocks.
\newblock Technical report, Los Alamos Scientific Lab., NM (USA), 1977.

\bibitem{Lazarus1981}
R.B.~Lazarus.
\newblock Self-similar solutions for converging shocks and collapsing cavities.
\newblock {\em SIAM J. Numer. Anal.}, 18(2):316--371, 1981.

\bibitem{meinsma1995elementary}
Gjerrit Meinsma.
\newblock Elementary proof of the Routh-Hurwitz test.
\newblock {\em Systems \& Control Letters}, 25(4):237--242, 1995.

\bibitem{MRRS2022}
Frank Merle, Pierre Rapha{\"e}l, Igor Rodnianski, and Jeremie Szeftel.
\newblock On blow up for the energy super critical defocusing nonlinear
  Schr\"odinger equations.
\newblock {\em Invent. Math.}, 227(1):247--413, 2022.

\bibitem{MRRS2022a}
Frank Merle, Pierre Rapha{\"e}l, Igor Rodnianski, and Jeremie Szeftel.
\newblock On the implosion of a compressible fluid {I}: {S}mooth self-similar
  inviscid profiles.
\newblock {\em Ann. of Math. (2)}, 196(2):567--778, 2022.

\bibitem{MRRS2022b}
Frank Merle, Pierre Rapha{\"e}l, Igor Rodnianski, and Jeremie Szeftel.
\newblock On the implosion of a compressible fluid {II}: {S}ingularity
  formation.
\newblock {\em Ann. of Math. (2)}, 196(2):779--889, 2022.

\bibitem{MeyerTerVehn1982}
J.~Meyer-ter Vehn and C.~Schalk.
\newblock Selfsimilar spherical compression waves in gas dynamics.
\newblock {\em Z. Naturforsch. A}, 37:955--969, 1982.

\bibitem{Morawetz1951}
Cathleen~Synge Morawetz.
\newblock {\em Contracting Spherical Shocks Treated by a Perturbation Method}.
\newblock PhD thesis, New York University, 1951.
\newblock Advised by K. O. Friedrichs.

\bibitem{NRSV24}
Isaac Neal, Calum Rickard, Steve Shkoller, and Vlad Vicol.
\newblock A new type of stable shock formation in gas dynamics.
\newblock {\em Communications on Pure and Applied Analysis}, 23(10):1423--1447,
  2024.

\bibitem{NealShkollerVicol25}
Isaac Neal, Steve Shkoller, and Vlad Vicol.
\newblock A characteristics approach to shock formation in 2d euler with
  azimuthal symmetry and entropy.
\newblock {\em Communications in Analysis and Mechanics}, 17(1):188--236, 2025.

\bibitem{NSV25}
Isaac Neal, Steve Shkoller, and Vlad Vicol.
\newblock Gradient catastrophes and an infinite hierarchy of h{\"o}lder
  cusp-singularities for 1d euler.
\newblock {\em Journal of the London Mathematical Society}, 112(2):e70261,
  2025.

\bibitem{ori1988simple}
Amos Ori and Tsvi Piran.
\newblock A simple stability criterion for isothermal spherical self-similar
  flow.
\newblock {\em Monthly Notices of the Royal Astronomical Society},
  234(4):821--829, 1988.

\bibitem{penston1969dynamics}
M.V.~Penston.
\newblock Dynamics of self-gravitating gaseous spheres---iii: Analytical
  results in the free-fall of isothermal cases.
\newblock {\em Monthly Notices of the Royal Astronomical Society},
  144(4):425--448, 1969.

\bibitem{Ramsey2012}
Scott~D. Ramsey, James~R. Kamm, and John~H. Bolstad.
\newblock The {G}uderley problem revisited.
\newblock {\em Int. J. Comput. Fluid Dyn.}, 26(2):79--99, 2012.

\bibitem{Sedov1946}
Leonid~Ivanovich Sedov.
\newblock Propagation of strong shock waves.
\newblock {\em Journal of Applied Mathematics and Mechanics}, 10:241--250,
  1946.

\bibitem{Sedov2018}
Leonid~Ivanovich Sedov.
\newblock {\em Similarity and dimensional methods in mechanics}.
\newblock CRC press, 10th edition, 2018.

\bibitem{ShaoWangWeiZhang2025}
Feng Shao, Shumao Wang, Dongyi Wei, and Zhifei Zhang.
\newblock Blow-up of the {3-D} compressible {N}avier-{S}tokes equations for
  monatomic gases.
\newblock {\em arXiv preprint arXiv:2501.15701}, 2025.

\bibitem{SWZ2024}
Feng Shao, Dongyi Wei, and Zhifei Zhang.
\newblock Self-similar imploding solutions of the relativistic {E}uler
  equations.
\newblock {\em arXiv preprint arXiv:2403.11471}, 2024.

\bibitem{shao2025blowup}
Feng Shao, Dongyi Wei, and Zhifei Zhang.
\newblock On blow-up for the supercritical defocusing nonlinear wave equation.
\newblock In {\em Forum of Mathematics, Pi}, volume~13, page e15. Cambridge
  University Press, 2025.

\bibitem{ShkollerVicol2024}
Steve Shkoller and Vlad Vicol.
\newblock The geometry of maximal development and shock formation for the
  {E}uler equations in multiple space dimensions.
\newblock {\em Invent. Math.}, 237:871--1252, 2024.

\bibitem{Stanyukovich1960}
K.~P. Stanyukovich.
\newblock {\em Unsteady Motion of Continuous Media}.
\newblock Pergamon Press, Oxford, 1960.
\newblock Translated from the Russian edition: Gostekhizdat, Moscow, 1955.

\bibitem{Taylor1950}
Geoffrey~Ingram Taylor.
\newblock The formation of a blast wave by a very intense explosion i.
  theoretical discussion.
\newblock {\em Proceedings of the Royal Society of London. Series A.
  Mathematical and Physical Sciences}, 201(1065):159--174, 1950.

\bibitem{teschl2012ordinary}
Gerald Teschl.
\newblock {\em Ordinary differential equations and dynamical systems}, volume
  140.
\newblock American Mathematical Soc., 2012.

\bibitem{vonNeumann1947}
John von Neumann.
\newblock The point source solution.
\newblock In Hans~A. Bethe, Klaus Fuchs, Joseph~O. Hirschfelder, John~L. Magee,
  Rudolph~E. Peierls, and John von Neumann, editors, {\em Blast Wave}, pages
  27--55. Los Alamos Scientific Laboratory Report, 1947.
\newblock Los Alamos Scientific Laboratory Report LA-2000, written August 1947,
  distributed March 27, 1958.

\bibitem{WuRoberts1996PRE}
C.~C. Wu and P.~H. Roberts.
\newblock Instability of converging shock waves and sonoluminescence.
\newblock {\em Phys. Rev. E}, 54(5):5004--5011, 1996.

\bibitem{WuRoberts1996}
C.~C. Wu and P.~H. Roberts.
\newblock Structure and stability of a spherical shock wave in a van der
  {W}aals gas.
\newblock {\em Quart. J. Mech. Appl. Math.}, 49(4):501--543, 1996.

\bibitem{ZelDovich2002}
Ya~B Zel'Dovich and Yu~P Raizer.
\newblock {\em Physics of shock waves and high-temperature hydrodynamic
  phenomena}.
\newblock Courier Corporation, 2002.

\end{thebibliography}

\end{document}